\documentclass[reqno,11pt,letterpaper]{amsart}

\usepackage{lipsum}
\usepackage{amsmath}
\usepackage{amssymb}
\usepackage{amsthm}
\usepackage{mathrsfs}
\usepackage{accents}
\usepackage{calc}
\usepackage{arydshln}
\usepackage{upgreek}
\usepackage{slashed}
\usepackage{xifthen}
\usepackage{graphicx}
\usepackage{longtable}
\usepackage[inline]{enumitem}

\usepackage{xcolor}
\definecolor{winered}{rgb}{0.6,0,0}
\definecolor{lessblue}{rgb}{0,0,0.7}

\usepackage[pdftex,colorlinks=true,linkcolor=winered,citecolor=lessblue,urlcolor=lessblue,breaklinks=true,bookmarksopen=true]{hyperref}

\makeatletter
\newcommand{\myitem}[3]{\item[#2]\def\@currentlabel{#3}\label{#1}}
\makeatother

\usepackage{titletoc}

\makeatletter

\def\@tocline#1#2#3#4#5#6#7{
\begingroup
  \par
    \parindent\z@ \leftskip#3 \relax \advance\leftskip\@tempdima\relax
                  \rightskip\@pnumwidth plus 4em \parfillskip-\@pnumwidth
    \ifcase #1 
       \vskip 0.6em \hskip 0em 
       \or
       \or \hskip 0em 
       \or \hskip 1em 
    \fi%
    #6
    \nobreak\relax{\leavevmode\leaders\hbox{\,.}\hfill}
    \hbox to\@pnumwidth {\@tocpagenum{#7}}
  \par
\endgroup
}

 \def\l@section{\@tocline{0}{0pt}{0pc}{}{}}

\renewcommand{\tocsection}[3]{%
  \indentlabel{\@ifnotempty{#2}{ 
    \ignorespaces\bfseries{#2. #3}}}
  \indentlabel{\@ifempty{#2}{\ignorespaces\bfseries{#3}}{}} 
    \vspace{1.5pt}}

\renewcommand{\tocsubsection}[3]{%
  \indentlabel{\@ifnotempty{#2}{
    \ignorespaces#2. #3}}
  \indentlabel{\@ifempty{#2}{\ignorespaces #3}{}}
    \vspace{1.5pt}}

\renewcommand{\tocsubsubsection}[3]{%
  \indentlabel{\@ifnotempty{#2}{
    \ignorespaces#2. #3}}
  \indentlabel{\@ifempty{#2}{\ignorespaces #3}{}}
    \vspace{1.5pt}}

\makeatother

\makeatletter
\def\@nomenstarted{0}
\newlength{\@nomenoldtabcolsep}

\newcommand{\nomenstart}
  {%
    \def\@nomenstarted{1}%
    \setlength{\@nomenoldtabcolsep}{\tabcolsep}%
    \setlength{\tabcolsep}{3.5pt}%
    \begin{longtable}{p{0.11\textwidth} p{0.86\textwidth}}
  }

\newcommand{\nomenitem}[2]{%
    \ifcase\@nomenstarted%
      \or 
      \or \\ 
    \fi%
    #1\,{\leavevmode\leaders\hbox{\,.}\hfill} & #2%
    \def\@nomenstarted{2}%
  }%
\newcommand{\nomenend}
  {\\%
      \end{longtable}%
      \setlength{\tabcolsep}{\@nomenoldtabcolsep}%
      \def\@nomenstarted{0}%
  }
\makeatother

\makeatletter
\newcommand{\bigish}{\bBigg@{0}}
\newcommand{\vast}{\bBigg@{4}}
\newcommand{\Vast}{\bBigg@{5}}
\newcommand{\VAST}[1]{\bBigg@{#1}}
\makeatother

\allowdisplaybreaks

\numberwithin{equation}{section}
\numberwithin{figure}{section}
\newtheorem{thm}{Theorem}[section]

\newtheorem{prop}[thm]{Proposition}
\newtheorem{lemma}[thm]{Lemma}
\newtheorem{cor}[thm]{Corollary}
\newtheorem*{thm*}{Theorem}
\newtheorem*{prop*}{Proposition}
\newtheorem*{cor*}{Corollary}
\newtheorem*{conj*}{Conjecture}

\theoremstyle{definition}
\newtheorem{definition}[thm]{Definition}
\newtheorem{notation}[thm]{Notation}

\theoremstyle{remark}
\newtheorem{rmk}[thm]{Remark}
\newtheorem{example}[thm]{Example}

\makeatletter
\newcommand{\fakephantomsection}{%
  \Hy@MakeCurrentHref{\@currenvir.\the\Hy@linkcounter}
  \Hy@raisedlink{\hyper@anchorstart{\@currentHref}\hyper@anchorend}%
  \Hy@GlobalStepCount\Hy@linkcounter%
}
\makeatother

\newcommand{\mc}{\mathcal}
\newcommand{\cA}{\mc A}
\newcommand{\cB}{\mc B}
\newcommand{\cC}{\mc C}
\newcommand{\cD}{\mc D}
\newcommand{\cE}{\mc E}
\newcommand{\cF}{\mc F}

\newcommand{\cI}{\mc I}

\newcommand{\cL}{\mc L}

\newcommand{\cO}{\mc O}

\newcommand{\cR}{\mc R}
\newcommand{\cS}{\mc S}
\newcommand{\cT}{\mc T}
\newcommand{\cU}{\mc U}
\newcommand{\cV}{\mc V}
\newcommand{\cW}{\mc W}
\newcommand{\cX}{\mc X}

\newcommand{\ms}{\mathscr}

\newcommand{\sC}{\ms C}

\newcommand{\scri}{\ms I}

\newcommand{\sP}{\ms P}

\newcommand{\sS}{\ms S}

\newcommand{\C}{\mathbb{C}}
\newcommand{\N}{\mathbb{N}}
\newcommand{\R}{\mathbb{R}}

\newcommand{\Sph}{\mathbb{S}}

\newcommand{\sfb}{\mathsf{b}}

\newcommand{\sfp}{\mathsf{p}}
\newcommand{\sfr}{\mathsf{r}}
\newcommand{\sfs}{\mathsf{s}}

\newcommand{\sfw}{\mathsf{w}}

\newcommand{\sfH}{\mathsf{H}}

\newcommand{\sfZ}{\mathsf{Z}}

\newcommand{\fm}{\mathfrak{m}}
\newcommand{\fp}{\mathfrak{p}}

\newcommand{\ft}{\mathfrak{t}}

\newcommand{\slg}{\slashed{g}{}}

\newcommand{\slDelta}{\slashed{\Delta}{}}

\newcommand{\End}{\operatorname{End}}

\renewcommand{\Re}{\operatorname{Re}}
\renewcommand{\Im}{\operatorname{Im}}
\newcommand{\Id}{\operatorname{Id}}

\newcommand{\supp}{\operatorname{supp}}

\newcommand{\tr}{\operatorname{tr}}

\newcommand{\diag}{\operatorname{diag}}

\newcommand{\eps}{\epsilon}

\newcommand{\hra}{\hookrightarrow}
\newcommand{\la}{\langle}

\newcommand{\ol}{\overline}
\newcommand{\pa}{\partial}
\newcommand{\dd}{{\mathrm d}}
\newcommand{\ra}{\rangle}
\newcommand{\spec}{\operatorname{spec}}
\newcommand{\specb}{\operatorname{spec}_\bop}

\newcommand{\ul}[1]{\underline{#1}{}}

\newcommand{\weakto}{\rightharpoonup}
\newcommand{\wh}{\widehat}
\newcommand{\wt}{\widetilde}
\newcommand{\xra}{\xrightarrow}
\newcommand{\ubar}[1]{\underaccent{\bar}#1}
\newcommand{\pfstep}[1]{$\bullet$\ \underline{\textit{#1}}}
\newcommand{\pfsubstep}[2]{{\bf#1}\ \textit{#2}}

\newcommand{\bop}{{\mathrm{b}}}
\newcommand{\cop}{{\mathrm{c}}}

\newcommand{\scop}{{\mathrm{sc}}}
\newcommand{\schop}{{\mathrm{sc,\semi}}}

\newcommand{\chop}{{\mathrm{c}\semi}}
\newcommand{\scbtop}{{\mathrm{sc}\text{-}\mathrm{b}}}

\newcommand{\ebop}{{\mathrm{e,b}}}
\newcommand{\eop}{{\mathrm{e}}}
\newcommand{\tbop}{{3\mathrm{b}}}

\newcommand{\cuop}{{\mathrm{cu}}}

\newcommand{\semi}{\hbar}

\newcommand{\ff}{\mathrm{ff}}

\newcommand{\cface}{{\mathrm{cf}}}

\newcommand{\scface}{{\mathrm{scf}}}
\newcommand{\sface}{{\mathrm{sf}}}

\newcommand{\tface}{{\mathrm{tf}}}

\newcommand{\ztface}{{\mathrm{ztf}}}
\newcommand{\sctface}{{\mathrm{sctf}}}
\newcommand{\zface}{{\mathrm{zf}}}

\newcommand{\res}{{\mathrm{res}}}

\newcommand{\cp}{{\mathrm{c}}}

\newcommand{\Diff}{\mathrm{Diff}}

\DeclareMathOperator{\Op}{Op}

\newcommand{\Vb}{\cV_\bop}
\newcommand{\Diffb}{\Diff_\bop}

\newcommand{\Psib}{\Psi_\bop}
\newcommand{\Diffch}{\Diff_\chop}
\newcommand{\Psich}{\Psi_\chop}
\newcommand{\Psisc}{\Psi_\scop}
\newcommand{\Psiscbt}{\Psi_\scbtop}
\newcommand{\Diffscbt}{\Diff_\scbtop}
\newcommand{\Diffsch}{\Diff_\schop}
\newcommand{\Psisch}{\Psi_{\scop,\semi}}

\newcommand{\etbop}{{\eop,3\bop}}
\newcommand{\betbop}{{\bop,\eop,3\bop}}
\newcommand{\Vbetb}{\cV_\betbop}
\newcommand{\Diffbetb}{\Diff_\betbop}
\newcommand{\Tbetb}{{}^\betbop T}

\newcommand{\Tetb}{{}^\etbop T}
\newcommand{\Setb}{{}^\etbop S}
\newcommand{\Diffetb}{\Diff_\etbop}
\newcommand{\Psietb}{\Psi_\etbop}

\newcommand{\Vtb}{\cV_\tbop}

\newcommand{\Difftb}{\Diff_\tbop}
\newcommand{\Psitb}{\Psi_\tbop}
\newcommand{\Veb}{\cV_\ebop}

\newcommand{\Diffeb}{\Diff_\ebop}
\newcommand{\Psieb}{\Psi_\ebop}

\newcommand{\Vscbt}{\cV_\scbtop}

\newcommand{\Vsc}{\cV_\scop}

\newcommand{\Diffsc}{\Diff_\scop}

\newcommand{\Vch}{\cV_\chop}

\newcommand{\WF}{\mathrm{WF}}

\newcommand{\Char}{\mathrm{Char}}

\newcommand{\WFb}{\WF_{\bop}}

\newcommand{\Omegab}{{}^{\bop}\Omega}

\newcommand{\Omegasc}{{}^{\scop}\Omega}

\newcommand{\WFtb}{\WF_{\tbop}}
\newcommand{\Elltb}{\mathrm{Ell_\tbop}}

\newcommand{\WFsc}{\WF_{\scop}}
\newcommand{\WFsch}{\WF_{\schop}}

\newcommand{\Tb}{{}^{\bop}T}
\newcommand{\Tcu}{{}^{\cuop}T}
\newcommand{\Tch}{{}^{\chop}T}
\newcommand{\Sch}{{}^{\chop}S}
\newcommand{\Tscbt}{{}^\scbtop T}
\newcommand{\Sscbt}{{}^\scbtop S}
\newcommand{\Tsc}{{}^\scop T}
\newcommand{\Tsch}{{}^\schop T}

\newcommand{\Teb}{{}^{\ebop}T}
\newcommand{\Ttb}{{}^{\tbop}T}

\newcommand{\Sb}{{}^{\bop}S}
\newcommand{\Ssc}{{}^{\scop}S}
\newcommand{\Seb}{{}^{\ebop}S}
\newcommand{\Stb}{{}^{\tbop}S}

\newcommand{\half}{{\tfrac{1}{2}}}

\newcommand{\sigmab}{{}^\bop\upsigma}
\newcommand{\sigmaeb}{{}^\ebop\upsigma}
\newcommand{\sigmasc}{{}^\scop\upsigma}

\newcommand{\sigmasch}{{}^{\scop,\semi}\upsigma}
\newcommand{\sigmatb}{{}^\tbop\upsigma}

\newcommand{\sigmach}{{}^\chop\upsigma}

\newcommand{\loc}{{\mathrm{loc}}}
\newcommand{\CI}{\cC^\infty}
\newcommand{\CIdot}{\dot\cC^\infty}

\newcommand{\CIc}{\cC^\infty_\cp}

\newcommand{\Hb}{H_{\bop}}

\newcommand{\Hbext}{\bar H_{\bop}}

\newcommand{\Hbsupp}{\dot H_{\bop}}

\newcommand{\Hext}{\bar H}

\newcommand{\Hsupp}{\dot H}

\newcommand{\Heb}{H_{\ebop}}

\newcommand{\Htb}{H_\tbop}

\newcommand{\Hsc}{H_{\scop}}

\newcommand{\bhm}{\fm}

\newcommand{\openbigpmatrix}[1]
  {%
    \def\@bigpmatrixsize{#1}%
    \addtolength{\arraycolsep}{-#1}%
    \begin{pmatrix}%
  }
\newcommand{\closebigpmatrix}
  {%
    \end{pmatrix}%
    \addtolength{\arraycolsep}{\@bigpmatrixsize}%
  }

\newlength{\enummargin}

\newcommand{\usref}[1]{{\upshape\ref{#1}}}

\DeclareGraphicsExtensions{.mps}

\makeatletter
\newcommand*{\fwbw}[1]{\expandafter\@fwbw\csname c@#1\endcsname}
\newcommand*{\@fwbw}[1]{\ifcase #1 \or {\rm fw}\or {\rm bw}\fi}
\AddEnumerateCounter{\fwbw}{\@fwbw}
\makeatother

\begin{document}

\title{Linear waves on asymptotically flat spacetimes. I}

\date{\today}

\begin{abstract}
  We introduce a novel framework for the analysis of linear wave equations on non-stationary asymptotically flat spacetimes, under the assumptions of mode stability and absence of zero energy resonances for a stationary model operator. Our methods apply in all spacetime dimensions and to tensorial equations, and they do not require any symmetry or almost-symmetry assumptions on the spacetime metrics or on the wave type operators. Moreover, we allow for the presence of terms which are asymptotically scaling critical at infinity, such as inverse square potentials. For simplicity of presentation, we do not allow for normally hyperbolic trapping or horizons.

  In the first part of the paper, we study stationary wave type equations, i.e.\ equations with time-translation symmetry, and prove pointwise upper bounds for their solutions. We establish a relationship between pointwise decay rates and weights related to the mapping properties of the zero energy operator. Under a nondegeneracy assumption, we prove that this relationship is sharp by extracting leading order asymptotic profiles at late times. The main tool is the analysis of the resolvent at low energies.

  In the second part, we consider a class of wave operators without time-translation symmetry which settle down to stationary operators at a rate $t_*^{-\delta}$ as an appropriate hyperboloidal time function $t_*$ tends to infinity. The main result is a sharp solvability theory for forward problems on a scale of polynomially weighted spacetime $L^2$-Sobolev spaces. The proof combines a regularity theory for the non-stationary operator with the invertibility of the stationary model established in the first part. The regularity theory is fully microlocal and utilizes edge-b-analysis near null infinity, as developed in joint work with Vasy, and 3b-analysis in the forward cone.
\end{abstract}


\author{Peter Hintz}
\address{Department of Mathematics, ETH Z\"urich, R\"amistrasse 101, 8092 Z\"urich, Switzerland}
\email{peter.hintz@math.ethz.ch}

\maketitle

\setlength{\parskip}{0.00in}
\tableofcontents
\setlength{\parskip}{0.05in}

\section{Introduction}
\label{SI}

We study pointwise and weighted $L^2$-bounds for solutions of wave equations on asymptotically flat spacetimes $(M,g)$ of dimension $n+1\geq 2$. The primary purpose is the development of a flexible analytic framework which
\begin{itemize}
\item allows for the spacetimes and wave operators to be non-stationary (though settling down to stationary models as time goes to infinity);
\item takes full advantage of spectral information about the stationary model operators;
\item is applicable to (tensorial) equations without any conditions on the existence of approximate symmetries or almost conserved energies. 
\end{itemize}

Spectral theory, resolvent analysis, and the Fourier transform in time cannot be applied directly for the analysis of an operator $P$ which is non-stationary, i.e.\ which breaks time translation invariance. Instead, insofar as they can be used to prove estimates for the stationary model operator $P_0$ (which the non-stationary operator $P$ is a decaying perturbation of), they provide the \emph{first} of \emph{two ingredients} for the global analysis of the wave operator $P$ on a non-stationary spacetime, namely, control of decay (in the forward timelike cone). We prove resolvent estimates (with a particular focus on low frequencies) for a general class of stationary wave operators and use them to prove pointwise decay estimates, which we moreover show to be sharp under certain nondegeneracy conditions.

The \emph{second ingredient}---control of an appropriate notion of regularity for solutions of the non-stationary equation $P u=f$---requires different techniques; in the present paper, we develop a fully microlocal framework for this purpose. The combination of control of decay and control of regularity of solutions of a linear partial differential equation (PDE) on a noncompact space implies the Fredholm property of the PDE via a Rellich-type compactness theorem. For the wave equations under study here, we can further improve the Fredholm property to invertibility by taking advantage of their hyperbolic character.

In this manner, we deduce the solvability of $P u=f$ together with quantitative estimates on spacetime Sobolev spaces with polynomially decaying weights, as well as pointwise estimates which follow from these via Sobolev embedding. We are able to treat wave type equations on spacetimes of arbitrary dimension, without symmetry assumptions, acting on sections of vector bundles, and with scaling critical terms (such as asymptotically inverse square potentials). We refer the reader to~\S\ref{SsIL} for a detailed comparison of our approach and results with the literature.

The present paper is entirely concerned with the theory of \emph{linear} wave equations. A natural place where non-stationary perturbations of stationary asymptotically flat spacetimes arise is the theory of \emph{quasilinear} wave equations; see~\S\ref{SsIL} for references. While the present work originated in the study of stability problems for black hole spacetimes in general relativity, a discussion of a number of features present already in linearizations of such problems (most notably the existence of stationary states, and the---however well-understood---analysis at trapping and horizons) is deferred to later work, as are applications of our theory to nonlinear wave equations. (Regarding the latter, the basic idea is to use iteration schemes on the precise spacetime function spaces developed here, much as in \cite{HintzVasySemilinear,HintzQuasilinearDS,HintzVasyQuasilinearKdS,HintzVasyKdSStability}.)  We hope that by restricting to nontrapping geometries and linear wave equations subject to mode stability and no-zero-energy-resonance conditions, the conceptual content of this work is more clearly visible.

\subsection{An example; overview}
\label{SsIEx}

We illustrate our results in the simple setting of a wave equation on Minkowski space $\R^{1+n}=\R_t\times\R^n_x$ with a non-stationary potential; we stress again that our results are vastly more general (see the discussion towards the end of~\S\ref{SsIEx}). We measure regularity (later called \emph{b-regularity}) on $\R^n_x$, resp.\ $\R^{1+n}_{t,x}$, using the collection of vector fields
\begin{equation}
\label{EqIZ}
\begin{alignedat}{3}
  Z_0 &= \{ \pa_{x^i},\ \Omega_{j k},\ S_0 \},&\qquad \Omega_{j k}&=x^j\pa_{x^k}-x^k\pa_{x^j},&\quad S_0&=\sum_{j=1}^n x^j\pa_{x^j},\quad \text{resp.} \\
  Z &= \{ \pa_t,\ \pa_{x^j},\ \Omega_{j k},\ \cL_j,\ S \}, &\qquad \cL_j&=x^j\pa_t+t\pa_{x^j},&\quad S&=t\pa_t + \sum_{j=1}^n x^j\pa_{x^j}.
\end{alignedat}
\end{equation}
We write $g_0=-\dd t^2+\sum_{j=1}^n (\dd x^j)^2$ and $\Box_{g_0}=-D_t^2+\Delta_{\R^n}$ for the Minkowski metric and the associated wave operator, where $\Delta_{\R^n}=\sum_{j=1}^n D_{x^j}^2$ and $D=\frac{1}{i}\pa$. We write
\[
  S^{-\alpha}(\R^n) = \cA^\alpha(\ol{\R^n}) = \{ a \colon \R^n\to\C \colon |Z_0^J a|\lesssim\la x\ra^{-\alpha}\ \text{for all multi-indices}\ J \}
\]
for the space of functions which are \emph{conormal} with weight $\alpha$. Finally, we shall work in the future causal cone
\[
  t_*\geq 0,\qquad t_*:=t-r-1,\quad r:=|x|.
\]
In this region, we define
\begin{equation}
\label{EqIRhos}
  \rho_{\!\scri} := \frac{\la t-r\ra}{t},\qquad
  \rho_+ := \frac{t}{\la r\ra\la t-r\ra},\qquad
  \rho_\cT := \frac{\la r\ra}{t}.
\end{equation}

\begin{figure}[!ht]
\centering
\includegraphics{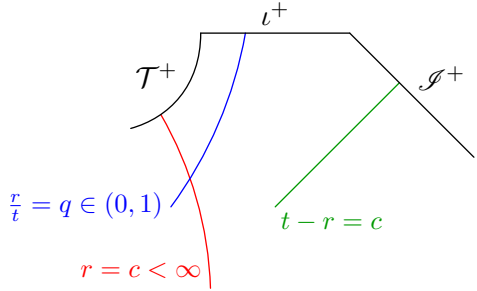}
\caption{Illustration of the region $r=|x|>1$, $t-r>0$ in $\R^{1+n}$, drawn in a compactified manner. Thus, $\cT^+$ (the future translation face), resp.\ $\iota^+$ (punctured future timelike infinity) is the closure of the endpoints at future infinity of curves along which $r$ remains bounded, resp.\ $\frac{r}{t}$ remains in a compact subset of $(0,1)$, while $\scri^+$ is the closure of the set of endpoints of future lightlike curves along which $t-r$ remains bounded. Also shown are level sets of $r$, $\frac{r}{t}$, and $t-r$. The functions $\rho_\cT$, $\rho_+$, and $\rho_{\!\scri}$ in~\eqref{EqIRhos} vanish simply at $\cT^+$, $\iota^+$, and $\scri^+$, respectively, and are positive elsewhere.}
\label{FigIRhos}
\end{figure}

\begin{rmk}[Weights]
\label{RmkIWeights}
  Roughly speaking, $\rho_{\!\scri}$ vanishes only at \emph{future null infinity} (denoted $\scri^+$); further, $\rho_+$ vanishes only in the transition region $\iota^+$ (called \emph{punctured future timelike infinity}) where $r=q t$, $q\in(0,1)$, at $t=\infty$; and $\rho_\cT$ vanishes only at future timelike infinity in spatially compact regions (denoted $\cT^+$ and called the \emph{future translation face}). See Figure~\ref{FigIRhos}. We will give precise meaning to these statements by identifying $\rho_{\!\scri},\rho_+,\rho_\cT$ with the defining functions of the boundary hypersurfaces of a compactification of $\R^{1+n}$ to a manifold with corners.
\end{rmk}

\begin{thm}[Decay for a wave equation with potential]
\label{ThmIV}
  Let $n\geq 3$. Let $\delta\in(0,1]$, and let $V_0\in\cA^{2+\delta}(\ol{\R^n})$ denote a complex-valued potential. Define the stationary operator $P_0:=\Box_{g_0}+V_0$. Assume that $P_0$ satisfies mode stability and has no resonance at zero energy.\footnote{We explain these notions after the statement of the Theorem.} Let moreover $\tilde V\in\CI(\R^{1+n})$ with\footnote{Thus, $\tilde V$ has $r^{-1-\delta}$ decay at null infinity, $t^{-\delta}$ decay in spatially compact sets, and $t^{-2-\delta}$ decay in the region $\eta t<r<(1-\eta)t$, $\eta>0$. As a special case, potentials $\tilde V$ satisfying $|Z^J\tilde V|\lesssim t_*^{-p}r^{-q}$ for $p>0$, $q>1$ with $p+q>2$ satisfy~\eqref{EqIVDecay} for $\delta=\min(1,p,q-1,p+q-2)$.}
  \begin{equation}
  \label{EqIVDecay}
    |Z^J\tilde V(t,x)| \lesssim \rho_{\!\scri}^{1+\delta}\rho_+^{2+\delta}\rho_\cT^\delta = t^{-\delta}\la r\ra^{-2}\frac{t}{\la t-r\ra}
  \end{equation}
  in $t_*\geq 0$ for all multi-indices $J$, and set $P:=\Box_{g_0}+V$ where $V:=V_0+\tilde V$. Let $f$ be a distribution on $\R^{1+n}$ with support in $t_*\geq 0$, and let $u_0,u$ denote the forward solutions (i.e.\ $t_*\geq 0$ on $\supp u_0$, $\supp u$) of
  \[
    P_0 u_0=f,\qquad P u=f.
  \]
  \begin{enumerate}
  \item\label{ItIVSchwartz}{\rm (Pointwise bounds for Schwartz forcing.)} Suppose that $f\in\sS(\R^{1+n})$. Then $|Z^J u|\lesssim(\rho_{\!\scri}\rho_+\rho_\cT)^{\frac{n-1}{2}-\eps}=t^{-\frac{n-1}{2}+\eps}$ for all $\eps>0$ and all multi-indices $J$.
  \item\label{ItIVL2}{\rm ($L^2$ bounds for forcing in weighted $L^2$ spaces.)} Let $\alpha_{\!\scri},\alpha_+,\alpha_\cT\in\R$, and suppose that
  \[
    \alpha_{\!\scri}<-\frac12,\qquad
    \alpha_+<-\frac12+\alpha_{\!\scri},\qquad
    \alpha_\cT \in \Bigl(\alpha_+-\frac{n}{2}+2,\ \alpha_++\frac{n}{2}\Bigr).
  \]
  Writing weighted spaces as $w L^2=\{ u \colon w^{-1}u\in L^2 \}$, we then have
  \begin{equation}
  \label{EqIVL2}
    Z^J f \in \rho_{\!\scri}^{\alpha_{\!\scri}+1}\rho_+^{\alpha_++2}\rho_\cT^{\alpha_\cT}L^2(\R^{1+n})\ \forall\,J \implies Z^J u \in \rho_{\!\scri}^{\alpha_{\!\scri}}\rho_+^{\alpha_+}\rho_\cT^{\alpha_\cT}L^2(\R^{1+n})\ \forall\,J.
  \end{equation}
  \item\label{ItIVLinfty}{\rm (Pointwise bounds for polynomially weighted forcing.)} Let $\gamma_{\!\scri},\gamma_+,\gamma_\cT\in\R$, and suppose that
  \[
    \gamma_{\!\scri}<\frac{n-1}{2},\qquad
    \gamma_+<\gamma_{\!\scri},\qquad
    \gamma_\cT\in(\gamma_+-n+2,\gamma_+).
  \]
  Let $\eps>0$ be arbitrary. Then
  \[
    Z^J f \in \rho_{\!\scri}^{\gamma_{\!\scri}+1}\rho_+^{\gamma_++2}\rho_\cT^{\gamma_\cT}L^\infty(\R^{1+n})\ \forall\,J \implies Z^J u \in \rho_{\!\scri}^{\gamma_{\!\scri}-\eps}\rho_+^{\gamma_+-\eps}\rho_\cT^{\gamma_\cT-\eps}L^\infty(\R^{1+n})\ \forall\,J.
  \]
  \item\label{ItIVStat}{\rm (Pointwise bounds and asymptotic profiles for the stationary problem.)} Suppose that $f\in\sS(\R^{1+n})$. Then $|Z^J u_0|\lesssim \rho_{\!\scri}^{\frac{n-1}{2}}(\rho_+\rho_\cT)^{n-1}=t^{-\frac{n-1}{2}}\la t-r\ra^{-\frac{n-1}{2}}$. Moreover, when $n$ is even, then $u_0$ has an asymptotic profile as $t_*\to\infty$: denoting by $a_\cT\in\cA^0(\ol{\R^n})$ the unique stationary (i.e.\ time-independent) solution of $P_0 a_\cT=0$ with $a_\cT-1\in\cA^\eta(\ol{\R^n})$ for some $\eta>0$, and letting $a_+(v)=(\frac{v}{v+2})^{\frac{n-1}{2}}$, we have, for some $\eta>0$,
  \begin{equation}
  \label{EqIVStat}
    \Bigl|u_0(t,x)-c\la t-r\ra^{-n+1} a_\cT(x)a_+\Bigl(\frac{\la t-r\ra}{\la r\ra}\Bigr)\Bigr|\lesssim \rho_{\!\scri}^{\frac{n-1}{2}}(\rho_+\rho_\cT)^{n-1+\eta}=t^{-\frac{n-1}{2}}\la t-r\ra^{-\frac{n-1}{2}-\eta}\,.
  \end{equation}
  The constant $c=c(f)\in\C$ is nonzero unless $f$ lies in a positive codimension subspace of $\{f\in\sS(\R^{1+n})\colon t_*\geq 0\ \text{on}\ \supp f\}$.
  \end{enumerate}
\end{thm}

Here, \emph{mode stability} is the statement that for $\sigma\neq 0$, $\Im\sigma\geq 0$, there does not exist a conormal function (with any weight) $u\neq 0$ on $\R_x^n$ with $\wh{P_0}(\sigma)(e^{i\sigma\la r\ra}u)=0$, where $\wh{P_0}(\sigma)=\Delta_{\R^n}+V_0-\sigma^2$ is the \emph{spectral family} of $P_0$. Moreover, $P_0$ not having a \emph{zero energy resonance} means that there does not exist $0\neq u\in\cA^\alpha(\ol{\R^n})$, $\alpha>0$, with $\wh{P_0}(0)u=0$. (These two assumptions together are equivalent to the requirement that the resolvent $\wh{P_0}(\sigma)^{-1}\colon L^2(\R^n)\to L^2(\R^n)$ extend analytically from $\Im\sigma\gg 1$ to $\Im\sigma>0$, and that it be continuous down to $\Im\sigma=0$ as a map $\CIc(\R^n)\to\CI(\R^n)$.) A simple class of potentials for which mode stability can easily be checked is given by nonnegative real-valued potentials $V_0\in\cA^{2+\delta}(\ol{\R^n})$ and small (complex-valued) perturbations thereof in $\cA^{2+\delta}(\ol{\R^n})$.

Theorem~\ref{ThmIV} is proved in~\S\ref{SsEW}; see Remark~\ref{RmkEWIntro}. Part~\eqref{ItIVLinfty} of Theorem~\ref{ThmIV} follows from part~\eqref{ItIVL2} via Klainerman--Sobolev inequalities \cite{KlainermanUniformDecay} (or equivalently via Sobolev embedding for b-Sobolev spaces on a suitable compactification of $\R^{1+n}$), with the $L^\infty$- and $L^2$-weights being related by $\gamma_{\!\scri}=\alpha_{\!\scri}+\frac{n}{2}$, $\gamma_+=\alpha_++\frac{n+1}{2}$, and $\gamma_\cT=\alpha_\cT+\frac12$. Part~\eqref{ItIVSchwartz} is a special case of part~\eqref{ItIVLinfty} where one chooses $\gamma_{\!\scri},\gamma_+,\gamma_\cT$ as large as possible.

\begin{rmk}[Pointwise and $L^2$-decay]
\label{RmkIDecay}
  The pointwise $\la t-r\ra^{-(n-1)}$ decay of the leading order term in~\eqref{EqIVStat} in $r<q t$, $q\in(0,1)$, is well-known for solutions of the wave equation on odd-dimensional Minkowski space, likewise for the $r^{-\frac{n-1}{2}}$ decay at null infinity. Note that the pointwise bounds which we prove for the solution $u$ of the non-stationary equation (part~\eqref{ItIVSchwartz}) are weaker than those for the solution $u_0$ of the stationary problem (part~\eqref{ItIVStat}); this is discussed in Remark~\ref{RmkWbCompStat}.
\end{rmk}

In this paper, we shall analyze spacetimes and operators which generalize the setting of Theorem~\ref{ThmIV} in the following ways.
\begin{enumerate}
\item\label{ItIStat0} The underlying spacetime manifold may have any dimension $\geq 2$ (under suitable spectral assumptions on the stationary operator).
\item\label{ItIStat1} The stationary operators may involve potentials, i.e.\ zeroth order terms, which to leading order have inverse square decay as $r\to\infty$, i.e.\ which are scaling critical in that they have the same homogeneity as the Laplacian on $\R^n$ with respect to dilations. We can also allow for general first order terms whose coefficients (relative to $\pa_t$, $\pa_{x^i}$) are, to leading order, homogeneous of degree $-1$ with respect to scaling. See Example~\ref{ExGSOPot}.
\item\label{ItIStat2} The stationary model metric, which in Theorem~\ref{ThmIV} is the Minkowski metric, can be very general, and need not even be asymptotic to Minkowski space as $r\to\infty$.\footnote{Natural examples include products of $\R_t$ with asymptotically conic Riemannian manifolds.} For simplicity of presentation, we impose a nontrapping condition here; see however Remark~\ref{RmkITrapping} below. The main assumptions on the stationary model operator are then mode stability and the absence of zero energy resonances (the precise definition of which, specifically regarding the pointwise decay as $r\to\infty$ of putative zero energy states, depends on the operator in question).
\item\label{ItITensor} Our methods apply directly to equations on vector bundles (without any requirements on symmetry or almost-symmetry, or the existence of positive definite fiber inner products with special properties). For example, Theorem~\ref{ThmIV} remains valid when $V_0$ and $\tilde V$ are valued in the space of complex $N\times N$ matrices, and correspondingly $f,u_0,u$ are valued in $\C^N$.
\item\label{ItINonstat1} The spacetime metric, which determines the principal part of the wave operators under consideration, roughly speaking only needs to settle down to a stationary metric at a rate $t_*^{-\delta}$, $\delta>0$. The class of metrics we can allow in $3+1$ dimensions includes those arising in nonlinear stability problems for asymptotically flat solutions of Einstein's field equations \cite{ChristodoulouKlainermanStability,KlainermanNicoloPeeling,LindbladAsymptotics,HintzVasyMink4,KlainermanSzeftelPolarized,DafermosHolzegelRodnianskiTaylorSchwarzschild,KlainermanSzeftelKerr}; see also \cite{HintzVasyScrieb,HintzMink4Gauge}. The non-stationary wave operators are subject to a linear version of a weak null condition at null infinity \cite{LindbladRodnianskiWeakNull}, \cite[Remark~1.7]{HintzVasyMink4}.
\item\label{ItIFinite} We prove finite regularity versions of Theorem~\ref{ThmIV}\eqref{ItIVL2}. More precisely, we present a sharp regularity theory for solutions of $P u=f$ on a scale of weighted spacetime $L^2$-Sobolev spaces; the notion of regularity for these spaces at future null infinity $\scri^+$ and at the future translation face $\cT^+$ is different from (indeed, weaker than) the b-regularity (i.e.\ regularity under $Z$) used in~\eqref{EqIVL2}. Concretely, we test for regularity in $t_*\geq 0$, $r>\frac14 t$, using \emph{edge-b-vector fields} $t\pa_t+r\pa_r$, $r(\pa_t+\pa_r)$, $\sqrt{\frac{t-r}{t}}\Omega_{j k}$ following \cite{HintzVasyScrieb}, and in $t_*\geq 0$, $r<\frac34 t$ using \emph{3-body/b-vector fields} (or \emph{3b-vector fields}) which in $r>\half$ are $r\pa_t$, $r\pa_r$, $\Omega_{j k}$, and which in $r<1$ are $\pa_t,\pa_x$. In the overlap region $\frac14 t<r<\frac34 t$, this gives rise to the same notion of regularity, and we shall speak globally of \emph{edge-3b-regularity}. Thus, we will prove estimates of the form
  \begin{equation}
  \label{EqIEst}
    \| \rho_{\!\scri}^{-\alpha_{\!\scri}}\rho_+^{-\alpha_+}\rho_\cT^{-\alpha_\cT}u \|_{H_\etbop^s} \leq C\| \rho_{\!\scri}^{-\alpha_{\!\scri}-1}\rho_+^{-\alpha_+-2}\rho_\cT^{-\alpha_\cT}P u \|_{H_\etbop^{s-1}}
  \end{equation}
  for suitable orders $s$ (which will typically need to be variable, see~\S\ref{SssIReg}) measuring the amount of regularity with respect to these vector fields. (Regularity under $Z$ will be proved as extra regularity on top of an appropriate amount of such b-edge-3-body-regularity.) Edge-3b-vector fields are very natural for the study of $P$ since up to an overall weight $P$ is to leading order a Lorentzian signature quadratic form in these; see~\S\ref{SsIN}.
\end{enumerate}
Points~\eqref{ItIStat0}--\eqref{ItITensor} are discussed in detail in~\S\ref{SsIS}, and points~\eqref{ItINonstat1}--\eqref{ItIFinite} in~\S\ref{SsIN}. We do \emph{not} handle:
\begin{itemize}
\item zero energy resonances;
\item trapping, ergoregions, or horizons, which are present on spacetimes which settle down to Schwarzschild \cite{SchwarzschildPaper} or Kerr \cite{KerrKerr} spacetimes at late times;
\item stationary operators for which mode stability fails, i.e.\ for which exponentially growing or purely oscillatory mode solutions exist;
\item geometric singularities, such as (asymptotically) stationary obstacles, or timelike curves of cone points;
\item massive waves, i.e.\ solutions of the Klein--Gordon equation $(\Box_{g_0}+m^2)u=0$ and generalizations (where $m\in\R\setminus\{0\}$).
\end{itemize}

\begin{rmk}[Trapping]
\label{RmkITrapping}
  The generalization of our methods to spacetimes with normally hyperbolic trapping (as well as ergoregions and horizons) is straightforward on a conceptual and technical level, but would require a number of notational modifications. See \cite[\S4]{HintzPrice} for stationary examples in $3+1$ dimensions; the general non-stationary case will be discussed in follow-up work \cite{HintzVasyNonstat2} by means of \cite{HintzPolyTrap}. For the purposes of this introduction, we shall thus discuss related work without regards to whether the underlying spacetimes have normally hyperbolic trapping or not.
\end{rmk}

Allowing for zero energy resonances is considerably more involved, and indeed will be the main advance of \cite{HintzVasyNonstat2}. The assumption of mode stability for nonzero frequencies on the other hand is crucial both here and in~\cite{HintzVasyNonstat2} (unlike in the work \cite{MetcalfeSterbenzTataruLED}, discussed further below). Allowing for some classes of singularities is certainly possible, see for example \cite{BaskinMarzuolaCone} for lines of cone points, but we shall not pursue such generalizations here. Moreover, with the sole exception of the asymptotic profiles for solutions of stationary equations as in Theorem~\ref{ThmIV}\eqref{ItIVStat}, we shall content ourselves with weighted $L^2$ memberships and pointwise upper bounds for solutions of non-stationary equations here. The existence and nature of partial asymptotic expansions of $u$ (including asymptotics at null infinity, the existence of radiation fields, and the existence of asymptotic profiles for solutions of non-stationary wave equations), for appropriate forcing terms $f$, will be discussed elsewhere.

Finally, Klein--Gordon type equations have a fundamentally different structure, and correspondingly their solutions $u$ have a fundamentally different asymptotic behavior, at null infinity and at future timelike infinity; see e.g.\ \cite[Theorem~7.2.7]{HormanderNonlinearLectures} for the case of the Klein--Gordon equation on Minkowski space. Roughly speaking, the qualitative asymptotic behavior of massive waves is determined by local information at null and future timelike infinity, unlike in the case of massless waves whose asymptotic behavior is global in character (see also \cite{BaskinVasyWunschRadMink}). Microlocal treatments of the Klein--Gordon equation can be found in~\cite{GellRedmanHaberVasyFeynman,SussmanKG}, and (asymptotically) stationary perturbations of the Klein--Gordon equation are studied in a 3-body-scattering framework (which takes into account spectral information at what we call $\cT^+$ here) in \cite{BaskinDollGellRedmanKG}.

Having laid out the scope of the paper, we now proceed to discuss the relationship of the results proved here with existing works in some detail; see~\S\ref{SsIL}. In~\S\S\ref{SsIS}--\ref{SsIN}, we describe the classes of spacetimes and operators of interest to us. In particular, \S\ref{SsIN} introduces the central ideas of the general framework for non-stationary wave equations developed in this paper. Finally, we provide an outline for the rest of the paper in~\S\ref{SsIO}.

\subsection{Prior literature on wave decay}
\label{SsIL}

Restricted to $n=3$ spatial dimensions, Theorem~\ref{ThmIV}\eqref{ItIVL2} is closely related to results by Metcalfe--Sterbenz--Tataru \cite{MetcalfeSterbenzTataruLED}. In the context of Theorem~\ref{ThmIV}, \cite[Theorem~2.12]{MetcalfeSterbenzTataruLED} (only assuming summable $r^{-2}$ decay of $V_0$ without derivatives) shows that mode stability and the absence of zero energy resonances together are equivalent to the validity of estimates capturing \emph{integrated local energy decay} (henceforth abbreviated \emph{ILED}), see \cite[Definition~1.6]{MetcalfeSterbenzTataruLED}. For the non-stationary operator $P$, \cite[Theorems~2.16 and 2.17]{MetcalfeSterbenzTataruLED} require the potential $V$ to be almost real\footnote{This almost symmetry assumption appears to be essential in \cite[\S7.3]{MetcalfeSterbenzTataruLED}, where it is used to upgrade a two point local energy decay estimate to bounds on the energy growth rate for the non-stationary problem.} and to have summable $r^{-2}$ decay (which at null infinity is almost a full order stronger than~\eqref{EqIVDecay}, and at $\iota^+$ it is also stronger than the pointwise $r^{-2}$ decay we can handle in this paper---cf.\ point~\eqref{ItIStat1} above). Under this restriction however, \cite[Theorem~2.16]{MetcalfeSterbenzTataruLED}, which proves a dichotomy between exponential growth of solutions and the validity of ILED, permits the existence of nonzero modes.\footnote{The work \cite{MetcalfeSterbenzTataruLED} also only requires $V$ to be slowly varying in time, which is a weaker condition than settling down to a stationary potential in $t_*\geq T\gg 1$. Our methods can likewise handle the presence of small time-dependent and non-decaying first and zeroth order terms; see Remarks~\ref{RmkWNonDecay} and \ref{RmkWbNonDecay}.} We also recall that Metcalfe--Tataru--Tohaneanu \cite{MetcalfeTataruTohaneanuPriceNonstationary}, still for $n=3$, showed that ILED\footnote{Weaker versions thereof which allow for derivative losses due to trapping are also sufficient. This was combined with \emph{stationary local energy decay} estimates in \cite{MetcalfeTataruTohaneanuPriceNonstationary,LooiDecay}, which are in essence elliptic-type estimates at zero frequency with the time derivative of the wave on the right hand side, which is then estimated using ILED.} for $u$ and its coordinate derivatives implies sharp pointwise $t^{-1}\la t-r\ra^{-2}$ decay when $\delta=1$ in Theorem~\ref{ThmIV}. For general $\delta$, Looi \cite{LooiDecay} proves $t^{-1}\la t-r\ra^{-1-\delta}$ bounds. (For Strichartz estimates under the assumption of ILED, see \cite{MetcalfeTataruGlobal}.) These pointwise decay rates (for $\delta=1$) are stronger than our bounds for $n=3$ by a factor of $\la t-r\ra^{-1}$, resp.\ $\la t-r\ra^{-2}$ in the stationary, resp.\ non-stationary setting; in the stationary setting, we explain in~\S\ref{SsIS} how this discrepancy arises from the fact that the wave operator on Minkowski space in $n+1$ dimensions with odd $n\geq 3$ violates a certain nondegeneracy condition. Compared to the local energy norms utilized in \cite{MetcalfeTataruGlobal,MetcalfeSterbenzTataruLED,LooiDecay}, we have a large degree of flexibility in the choice of spatial and temporal weights for the forcing and the solution in part~\eqref{ItIVL2}; at the level of generality of the present paper, this flexibility is essential. In particular, ILED or uniform energy bounds do not hold in general here (and even in settings where they do hold, our results do not imply them).

\begin{rmk}[Local energy spaces]
\label{RmkILLE}
  In $n=3$ spatial dimensions, we have the following rough relationships of the spaces $LE^1$ and $LE^*$ defined in \cite[\S1.4]{MetcalfeTataruTohaneanuPriceNonstationary} to weighted spacetime $L^2$-spaces, valid for any $\eps>0$:
  \begin{alignat*}{2}
    \bigl\{ u\in \la r\ra^{-\frac32}L^2(\R^4) \colon Z u \in \la r\ra^{-\frac32}L^2(\R^4) \bigr\} &\subset LE^1 &\,\subset\,& \la r\ra^{-\frac32+\eps}L^2(\R^4), \\
    \la r\ra^{\frac12-\eps}L^2(\R^4) &\subset LE^* &\,\subset\,& \la r\ra^{\frac12}L^2(\R^4).
  \end{alignat*}
  A typical ILED estimate being $\|u\|_{LE^1}\lesssim\|f\|_{LE^*}$ when $u$ has trivial Cauchy data and $f=P u$ (see \cite[Definition~1.2]{MetcalfeTataruTohaneanuPriceNonstationary}), the space $LE^1$ for $u$, resp.\ $LE^*$ for $f$ therefore corresponds roughly (i.e.\ ignoring $\eps$-losses in the weights and the precise notion of regularity) to the spaces in~\eqref{EqIVL2} with $(\alpha_\cT,\alpha_+,\alpha_{\!\scri})=(0,-\frac32,-\frac32)$, resp.\ $(\alpha_\cT,\alpha_+,\alpha_{\!\scri})=(0,-\frac32,-\frac12)$.
\end{rmk}

These papers are developments of Tataru's seminal work \cite{TataruDecayAsympFlat} on Price's law in the stationary setting; the key idea of \cite{TataruDecayAsympFlat} to upgrade ILED estimates to stronger bounds, including sharp pointwise bounds, has predecessors in \cite{TataruParametrices,MetcalfeTataruGlobal}. (We remark that these works require that the lower order terms of the operators be scaling subcritical, i.e.\ better than $r^{-2}$, resp.\ $r^{-1}$ for zeroth, resp.\ first order terms, rendering them very short range.)

Metcalfe--Tataru--Tohaneanu \cite{MetcalfeTataruTohaneanuMaxwellSchwarzschild} upgraded ILED to strong pointwise decay estimates for the Maxwell field $F$ (a 2-form) on a class of non-stationary asymptotically flat spacetimes, including Schwarzschild or Kerr type black hole spacetimes, by working directly with the \emph{tensorial} equation (the first order Maxwell system). Moreover, while there may be zero energy bound states (Coulomb solutions) corresponding to global electric or magnetic charges, their emergence at late times is suppressed via sufficient decay of the forcing. (Earlier work by Sterbenz--Tataru \cite{SterbenzTataruMaxwellSchwarzschild} on spherically symmetric and stationary spacetimes, including Schwarzschild, allows for nonzero charges.)

The proofs of (versions of) Price's law on non-stationary spacetimes in \cite{MetcalfeTataruTohaneanuPriceNonstationary,MetcalfeTataruTohaneanuMaxwellSchwarzschild,LooiDecay} crucially use properties of the scalar wave equation on Minkowski space which are specific to the $(3+1)$-dimensional case, such as the positivity of its fundamental solution as well as its relationship to the $(1+1)$-dimensional wave equation in spherical symmetry. This is by contrast with Tataru's work \cite{TataruDecayAsympFlat} on Price's law, which (like \cite{HintzPrice}) mainly relies on the Fourier transform in time and spectral theory in the form of resolvent estimates. The present work, much like \cite{MetcalfeSterbenzTataruLED}, chooses a middle ground in which decay estimates for stationary operators, proved via spectral theory, are used to estimate certain error terms arising in estimates for non-stationary problems. We also mention the conceptually related work for time-dependent Schr\"odinger equations on hyperbolic space by Lawrie--L\"uhrmann--Oh--Shahshahani \cite{LawrieLuehrmannOhShahshahaniHypSmoothing}.

Local energy decay estimates and ILED originate in work by Morawetz \cite{MorawetzIBVP,MorawetzIdentities,MorawetzExponentialDecay,MorawetzNonlinearKG} on decay estimates for solutions of wave and Klein--Gordon equations on Minkowski space. Subsequent works established similar estimates on perturbations of Minkowski space. An influential advance was the vector field method introduced by Klainerman \cite{KlainermanUniformDecay}, a significant geometric generalization of which led to proofs of the stability of Minkowski space in $3+1$ dimensions by Christodoulou--Klainerman \cite{ChristodoulouKlainermanStability} and Lindblad--Rodnianski \cite{LindbladRodnianskiGlobalStability}; see also \cite{HintzVasyMink4} and \cite[\S6]{HormanderNonlinearLectures}. Further developments include the proof of ILED for small rough non-stationary perturbations of $\Box_g$ on Minkowski space by Metcalfe--Tataru \cite{MetcalfeTataruGlobal} and for (stationary) slowly rotating Kerr spacetimes by Tataru--Tohaneanu \cite{TataruTohaneanuKerrLocalEnergy}. Closely related Keel--Smith--Sogge type estimates \cite{KeelSmithSoggeSemilinearAlmost} were proved on curved backgrounds by Alinhac \cite{AlinhacMorawetzKSSCurved}; in the exterior of star-shaped obstacles, estimates for linear wave equations were proved by Morawetz--Ralston--Strauss \cite{MorawetzRalstonStrauss}, and (almost) global existence results for quasilinear wave equations by Metcalfe--Sogge \cite{MetcalfeSoggeExteriorQuasiHigh,MetcalfeSoggeExteriorQuasi}. For Keel--Smith--Sogge estimates on product spacetimes with asymptotically Euclidean spatial metrics, see Bony--H\"afner \cite{BonyHaefnerSemilinear}. We also recall that ILED can be combined with the $r^p$-method of Dafermos--Rodnianski \cite{DafermosRodnianskiRp,MoschidisRp} to yield stronger decay estimates in spacetime dimensions $n+1\geq 4$.

Energy methods have been successfully applied to asymptotically flat spacetimes arising in General Relativity, in particular black hole spacetimes. Early results on Schwarzschild spacetimes showed the boundedness of waves \cite{WaldSchwarzschild,KayWaldSchwarzschild} and non-quantitative pointwise decay \cite{TwainySchwarzschild}; sharp decay in spherical symmetry but in a nonlinear setting was proved in \cite{DafermosRodnianskiPrice}. Sharp $t_*^{-3}$ decay, known as Price's law \cite{PriceLawI,PriceLawII,PriceBurkoLaw}, was proved on a general class of spacetimes, including Schwarzschild and subextremal Kerr black holes, by Tataru \cite{TataruDecayAsympFlat} (using resolvent estimates) under the assumption that ILED holds; \cite{TataruTohaneanuKerrLocalEnergy} verified this for small angular momenta, whereas the full subextremal range was treated by Dafermos--Rodnianski--Shlapentokh-Rothman \cite{DafermosRodnianskiShlapentokhRothmanDecay} using a combination of energy and Fourier methods. (See \cite{HolzegelKauffmannKerrFirstOrder} for a recent extension which allows for the addition of small stationary first order terms.) Earlier results by Andersson--Blue \cite{AnderssonBlueHiddenKerr} and Dafermos--Rodnianski \cite{DafermosRodnianskiRedShift,DafermosRodnianskiKerrBoundedness} treat the case of small angular momenta. The work \cite{AnderssonBlueMaxwellKerr} by Andersson--Blue on decay to Coulomb solutions for the Maxwell equations on Kerr spacetimes, following earlier work by Blue \cite{BlueMaxwellSchwarzschild}, relies on energy estimates for scalar equations (with complex potentials) obtained after a suitable separation. Applications to semilinear equations include \cite{BlueSterbenzSemilinear,TohaneanuKerrStrichartz,LukKerrNonlinear,LindbladMetcalfeSoggeTohaneanuWangStrauss,StoginKerrWaveMap}. The first quasilinear results in asymptotically flat black hole settings were obtained by Lindblad--Tohaneanu \cite{LindbladTohaneanuSchwarzschildQuasi,LindbladTohaneanuKerrQuasi}; see \cite{LooiTohaneanuDecayNull,LooiDecayEnergyCritical,LooiDecayCubic,LooiDecayQuasilinearImproved} for further nonlinear results. More recently, Dafermos--Holzegel--Taylor--Rodnianski \cite{DafermosHolzegelRodnianskiTaylorSchwarzschild} proved the codimension $3$ nonlinear stability of the Schwarzschild family of black holes as solutions of the Einstein vacuum equations; while this is a tensorial equation, \cite{DafermosHolzegelRodnianskiTaylorSchwarzschild} phrases it is a coupled system for a large number of unknowns and utilizes a delicate hierarchical structure of the equations in which decay estimates for carefully chosen scalar quantities play a key role. Klainerman--Szeftel \cite{KlainermanSzeftelKerr} similarly base their proof of the nonlinear stability of slowly rotating Kerr spacetimes on estimates for a scalar quantity (the Teukolsky scalar). For a recent result in which a (degenerate) local energy decay estimate for a stationary problem is a key ingredient in the proof of \emph{quasilinear} existence results in a very general geometric setting, see Dafermos--Holzegel--Rodnianski--Taylor \cite{DafermosHolzegelRodnianskiTaylorQuasilinear}.

In part, the popularity and success of energy methods stem from the fact that flexible arguments (meaning: not requiring special algebraic properties of the operator) yielding energy decay on stationary backgrounds which are entirely based on physical space techniques typically apply on time-dependent perturbations with only small modifications. On the other hand, spectral methods are typically restricted to time-translation invariant settings, unless they are complemented by additional techniques, as done here or in \cite{MetcalfeSterbenzTataruLED}.

We recall in this context that the proofs of the nonlinear stability of slowly rotating Kerr--de~Sitter and Kerr--Newman--de~Sitter black holes by the author and Vasy \cite{HintzVasyKdSStability,HintzKNdSStability} (see also \cite{FangKdS}) are of this latter kind: rough exponential bounds for waves on asymptotically Kerr--Newman--de~Sitter spacetimes are obtained using a simple energy estimate, and high regularity in such exponentially weighted spaces is proved by microlocal means (i.e.\ propagation of regularity in phase space). Spectral theory on the exact black hole spacetimes is then used to extract precise asymptotic expansions of waves by putting the non-stationary terms (which have \emph{exponential} decay) of the wave operators, regarded as error terms, on the right hand side of an equation that involves only the stationary wave operator on the left hand side. (See also Remark~\ref{RmkI3bSolv}.) In the context of wave equations on cosmological spacetimes, this method has its origins in the works \cite{VasyWaveOndS,VasyMicroKerrdS,HintzVasySemilinear}, which are in turn inspired by the theory of (singular) elliptic PDE in the form developed by Melrose, Mazzeo, and others \cite{MelroseAPS,MazzeoMelroseHyp,MazzeoEdge}.

By contrast, in the asymptotically flat setting under study here, a combination of spectral methods and spacetime regularity estimates as in \cite{HintzVasyKdSStability,HintzKNdSStability} is significantly more delicate to implement, and it is not clear at all how a useful initial energy estimate could be proved: such an energy estimate would need to take place on polynomially weighted spaces so that the (merely \emph{polynomially} decaying) non-stationary terms could be regarded as error terms when proving decay. Rather, the key advance of the present paper is that we can prove the solvability for a large class of non-stationary wave type equations in polynomially bounded (in fact, decaying) spaces at all (through a combination of spectral methods and microlocal regularity estimates); see~\S\ref{SsIN}.

A novel perspective on waves on asymptotically Minkowskian spacetimes was put forth by Baskin--Vasy--Wunsch \cite{BaskinVasyWunschRadMink,BaskinVasyWunschRadMink2} following the works \cite{VasyMicroKerrdS,VasyMinkDSHypRelation}: asymptotic homogeneity of degree $-2$ with respect to the scaling vector field $S=t\pa_t+r\pa_r$ becomes the structure of main interest.\footnote{Approximate time translation invariance on the other hand plays no role anymore, and indeed is incompatible with the setup of \cite{BaskinVasyWunschRadMink} except in special cases, such as exact Minkowski space. The classes of spacetimes covered by \cite{BaskinVasyWunschRadMink,BaskinVasyWunschRadMink2} on the one hand and those covered by most other papers on asymptotically flat spacetimes which focus on approximate time-translation symmetry on the other hand thus have a very small intersection.} For spacetimes with this approximate homogeneity for large $|t|+|x|$ which near the light cone $\frac{t}{|x|}=\pm 1$ have a geometry similar to Minkowski space, \cite{BaskinVasyWunschRadMink,BaskinVasyWunschRadMink2} obtain full compound asymptotics for solutions of the wave equation on a suitable compactification of $\R^{1+n}$, $n\geq 1$, to a manifold with corners. This is one of two asymptotically flat settings in which \emph{full} asymptotics have been proved to date, the other setting being the nonlinear stability of Minkowski space as analyzed by the author with Vasy \cite{HintzVasyMink4}. The proofs by Baskin--Vasy--Wunsch combine regularity estimates (proved by microlocal means in Melrose's b-calculus \cite{MelroseAPS}) for the full wave operator with spectral theoretic information for the exactly homogeneous model problem at infinity of the spacetime. The spectral theoretic information is encoded by \emph{resonances} and \emph{resonant states} which here are solutions of the model equation which are (quasi)homogeneous with respect to spacetime dilations and which are supported in the future causal cone. In the present paper, approximate dilation-invariance will play an important role in the intermediate region $\eta t<r<(1-\eta)t$, $\eta>0$, which is the interior of $\iota^+$ in Figure~\ref{FigIRhos}; see~\S\ref{SsIN}.

Turning fully to \emph{stationary} (i.e. time-translation invariant) wave equations now, we recall further works in the context of General Relativity. The linear stability of slowly rotating Kerr spacetimes was proved by Andersson--B\"ackdahl--Blue--Ma \cite{AnderssonBackdahlBlueMaKerr} using vector field methods and by the author with H\"afner and Vasy \cite{HaefnerHintzVasyKerr} using resolvent estimates; due to the presence of a second order resolvent pole at zero frequency in \cite{HaefnerHintzVasyKerr}, generalizations of the estimates in \cite{HaefnerHintzVasyKerr} to non-stationary asymptotically Kerr spacetimes must be deferred to follow-up work. See \cite{HeLinearKerrNewman} for the linear stability of mildly charged and slowly rotating black holes. Previously, the linear stability of the Schwarzschild spacetime had been shown by Dafermos--Holzegel--Rodnianski \cite{DafermosHolzegelRodnianskiSchwarzschildStability}; see also \cite{HungKellerWangSchwarzschild,HungSchwarzschildOdd,HungSchwarzschildEven,JohnsonSchwarzschild}. In the special case of the Teukolsky equation on Kerr spacetimes, sharp asymptotics were derived by Ma--Zhang \cite{MaZhangTeukolsky} (building on \cite{MaGravityKerr}) using methods specific to the algebraic structure of the equation; see \cite{MaZhangSchwarzschildDiracSharp} for the case of the Dirac equation on Schwarzschild spacetimes, and \cite{MaMaxwellAlmost,MaMaxwellKerr} for results for the Maxwell equation. For further results, see \cite{FinsterSmollerKerrStability,DafermosHolzegelRodnianskiTeukolsky,ShlapentokhRothmanTeixeiradCTeukolskyI,ShlapentokhRothmanTeixeiradCTeukolskyII,MilletTeukolskyDecay}, with mode stability proved in \cite{WhitingKerrModeStability,ShlapentokhRothmanModeStability,AnderssonMaPaganiniWhitingModeStab,CasalsTeixeiradCModes}.

Angelopoulos, Aretakis, and Gajic have been developing an approach for extracting asymptotic profiles, including subleading terms, which is based entirely on physical space methods, with the works \cite{AngelopoulosAretakisGajicLate,AngelopoulosAretakisGajicVF} (sharp vector field methods in spherical symmetry), \cite{AngelopoulosAretakisGajicRNPrice} (Price's law on Reissner--Nordstr\"om spacetimes), and \cite{AngelopoulosAretakisGajicLog} (subleading logarithmic terms at late times) leading up to \cite{AngelopoulosAretakisGajicKerr} (asymptotic profiles and subleading terms for scalar waves on subextremal Kerr spacetimes); the leading order asymptotic profile was independently obtained previously by the author in \cite{HintzPrice}. (The spacetimes and operators in all these works violate the nondegeneracy condition discussed in~\S\ref{SsIS}, leading to stronger decay than what generically holds in the general setting of the present paper.) Gajic \cite{GajicInverseSquare} obtained asymptotic profiles for wave equations on the Schwarzschild spacetime coupled to stationary potentials with inverse square decay at infinity; not taking into account the minor issue of trapping, we are able to generalize \cite{GajicInverseSquare} to all dimensions and to complex- or matrix-valued potentials; see~\S\ref{SsEV}. An important ingredient in the physical space approach is a hierarchy of $r^p$-weighted estimates, which is related to the iterative construction of the Taylor series of the low energy resolvent in the present paper; see~\S\ref{SsIS} and \S\S\ref{SsStCo}--\ref{SAS} for further details. We also mention the work by Burq--Planchon--Stalker--Tahvildar-Zadeh \cite{BurqPlanchonStalkerTahvildarZadehInvSq} on Morawetz estimates for Schr\"odinger and wave equations coupled to \emph{exact} inverse square potentials---including the singularity at $r=0$. Similarly, Baskin--Gell-Redman--Marzuola \cite{BaskinGellRedmanMarzuolaPriceInvSq}, working in the framework of \cite{BaskinVasyWunschRadMink}, prove sharp decay and asymptotics for exact inverse square potentials on $(3+1)$-dimensional Minkowski space; besides the critical inverse quadratic decay at infinity, this work also treats the inverse square singularity at the spatial origin (relying on \cite{BaskinMarzuolaCone}). (We already mention here that even an \emph{asymptotically} inverse square potential as $r\to\infty$ gives rise to an inverse square \emph{singularity} at $\frac{r}{t}=0$ from the perspective of a dilation-invariant model problem at punctured future timelike infinity; see~\S\ref{SsIN} for further details.)

Spectral theoretic methods for Schr\"odinger and wave equations on stationary backgrounds have a long history as well, and are the main tool for the analysis of stationary operators in the present paper. Jensen--Kato \cite{JensenKatoResolvent} derived late time expansions for the Schr\"odinger flow on Euclidean space with very short range potentials in $3$ spatial dimensions from the singularity structure of the resolvent near zero energy; see also \cite[\S3.3]{DyatlovZworskiBook}, further \cite{StrohmaierWatersHodge} for detailed results for the Laplacian on differential forms, and \cite{MullerStrohmaierResolvent} for convergent low energy expansions. Work done on product spacetimes with curved spatial geometry includes the paper by Guillarmou--Hassell--Sikora \cite{GuillarmouHassellSikoraResIII} on asymptotically conic manifolds of general dimension (which moreover allows for real asymptotically inverse square potentials), which gives a leading order asymptotic profile for the Schr\"odinger and wave evolution under a nondegeneracy condition (see \cite[Corollary~1.3]{GuillarmouHassellSikoraResIII}) which is a special case of the one in the present paper. Vasy--Wunsch \cite{VasyWunschMorawetz} prove Morawetz estimates on scattering manifolds. Work by Bouclet--Burq \cite{BoucletBurqSharpDecay} gives sharp low energy resolvent bounds in any dimension for long range perturbations of the Euclidean metric, improving upon earlier work by Bony--H\"afner \cite{BonyHaefnerLED} and Bouclet \cite{BoucletLEDEuclidean}; see also \cite{BurqDecay}. Schlag--Soffer--Staubach \cite{SchlagSofferStaubachConicI,SchlagSofferStaubachConicII} treat the case of asymptotically conic warped product manifolds with hyperbolic trapping in a compact set. Building also on the earlier \cite{CostinSchlagStaubachTanveerInvSq,DonningerSchlag1dDecay,DonningerKriegerVectorField} on inverse square and inverse power law potentials in one dimension, Donninger--Schlag--Soffer \cite{DonningerSchlagSofferSchwarzschild,DonningerSchlagSofferPrice} obtained Price's law on Schwarzschild spacetimes via low energy resolvent estimates (for the late time tail) and high energy estimates (to deal with trapping); see also the survey articles \cite{SchlagSchrodingerSurvey,SchlagDecaySurvey}. On stationary spacetimes which asymptote to $(3+1)$-dimensional Minkowski space at fast inverse polynomial rates as $r\to\infty$, Morgan and Wunsch \cite{MorganDecay,MorganWunschPrice} proved fast inverse polynomial time decay using low energy resolvent estimates, thus giving an independent proof (with an $\eps$-loss) of the stationary special cases of Looi's results \cite{LooiDecay}.

Of central importance to the present paper (and the earlier \cite{HaefnerHintzVasyKerr,HintzPrice}) is work by Vasy \cite{VasyLAPLag,VasyLowEnergyLag} on limiting absorption principles and low energy resolvent estimates, which applies also to spectral families of stationary spacetimes with general (not product-type) asymptotically flat metrics. An important insight is that precise estimates for the limiting resolvent on the real axis, including at low frequencies, are more readily accessible if one uses the Fourier transform in $t_*\approx t-r$ instead of $t$.\footnote{The level sets of $t_*$ are transversal to future null infinity, and thus the behavior of waves as $t_*\to\infty$ is directly encoded in the behavior at $\sigma=0$ of the resolvent defined using $t_*$. On the other hand, the behavior of the resolvent defined using $t$ does not, upon taking the inverse Fourier transform, distinguish null infinity and future timelike infinity cleanly. See also \cite[\S1.3]{HintzPrice}, which builds on \cite{VasyLowEnergyLag}.} We also mention Sussman's recent work \cite{SussmanCoulomb} (and references therein) on low energy resolvent asymptotics in the presence of very long range (Coulomb like) potentials which have an altogether different low energy spectral theory than the inverse square type potentials allowed for here.

\subsection{Stationary operators}
\label{SsIS}

We now give a more detailed description of our general setup in the stationary case. (The non-stationary operators will be decaying (as $t_*\to\infty$) perturbations of the stationary operators introduced here; see~\S\ref{SsIN} below.) While the setup here is significantly more general than that in \cite{HintzPrice} on Price's law (ignoring horizons and trapping, which are treated in \cite{HintzPrice}) and \cite{VasyLowEnergyLag} on low energy resolvent estimates, the proof strategies for obtaining resolvent bounds and low energy resolvent expansions are similar.

\subsubsection{Metrics and operators}

On $\R^{1+n}$, $n\geq 1$, with coordinates $z=(t,x)$, we consider stationary metrics of the form
\[
  g_0 = -\dd t^2+\dd x^2 + \tilde g_0,\qquad \tilde g_0=\sum \tilde g_{0,\mu\nu}\,\dd z^\mu\,\dd z^\nu,\quad \tilde g_{0,\mu\nu}\in\cA^\delta(\ol{\R^n}),
\]
for some $\delta\in(0,1]$. Letting $r=|x|$, $\omega=\frac{x}{|x|}\in\Sph^{n-1}$, we assume that there exists a smooth function $t_*=t-\breve F(x)$ which for $r>1$ is equal to $t_*=t-r-\tilde F(r)$ with $\tilde F\in\cA^{-1+\delta}(\ol{\R^n})$, so that $\dd t_*$ is globally timelike or null and, in the coordinates $(t_*,r,\omega)$,
\[
  g_0^{-1}(\dd t_*,-) \equiv -\pa_r \bmod \cA^{1+\delta}.
\]
One should think of $t_*$ as an approximate null coordinate for large $r$, and of $\pa_r$ as an approximate null generator of the level sets of $t_*$. We require that $g_0$ be nontrapping, i.e.\ $r\to\infty$ along all future null-geodesics. The Minkowski metric, with $\tilde g_0=0$ and $t_*=t-r$ in $r>1$ (appropriately extended to $r\leq 1$), is the simplest example. Metrics which are asymptotic to a $(3+1)$-dimensional Schwarzschild or Kerr metric as $r\to\infty$ are also of this form; see Example~\ref{ExGSGKerr}.

For the sake of notational simplicity, we restrict to the case $\delta=1$ here and write $\cO(r^{-\nu})$ instead of $\cA^\nu(\ol{\R^n})$. Working in $(t_*,r,\omega)$ coordinates, the dual metric is of the form
\[
  g_0^{-1} = -2\pa_{t_*}\otimes_s\pa_r + \pa_r^2 + r^{-2}\slg^{-1} + \tilde G,
\]
where $\tilde G=\cO(r^{-1})$ (with respect to $\pa_{t_*}$, $\pa_r$, $r^{-1}\pa_\omega$) and $\tilde G(\dd t_*,-)=\cO(r^{-2})$. Here $\slg$ is the standard metric on $\Sph^{n-1}$.\footnote{In the main part of the paper, the metric $\slg$ may be an arbitrary Riemannian metric on $\Sph^{n-1}$.}

We then consider wave type operators corresponding to such metrics,
\begin{equation}
\label{EqISOp}
  P_0 = 2\pa_{t_*}\Bigl(\pa_r+\frac{n-1}{2 r}+r^{-1}S\Bigr) + \wh{P_0}(0)
\end{equation}
with $S=S(r,\omega)$ a smooth function of $(r^{-1},\omega)$ which roughly speaking accounts for possible modifications of the decay rate at null infinity from $r^{-\frac{n-1}{2}}$ to $r^{-\frac{n-1}{2}-S}$.\footnote{Such modifications arise e.g.\ in the presence of weak damping; see \cite[\S1.1]{HintzMink4Gauge}.} (For the purpose of this introduction we do not include terms which are of lower order in the sense that they do not require additional effort to take into account; also, the smoothness requirement on $S$ can be relaxed. See Definition~\ref{DefGSO} for the full setup.) The principal part of the (elliptic) \emph{zero energy operator} $\wh{P_0}(0)\in\Diff^2(\R^n)$ is asymptotic to the Euclidean Laplacian as $r\to\infty$. More precisely, we require that $\wh{P_0}(0)$ be asymptotically homogeneous of degree $-2$ with respect to dilations, so $\wh{P_0}(0)$ is equal to
\begin{equation}
\label{EqISP0}
  r^{-2} P_{(0)}(r D_r,\omega,D_\omega) = D_r^2 + r^{-2}\slDelta + a(\omega)r^{-1}D_r + b(\omega)r^{-2}D_\omega + c(\omega)r^{-2}
\end{equation}
plus lower order terms (i.e.\ operators built out of $D_r^2$, $r^{-1}D_r D_\omega$, $r^{-2}D_\omega D_\omega$, $r^{-1}D_r$, $r^{-2}D_\omega$, $r^{-2}$ with $\cO(r^{-1})$ coefficients); here $\slDelta=\Delta_\slg$ is the (nonnegative) Laplacian on $(\Sph^{n-1},\slg)$. The coefficients $a,b,c$ of the terms of differential order $\leq 1$ may be complex.

The wave operator on Minkowski space is of this form, with $S=0$ and $\wh{P_0}(0)=D_r^2-\frac{i(n-1)}{r}D_r+r^{-2}\slDelta$. To add an inverse square potential, one may take $c(\omega)\neq 0$ in~\eqref{EqISP0}. The coefficients $a(\omega)$ and $b(\omega)$ encode first order terms with scaling-critical behavior at infinity. See also Example~\ref{ExGSOPot}.

\subsubsection{Spectral information and pointwise decay}

A key piece of data about $\wh{P_0}(0)$ is its \emph{boundary spectrum} \cite{MelroseAPS}
\[
  \specb(\wh{P_0}(0)) \subset \C,
\]
consisting of all $\lambda\in\C$ for which there exists a nonzero function of the form $r^{-\lambda}u(\omega)$ which (in $r>0$) is annihilated by $P_{(0)}$. We require that $\wh{P_0}(0)$ have a nonempty \emph{indicial gap}, which is an interval
\[
  (\beta^-,\beta^+)\subset\R
\]
so that for all $\beta\in(\beta^-,\beta^+)$, the operator
\begin{equation}
\label{EqINP0}
  \wh{P_0}(0)\colon\cA^\beta(\ol{\R^n})\to\cA^{\beta+2}(\ol{\R^n})
\end{equation}
is invertible;\footnote{Later on, see in particular Lemma~\ref{LemmaStEst0}, we mainly work with $L^2$-based spaces instead. We note though that the invertibility statement here follows from the invertibility on the b-Sobolev spaces in Lemma~\ref{LemmaStEst0}. This is a consequence of the mapping properties of an elliptic parametrix for $\wh{P_0}(0)$ in the large b-pseudodifferential calculus \cite{MelroseAPS}.} this is the appropriate generalization of the notion of \emph{absence of zero energy resonances} from the case of Euclidean potential scattering. On $(3+1)$-dimensional Minkowski space, the indicial gap is $(0,1)$; see Example~\ref{ExIMink} below. The largest such interval necessarily has the property that $\beta^\pm=\Re\lambda^\pm$ for some $\lambda^\pm\in\specb(\wh{P_0}(0))$ (see Lemma~\ref{LemmaStEst0}). Up to possible logarithmic corrections, $\beta^+$ is the decay rate\footnote{in the generalized sense that negative decay rates $\beta^+<0$ are growth rates with exponent $-\beta^+$} of the Green's function of $\wh{P_0}(0)$ (with second argument being an arbitrary point in $\R^n$), while $\beta^-$ is the decay rate of the smallest (in a pointwise sense as $r\to\infty$) element in the nullspace of $\wh{P_0}(0)$ (cf.\ $a_\cT$ in Theorem~\ref{ThmIV}\eqref{ItIVStat}).

\begin{example}[Minkowski space]
\label{ExIMink}
  In the case $\wh{P_0}(0)=\Delta_{\R^n}$, we have $\specb(\Delta_{\R^n})=(-\N_0)\cup(n-2+\N_0)$. The invertibility assumption~\eqref{EqINP0} is satisfied for $n\geq 3$ for $\beta$ in the indicial gap $(0,n-2)$. For $n\in\{1,2\}$, there is no indicial gap, as can be verified using the properties of the Green's function (which for $n=1,2$ does not decay at infinity): the operator $\wh{P_0}(0)=\Delta_{\R^n}$ on $\R^n$ is injective on $\cA^\beta(\ol\R)$ only when $\beta>0$ (to exclude constants), but typically solutions of $\Delta_{\R^n}u=f$, even for compactly supported $f$, behave at infinity like the Green's function (so grow linearly when $n=1$, or logarithmically when $n=2$), i.e.\ they are \emph{not} bounded and thus do not lie in $\cA^\beta(\ol{\R^n})$ when $\beta>0$. --- In summary, the scalar wave operator on Minkowski space of dimension $n+1$ fits into the setting of the present paper if and only if $n\geq 3$. We have $S=0$ in~\eqref{EqISOp}, and hence the decay at null infinity, and the decay of outgoing spherical waves, is $\cO(r^{-\frac{n-1}{2}})$.
\end{example}

\begin{example}[Inverse square potentials]
\label{ExIMinkSq}
  If one adds to $\Delta_{\R^n}$ an inverse square potential, say with $c(\omega)=\alpha\in\C\setminus(-\infty,-(\frac{n-2}{2})^2]$ in~\eqref{EqISP0}, then the boundary spectrum becomes $\frac{n-2}{2}\pm(\nu_0+\N_0)$ where $\nu_0=\sqrt{(\frac{n-2}{2})^2+\alpha}$. The indicial gap is contained in $(\frac{n-2}{2}-\Re\nu_0,\frac{n-2}{2}+\Re\nu_0)$, though the invertibility of~\eqref{EqINP0} (and thus the existence of a nonempty indicial gap) depends on the precise expression for $\wh{P_0}(0)$. We again have $S=0$ in~\eqref{EqISOp}.
\end{example}

The assumptions of mode stability, the absence of zero energy resonances, the nontrapping assumption, and the timelike nature of $\dd t_*$ suffice for pointwise and high energy estimates for the resolvent family
\begin{equation}
\label{EqINSpecFam}
  \wh{P_0}(\sigma)^{-1},\qquad \wh{P_0}(\sigma)=-2 i\sigma\Bigl(\pa_r+\frac{n-1}{2 r}+r^{-1}S\Bigr) + \wh{P_0}(0),
\end{equation}
for $\Im\sigma\geq 0$ on $L^2$-analogues of the spaces $\cA^\alpha(\ol{\R^n})$ for appropriate weights $\alpha$: for $\sigma=0$, one needs $\alpha\in(\beta^-,\beta^+)$; for $\sigma\neq 0$, we note that $\frac{n-1}{2}+S$ is the expected decay rate of outgoing spherical waves (being the power of $r^{-1}$ which is annihilated by the first term in~\eqref{EqINSpecFam}), whence one needs $\alpha<\frac{n-1}{2}+\ubar S$ where\footnote{If one considers operators acting on sections of vector bundles, then $\ubar S$ is defined in terms of the spectrum of $S(\infty,\omega)$; see Definition~\ref{DefStEstThr}.}
\[
  \ubar S=\min_{\omega\in\Sph^{n-1}}\Re S(\infty,\omega).
\]
The final assumption we impose on $P_0$ in order to get uniform low energy estimates concerns the \emph{transition face model operator} $N_\tface^\pm(P_0)$.\footnote{This was introduced in \cite{VasyLowEnergyLag} and played an important role in \cite{HintzPrice}, where its invertibility was used to obtain uniform estimates on the low energy resolvent, and its inversion produced the leading order logarithmic term at $\sigma=0$ in the low energy resolvent expansion; in these works, this operator is the spectral family of the Laplacian on an exact conic manifold at energy $1$, whereas here it may be more general. If $P_0$ is the scalar wave operator on an asymptotically Minkowski spacetime, the operator $N_\tface^+(P_0)$ is equal to $\wt\Box(1)$ in the notation of \cite[Definition~2.20]{HintzPrice}.} The transition face model operator
\begin{equation}
\label{EqINNtf}
  N_\tface^\pm(P_0) = \mp 2 i\Bigl(\pa_{\hat r}+\frac{n-1}{2\hat r} + \hat r^{-1}S(\infty,\omega)\Bigr) + \hat r^{-2}P_{(0)}(\hat r D_{\hat r},\omega,D_\omega)
\end{equation}
arises by considering the limit of the operator $\sigma^{-2}\wh{P_0}(\sigma)$ as $\pm\sigma\searrow 0$ while $\hat r=|\sigma|r\in(0,\infty)$ remains bounded. We then require $N_\tface^\pm(P_0)$ to be invertible on function spaces with appropriate weights at $\hat r=0$ (matching the weights for the zero energy problem) and $\hat r=\infty$ (matching the weights for nonzero frequencies). See Definition~\ref{DefGSOSpec} for the full set of spectral assumptions on $P_0$.

The first main result of the paper is a decay estimate for forward solutions of $P_0$. We state this in terms of the functions
\[
  \rho_{\!\scri}=\frac{\la t_*\ra}{t},\qquad
  \rho_+ = \frac{t}{\la r\ra\la t_*\ra},\qquad
  \rho_\cT = \frac{\la r\ra}{t},
\]
which on Minkowski space and for $t_*=t-r$ are precisely the functions~\eqref{EqIRhos}.

\begin{thm}[Decay for stationary problems]
\label{ThmINDecay}
  Let $f\in\sS(\R^{1+n})$, with $t_*\geq 0$ on $\supp f$. Then the forward solution $u_0$ of $P_0 u_0=f$, i.e.\ the unique solution with $t_*\geq 0$ on $\supp u_0$, satisfies for all $\eps>0$ the pointwise bounds
  \begin{equation}
  \label{EqINDecay}
  \begin{split}
    | Z^J u_0 | &\lesssim \rho_{\!\scri}^{\frac{n-1}{2}+\ubar S-\eps}\rho_+^{\beta^++1-\eps}\rho_\cT^{\beta^+-\beta^-+1-\eps} \\
      & = t^\eps\la t_*\ra^{-(\beta^+-\beta^-+1)+\eps}\la r\ra^{-\beta^-}\Bigl(\frac{\la t_*\ra}{t}\Bigr)^{\frac{n-1}{2}+\ubar S-\beta^-}
  \end{split}
  \end{equation}
  in $t_*\geq 0$ for all multi-indices $J$. In particular, the decay rate as $t_*\to\infty$ in spatially compact regions is equal to $1$ plus the size of the indicial gap, and the decay rate in $\eta t<r<(1-\eta)r$, $\eta>0$, is $1$ plus the upper bound of the indicial gap (up to $\eps$-losses).
\end{thm}

\begin{rmk}[Initial value problems and regularity/decay assumptions]
\label{RmkINDecayIVP}
  It is standard to convert estimates for forward problems into estimates for initial value problems, though we shall not do this explicitly in this paper. In the context of Theorem~\ref{ThmINDecay}, the estimates~\eqref{EqINDecay} hold for the solution $u_0$ of $P_0 u_0=0$ with $(u_0,\pa_{t_*}u_0)|_{t_*=0}=(u^{(0)},u^{(1)})$ when $t_*^{-1}(0)$ is spacelike and the initial data $u^{(0)},u^{(1)}$ are smooth and compactly supported. The precise regularity theory developed in the second part of this paper allows one to relax the regularity and decay assumptions of forcing terms or initial data considerably.
\end{rmk}

In the context of Example~\ref{ExIMink}, the bound~\eqref{EqINDecay} is almost sharp (i.e.\ sharp up to the $\eps$-loss) on Minkowski space of odd spacetime dimension $n+1\geq 5$, whereas in even spacetime dimensions $n+1\geq 4$, where the sharp Huygens principle holds, it is not. (See also \cite{MorganDecay,MorganWunschPrice,LooiDecay} for precise statements for $n+1=4$.) But in \emph{all} dimensions $n+1\geq 2$, the bound~\eqref{EqINDecay} is almost sharp in the presence of an inverse square potential $\sim \alpha r^{-2}$ (subject to mode stability and absence of zero energy resonances) for an open and dense subset of $\alpha\in\C\setminus(-\infty,-(\frac{n-2}{2})^2]$. See Propositions~\ref{PropESSharp} and \ref{PropEV} for precise statements.

\begin{figure}[!ht]
\centering
\includegraphics{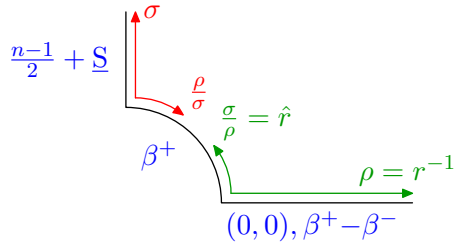}
\caption{Resolved space for the low energy analysis in $\sigma\geq 0$ (following \cite{VasyLowEnergyLag}). The horizontal axis is (the lift of) $\sigma=0$. Shown are local coordinates near the two corners, expressed in terms of $\rho=r^{-1}$. Indicated in blue are the pointwise decay orders of $\wh{P_0}(\sigma)^{-1}\hat f(\sigma,r\omega)$ (near the top left corner meaning a bound by $\sigma^{\beta^+}(\frac{\rho}{\sigma})^{\frac{n-1}{2}+{\underline{S}}}$) when $f$ is Schwartz; by `$(0,0),\beta^+{-}\beta^-$' we mean smoothness up to errors with decay order $\beta^+{-}\beta^-$ (ignoring possible logarithmic terms).}
\label{FigINRes}
\end{figure}

To prove Theorem~\ref{ThmINDecay}, we express the Fourier transform (with the sign convention commonly used in spectral theory) of $u_0$ in $t_*$,
\begin{equation}
\label{EqIFT}
  \wh{u_0}(\sigma,x) = \int_\R e^{i\sigma t_*}u(t_*,x)\,\dd t_*,
\end{equation}
in terms of $f$ via $\wh{u_0}(\sigma,x)=\wh{P_0}(\sigma)^{-1}\hat f(\sigma,x)$ for $\Im\sigma=C\gg 1$. One can shift the contour in $u_0(t_*,x)=(2\pi)^{-1}\int_{\Im\sigma=C} e^{-i\sigma t_*}\wh{P_0}(\sigma)^{-1}\hat f(\sigma,x)\,\dd\sigma$ to $\R$ using the following inputs.
\begin{itemize}
\item The nontrapping assumption implies high energy estimates for the resolvent $\wh{P_0}(\sigma)^{-1}$ in $\Im\sigma\geq 0$, $|\Re\sigma|\gg 1$; see Proposition~\ref{PropStEstHi}\eqref{ItStEstHib}. This builds on work by Vasy--Zworski \cite{VasyZworskiScl} in the form given by Vasy \cite[\S5]{VasyLAPLag}; see also \cite{BurqNontrapping,VasyMinicourse}.
\item Mode stability implies the existence of, and uniform bounds for, $\wh{P_0}(\sigma)^{-1}$ and its $\sigma$-derivatives also in compact subsets of $\Im\sigma\geq 0$, $\sigma\neq 0$; see Proposition~\ref{PropStEstNz} and \cite{MelroseEuclideanSpectralTheory,VasyLAPLag}.
\item The absence of zero energy resonances and the invertibility of the transition face model operators $N_\tface^\pm(P_0)$ imply uniform bounds on the resolvent for small $\sigma$; see Proposition~\ref{PropStEstLo}\eqref{ItStEstLob} and \cite{VasyLowEnergyLag}. These estimates take place on function spaces which keep track of weights and asymptotics in three regimes, depending on the relative size of $|\sigma|$ and $\rho=r^{-1}$. See Figure~\ref{FigINRes} (though at this point we do not prove an expansion at $\sigma=0$ yet). The blown-up corner, parameterized by $\hat r=|\sigma|r=\frac{|\sigma|}{\rho}$, is the \emph{low energy transition face}: this is the place where $N_\tface^\pm(P_0)$ lives.
\end{itemize}

Together with the conormality of the resolvent near low energies (regularity under repeated application of $\sigma\pa_\sigma$, see Lemma~\ref{LemmaStCo}), this is sufficient to give rough pointwise decay estimates; the only contribution to $(2\pi)^{-1}\int_\R e^{-i\sigma t_*}\wh{P_0}(\sigma)^{-1}\hat f(\sigma,x)\,\dd\sigma$ which is not rapidly decaying in $t_*$ arises from $\sigma\approx 0$.

Stronger decay estimates are equivalent to the existence of a partial Taylor expansion of the resolvent at (the lift to the resolved space in Figure~\ref{FigINRes} of) $\sigma=0$. We obtain this expansion by adapting the algorithmic procedure introduced in \cite{HintzPrice}. Dropping possible $\eps$- or logarithmic losses for the sake of simplicity, this proceeds as follows.
\begin{enumerate}
\item We first approximate the solution $\wh{u_0}(\sigma)$ of the spectral problem $\wh{P_0}(\sigma)\wh{u_0}(\sigma)=\hat f(\sigma)$ using the zero energy inverse $\wh{P_0}(0)^{-1}\hat f(0)\in\cA^{\beta^+}(\ol{\R^n})$; this solves the spectral problem up to an error term (produced by the first term in~\eqref{EqINSpecFam}) of class $\sigma\cA^{\beta^++1}$.
\item If $\beta^-<\beta^+-1$, we can solve this error term away by applying $\wh{P_0}(0)^{-1}$ to it, producing a $\sigma\cA^{\beta^+-1}$ (cf.\ \eqref{EqINP0}) correction to the solution of the spectral problem, and leaving an error in $\sigma^2\cA^{\beta^+}$. We can continue in this fashion as long as we remain in the invertible range of weights for $\wh{P_0}(0)$, thus for a total of roughly $\beta^+-\beta^-$ steps.
\item The final error term (of class $\sigma^{\beta^+-\beta^-}\cA^{\beta^-+2}$ when $\beta^+-\beta^-$ is an integer) is solved away by applying the true resolvent; the resulting conormal $\cO(\sigma^{\beta^+-\beta^-})$ behavior of $\wh{u_0}(\sigma)$ gives rise to $t_*^{-(\beta^+-\beta^-+1)}$ decay of $u_0$ in spatially compact sets upon inverse Fourier transforming.
\item Moreover, since $t_*/r$ is a coordinate along punctured future timelike infinity $\iota^+$, the behavior of $u_0$ at $\iota^+$ can be determined by taking the Fourier transform of $\wh{u_0}(\sigma)$ in $\sigma r$---which is precisely the coordinate along the low energy transition face. The $r^{-\beta^+}$ bound on $\wh{u_0}(\sigma)$ for $|\sigma r|=|\frac{\sigma}{\rho}|\lesssim 1$ translates into the $r^{-1-\beta^+}$ bound at $\iota^+$ in~\eqref{EqINDecay} (noting that $t_*\sim r$ there).
\end{enumerate}
See Proposition~\ref{PropStCoReg} and Theorem~\ref{ThmStCo} for details.

This procedure is closely related to the physical space method developed by Angelopoulos, Aretakis, and Gajic. Specifically, in Step~3 in \cite[\S4.2]{GajicInverseSquare}, the wave equation under study is rewritten so as to feature the zero energy operator $\wh{P_0}(0)$, denoted $\cL$ in the reference, on the left hand side, and the remaining terms (all of which involve at least one time derivative), denoted $F\circ T$ (with $T=\pa_{t_*}$) in the reference, on the right hand side. More precisely, this is done at the level of initial data: the data of $T^{-n}u$ ($n$-fold integration from $t_*=\infty$) are put in relation to those of $T^{-(n-1)}u$ for $n=1,\ldots,N$, with a hierarchy of $r^p$-weighted energy estimates providing ILED estimates for $T^{-N}u$. On the spectral side, this roughly corresponds to starting with good regularity estimates for $\sigma^N\hat u(\sigma)$ at $\sigma=0$, where $\hat u(\sigma)=\wh{P_0}(\sigma)^{-1}\hat f(\sigma)$, and extracting an expansion for $\hat u(\sigma)$ itself at (the lift of) $\sigma=0$ by iterating the formula $\hat u(\sigma)=\wh{P_0}(0)^{-1}(\hat f(\sigma)-(\wh{P_0}(\sigma)-\wh{P_0}(0))\hat u(\sigma))$, where $\wh{P_0}(\sigma)-\wh{P_0}(0)$ corresponds to the conjugation of $F\circ T$ by the Fourier transform in $t_*$.

Finally, we explain the aforementioned nondegeneracy condition under which the decay estimate~\eqref{EqINDecay} is sharp up to the $\eps$-loss, and indeed so that $u_0$ has a leading order asymptotic profile as $t_*\to\infty$ which is consistent with~\eqref{EqINDecay} with $\eps=0$ (as in Theorem~\ref{ThmIV}\eqref{ItIVStat}). Namely, in the above procedure, one can keep track of the leading terms of the asymptotic expansions of all terms in the Taylor expansion. In the last step one can extract a leading order singular term in $\sigma$, which gives a leading order asymptotic profile for $u_0$ at $\cT^+$; to get the asymptotic profile at $\iota^+$, one also needs to compute the leading order term at the low energy transition face. (By contrast, these two asymptotic profiles are extracted in a single step in Step 4 of \cite[\S4.2]{GajicInverseSquare}.) The nondegeneracy condition now is that \emph{none of these leading order terms vanish}. This requires in the first step that $\wh{P_0}(0)^{-1}\hat f(0)$ have a nonvanishing size $r^{-\beta^+}$ leading order term, which is true for generic $f$; and in later steps this requires that the term $\pa_r+\frac{n-1}{2 r}+r^{-1}S$ in~\eqref{EqINSpecFam} not annihilate the leading order terms of the Taylor coefficients of $\wh{u_0}(\sigma)$ in $\sigma$. This thus amounts to excluding integer coincidences between the decay rate of the Green's function of the zero energy operator on the one hand and the decay rate of outgoing spherical waves on the other hand. Correspondingly, the computations in Examples~\ref{ExIMink}--\ref{ExIMinkSq} explain why the wave equation on $(n+1)$-dimensional Minkowski space with $n+1\geq 5$ odd, and most inverse square potentials regardless of the spacetime dimension, satisfy this nondegeneracy condition, while Minkowski space with $n+1\geq 4$ even violates it. This applies to the settings of \cite{GajicInverseSquare} and \cite[Corollary~1.3]{GuillarmouHassellSikoraResIII}. See~\S\ref{SAS} for details, in particular Remarks~\ref{RmkASCounterex1} and \ref{RmkASCounterex2} as well as assumptions~\eqref{EqASAssmSimple}, \eqref{EqASAssmfLead}, and~\eqref{EqASAssmtfLead}. The rough result thus reads:

\begin{thm}[Sharp asymptotics for stationary problems]
\label{ThmISharp}
  Let $P_0$ denote a stationary wave type operator which satisfies the above spectral assumptions and, in addition, a (conjecturally generic) nondegeneracy condition. Then the decay rates of forward solutions $u_0$ of $P_0 u_0=f$ for generic Schwartz $f$ with support in $t_*\geq 0$ in~\eqref{EqINDecay} are sharp at $\iota^+\cup\cT^+$ up to the $\eps$-loss. In fact in $\frac{r}{t_*}<C<\infty$, $u_0$ is equal to a leading order term $\rho_+^{\beta^++1}\rho_\cT^{\beta^+-\beta^-+1}a(t_*,x)$ plus a more decaying remainder, where the limits $a_\cT(x):=\lim_{t_*\to\infty}a(t_*,x)$ (with $x\in\R^n$) at $\cT^+$ and $a_+(X)=\lim_{t_*\to\infty}a(t_*,t_* X)$ (with $0\neq X\in\R^n$) at $\iota^+$ exist, agree at the corner $\iota^+\cap\cT^+$ (i.e.\ $\lim_{r\to\infty}a_\cT(r\omega)=\lim_{R\to 0}a_+(R\omega)$, $\omega\in\Sph^{n-1}$), and are not identically $0$.
\end{thm}

Here, $a_\cT$ is a smooth multiple of a non-trivial element of $\ker\wh{P_0}(0)$ with minimal growth at infinity, and $a_+$ can similarly be found by solving a certain PDE. See Theorem~\ref{ThmAS} for the precise statement, which includes uniform bounds near $\scri^+\cap\iota^+$. Since in the generality at which we work in this paper the nondegeneracy condition is often verified, we restrict ourselves here to the basic decay estimate in Theorem~\ref{ThmINDecay}, and do not study the problem of obtaining sharp bounds also in degenerate situations (with \cite{HintzPrice} being a special case).

\subsection{Nonstationary operators}
\label{SsIN}

We now turn to the class of non-stationary operators $P$ we shall investigate in the second part of this paper. We will explain the (microlocal) regularity theory for $P$ and how to combine it with decay estimates for its stationary model operator $P_0$ in order to get a solvability theory for $P$. For an illustration of some aspects of our approach in a simple ordinary differential equation (ODE) setting, we refer the reader to Appendix~\ref{SODE}.

With $g_0$ denoting a stationary metric as in~\S\ref{SsIS}, we consider non-stationary metrics
\[
  g = g_0 + \tilde g.
\]
In $r=|x|<\half t$, we require $|Z^J\tilde g_{\mu\nu}|\lesssim t_*^{-\delta}$ for the coefficients $\tilde g_{\mu\nu}$ of $\tilde g$ in the coordinates $z=(t,x)\in\R^{1+n}$. For $\frac{r}{t}\in(\half,2)$ and large $r$, i.e.\ near null infinity, the precise decay requirements for different components of $\tilde g$ differ. As a simple example, the decay $|Z^J\tilde g_{\mu\nu}|\lesssim r^{-1-\delta}t_*^{-\delta}$ in $t_*\geq 1$ is acceptable (see Example~\ref{ExGAGEx} for a discussion); but the ability to handle metric perturbations with weaker decay at $\scri^+$ and $\iota^+$ is important in applications. See Definition~\ref{DefGAG} and Lemmas~\ref{LemmaGAGPert}--\ref{LemmaGAGPertMet} for the precise assumptions, which in particular allow for the type of behavior arising in nonlinear stability problems in $3+1$ spacetime dimensions.

We then consider a class of non-stationary operators
\[
  P = P_0 + \tilde P
\]
which are wave operators with respect to the metric $g$ (i.e.\ the principal symbol of $P$ is scalar and equal to that of $\Box_g$); here, $P_0$ is a stationary operator with respect to $g_0$ as in~\S\ref{SsIS}, and $\tilde P$ has decaying coefficients (except for the possibility of certain leading order contributions at null infinity) relative to $P_0$ when written in terms of an appropriate frame of vector fields on $\R^{1+n}$. See Examples~\ref{ExGAWWave}--\ref{ExGAWFirst} for various classes of operators $P$ of interest in applications; this includes wave operators associated with non-stationary metrics such as $g$ and their couplings to potentials $V_0+\tilde V$ as (or even more general than) in Theorem~\ref{ThmIV} as well as first order perturbations. See Definition~\ref{DefGAW} (and also Lemma~\ref{LemmaGAWStatbetb}) for the general setting.

Roughly speaking, we require that the coefficients of $\tilde P$, relative to those of $P_0$, decay at null infinity and in $t_*=t-r>1$, $r>\frac14 t$ when expressing $P_0,\tilde P$ in terms of the vector fields\footnote{An equivalent collection of vector fields, for present purposes, is $t\pa_t+r\pa_r$, $r(\pa_t+\pa_r)$, $\sqrt{\frac{t-r}{t}}\pa_\omega$.}
\begin{equation}
\label{EqIeb}
  (t-r)(\pa_t-\pa_r),\qquad
  r(\pa_t+\pa_r),\qquad
  \sqrt{\frac{t-r}{t}}\pa_\omega
\end{equation}
(writing $\pa_\omega$ for a spherical vector field), which span the space of edge-b-vector fields on a compactification of $\R^{1+n}$ to a manifold with corners introduced in~\cite{HintzVasyScrieb}. (Notions of edge-b-geometry are recalled in~\S\ref{SsMeb}, following \cite{MelroseVasyWunschDiffraction,HintzVasyScrieb}.) Near the future translation face $\cT^+$, the relevant class of vector fields is given by the 3b-vector fields introduced in~\cite{Hintz3b}, which in the coordinates $(t_*,r,\omega)$ and away from $r=0$ are
\begin{equation}
\label{EqI3b}
  r\pa_{t_*},\qquad r\pa_r,\qquad \pa_\omega.
\end{equation}

The utility of edge-b-vector fields near $\scri^+$, resp.\ 3b-vector fields near $\cT^+$, is explained in~\cite{HintzVasyScrieb}, resp.\ \cite{Hintz3b}. Namely, up to an overall weight, $g$ has Lorentzian signature with respect to these vector fields \emph{uniformly} in $t_*\geq 1$, i.e.\ in all asymptotic regimes (so near $\scri^+$, $\iota^+$, and $\cT^+$); and correspondingly the principal symbol of $P$ is, up to an overall weight, a Lorentzian signature quadratic form in these vector fields. We illustrate this only in the case of the Minkowski metric $g_M=-\dd t^2+\dd r^2+r^2\slg$ here. The dual metric is $g_M^{-1}=-\pa_t^2+\pa_r^2+r^{-2}\slg^{-1}$; we schematically write $\slg^{-1}=\pa_\omega\otimes_s\pa_\omega$. In $r>\frac14 t$, $t-r>1$, we then note that
\begin{subequations}
\begin{equation}
\label{EqINeb}
  r(t-r) g_M^{-1} = -r(\pa_t+\pa_r)\otimes_s (t-r)(\pa_t-\pa_r) + \sqrt{\frac{t-r}{r}}\pa_\omega\otimes_s\sqrt{\frac{t-r}{r}}\pa_\omega
\end{equation}
is a Lorentzian signature quadratic form in the aforementioned edge-b-vector fields, uniformly as $r\to\infty$. In $t>0$, $1<r<\frac34 t$, the rescaling
\begin{equation}
\label{EqIN3b}
  r^2 g_M^{-1} = -(r\pa_t)^2 + (r\pa_r)^2 + \pa_\omega^2
\end{equation}
\end{subequations}
is a Lorentzian signature quadratic form in the 3b-vector fields~\eqref{EqI3b}.\footnote{A uniform description down to $r=0$ simply uses $\la x\ra$ instead of $r$, so $\la x\ra^2 g_M^{-1}=-(\la x\ra\pa_t)^2+(\la x\ra\pa_x)^2$.} In the overlap region where $\frac14 t<r<\frac34 t$, we have $r(t-r)\sim r^2$, i.e.\ the two weights in~\eqref{EqINeb} and \eqref{EqIN3b} are comparable.

\begin{rmk}[Other frames of vector fields]
\label{RmkINOther}
  The Minkowski dual metric $g_M^{-1}$ is a (weighted) Lorentzian signature quadratic form, uniformly in all asymptotic regimes, with respect to a variety of other frames of vector fields as well.
  \begin{enumerate}
  \item It is a Lorentzian signature quadratic form in the spacetime translation-invariant vector fields\footnote{These are \emph{scattering vector fields} on the radial compactification $\ol{\R^{1+n}}$ in the terminology of \cite{MelroseEuclideanSpectralTheory}, and their span with $\CI(\ol{\R^{1+n}})$-coefficients defines the Lie algebra $\Vsc(\ol{\R^{1+n}})$ of scattering vector fields.} $\pa_t$, $\pa_{x^j}$ ($j=1,\ldots,n$) uniformly as $|t|+|x|\to\infty$ (its coefficients being constant). However, estimates for $\Box_{g_M}$ that only use these vector fields (but not weighted versions thereof) are \emph{very weak}, and are not suitable for obtaining decay estimates. More subtly, $\Box_{g_M}$ is \emph{degenerate} (its scattering principal symbol vanishes quadratically at the 0-section over $\pa\ol{\R^{1+n}}$) and thus does not have good mapping properties between (weighted) Sobolev spaces defined with respect to these vector fields.
  \item Similarly, if $z=\sqrt{1+t^2+|x|^2}$, then $z^2 g_M^{-1}=-(z\pa_t)^2+(z\pa_x)^2$ is a Lorentzian signature quadratic form in the approximately scaling-invariant vector fields $z\pa_t$, $z\pa_{x^j}$.\footnote{These are \emph{b-vector fields} on the radial compactification $\ol{\R^{1+n}}$ in the terminology of \cite{MelroseAPS}, and their $\CI(\ol{\R^{1+n}})$-span is the Lie algebra $\Vb(\ol{\R^{1+n}})$ of b-vector fields.} This is the perspective of \cite{BaskinVasyWunschRadMink,BaskinVasyWunschRadMink2}. These vector fields are \emph{too strong}, both at null infinity and in spatially compact regions, if instead of $g_M$ one considers perturbations $g$ which are (asymptotically) stationary, or which have a delicate structure near null infinity: the coefficients of $g$ are not regular with respect to these strong vector fields (e.g.\ the derivative of a stationary spatially compactly compactly supported term $\chi(x)\dd t^2$ along $\la z\ra\pa_x$ loses a power of $\la z\ra$), and correspondingly neither are solutions of corresponding wave equations. However, this perspective \emph{is} adequate in the present paper in the intermediate region $q t<r<(1-q)t$, $q\in(0,\half)$, where it matches the edge-b- and 3b-perspectives.
  \end{enumerate}
  In a nutshell, one must strike a balance between using sufficiently weak vector fields (i.e.\ weighted versions of $\pa_t$, $\pa_x$ where the weights are not very strong at infinity) so that the coefficients of the wave operators are regular with respect to these, while one wants to use vector fields which are as strong as possible so that one obtains estimates on function spaces encoding strong weighted regularity. The edge-3b-perspective put forth in the present work satisfies both requirements.
\end{rmk}

Lastly, we remark that if $P$ acts on sections of a vector bundle, then we need to require its leading order part at null infinity to have a certain lower triangular structure (see Definition~\ref{DefGAW}\eqref{ItGAWEdgeN}) which is reminiscent of linearized versions of weak null conditions, as explained in \cite[Remark~6.3]{HintzVasyScrieb} and \cite[Remark~1.7]{HintzVasyMink4}.

\subsubsection{Regularity via microlocal analysis}
\label{SssIReg}

Here, we shall only discuss the microlocal perspective on analyzing the (propagation of) regularity of solutions of $P u=f$ in the timelike cone $r<\half t$. The microlocal analysis of operators near null infinity is the subject of \cite{HintzVasyScrieb}, to which we refer the reader for a detailed discussion. For the sake of simplicity, we restrict in this introduction to the case that $g=g_0$ is the Minkowski metric, further
\[
  P_0=\Box_{g_0}+V_0
\]
is the wave operator on $(n+1)$-dimensional Minkowski space with a stationary potential $V_0\in\cA^{2+\delta}(\ol{\R^n})$, so in the coordinates $(t_*,r,\omega)=(t-r,r,\omega)$
\begin{equation}
\label{EqI3bP0}
  r^2 P_0 = 2 r\pa_{t_*}\Bigl(r\pa_r + \frac{n-1}{2}\Bigr) + (r D_r)^2 - i(n-2)r D_r + \slDelta + r^2 V_0(r,\omega),
\end{equation}
and finally $\tilde P$ is a zeroth order term, i.e.\ $\tilde P=\tilde V$ is a (complex-valued) potential with $|r^2 Z^J\tilde P|\lesssim t_*^{-\delta}$ and support in $r<\half t$.

\begin{figure}[!ht]
\centering
\includegraphics{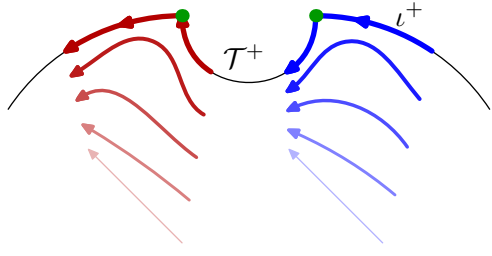}
\caption{Projections of null-bicharacteristics of $P$ in the 3b-cotangent bundle to the base manifold $M$---that is, null-geodesics (and limits of families of null-geodesics)---in $r<\half t$. The null-geodesics are $\gamma_{t_0}\colon s\mapsto (t,x)=(t_0+s,x_0+s v_0)$ where $x_0,v_0$ with $|v_0|<1$ are fixed. \textit{In blue:} incoming null-geodesics $\gamma_{t_0}(s)$ where $-\frac{t_0}{2}\leq s\leq -1$ and $t_0\nearrow\infty$. \textit{In red:} outgoing null-geodesics $\gamma_{t_0}(s)$ where $1\leq s\leq\frac{t_0}{2}$ and $t_0\nearrow\infty$. \textit{In green:} the projection to $M$ of the two radial sets over $\cT^+$. The incoming and outgoing radial sets are bundles of lines (with the origin removed) in $\Ttb^*M$ over $\cT^+$.}
\label{FigIBichar}
\end{figure}

Note then that $r^2 P=r^2 P_0+r^2\tilde P$ enjoys an approximate invariance under time translations $(t_*,r,\omega)\mapsto(t_*+s,r,\omega)$ (broken by the presence of $\tilde P$) and under spacetime scalings $(t_*,r,\omega)\mapsto(\lambda t_*,\lambda r,\omega)$ (broken by the presence of $V_0$ and $\tilde P$). This is precisely the setting for which the framework of 3b-analysis---microlocal analysis of (pseudo)differential operators built out of the vector fields~\eqref{EqI3b}---was developed by the author in \cite{Hintz3b}; see also~\S\ref{SsM3b}. We explain this in the present context by working in a coordinate chart
\begin{equation}
\label{EqI3bChart}
  [0,1)_{\rho_\cT}\times[0,1)_{\rho_+}\times\Sph^{n-1},\qquad
  \rho_\cT=\frac{r}{t_*},\qquad
  \rho_+=\frac{1}{r},
\end{equation}
on a compactification $M$ of $\R^{1+n}$ to a manifold with corners; this chart covers a neighborhood of the intersection of $\iota^+$ and $\cT^+\subset M$ (cf.\ Figure~\ref{FigIRhos}), with $\R^{1+n}$ locally given by $\rho_\cT,\rho_+>0$. Over this chart, we define the 3b-cotangent bundle $\Ttb^*M\to M$ to have as a smooth frame \emph{down to} $\iota^+=\rho_+^{-1}(0)$ and $\cT^+=\rho_\cT^{-1}(0)$ the 1-forms
\[
  \frac{\dd t_*}{r},\qquad \frac{\dd r}{r},\qquad \dd\omega
\]
dual to~\eqref{EqI3b}. By~\eqref{EqI3bP0} and the assumptions on $\tilde P$, the principal symbol of $r^2 P$ is a Lorentzian signature quadratic form on $\Ttb^*M$. The null-bicharacteristic flow of $P$, which here is the null-geodesic flow on Minkowski space lifted to phase space $T^*\R^{1+n}$, extends to the flow of a vector field $H$ on $\Ttb^*M$. As such, $H$ has radial sets, i.e.\ it is fiber-radial at two submanifolds of $\Ttb^*M$ which lie over the corner $\iota^+\cap\cT^+$. These radial sets are the links between the near field region ($r$ bounded, $t_*=\infty$, i.e.\ the interior of $\cT^+$) and the intermediate region ($\frac{r}{t}\in(\eta,1-\eta)$, $t_*=\infty$, i.e.\ the interior of $\iota^+$). Null-bicharacteristics over $t_*^{-1}(\infty)=\iota^+\cap\cT^+$ entering the near field region pass through an incoming radial set, and those leaving the near field region pass through the outgoing radial set; these null-bicharacteristics are limits of reparameterized lifts of families of null-geodesics translated by increasingly large amounts in $t_*$. See Figure~\ref{FigIBichar}.

Using standard positive commutator arguments, we then prove microlocal radial point estimates (see \cite{MelroseEuclideanSpectralTheory}, \cite[\S2]{VasyMicroKerrdS}, \cite[Appendix~E]{DyatlovZworskiBook}) for the non-stationary operator $P$ in weighted 3b-Sobolev spaces, which are $L^2$-based weighted spacetime Sobolev spaces in which regularity is measured with respect to the vector fields~\eqref{EqI3b}; see~\S\ref{SsWT}. Since we work in nontrapping geometries, microlocal elliptic and real principal type propagation estimates (to deduce regularity along outgoing null-bicharacteristics from regularity along incoming null-bicharacteristics) suffice for all other parts of phase space in $r<\half t$.

Altogether, we can therefore propagate microlocal regularity of $u$ solving $P u=f$ passing from $r>\half t$ into the timelike cone $r<\half t$ through this entire cone and back out again. That is, we obtain a schematic estimate
\begin{equation}
\label{EqI3bNearT}
  \|\chi u\|\leq\|E u\|+\|\chi P u\|,
\end{equation}
where $\chi$ localizes to $r<\half t$, and $E$ is a pseudodifferential operator which is elliptic at (and localized near) incoming momenta in the region $\frac{r}{t}\in(\frac14,\frac34)$.\footnote{More precisely, $\|\chi u\|_{\Htb^{s,(\alpha_+,\alpha_\cT)}}\lesssim\|E u\|_{\Htb^{s,(\alpha_+,\alpha_\cT)}}+\|\tilde\chi P u\|_{\Htb^{s-1,(\alpha_++2,\alpha_\cT)}}+\|\tilde\chi u\|_{\Htb^{s_0,(\alpha_+,\alpha_\cT)}}$, where $E\in\Psitb^0$ is elliptic on the stated set. The orders $s_0<s$ and $\alpha_+,\alpha_\cT$ are subject to threshold conditions at the incoming and outgoing radial sets over $\pa\cT^+$; and $\tilde\chi$ equals $1$ near $\supp\chi$. In the main part of the paper, we do not state such an `intermediate' regularity estimate, and instead combine the microlocal elliptic, real principal type, and radial point estimates, including those near null infinity discussed below, all at once.} In this manner, the 3b-perspective handles scattering off the asymptotically stationary potential $V_0+\tilde V$ both geometrically and analytically.

\begin{rmk}[Variable orders]
\label{RmkI3bVar}
  The microlocal 3b-estimates only use the principal symbol of $P$ in $\Tb^*M$ and the corresponding null-bicharacteristic flow, and its subprincipal symbol at the radial sets. No spectral assumptions enter yet; thus, these estimates in principle allow for the possibility e.g.\ of embedded eigenvalues $0\neq\sigma_0\in\R$ leading to $e^{-i\sigma_0 t_*}$ behavior of $u$. Since differentiation of this along $r\pa_{t_*}$ increases growth in $r$, there is a definite upper bound on the amount of 3b-regularity (relative to an $L^2$-space with fixed weights) we can control in the outgoing directions using these estimates. On the flipside, propagation \emph{into} the near field region imposes a lower bound on 3b-regularity, roughly corresponding to requiring the absence of incoming spherical waves for resolvent estimates for the stationary operator. To accommodate both threshold requirements at the same time, we need to use 3b-Sobolev spaces with suitable variable differential orders.\footnote{This is already the case for the wave operator on Minkowski space, where it is closely related to the need for variable decay orders in sharp estimates for stationary scattering as in \cite[Proposition~4.13]{VasyMinicourse}.} Roughly speaking, the 3b-regularity order can be arbitrarily high, except in the region $1\ll r\leq t=t_*+r$ near outgoing momenta $\pa_r^\flat\approx\dd t_*$ (in $t_*,r$ coordinates) where one needs to impose an upper bound. Such Sobolev spaces, and corresponding classes of pseudodifferential operators, go back to work of Unterberger \cite{UnterbergerVariable} and feature also in \cite{BaskinVasyWunschRadMink} (albeit for different reasons). See Lemma~\ref{LemmaWROrder}.
\end{rmk}

In order to control regularity of $u$ \emph{globally} in $t_*\geq 0$, we combine the estimate~\eqref{EqI3bNearT} with microlocal estimates in the remaining regions of spacetimes in the usual modular microlocal fashion. Concretely, one needs to control the term $E u$; this is feasible since backwards null-geodesics starting in the region where $E$ is localized remain outside of $r<\frac14 t$, i.e.\ control of $E u$ follows from the propagation of regularity in $r>\frac14 t$ (which in particular includes a neighborhood of null infinity) on an appropriate scale of Sobolev spaces. The microlocal framework near null infinity, based on edge-b-analysis, is the subject of \cite{HintzVasyScrieb}. Upon concatenating all radial point and real principal type propagation estimates, one can propagate regularity of $u$ in $t_*<0$ (where we require $u$ to vanish when $f$ does) along future null-bicharacteristics and through radial sets to obtain global control; see Figure~\ref{FigWFlow}. Schematically, the resulting estimate is\footnote{The factor of $2$ in the weight of the Sobolev spaces here is introduced as a matter of convention: in the main part of the paper, we use powers of $x_{\!\scri}=\rho_{\!\scri}^{1/2}$ for weights at $\scri^+$, as $x_{\!\scri}$ is the correct defining function of $\scri^+$ from the geometric (edge-b-)perspective. The number $\alpha_{\!\scri}$ refers to decay of $u$ as measured (for bounded $|t-r|$) in powers of $r^{-1}$.}
\begin{equation}
\label{EqI3bRegEst}
  \|u\|_{H_{\etbop}^{s,(2\alpha_{\!\scri},\alpha_+,\alpha_\cT)}} \lesssim \|P u\|_{H_{\etbop}^{s-1,(2(\alpha_{\!\scri}+1),\alpha_++2,\alpha_\cT)}} + \|u\|_{H_{\etbop}^{s_0,(2\alpha_{\!\scri},\alpha_+,\alpha_\cT)}},
\end{equation}
where the norm on the left is the $L^2$-norm of $\rho_{\!\scri}^{-\alpha_{\!\scri}}\rho_+^{-\alpha_+}\rho_\cT^{-\alpha_\cT}u$ and its derivatives of order $s$ along edge-b-vector fields in $\frac14 t\leq r<2 t$ and 3b-vector fields in $r<\frac34 t$ (ignoring the technical issue of having to work with variable, thus microlocally defined, Sobolev spaces), and the orders $s,\alpha_{\!\scri},\alpha_+,\alpha_\cT$ are subject to certain threshold conditions required by the radial point estimates;\footnote{If $P$ acts on sections of a vector bundle, the precise form of these threshold conditions depends on a choice of fiber inner product. Due to the (modular and) microlocal character of the radial point estimates, one can optimize the choice of fiber inner product individually near each radial set. We shall use Proposition~\ref{PropLA} to find inner products for which the threshold conditions are as weak as possible.} lastly, $s_0<s$. See Proposition~\ref{PropWR} for the precise statement. Importantly for the solvability theory for $P$, we also prove an estimate for solutions of the adjoint equation $P^*\tilde u=\tilde f$ on the dual function spaces.

\begin{rmk}[Edge-3b-regularity vs.\ b-regularity]
\label{RmkI3bRegCompare}
  In $t_*\geq 1$, $r>\frac14 t$, one can test for b-regularity (i.e.\ regularity with respect to $Z$ in~\eqref{EqIZ}) by differentiating along the vector fields $(t-r)(\pa_t-\pa_r)$, $r(\pa_t+\pa_r)$, $\pa_\omega$, which are stronger than the edge-b-vector fields in~\eqref{EqIeb} in that the spherical derivatives do not come with the decaying (at $\scri^+$) weight $\sqrt{\frac{t-r}{t}}$. Testing for b-regularity in $1<r<\frac34 t$ can be accomplished using the vector fields $t_*\pa_{t_*}$, $r\pa_r$, $\pa_\omega$, which are stronger than the 3b-vector fields in~\eqref{EqI3b} in that the 3b-vector field $r\pa_{t_*}=\frac{r}{t_*}t_*\pa_{t_*}$ has a vanishing (at $\cT^+$) factor $\frac{r}{t_*}$ in front compared to $t_*\pa_{t_*}$.
\end{rmk}

These (a priori) regularity estimates are rather permissive, as they are purely symbolic and do not rely on any spectral assumptions, such as mode stability, on the stationary model $P_0$, and indeed they do not provide any control on the \emph{decay} of $u$. This is where the inverses of model operators of $P$ (and thus the spectral assumptions on $P_0$) enter.

Moreover, there is an upper bound on the amount of regularity which these propagation results can control, cf.\ Remark~\ref{RmkI3bVar}; one can remove this upper bound only after one has exploited the spectral assumptions on $P_0$ (see below).

\subsubsection{Decay via inversion of model operators}

Corresponding to the three decay orders of the Sobolev spaces appearing in~\eqref{EqI3bRegEst}, there are three model operators of the edge-3b-operator $P$. We first discuss the model operator at the future translation face $\cT^+$---which is precisely the stationary model $P_0$. An important aspect of the (variable order) 3b-spaces is that (for a special value of $\alpha_\cT$) under the Fourier transform in $t_*$ they are isometric to $L^2$-spaces in the spectral parameter $\sigma$ with values in \emph{precisely} those function spaces on which the spectral theory of $\wh{P_0}(\sigma)$ takes place, i.e.\ scattering Sobolev spaces with variable scattering decay order for $\sigma\neq 0$, b-Sobolev spaces for $\sigma=0$, and scattering-b-transition Sobolev spaces (see \cite{GuillarmouHassellResI}, \cite[\S2.4]{Hintz3b}) for uniformity near $\sigma=0$. \emph{This perfect match lends strong support to the 3b-perspective near $\cT^+$, and is indeed our motivation to introduce and use it in this paper for the analysis of waves on non-stationary spacetimes.} See also \cite[\S\S1.1, 4.4]{Hintz3b}. By replacing $P$ with its $\cT^+$-model operator $P_0$ and using the invertibility of the latter, one can improve the order $\alpha_\cT$ of the final, error, term in~\eqref{EqI3bRegEst} to $\alpha'_\cT<\alpha_\cT$.

Likewise, by inverting the normal operator at $\scri^+$, one can improve the order $\alpha_{\!\scri}$ of the error term (see \cite[Definition~3.5, \S7]{HintzVasyScrieb} and Definition~\ref{DefGAW}\eqref{ItGAWEdgeN} for definitions and proofs). The final model operator of $r^2 P$ is its model at $\iota^+$ which is exactly homogeneous under spacetime dilations. Using the coordinates $R=\frac{r}{t_*}$, $\omega\in\Sph^{n-1}$, note that punctured future timelike infinity $\iota^+$ is locally near $r/t_*=0$ diffeomorphic to $[0,1)_R\times\Sph^{n-1}$, i.e.\ the origin $R=0$ of polar coordinates is resolved; this corresponds to the fact that the $\iota^+$-normal operator has a conic singularity at $R=0$. In the present special case, this normal operator is equal to $r^2\Box_{g_0}$, which in the coordinates $T=t_*^{-1}$, $R$, $\omega$ is given by
\begin{equation}
\label{EqI3bIplusNorm}
  R^2 \Bigl( D_R^2 - \frac{i(n-1)}{R}D_R + R^{-2}\slDelta - 2 R^{-1}(T\pa_T+R\pa_R)(R \pa_R + n-1) \Bigr)
\end{equation}
(cf.\ \eqref{EqI3bP0}). The Mellin-transformed normal operator in $T$ lies in the algebra of b-operators near $R=0$, and in the high frequency regime in the algebra of semiclassical cone operators \cite{HintzConicPowers,HintzConicProp}. This is discussed in the general 3b-setting in \cite[\S\S1.1, 3.3]{Hintz3b}. Under the Mellin-transform, the function spaces in~\eqref{EqI3bRegEst} transform to Sobolev spaces corresponding to these algebras \cite[\S4.4]{Hintz3b}. Analytically, the conic singularity necessitates allowing for a non-trivial relative decay order $\alpha_\cT-\alpha_+$ at $\cT^+$ and $\iota^+$. (To connect this with the discussion in~\S\ref{SsIL}, we note that if the potential $V_0$ included an inverse square term $c r^{-2}$ in $r\gg 1$, there would be an extra term $R^2\cdot c R^{-2}$ in~\eqref{EqI3bIplusNorm}, i.e.\ an inverse square singularity at $R=0$.) For the class of non-stationary operators considered here, the Mellin-transformed $\iota^+$-normal operator family is equal to the action of the transition face model operators $N_\tface^\pm(P_0)$ (acting on suitable positive and negative frequency parts) under a sequence of transformations (see~\S\ref{SsWip}).\footnote{Conceptually, notice the similarity of the transition from zero to nonzero energies and the decay rates of the Green's functions on the one hand, and the transition from the near-field (where the zero energy operator dominates the asymptotic behavior) to the far-field (null infinity, where the decay is inherited from outgoing spherical waves) on the other hand.} This explains why, in contrast to \cite{BaskinVasyWunschRadMink}, we do not need to make separate assumptions on the absence of resonances at $\iota^+$ in suitable half spaces of the Mellin-dual variable for our decay results for non-stationary operators: they are automatic.\footnote{For stationary wave type operators, spectral methods such as those in \cite{HintzPrice} or~\S\ref{SSt} suffice for obtaining sharp decay results as $t_*\to\infty$ (or more precisely at $\iota^+\cup\cT^+$). It is thus possible that in the present work one could dispense of the analysis of the $\iota^+$-normal operator altogether if one uses the stationary model $P_0$ as a good approximation of $P$ globally near $\iota^+\cup\cT^+$. However, it is not clear how to prove estimates for $P_0$ on the same spaces on which the (microlocal) analysis of $P$ takes place without dealing with $\iota^+$ separately after all.}

\begin{rmk}[Hierarchies of model operators]
\label{RmkI3bModel}
  The non-stationary operator $P$ does not possess \emph{any} global symmetries; but it does possess approximate symmetries in the three aforementioned asymptotic regimes: approximate time-translation invariance near $\cT^+$, approximate degree $-2$ homogeneity under spacetime dilations near $\iota^+$, and a more subtle approximate scaling invariance near $\scri^+$ described in \cite[Remark~1.4]{HintzVasyScrieb}. The stationary model $P_0$ is the \emph{exactly} time-translation invariant operator that $P$ is equal to at $\cT^+$ to leading order. (In this sense, the stationary analysis of~\S\S\ref{SsIS} and \ref{SSt}, which builds on \cite{HintzPrice}, merely concerns the inversion of this model operator.) Exploiting the time-translation invariance of $P_0$ by passing to its spectral family, one obtains a family $\sigma\mapsto\wh{P_0}(\sigma)$ of operators acting (parametrically) on a lower-dimensional space; this family in turn has approximate invariances in various regimes (e.g.\ under scalings $(\sigma,\rho)\mapsto(\lambda\sigma,\lambda\rho)$, $\rho=r^{-1}$, near $(\sigma,\rho)=(0,0)$) in which one can then in turn define \emph{exactly} invariant model operators (e.g.\ the transition face normal operator); and so on. Analogous considerations apply to the exactly dilation-homogeneous model operator at $\iota^+$. Figure~\ref{FigWMFred} summarizes some of the model operators involved, and their interrelationships.
\end{rmk}

With all three decay orders in~\eqref{EqI3bRegEst} improved, and upon proving analogous estimates for the adjoint operator $P^*$ on dual function spaces, we conclude that $P$ acts as a Fredholm operator between spacetime 3b-Sobolev spaces with suitable weights and differentiability orders. While $P$, as a wave operator, is clearly injective on functions supported in $t_*\geq 0$, its surjectivity, i.e.\ the solvability of $P u=f$ (for $f$ Schwartz, say) \emph{in polynomially weighted spaces} is not obvious at all. We prove the surjectivity of $P$ indirectly by showing the injectivity of $P^*$ on the dual function spaces; we do this in two steps.
\begin{enumerate}
\item First, we apply the previous arguments to the stationary model $P_0$. Since by Theorem~\ref{ThmINDecay} Schwartz forcing terms are in the range of $P_0$ acting on a space of polynomially bounded functions, and since the range of $P_0$ is closed in a weighted 3b-Sobolev space, this suffices to obtain the surjectivity of $P_0$ acting between these spaces.
\item Rephrasing this as a quantitative injectivity statement for $P_0^*$, i.e.\ an estimate
  \[
    \|v\|_{H_\etbop^{-s+1,(-2(\alpha_{\!\scri}+1),-\alpha_+-2,-\alpha_\cT)}}\lesssim\|P_0^* v\|_{H_\etbop^{-s,(-2\alpha_{\!\scri},-\alpha_+,-\alpha_\cT)}},
  \]
  we then prove the injectivity of $P^*$ by upgrading the Fredholm estimate for $P^*$ via comparing $P^*$ with $P_0^*$ and absorbing the error $P^*-P_0^*$, which is small for late times.
\end{enumerate}
(We remind the reader of Appendix~\ref{SODE}, where this procedure is carried out in a simple ODE setting.) See Theorem~\ref{ThmWM} for the resulting detailed statement on the solvability of $P u=f$.

\begin{rmk}[Solvability]
\label{RmkI3bSolv}
  Here as well as in the de~Sitter type settings considered in \cite{HintzVasySemilinear,HintzQuasilinearDS,HintzVasyQuasilinearKdS,HintzVasyKdSStability,HintzKNdSStability}, the existence of forward solutions of $P u=f$ obeying some \emph{exponential} upper bound follows from a simple energy estimate. Unlike in the de~Sitter type settings, however, it is impossible to upgrade such an exponential bound to polynomial bounds solely by using the mapping properties of $P_0$, unless the non-stationary perturbation $\tilde P=P-P_0$ of the stationary operator $P_0$ has exponentially decaying coefficients---an assumption which is much too restrictive in asymptotically flat settings.
\end{rmk}

Higher b-regularity of solutions (that is, regularity under repeated application of elements of $Z$ in~\eqref{EqIZ}) is proved by commuting suitable vector fields related to $Z$ through $P$; see~\eqref{EqWbVF}. It is important here that there is a sufficient supply of vector fields which, at punctured future timelike infinity $\iota^+$, approximately commute with the principal part of $P$; this is related to Klainerman's vector field method \cite{KlainermanUniformDecay}, though the requirements on the vector fields are significantly less stringent. In this manner, one can trade edge-3b-derivatives for b-derivatives (which are stronger); one can then regain edge-3b-regularity (while retaining the b-regularity improvement) using a microlocal regularity theory on Sobolev spaces which mix microlocal edge-3b-regularity with integer order b-regularity, akin to similar mixed spaces employed for instance in \cite{VasyPropagationCorners} or implicitly in \cite[\S12]{MelroseEuclideanSpectralTheory}. In this argument, it is crucial that we can solve (wave type equations related to) $P u=f$ on polynomially weighted spaces with sharp control on decay. Therefore, the ability to obtain higher b-regularity is directly tied to the validity of the spectral assumptions (in particular: mode stability) of the stationary model $P_0$ (cf.\ the discussion of embedded eigenvalues in Remark~\ref{RmkI3bVar}). Crucially, the amount of b-regularity one can prove in this manner is \emph{unlimited}, which thus allows one to eliminate the upper bound on regularity in the outgoing part of phase space which was previously imposed by the symbolic regularity theory (see again Remark~\ref{RmkI3bVar}). See Theorem~\ref{ThmWb}.

 A rough version of Theorems~\ref{ThmWM} and \ref{ThmWb} for the general class of wave type operators $P$ we study in the main part of the paper reads as follows (using the notation of~\eqref{EqI3bRegEst}, cf.\ also the estimate \eqref{EqIEst}):

\begin{thm}[Solvability and decay for non-stationary wave type equations]
\label{ThmIWM}
  Denote by $\ubar p_1\in\R$ a certain quantity controlling decay at $\scri^+$.\footnote{This quantity can be read off from the coefficients of the wave type operator $P$, and is essentially equal to a possible shift of $\frac{n-1}{2}$ in~\eqref{EqI3bP0}; see equation~\eqref{EqWRubarp1} in the notation of Definition~\usref{DefGAW}\eqref{ItGAWEdgeN} for the precise definition.} Assume that the stationary model $P_0$ of $P$ satisfies all spectral assumptions (mode stability, no zero energy resonances, invertibility of transition face model operators). Let $(\beta^-,\beta^+)$ denote the indicial gap of $P_0$. Let $\alpha_{\!\scri},\alpha_+,\alpha_\cT\in\R$ with $\alpha_{\!\scri}<-\half+\ubar p_1$, $\alpha_+<\min(-\half+\alpha_{\!\scri},-\frac{n-1}{2}+\beta^+)$, and $\alpha_\cT\in(\alpha_++\frac{n}{2}-\beta^+,\alpha_++\frac{n}{2}-\beta^-)$. Let $f\in H_\etbop^{s-1,(2(\alpha_{\!\scri}+1),\alpha_++2,\alpha_\cT)}$ be supported in $t_*\geq 0$, where $s$ is a suitable differential order. Then the forward solution of
  \[
    P u=f
  \]
  satisfies $u\in H_\etbop^{s,(2\alpha_{\!\scri},\alpha_+,\alpha_\cT)}$. If $f$ in addition enjoys $k$ degrees of b-regularity (i.e.\ up to $k$ b-derivatives of $f$ remain in the stated space), then $u$ also enjoys $k$ additional degrees of b-regularity.
\end{thm}

Infinite regularity $L^2$-based estimates for forward solutions, and pointwise estimates which follow from Sobolev embedding, are immediate consequences; see Corollaries~\ref{CorWbCon} and \ref{CorWbPointwise}.

\subsection{Outline}
\label{SsIO}

The regularity theory for non-stationary wave equations developed in this paper utilizes a sizeable arsenal of pseudodifferential algebras: edge- and b-algebras (and their combination) near null infinity, the 3b-algebra near the translation face, and algebras related to these two via normal operator maps and associated Fourier- or Mellin-transformed spectral or normal operator families---including b-, scattering, 0-, and semiclassical cone algebras. These algebras and their relationships are recalled in~\S\ref{SM}.

In~\S\ref{SG}, we describe in detail the metrics and operators to be studied in the remainder of the paper, both in the stationary (\S\ref{SsGS}) and in the non-stationary case (\S\ref{SsGA}).

The focus of~\S\ref{SSt} is on spectral theory and forward solutions of stationary operators. In~\S\ref{SsStEst}, we prove resolvent estimates for the spectral family which cover all spectral parameters $\sigma\in\C$ with $\Im\sigma\geq 0$. We assemble these in~\S\ref{SsStCo} to deduce pointwise decay estimates for waves (cf.\ Theorem~\ref{ThmIV}\eqref{ItIVSchwartz}) via the conormal regularity of the resolvent and a partial expansion at $\sigma=0$. In~\S\ref{SAS}, we obtain leading order asymptotic profiles of solutions of stationary wave type operators subject to a nondegeneracy condition (cf.\ Theorem~\ref{ThmIV}\eqref{ItIVStat}).

In~\S\ref{SW}, we turn to the analysis of non-stationary operators. In~\S\S\ref{SsWT}--\ref{SsWR}, we develop the microlocal regularity theory based on 3b-analysis near $\cT^+$ and edge-b-analysis (from \cite{HintzVasyScrieb}) near $\scri^+$; in~\S\ref{SsWip}, we study the $\iota^+$-normal operator and its relationship to the transition face model operators in the low energy spectral theory. In~\S\ref{SsWM}, we combine these estimates with the basic decay result for the stationary model operator from~\S\ref{SsStCo} to prove the solvability of non-stationary equations with precise mapping properties on polynomially weighted edge-3b-Sobolev spaces. In~\S\ref{SsWb}, we prove higher b-regularity in the non-stationary setting (cf.\ Theorem~\ref{ThmIV}\eqref{ItIVL2}) and deduce pointwise decay estimates (cf.\ Theorem~\ref{ThmIV}\eqref{ItIVLinfty}).

Simple concrete examples for the main results of this paper---Theorems~\ref{ThmStCo}, \ref{ThmAS}, \ref{ThmWM}, and \ref{ThmWb}, and Corollaries~\ref{CorWbCon} and \ref{CorWbPointwise}---are given in~\S\ref{SE}.

The reader only interested in the theory of stationary wave equations (with rapidly decaying forcing or initial data) may skip~\S\ref{SW} entirely; moreover, the microlocal background required for~\S\ref{SSt} does not involve edge-b-, 0-, or 3b-analysis, and indeed is largely restricted to the b-algebra (though in some proofs we do make use of a second microlocal algebra, following~\cite{VasyLowEnergy,VasyLowEnergyLag}).

Readers interested in the non-stationary theory may take Theorem~\ref{ThmStCo} as a black box and proceed to~\S\ref{SW}. Only the estimates for the $\cT^+$-normal operator, i.e.\ the stationary model operator, in the proof of spacetime estimates on weighted 3b-Sobolev spaces rely on the estimates on the spectral family proved in~\S\ref{SsStEst}. (We mention here that many statements in~\S\ref{SsStEst} have two versions, e.g.\ Proposition~\ref{PropStEstNz}: part~\eqref{ItStEstNzsc} takes place on variable order spaces and is used in the variable order 3b-regularity theory in~\S\ref{SW}, whereas part~\eqref{ItStEstNzb} on b-Sobolev spaces is only used in the proof of Theorem~\ref{ThmStCo}.)

In Appendix~\ref{SODE}, we illustrate the interplay of estimates for non-stationary operators, control of decay for stationary models via the Fourier transform, and functional analytic techniques for the solution of non-stationary equations in the simple setting of ODEs on the positive half line. Appendix~\ref{SLA} proves a linear algebra lemma on near-optimal fiber inner products on vector bundles. Appendix~\ref{SWROrder} contains the proof of a technical result from~\S\ref{SsWR}. Appendix~\ref{SF} gives a pictorial summary of the function spaces used in the paper.

\subsection*{Acknowledgments}

I am very grateful to Andr\'as Vasy for countless discussions over the years on wave equations on asymptotically flat spacetimes which strongly influenced many parts of this paper. I also wish to express my gratitude to Shi-Zhuo Looi for many helpful comments and suggestions. Part of this research was conducted during the time I served as a Clay Research Fellow. Furthermore, I gratefully acknowledge support from the U.S.\ National Science Foundation (NSF) under Grant No.\ DMS-195514 and from a Sloan Research Fellowship. This material is based upon work supported by the NSF under Grant No.\ DMS-1440140 while I was in residence at MSRI in Berkeley, California, during the Fall 2019 semester.

\section{Microlocal toolkit}
\label{SM}

We recall a variety of Lie algebras of vector fields on manifolds with (fibered) boundaries which play a role in the present paper. Many of these algebras are well-known; accounts for b-analysis are given in \cite{MelroseAPS,GrieserBasics,MelroseDiffOnMwc}, and the 0-, edge, semiclassical scattering, b-edge, and scattering-b-transition algebras are discussed in the original papers \cite{MazzeoMelroseHyp,MazzeoEdge,VasyZworskiScl,MelroseVasyWunschEdge,GuillarmouHassellResI}. The 3b- and semiclassical cone algebras were defined in \cite{Hintz3b,HintzConicPowers,HintzConicProp}. For the present paper, the detailed presentations in \cite[\S2]{HintzVasyScrieb} and \cite[\S2]{Hintz3b} are particularly relevant. The main novel material in this section concerns commutator b-vector fields in the 3b-algebra, discussed in~\S\S\ref{SssM3bDil}--\ref{SssM3bb}.

We denote by $X$ a smooth manifold with boundary and by $M$ a smooth manifold with corners; all boundary hypersurfaces are assumed to be embedded. We write $x\in\CI(X)$ for a boundary defining function, and $x,z\in\CI(M)$ for defining functions of two intersecting boundary hypersurfaces $H_1,H_2$ of $M$; local coordinates in $\pa X$ or in the interior of $Y:=H_1\cap H_2$ are denoted $y\in\R^m$ (where $m=\dim X-1$, resp.\ $m=\dim M-2$). (We shall omit explicit descriptions near codimension $\geq 3$ corners of $M$.) We identify a collar neighborhood of $\pa X$, resp.\ $H_1\cap H_2$, with $[0,1)_x\times\pa X$, resp.\ $[0,1)_x\times[0,1)_z\times Y$. The space of smooth vector fields on $M$ is denoted $\cV(M)=\CI(M;T M)$, and $\Vb(M)\subset\cV(M)$ is the Lie algebra of \emph{b-vector fields} (i.e.\ vector fields tangent to the boundary); thus, $\Vb(M)$ is the space of smooth sections of the b-tangent bundle $\Tb M\to M$ with local frame $x\pa_x$, $z\pa_z$, $\pa_y$. Similarly, $\Vb(X)$ is spanned over $\CI(X)$ by $x\pa_x$, $\pa_y$. Restriction to a hypersurface $H\subset M$ defines a surjective map $\Vb(M)\to\Vb(H)$ and correspondingly a surjective bundle map $\Tb_H M\to\Tb H$. The space of b-differential operators is denoted $\Diffb^m(M)$ (locally finite sums of up to $m$-fold compositions of b-vector fields). If $w\in\CI(M^\circ)$ is a product of real powers of boundary defining functions, we write $w^{-1}\Diffb^m(M)=\{w^{-1}A\colon A\in\Diffb^m(M)\}$. We make analogous definitions on $X$, and write $\Diffb^{m,\alpha}(X)=x^{-\alpha}\Diffb^m(X)$. If $M$ has $N$ boundary hypersurfaces $H_1,\ldots,H_N$ and we assign a weight $\alpha_1,\ldots,\alpha_N\in\R$ to each of them, we shall also write $\Diffb^{m,(\alpha_1,\ldots,\alpha_N)}(M)=w^{-1}\Diffb^m(M)$ where $w=\prod_{i=1}^N\rho_i^{\alpha_i}$, with $\rho_i\in\CI(M)$ a defining function of $H_i$ for $i=1,\ldots,N$. We moreover write
\[
  \cA^\alpha(M) = \{ u\in\CI(M^\circ) \colon w^{-1}P u\in L^\infty_\loc(M)\ \forall P\in\Diffb(M) \}
\]
for spaces of ($L^\infty$-)\emph{conormal functions}; we put $\cA(M)=\bigcup_{\alpha\in\R^N}\cA^\alpha(M)$; we allow for $\alpha_i=\infty$ for $i$ in a subset $I\subset\{1,\ldots,N\}$, in which case $\cA^\alpha(M)$ is the intersection of all $\cA^\beta(M)$ where $\beta_j=\alpha_j$ for $j\notin I$, while $\beta_i$ is arbitrary for $i\in I$. (Thus, for example, $\cA^{(\infty,\ldots,\infty)}(M)=\CIdot(M)$ is the space of smooth functions vanishing to infinite order at all boundary hypersurfaces.) Generalizing $\Diffb^{m,\alpha}(M)$, we have the space $\cA^{-\alpha}\Diffb^m(M)$ of operators with conormal coefficients, consisting of locally finite sums of operators of the form $a P$ where $a\in\cA^{-\alpha}(M)$ and $P\in\Diffb^m(M)$.

We shall also consider spaces with (partial) expansions. Thus, for $\alpha<\beta\in\R$ we write
\[
  \cA^{(\alpha,0)}(X) = x^\alpha\CI(X),\qquad \cA^{(\alpha,0),\beta}(X) := x^\alpha\CI(X) + \cA^\beta(X).
\]
On $M$, given weights $\alpha_i<\beta_i$ for $i=1,\ldots,N$, we define
\begin{equation}
\label{EqMphgcon}
  \cA^{(((\alpha_1,0),\beta_1),\ldots,((\alpha_N,0),\beta_N))}(M)
\end{equation}
by induction over the dimension of $M$ as the space of smooth functions on $M$ which in a collar neighborhood $[0,1)_{\rho_i}\times H_i$ of $H_i$ are the sum of an element, valued in the space $\cA^{(((\alpha_{I(i,1)},0),\beta_{I(i,1)}),\ldots,((\alpha_{I(i,N(i))},0),\beta_{I(i,N(i))}))}(H_i)$, of $\rho_i^{\alpha_i}\CI([0,1))$ and an element of $\cA^{\beta_i}([0,1))$, where $H_{I(i,j)}$, $1\leq j\leq N(i)$, are the boundary hypersurfaces which intersect $H_i$ non-trivially. When $\alpha_i$, resp.\ $\beta_i$ is omitted, this means that the term $\rho_i^{\alpha_i}\CI([0,1))$, resp.\ $\cA^{\beta_i}([0,1))$ is absent; thus, for example,
\[
  \cA^{(((1,0),2),(3,0)}([0,1)_{\rho_1}\times[0,1)_{\rho_2}) = \rho_1^1\rho_2^3\CI([0,1)\times[0,1)) + \cA^2\bigl([0,1)_{\rho_1};\rho_2^3\CI([0,1)_{\rho_2})\bigr).
\]

One can define spaces $\Diffb^m(M;E):=\Diffb^m(M;E,E)$ and $\Diffb^m(M;E,F)$ of operators acting on sections of smooth vector bundles $E,F\to M$ in the usual manner, and one can then also define $\cA^\alpha(M;E)$ (which are tuples of elements of $\cA^\alpha(M)$ in local trivializations of $E$) etc. We shall not explicitly state generalizations to bundles for the other algebras of operators recalled below.

For $P\in\Diffb^m(X)$, and in the collar neighborhood $[0,1)_x\times\pa X$, we can write $P=\sum_{j=0}^m P_j(x,y,D_y)(x D_x)^j$, where $P_j$ is a smooth family (in $x$) of elements of $\Diff^{m-j}(\pa X)$. The \emph{b-normal operator} of $P$ (also called the normal operator of $P$ at $\pa X$) is defined as $N(P)=\sum_{j=0}^m P_j(0,y,D_y)(x D_x)^j$, and the \emph{Mellin-transformed normal operator family} is
\[
  \wh N(P,\lambda) = \sum_{j=0}^m P_j(0,y,D_y)\lambda^j,\qquad \lambda\in\C.
\]
When $P$ is a b-differential operator on $M$, it is convenient to indicate the relevant boundary hypersurface $H\subset M$ as a subscript; we thus write $N_H(P)$ and $\wh{N_H}(P,\lambda)$ for the normal operator and the Mellin-transformed normal operator family in this case. Returning to the manifold with boundary $X$, suppose that $P\in\Diffb^m(X)$ is elliptic, i.e.\ its principal symbol $\sigmab^m(P)\in(S^m/S^{m-1})(\Tb^*X)$ is invertible; here $\sigmab^m(P)$ is in fact a homogeneous polynomial of degree $m$ in the fibers, with value at $\lambda\frac{\dd x}{x}+\eta$ given in terms of the principal symbols of the $P_j$ by $\sum_{j=0}^m \upsigma^{m-j}(P_j)(0,y,\eta)\lambda^j$. Then $\wh N(P,\lambda)\in\Diff^m(\pa X)$ is a holomorphic family (in $\lambda\in\C$) of elliptic operators on $\pa X$; we define the \emph{boundary spectrum} of $P$ as
\begin{equation}
\label{EqMspecb}
  \specb(P) := \{ i\lambda \colon \wh N(P,\lambda)\ \text{is not invertible on}\ \CI(\pa X) \}.
\end{equation}
This is a discrete subset of $\C$.

Continuing the discussion of Lie algebras of vector fields on $X$, we write $\cV_0(X)=x\cV(X)$ and $\Vsc(X)=x\Vb(X)$ for the spaces of \emph{0-} and \emph{scattering vector fields}, respectively; the corresponding tangent bundles are ${}^0 T X$, with local frame
\[
  x\pa_x,\qquad x\pa_y,
\]
and $\Tsc X$, with local frame
\[
  x^2\pa_x,\qquad x\pa_y.
\]
For example, define the \emph{radial compactification}
\[
  \ol{\R^n} := \Bigl(\R^n \sqcup \bigl( [0,\infty)_\rho\times\Sph^{n-1}\bigr) \Bigr) / \sim
\]
where $r\omega\sim(\rho,\omega)=(r^{-1},\omega)$ in polar coordinates $r,\omega$ on $\R^n$ when $r\neq 0$; then the space $\Vsc(\ol{\R^n})$ is spanned over $\CI(\ol{\R^n})$ (which is the space of classical symbols of order $0$ on $\R^n$) by coordinate vector fields. The corresponding spaces of differential operators are denoted $\Diff_0^m(X)$ and $\Diffsc^m(X)$, and weighted versions of these are defined as in the b-setting. On $M$, one can consider Lie algebras of vector fields featuring different types of degenerations at different hypersurfaces; as a simple example, if $M=[0,1]\times Y$ where $Y$ has no boundary, one can consider the space $\cV_{0,\bop}(M)$ of 0-b-vector fields, which on $[0,1)\times Y$ are 0-vector fields and on $(0,1]\times Y$ b-vector fields.

For each of these algebras, we also consider a semiclassical version. In the b-setting, we thus consider on $[0,1)_h\times X$ the Lie algebra $\cV_{\bop\semi}(X)$ of \emph{semiclassical b-vector fields} consisting of those $V\in\Vb([0,1)_h\times X)$ which are tangent to the level sets of $h$ and vanish at $h=0$; a local frame is $h x\pa_x$, $h\pa_y$, and this is also a frame for the semiclassical b-tangent bundle ${}^{\bop\semi}T X\to [0,1)\times X$. The semiclassical scattering algebra (with tangent bundle denoted $\Tsch X$) and the semiclassical 0-algebra (with tangent bundle ${}^{0\semi}X$) are defined analogously. We commit an abuse of notation here: a section of $\Tsch X$ is really an $h$-dependent vector field on $X$ for $h\in(0,1)$ which depends smoothly on $h$ and degenerates in a particular manner at $h=0$ and at $\pa X$.

\begin{figure}[!ht]
\centering
\includegraphics{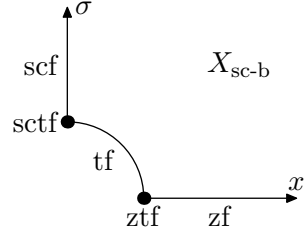}
\caption{The $\scbtop$-single space $X_\scbtop$ of $X$ (with only the $x$-coordinate shown), and its boundary hypersurfaces.}
\label{FigMscbt}
\end{figure}

We furthermore recall from \cite{GuillarmouHassellResI,HintzKdSMS,MelroseSaBarretoLow} the Lie algebra of \emph{scattering-b-tran\-si\-tion vector fields}, denoted $\Vscbt(X)$. To this end, we define the $\scbtop$-single space as the real blow-up
\[
  X_\scbtop := \bigl[ [0,1)_\sigma\times X; \{0\}\times\pa X \bigr].
\]
See Figure~\ref{FigMscbt}. Thus,
\[
  \rho_\scface=\frac{x}{x+\sigma},\qquad
  \rho_\tface=x+\sigma,\qquad
  \rho_\zface=\frac{\sigma}{x+\sigma}
\]
are smooth functions on $X_\scbtop$, and indeed they are defining functions of the scattering face $\scface$ (the lift of $[0,1)\times\pa X$), the transition face $\tface$ (the front face), and the zero face $\zface$ (the lift of $\{0\}\times X$), respectively. See \cite{MelroseDiffOnMwc} for a detailed discussion of blow-ups. Here we recall that $X_\scbtop\setminus\tface=([0,1)\times X)\setminus(\{0\}\times\pa X)$, and $X_\scbtop$ arises by replacing $\{0\}\times\pa X$ with its inward pointing spherical normal bundle $\tface=[0,1]_{\rho_\scface}\times\pa X$; and a collar neighborhood of $\tface$ is $[0,1)_{\rho_\tface}\times[0,1]_{\rho_\scface}\times\pa X$. Moreover, $\tface$ is naturally diffeomorphic to the radially compactified inward pointing normal bundle $\ol{{}^+N}\pa X$. Then $\Vscbt(X)$ consists of all smooth b-vector fields on $X_\scbtop$ which are tangent to the level sets of $\sigma$ and which vanish at $\scface$; they are thus spanned over $\CI(X_\scbtop)$ by the vector fields
\[
  \rho_\scface x\pa_x,\qquad \rho_\scface\pa_y,
\]
which form a frame of $\Tscbt X\to X_\scbtop$. Restriction to $\tface$, resp.\ $\zface$ induces a surjective map $\Vscbt(X)\to\cV_{\scop,\bop}(\tface)$ (with scattering, resp.\ b-behavior at $\sctface=\tface\cap\scface$, resp.\ $\ztface=\tface\cap\zface$), resp.\ $\Vb(\zface)=\Vb(X)$ with kernel $\rho_\tface\Vscbt(X)$, resp.\ $\rho_\zface\Vscbt(X)$, and thus a bundle isomorphism
\[
  \Tscbt_\tface X \cong {}^{\scop,\bop}T\tface,\quad\text{resp.}\quad
  \Tscbt_\zface X \cong \Tb X.
\]
The corresponding spaces of weighted $\scbtop$-operators are denoted
\begin{equation}
\label{EqMscbtDiff}
  \Diffscbt^{m,r,l,b}(X)=\rho_\scface^{-r}\rho_\tface^{-l}\rho_\zface^{-b}\Diffscbt^m(X).
\end{equation}

Finally, we recall from \cite{HintzConicPowers,HintzConicProp} the \emph{semiclassical cone algebra}. We define
\[
  X_\chop := \bigl[ [0,1)_h\times X; \{0\}\times\pa X \bigr].
\]
Its boundary hypersurfaces are the cone face $\cface$ (the lift of $[0,1)\times\pa X$), the transition face $\tface$ (the front face), and the semiclassical face $\sface$ (the lift of $\{0\}\times X$). Defining functions of these boundary hypersurfaces are
\[
  \rho_\cface=\frac{x}{x+h},\qquad
  \rho_\tface=x+h,\qquad
  \rho_\sface=\frac{h}{x+h}\,.
\]
We define $\Vch(X)$ to consist of all smooth b-vector fields on $X_\chop$ which vanish at $\sface$; a local frame of the corresponding tangent bundle $\Tch X\to X_\chop$ is thus
\[
  \rho_\sface x\pa_x,\qquad
  \rho_\sface\pa_y.
\]
Restriction to $\tface$ induces a surjective map $\Vch(X)\to\cV_{\bop,\scop}(\tface)$ with kernel $\rho_\tface\Vch(X)$, and thus a bundle isomorphism $\Tch_\tface X\cong{}^{\bop,\scop}T\tface$.\footnote{Since $\Vch(X)$ will be used for high frequency and $\Vscbt(X)$ for low frequency analysis, the notation $\tface$ will be unambiguous, given the context in which it is used.} Spaces of weighted $\chop$-differential operators are denoted
\begin{equation}
\label{EqMchDiff}
  \Diffch^{m,l,\alpha,b}(X) := \rho_\cface^{-l}\rho_\tface^{-\alpha}\rho_\sface^{-b}\Diffch^m(X).
\end{equation}

We next turn to symbols of operators. The principal symbol of $V\in\cV(X^\circ)$ is the map $T^*X^\circ\ni\zeta\mapsto i \zeta(V)$. For $V\in\Vb(X)$, this is a smooth linear function on $\Tb^*X$; for $V\in\Vch(X)$ it is a smooth linear function on $\Tch^*X$; and so on. Principal symbols of differential operators are then defined by linearity and multiplicativity; if $w\in\CI(X)$ (or $w\in\CI([0,1)\times X)$ in the case of semiclassical algebras, or $w\in\CI(X_\scbtop)$ or $w\in\CI(X_\chop)$ in the case of the $\scbtop$- or $\chop$-algebras) is a product of boundary defining functions so that for the Lie algebra $\cV_\bullet(X)$ under consideration one has $[\cV_\bullet(X),\cV_\bullet(X)]\subset w\cV_\bullet(X)$, then the principal symbol ${}^\bullet\upsigma^m(P)$ of an operator $P\in\Diff_\bullet^m(X)$ (finite linear combinations of up to $m$-fold compositions of elements of $\cV_\bullet(X)$) is well-defined as an element of $(P^m/w P^{m-1})({}^\bullet T^*X)$, where $P^m$ denotes the space of smooth functions on ${}^\bullet T^*X$ which are polynomials of degree $m$ in the fibers. The commutator of $P_j\in\Diff_\bullet^{m_j}(X)$, $j=1,2$, lies in $w\Diff_\bullet^{m_1+m_2-1}(X)$, and has principal symbol given by
\begin{equation}
\label{EqMSymbComm}
  {}^\bullet\upsigma^{m_1+m_2-1}(i[P_1,P_2]) = H_{p_1}p_2,\qquad p_j={}^\bullet\upsigma^{m_j}(P_j),
\end{equation}
where $H_{p_1}$ is the smooth extension from $T^*X^\circ$ of the Hamiltonian vector field of (a representative) $p_1\in\CI(T^*X^\circ)$ to a vector field on ${}^\bullet T^*X$. For the Lie algebras of vector fields on $X$ (i.e.\ without dependence on a semiclassical parameter), note that the phase space ${}^\bullet T^*X$ has a boundary hypersurface ${}^\bullet T^*_{\pa X}X$, and hence one can consider $\bullet$-vector fields on ${}^\bullet T^*X$; one then has $H_p\in\cV_\bullet({}^\bullet T^*X)$ for $p\in\CI({}^\bullet T^*X)$, and indeed the map $p\mapsto H_p$ is an element of $\Diff_\bullet^1({}^\bullet T^*X;\ubar\C,{}^\bullet T({}^\bullet T^*X))$ where $\ubar\C\to{}^\bullet T^*X$ is the trivial bundle.

Two more algebras, which are of central interest in this paper, are recalled in~\S\S\ref{SsMeb}--\ref{SsM3b}. We discuss the corresponding spaces of pseudodifferential operators in~\S\ref{SsMO}. Finally, we define the corresponding scales of $L^2$-based Sobolev spaces in~\S\ref{SsMF}.

\subsection{Edge-b-algebra}
\label{SsMeb}

We work on a manifold with corners $M$, and specifically near the interior of the intersection $Y=H_1\cap H_2$ of two boundary hypersurfaces; we write $x$, resp.\ $z\in\CI(M)$ for a defining function of $H_1$, resp.\ $H_2$ as above. Suppose $H_1$ is the total space of a fibration which in local coordinates $(x,z,y)$ on $M$ is given by $(z,y)\mapsto y$, then we write $\Veb(M)$ for the subspace of $\Vb(M)$ consisting of those $V\in\Vb(M)$ which over $H_1$ are tangent to the fibers. Thus, in a collar neighborhood $[0,1)_x\times[0,1)_z\times Y$ of $Y$, such edge-b-vector fields are smooth linear combinations of the vector fields $x\pa_x$, $z\pa_z$, $x\pa_y$, which form a frame of $\Teb M$.

Restriction to $H_2=z^{-1}(0)$ induces a surjective map $\Veb(M)\to\cV_0(H_2)$ and thus a surjective bundle map $\Teb_{H_2}M\to{}^0 T H_2$; the adjoint is the inclusion map ${}^0 T^*H_2\hra\Teb^*_{H_2}M$. If $P\in\Diffeb^m(M)$, then the pullback of $\sigmaeb^m(P)$ under this inclusion map is ${}^0\upsigma(\wh{N_{H_2}}(P,0))$, as follows from the fact that $\wh{N_{H_2}}(P,0)\in\Diff_0^m(H_2)$ is the restriction of $P$ to $H_2$, i.e.\ the multiplicative extension of the restriction of vector fields. More generally, we define for $\lambda\in\C$ the operator $\wh{N_{H_2}}(P,\lambda)$ as in the b-setting upon fixing a collar neighborhood $[0,1)_z\times H_2$ of $H_2\subset M$. Since $\wh{N_{H_2}}(P,\lambda)-\wh{N_{H_2}}(P,0)\in\Diff_0^{m-1}(M)$, the principal symbol of $\wh{N_{H_2}}(P,\lambda)$ is independent of $\lambda$; for $\lambda\in\R$, this is geometrically due to the fact that the closure of $\lambda\frac{\dd z}{z}+{}^0 T^*H_2$ in the fiber-wise radial compactification $\ol{\Teb^*_{H_2}}M$ intersects fiber infinity $\Seb^*_{H_2}M$ in the $\lambda$-independent submanifold ${}^0 S^*H_2\hra\Seb^*_{H_2}M$. When we allow $\Re\lambda$ to become large, the situation changes; see \cite[\S2.6]{HintzVasyScrieb}. Concretely, fixing $\mu\in\R$, the two families
\begin{equation}
\label{EqMebScl}
  P_{\mu,h} \colon (0,1) \ni h \mapsto \wh{N_{H_2}}(P,-i\mu\pm h^{-1})
\end{equation}
are elements of $\Diff_{0\semi}^{m,0,m}(H_2)=h^{-m}\Diff_{0\semi}^m(H_2)$ (the order `$0$' referring to the weight at the boundary $\pa H_2$), with smooth dependence on $\mu$. This again follows by multiplicativity from the corresponding fact for vector fields: for $P=x D_x$, $z D_z$, $x D_y$, the families~\eqref{EqMebScl} are given by $h^{-1}h x D_x$, $h^{-1}(-i\mu h\pm 1)$, $h^{-1}h x D_y$. The principal symbol of $P_{\mu,\semi}:=(P_{\mu,h})_{h\in(0,1)}$ is the equivalence class of
\begin{equation}
\label{EqMebSclSymb}
  {}^{0\semi}T^*H_2 \ni (h,\varpi) \mapsto \sigmaeb^m(P)\Bigl(\pm h^{-1}\frac{\dd z}{z}+\varpi\Bigr)
\end{equation}
in $(h^{-m}S^m/h^{-(m-1)}S^{m-1})({}^{0\semi}T^*H_2)$, in other words the pullback of $\sigmaeb^m(P)$ under the map
\begin{equation}
\label{EqMebSclMap}
  {}^{0\semi}\iota \colon {}^{0\semi}T^*H_2\ni(h,\varpi)\mapsto\pm h^{-1}\frac{\dd z}{z}+\varpi\in\Teb^*_{H_2}M,\qquad h>0.
\end{equation}

One may think of the principal symbol of $P$ as an object defined at $\Seb^*M\subset\ol{\Teb^*}M$ (indeed it is a section of an appropriate line bundle over $\Seb^*M$), which is the place where edge-b-operators are commutative to leading order; likewise, the principal symbol~\eqref{EqMebSclSymb} of $P_{\mu,\semi}=(P_{\mu,h})_{h\in(0,1)}$ is an object defined at the two boundary hypersurfaces of $\ol{{}^{0\semi}T^*}H_2$ at which semiclassical 0-operators are commutative to leading order, i.e.\ at ${}^{0\semi}S^*H_2=[0,1)_h\times{}^0 S^*H_2$ and $\ol{{}^{0\semi}T^*_{h^{-1}(0)}}H_2$. Underlying the relationship~\eqref{EqMebSclMap} is then the geometric fact that the map~\eqref{EqMebSclMap} can be extended to $\ol{{}^{0\semi}T^*}H_2\to\ol{\Teb^*_{H_2}}M$ and then restricted to these two boundary hypersurfaces to give a pair of smooth maps
\begin{equation}
\label{EqMebSclMaps}
  {}^{0\semi}S^*H_2 \to \Seb^*_{H_2}M,\qquad
  \ol{{}^{0\semi}T^*_{h^{-1}(0)}}H_2 \to \Seb^*_{H_2}M.
\end{equation}
The first map is $h$-independent and for each $h$ given by the inclusion ${}^0 S^*H_2\hra\Seb^*_{H_2}M$; the second map depends on the choice of sign in~\eqref{EqMebSclMap}. As a consequence, if $\Sigma\subset\Seb^*M$ denotes the characteristic set of $P$ (which as a conic subset of $\Teb^*M\setminus o$ we may regard as a subset of $\Seb^*M\cong(\Teb^*M\setminus o)/\R_+$ indeed), then the characteristic set of $P_{\mu,\semi}$ as a semiclassical 0-operator is the union of the preimages of $\Sigma$ under the two maps~\eqref{EqMebSclMaps}.

Moreover, since $N_{H_2}(P)\in\Diffeb^m([0,\infty)_z\times H_2)$ is dilation-invariant in $z$, the Hamiltonian vector field $H_p$ of $p=\sigmaeb^m(P)$ restricted to $\Teb^*_{H_2}M$ conserves the dual momentum to $z$. Therefore, $H_p$ is tangent to the images of the $h$-level sets of ${}^{0\semi}T^*H_2$, and the Hamiltonian vector field of ${}^{0\semi}\upsigma^{m,0,m}(P_h)$ is equal to the pushforward of $H_p$ under~\eqref{EqMebSclMaps}.

A relationship of $P$ and $P_{\mu,\semi}$ via~\eqref{EqMebSclSymb} exists also on the level of subprincipal terms. Concretely, suppose that $P$ is symmetric to leading order, i.e.\ $\Im P:=\frac{P-P^*}{2 i}\in\Diffeb^{m-1}(M)$; here we define the $L^2$-adjoint of $P$ with respect to a density on $M$ of the form $|\frac{\dd z}{z}|\nu$ where $\nu$ is a weighted b-density on $H_2$, i.e.\ $\nu=a(x,y)x^\mu|\frac{\dd x}{x}\dd y|$ for some weight $\mu\in\R$ and a smooth function $a>0$. (Working with such densities ensures that $P^*\in\Diffeb^m(M)$.) Since $\wh{N_{H_2}}(P,\lambda)^*=\wh{N_{H_2}}(P^*,\bar\lambda)$ (the adjoint on the left taken with respect to $\nu$), the semiclassical 0-principal symbol of $\Im\wh{N_{H_2}}(P,\lambda)=\wh{N_{H_2}}(\Im P)(\bar\lambda)+\frac{1}{2 i}(\wh{N_{H_2}}(P,\lambda)-\wh{N_{H_2}}(P,\bar\lambda))$ for $\lambda=-i\mu\pm h^{-1}$ is equal to
\begin{equation}
\label{EqMebSclSubpr}
  {}^{0\semi}\sigma^{m-1,0,m-1}\bigl(\Im\wh{N_{H_2}}(P,-i\mu\pm h^{-1})\bigr) = {}^{0\semi}\iota^*\Bigl(\sigmaeb^{m-1}(\Im P) - (\Im\mu)z^{-1}H_p z\Bigr).
\end{equation}
Here, the second summand arises from the fact that for $P=(z D_z)^k$ and $p=\sigmaeb^k(P)$ one has $\wh{N_{H_2}}(P,\lambda)=\lambda^k=p(\lambda\frac{\dd z}{z})$ and $\Im(\lambda^k)=k(\Im\lambda)\lambda^{k-1}+\cO(|\lambda|^{k-2})$ for bounded $|\Im\lambda|$; note then that $z^{-1}H_p z=\pa_\zeta p=\pa_\zeta(\zeta^k)=k\zeta^{k-1}$ where $\zeta$ is the edge-b-momentum variable dual to $z$. See \cite[Lemma~2.12]{HintzVasyScrieb} for a direct proof of~\eqref{EqMebSclSubpr}.

From~\cite[\S5.1]{HintzVasyScrieb}, we finally recall that the space of \emph{commutator b-vector fields}
\[
  \cV_{[\bop]}(M) = \{ V\in\Vb(M) \colon [V,W]\in\Veb(M)\ \forall\,W\in\Veb(M) \}
\]
spans $\Vb(M)$ over $\CI(M)$. Given a vector bundle $E\to M$, the space $\Diff_{[\bop]}^1(M;E)$ of \emph{commutator b-operators} consists of operators with scalar principal part given by a commutator b-vector field. For $A\in\Diffeb^m(M;E)$ and $X\in\Diff_{[\bop]}^1(M;E)$, one then has $[X,A]\in\Diffeb^m(M;E)$. The analogous statement remains true also for pseudodifferential $A$, and when the coefficients of $A$ are merely conormal at $\pa M$; see \cite[Lemma~5.6]{HintzVasyScrieb}.

\subsection{3b-algebra}
\label{SsM3b}

Following \cite{Hintz3b}, we consider a smooth manifold $M_0$ with boundary, a point $\fp\in\pa M_0$, and the blow-up $M:=[M_0;\{\fp\}]$. We write $\rho_\cD$, resp.\ $\rho_\cT\in\CI(M)$ for a defining function of $\cD$ (the lift of $\pa M_0$) and the front face $\cT$. We write $\rho_0\in\CI(M_0)$ for a boundary defining function on $M_0$. Denoting the blow-down map by $\upbeta\colon M\to M_0$, the function $\upbeta^*\rho_0$ is a total boundary defining function on $M$; for brevity, we denote it by $\rho_0$ as well. The space $\Vtb(M)$ of 3b-vector fields is then the span over $\CI(M)$ of $\rho_\cD^{-1}\upbeta^*\Vsc(M_0)$; equivalently, a vector field $V\in\Vb(M)$ lies in $\Vtb(M)$ if and only if $V\rho_0\in\rho_0\rho_\cT\CI(M)$. In local coordinates $T\geq 0$, $X\in\R^n$ near $\fp$, with $\fp=(0,0)$, and letting $t=T^{-1}$, $x=\frac{X}{T}$, the space $\Vtb(M)$ is spanned in $r=|x|\lesssim 1$ by $\pa_t$, $\pa_x$ and in $r\gtrsim 1$ by
\begin{equation}
\label{EqM3btr}
  r\pa_t,\quad r\pa_r,\quad \pa_\omega,
\end{equation}
where $\pa_\omega$ is schematic notation for spherical derivatives in $x$ such as $x^i\pa_{x^j}-x^j\pa_{x^i}$. The tangent bundle corresponding to $\Vtb(M)$ is denoted $\Ttb M$, and we write $\Difftb^m(M)$ and $\Difftb^{m,(\alpha_\cD,\alpha_\cT)}(M)=\rho_\cD^{-\alpha_\cD}\rho_\cT^{-\alpha_\cT}\Difftb^m(M)$ for the corresponding spaces of (weighted) 3b-differential operators.

Since $\Vtb(M)\subset\Vb(M)$, any $P\in\Difftb^m(M)$ has a normal operator at $\cD=[\pa M_0;\{\fp\}]$; by contrast to the b-setting however, one regards it dilation-invariant operator on the resolution
\[
  {}^+N_\tbop\cD := \bigl[ [0,\infty)_{\rho_0}\times\pa M_0; [0,\infty)_{\rho_0}\times\{\fp\} \bigr] = [0,\infty)_{\rho_0}\times\cD
\]
of the lifted inward pointing normal bundle of $\pa M_0$ (rather than an operator on the inward pointing normal bundle of $\cD\subset M$); then $N_\cD(P)$ is a dilation-invariant edge-b-operator on ${}^+N_\tbop\cD$ with respect to the boundary fibration $[0,\infty)\times\pa\cD\to[0,\infty)$. Concretely, working in coordinates $t,r,\omega$ as before, and writing $T=\frac{1}{t}$, $R=\frac{r}{t}$, then given a 3b-differential operator
\begin{align}
\label{EqM3bOp}
  P &= \sum_{j+k+|\alpha|\leq m} a_{j k\alpha}(r D_t)^j (r D_r)^k D_\omega^\alpha \\
    &= \sum_{j+k+|\alpha|\leq m} a_{j k\alpha}(-R T D_T-R^2 D_R)^j(R D_R)^k D_\omega^\alpha \nonumber
\end{align}
with $a_{j k\alpha}\in\CI(M)$, the operator $N_\cD(P)$ is given by freezing the coefficients $a_{j k\alpha}$ at $T=0$; thus $a_{j k\alpha}|_\cD=:b_{j k\alpha}\in\CI(\cD)$, in local coordinates, is a smooth function of $R\geq 0$ and $\omega\in\Sph^{n-1}$. (We also note that 3b-vector fields are locally spanned by $R T\pa_T$, $R\pa_R$, $\pa_\omega$ in the $T,R,\omega$ coordinates.) The Mellin-transformed normal operator family is defined with respect to $\rho_0$, so
\begin{equation}
\label{EqM3bNDMT}
  \wh{N_\cD}(P,\lambda) = \sum_{j+k+|\alpha|\leq m} b_{j k\alpha}(R,\omega) (-R\lambda-R^2 D_R)^j (R D_R)^k D_\omega^\alpha.
\end{equation}
Its b-normal operator at $\pa\cD=R^{-1}(0)$ is independent of $\lambda$. For $\mu\in\R$, the high energy family
\begin{equation}
\label{EqM3bScl}
  P_{\mu,\semi} \colon (0,1)\ni h \mapsto \wh{N_\cD}(P,-i\mu\pm h^{-1})
\end{equation}
defines an element of $\Diffch^{m,0,0,m}(\cD)$; for $m=1$ this follows from the fact that for
\[
  P=-R T D_T-R^2 D_R,\qquad R D_R,\qquad D_\omega,
\]
the family~\eqref{EqM3bScl} is given by
\begin{align*}
  &i R\mu\mp R h^{-1}-R^2 D_R=\frac{h+R}{h}\Bigl(\mp\frac{R}{R+h}+i\frac{h R}{h+R}\mu-R\frac{h}{h+R}R D_R\Bigr), \\
  &\qquad \frac{h+R}{h}\frac{h}{h+R}R D_R,\qquad
  \frac{h+R}{h}\frac{h}{h+R}D_\omega,
\end{align*}
with $\Vch(\cD)$ spanned over $\CI(\cD_\chop)$ by $\frac{h}{h+R}R D_R$ and $\frac{h}{h+R}D_\omega$. Letting $\tilde R=\frac{R}{h}$, we write
\begin{equation}
\label{EqM3bNDtf}
  N_{\cD,\tface}^\pm(P)=\sum_{j+k+|\alpha|\leq m} b_{j k\alpha}(0,\omega)(\mp\tilde R)^j(\tilde R D_{\tilde R})^k D_\omega^\alpha \in \Diff_{\bop,\scop}^{m,0,m}(\tface)
\end{equation}
for the normal operators of~\eqref{EqM3bScl} at $[0,\infty]_{\tilde R}\times\Sph^{n-1}=\tface\subset\cD_\chop$.

\begin{figure}[!ht]
\centering
\includegraphics{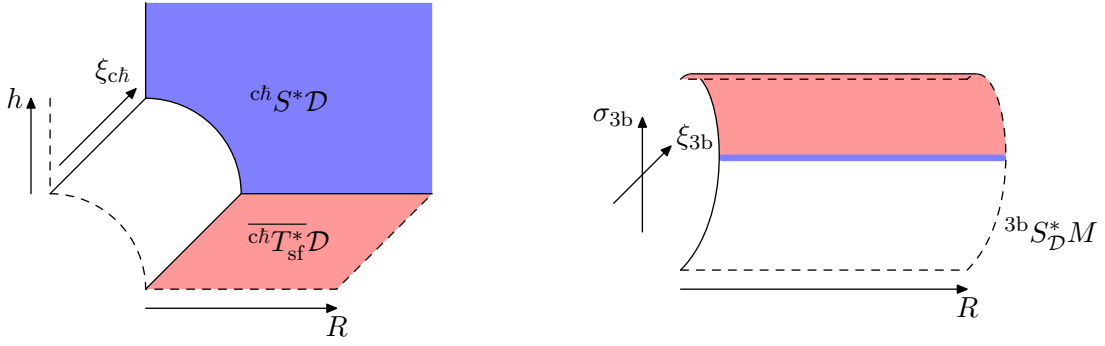}
\caption{Illustration of the maps~\eqref{EqM3bSclMaps} (for the `$+$' sign in~\eqref{EqM3bSclMap}). \textit{On the left:} the semiclassical cone phase space, with $\xi_\chop=\frac{h}{h+R}R\pa_R(\cdot)$. \textit{On the right:} the 3b-phase space, with $\sigma_\tbop=R T\pa_T(\cdot)$ and $\xi_\tbop=R\pa_R(\cdot)$. The domain and codomain of the first, resp.\ second map in~\eqref{EqM3bSclMaps} is drawn in blue, resp.\ red. We only show the region where $\xi_\chop,\xi_\tbop>0$.}
\label{FigM3bSclMaps}
\end{figure}

We recall the relationships between principal symbols from \cite[Proposition~3.22, Lemma~2.31]{Hintz3b}. Namely, $\sigmab^m(\wh{N_\cD}(P,\lambda))$ is the pullback of $\sigmatb^m(P)$ along the inclusion $\Tb^*\cD\hra\Ttb^*_\cD M$ which is the dual of the map of tangent bundles induced by the restriction $\Vtb(M)\to\Vb(\cD)$ of vector fields; and $\sigmach^{m,0,0,m}(P_{\mu,\semi})$ is the pullback along
\begin{equation}
\label{EqM3bSclMap}
  {}^\chop\iota \colon \Tch^*\cD \in (h,\varpi) \mapsto \pm h^{-1}\frac{\dd\rho_0}{\rho_0} + \varpi \in \Ttb^*_\cD M,\qquad h>0.
\end{equation}
(Note here that in local coordinates, $\frac{\dd\rho_0}{\rho_0}=R\frac{\dd T}{R T}$ is a 3b-covector vanishing simply at $R=0$.) In analogy with~\eqref{EqMebSclMaps}, we record that the smooth extension $\ol{\Tch^*}\cD\to\ol{\Ttb^*_\cD}M$ of the map~\eqref{EqM3bSclMap} restricts to a pair of smooth maps
\begin{equation}
\label{EqM3bSclMaps}
  \Sch^*\cD \to \Stb^*_\cD M,\qquad
  \ol{\Tch^*_\sface}\cD \to \Stb^*_\cD M,
\end{equation}
each of which depends on the choice of sign in~\eqref{EqM3bSclMap}. This is not explicitly stated in the reference; but the smoothness of the maps~\eqref{EqM3bSclMaps} follows from a direct calculation in local coordinates. For example, near $\sface\subset\cD_\chop$, we may work with the coordinates $\tilde h=\frac{h}{R}$, $R$, $\omega$ in the base and $\xi_\chop,\eta_\chop$ in the fibers of $\Tch^*\cD$, where we write $\chop$-covectors as $\xi_\chop\frac{R}{h}\frac{\dd R}{R}+\eta_\chop\frac{R}{h}\dd\omega$; then ${}^\chop\iota$ maps $(\tilde h,R,\omega,\xi_\chop,\eta_\chop)\mapsto\tilde h^{-1}(\pm\frac{\dd T}{R T}+\xi_\chop\frac{\dd R}{R}+\eta_\chop\,\dd\omega)$, a 3b-covector at the point $(0,R,\omega)\in\cD$ (the first coordinate being $T=0$). Thus, $\ol{\Tch^*_\sface}\cD$ (where $\tilde h=0$) gets mapped to the set of endpoints at fiber infinity of $\ol{\Ttb^*_\cD}M$ of the rays $\R_+(\pm\frac{\dd T}{R T}+\xi_\chop\frac{\dd R}{R}+\eta_\chop\,\dd\omega)$. On the other hand, the endpoint at fiber infinity of $\ol{\Tch^*}\cD$ of the ray $\R_+(\xi_\chop\frac{\dd R}{R}+\eta_\chop\,\dd\omega)$ with $(\xi_\chop,\eta_\chop)\neq(0,0)$ gets mapped to the endpoint of $\R_+(\xi_\chop\frac{\dd R}{R}+\eta_\chop\,\dd\omega)$ (independently of the value of $\tilde h$). See Figure~\ref{FigM3bSclMaps}.

Turning to subprincipal symbols, suppose that $P\in\Difftb^m(M)$ satisfied $\Im P=\frac{P-P^*}{2 i}\in\Difftb^{m-1}(M)$, where the adjoint is defined with respect to the lift to $M$ of a density $\frac{\dd\rho_0}{\rho_0}\nu$ on $M_0$, where $\nu$ is a positive density on $\pa M_0$. Then
\begin{equation}
\label{EqM3bSubpr}
  \sigmach^{m-1,0,0,m-1}\bigl(\Im\wh{N_\cD}(P,-i\mu\pm h^{-1})\bigr) = {}^\chop\iota^*\Bigl(\sigmatb^{m-1}(\Im P) - (\Im\mu)\rho_0^{-1}H_p\rho_0\Bigr),
\end{equation}
where the adjoint of $\wh{N_\cD}(P,-)$ is defined with respect to the lift of the density $\nu$ to $\cD$. This relationship follows from the same arguments as~\eqref{EqMebSclSubpr}, and can also be checked explicitly using~\eqref{EqM3bNDMT}.

Following~\cite[\S3.2]{Hintz3b}, we now turn to the second, time translation invariant normal operator of $P\in\Difftb^m(M)$. Writing $P$ in $r>1$ as in~\eqref{EqM3bOp}, we freeze the coefficients $a_{j k\alpha}$ at $\cT$; thus, with $c_{j k\alpha}(r,\omega)$ denoting the restriction $a_{j k\alpha}|_\cT$, we set
\[
  N_\cT(P) = \sum_{j+k+|\alpha|\leq m} c_{j k\alpha}(r,\omega) (r D_t)^j (r D_r)^k D_\omega^\alpha.
\]
(For $r\lesssim 1$ and in terms of $x=r\omega$, the operator $N_\cT(P)$ can be written in the form $\sum_{j+|\alpha|\leq m} d_{j\alpha}(x) D_t^j D_x^\alpha$.) The spectral family is
\[
  \wh{N_\cT}(P,\sigma) = \sum_{j+k+|\alpha|\leq m} c_{j k\alpha} (-r\sigma)^j (r D_r)^k D_\omega^\alpha.
\]
This is formally the conjugation of $N_\cT(P)$ by the Fourier transform $\cF_{t\to\sigma}$ in $t$, using the convention $\cF_{t\to\sigma}u=\int e^{i\sigma t}u(t)\,\dd t$. Thus, $\wh{N_\cT}(P,0)\in\Diffb^m(\cT)$ and $\wh{N_\cT}(P,\sigma)\in\Diffsc^{m,m}(\cT)$ for $\sigma\neq 0$. Passing to $\rho_\cD=r^{-1}$, which is a local defining function of the boundary $\pa\cT=\cT\cap\cD$ of $\cT$, one finds that
\begin{equation}
\label{EqM3bscbt}
  \pm[0,1)\ni\sigma\mapsto\wh{N_\cT}(P,\sigma) = \sum_{j+k+|\alpha|\leq m} c_{j k\alpha} (-\sigma\rho_\cD^{-1})^j(-\rho_\cD D_{\rho_\cD})^k D_\omega^\alpha
\end{equation}
is an element of $\Diffscbt^{m,m,0,0}(\cT)$. (For $m=1$ and setting $\rho_\scface=\frac{\rho_\cD}{|\sigma|+\rho_\cD}$, this follows from $\sigma\rho_\cD^{-1}=\rho_\scface^{-1}\frac{\sigma}{|\sigma|+\rho_\cD}\in\Diffscbt^{0,1,0,-1}(\cT)\subset\Diffscbt^{1,1,0,0}(\cT)$ and $\rho_\cD D_{\rho_\cD}=\rho_\sface^{-1}\rho_\sface\rho_\cD D_{\rho_\cD}$, $D_\omega=\rho_\sface^{-1}\rho_\sface D_\omega\in\Diffscbt^{1,1,0,0}(\cT)$.) In terms of $\hat\rho_\cD=\frac{\rho_\cD}{|\sigma|}=\pm\frac{\rho_\cD}{\sigma}$, the transition face model operator of~\eqref{EqM3bscbt} is
\begin{equation}
\label{EqMtbTtf}
  N_{\cT,\tface}^\pm(P) = \sum_{j+k+|\alpha|\leq m} b_{j k\alpha}(0,\omega) (\mp\hat\rho_\cD^{-1})^j (-\hat\rho_\cD D_{\hat\rho_\cD})^k D_\omega^\alpha \in \Diff_{\scop,\bop}^{m,m,0}(\tface)
\end{equation}
where $b_{j k\alpha}$ is the same function as appearing already in~\eqref{EqM3bNDtf}; note indeed that $b_{j k\alpha}$ is the restriction of $a_{j k\alpha}$ to $\cT\cap\cD$ and thus of $c_{j k\alpha}$ to $\pa\cT$. Pullback under the diffeomorphism $(\hat\rho_\cD,\omega)\mapsto(\tilde R,\omega)=(\hat\rho_\cD^{-1},\omega)$ between the $\scbtop$- and $\chop$-transition faces identifies $N_{\cT,\tface}^\pm(P)$ and $N_{\cD,\tface}^\pm(P)$; see also~\cite[Proposition~3.28]{Hintz3b}. In the high energy regime on the other hand, we recall that the operator families
\begin{equation}
\label{EqM3bsch}
  (0,1)\ni h\mapsto \wh{N_\cT}(P,\pm h^{-1})
\end{equation}
define elements of $\Diffsch^{m,m,m}(\cT)=\rho_\cD^{-m}h^{-m}\Diffsch^m(\cT)$.

\begin{figure}[!ht]
\centering
\includegraphics{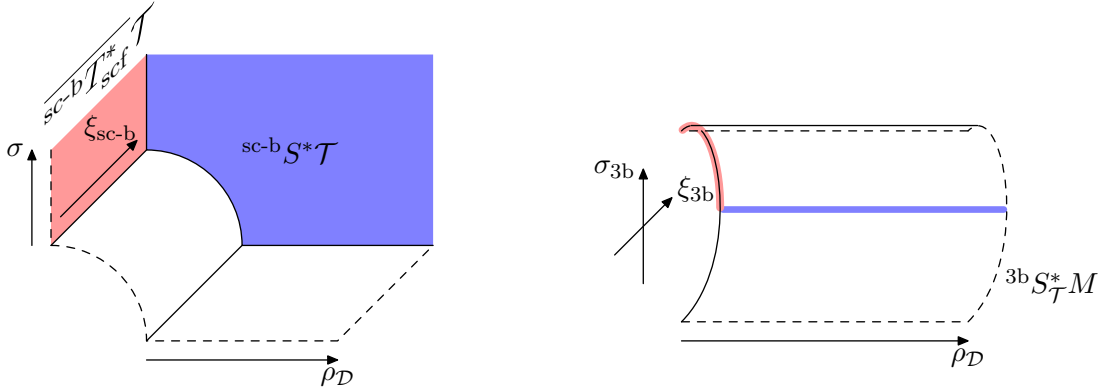}
\caption{Illustration of the maps~\eqref{Eq3bLoMaps} corresponding to positive frequencies. \textit{On the left:} the scattering-b-transition phase space. \textit{On the right:} the 3b-phase space, with $\sigma_\tbop=-r\pa_t(\cdot)$ and $\xi_\tbop=r\pa_r(\cdot)$. The domain and codomain of the first, resp.\ second map in~\eqref{Eq3bLoMaps} is drawn in blue, resp.\ red.}
\label{FigM3bLoMaps}
\end{figure}

The various phase spaces are related as follows: we have an inclusion $\Tb^*\cT\hra\Ttb^*_\cT M$ which is dual to the map induced by the surjective restriction map $\Vtb(M)\to\Vb(\cT)$ of spaces of vector fields; the b-principal symbol of $\wh{N_\cT}(P,0)$ is the pullback of $p=\sigmatb^m(P)$ under this map. For $\sigma\neq 0$, the scattering principal symbol of $\wh{N_\cT}(P,\sigma)$ is the pullback of $p$ under the map
\[
  T^*\cT^\circ \ni \varpi \mapsto -\sigma\,\dd t + \varpi \in \Ttb^*_\cT M.
\]
For finite $\sigma\neq 0$, this map extends to $\Tsc^*\cT$ only if one compactifies the fibers of the codomain; we may then also compactify the domain and obtain the map
\begin{equation}
\label{EqM3biotasigma}
  \iota(\sigma,-) \colon \ol{\Tsc^*}\cT \ni \varpi \mapsto -\sigma\,\dd t+\varpi \in \ol{\Ttb^*_\cT}M.
\end{equation}
Writing scattering covectors as $\varpi=\xi_\scop\,\dd r+\eta_\scop\,r^{-1}\dd\omega$ and using $\rho_\cD=r^{-1}$, this map takes the form
\[
  \iota(\sigma,\xi_\scop\,\dd r+\eta_\scop\,r^{-1}\dd\omega) = \rho_\cD^{-1}\Bigl(-\sigma\frac{\dd t}{r}+\xi_\scop\frac{\dd r}{r}+\eta_\scop\,\dd\omega\Bigr);
\]
the smoothness of~\eqref{EqM3biotasigma} follows easily from this. Moreover, regarding
\[
  {}^\scbtop\iota:=(\iota(\sigma,-))_{\sigma\in\pm[0,1)}
\]
as a single map, the map ${}^\scbtop\iota\colon\ol{\Tscbt^*}\cT\to\ol{\Ttb^*_\cT}M$ is smooth, and it restricts to smooth inclusion maps
\begin{equation}
\label{Eq3bLoMaps}
  \Sscbt^*\cT \to \Stb^*_\cT M,\qquad
  \ol{\Tscbt^*_\scface}\cT \to \Stb^*_{\pa\cT}M.
\end{equation}
For instance, near $\scface\subset\cT_\scbtop$, with coordinates $\sigma$, $\hat\rho_\cD=\frac{\rho_\cD}{|\sigma|}=\pm\frac{\rho_\cD}{\sigma}$, $\omega$, and writing $\scbtop$-covectors as $\xi_\scbtop\frac{\dd\rho_\cD}{\hat\rho_\cD\rho_\cD}+\eta_\scbtop\frac{\dd\omega}{\hat\rho_\cD}$, the map ${}^\scbtop\iota$ takes the form $(\sigma,\hat\rho_\cD,\omega,\xi_\scbtop,\eta_\scbtop)\mapsto\hat\rho_\cD^{-1}(\mp\frac{\dd t}{r}-\xi_\scbtop\frac{\dd r}{r}+\eta_\scbtop\,\dd\omega)$, the output being a 3b-covector over the point on $\cT$ with coordinates $\rho_\cD=|\sigma|\hat\rho_\cD$ and $\omega$; this easily gives~\eqref{Eq3bLoMaps} in the region $\hat\rho_\cD\lesssim 1$. See Figure~\ref{FigM3bLoMaps} (and also \cite[Figure~3.3]{Hintz3b}).

At high positive or negative frequencies, the map ${}^\schop\iota:=(\iota(\pm h^{-1},-))_{h\in(0,1)}$ extends to a smooth map $\ol{\Tsch^*}\cT\to\ol{\Ttb^*_\cT}M$ which then restricts to smooth inclusion maps
\begin{equation}
\label{Eq3bHiMaps}
  {}^\schop S^*\cT \to \Stb^*_\cT M,\qquad
  \ol{\Tsch^*_{[0,1)\times\pa\cT}}\cT \to \Stb^*_{\pa\cT}M,\qquad
  \ol{\Tsch^*_{h^{-1}(0)}}\cT \to \Stb^*_\cT M.
\end{equation}
The principal symbols of~\eqref{EqM3bscbt} and \eqref{EqM3bsch} are given via pullback of $\sigmatb^m(P)$ along the maps ${}^\scbtop\iota$ and ${}^\schop\iota$, respectively; the corresponding Hamiltonian vector fields are related in the same manner. We leave the statements of relationships of subprincipal symbols to the interested reader.

Finally, for later use, we note that if $T\in\CI(M_0)$ denotes a boundary defining function (and hence its lift to $M$ is a particular total boundary defining function), then
\begin{equation}
\label{Eq3bConjT}
  P\in\Difftb^m(M) \implies T^\alpha[P,T^{-\alpha}] \in \rho_\cT\Difftb^{m-1}(M).
\end{equation}
Indeed, this follows from the case $P\in\Vtb(M)$ and the observation that, in the coordinates $t,r$ used in~\eqref{EqM3btr}, one has $t^{-\alpha}[r\pa_t,t^\alpha]=\alpha\frac{r}{t}\in\rho_\cT\CI(M)$, whereas the other vector fields in~\eqref{EqM3btr} commute with $t^\alpha$.

\subsubsection{Time dilations}
\label{SssM3bDil}

In~\cite{Hintz3b}, the b-regularity of solutions of 3b-equations was proved in the fully elliptic setting via a parametrix construction. In the wave equation context of the present paper, we will instead proceed via testing with suitable vector fields.

\begin{definition}[Time dilation vector field]
\label{DefM3bb}
  We call $V\in\Vb(M)$ a \emph{time dilation vector field} if and only if $V$ is a b-normal vector field at $\cT$, i.e.\ $V-\rho_\cT\pa_{\rho_\cT}\in\rho_\cT\Vb(M)$. If $E\to M$ is a vector bundle, then $X\in\Diffb^1(M;E)$ is a \emph{time dilation operator} if the principal part of $X$ is scalar and a time dilation vector field, and if moreover $X\colon\CI(M;E)\to\rho_\cT\CI(M;E)$.
\end{definition}

Here, we recall the general fact that if $X$ is a manifold with boundary and $[0,1)_x\times\pa X$ is a collar neighborhood of $\pa X$, then $x\pa_x$ is independent, modulo $x\Vb(X)$, of the choice of the collar neighborhood. Note that $-t\pa_t$ is a time dilation vector field; hence the terminology. Moreover, any element of $\Vb(M)$ is the sum of a smooth multiple of $-t\pa_t$ and a 3b-vector field, cf.\ \eqref{EqM3btr}. Lastly, given a time dilation vector field $V$, we can construct a corresponding time dilation operator $X$ locally in a local trivialization of $E$ as a diagonal matrix with diagonal entries $V$, and globally via a partition of unity.

\begin{lemma}[Commutators I: vector fields]
\label{LemmaM3bbComm}
  Let $V$ be a time dilation vector field, and let $W\in\Vtb(M)$. Then there exist $w^\flat\in\CI(M)$ and $W^\sharp\in\Vtb(M)$ so that
  \begin{equation}
  \label{EqM3bbComm}
    [V,W] = \rho_\cT w^\flat V + \rho_\cT W^\sharp.
  \end{equation}
  In particular, $[V,W]\in\Vtb(M)$.
\end{lemma}
\begin{proof}
  We verify this in the coordinates $t,r,\omega$ used in~\eqref{EqM3btr}, in which we can write
  \[
    W = a\Bigl(\frac{r}{t},\frac{1}{r},\omega\Bigr) r\pa_t + b(\frac{r}{t},\frac{1}{r},\omega\Bigr) r\pa_r + \sum_{j=1}^{n-1} c^j(\frac{r}{t},\frac{1}{r},\omega\Bigr)\pa_{\omega^j},
  \]
  where $a,b,c^j\in\CI([0,\infty)\times[0,\infty)\times\R^{n-1})$ are smooth functions of $\rho_\cT=\frac{r}{t}$, $\rho_\cD=\frac{1}{r}$, $\omega$. We have
  \[
    V=-t\pa_t+\tilde V,\qquad
    \tilde V=\frac{r}{t}(\tilde a t\pa_t+\tilde b r\pa_r+\sum\tilde c^j\pa_{\omega^j}),
  \]
  with $\tilde a,\tilde b,\tilde c^j$ smooth. We compute
  \begin{align*}
    [-t\pa_t,W] &= a r\pa_t - a'r\pa_t-b'r\pa_r-\sum_{j=1}^{n-1} (c^j)'\pa_{\omega^j},
  \end{align*}
  where the terms involving $a':=-\rho_\cT\pa_{\rho_\cT}a\in\rho_\cT\CI(M)$ and $b':=-\rho_\cT\pa_{\rho_\cT}b$, $(c^j)':=-\rho_\cT\pa_{\rho_\cT}c^j\in\rho_\cT\CI(M)$ can be absorbed into $W^\sharp$ in~\eqref{EqM3bbComm}, while $a r\pa_t=\rho_\cT a t\pa_t$ can be absorbed into the first. On the other hand, since $r\pa_r$, $\pa_{\omega^j}$ are 3b-vector fields, we have $[\frac{r}{t}\tilde b r\pa_r,W]$, $[\frac{r}{t}\tilde c^j\pa_{\omega^j},W]\in\rho_\cT\Vtb(M)$; moreover, the first term of $[\frac{r}{t}\tilde a t\pa_t,W]=\rho_\cT\tilde a[t\pa_t,W]-W(\frac{r}{t}\tilde a)t\pa_t$ is of the form~\eqref{EqM3bbComm} by what we have already shown, and so is the second term since $W(\rho_\cT\tilde a)\in\rho_\cT\CI(M)$.
\end{proof}

\begin{rmk}[Equivalent characterizations]
\label{RmkM3bbEquiv}
  Let $V_\cT$ be a time dilation vector field. The space $\cV_{[\bop]}(M)$ of $V\in\Vb(M)$ so that $[V,\Vtb(M)]\subset\Vtb(M)$ contains $\Vtb(M)$; since every $V\in\Vb(M)$ can be written as $V=a V_\cT+V'$ with $a\in\CI(M)$ and $V'\in\Vtb(M)$, we infer that $[a V_\cT,W]=-(W a)V_\cT+a[V_\cT,W]\in\Vtb(M)$ for all $W\in\Vtb(M)$, so necessarily $W a\in\rho_\cT\CI(M)$, since by Lemma~\ref{LemmaM3bbComm} the second summand lies in $\Vtb(M)$ already. Therefore, $a|_\cT$ is constant. This shows that $\Vtb(M)$ and $V_\cT$ generate the space of all `commutator b-vector fields' $\cV_{[\bop]}(M)$ over $\CI_\cT(M)=\{a\in\CI(M)\colon a|_\cT\ \text{is constant}\}$. The space $\cV_{[\bop]}(M)$ can equivalently be characterized as the $\CI_\cT(M)$-span of all pullbacks $\upbeta^*V_0$ where $V_0\in\Vb(M_0)$ is tangent to $\fp$.
\end{rmk}

\begin{lemma}[Commutators II: differential operators]
\label{LemmaM3bbComm2}
  Let $X\in\Diffb^1(M;E)$ be a time dilation operator, and let $A\in\Difftb^m(M;E)$. Then there exist $A^\flat\in\Difftb^{m-1}(M;E)$ and $A^\sharp\in\Difftb^m(M;E)$ so that
  \begin{equation}
  \label{EqM3bbComm2}
    [X,A] = \rho_\cT A^\flat X + \rho_\cT A^\sharp.
  \end{equation}
  In particular, $[X,A]\in\Difftb^m(M;E)$. For $A\in\cA^{(\alpha_\cD,\alpha_\cT)}\Difftb^m(M;E)$, we have $[X,A]\in\cA^{(\alpha_\cD,\alpha_\cT)}\Difftb^m(M;E)$.
\end{lemma}
\begin{proof}
  For $m=0$, so $A\in\CI(M;\End(E))$, the conclusion (with $A^\flat=0$ and $A^\sharp=\rho_\cT^{-1}X(A)$) follows from Definition~\ref{DefM3bb}. For $m=1$, the claim follows from Lemma~\ref{LemmaM3bbComm}. For $m\geq 2$, we write $A=\sum_i Y_i A_i$ with $Y_i\in\Difftb^1(M;E)$, $A_i\in\Difftb^{m-1}(M;E)$; then for suitable $A_i^\flat\in\Difftb^{m-2}$, $A_i^\sharp\in\Difftb^{m-1}$ and $Y_i^\flat\in\CI$, $Y_i^\sharp\in\Difftb^1$, we compute
  \[
    [X,Y_i A_i] = [X,Y_i]A_i + Y_i[X,A_i] = \rho_\cT Y_i^\flat X A_i + \rho_\cT Y_i^\sharp A_i + Y_i \rho_\cT A_i^\flat X + Y_i\rho_\cT A_i^\sharp.
  \]
  In the first term, we write $X A_i=A_i X+[X,A_i]$ and use the inductive hypothesis on $[X,A_i]$. The second term is already of the desired form. In third and fourth terms, we write $Y_i\rho_\cT=\rho_\cT Y_i+[Y_i,\rho_\cT]$ with $[Y_i,\rho_\cT]\in\rho_\cT\CI$. This establishes the decomposition~\eqref{EqM3bbComm2}.

  For the case of conormal coefficients, it suffices to consider an operator $A=f B$ with $f\in\cA^{(\alpha_\cD,\alpha_\cT)}(M)$, $B\in\Difftb^m(M;E)$. Then
  \begin{align*}
    [X,A]=f[X,B]+[X,f]B &\in f\Difftb^m(M;E) + (X f)\Difftb^m(M;E) \\
      &\subset \cA^{(\alpha_\cD,\alpha_\cT)}\Difftb^m(M;E),
  \end{align*}
  as desired.
\end{proof}

\subsubsection{Mixed \texorpdfstring{$\tbop;\bop$}{3b;b}-operators}
\label{SssM3bb}

For the purpose of tracking b-regularity relative to 3b-reg\-u\-lar\-i\-ty, we now introduce a mixed algebra in the spirit of similar earlier constructions e.g.\ in \cite{VasyPropagationCorners,VasyWaveOnAdS,MelroseVasyWunschEdge,MelroseVasyWunschDiffraction,HintzVasyScrieb}.

\begin{definition}[Mixed $\tbop;\bop$-differential operators]
\label{DefM3bbOp}
  Let $m,k\in\N_0$. Then the space $\Diffb^k\Difftb^m(M)$ consists of all operators which are locally finite sums $\sum_i Q_i P_i$ where $Q_i\in\Diffb^k(M)$ and $P_i\in\Difftb^m(M)$.
\end{definition}

With purely notational changes, one can also consider spaces of weighted operators, $\Diffb^k\Difftb^{m,(\alpha_\cD,\alpha_\cT)}(M)=\Diffb^{k,(\alpha_\cD,\alpha_\cT)}\Difftb^m(M)$, and also versions with conormal coefficients, $\cA^{(\alpha_\cD,\alpha_\cT)}\Diffb^k\Difftb^m(M)$, as well as operators acting between sections of bundles. We do not spell out these purely notational generalizations.

\begin{lemma}[Commuting the two factors]
\label{LemmaM3bbF}
  Let $Q\in\Diffb^k(M)$ and $P\in\Difftb^m(M)$. Then there exist finitely many operators $Q_j^\flat,Q_j^\sharp\in\Diffb^k(M)$ and $P_j^\flat,P_j^\sharp\in\Difftb^m(M)$ so that
  \begin{equation}
  \label{EqM3bbF}
    Q P = \sum_j P_j^\flat Q_j^\flat,\qquad
    P Q = \sum_j Q_j^\sharp P_j^\sharp.
  \end{equation}
\end{lemma}
\begin{proof}
  The case $k=0$ is trivial. For $k=1$, we note that $\Diffb^1(M)$ is spanned over $\CI(M)$ by $\Difftb^1(M)$ and any time dilation vector field $V\in\Vb(M)$; it thus suffices to consider the case $Q=a V$ where $a\in\CI(M)$. In this case, we have $Q P=P Q-[P,Q]$ where $[P,Q]=[P,a]V+a[P,V]$, with the first term in $\Difftb^{m-1}\Diffb^1$ and the second in $\Difftb^m$ by Lemma~\ref{LemmaM3bbComm2}. Similarly, $P Q=Q P+[P,Q]$ where we now further write $[P,Q]=V[P,a]+[[P,a],V]+a[P,V]$ and apply Lemma~\ref{LemmaM3bbComm2} to the second term on the right.

  If~\eqref{EqM3bbF} is established for b-differential orders $\leq k-1$, then it follows for orders $\leq k$ by writing $Q\in\Diffb^k(M)$ as a finite sum of terms $Q_1 Q_2$ with $Q_1\in\Diffb^1$ and $Q_2\in\Diffb^{k-1}$. Then $Q_2 P=\sum P_{2 j}^\flat Q_{2 j}^\flat$ with $Q_{2 j}^\flat\in\Diffb^{k-1}$ and $P_{2 j}^\flat\in\Difftb^m$ by the inductive hypothesis. Using the case $k=1$, we can further write $Q_1 P_{2 j}^\flat=\sum P_{2 j k}^\flat Q_{2 j k}^\flat$ with $Q_{2 j k}^\flat\in\Diffb^1$ and $P_{2 j k}^\flat\in\Difftb^m$. This proves the first part of~\eqref{EqM3bbF}. The proof of the second part is similar.
\end{proof}

\subsection{Pseudodifferential operators}
\label{SsMO}

We now turn to a brief description of the main features of the classes of pseudodifferential operators corresponding to the Lie algebras of vector fields recalled above. We give a unified description in the same manner as around~\eqref{EqMSymbComm}; we write $X$ for the manifold which, depending on the algebra under consideration, is a manifold with boundary or with corners. Thus, we consider a Lie algebra $\cV_\bullet(X)$ of vector fields on the appropriate single space $X_\bullet$. This means that $X_\bullet=X$ for $\bullet=\bop,0,\scop,(\ebop),\tbop$, further $X_\bullet=[0,1)\times X$ when $\bullet=\bop\semi,0\semi$, or $\schop$, and finally $X_\bullet=X_\scbtop$ or $X_\bullet=X_\chop$ for $\bullet=\scbtop$ or $\chop$. We denote by $w\in\CI(X_\bullet)$ the product of the largest powers of defining functions of boundary hypersurfaces of $X_\bullet$ so that $[\cV_\bullet(X),\cV_\bullet(X)]\subset w\cV_\bullet(X)$; thus we can take $w=1$ for $\bullet=\bop,0,(\ebop),\tbop$, while $w=x$ (a boundary defining function of $X$) for $\bullet=\scop$, further $w=h$ for $\bullet=\bop\semi,0\semi$ and $w=h x$ for $\bullet=\schop$, and finally $w=\rho_\scface$, resp.\ $w=\rho_\sface$ for $\bullet=\scbtop$, resp.\ $\chop$. Recall that differential operators $A\in\Diff_\bullet^m(X)$ have principal symbols in $(P^m/w P^{m-1})({}^\bullet T^*X)$; for $s\in\R$, the space $\Psi_\bullet^s(X)$ of $\bullet$-pseudodifferential operators then consists of quantizations of symbols in $S^s({}^\bullet T^*X)$, and the principal symbol map ${}^\bullet\upsigma^s$ gives rise to a short exact sequence
\[
  0 \to w\Psi_\bullet^{s-1}(X) \hra \Psi_\bullet^s(X) \xra{{}^\bullet\upsigma^s} (S^s/w S^{s-1})({}^\bullet T^*X)\to 0.
\]
See below for explicit quantization maps in the various settings of interest. Given an element $a\in S^s({}^\bullet T^*X)$, one can define $A=\Op_\bullet(a)$ via a partition of unity in local coordinates by means of suitable variants of the standard quantization, and the space $\Psi_\bullet^s(X)$ is the space spanned by all such quantizations of elements of $S^s({}^\bullet T^*X)$ plus a space $w^\infty\Psi_\bullet^{-\infty}(X)$ of residual operators. Moreover, the space $\Psi_\bullet(X)=\bigcup_{s\in\R}\Psi_\bullet^s(X)$ is an algebra under composition, the principal symbol map is multiplicative, and the principal symbol of $i$ times a commutator is the Poisson bracket of principal symbols, just as in the case~\eqref{EqMSymbComm} of differential operators. All algebras under consideration here are invariant under conjugation by real powers of boundary defining functions of $X_\bullet$.

Now, ${}^\bullet\upsigma^s$ captures an element of $\Psi_\bullet^s(X)$ to leading order in the sense that its vanishing implies an extra amount $w$ of decay of its coefficients and also a gain in the differential order. In order to unify the two gains, fix on $\ol{{}^\bullet T^*}X$ a defining function $\rho_\infty\in\CI(\ol{{}^\bullet T^*}X)$ of fiber infinity ${}^\bullet S^*X$, and let ${}^\bullet w=(\pi^*w)\rho_\infty$ where $\pi\colon\ol{{}^\bullet T^*}X\to X$ is the projection. In the settings we consider here, ${}^\bullet w$ is thus the product of boundary defining functions of boundary hypersurfaces ${}^\bullet S^*X$ and ${}^\bullet W_1,\ldots,{}^\bullet W_K\subset\ol{{}^\bullet T^*}X$, and we say that the algebra $\Psi_\bullet(X)=\bigcup\Psi_\bullet^s(X)$ is symbolic to leading order at ${}^\bullet W_k$, $k=1,\ldots,K$. If we write $\rho_k\in\CI(X_\bullet)$ for defining functions of ${}^\bullet W_k$ for $k=1,\ldots,K$, and if $s\in\R$ and $\alpha=(\alpha_1,\ldots,\alpha_K)\in\R^K$ denote orders, we consider the space $\Psi_\bullet^{s,\alpha}(X)=(\prod_{k=1}^K\rho_k^{-\alpha_k})\Psi_\bullet^s(X)$ of weighted operators, and thus
\begin{equation}
\label{EqMOSES}
  0 \to \Psi_\bullet^{s-1,\alpha-(1,\ldots,1)}(X) \hra \Psi_\bullet^{s,\alpha}(X) \xra{{}^\bullet\upsigma^{s,\alpha}} (S^{s,\alpha}/{}^\bullet w S^{s,\alpha})({}^\bullet T^*X) \to 0.
\end{equation}
(We can also introduce weights $\beta$ at the remaining boundary hypersurfaces of $X_\bullet$, which amounts to adding an extra order $\beta$ to each of the terms in this exact sequence: the vanishing of the principal symbol does not imply additional decay at those boundary hypersurfaces.)

In the settings considered here, one may generalize the spaces $\Psi_\bullet^{s,\alpha}(X)$ by allowing the orders $s,\alpha$ to be variable along the respective boundary hypersurfaces. Thus, consider variable order functions $\sfs\in\CI({}^\bullet S^*X)$ and $\upalpha_k\in\CI({}^\bullet W_k)$, $k=1,\ldots,K$; denote arbitrary smooth extensions of these to functions on $\ol{{}^\bullet T^*}X$ by the same letters. One can then consider the symbol class
\begin{equation}
\label{EqMOVar}
  S^\sfw(\ol{{}^\bullet T^*}X)=\left(\rho_\infty^{-\sfs}\prod_{k=1}^K\rho_k^{-\upalpha_k}\right)\bigcap_{\delta\in(0,\frac12)}S^0_{1-\delta,\delta}(\ol{{}^\bullet T^*}X),
\end{equation}
where we follow the standard notation for $(\rho,\delta)$-symbol classes (with $\rho=1-\delta$) from \cite{HormanderFIO1}. Concretely, a smooth function $a$ on the interior of $\ol{{}^\bullet T^*}X$ lies in $S^0_{1-\delta,\delta}(\ol{{}^\bullet T^*}X)$ if and only if $V_1\cdots V_N a\in \rho_\infty^{-N\delta}(\prod_{k=1}^K w_k^{-N\delta})L^\infty_\loc(\ol{{}^\bullet T^*}X)$ for all $N\in\N_0$ and all $V_j\in\Vb(\ol{{}^\bullet T^*}X)$. Allowing for $\delta>0$ in~\eqref{EqMOVar} ensures that the definition is independent of the choice of the extensions of $\sfs,\upalpha_k$ to smooth functions on $\ol{{}^\bullet T^*}X$. We only consider orders $\sfs,\upalpha_k$ which are bounded above; thus, if $s_0\geq\sup\sfs$ and $\alpha_0=(\alpha_{1,0},\ldots,\alpha_{K,0})$ where $\alpha_{k,0}\geq\sup\upalpha_k$, then
\[
  S^\sfw(\ol{{}^\bullet T^*}X) \subset \bigcap_{\delta\in(0,\frac12)} S_{1-\delta,\delta}^{s_0,\alpha_0}({}^\bullet T^*X).
\]
For quantizations of symbols of class $(1-\delta,\delta)$, we have a short exact sequence as in~\eqref{EqMOSES} but with $s-(1-2\delta)$ and $\alpha-(1-2\delta,\ldots,1-2\delta)$ on the left, and with ${}^\bullet w^{1-2\delta}$ on the right. Defining $\Psi_\bullet^{\sfs,\upalpha}(X)$ (modulo the space $\Psi_\bullet^{-\infty,-\infty}(X)$ of residual operators) via quantizations in local coordinates of elements of $S^\sfw(\ol{{}^\bullet T^*}X)$, we then obtain a short exact sequence
\begin{align*}
  0 &\to \bigcap_{\delta\in(0,\frac12)}\Psi_\bullet^{\sfs-(1-2\delta),\upalpha-(1-2\delta,\ldots,1-2\delta)}(X) \hra \Psi_\bullet^{\sfs,\upalpha}(X) \\
    &\quad\hspace{6em} \xra{{}^\bullet\upsigma^{\sfs,\upalpha}} \biggl(S^\sfw/\bigcap_{\delta\in(0,\frac12)}S^{\sfw-(1-2\delta,(1-2\delta,\ldots,1-2\delta))}\biggr)(\ol{{}^\bullet T^*}X) \to 0.
\end{align*}

We now make this concrete. In the b-setting on $X$, and in local coordinates $x\geq 0$, $y\in\R^{n-1}$ near a boundary point, we can introduce fiber-linear coordinates on $\Tb^*X$ by writing b-covectors as $\xi_\bop\frac{\dd x}{x}+\eta_\bop\,\dd y$; the quantization of a symbol $a\in S^s(\Tb^*X)$, in local coordinates taking the form $a=a(x,y,\xi_\bop,\eta_\bop)$, with support in the coordinate chart, is then
\begin{equation}
\label{EqMOOpb}
\begin{split}
  (\Op_\bop(a)u)(x,y) &= (2\pi)^{-n}\iiiint_{\R\times\R^{n-1}\times[0,\infty)\times\R^{n-1}}\exp\Bigl[i\Bigl(\frac{x-x'}{x}\xi_\bop+(y-y')\cdot\eta_\bop\Bigr)\Bigr] \\
    &\qquad \times \chi\Bigl(\log\frac{x}{x'}\Bigr)\chi(|y-y'|)a(x,y,\xi_\bop,\eta_\bop)u(x',y')\,\dd\xi_\bop\,\dd\eta_\bop\frac{\dd x'}{x'}\dd y',
\end{split}
\end{equation}
where $\chi\in\CIc(\R)$ is identically $1$ near $0$ and supported in a small neighborhood of $0$. One can also consider weighted conormal symbols $a\in S^{s,\beta}(\Tb^*X)$, $s,\beta\in\R$, which in local coordinates satisfy the bound $|a|\lesssim x^\beta(1+|\xi_\bop|+|\eta_\bop|)^s$ together with all derivatives along $x\pa_x$, $\pa_y$, $\pa_{\xi_\bop}$, $\pa_{\eta_\bop}$, gaining a power of $(1+|\xi_\bop|+|\eta_\bop|)^{-1}$ upon differentiating in $\xi_\bop,\eta_\bop$. The b-differential order $s$ can be variable, i.e.\ an element of $\CI(\Sb^*X)$, with the quantization of such variable order symbols being given by the above formula still. See \cite[Appendix~A]{BaskinVasyWunschRadMink} for more on variable order b-operators. The quantization of a symbol $a=a(x,y,\xi_0,\eta_0)$ on ${}^0 T^*X$, where we write 0-covectors as $\xi_0\frac{\dd x}{x}+\eta_0\frac{\dd y}{x}$, is given by the formula~\eqref{EqMOOpb} with subscripts `$0$' instead of `$\bop$' and with $y-y'$ replaced by $\frac{y-y'}{x}$ and $\dd y'$ by $\frac{\dd y'}{x'^{n-1}}$.

Scattering ps.d.o.s on $\ol{\R^n}$, including with variable orders, are discussed in detail in \cite{VasyMinicourse}; see also \cite[Appendix~A.1]{HintzKdSMS} for an explicit quantization map. Thus, on a manifold $X$ with boundary, we have spaces $\Psisc^{s,r}(X)=x^{-r}\Psisc^s(X)$ of scattering ps.d.o.s; since $[\Vsc(X),\Vsc(X)]\subset x\Vsc(X)$, both the differential order $s$ and the scattering decay order $r$ can be taken to be variable, i.e.\ they can be replaced by $\sfs\in\CI(\Ssc^*X)$, $\sfr\in\CI(\ol{\Tsc^*_{\pa X}}X)$, giving rise to spaces $\Psisc^{\sfs,\sfr}(X)$. For scattering theory at nonzero real energies on asymptotically Euclidean (or more generally conic) spaces, one typically takes $\sfs$ to be constant, whereas $\sfr$ is variable and is above, resp.\ below some threshold at incoming, resp.\ outgoing momenta (matching the chosen energy). See \cite[Proposition~4.13]{VasyMinicourse}.

Semiclassical scattering ps.d.o.s\ were introduced by Vasy--Zworski \cite{VasyZworskiScl}; see again \cite[Appendix~A.1]{HintzKdSMS} for an explicit quantization map. In this setting, all three orders---the differential order $\sfs$, scattering decay order $\sfr$, and semiclassical order $\sfb$---can be variable, giving spaces of ps.d.o.s
\[
  \Psisch^{\sfs,\sfr,\sfb}(X),\qquad \sfs\in\CI({}^\schop S^*X),\quad \sfr\in\CI(\ol{{}^\schop T^*_{[0,1)\times\pa X}}X),\quad \sfb\in\CI(\ol{{}^\schop T^*_{h^{-1}(0)}}X).
\]
In this paper, only the case of variable $\sfr,\sfb$, but constant $\sfs$, is relevant.

The scattering-b-transition case is similar, see \cite[Appendix~A.3]{HintzKdSMS}, leading to spaces of operators
\begin{equation}
\label{EqMscbtPsi}
  \Psiscbt^{\sfs,\sfr,l,b}(X),\qquad \sfs\in\CI(\Sscbt^*X),\quad\sfr\in\CI(\ol{\Tscbt^*_\scface}X),\quad l,b\in\R,
\end{equation}
where $l,b$ are the orders at $\tface$ and $\zface$, as in the case of differential operators in~\eqref{EqMscbtDiff}. For a detailed description of the semiclassical cone setting, see \cite[\S3.2]{HintzConicProp}, where variable semiclassical orders are considered; this suffices for the present paper, though allowing for variable differential orders causes no additional issues. Thus, we have spaces of variable order operators
\begin{equation}
\label{EqMchPsi}
  \Psich^{\sfs,l,\alpha,\sfb}(X),\qquad \sfs\in\CI(\Sch^*X),\quad l,\alpha\in\R,\quad \sfb\in\CI(\ol{\Tch_\sface^*}X),
\end{equation}
where $l,\alpha$ are the orders at $\cface$ and $\tface$, as in the case of differential operators in~\eqref{EqMchDiff}.

The definition of edge-b-pseudodifferential operators depends in general on the properties of the boundary fibration. We only consider the case described in~\S\ref{SsMeb}. (Another case, of importance in the 3b-calculus, is discussed in detail in \cite[\S2.7]{Hintz3b}.) See \cite[Appendix~B]{MelroseVasyWunschDiffraction} for a discussion of the general case. In local coordinates $x\geq 0$, $z\geq 0$, $y\in\R^{n-2}$ on $M$, with the fibration of $x^{-1}(0)$ given by $(z,y)\mapsto y$, we write edge-b-covectors as $\xi_\ebop\frac{\dd x}{x}+\eta_\ebop\frac{\dd y}{x}+\zeta_\ebop\frac{\dd z}{z}$. The quantization of a symbol $a=a(x,y,z,\xi_\ebop,\eta_\ebop,\zeta_\ebop)$ with support in the local coordinate chart is then
\begin{equation}
\label{EqMebQuant}
\begin{split}
  &(\Op_\ebop(a)u)(x,y,z) = (2\pi)^{-n}\iiint_{[0,\infty)\times\R^{n-2}\times[0,\infty]}\iiint_{\R\times\R^{n-2}\times\R} \\
  &\qquad \exp\Bigl[i\Bigl(\frac{x-x'}{x}\xi_\ebop+\frac{y-y'}{x}\eta_\ebop+\frac{z-z'}{z}\zeta_\ebop\Bigr)\Bigr]\chi\Bigl(\log\frac{x}{x'}\Bigr)\chi\Bigl(\Bigl|\frac{y-y'}{x}\Bigr|\Bigr)\chi\Bigl(\log\frac{z}{z'}\Bigr) \\
  &\qquad\hspace{8em} \times a(x,y,z,\xi_\ebop,\eta_\ebop,\zeta_\ebop)u(x',y',z')\,\dd\xi_\ebop\,\dd\eta_\ebop\,\dd\zeta_\ebop\,\frac{\dd x'}{x'}\frac{\dd y'}{x'{}^{n-2}}\frac{\dd z'}{z'},
\end{split}
\end{equation}
with $\chi\in\CIc(\R)$ identically $1$ near $0$ and supported near $0$. When $|a|\lesssim x^{-\alpha_1}z^{-\alpha_2}(1+|\xi_\ebop|+|\eta_\ebop|+|\zeta_\ebop|)^s$, with the same bound for derivatives along $x\pa_x$, $\pa_y$, $z\pa_z$, as well as along $\xi_\ebop$, $\eta_\ebop$, $\zeta_\ebop$ with a gain of one order in $s$ for each of the latter three derivatives, this is a local description of elements of $\Psieb^{s,\alpha_1,\alpha_2}(M)$. The edge-b-differential order $s$ may be variable, i.e.\ one may replace it with $\sfs\in\CI(\Seb^*M)$. 

When $a$ is independent of $z$, the operator defined by~\eqref{EqMebQuant} is dilation-invariant in $z$. Upon appropriately modifying the symbol $a$, we may replace the factor $\exp(i\frac{z-z'}{z}\zeta_\ebop)$ by $(\frac{z}{z'})^{i\zeta_\ebop}$, and we can moreover drop the localizer $\chi(\log\frac{z}{z'})$ (via replacing $a$ by the convolution of $a$ in $\zeta_\ebop$ with the Fourier transform of $\chi$, and hence $a$ is entire in $\zeta_\ebop\in\C$ with suitable Paley--Wiener type bounds as $|\Im\zeta_\ebop|\to\infty$); then the Mellin-transformed normal operator family of $A=\Op_\ebop(a)$ at $H_2=z^{-1}(0)$ has Schwartz kernel
\begin{align*}
  \wh{N_{H_2}}(A,\lambda) &= (2\pi)^{-(n-1)}\iint_{\R\times\R^{n-2}} \exp\Bigl[i\Bigl(\frac{x-x'}{x}\xi_0+\frac{y-y'}{x}\eta_0\Bigr)\Bigr] \\
    &\quad\hspace{4em} \times \chi\Bigl(\log\frac{x}{x'}\Bigr)\chi\Bigl(\Bigl|\frac{y-y'}{x}\Bigr|\Bigr) a(x,y,0,\xi_0,\eta_0,\lambda)\,\dd\xi_0\,\dd\eta_0\,\Bigl|\frac{\dd x'}{x'}\frac{\dd y'}{x'{}^{n-2}}\Bigr|,
\end{align*}
where we relabeled $\xi_\ebop,\eta_\ebop$ as $\xi_0,\eta_0$. This shows that $\wh{N_{H_2}}(A,\lambda)\in\Psi_0^m(H_2)$, with principal symbol given by pullback of $\sigmaeb^m(A)$ along the inclusion ${}^0 T^*H_2\hra\Teb^*_{H_2}M$ as in the case of differential operators. In the high frequency regime $\lambda=-i\mu\pm h^{-1}$, we pass to semiclassical 0-momentum variables by writing covectors as $\xi_{0\semi}\frac{\dd x}{h x}+\eta_{0\semi}\frac{\dd y}{h x}$; the Schwartz kernel of $\wh{N_{H_2}}(A,-i\mu\pm h^{-1})$ is then
\begin{align*}
  &(2\pi)^{-(n-1)}\iint_{\R\times\R^{n-2}} \exp\Bigl[i\Bigl(\frac{x-x'}{h x}\xi_{0\semi}+\frac{y-y'}{h x}\eta_{0\semi}\Bigr)\Bigr]\chi\Bigl(\log\frac{x}{x'}\Bigr)\chi\Bigl(\Bigl|\frac{y-y'}{x}\Bigr|\Bigr) \\
  &\hspace{6em} \times a\bigl(x,y,0,h^{-1}\xi_{0\semi},h^{-1}\eta_{0\semi},h^{-1}(\pm 1-i h\mu)\bigr)\,\dd\xi_{0\semi}\,\dd\eta_{0\semi}\,\Bigl|\frac{\dd x'}{h x'}\frac{\dd y'}{(h x')^{n-2}}\Bigr|.
\end{align*}
This is the semiclassical 0-quantization of the symbol
\[
  (h,x,y,\xi_{0\semi},\eta_{0\semi}) \mapsto a\bigl(x,y,0,h^{-1}\xi_{0\semi},h^{-1}\eta_{0\semi},h^{-1}(\pm 1-i h\mu)\bigr),
\]
which in the case $a\in S^s(\Teb^*M)$ is an element of $S^{s,0,s}(\ol{{}^{0\semi}T^*}H_2)=h^{-s}S^s(\ol{{}^{0\semi}T^*}H_2)$, where the first order is the semiclassical 0-differential order, the second order is the order at the boundary $\ol{{}^{0\semi}T^*_{\pa H_2}}H_2$, and the third order is the semiclassical order, i.e.\ the order at $\ol{{}^{0\semi}T^*_{h^{-1}(0)}}H_2$. The principal symbols are related as in the case of differential operators, see~\eqref{EqMebSclSymb}. In particular, when $A=\Op_\ebop(a)$ is elliptic as an edge-b-operator at $H_2$, then $\wh{N_{H_2}}(A,-i\mu\pm h^{-1})$ is elliptic as a semiclassical 0-operator. This continues to hold also in the case that the edge-b-differential order is variable: if $\sfs\in\CI(\Seb^*M)$ and $A\in\Psieb^\sfs(M)$, then
\begin{equation}
\label{EqMebMTscl}
  \Bigl((0,1)\ni h \mapsto \wh{N_{H_2}}(A,-i\mu\pm h^{-1})\Bigr) \in \Psi_{0\semi}^{\sfs_\infty,0,\sfs_\semi}(H_2),
\end{equation}
where $\sfs_\infty\in\CI({}^{0\semi}S^*H_2)$ and $\sfs_\semi\in\CI(\ol{{}^{0\semi}T^*_{h^{-1}(0)}}H_2)$ are the pullbacks of $\sfs$ along the maps~\eqref{EqMebSclMaps}; and this operator family inherits ellipticity from $A$.

Finally, we recall the 3b-pseudodifferential algebra which was introduced in \cite{Hintz3b}. In the notation of~\S\ref{SsM3b}, we note that for $\cU\subset M$ with $\cU\cap\cT=\emptyset$, we have $\Ttb^*_\cU M=\Tb^*_\cU(M\setminus\cT)$, and correspondingly one can quantize 3b-symbols with support in $\cU$ using the above b-quantization map. When $\cU\cap\cD=\emptyset$ on the other hand, then $\Ttb^*_\cU M=\Tcu^*_\cU(M\setminus\cD)$ is the cusp cotangent bundle; in coordinates $T\geq 0$ and $x\in\R^n$ near $\cT^\circ$, this is the bundle dual to the cusp tangent bundle which is spanned by $T^2\pa_T$, $\pa_x$. Quantizations of cusp symbols, i.e.\ symbols on $\Tcu^*_\cU(M\setminus\cD)$, are described in detail in \cite{MazzeoMelroseFibred}, and an explicit quantization map is given by
\begin{align*}
  \Op_\cuop(a) &= (2\pi)^{-(n+1)}\iint_{\R\times\R^n} \exp\Bigl[i\Bigl(\frac{T-T'}{T^2}\sigma_\cuop+(x-x')\cdot\xi_\cuop\Bigr)\Bigr] \\
    &\hspace{8em}\times\chi\Bigl(\frac{T-T'}{T^2}\Bigr)\chi(|x-x'|)a(T,x,\sigma_\cuop,\xi_\cuop)\,\dd\sigma_\cuop\,\dd\xi_\cuop\,\Bigl|\frac{\dd T'}{T'^2}\dd x'\Bigr|
\end{align*}
where we write cusp covectors as $\sigma_\cuop\frac{\dd T}{T^2}+\xi_\cuop\,\dd x$. In both the b- and the cusp settings, the differential order (growth in $(\xi_\bop,\eta_\bop)$, resp.\ $(\sigma_\cuop,\xi_\cuop)$), may be variable. Finally, in a neighborhood of $\cT\cap\cD$ where $\rho_\cT=R=\frac{r}{t}\in[0,1)$ and $\rho_\cD=\frac{T}{R}=\frac{1}{r}\in[0,1)$ (with $T=t^{-1}$), we may write 3b-covectors (dually to~\eqref{EqM3btr}) as $\sigma_\tbop\frac{\dd t}{r}+\xi_\tbop\frac{\dd r}{r}+\eta_\tbop\,\dd\omega$ where $\omega\in\R^{n-1}$ denotes coordinates on $\cT\cap\cD\cong\Sph^{n-1}$. Given a symbol $a=a(\rho_\cT,\rho_\cD,\omega,\sigma_\tbop,\xi_\tbop,\eta_\tbop)$, of possibly variable differential order $\sfs\in\CI(\Stb^*M)$, with support in this chart, we define its quantization by
\begin{align*}
  &\Op_\tbop(a) = (2\pi)^{-n-1}\iiint_{\R\times\R\times\R^{n-1}} \exp\Bigl[i\Bigl(\frac{t-t'}{r}\sigma_\tbop+\frac{r-r'}{r}\xi_\tbop+(\omega-\omega')\cdot\eta_\tbop\Bigr)\Bigr] \\
  &\quad\hspace{4em} \times \chi\Bigl(\frac{t-t'}{r}\Bigr)\chi\Bigl(\frac{r-r'}{r}\Bigr) a\Bigl(\frac{r}{t},\frac{1}{r},\omega,\sigma_\tbop,\xi_\tbop,\eta_\tbop\Bigr)\,\dd\sigma_\tbop\,\dd\xi_\tbop\,\dd\eta_\tbop\,\Bigl|\frac{\dd t'}{r'}\frac{\dd r'}{r'}\dd\omega'\Bigr|.
\end{align*}

Alternatively, one can define $\Op_\tbop(a)$ as an edge-b-quantization of a conormal symbol; indeed, writing edge-b-covectors on $[0,1)_T\times[0,1)_R\times\Sph^{n-1}$ as $\sigma_\ebop\frac{\dd T}{T R}+\xi_\ebop\frac{\dd R}{R}+\eta_\ebop\,\dd\omega$, the expression of $a$ in these edge-b-coordinates is
\[
  a_\ebop \colon (T,R,\omega,\sigma_\ebop,\xi_\ebop,\eta_\ebop)\mapsto a\Bigl(R,\frac{T}{R},\omega,-\sigma_\ebop-R\xi_\ebop,\xi_\ebop,\eta_\ebop\Bigr),
\]
where we use $\sigma_\ebop\frac{\dd T}{T R}+\xi_\ebop\frac{\dd R}{R}+\eta_\ebop\,\dd\omega=\sigma_\tbop\frac{\dd t}{r}+\xi_\tbop\frac{\dd r}{r}+\eta_\tbop\,\dd\omega$ for $\sigma_\tbop=-R\xi_\ebop-\sigma_\ebop$, $\xi_\tbop=\xi_\ebop$, $\eta_\tbop=\eta_\ebop$ (and conversely $\sigma_\ebop=-\sigma_\tbop-R\xi_\tbop$). Note that typically this edge-b-symbol is only conormal near $T=R=0$ even if $a$ itself is a smooth symbol. (The close relationship between 3b-geometry near $\cT\cap\cD$ and edge-b-geometry on $[0,1)_T\times[0,1)_R\times\Sph^{n-1}$ is key also to the analysis of the $\cD$-normal operator; see \cite[\S3.3]{Hintz3b} for a detailed discussion.) That is, we may alternatively define $\Op_\tbop(a)$ to be equal to
\begin{align*}
  &(2\pi)^{-n-1}\iiint_{\R\times\R\times\R^{n-1}} \exp\Bigl[i\Bigl(\frac{T-T'}{T R}\sigma_\ebop+\frac{R-R'}{R}\xi_\ebop+(\omega-\omega')\cdot\eta_\ebop\Bigr)\Bigr] \\
  &\hspace{6em} \times \chi\Bigl(\frac{T-T'}{T R}\Bigr)\chi\Bigl(\frac{R-R'}{R}\Bigr)a_\ebop(T,R,\omega,\sigma_\ebop,\xi_\ebop,\eta_\ebop)\,\dd\sigma_\ebop\,\dd\xi_\ebop\,\dd\eta_\ebop.
\end{align*}
Both definitions produce operators with principal symbol $a$, though the two quantization maps typically differ already at the subprincipal level. Either map can be used to define quantizations of 3b-symbols with variable differential order. Quantizations of weighted symbols of class $\rho_\cD^{-\alpha_\cD}\rho_\cT^{-\alpha_\cT}S^s(\Ttb^*M)$ (or more generally $\cA^{(-\alpha_\cD,-\alpha_\cT)}S^s(\Ttb^*M)$) give rise to spaces of weighted 3b-ps.d.o.s
\[
  \Psitb^{s,(\alpha_\cD,\alpha_\cT)}(M).
\]

\subsubsection{Mixed operators}
\label{SssMOM}

In \cite[\S5]{HintzVasyScrieb}, spaces of mixed b-differential edge-b-pseu\-do\-dif\-fer\-en\-tial operators $\Diffb^k\Psieb^{s,\alpha}$ were defined on manifolds with corners; such operators are finite sums of compositions $Q P$ where $Q\in\Diffb^k$ and $P\in\Psieb^{s,\alpha}$. Since we shall not work with such operators directly here, we refer the reader to \cite[\S5]{HintzVasyScrieb} for a detailed discussion. Here, we introduce a similar mixed algebra in the 3b-setting, following \cite{HintzVasyScrieb,VasyPropagationCorners,VasyWaveOnAdS,MelroseVasyWunschDiffraction}. Thus, we work on the $(n+1)$-dimensional manifold $M=[M_0;\{\fp\}]$ as in~\S\ref{SsM3b}, with the blow-down map denoted $\upbeta\colon M\to M_0$. First, we consider commutators of 3b-ps.d.o.s with time dilation vector fields or time dilation operators; this is the pseudodifferential generalization of Lemmas~\ref{LemmaM3bbComm} and \ref{LemmaM3bbComm2}.

\begin{lemma}[Commutators III: pseudodifferential operators]
\label{LemmaMOMComm}
  Let $V$ be a time dilation vector field, and let $A\in\Psitb^s(M)$, $s\in\R$. Then there exist $A^\flat\in\Psitb^{s-1}(M)$ and $A^\sharp\in\Psitb^s(M)$ so that
  \begin{equation}
  \label{EqMOMComm}
    [V,A] = \rho_\cT A^\flat V + \rho_\cT A^\sharp.
  \end{equation}
  In particular, $[V,A]\in\Psitb^s(M)$. If $V\in\Diffb^1(M;E)$ is a time dilation operator acting on sections of a bundle $E\to M$, and $A\in\Psitb^s(M;E)$, the same conclusion holds, now with $A^\flat,A^\sharp$ acting on sections of $E$. If $A\in\cA^{(-\alpha_\cD,-\alpha_\cT)}\Psitb^s(M)$, we have $[V,A]\in\cA^{(-\alpha_\cD,-\alpha_\cT)}\Psitb^s(M)$, similarly when vector bundles are present. When $\sfs\in\CI(\Stb^*M)$ is a variable order and $A\in\Psitb^\sfs(M)$, then $[V,A]\in\Psitb^{\sfs+\eps}(M)$ for all $\eps>0$, similarly with conormal coefficients and vector bundles.
\end{lemma}
\begin{proof}
  Unlike in the case of 3b-differential operators, we cannot argue by induction on the order of $A$ here; we instead proceed directly using the $\cT$-normal operator and the Schwartz kernel $K$ of $A$. We begin with the case $A\in\Psitb^s(M)$. Let $\chi\in\CI(M)$ be a cutoff function which is identically $1$ near $\cT$ and supported in the preimage under $\upbeta$ of a coordinate chart $T\in[0,1)$, $X\in\R^{n-1}$, $|X|<1$. Then $(1-\chi)A\in\rho_\cT^N\Psitb^s(M)$ for any $N$; we take $N=2$. Since $\rho_\cT V\in\Vtb(M)$, we conclude that $[V,(1-\chi)A]\in\rho_\cT^{N-1}\Psitb^s(M)=\rho_\cT\tilde A^\sharp$ where $\tilde A^\sharp\in\rho_\cT^{N-2}\Psitb^s(M)=\Psitb^s(M)$. Thus, it suffices to prove the Lemma for $\chi A$ instead of $A$. A similar argument (multiplying by $1-\chi$ on the right) reduces this further to the study of $\chi A\chi$, which we relabel as $A$; now the Schwartz kernel $K$ of $A$ is supported in both factors in the coordinate chart.

  Next, in the coordinates $t=T^{-1}$ and $x=\frac{X}{T}$, and setting $\rho_\cT=\frac{\la x\ra}{t}$ (which is a local defining function of $\cT$), we can decompose $V=V_0+\tilde V$ where $V_0=-t\pa_t$ and $\tilde V$ is a b-vector field vanishing at $\cT$, so
  \[
    \tilde V=\rho_\cT \tilde a V_0 + \rho_\cT W,\qquad \tilde a\in\CI(M),\quad W\in\Vtb(M).
  \]
  Since $[\rho_\cT W,\Psitb^s(M)]\subset\rho_\cT\Psitb^s(M)$, and since $[\rho_\cT\tilde a V_0,A]=\rho_\cT\tilde a[V_0,A]-[A,\rho_\cT\tilde a]V_0$ with $[A,\rho_\cT\tilde a]\in\rho_\cT\Psitb^{s-1}(M)$, we conclude that it suffices to prove the Lemma for the time-dilation vector field $V_0=-t\pa_t$.

  Note first that $-t\pa_t=-\rho_\cT^{-1}\la x\ra\pa_t\in\rho_\cT^{-1}\Vtb(M)$, and hence $[-t\pa_t,A]\in\rho_\cT^{-1}\Psitb^s(M)$. As a first step, we improve this to
  \begin{equation}
  \label{EqMOMCommPf}
    [-t\pa_t,A] \in \Psitb^s(M).
  \end{equation}
  The Schwartz kernel $K_1$ of $[-t\pa_t,A]$ is equal to $-t\pa_t K+(t'\pa_{t'})^*K$ where $t$ and $t'$ denote the lifts of $t$ to the left and right factor of $M^2$, respectively, and $(t'\pa_{t'})^*$ is the adjoint of $t'\pa_{t'}$, which thus acts on right densities. Trivializing the right 3b-density bundle using $\la x'\ra^{n-1}|\dd t'\,\dd x'|$, we have $(t'\pa_{t'})^*=-t'\pa_{t'}-1$, and therefore
  \begin{equation}
  \label{EqMOMCommPfPf}
    K_1 = -K + R_\fp K,\qquad R_\fp := -(t\pa_t+t'\pa_{t'}).
  \end{equation}
  Observe then that passing to the coordinates $T=t^{-1}$, $X=\frac{x}{t}$, $T'=t'{}^{-1}$, $X'=\frac{x'}{t'}$ on $(M^\circ)^2$, the vector field
  \[
    R_\fp = T\pa_T+X\pa_X+T'\pa_{T'}+X'\pa_{X'}
  \]
  is the generator of dilations in $(T,X,T',X')$; using the terminology of \cite[Definitions~4.2 and 4.3]{Hintz3b}, this is tangent to the 3b-diagonal $\diag_{\tbop,\flat}$ inside the 3b-double space $M^2_\tbop$, and it is the b-normal vector field at $\ff_{\cT,\flat}\subset M^2_\tbop$. Therefore,
  \begin{equation}
  \label{EqMOMCommPfR}
     R_\fp \colon \Psitb^s(M) \to \rho_\cT\Psitb^s(M),
  \end{equation}
  since the space (of Schwartz kernels of elements of) $\Psitb^s(M)$ consists of conormal distributions to $\diag_{\tbop,\flat}$ which are smooth down to $\ff_{\cT,\flat}$, the gain of $\rho_\cT$ arising from the smoothness. Hence,~\eqref{EqMOMCommPfPf} implies~\eqref{EqMOMCommPf}.

  Finally, we improve~\eqref{EqMOMCommPf} to~\eqref{EqMOMComm}. To this end, we write
  \begin{equation}
  \label{EqMOMCommPf2}
    [-t\pa_t,A] = -t[\pa_t,A] + ([A,t]t^{-1})t\pa_t.
  \end{equation}
  Since $(\pa_{t'})^*=-\pa_{t'}$, the Schwartz kernel of $-t[\pa_t,A]$ is
  \[
    -t(\pa_t+\pa_{t'})K=R_\fp K+\frac{t'-t}{t'}t'\pa_{t'}K = R_\fp K - \frac{t'-t}{t'} (t'\pa_{t'})^* K - \frac{t'-t}{t'}K.
  \]
  The first term lies in $\rho_\cT\Psitb^s(M)$ by~\eqref{EqMOMCommPfR}. In the third term, we note that $q:=\frac{t'-t}{t'}=\frac{T-T'}{T}$ vanishes at $\ff_{\cT,\flat}$ (in the notation of \cite[\S4]{Hintz3b}), hence this term lies in $\rho_\cT\Psitb^s(M)$ as well. The second term is $q$ times the Schwartz kernel of $A\circ(-t\pa_t)$ and thus of the form of the first term in~\eqref{EqMOMComm}.

  In the second term of~\eqref{EqMOMCommPf2}, we note that the Schwartz kernel of $[A,t]$ is $(t-t')K=(\frac{1}{T}-\frac{1}{T'})K=T^{-1}\frac{T'-T}{T'}K$. Since, as noted before, $\frac{T'-T}{T'}$ vanishes at $\ff_{\cT,\flat}$, we have $T^{-1}\cdot\frac{T'-T}{T'}=\rho_{\ff_{\cD,\flat}}^{-1}\rho_{\ff_{\cT,\flat}}^{-1}\cdot\rho_{\ff_{\cT,\flat}}\cdot a$ where $a$ is smooth down to $\ff_{\cT,\flat}\cup\ff_{\cD,\flat}$ and blows up at most simply at the other boundary hypersurfaces of $M^2_{\tbop,\flat}$, at which however $K$ vanishes to infinite order. Therefore, $[A,t]\in\rho_\cD^{-1}\Psitb^{s-1}(M)$. In view of $t^{-1}\in\rho_\cT\rho_\cD\CI(M)$, this shows that $[A,t]t^{-1}\in\rho_\cT\Psitb^{s-1}(M)$; composed with $t\pa_t$, this is of the same class as the first term in~\eqref{EqMOMComm}. This completes the proof of the Lemma when $A\in\Psitb^s(M)$. The proof requires only minimal adaptations in the presence of a vector bundle.

  When $A\in\cA^{(-\alpha_\cD,-\alpha_\cT)}\Psitb^s(M)$, the Schwartz kernel of $[-A,t\pa_t]$ is again given by~\eqref{EqMOMCommPfPf}, and the conclusion follows from the fact that $R_\fp$ maps $\cA^{(-\alpha_\cD,-\alpha_\cT)}\Psitb^s(M)$ into itself. In the variable order case, application of the vector field $R_\fp$ increases the order of a variable order conormal distribution (such as the Schwartz kernel of $A$) by a logarithmic amount, and thus a fortiori by less than any positive amount $\eps>0$.
\end{proof}

\begin{definition}[Mixed $\tbop;\bop$-pseudodifferential operators]
\label{DefMOM}
  Let $s\in\R$ or $s\in\CI(\Stb^*M)$ and $k\in\N_0$. Then $\Diffb^k\Psitb^s(M)$ consists of all operators which are locally finite sums $\sum_i Q_i P_i$ where $Q_i\in\Diffb^k(M)$ and $P_i\in\Psitb^s(M)$.
\end{definition}

More general spaces such as $\Diffb^k\Psitb^{s,(\alpha_\cD,\alpha_\cT)}(M)$ and $\cA^{(-\alpha_\cD,-\alpha_\cT)}\Diffb^k\Psitb^s(M)$ are defined similarly.

\begin{lemma}[Commuting the two factors]
\label{LemmaMOMF}
  Lemma~\usref{LemmaM3bbF} remains valid if one replaces $\Difftb$ by $\Psitb$ throughout, including in the case of variable 3b-differential orders.
\end{lemma}
\begin{proof}
  In the constant order case and in view of Lemma~\ref{LemmaMOMComm}, the proof of Lemma~\ref{LemmaM3bbF} applies verbatim if we replace $\Difftb$ by $\Psitb$. If $P\in\Psitb^\sfs(M)$ is a variable order operator, we need to revisit aspects of the proof. Let $Q\in\Diffb^k(M)$. In the case $k=1$, it suffices to consider $Q=a V$ where $a\in\CI(M)$ and $V$ is a time dilation vector field, in which case $[P,Q]=[P,a]V+a[P,V]\in\Psitb^{\sfs+\eps}(M)$ by Lemma~\ref{LemmaMOMComm}; but we can write any element of $\Psitb^{\sfs+\eps}(M)$ as a sum of terms $P'Q'$ where $P'\in\Psitb^{\sfs-1+\eps}(M)\subset\Psitb^\sfs(M)$ and $Q'\in\Difftb^1(M)\subset\Diffb^1(M)$. This establishes the desired decomposition of $Q P=P Q-[P,Q]$. For the decomposition of $P Q=Q P+[P,Q]$, one further writes $[P,Q]=V[P,a]+[[P,a],V]+a[P,V]$ and notes that $[[P,a],V]\in[\Psitb^{\sfs-1+\eps}(M),V]\subset\Psitb^{\sfs-1+2\eps}(M)\subset\Psitb^\sfs(M)$ by Lemma~\ref{LemmaMOMComm}. We leave the minor modifications required for adapting the inductive argument for $Q\in\Diffb^k(M)$ with $k\geq 2$ to the reader.
\end{proof}

\subsection{Function spaces}
\label{SsMF}

Corresponding to each of the algebras of differential operators introduced above, there is a corresponding scale of weighted Sobolev spaces with nonnegative integer regularity; using pseudodifferential operators, one can more generally define spaces with variable orders. Consider first weighted b-Sobolev spaces on the manifold $X$ with boundary; assume that $X$ is compact. Fix on $X$ a positive weighted b-density $\nu$, i.e.\ $\nu=x^\mu\nu_0$ where $\mu\in\R$ and $0<\nu_0\in\CI(X;\Omegab X)$ where $\Omegab X\to X$ is the density bundle associated with $\Tb X$; that is, in local coordinates $x\geq 0$ and $y\in\R^{n-1}$, we have $\nu_0=a(x,y)|\frac{\dd x}{x}\dd y|$ where $0<a\in\CI$. We then set $\Hb^0(X)=L^2(X,\nu)$. (When the need arises to emphasize the choice of volume density, we shall write $\Hb^0(X,\nu)$.) For $\alpha\in\R$ and $s\in\N_0$, we set
\[
  \Hb^{s,\alpha}(X) = x^\alpha\Hb^s(X) = \{ x^\alpha u \colon u\in\Hb^s(X) \},
\]
where $\Hb^s(X)=\{u\in L^2(X)\colon A u\in L^2(X)\}$ for all $A\in\Diffb^s(X)$. For even $s$, this can be given the structure of a Hilbert space by fixing $A\in\Diffb^s(X)$ with an elliptic principal symbol and setting
\begin{equation}
\label{EqMFHb}
  \| u \|_{\Hb^{s,\alpha}(X)}^2 := \| x^{-\alpha}u \|_{L^2(X)}^2 + \| x^{-\alpha}A u \|_{L^2(X)}^2;
\end{equation}
for general $s\in\N$, one takes $\|x^{-\alpha}u\|_{L^2(X)}^2+\sum\|x^{-\alpha}A u\|_{L^2(X)}^2$ where one sums over a finite collection of $A$ which span $\Diffb^s(X)$ over $\CI(X)$. More generally, when $s\in\R$ is nonnegative, or when $\sfs\in\CI(\Sb^*X)$ is a variable order function with $\sfs\geq 0$, we define $\Hb^{s,\alpha}(X)$ and $\Hb^{\sfs,\alpha}(X)$ in the same manner, except we work with (elliptic) $A\in\Psib^s(X)$ and $A\in\Psib^\sfs(X)$, respectively. For real $s<0$, we can define $\Hb^{s,\alpha}(X)$ as the dual space, with respect to the $L^2(X)$ inner product, of the already defined space $\Hb^{-s,-\alpha}(X)$; in particular, it is a Hilbert space. Equivalently, and more directly, the space $\Hb^{s,\alpha}(X)$ consists of all distributions on $X^\circ$ of the form $u=u_0+A u_1$ where $u_0,u_1\in\Hb^{0,\alpha}(X)$ and $A\in\Psib^{-s}(X)$; if we fix $A$ to be elliptic, an equivalent norm of $u\in\Hb^{s,\alpha}(X)$ is then $\inf(\|u_0\|_{\Hb^{0,\alpha}(X)}+\|u_1\|_{\Hb^{0,\alpha}(X)})$, where the infimum is taken over all $u_0,u_1$ with $u=u_0+A u_1$. The case of negative variable orders is analogous. More generally, given $\sfs\in\CI(\Sb^*X)$, set $s_0:=\inf\sfs$; then $\Hb^{\sfs,\alpha}(X)$ consists of all $u\in\Hb^{s_0,\alpha}(X)$ so that $A u\in\Hb^{0,\alpha}(X)$ for an elliptic element of $\Psib^\sfs(X)$. (See \cite[Appendix~A]{MelroseVasyWunschDiffraction} for functional analytic background.)

If $E\to X$ is a vector bundle, then the space $\Hb^{s,\alpha}(X;E)$ can be defined to consist of tuples of elements of $\Hb^{s,\alpha}(X)$ in local trivializations of $E$; an equivalent definition can be given via testing with (elliptic) b-differential operators acting on sections of $E$. Elements of $\Psib^{s,\alpha}(M;E,F)$ then define continuous linear maps $\Hb^{t,\beta}(M;E)\to\Hb^{t-s,\beta-\alpha}(M;F)$ for all $t,\beta\in\R$; this remains true also when $s,t$ are variable orders.

A fundamental result for elliptic b-differential (or b-pseudodifferential) operators $P\in\Diffb^m(M;E,F)$ states that if $\alpha\notin\Re\specb(P)$, then $P\colon\Hb^{s,\alpha}(M;E)\to\Hb^{s-m,\alpha}(M;F)$ is a Fredholm operator for all $s\in\R$, whereas this fails when $\alpha\in\Re\specb(P)$. See \cite[\S5.17]{MelroseAPS}, and \cite[\S2]{HintzUnDet} for a quick summary of the proof. Moreover, the Fredholm indices $I_\alpha$, $I_\beta$ of $P$ for different choices of weights $\alpha,\beta$ are related by a relative index theorem \cite[\S6.2]{MelroseAPS}; in particular, if there exists a point in $\Re\specb(P)$ between two weights $\alpha<\beta$ not in $\Re\specb(P)$, then $I_\alpha>I_\beta$.

Next, if $\Omega\subset X$ is relatively open, we define the closed subspace
\begin{equation}
\label{EqMFHbsupp}
  \Hsupp_\bop^{s,\alpha}(\ol\Omega) := \{ u\in\Hb^{s,\alpha}(X) \colon \supp u\subset\ol\Omega \}
\end{equation}
of supported distributions, and the quotient space
\begin{equation}
\label{EqMFHbext}
  \Hext_\bop^{s,\alpha}(\Omega) := \{ u|_\Omega \colon u\in\Hb^{s,\alpha}(X) \} \cong \Hb^{s,\alpha}(X) / \Hbsupp^{s,\alpha}(X\setminus\Omega)
\end{equation}
of extendible distributions, using the terminology of H\"ormander \cite[Appendix~B]{HormanderAnalysisPDE3}. Finally, when $X$ is not compact, one can define $H_{\bop,\cop}^{s,\alpha}(X)$, $s\geq 0$, to consist of distributions on $X^\circ$ whose support has compact closure in $X$ and for which the expression~\eqref{EqMFHb} is finite (where $A$ is required to be properly supported), and $H_{\bop,\loc}^{s,\alpha}(X)$ consists of all distributions on $X^\circ$ so that $\phi u\in H_{\bop,\cp}^{s,\alpha}(X)$ for all $\phi\in\CIc(X)$. For negative $s$, these spaces can be defined via duality.

Returning to the case that $X$ is compact, Rellich's compactness theorem states that the inclusion
\begin{equation}
\label{EqMFHbRellich}
  \Hb^{s,\alpha}(X) \hra \Hb^{s',\alpha'}(X)
\end{equation}
is compact if (and only if) $s>s'$ and $\alpha>\alpha'$. Similarly, in the notation of~\eqref{EqMFHbsupp}--\eqref{EqMFHbext}, the inclusions $\Hbsupp^{s,\alpha}(\ol\Omega)\hra\Hbsupp^{s',\alpha'}(\ol\Omega)$ and $\Hbext^{s,\alpha}(\Omega)\hra\Hbext^{s',\alpha'}(\Omega)$ are compact. These are special cases of a general result for compact inclusions of weighted Sobolev spaces on manifolds with bounded geometry (see e.g.\ \cite[Corollary~4.10]{AlbinLectureNotes}), which states that an inclusion of two weighted Sobolev spaces is compact if both the differential order and the weight at infinity of the domain are stronger than those of the codomain. (See also \cite[Lemma~4.23]{Hintz3b} in the 3b-setting, and \cite[\S4.5]{Hintz3b} for a brief discussion of the case of variable differential orders.)

The wave front set $\WFb^{s,\alpha}(u)\subset\Sb^*X$ (equivalently, a conic subset of $\Tb^*X\setminus o$, with $o$ denoting the zero section) is defined for $u\in\Hb^{-\infty,\alpha}(X)=\bigcup_{s\in\R}\Hb^{s,\alpha}(X)$ (or $H_{\bop,\loc}^{-\infty,\alpha}(X)$ in the noncompact case) as the complement of the set of $\varpi\in\Sb^*X$ for which there exists an operator $A\in\Psib^0(X)$ which is elliptic at $\varpi$ and for which $A u\in\Hb^{s,\alpha}(X)$. Thus, $\WFb^{s,\alpha}(u)=\emptyset$ implies $u\in\Hb^{s,\alpha}(X)$. Note that $u$ needs to have weight $\alpha$ for $\WFb^{s,\alpha}(u)$ to be well-defined; by contrast, the b-regularity of $u$ may be arbitrarily negative. The reason is that the algebra $\bigcup_{s,\alpha\in\R}\Psib^{s,\alpha}(X)$ is commutative to leading order only in the b-differential sense but not in the weight; consequently, one can construct symbolic elliptic parametrices with errors which are trivial (order $-\infty$) in the b-differential order sense, but not in the decay order.

We make the same definitions for all of the other algebras of degenerate (pseudo)dif\-fer\-en\-tial operators on $X$ or appropriate manifolds with corners $M$ discussed previously. Thus, we have spaces
\[
  \Hsc^{s,r}(X),\quad
  H_{\scop,h}^{s,r,b}(X).
\]
For $r\in\R$, we have $\Hsc^{s,r}(X)=x^r\Hsc^s(X)$, which in the case $X=\ol{\R^n}$ is the standard weighted Sobolev space $\la z\ra^{-r}H^s(\R^n_z)$; and for $r,b\in\R$, we have $H_{\scop,h}^{s,r,b}(X)=x^r h^b H_{\scop,h}^s(X)$, which in the case $X=\ol{\R^n}$ is the semiclassical weighted Sobolev space on $\R^n$, with $h$-dependent norm
\[
  \|u\|_{H_{\scop,h}^{s,r,b}(\ol{\R^n})}^2 = \| \la h D_z\ra^s \la z\ra^{-r} h^{-b} u \|_{L^2(\R^n)}^2.
\]
In the case that $s,r,b$ are variable orders (as discussed in~\S\ref{SsMO}), these spaces can only be defined microlocally, and their norms need to be defined by means of testing with variable order pseudodifferential operators. Note that $H_{\scop,h}^{s,r,b}(X)=\Hsc^{s,r}(X)$ as a set for any $h>0$, but the norms are not uniformly equivalent as $h\searrow 0$ (unless $s=0$). For $u\in\Hsc^{-\infty,-\infty}(X)$, the set $\WFsc^{s,r}(u)$ is well-defined, and it is a subset of the union of all boundary hypersurfaces of $\ol{\Tsc^*}X$ at which scattering ps.d.o.s are commutative to leading order; thus,
\[
  \WFsc^{s,r}(u) \subset \Ssc^*X \cup \ol{\Tsc^*_{\pa X}}X.
\]
In the present paper, the wave front set at a point $\varpi\in\Tsc^*_{\pa X}X$ are of particular importance; the absence of such wave front set signals decay at $\pa X$ of order $r$ when localizing $u$ at frequency $\varpi$. In the semiclassical setting, for $u\in H_{\scop,h}^{-\infty,-\infty,-\infty}(X)$,
\[
  \WFsch^{s,r,b}(u) \subset {}^\schop S^*X \cup \ol{\Tsch^*_{[0,1)\times\pa X}}X \cup \ol{\Tsch^*_{h^{-1}(0)}}X.
\]
The wave front set over $h^{-1}(0)$ measures, microlocally, the size of $u$ relative to $h^b$, and its absence signals smallness of microlocalizations of $u$.

We similarly have weighted $\scbtop$-Sobolev spaces
\begin{equation}
\label{EqMFscbt}
  H_{\scbtop,\sigma}^{s,r,l,b}(X),\qquad \sigma\in\pm[0,1),
\end{equation}
where the orders $s,r$ may be variable as in~\eqref{EqMscbtPsi}. For any $\sigma\neq 0$, this space is equal to $\Hsc^{s,r}(X)$ as a set, but its norm, which is defined via testing with a fixed elliptic operator in $\Psiscbt^{s,r,l,b}(X)$, is not uniformly equivalent to the $\Hsc^{s,r}(X)$-norm as $\sigma\to 0$. In fact, for $b=0$, the norm on the space~\eqref{EqMFscbt} for $\sigma=0$ is the $\Hb^{s,l}(X)$-norm (and is thus in particular not the norm on a scattering Sobolev space anymore). Furthermore, semiclassical cone Sobolev spaces
\begin{equation}
\label{EqMFch}
  H_{\cop,h}^{s,l,\alpha,b}(X)
\end{equation}
possibly with variable orders $s,b$ as in~\eqref{EqMchPsi}, can be defined analogously. We also have
\[
  \Heb^{s,\alpha}(M),\qquad
  \Htb^{s,(\alpha_\cD,\alpha_\cT)}(M)
\]
in the edge-b- and 3b-settings discussed in~\S\S\ref{SsMeb}--\ref{SsM3b}; the differential order $s$ may be variable in both cases. We shall also use (semiclassical) 0-Sobolev spaces
\begin{equation}
\label{EqMF0}
  H_0^{s,\gamma}(X),\qquad
  H_{0,h}^{s,\gamma,b}(X),
\end{equation}
where $\gamma\in\R$, while the (semiclassical) 0-differential order $s$ and the semiclassical order $b$ may be variable. In local coordinates $x\geq 0$ and $y\in\R^{n-1}$, the norms on these spaces (in the case of $s\in\N_0$) are defined with respect to testing with the vector fields $x\pa_x$, $x\pa_y$, and $h x\pa_x$, $h x\pa_y$, respectively. See \cite{MazzeoMelroseHyp,HintzVasyScrieb} for further details.

Corresponding to the normal operator homomorphisms from the various (pseudo)dif\-fer\-en\-tial algebras into model algebras---such as the normal operator map $\Diffscbt^{s,r,0,b}(X)\to\Diff_{\scop,\bop}^{s,r,b}(\tface)$ at the transition face $\tface\subset X_\scbtop$---there are norm equivalences for elements of the corresponding Sobolev spaces---such as \cite[Proposition~2.21]{Hintz3b}. For b-Sobolev spaces, this is discussed in \cite[\S3.1]{VasyMicroKerrdS} and \cite[\S2.1]{Hintz3b}; for $\scbtop$-spaces, see \cite[Proposition~2.21]{Hintz3b}; for $\chop$-Sobolev spaces, see \cite[Proposition~2.29]{Hintz3b} and \cite[Corollary~3.7]{HintzConicProp}; and for 3b-Sobolev spaces, see \cite[Propositions~4.24 and 4.26]{Hintz3b} for such norm equivalences at $\cT$ and $\cD$, respectively, and also \cite[\S4.5]{Hintz3b} for the variable order case.

In the present paper, we need a further such equivalence statement, namely in the edge-b-setting of~\S\ref{SsMeb}, which generalizes the constant order result given in \cite[Lemma~2.10]{HintzVasyScrieb}. Thus, $M$ is compact, with two boundary hypersurfaces $H_1,H_2$ intersecting at $Y$. We work locally on
\[
  M = [0,1)_x\times[0,1)_z\times\R^m_y,\qquad \mu=\Bigl|\frac{\dd z}{z}\Bigr|\nu,
\]
where $\nu$ is a weighted positive b-density on $H_2=z^{-1}(0)$; and $H_1=x^{-1}(0)$ is fibered via the map $(z,y)\mapsto y$.

\begin{prop}[Equivalence of norms]
\label{PropMFebEq0}
   Let $\sfs\in\CI(\Seb^*M)$, $\alpha_1,\alpha_2\in\R$. Let $\sfs_\infty\in\CI({}^{0\semi}S^*H_2)$ and $\sfs_\semi\in\CI(\ol{{}^{0\semi}T^*_{h^{-1}(0)}}H_2)$ denote the pullbacks of $\sfs$ along the maps~\eqref{EqMebSclMaps}. Define the Mellin transform of $u=u(x,y,z)$ in $z$ by $\hat u((x,y),\lambda)=\int_0^\infty z^{-i\lambda}u(x,y,z)\,\frac{\dd z}{z}$. Recalling that $\sfs_\infty$ is an $h$-independent smooth function on ${}^0 S^*H_2$, put
   \begin{align*}
     \|u\|_{\sfs,(\alpha_1,\alpha_2)}^2 &:= \int_{[-1,1]} \| \hat u(-,-i\alpha_2+\lambda) \|_{H_0^{\sfs_\infty,\alpha_1}(H_2,\nu)}^2 \,\dd\lambda \\
       &\qquad + \sum_\pm\int_{[1,\infty)} \| \hat u(-,-i\alpha_2\pm\lambda) \|_{H_{0,|\lambda|^{-1}}^{\sfs_\infty,\alpha_1,\sfs_\semi}(H_2,\nu)}^2\,\dd\lambda.
   \end{align*}
   Let $\eps>0$. Then there exists $z_0>0$ so that for all $\chi\in\CIc(M)$ with $\supp\chi\subset\{z\leq z_0\}$, there exists a constant $C>1$ so that for all $u$,
   \begin{equation}
   \label{EqMFebEq0}
     C^{-1}\| \chi u \|_{\Heb^{\sfs-\eps,(\alpha_1,\alpha_2)}(M,\mu)} \leq \|\chi u\|_{\sfs,(\alpha_1,\alpha_2)} \leq  C\|\chi u\|_{\sfs+\eps,(\alpha_1,\alpha_2)}.
   \end{equation}
   Here, $\Heb^{\sfs,(\alpha_1,\alpha_2)}(M,\mu)=x^{\alpha_1}z^{\alpha_2}\Heb^\sfs(M,\mu)$. If $\sfs$ is dilation-invariant in $z$ on $\supp\chi$, then the estimates~\eqref{EqMFebEq0} hold also for $\eps=0$.
\end{prop}
\begin{proof}
  By definition of the Mellin transform, we may reduce to the case $\alpha_2=0$ by considering $z^{-\alpha_2}u$ instead of $u$; we may also divide by $x^{\alpha_1}$ (which commutes with the Mellin transform). For $\sfs=0$, the claim is a restatement of Plancherel's theorem. Since the difference of $\sfs$ and its $z$-dilation invariant extension is smaller (in $L^\infty$) than any $\eps>0$ in the region $z\leq z_0$ when $z_0>0$ is small enough (depending on $\eps$), we only need to show~\eqref{EqMFebEq0} when $\sfs$ is dilation-invariant and $\eps=0$.

  Consider first the case $\sfs\geq 0$. Fix an elliptic operator $A\in\Psieb^\sfs(M)$ which is dilation-invariant on $\supp\chi$; then the edge-b-Sobolev norm is
  \begin{equation}
  \label{EqMFebEq0Pf}
    \|\chi u\|_{\Heb^\sfs(M)}^2 = \|\chi u\|_{L^2(M)}^2 + \|A(\chi u)\|_{L^2(M)}^2.
  \end{equation}
  (Any two choices of $A$ give equivalent norms.) We pass to the Mellin transform in both terms on the right hand side; in the second term, this gives the squared $L^2(\R_\lambda)$-norm of $\|\wh{N_{H_2}}(A,\lambda)\wh{\chi u}(-,\lambda)\|_{L^2(H_2)}^2$. But $\wh{N_{H_2}}(A,\lambda)\in\Psi_0^{\sfs_\infty}(H_2)$ for bounded $\lambda\in\R$, and $(0,1)\ni h\mapsto\wh{N_{H_2}}(A,\pm h^{-1})\in\Psi_{0\semi}^{\sfs_\infty,0,\sfs_\semi}(H_2)$ in the semiclassical regime. Moreover, these operators are elliptic since $A$ is; see the discussion leading up to~\eqref{EqMebMTscl}. But this means that~\eqref{EqMFebEq0Pf} is precisely given by $\|u\|_{\sfs,(\alpha_1,\alpha_2)}^2$ (up to equivalence of norms due to the arbitrariness of the choice of $A$).

  By duality, we now obtain~\eqref{EqMFebEq0} also for $\sfs\leq 0$. When $\sfs$ does not have a sign, we pick $s_0<\inf_{\supp\chi}\sfs$ and note that $\|\chi u\|_{\Heb^\sfs(M)}^2$ is equivalent to $\|\chi u\|_{\Heb^{s_0}(M)}^2+\|A(\chi u)\|_{L^2(M)}$ where $A\in\Psieb^\sfs(M)$; the estimate~\eqref{EqMFebEq0} (with $\eps=0$) applies to the first term since $s_0$ is constant, and the previous arguments apply to the second term without change.
\end{proof}

We recall the Sobolev spaces corresponding to the mixed $\ebop;\bop$-algebra in the setting of~\S\ref{SsMeb}; see \cite[\S5.2]{HintzVasyScrieb}. In the case that $M$ is compact, these spaces are denoted
\[
  H_{\ebop;\bop}^{(s;k),\alpha}(M),
\]
where $s\in\R$ or $s\in\CI(\Seb^*M)$, $k\in\N_0$, and $\alpha$ is a vector of weights for each boundary hypersurface of $M$. (In \cite{HintzVasyScrieb} only the case of constant orders is discussed. The case of variable orders requires only minor modifications similar to those in~\S\ref{SssMOM}; we omit the details here.) The space $H_{\ebop;\bop}^{s,k;\alpha}(M)$ consists of all $u\in\Heb^{s,\alpha}(M)$ so that $A u\in\Heb^{s,\alpha}(M)$ for all $A\in\Diffb^k(M)$. On noncompact $M$, analogous local spaces or spaces of compactly supported distributions can be defined in the usual manner; and the addition of vector bundles is routine as well.

In the 3b-setting considered in~\S\ref{SsM3b}, we analogously define
\begin{equation}
\label{EqMFtbb}
  H_{\tbop;\bop}^{(s;k),\alpha}(M)
\end{equation}
for $s\in\R$ or $s\in\CI(\Stb^*M)$, $k\in\N_0$, and $\alpha=(\alpha_\cD,\alpha_\cT)\in\R^2$ to consist of all $u\in\Htb^{s,\alpha}(M)$ so that $A u\in\Htb^{s,\alpha}(M)$ for all $A\in\Diffb^k(M)$. By virtue of the commutation result Lemma~\ref{LemmaMOMF}, elements of $\Diffb^{k'}\Psitb^{s',\alpha'}(M)$ map the space~\eqref{EqMFtbb} into the analogous space with orders reduced by $s'$, $k'$, $\alpha'$. For $u\in H_{\tbop;\bop}^{(-\infty;k),\alpha}(M)=\bigcup_{s\in\R} H_{\tbop;\bop}^{(s;k),\alpha}(M)$, we define
\begin{equation}
\label{EqMFtbbWF}
  \WF_{\tbop;\bop}^{(\sfs;k),\alpha}(u) \subset \Stb^*M
\end{equation}
as the complement of the set of $\varpi\in\Stb^*M$ for which there exists an operator $A\in\Psitb^\sfs(M)$ which is elliptic at $\varpi$ and for which $A u\in H_\tbop^{\sfs,\alpha}(M)$. Note that $u$ needs to possess $k$ degrees of b-regularity in order for the set~\eqref{EqMFtbbWF} to be well-defined.

\section{Geometric and analytic setup}
\label{SG}

In~\S\ref{SsGS}, we describe the setup for the analysis of stationary (time translation invariant) wave type equations, starting with a description of the underlying spacetime manifold and the class of Lorentzian metrics of interest on it (Definition~\ref{DefGSG}) in~\S\ref{SssGSG}, followed by the definition of the class of stationary wave type operators we shall consider (Definition~\ref{DefGSO}) and the spectral assumptions placed on them (Definition~\ref{DefGSOSpec}) as required for our subsequent analysis in~\S\ref{SSt}.

The compactification of the spacetime manifold on which the analysis of non-stationary wave type operators will take place is more involved; it is described at the beginning of~\S\ref{SsGA} (Definition~\ref{DefGAMfd}). In~\S\ref{SssGAG} then, we introduce the class of (non-stationary) admissible asymptotically flat metrics (Definition~\ref{DefGAG}), and in~\S\ref{SssGAW} the class of (non-stationary) wave type operators to which the sharp solvability theory developed in~\S\ref{SW} applies.

Readers only interested in our main results on stationary wave operators obtained in~\S\ref{SSt} (the pointwise bounds of Theorems~\ref{ThmStCo} and the asymptotic profiles of Theorem~\ref{ThmAS}) may skip~\S\ref{SsGA} and continue straight to~\S\ref{SSt}. We stress, however, that the proof of sharp mapping properties on weighted $L^2$-type Sobolev spaces \emph{even for stationary operators} in~\S\ref{SW} depends in a crucial manner on the geometric and analytic structures introduced in~\S\ref{SsGA}.

\subsection{The stationary model}
\label{SsGS}

Let $n\geq 1$. The spatial manifold is
\[
  X := \ol{\R^n}.
\]
We denote the standard coordinates on $X^\circ=\R^n$ by $x=(x^1,\ldots,x^n)\in\R^n$; we shall also frequently use polar coordinates $x=r\omega$, $r^2=\sum_{j=1}^n (x^j)^2$, $\omega\in\Sph^{n-1}$, and inverse polar coordinates $\rho=r^{-1}$. On $X$, we have the scattering cotangent bundle $\Tsc^* X\to X$ as well as the bundle
\[
  \wt\Tsc{}^*X := \ul\R\oplus\Tsc^*X,
\]
where we write $\ul\R:=X\times\R\to X$ for the trivial bundle. On the product
\begin{equation}
\label{EqGM0}
  M_0 := \ol\R \times X,
\end{equation}
we write $t_*\in\R$ for the coordinate in the first factor of $\R\times X\subset M_0$. Moreover, $M_0$ is equipped with a projection map $\pi_0\colon M_0\to X$. We define the bundle
\[
  \wt T^*M_0 := \pi_0^*\bigl(\wt\Tsc{}^*X\bigr).
\]
Over the interior $M_0^\circ=\R\times\R^n$, we identify this bundle with $T^*M_0^\circ$ via the map $\R\oplus T^*_x\R^n\ni(\sigma,\xi)\mapsto\sigma\,\dd t_*+\xi\in T^*\R^{n+1}$.

We furthermore fix a stationary vector bundle
\[
  E\to M_0,
\]
by which we mean that $E=\pi_0^*(E_X)$ for a vector bundle $E_X\to X$. The $(\R,+)$-action on $M_0$ given by translations in the first factor induces a translation action on any stationary bundle. Stationary sections of $E$ are then exactly those which are invariant under this action. These notions apply in particular to tensor bundles built from $\wt T^*M_0$ and its dual $\wt T M_0$. Lastly, we say that a differential operator $P\in\Diff^m(M_0^\circ;E)$ is \emph{stationary} if it commutes with time translations; equivalently, we have
\[
  P = \sum_{j=0}^m P_j\pa_{t_*}^j,\qquad P_j\in\Diff^{m-j}(X^\circ;E_X),
\]
where we identify a `spatial operator' $P_j$ with an operator $\tilde P_j\in\Diff^{m-j}(M_0^\circ;E)$ via $(\tilde P_j u)(t_*,x)=(P_j u(t_*,-))(x)$. For such a stationary operator $P$, we define its \emph{spectral family} as
\[
  \hat P(\sigma) := \sum_{j=0}^m (-i\sigma)^j P_j,\qquad \sigma\in\C.
\]
This is the formal conjugation of $P$ by the Fourier transform~\eqref{EqIFT}. We shall also fix a stationary positive definite fiber inner product on $E$, with respect to which adjoints of operators acting on sections of $E$ are defined. Finally, we denote by $E|_{\pa X}\to\pa X$ the restriction of $E_X$ to $\pa X$.

\begin{rmk}[Fibered cusp perspective]
\label{RmkGphi}
  With respect to the lift under $\pi_0$ of a boundary defining function of $X$, and with respect to the fibration $\ol\R\times\pa X\to\pa X$ induced by $\pi_0$, the bundle $\wt T^*M_0$ is isomorphic to the \emph{fibered cusp} (or $\Phi$-) \emph{cotangent bundle} \cite{MazzeoMelroseFibred} over $\R\times\pa X$, and to the $\Phi$-cusp cotangent bundle globally. However, since the translation action on $M_0$ plays a central role for us, we stick to the above more specific setting.
\end{rmk}

\subsubsection{Geometry}
\label{SssGSG}

For the purposes of the present paper, we use the following variant of \cite[Definition~2.3]{HintzPrice}:

\begin{definition}[Stationary and asymptotically flat metrics]
\label{DefGSG}
  We call a Lorentzian metric $g$ on $M_0^\circ=\R_{t_*}\times\R^n$ \emph{stationary and asymptotically flat} (on $M_0$) if it is stationary (i.e.\ $\pa_{t_*}$ is a Killing vector field), and if moreover
  \begin{enumerate}
  \item\label{ItGSGTime} $\dd t_*$ is everywhere past timelike, i.e.\ $g^{0 0}<0$ in the notation below;
  \item\label{ItGSGTime2} $\pa_{t_*}$ is (future) timelike;
  \item\label{ItGSGMetric} the dual metric of $g$ takes the form
    \begin{equation}
    \label{EqGSGMetric}
      g^{-1} = g^{0 0}\pa_{t_*}^2 + 2\pa_{t_*}\otimes_s g^{0 X} + g^{X X},
    \end{equation}
    where for some $\delta\in(0,1]$, the coefficients are of the form
    \begin{equation}
    \label{EqGSGMetricCoeff}
    \begin{split}
      g^{0 0} &\in \cA^{1+\delta}(X), \\
      g^{0 X} &\equiv -\pa_r \bmod \cA^{1+\delta}(X;\Tsc X), \\
      g^{X X} &\equiv \pa_r^2 + r^{-2}\slg^{-1} \bmod \cA^\delta(X;S^2\,\Tsc X),
    \end{split}
    \end{equation}
    with $\slg\in\CI(\pa X;S^2 T\pa X)$ a Riemannian metric.
  \end{enumerate}
  We say that $g$ is \emph{nontrapping} if, moreover, the following condition is satisfied:
  \begin{enumerate}
  \setcounter{enumi}{3}
  \item\label{ItGSGNontrap} let $\gamma\colon I\subseteq\R\to M_0^\circ$ denote a maximally extended null-geodesic.\footnote{In particular, $\dd t_*(\gamma'(0))\neq 0$ since $\ker\dd t_*$ is spacelike by assumption~\eqref{ItGSGTime}.} Then, as $s\searrow\inf I$ or $s\nearrow\sup I$, we have $\rho(\gamma(s))\searrow 0$, where $\rho=r^{-1}=|x|^{-1}$.
  \end{enumerate}
\end{definition}

In applications, the metric coefficients are often smooth in $(r^{-1},\omega)$, i.e.\ we have $g^{0 0}\in\rho^2\CI(X)$, $g^{0 X}+\pa_r\in\rho^2\CI(X;\Tsc X)$, and $g^{X X}-(\pa_r^2+r^{-2}\slg^{-1})\in\rho\CI(X;S^2\,\Tsc X)$, but such strong regularity is not needed in our analysis. We also note that assumption~\eqref{ItGSGTime2} is not present in \cite{HintzPrice}; in the present paper, it is used to obtain high frequency resolvent estimates for associated wave operators, which is assumed as part of a spectral admissibility condition in \cite[Definition~2.9]{HintzPrice} (which differs from the spectral admissibility condition introduced below in Definition~\ref{DefGSOSpec}).

\begin{rmk}[Timelike nature of $\dd t_*$]
\label{RmkGSGTimelike}
  Given a metric $g^{-1}$ of the form~\eqref{EqGSGMetric}--\eqref{EqGSGMetricCoeff}, condition~\eqref{ItGSGTime} is automatic in a neighborhood of $\pa X$ via replacing $t_*$ by $\tilde t_*:=t_*-C r^{-\delta}$ for sufficiently large $C$; indeed,
  \begin{align*}
    g^{-1}(\dd\tilde t_*,\dd\tilde t_*) &= g^{0 0} + C^2\delta^2 r^{-2-2\delta} g^{X X}(\dd r,\dd r) + 2 C\delta r^{-1-\delta} g^{0 X}(\dd r) \\
      &= -2 C\delta r^{-1-\delta} + g^{0 0} + \cO(r^{-2-2\delta})
  \end{align*}
  is then negative.
\end{rmk}

A more familiar form of the metric $g$ arises via a change of the time coordinate:
\begin{lemma}[Other time coordinates]
\label{LemmaGSGTime}
  Let $g$ be a stationary metric of the form~\eqref{EqGSGMetric}--\eqref{EqGSGMetricCoeff}. Let $t=t_*+F$, where $F\in\CI(X^\circ)$ satisfies $F-r\in\cA^{-1+\delta}(X)$.\footnote{An example is $F(x)=\la r\ra$ where $r=|x|$, or simply $F(x)=r$ for $r>1$, extended smoothly to $r\leq 1$.} Then
  \begin{equation}
  \label{EqGSGTimeStd}
  \begin{split}
    g^{-1} &\equiv -\pa_t^2 + \pa_r^2 + r^{-2}\slg^{-1} \\
      &\qquad \bmod \cA^\delta(X)\pa_t^2 + \pa_t\otimes_s \cA^\delta(X;\Tsc X) + \cA^\delta(X;S^2\,\Tsc X).
  \end{split}
  \end{equation}
\end{lemma}

We caution however that, conversely, the membership~\eqref{EqGSGTimeStd} does not imply \eqref{EqGSGMetricCoeff}, as it does not capture the additional orders of decay of $g^{0 0}$ and of the remainder term of $g^{0 X}$ in~\eqref{EqGSGMetricCoeff}.

\begin{proof}[Proof of Lemma~\usref{LemmaGSGTime}]
  Write $F=r+\tilde F(r,\omega)$ near $\pa X$. Changing coordinates in the vector fields $\pa_{t_*}$, $\pa_r$, $\pa_\omega$ to $(t,r,\omega)$ gives $\pa_t$, $\pa_r+(1+\pa_r\tilde F)\pa_t$, $\pa_\omega+(\pa_\omega\tilde F)\pa_t$, which implies~\eqref{EqGSGTimeStd}.
\end{proof}

\begin{example}[Minkowski type metrics]
\label{ExGSGMink}
  If $\slg$ is the standard metric on $\Sph^{n-1}$, then the Minkowski dual metric $g_M^{-1}=-\pa_t^2+\pa_r^2+r^{-2}\slg^{-1}$,  expressed in terms of $t_*=t-r$ in $r>0$, is equal to
  \[
    g_M^{-1}=-2\pa_{t_*}\otimes_s\pa_r+\pa_r^2+r^{-2}\slg^{-1}.
  \]
  Any dual metric $g^{-1}$ of the form~\eqref{EqGSGMetric}--\eqref{EqGSGMetricCoeff} thus satisfies
  \begin{equation}
  \label{EqGSGMetMink}
  \begin{split}
    g^{-1} - g_M^{-1} &\in \cA^{1+\delta}(X)\pa_{t_*}^2 + 2\pa_{t_*}\otimes_s \cA^{1+\delta}(X;\Tsc X) + \cA^\delta(X;S^2\,\Tsc X) \\
      &\subset \cA^\delta(X;S^2\,\wt{\Tsc}X).
  \end{split}
  \end{equation}
  More generally, we call $g_M^{-1}$ a \emph{Minkowski type metric} when $\slg$ is a general Riemannian metric on $\Sph^{n-1}$.
\end{example}

\begin{example}[Asymptotically Kerr metrics]
\label{ExGSGKerr}
  Black hole metrics violate Definition~\ref{DefGSG} in three ways: they are not defined on the entire spatial manifold $X^\circ=\R^n$; they violate the nontrapping condition; and they typically feature ergoregions, where $\pa_{t_*}$ fails to be timelike. We discuss the required modifications in the companion paper \cite{HintzVasyNonstat2}. Kerr metrics do satisfy assumptions~\eqref{ItGSGTime}--\eqref{ItGSGMetric} of Definition~\ref{DefGSG}, with $\delta=1$, for large $r$; this can be read off from the principal symbol of \cite[Lemma~4.2]{HaefnerHintzVasyKerr}. We only explicitly describe the case of a Schwarzschild black hole with mass $\bhm$, where this follows from the expression $g^{-1}=-2\pa_{t_*}\otimes_s\pa_r+(1-\frac{2\bhm}{r})\pa_r^2+r^2\slg^{-1}$ for the dual metric in outgoing Eddington--Finkelstein coordinates.
\end{example}

\begin{lemma}[Volume density]
\label{LemmaGSGVol}
  Let $g$ be as in Definition~\usref{DefGSG}. Then the volume density $|\dd g|\in\CI(M_0^\circ;\Omega M_0^\circ)$ takes the form $|\dd g|=|\dd t_*||\dd g_X|$, where
  \[
    |\dd g_X| \in (\CI+\cA^\delta)(X;\Omegasc X),\qquad |\dd g_X|\equiv r^{n-1}|\dd r\,\dd\slg|\bmod\cA^\delta(X;\Omegasc X).
  \]
\end{lemma}
\begin{proof}
  This follows directly from the definition via a simple calculation in the coordinates $(t_*,r,\omega)$ on $M_0^\circ$, or more easily still using Lemma~\ref{LemmaGSGTime}.
\end{proof}

\emph{We shall use $|\dd g|$, resp.\ $|\dd g_X|$ to define adjoints of operators on $M_0^\circ$, resp.\ $X^\circ$.}

\subsubsection{Wave type operators}
\label{SssGSO}

The class of stationary models we consider in this paper is the following:

\begin{definition}[Stationary wave type operators]
\label{DefGSO}
  Let $g$ be a stationary and asymptotically flat Lorentzian metric on $M_0$, and let $E\to M_0$ be a stationary vector bundle. Then a stationary operator $P\in\Diff^2(M_0^\circ;E)$ is called a \emph{stationary wave type operator} (with respect to $g$) if:
  \begin{enumerate}
  \item\label{ItGSOSymb} the principal symbol of $P$ is scalar and equal to the dual metric function $\zeta\mapsto g^{-1}(\zeta,\zeta)$ of $g$;
  \item\label{ItGSOStruct} upon writing $P$ near $\rho=r^{-1}=0$ in the form
    \begin{equation}
    \label{EqGSOStruct}
      P = -2\pa_{t_*}\rho\Bigl(\rho\pa_\rho-\frac{n-1}{2}-S\Bigr) + \hat P(0) + Q\pa_{t_*} - g^{0 0}\pa_{t_*}^2,
    \end{equation}
    we have, for some $\delta\in(0,1]$, the membership $g^{0 0}\in\cA^{1+\delta}(X)$ from~\eqref{EqGSGMetricCoeff} and
    \begin{align*}
      \hat P(0) &\in \rho^2\Diffb^2(X;E_X) + \cA^{2+\delta}\Diffb^2(X;E_X), \\
      S &\in \CI(X;\End(E_X)) + \cA^\delta(X;\End(E_X)), \\
      Q &\in \cA^{2+\delta}\Diffb^1(X;E_X).
    \end{align*}
  \end{enumerate}
\end{definition}

\begin{rmk}[Comments on Definition~\ref{DefGSO}]
\label{RmkGSOStruct}
  Given condition~\eqref{ItGSOSymb}, the operator $P$ is necessarily of the form~\eqref{EqGSOStruct} near $\rho=0$ for some $S\in\CI(X^\circ;\End(E_X))$, and also the principal parts of $\hat P(0)$ and $Q$ are determined. The new pieces of information in part~\eqref{ItGSOStruct} are thus the smoothness (up to a decaying conormal error) of $S$ down to $\pa X$ (the term $\frac{n-1}{2}$ in~\eqref{EqGSOStruct} is merely a convenient normalization of $S$) as well as the membership of the subprincipal terms of $\hat P(0)$ and $Q$ in $\rho^2\Diffb^1(X;E_X)+\cA^{2+\delta}\Diffb^1(X;E_X)$ and $\cA^{2+\delta}(X;\End(E_X))$, respectively. One may allow for $Q$ to contain also terms in $\cA^{1+\delta}\Diffb^1(X;E_X)$; however, this does not enlarge the class of operators considered, since such a term can equivalently be regarded as a contribution $\half\rho^{-1}\tilde Q\in\cA^\delta(X;\End(E_X))$ to $S$.
\end{rmk}

\begin{example}[Wave operator]
\label{ExGSO}
  The wave operator $\Box_g$ of a stationary and asymptotically flat Lorentzian metric, such as the Minkowski metric, is an example of a stationary wave type operator, with $E_X$ (and thus $E$) the trivial bundle and with $S|_{\pa X}=0$. See Proposition~\ref{PropES}.
\end{example}

\begin{example}[Coupling with potentials or first order terms]
\label{ExGSOPot}
  If $g$ is stationary and asymptotically flat, then $\Box_g+V_0$ is a stationary wave type operator for all potentials $V_0\in\la x\ra^{-2}\CI(X;E_X)+\cA^{2+\delta}(X;E_X)$ on $X=\ol{\R^n}$. These are precisely those potentials which (in local trivializations of $E_X$) have a leading order term at infinity with inverse square decay $|x|^{-2}$ plus a remainder which is very short range, namely of size $|x|^{-2-\delta}$ together with all b-derivatives (derivatives along $\la x\ra\pa_x$). Indeed, such $V_0$ can be absorbed into $\hat P(0)$. Allowed first order terms are of the schematic form $a(x)\pa_{t_*}$ or $a(x)\pa_x$ where $a\in\la x\ra^{-1}\CI+\cA^{1+\delta}$ (i.e.\ inverse linear decay plus lower order terms): the first type of term can be absorbed into $S$, and the second type of term into $\hat P(0)$.
\end{example}

The spectral family of $P$ takes the form
\begin{equation}
\label{EqGSOSpecFam}
  \hat P(\sigma) = 2 i\sigma\rho\Bigl(\rho\pa_\rho-\frac{n-1}{2}-S\Bigr) + \hat P(0) - i\sigma Q + g^{0 0}\sigma^2.
\end{equation}

Therefore, $\hat P(0)\in\rho^2\Diffb^2(X;E_X)$. For general $\sigma\in\C$, the operator~\eqref{EqGSOSpecFam} lies in\footnote{The coefficients are smooth plus decaying conormal; we shall typically not write this out explicitly.} $\Diffsc^2(X;E_X)$ since $\rho^j\Diffb^j\subset\Diffsc^j$ for $j\in\N_0$. (We shall see that for $\sigma\neq 0$, the operator $\hat P(\sigma)$ is, in a suitable sense, nondegenerate as a scattering differential operator, whereas for $\sigma=0$ it is important to use the sharper b-nature of the operator. See already Lemma~\ref{LemmaStEst0} and Proposition~\ref{PropStEstNz}.)

\begin{lemma}[Low energy behavior]
\label{LemmaGSOLo}
  Let $P$ be a stationary wave type operator, and let $\theta\in[0,\pi]$. Then the low energy spectral family $[0,1)e^{i\theta}\ni\sigma\mapsto\hat P(\sigma)$ defines an element of $\Diffscbt^{2,0,-2,0}(X;E_X)$.
\end{lemma}
\begin{proof}
  Recall that $\sigma$-independent elements of $\Diffb^1(X;E_X)$ (such as $\rho\pa_\rho$) are elements of $\Diffscbt^{1,1,0,0}(X;E_X)$; and moreover $\sigma\in\Diffscbt^{0,0,-1,-1}(X;E_X)$ and $\rho\in\Diffscbt^{0,-1,-1,0}(X;E_X)$. This implies the claim.
\end{proof}

\begin{definition}[Transition face normal operator]
\label{DefGSOtf}
  For a stationary wave type operator $P$, we denote by
  \[
    N_\tface^\theta(P)\in\Diff_{\scop,\bop}^{2,(0,2)}(\tface;\upbeta_\tface^*E|_{\pa X}),\qquad \theta\in[0,\pi],
  \]
  the transition face normal operator of $[0,1)e^{i\theta}\ni\sigma\mapsto\sigma^{-2}\hat P(\sigma)$; here $\upbeta_\tface\colon\tface\to\pa X$ is the projection (or blow-down) map. We also write
  \[
    N_\tface^+(P) := N_\tface^0(P),\qquad
    N_\tface^-(P) := N_\tface^\pi(P)
  \]
  for the transition face normal operators of $\pm[0,1)\ni\sigma\mapsto\sigma^{-2}\hat P(\sigma)$.
\end{definition}

Concretely, fixing a collar neighborhood $[0,\eps)_\rho\times\pa X$ and a bundle isomorphism of $E_X$ with the pullback of $E|_{\pa X}$ along the projection $[0,\eps)\times\pa X\to\pa X$, and for $u\in\CIc((0,\infty)_{\hat\rho}\times\Sph^{n-1};E|_{\pa X})$, we set $u_\sigma(\rho,\omega)=u(\frac{\rho}{|\sigma|},\omega)$ and $N_\tface^\theta(P)u(\hat\rho,\omega) := \lim_{\sigma_0\searrow 0}\bigl( (\sigma^{-2}\hat P(\sigma))u_\sigma \bigr) \big|_{(|\sigma|\hat\rho,\omega)}$ where $\sigma=e^{i\theta}\sigma_0$. Formally, we obtain $N_\tface^\theta(P)$ from $P$ by replacing $\pa_{t_*}$ and $\rho$ by $-i e^{i\theta}$ and $\hat\rho$, dropping the terms in~\eqref{EqGSOStruct} involving $Q$ and $g^{0 0}$, and restricting $S$ and the coefficients of $\rho^{-2}\hat P(0)$ (as a b-operator) to $\pa X$. That is, writing
\begin{equation}
\label{EqGSOP0}
  \hat P(0) = \rho^2 P_{(0)}(\rho,\omega,\rho D_\rho,D_\omega),
\end{equation}
we have
\begin{equation}
\label{EqGSONtf}
\begin{split}
  N_\tface^\theta(P) &= 2 i e^{i\theta}\hat\rho\Bigl(\hat\rho\pa_{\hat\rho}-\frac{n-1}{2}-S|_{\pa X}\Bigr) + \hat\rho^2 P_{(0)}(0,\omega,\hat\rho D_{\hat\rho},D_\omega), \\
  N_\tface^\pm(P) &= \pm 2 i\hat\rho\Bigl(\hat\rho\pa_{\hat\rho}-\frac{n-1}{2}-S|_{\pa X}\Bigr) + \hat\rho^2 P_{(0)}(0,\omega,\hat\rho D_{\hat\rho},D_\omega).
\end{split}
\end{equation}
In view of Lemma~\ref{LemmaGSOLo}, or simply by direct inspection, we have $N_\tface^\theta(P)\in\Diff_{\scop,\bop}^{2,(0,2)}(\tface;E|_{\pa X})$, where $\tface=[0,\infty]_{\hat\rho}\times\Sph^{n-1}$.

\begin{definition}[Spectral admissibility]
\label{DefGSOSpec}
  Let $P$ be a stationary wave type operator associated with the stationary and asymptotically flat metric $g$ on $M_0$, and denote by $\hat P(\sigma)$ its spectral family. Let $\beta^-<\beta^+$. We then say that $P$ is \emph{spectrally admissible (with indicial gap $(\beta^-,\beta^+)$)} if the following assumptions are satisfied.
  \begin{enumerate}
  \item\label{ItGSOSpec}{\rm (Mode stability away from zero energy.)} For all $\sigma\in\C$ with $\Im\sigma\geq 0$, $\sigma\neq 0$, every $u\in\cA(X;E_X)$ with $\hat P(\sigma)u=0$ must vanish identically.
  \item\label{ItGSOSpec0}{\rm (No zero energy resonance.)} For all $\beta\in(\beta^-,\beta^+)$, the kernel of $\hat P(0)$ on $\cA^\beta(X;E_X)$ and the kernel of $\hat P(0)^*$ on $\cA^{n-2-\beta}(X;E_X)$ are trivial.
  \item\label{ItGSOSpectf}{\rm (Invertibility of the transition face normal operators.)} For all $\theta\in[0,\pi]$, the kernels of $N_\tface^\theta(P)$ on $\cA^{(\alpha,-\beta)}(\tface;\upbeta_\tface^*E|_{\pa X})$ are trivial for all $\alpha\in\R$, $\beta\in(\beta^-,\beta^+)$, and the kernels of $N_\tface^\pm(P)^*$ on\footnote{Observe that $\frac{\hat\rho}{1+\hat\rho}$ is a defining function of $\hat\rho^{-1}(0)\subset\tface$.} $e^{\mp 2 i(1+\hat\rho)/\hat\rho}\cA^{(\alpha,-n+2+\beta)}(\tface;\upbeta_\tface^*E|_{\pa X})$ are trivial for all such $\alpha,\beta$ as well. (The weights refer, in this order, to decay at $\hat\rho^{-1}(0)$ and $\hat\rho^{-1}(\infty)\subset\tface$..)
  \item\label{ItGSOOptimal}{\rm (Optimality of weights.)} The interval $(\beta^-,\beta^+)$ is the largest open interval of weights $\beta$ for which conditions~\eqref{ItGSOSpec0} and \eqref{ItGSOSpectf} hold.
  \end{enumerate}
\end{definition}

This definition is related to \cite[Definition~2.9]{HintzPrice}, though in the reference only scalar operators are studied (in which case the transition face normal operators are automatically invertible, see \cite[Theorem~2.22]{HintzPrice} for the case $n=3$). Moreover, as we recall in Proposition~\ref{PropStEstHi}, the nontrapping assumption on $g$ in Definition~\ref{DefGSG}\eqref{ItGSGNontrap} implies high energy estimates, which we thus do not need to assume separately here.

\begin{rmk}[Outgoing solutions]
\label{RmkGSOSpecOut}
  Condition~\eqref{ItGSOSpectf} demands the nonexistence of purely outgoing solutions of $N_\tface^\theta(P)$ (which in the present, conjugated, perspective do not have any oscillatory behavior at $\hat\rho=0$) and of purely incoming solutions of $N_\tface^\pm(P)^*$. (The triviality of the kernel of $N_\tface^\theta(P)^*$, $\theta\in(0,\pi)$, on appropriate function spaces does not need to verified separately, as it follows from Fredholm index arguments; see Lemma~\ref{LemmaStEsttf}.)
\end{rmk}

\begin{rmk}[Zero energy weights: I]
\label{RmkGSOBetapm}
  The values of $\beta^-$ and $\beta^+$ can be read off from the boundary spectrum of $\hat P(0)$; see Lemma~\ref{LemmaStEst0} and Remark~\ref{RmkStEst0Weights}.
\end{rmk}

\subsection{Asymptotically flat spacetimes}
\label{SsGA}

Starting with the product manifold $M_0=\ol\R\times X$ from~\eqref{EqGM0}, we now define
\[
  M_1 := [M_0; \pa\ol\R\times\pa X] = [M_0; \{+\infty\}\times\pa X, \{-\infty\}\times\pa X],\qquad
  \wt T^*M_1 := \upbeta_1^*\bigl(\wt T^*M_0\bigr),
\]
where $\upbeta_1\colon M_1\to M_0$ is the total blow-down map. The following result explains the sense in which $M_1$ fits into the framework of 3b-analysis near the lift of $\{+\infty\}\times X$ (see also Definition~\ref{DefGAGVF} below):

\begin{lemma}[Alternative description of $M_1$]
\label{LemmaGA}
  Let $t=t_* + F$ where $F\in\CI(X^\circ)$ (regarded as a stationary function on $\R_{t_*}\times X^\circ$) and $F(x)-\la x\ra\in\CI(X)$. Then the map
  \begin{equation}
  \label{EqGAMap}
    \R_{t_*}\times\R^n\ni (t_*,x) \mapsto (t,x) = (t_*+F(x),x) \in \R_t\times\R^n_x=\R^{n+1}
  \end{equation}
  extends to a diffeomorphism of $M_1$ and the manifold $M'_1$ obtained from $\ol{\R^{n+1}}=\ol{\R_t\times\R_x^n}$ by blowing up the `north' and `south' poles $N=\{(+\infty,0)\}$ and $S=\{(-\infty,0)\}$ as well as the light cone at future infinity $Y_+$, defined as the intersection of the closure of $t=r$ with $\pa\ol{\R^{n+1}}$.

  Moreover, the map induced by~\eqref{EqGAMap} on the cotangent bundles extends to a bundle isomorphism $\wt T^*M_1\cong(\upbeta_1')^*(\Tsc^*\ol{\R^{n+1}})$, where $\upbeta_1'\colon M'_1\to\ol{\R^{n+1}}$ is the blow-down map.
\end{lemma}

The asymmetric definitions of $M_1$ and $M'_1$ arise because we choose to work with the retarded time $t_*=t-F\approx t-r$ in one case and the standard time $t$ in the other. Note also that in $M'_1$ the light cone at past null infinity is not resolved, since we are interested only in the region $t\geq 0$ here.\footnote{For completeness, we note that if one equips $M_1$ with a stationary, asymptotically flat (with mass $\bhm\in\R$) metric in the sense of \cite[Definition~2.3]{HintzPrice}, then past null infinity, given as the intersection of the closure of $t=-r$ with $M_1$, is correctly resolved only if one changes the smooth structure near the light cone at past infinity, due to the possible presence of long range mass terms in the metric; cf.\ \cite[\S2.1]{HintzVasyMink4}.}

\begin{proof}[Proof of Lemma~\usref{LemmaGA}]
  Let $\tilde F:=F(x)-\la x\ra\in\CI(X)$. It is easy to check that $(t_*,x)\mapsto(t_*+\tilde F(x),x)$ induces a diffeomorphism of $M_1$. Therefore, we may assume that $F(x)=\la x\ra$.

  Since $N$ is defined by $\frac{1}{t}=0$, $\frac{x}{t}=0$, a collar neighborhood of the lift of $N$ to $M'_1$ is given by $[0,1)_{\rho'_\cT}\times\ol{\R^n_x}$ where $\rho'_\cT=((\frac{1}{t})^2+|\frac{x}{t}|^2)^{1/2}=\frac{\la x\ra}{t}$. On the other hand, since $t_*^{-1}$ and $\la x\ra^{-1}$ are local defining functions of $\{+\infty\}\times X$ and $\ol\R\times\pa X\subset M_0$, a collar neighborhood of the lift of $\{+\infty\}\times X$ to $M_1$ is given by $\ol{\R^n_x}\times[0,\infty)_{\rho_\cT}$, where $\rho_\cT=\frac{t_*^{-1}}{\la x\ra^{-1}}=\frac{\la x\ra}{t_*}$. Since $\rho_\cT'=\rho_\cT(1+\rho_\cT)^{-1}$, and conversely $\rho_\cT=\rho'_\cT(1-\rho'_\cT)^{-1}$, this shows that the map~\eqref{EqGAMap} restricts to a diffeomorphism of these two collar neighborhoods.

  Similar arguments show that a collar neighborhood of the lift of $S$ to $M'_1$, resp.\ of the lift of $\{-\infty\}\times X$ to $M_1$, are naturally diffeomorphic.

  Next, since $Y_+$ is defined by $\frac{t_*}{r}=0$, $r^{-1}=0$, a collar neighborhood of the lift of $Y_+$ to $M'_1$ is given by $\ol{\R_{t_*}}\times[0,2)_{\rho_{\!\scri}}\times\pa X$ where $\rho_{\!\scri}=((\frac{t_*}{r})^2+(\frac{1}{r})^2)^{1/2}=\frac{\la t_*\ra}{r}$ . With $\ol{\R_{t_*}}\times[0,2)_{r^{-1}}\times\pa X$ being a collar neighborhood of $\ol\R\times\pa X\subset M_0$ and $\la t_*\ra^{-1}$ being a defining function of $\{\pm\infty\}\times X\subset M_0$, the chart $\ol{\R_{t_*}}\times[0,2)_{\rho_{\!\scri}}\times\pa X$ is also collar neighborhood of the lift of $\ol\R\times\pa X$ to $M_1$; note indeed that $\frac{r^{-1}}{\la t_*\ra^{-1}}=\rho_{\!\scri}$.

  The bundle isomorphism of $\wt T^*M_1$ with the pullback of $\Tsc^*\ol{\R^{n+1}}$ follows directly from the definitions: both bundles have smooth frames given by the coordinate differentials $\dd t_*$, $\dd x^j$ ($j=1,\ldots,n$).
\end{proof}

\begin{definition}[Spacetime manifold]
\label{DefGAMfd}
  We define the manifold with corners $M$ as the square root blow-up of $M_1=[\ol\R\times X;\pa\ol\R\times\pa X]$ at the lift of $\ol\R\times\pa X$; that is, $M=M_1$ as sets, and $\CI(M)$ is the smallest algebra containing $\CI(M_1)$ and the square root of a defining function of the lift of $\ol\R\times\pa X$.\footnote{Explicit charts are given below.} We label the boundary hypersurfaces of $M$ and $M_1$ as follows.
  \begin{itemize}
  \item{\rm (Future null infinity).} $\scri^+$ is the lift of $\ol\R\times\pa X$.
  \item{\rm (Punctured future timelike infinity.)} $\iota^+$ is the lift of $\{+\infty\}\times\pa X$.
  \item{\rm (Future translation face.)} $\cT^+$ is the lift of $\{+\infty\}\times X$.
  \item{\rm (Spacelike infinity.)} $I^0$ is the lift of $\{-\infty\}\times\pa X$.\footnote{The terminology is correct only in $t\geq 0$, which is the only region we care about here.}
  \end{itemize}
  We write $\upbeta\colon M\to M_0$ for the total blow-down map, and write
  \[
    \wt T^*M := \upbeta^*\bigl(\wt T^*M_0\bigr).
  \]
  Lastly, we denote by $\rho_0$, $x_{\!\scri}$, $\rho_+$, and $\rho_\cT\in\CI(M)$ defining functions of $I^0$, $\scri^+$, $\iota^+$, and $\cT^+$, respectively; and we set $\rho_{\!\scri}=x_{\!\scri}^2$ (which is a defining function of the lift of $\ol\R\times\pa X$ to $M_1$). Their definitions may change throughout the paper; we moreover use this notation also for defining functions in local coordinates.
\end{definition}

Concretely, a collar neighborhood of $\cT^+$ (where $r\ll t_*$) is given by $[0,\infty)_{\rho_\cT}\times\ol{\R^n_x}$ where $\rho_\cT:=\frac{\la x\ra}{t_*}\in[0,\infty)$. Near $\pa\cT^+=\cT^+\cap\iota^+$ (where $1\lesssim r\ll t_*$), we have a local chart $[0,\infty)_{\rho_\cT}\times[0,\infty)_{\rho_+}\times\Sph^{n-1}$ where now $\rho_\cT=\frac{|x|}{t_*}$ and $\rho_+=\frac{1}{|x|}$. A chart near $\iota^+\setminus\cT^+$ (where $1\lesssim t_*\lesssim r$ and $r\gg 1$) is given by $[0,\infty)_{\rho_+}\times[0,\infty)_{x_{\!\scri}}\times\Sph^{n-1}$ where now $\rho_+=\frac{1}{t_*}$ and $x_{\!\scri}=(\frac{t_*}{r})^{\frac12}$. A chart near $\scri^+\cap I^0$ (where $1\lesssim -t_*<r$ and $r\gg 1$) is $[0,\infty)_{\rho_0}\times[0,\infty)_{x_{\!\scri}}\times\Sph^{n-1}$ where $\rho_0=\frac{1}{|t_*|}$ and $x_{\!\scri}=(\frac{|t_*|}{r})^{\frac12}$. (See also Lemma~\ref{LemmaGAGeb} below.) Global choices in $t_*+\la x\ra>-\la x\ra$ for the defining functions can be obtained by suitably gluing these local definitions; explicit examples are smoothed out versions of
\[
  \rho_0 = \frac{1}{1+(t_*)_-},\qquad
  x_{\!\scri} = \sqrt{\frac{\la t_*\ra}{t_*+2 r}},\qquad
  \rho_+ = \frac{t_*+2 r}{r(1+(t_*)_+)},\qquad
  \rho_\cT = \frac{r}{t_*+2 r}.
\]
In $t_*\geq 1$, one may take $\rho_\cT$, $\rho_+$, and $x_{\!\scri}:=\rho_{\!\scri}^{\frac12}$ in terms of~\eqref{EqIRhos} (with $r$ replaced by $\la x\ra$ to ensure smoothness at $r=0$).

\begin{figure}[!ht]
\centering
\includegraphics{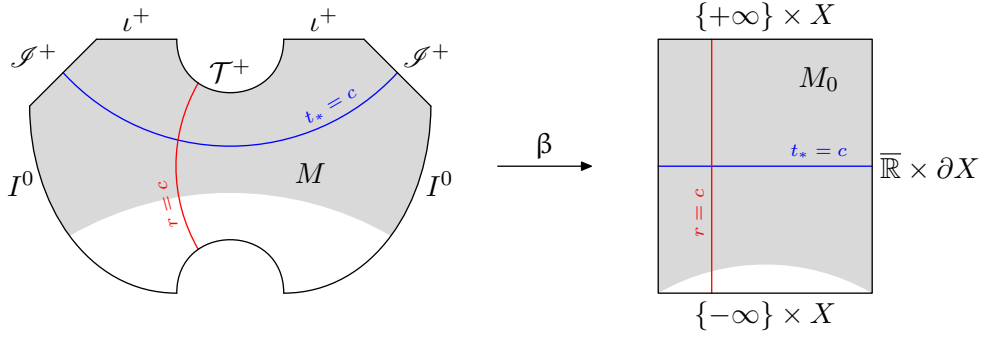}
\caption{\textit{On the left:} the spacetime manifold $M$ from Definition~\ref{DefGAMfd}, and its boundary hypersurfaces. The shaded region is the set where $t_*+r\geq-\half r-1$. \textit{On the right:} the product manifold $M_0=\ol{\R_{t_*}}\times X$ where $X=\ol{\R^n_x}$. Shown on both sides are level sets of $t_*$ and $r=|x|$.}
\label{FigGAMfd}
\end{figure}

The part of $M$ on which we shall consider wave equations is contained in $\{t_*+r\geq-\half r-1\}$ (or more precisely its closure in $M$). See Figure~\ref{FigGAMfd}.

\subsubsection{Geometry}
\label{SssGAG}

Null infinity $\scri^+\subset M$ is the total space of a fibration $\upbeta|_{\scri^+}\colon\scri^+=\ol\R\times\pa X\to\pa X=\Sph^{n-1}$. Moreover, $\cT^+\subset M$ is the front face of the blow-up of $\ol{\R^{n+1}}$ at the `north pole' $N$ in the notation of Lemma~\ref{LemmaGA}. The Lie algebra of vector fields on $M$ which will be the basis of our analysis is then:

\begin{definition}[b-edge-3b-vector fields]
\label{DefGAGVF}
  We write $\Vbetb(M)$ for the space of smooth vector fields on $M$ which are tangent to (i.e.\ `b at') all boundary hypersurfaces, tangent also to the fibers of (i.e.\ `edge at') $\scri^+$, and which near $\cT^+$ and in the notation of Lemma~\usref{LemmaGA} are moreover of the form $\rho_+^{-1}(\upbeta'_1)^*W$ where $\rho_+\in\CI(M)$ is a defining function of $\iota^+$ and $W\in\Vsc(\ol{\R^{n+1}})$, i.e.\ they are `3b (3-body/b-)vector fields near $\cT^+$'. We write $\Tbetb M\to M$ for the corresponding tangent bundle so that $\Vbetb(M)=\CI(M;\Tbetb M)$.
\end{definition}

Concretely, in the coordinates $t_*$ and $x=r\omega$ on $\R$ and $X^\circ=\R^n$, elements of $\Vbetb(M)$ are smooth (on $M$) linear combinations of the following vector fields:
\begin{itemize}
\item near $\cT^+$: $\la x\ra\pa_{t_*}$, $\la x\ra\pa_x$ (or $r\pa_{t_*}$, $r\pa_r$, $\pa_\omega$ away from the closure of $r=0$);
\item near $\scri^+$ and in $t_*>0$, or more generally in $t_*>T$, $T\in\R$, upon replacing $t_*$ in the following expressions by $t_*-T$: $t_*\pa_{t_*}$, $r\pa_r$, $\sqrt{\frac{t_*}{r}}\pa_\omega$ (and likewise in $t_*<0$ upon replacing $t_*$ by $|t_*|$). In terms of $t=t_*+r$, the first two vector fields can be replaced by $r(\pa_t+\pa_r)$, $r\pa_r+t\pa_t$.
\end{itemize}
On $M\setminus(\scri^+\cup\cT^+)$, the bundle $\Tbetb M$ is isomorphic to $\Tb M$. That is, the spaces of b-edge-3b-vector fields on $M$ and b-vector fields on $M$ with support disjoint from $\scri^+\cup\cT^+$ are equal; furthermore, by Lemma~\ref{LemmaGA}, such vector fields are lifts of smooth b-vector fields on $\ol{\R^{n+1}}$. Since $\la r\ra^{-1}$ is a defining function of $I^0\cup\iota^+\setminus(\scri^+\cup\cT^+)$ in $M\setminus(\scri^+\cup\cT^+)$, we thus conclude that away from $\scri^+\cup\cT^+$ the coordinate vector fields $\pa_{t_*}$, $\pa_x$ are a local frame of $\la r\ra^{-1}\,\Tbetb M$; this gives a bundle isomorphism of (the pullback of) $\Tsc\ol{\R^{n+1}}$ and $\rho_0\rho_+\,\Tbetb M$ over $M\setminus(\scri^+\cup\cT^+)$.

Let now $g_0$ be stationary and asymptotically flat as in Definition~\ref{DefGSG}. Thus, we have $g_0\in\CI(M_0;S^2\wt T^*M_0)$, and upon pulling $g_0$ back to $M$ also $g_0\in\CI(M;S^2\wt T^*M)$ in the notation of Definition~\ref{DefGAMfd}. As already mentioned in~\S\ref{SsIN} and discussed in detail in~\cite{HintzVasyScrieb}, the perturbations of $g_0$ we shall allow for near $\scri^+$ lie in an appropriate class of weighted edge-b-metrics. As a preparation, we first show:

\begin{lemma}[$g_0$ as a weighted edge-b-metric near $\scri^+$]
\label{LemmaGAGeb}
  Let $T\in\R$, and use inverse polar coordinates $\rho=r^{-1}$, $\omega\in\Sph^{n-1}$ on $X$ near $\pa X$.
  \begin{enumerate}
  \item\label{ItGAGeb0}{\rm (Coordinates near $\scri^+\setminus\iota^+$.)} In $t_*<T$, define the smooth coordinates
    \begin{equation}
    \label{EqGAGeb0}
      \rho_0 := \frac{1}{T-t_*},\qquad
      x_{\!\scri} := \sqrt{(T-t_*)\rho}
    \end{equation}
    on $M$. (Thus, $[0,\infty)_{\rho_0}\times[0,\infty)_{x_{\!\scri}}\times\Sph^{n-1}_\omega$ is a local coordinate chart on $M$ near $\scri^+\cap I^0$.) Write $\Teb M$ for the edge-b-bundle in this coordinate chart, corresponding to b-behavior at $I^0$ and edge behavior at $\scri^+$. Then
    \begin{equation}
    \label{EqGAGeb0Metric}
      g_0^{-1} \equiv \frac12 \rho_0^2 x_{\!\scri}^2\Bigl( -x_{\!\scri}\pa_{x_{\!\scri}} \otimes_s (x_{\!\scri}\pa_{x_{\!\scri}}-2\rho_0\pa_{\rho_0}) + 2 x_{\!\scri}^2\slg^{-1} \Bigr) \bmod \cA^{(((2,0),2+\delta),2+2\delta)}(M;S^2\,\Teb M).
    \end{equation}
    (The space on the right consists of symmetric 2-cotensors whose coefficients in the frame $\frac{\dd\rho_0}{\rho_0}$, $\frac{\dd x_{\!\scri}}{x_{\!\scri}}$, and $x_{\!\scri}^{-1}T^*\Sph^{n-1}$ lie in $(\rho_0^2\CI+\cA^{2+\delta})([0,\infty)_{\rho_0};\cA^{2+2\delta}([0,\infty)_{x_{\!\scri}}\times\pa X))$. See~\eqref{EqMphgcon} for the general definition.)
  \item\label{ItGAGebp}{\rm (Coordinates near $\scri^+\setminus I^0$.)} In $t_*>T$, and in the coordinates
    \begin{equation}
    \label{EqGAGebp}
      \rho_+ := \frac{1}{t_*-T},\qquad
      x_{\!\scri} := \sqrt{(t_*-T)\rho}
    \end{equation}
    on $M$, we have (now with edge behavior at $\scri^+$ and b-behavior at $\iota^+$)
    \begin{equation}
    \label{EqGAGebpMetric}
    \begin{split}
      g_0^{-1} &\equiv \frac12 x_{\!\scri}^2\rho_+^2\Bigl( x_{\!\scri}\pa_{x_{\!\scri}} \otimes_s (x_{\!\scri}\pa_{x_{\!\scri}}-2\rho_+\pa_{\rho_+}) + \half x_{\!\scri}^2(x_{\!\scri}\pa_{x_{\!\scri}})^2 + 2 x_{\!\scri}^2\slg^{-1}\Bigr) \\
        &\quad\hspace{14em} \bmod \cA^{(2+2\delta,2+\delta)}(M;S^2\,\Teb M).
    \end{split}
    \end{equation}
    (The orders on the error term here refer to the powers of $x_{\!\scri}$ and $\rho_+$, in this order.)
  \end{enumerate}
\end{lemma}
\begin{proof}
  For part~\eqref{ItGAGeb0}, we compute $r^{-1}=\rho=\rho_0 x_{\!\scri}^2$, and
  \begin{equation}
  \label{EqGAGeb0VF}
  \begin{alignedat}{2}
    \pa_{t_*} &= \rho_0\bigl(\rho_0\pa_{\rho_0}-\half x_{\!\scri}\pa_{x_{\!\scri}}\bigr) &&\in \rho_0\Veb(M), \\
    \pa_r &= -\half\rho_0 x_{\!\scri}^3\pa_{x_{\!\scri}} &&\in \rho_0 x_{\!\scri}^2\Veb(M), \\
    r^{-1}\pa_\omega &= \rho_0 x_{\!\scri}^2\pa_\omega &&\in \rho_0 x_{\!\scri}\Veb(M).
  \end{alignedat}
  \end{equation}
  where we schematically write $\pa_\omega$ for a vector field on $\Sph^{n-1}$. Thus, in this coordinate patch we have $\CI(X;\Tsc X)\subset\rho_0 x_{\!\scri}\CI(M;\Teb M)$. Plugging these expressions into~\eqref{EqGSGMetric}--\eqref{EqGSGMetricCoeff}, and noting that $\cA^\alpha(X)=\rho^\alpha\cA^0(X)$ regarded as a space of stationary functions on $M$ is a subset of $\rho_0^\alpha x_{\!\scri}^{2\alpha}\cA^{(0,0)}(M)=\cA^{(\alpha,2\alpha)}(M)$ in the coordinate patch, gives~\eqref{EqGAGeb0Metric}.

  For part~\eqref{ItGAGebp}, we similarly compute $r^{-1}=\rho=x_{\!\scri}^2 \rho_+$ and
  \begin{equation}
  \label{EqGAGebpVF}
  \begin{alignedat}{2}
    \pa_{t_*} &= -\rho_+\bigl(\rho_+\pa_{\rho_+}-\half x_{\!\scri}\pa_{x_{\!\scri}}\bigr) &&\in \rho_+\Veb(M), \\
    \pa_r &= -\half\rho_+ x_{\!\scri}^3\pa_{x_{\!\scri}} &&\in x_{\!\scri}^2\rho_+\Veb(M), \\
    r^{-1}\pa_\omega &= x_{\!\scri}^2\rho_+\pa_\omega &&\in x_{\!\scri}\rho_+\Veb(M),
  \end{alignedat}
  \end{equation}
  which implies~\eqref{EqGAGebpMetric}.
\end{proof}

\begin{lemma}[$g_0$ as a weighted $\betbop$-metric on $M$]
\label{LemmaGAGbetb}
  Then any stationary and asymptotically flat metric $g_0$ is a nondegenerate weighted b-edge-3b-metric on $\{t_*+r\geq-\half r-1\}$; more precisely, on this set we have
  \begin{equation}
  \label{EqGAGbetbMetric}
  \begin{split}
    g_0^{-1} &\in \cA^{\bigish(((2,0),2+\delta),\ ((2,0),2+2\delta),\ ((2,0),2+\delta),\ (0,0)\bigish)}(M;S^2\,\Tbetb M), \\
    g_0 &\in \cA^{\bigish(((-2,0),-2+\delta),\ ((-2,0),-2+2\delta),\ ((-2,0),-2+\delta),\ (0,0)\bigish)}(M;S^2\,\Tbetb^*M).
  \end{split}
  \end{equation}
  Moreover, $g_0^{-1}$ is equal to the dual Minkowski type metric $g_M^{-1}$ (see Example~\usref{ExGSGMink}) to leading order at $I^0$, $\scri^+$, and $\iota^+$; that is, on $\{t_*+r\geq-\half r-1\}$ we have
  \begin{equation}
  \label{EqGAGbetbMetMink}
    g_0^{-1} - (-2\pa_{t_*}\otimes_s\pa_r + \pa_r^2 + r^{-2}\slg^{-1}) \in \cA^{(2+\delta,2+2\delta,2+\delta,0)}(M;S^2\,\Tbetb M).
  \end{equation}
\end{lemma}

In the case that $g_0^{-1}$ is smooth, the memberships~\eqref{EqGAGbetbMetric}--\eqref{EqGAGbetbMetMink} read
\[
  g_0^{-1}\in\rho_0^2 x_{\!\scri}^2 \rho_+^2\CI(M;S^2\,\Tbetb M),\qquad
  g_0\in\rho_0^{-2}x_{\!\scri}^{-2}\rho_+^{-2}\CI(M;S^2\,\Tbetb^*M),
\]
and $g_0^{-1}-(-2\pa_{t_*}\otimes\pa_r+\pa_r^2+r^{-2}\slg^{-1})\in\rho_0^3 x_{\!\scri}^3\rho_+^3\CI(M;S^2\,\Tbetb M)$, respectively.

\begin{proof}[Proof of Lemma~\usref{LemmaGAGbetb}]
  Near $\scri^+$, the memberships~\eqref{EqGAGbetbMetric} are a consequence of Lemma~\ref{LemmaGAGeb} and the fact that the explicit leading order terms in~\eqref{EqGAGeb0Metric} and \eqref{EqGAGebpMetric} are nondegenerate Lorentzian signature quadratic forms (in the edge-b-vector fields $\rho_0\pa_{\rho_0}$, resp.\ $\rho_+\pa_{\rho_+}$, and $x_{\!\scri}\pa_{x_{\!\scri}}$, $x_{\!\scri} T\Sph^{n-1}$) upon factoring out $\rho_0^2 x_{\!\scri}^2$, resp.\ $x_{\!\scri}^2\rho_+^2$. Away from $\scri^+$, the claim follows from the fact---a consequence of Lemma~\ref{LemmaGA}---that $g_0$ is a nondegenerate Lorentzian signature section of the pullback $S^2\wt T^*M_1\to M_1$ of the bundle of scattering 2-tensors on $\ol{\R^{n+1}}$ away from the light cone at future infinity, and therefore a weighted nondegenerate b-metric away from $\scri^+\cup\cT^+$, and (by definition) a 3b-metric near $\cT^+$; since $M_1\setminus\scri^+=M\setminus\scri^+$ as smooth manifolds, this implies the claim. One can also argue using the explicit local frames of $\Tbetb M$ given after Definition~\ref{DefGAGVF}.

  For the final claim, it suffices to note that $g_M$ is stationary and asymptotically flat, and hence~\eqref{EqGAGbetbMetric} holds for $g_M$ in place of $g_0$; and moreover, the difference $g_0^{-1}-g_M^{-1}$, as given in~\eqref{EqGSGMetMink}, vanishes (and indeed is conormal with weights $\delta,2\delta,\delta,0$ at $I^0,\scri^+,\iota^+,\cT^+$) as a scattering 2-cotensor at $M\cap\rho^{-1}(0)=I^0\cup\scri^+\cup\iota^+$.
\end{proof}

\begin{definition}[Admissible asymptotically flat metrics]
\label{DefGAG}
  Let $g_0$ be a stationary and asymptotically flat Lorentzian metric on $M_0=\ol{\R_{t_*}}\times X$. Let $\delta\in(0,1]$ be as in Definition~\usref{DefGSG}. Let $\ell_0,\ell_+\in(0,\delta]$, $\ell_{\!\scri}\in(0,\frac{\delta}{2}]$, and $\ell_\cT>0$. A Lorentzian metric $g\in\CI(M^\circ;S^2 T^*M^\circ)$ is then called an \emph{$(\ell_0,2\ell_{\!\scri},\ell_+,\ell_\cT)$-admissible asymptotically flat metric (relative to $g_0$)} if
  \begin{equation}
  \label{EqGAGMetric}
    g^{-1} - g_0^{-1} \in \rho_0^2 x_{\!\scri}^2\rho_+^2 \cA^{(\ell_0,2\ell_{\!\scri},\ell_+,\ell_\cT)}(M;S^2\,\Tbetb M)
  \end{equation}
  on $\{t_*+r\geq-\half r-1\}$ and if, moreover, there exists a smooth function $\ft_*\in\CI(M^\circ)$ with $\ft_*-t_*\in\cA^{2\ell_{\!\scri}}(X)$ so that $\dd\ft_*$ is past timelike in $M^\circ\cap\{\ft_*\geq -1\}$, and so that $t_*+r\geq-\half r-1$ in $\{\ft_*\geq -1\}$.
\end{definition}

Near $\cT^+$, we thus require $g$ to asymptote to $g_0$. By Lemma~\ref{LemmaGAGbetb}, we can replace $g_0$ in~\eqref{EqGAGMetric} near $\scri^+$, or indeed on $M\setminus\cT^+$, by the Minkowski type metric $g_M$ (which is smooth, rather than merely conormal). We then recall from \cite[Lemma~3.4]{HintzVasyScrieb} that functions $\ft_*$ of the required form exist locally near $\scri^+$; one can take $\ft_*=t_*-C x_{\!\scri}^{2\ell_{\!\scri}}$ where $C>0$ is sufficiently large. Since $g_0$ is nondegenerate, a condition equivalent to~\eqref{EqGAGMetric} is
\[
  g-g_0\in\rho_0^{-2}x_{\!\scri}^{-2}\rho_+^{-2}\cA^{(\ell_0,2\ell_{\!\scri},\ell_+,\ell_\cT)}(M;S^2\,\Tbetb^*M).
\]
We stress that $g-g_0$ is required to be of lower order relative to $g_0$ also at $I^0$ and $\iota^+$; thus, the definition of admissibility given here is more restrictive than \cite[Definition~3.2]{HintzVasyScrieb}.\footnote{The present definition ensures that the null-bicharacteristic flow over $\iota^+$ is the same for $g$ as for $g_0$. A similar assumption in Definition~\ref{DefGAW} below for associated wave type operators ensures that the leading order behavior of waves at $\iota^+$ can be inferred from properties of the stationary model. More general settings can in principle be considered, but appear to be of little importance in applications.}

\begin{rmk}
\label{RmkGAGWeights}
  For the purposes of the present paper, no strength of any result is lost if one takes $\ell_0=2\ell_{\!\scri}=\ell_+=\ell_\cT=\delta\in(0,1]$ to be all equal. In the presence of zero energy resonances, as in \cite{HintzVasyNonstat2}, it will be important to keep the weight $\ell_\cT$ separate from the others; and the full set of weights is needed only in applications to quasilinear wave equations where the metric decay and the decay or solutions of linear wave equations are coupled. We use the present notation with such applications in mind.
\end{rmk}

For completeness, we end with supplementary results (not used in the remainder of the paper) in which we make Definition~\ref{DefGAG} concrete in local coordinates. This is useful in applications for checking admissibility in concrete cases; see e.g.\ Example~\ref{ExGAGEx}.

\begin{lemma}[Allowed perturbations of the dual metric]
\label{LemmaGAGPert}
  Let $g_0$ be stationary and asymptotically flat. Let $\chi_{\cT^+}$, resp.\ $\chi_{\scri^+}\in\CI(M)$ be identically $1$ near $\cT^+$, resp.\ $\scri^+$, supported in a collar neighborhood of $\cT^+$, resp.\ $\scri^+$, and with $t_*\geq 1$ on $\supp\chi_{\cT^+}$ and $r\geq 1$ on $\supp\chi_{\scri^+}$. Use the coordinates $t_*,x=(x^1,\ldots,x^n)$ on $M^\circ$, and also spatial polar coordinates $x=r\omega$, $\omega\in\Sph^{n-1}$; we work in $\{t_*+r\geq-\half r-1\}$. Then a Lorentzian metric $g$ is {$(\ell_0,2\ell_{\!\scri},\ell_+,\ell_\cT)$-admissible asymptotically flat metric relative to $g_0$} if and only if the following conditions hold.
  \begin{enumerate}
  \item\label{ItGAGPertT}{\rm (Near $\cT^+$.)} Letting $\rho_\cT=\frac{\la r\ra}{t_*}$ and $\rho_+=\la r\ra^{-1}$, symmetric 2-cotensor $\chi_{\cT^+}(g^{-1}-g_0^{-1})$ is a linear combination of
    \begin{equation}
    \label{EqGAGPertTSections}
      \pa_{t_*}^2,\qquad \pa_{t_*}\otimes_s\pa_{x^j},\qquad \pa_{x^i}\otimes_s\pa_{x^j}\qquad (1\leq i,j\leq n),
    \end{equation}
    with coefficients which lie in $\rho_\cT^{\ell_\cT}\rho_+^{\ell_+}L^\infty(M)$ together with all conormal derivatives, i.e.\ derivatives along all finite products of $t_*\pa_{t_*}$, $\la x\ra\pa_x$.
  \item\label{ItGAGPertScri}{\rm (Near $\scri^+$.)} Letting $\rho_0=\frac{1}{1+(t_*)_-}$, $\rho_+=\frac{1}{1+(t_*)_+}$, and $x_{\!\scri}^2=\frac{r^{-1}}{\rho_0\rho_+}$, the cotensor $\chi_{\scri^+}(g^{-1}-g_0^{-1})$ is a linear combination of the symmetric 2-cotensors\footnote{Here, as before, $\pa_\omega$ is a schematic notation for a spherical vector field; two occurrences $\pa_\omega$ need not denote the same vector field.}
    \begin{equation}
    \label{EqGAGPertScriSections}
    \begin{alignedat}{3}
      &x_{\!\scri}^2\pa_{t_*}^2, &\qquad& \pa_{t_*}\otimes_s\pa_r, &\qquad& x_{\!\scri}\pa_{t_*}\otimes_s r^{-1}\pa_\omega, \\
      &x_{\!\scri}^{-2}\pa_r^2, &\qquad& x_{\!\scri}^{-1}\pa_r\otimes_s r^{-1}\pa_\omega, &\qquad& r^{-1}\pa_\omega\otimes r^{-1}\pa_\omega,
    \end{alignedat}
    \end{equation}
    with coefficients which lie in $\rho_0^{\ell_0}x_{\!\scri}^{2\ell_{\!\scri}}\rho_+^{\ell_+}L^\infty(M)$ together with all conormal derivatives, i.e.\ derivatives along all finite products of $t\pa_t+r\pa_r$ (scaling), $r(\pa_t+\pa_r)$ (weighted derivative along outgoing null cones), $\pa_\omega$, where $t=t_*+r$.
  \item\label{ItGAGPertOther}{\rm (Everywhere else.)} Letting $\rho_0=\frac{1+(t_*)_+}{\la r\ra}$ and $\rho_+=\frac{1}{1+(t_*)_+}$, the cotensor $(1-\chi_{\cT^+}-\chi_{\scri^+})(g^{-1}-g_0^{-1})$ is a linear combination of~\eqref{EqGAGPertTSections} with coefficients which lie in $\rho_0^{\ell_0}\rho_+^{\ell_+}L^\infty(M)$ together with all conormal derivatives, i.e.\ derivatives along all finite products of $\la t_*\ra\pa_{t_*}$, $\la x\ra\pa_x$.
  \end{enumerate}
\end{lemma}
\begin{proof}
  The cotensors~\eqref{EqGAGPertTSections} are a smooth (on $M$) frame of the bundle $S^2\wt T M$. Thus, part~\eqref{ItGAGPertT} follows at once from the fact that $(\frac{t_*}{\la r\ra})^{-1}=\frac{\la r\ra}{t_*}$, resp.\ $\la r\ra^{-1}$ is a local defining function of $\cT^+$, resp.\ $\iota^+$ near $\supp\chi_{\cT^+}$. For part~\eqref{ItGAGPertScri}, we recall the calculations~\eqref{EqGAGeb0VF}--\eqref{EqGAGebpVF}, which show that the 2-cotensors listed in~\eqref{EqGAGPertScriSections} are a frame of $\rho_0^2 x_{\!\scri}^2\rho_+^2\CI(M;S^2\,\Teb M)$ near $\supp\chi_{\scri^+}$. For part~\eqref{ItGAGPertOther} finally, we use the relationship between $\Tbetb M$ and $\Tsc\ol{\R^{n+1}}$ described after Definition~\ref{DefGAGVF}.
\end{proof}

\begin{lemma}[Allowed perturbations of the metric]
\label{LemmaGAGPertMet}
  In the notation of Lemma~\usref{LemmaGAGPert}, the Lorentzian metric $g$ is $(\ell_0,2\ell_{\!\scri},\ell_+,\ell_\cT)$-admissible relative to $g_0$ if and only if, in $\{t_*+r\geq-\half r-1\}$, the tensor $\chi_{\cT^+}(g-g_0)$ is a linear combination of $\dd t_*^2$, $\dd t_*\otimes_s\dd x^j$, $\dd x^i\otimes_s\dd x^j$ with coefficients in the spaces stated in Lemma~\usref{LemmaGAGPert}\eqref{ItGAGPertT} and \eqref{ItGAGPertOther}, and if in the notation of Lemma~\usref{LemmaGAGPert}\eqref{ItGAGPertScri}, $\chi_\scri^+(g-g_0)$ is a linear combination (with coefficients as there) of\footnote{Here $\dd\omega$ is a spherical 1-form, and no two $\dd\omega$ need to be the same 1-forms.}
  \begin{equation}
  \label{EqGAGPertMet0}
  \begin{alignedat}{3}
    &x_{\!\scri}^{-2}\,\dd t_*^2,&\qquad& \dd t_*\otimes_s\dd r, &\qquad& x_{\!\scri}^{-1}\dd t_*\otimes_s r\,\dd\omega, \\
    &x_{\!\scri}^2\,\dd r^2, &\qquad& x_{\!\scri}\,\dd r\otimes_s r\,\dd\omega, &\qquad& r\,\dd\omega\otimes_s r\,\dd\omega,
  \end{alignedat}
  \end{equation}
  or equivalently of
  \begin{equation}
  \label{EqGAGPertMet}
  \begin{alignedat}{3}
    &x_{\!\scri}^2(\dd x^0)^2,&\qquad& \dd x^0\otimes_s\dd x^1, &\qquad& x_{\!\scri}\,\dd x_0\otimes_s r\,\dd\omega, \\
    &x_{\!\scri}^{-2}(\dd x^1)^2, &\qquad& x_{\!\scri}^{-1}\dd x^1\otimes_s r\,\dd\omega, &\qquad& r\,\dd\omega\otimes_s r\,\dd\omega,
  \end{alignedat}
  \end{equation}
  where $x^0=t+r=t_*+2 r$ and $x^1=t_*$.\footnote{These are null coordinates for any Minkowski type metric, which is equal to $-\dd x^0\,\dd x^1+r^2\slg$.}
\end{lemma}
\begin{proof}
  The symmetric 2-tensors in~\eqref{EqGAGPertMet0} are the tensors dual to~\eqref{EqGAGPertScriSections}. Expressing them in terms of $x^0$ and $x^1$ via $\dd t_*=\dd x^1$ and $2\dd r=\dd x^0-\dd x^1$ and taking suitable linear combinations gives~\eqref{EqGAGPertMet}.
\end{proof}

We include the description~\eqref{EqGAGPertMet} here since double null coordinates are particularly well-suited to computations near $\scri^+$, see e.g.\ \cite[\S3]{HintzMink4Gauge}.

\begin{example}[Examples of admissible asymptotically flat metrics]
\label{ExGAGEx}
  Restricting to $t_*\geq 1$ for notational simplicity, and taking $g_0$ to be the Minkowski metric or any stationary and asymptotically flat metric (such as a metric asymptoting as $r\to\infty$ to a Kerr metric as in Example~\ref{ExGSGKerr}), examples of admissible asymptotically flat metrics are metrics $g$ for which $g-g_0$ is a linear combination of $\dd t^2$, $\dd t\otimes_s\dd x^i$, $\dd x^i\otimes_s\dd x^j$ ($1\leq i,j\leq n$) with coefficients lying in $r^{-1-\delta}t_*^{-\delta}L^\infty=x_{\!\scri}^{2+2\delta}\rho_+^{1+2\delta}\rho_\cT^\delta L^\infty$ together with all derivatives along $t_*\pa_{t_*}$ and $\la x\ra\pa_x$. (The extra order of vanishing at $r^{-1}=0$ ensures the order $\delta$ vanishing when expanding into the tensors in~\eqref{EqGAGPertMet0}.) Such metrics are thus $(\delta,\delta,\delta,\delta)$-admissible. In $n+1\geq 5$ spacetime dimensions, such metrics arise as perturbations of the Minkowski metric in the context of the stability problem \cite{ChoquetBruhatChruscielLoizeletEinsteinMaxwell,AndersonChruscielSimple,LindbladRodnianskiGlobalStability}. For $n+1=4$, the $r^{-1-\delta}$ vanishing of such metrics near $(\scri^+)^\circ$ is not valid in the context of the stability problem \cite{LindbladRodnianskiGlobalStability,HintzVasyMink4}; this is when the precision of Lemma~\ref{LemmaGAGPertMet} (which e.g.\ allow for less vanishing of the $(\dd x^1)^2$ and $r\,\dd\omega\otimes_s r\,\dd\omega$ coefficients which play a crucial role in \cite[Theorem~8.14]{HintzVasyMink4}) become important. See \cite{HintzMink4Gauge} for details near $I^0\cap\scri^+$. Further examples, arising from quasilinear wave equations, will be discussed elsewhere.
\end{example}

\subsubsection{Wave type operators}
\label{SssGAW}

We first relate our class of stationary wave type operators to $\betbop$-differential operators. Lemma~\ref{LemmaGAGbetb} suggests the following result:

\begin{lemma}[Stationary wave type operators as $\betbop$-operators]
\label{LemmaGAWStatbetb}
  Let $P_0$ be a stationary wave type operator (relative to $g_0$, see Definition~\usref{DefGSO}) acting on sections of the stationary vector bundle $E\to M_0$. Then
  \begin{equation}
  \label{EqGAWStatbetb}
    P_0\in\bigl(\rho_0^2 x_{\!\scri}^2\rho_+^2\CI+\cA^{(2+\delta,\ 2+2\delta,\ 2+\delta,\ (0,0))}\bigr)\Diffbetb^2(M;\upbeta^*E).
  \end{equation}
  (When the conormal terms in Definition~\usref{DefGSO} are in fact smooth, then we simply have $P_0\in\rho_0^2 x_{\!\scri}^2\rho_+^2\Diffbetb^2(M;\upbeta^*E)$.)
\end{lemma}
\begin{proof}
  For bounded $r=|x|$, the membership~\eqref{EqGAWStatbetb} asserts smoothness of the coefficients of $P_0$ in terms of $\pa_{t_*},\pa_x$ as functions of $(t_*^{-1},x)$, which is indeed true.

  Note now that $\cA^\delta(X)$ lifts (by stationarity) to $M$ as a subspace of $\cA^{(\delta,2\delta,\delta,(0,0))}(M)$, and $\rho=r^{-1}$ lifts to an element of $\rho_0 x_{\!\scri}^2\rho_+\CI(M)$. Similarly, stationary lifts of elements of $\Vb(X)$ (which over $r\gtrsim 1$ are spanned over $\CI(X)$ by $\rho\pa_\rho$, $\pa_\omega$) to $M$ lie in $\Vtb(M)$ near $\cT^+$, and in $x_{\!\scri}^{-1}\Veb(M)$ near $\scri^+$; cf.\ the explicit generators of $\Vbetb(M)$ listed after Definition~\ref{DefGAGVF}. The claim then follows from upon inspection of~\eqref{EqGSOStruct} from
  \begin{equation}
  \label{EqGAWStatbetbPf}
    \Vb(X) \subset x_{\!\scri}^{-1}\Vbetb(M),\qquad
    \pa_{t_*} \in \rho_0\rho_+\Vbetb(M).\qedhere
  \end{equation}
\end{proof}

The main goal of the second part of this paper is then the study of the following class of operators:

\begin{definition}[Wave type operators]
\label{DefGAW}
  Let $g$ be an $(\ell_0,2\ell_{\!\scri},\ell_+,\ell_\cT)$-admissible asymptotically flat metric (relative to the stationary and asymptotically flat metric $g_0$) on $M$ (see Definitions~\usref{DefGSG}, \usref{DefGAMfd}, \usref{DefGAG}). Let $P_0$ be a stationary wave type operator (relative to $g_0$, see Definition~\usref{DefGSO}) acting on sections of the stationary vector bundle $E\to M_0$. Then an operator $P\in\Diff^2(M^\circ;E)$ is called an \emph{admissible wave type operator} (with respect to $g$ and $P_0$) if:
  \begin{enumerate}
  \item\label{ItGAWSpec} the \emph{stationary model} $P_0$ is spectrally admissible (Definition~\usref{DefGSOSpec});
  \item\label{ItGAWSymb} the principal symbol of $P$ is scalar, and equal to the dual metric function of $\zeta\mapsto g^{-1}(\zeta,\zeta)$ of $g$;
  \item\label{ItGAWStruct} we have $P=P_0+\tilde P$, where $\tilde P=\tilde P_1+\tilde P_2$ with
    \begin{equation}
    \label{EqGAWStruct}
    \begin{split}
      \tilde P_1 &\in \rho_0^2 x_{\!\scri}^2\rho_+^2\cA^{(\ell_0,\;((0,0),2\ell_{\!\scri}),\;\ell_+,\;\infty)}\Diffbetb^1(M;\upbeta^*E), \\
      \tilde P_2 &\in \rho_0^2 x_{\!\scri}^2\rho_+^2\cA^{(\ell_0,2\ell_{\!\scri},\ell_+,\ell_\cT)}\Diffbetb^2(M;\upbeta^*E);
    \end{split}
    \end{equation}
  \item\label{ItGAWEdgeN} there exist $p_0\in\cA^{(\ell_0+1,\ell_++1)}(\scri^+;\upbeta^*\End(E))$ and $\tilde p_1\in\cA^{(\ell_0,\ell_+)}(\scri^+;\upbeta^*\End(E))$ so that in a collar neighborhood of $\scri^+$  and in the coordinates~\eqref{EqGAGeb0} (with $T=T^0$), resp.\ \eqref{EqGAGebp} (with $T=T^+<T^0$), and upon setting $p_0^0=\rho_0^{-1}p_0\in\cA^{\ell_0}(\scri^+\setminus I^+;\End(E))$, $p_0^+=\rho_+^{-1}p_0\in\cA^{\ell_+}(\scri^+\setminus I^0;\End(E))$, and $p_1:=S|_{\pa X}+\tilde p_1$ where $S$ is defined in~\eqref{EqGSOStruct}, we have
    \begin{equation}
    \label{EqGAWEdgeN}
    \begin{split}
      2\rho_0^{-2}x_{\!\scri}^{-2}P &\equiv -\Bigl(x_{\!\scri} D_{x_{\!\scri}}-2 i^{-1}\Bigl(\frac{n-1}{2}+p_1\Bigr)\Bigr)(x_{\!\scri} D_{x_{\!\scri}}-2\rho_0 D_{\rho_0}) + 2 x_{\!\scri}^2\slDelta + p_0^0, \\
      2x_{\!\scri}^{-2}\rho_+^{-2}P &\equiv \Bigl(x_{\!\scri} D_{x_{\!\scri}}-2 i^{-1}\Bigl(\frac{n-1}{2}+p_1\Bigr)\Bigr)(x_{\!\scri} D_{x_{\!\scri}}-2\rho_+ D_{\rho_+}) + 2 x_{\!\scri}^2\slDelta + p_0^+,
    \end{split}
    \end{equation}
    respectively, modulo $\cA^{(((0,0),\ell_0),\;2\ell_{\!\scri},\;((0,0),\ell_+))}\Diffeb^2(M;\upbeta^*E)$ (with edge behavior at $\scri^+$ and b-behavior at $I^0$ and $\iota^+$);\footnote{The zeroth order terms in~\eqref{EqGAWEdgeN} are consistent; multiplying the first, resp.\ second equation by $\rho_0$, resp.\ $\rho_+$, and noting that $\rho_0^{-1}x_{\!\scri}^{-2}=r^{-1}=x_{\!\scri}^{-2}\rho_+^{-1}$, this follows from $\rho_0 p_0^0=p_0=\rho_+ p_0^+$.}
  \item\label{ItGAWSpecial} there exists a bundle splitting $E|_{\pa X}=\oplus_{j=1}^J E_j$ so that $p_1$ is lower triangular, with diagonal entries $p_{1,j j}\in\End(E_j)$ having real spectrum, and so that $p_0$ is strictly lower triangular everywhere on $\scri^+$.
  \end{enumerate}
\end{definition}

\begin{rmk}[Comments on Definition~\usref{DefGAW}]
\fakephantomsection
\label{RmkGAWComments}
  \begin{enumerate}
  \item In~\eqref{EqGAWEdgeN}, the operator $\slDelta$ is any element of $\Diffb^2(\scri^+;E|_{\pa X})$ (extended to the collar neighborhood $[0,\eps)_{x_{\!\scri}}\times\scri^+$ of $\scri^+$ to be independent of $x_{\!\scri}$) which acts on fiber-constant sections on $\scri^+=[0,\infty]\times\Sph^{n-1}$ via any fixed second order operator on $\Sph^{n-1}$ whose principal symbol is scalar and equal to the dual metric function of $\slg$. Since the difference of any two such operators is a second order b-differential operator on $\scri^+$ which is vertical, i.e.\ involves derivatives only along the fibers of $\scri^+$, the operator $x_{\!\scri}^2\slDelta$ on the collar neighborhood is well-defined modulo $x_{\!\scri}\Diffbetb^2(M;\upbeta^*E)$.
  \item\label{ItGAWCommentsScriEB} An admissible wave type operator is admissible also in the sense of \cite[Definition~3.5]{HintzVasyScrieb} (which concerns only the structure of $P$ near $\scri^+$).\footnote{Regarding the requirements on the bundle in \cite[Definition~3.5]{HintzVasyScrieb}, we note that the pullback $\upbeta^*E\to M$ of the stationary bundle $E=\pi_0^*(E_X)$ (using the notation of~\S\ref{SsGS}) is, by Lemma~\ref{LemmaGA}, in a neighborhood of $\scri^+$ equal to the pullback of a bundle $\wt E\to\ol{\R^{n+1}}$ (defined near $Y_+$ in the notation of Lemma~\ref{LemmaGA}) to $M$. Indeed, one may take $\wt E$ to be the pullback of $E_X$ along the map $(t_*,x)\mapsto x$, which extends to a smooth map from the closure of $\{|x|>1,\ \frac{t_*}{|x|}<\half\}$ in $\ol{\R^{n+1}}$ to the closure of $\{|x|>1\}$ in $X$.} Our present definition is stronger and requires $P$ to be equal to the stationary model $P_0$ to leading order also at $I^0$ and $\iota^+$, and it moreover requires $p_0,p_1$ to vanish at $\pa\scri^+=(\scri^+\cap I^0)\cup(\scri^+\cap\iota^+)$. The leading order equality of $P$ and $P_0$ at $\cT^+$ is essential for our analysis. The equality at $\iota^+$ on the other hand, while frequently used in our analysis, can in principle be relaxed, at the expense of necessitating additional hypotheses on the structure of the null-bicharacteristic flow over $\iota^+$ (cf.\ the comments following Definition~\ref{DefGAG}) and the absence of resonances for the Mellin-transformed normal operator family in suitable half spaces; see in particular Lemma~\ref{LemmaWFlowI}, Proposition~\ref{PropWFlow}, and~\S\ref{SsWip} for the results whose validity would need to be \emph{assumed} if one did not demand the equality of $P$ and $P_0$ at $\iota^+$. This is also closely related to \cite{BaskinVasyWunschRadMink,BaskinVasyWunschRadMink2}.
  \item Condition~\eqref{ItGAWEdgeN} can be formulated equivalently using the edge normal operator, see \cite[Definition~3.5]{HintzVasyScrieb}. Note also that the principal part of the right hand side of~\eqref{EqGAWEdgeN} is determined by the dual metric function of $g$, which to leading order at $\scri^+$ is equal to the dual metric function of $g_0$ and can thus be read off from the expressions~\eqref{EqGAGeb0Metric} and \eqref{EqGAGebpMetric}. The new information in~\eqref{EqGAWEdgeN} is therefore the structure of the lower order terms, encoded by $p_0,p_1$. The terms $p_0$ and $\tilde p_1=p_1-S|_{\pa X}$ encode the allowed leading order behavior (in the sense of decay) of $\tilde P_1$ at $\scri^+$.
  \item Condition~\eqref{ItGAWSpecial} is the special admissibility condition of \cite[Definition~6.2]{HintzVasyScrieb}. Though we shall not do so in this paper, it can be relaxed at the expense of modifications to the weights at $I^0$, $\scri^+$, and $\iota^+$ of the function spaces in which the analysis in~\S\ref{SW} takes place. See \cite[Remark~6.5]{HintzVasyScrieb}.
  \end{enumerate}
\end{rmk}

We proceed to give a number of examples.

\begin{lemma}[Stationary operators]
\label{LemmaGAWStat}
  Let $P_0$ be stationary wave type operator (relative to $g_0$) that is spectrally admissible; assume moreover that there exists a bundle splitting $E|_{\pa X}=\bigoplus_{j=1}^J E_j$ so that $S|_{\pa X}$ is lower triangular, with diagonal entries having real spectrum. Then $P_0$ is an admissible wave type operator (with respect to $g_0$ and $P_0$), with $p_0=\tilde p_1=0$ in the notation of Definition~\usref{DefGAW}\eqref{ItGAWEdgeN}.
\end{lemma}
\begin{proof}
  Indeed, conditions~\eqref{ItGAWSpec} and \eqref{ItGAWSymb} of Definition~\ref{DefGAW} are satisfied by assumption. Condition~\eqref{ItGAWStruct} holds with $\tilde P_1=\tilde P_2=0$. Upon plugging the expressions~\eqref{EqGAGeb0VF} and \eqref{EqGAGebpVF}, respectively, into~\eqref{EqGSOStruct}, one obtains the two leading order expressions in~\eqref{EqGAWEdgeN} with $p_0=0$ and $\tilde p_1=0$, i.e.\ $p_1=S|_{\pa X}$; here one uses that $\hat P(0)\equiv D_r^2+r^{-2}\slDelta\bmod\rho^2\Diffb^1(X;E_X)+\cA^{2+\delta}\Diffb^2(X;E_X)$, as well as the $r\pa_r=-\rho\pa_\rho=-\half x_{\!\scri}\pa_{x_{\!\scri}}$ and the memberships~\eqref{EqGAWStatbetbPf}. Condition~\eqref{ItGAWSpecial} is then the same as the present assumption on $S$.
\end{proof}

This result applies to Examples~\ref{ExGSO} and \ref{ExGSOPot}. In the non-stationary setting, we note:

\begin{example}[Nonstationary wave operators coupled with stationary potentials or first order terms]
\label{ExGAWWave}
  The wave operator $\Box_g$ of an admissible asymptotically flat metric (relative to a stationary and asymptotically flat metric $g_0$) satisfies conditions~\eqref{ItGAWSymb}--\eqref{ItGAWSpecial} of Definition~\ref{DefGAW}, with $p_0=0$ and $\tilde p_1=0$. If $\Box_{g_0}$ is spectrally admissible (that is, condition~\eqref{ItGAWSpec} is also satisfied), then $\Box_g$ is an admissible wave type operator. (We prove this in Proposition~\ref{PropES2}.) More generally, these statements remain valid, mutatis mutandis, for $\Box_g+V_0$ and $\Box_{g_0}+V_0$ where $V_0\in\la x\ra^{-2}\CI(X;E_X)$ is a stationary potential with (approximately) inverse square decay. One can more generally allow $V_0\in\la x\ra^{-2}\CI(X;E_X)+\cA^{2+\delta}(X;E_X)$, $\delta>0$, or first order terms as in Example~\ref{ExGSOPot}.
\end{example}

\begin{example}[Coupling with non-stationary potentials]
\label{ExGAWPotential}
  If $P$ is an admissible wave type operator, then so is $P+\tilde V$ for all $\tilde V\in\rho_0^2 x_{\!\scri}^2\rho_+^2\cA^{(\ell_0,2\ell_{\!\scri},\ell_+,\ell_\cT)}(M;\End(\upbeta^*E))$. In $t_*\geq 1$, and for $\ell_0=\ell_{\!\scri}=\ell_+=\ell_\cT=\delta$, this means that the matrix elements of $V$ in local trivializations of $E_X$ lie in $t_*^{-\delta}\la r\ra^{-2}\frac{t}{\la t-r\ra}L^\infty$ together with all derivatives along $t_*\pa_{t_*}$ and $\la x\ra\pa_x$; this is precisely the condition~\eqref{EqIVDecay}. We recall that this class of potentials includes those which are conormal relative to $t_*^{-p}r^{-q}L^\infty$ with $p>0$, $q>1$, $p+q>2$ as a special case (upon reducing $\ell_0,\ell_{\!\scri},\ell_+,\ell_T>0$ if necessary). One can also allow for additional zeroth order terms which have a $r^{-1}$ leading order term at $\scri^+$ and an overall $r^{-1}\la t_*\ra^{-1-\delta}$ decay rate (cf.\ the term $p_0$ in~\eqref{EqGAWEdgeN}), which arise in stability problems for the Einstein equations in $3+1$ spacetime dimensions; see \cite[Lemma~3.8]{HintzVasyMink4} and \cite[Proposition~3.29]{HintzMink4Gauge} (where different notation is used).
\end{example}

\begin{example}[Coupling with first order terms]
\label{ExGAWFirst}
  If $P$ is an admissible wave type operator, say for simplicity with $p_1=0$ in~\eqref{EqGAWEdgeN}, then so is $P+\chi_\scri r^{-1}\tilde p_1(t_*,\omega)\pa_{t_*}$ where $\chi_\scri$ is a cutoff to a small neighborhood of $\scri^+$, and $|\tilde p_1|\lesssim\la t_*\ra^{-\delta}$ and likewise for all derivatives along $\la t_*\ra\pa_{t_*}$ and $\pa_\omega$, provided $\tilde p_1$ is lower triangular and has real spectrum. (Note here that $4 r^{-1}\pa_{t_*}=2 i^{-1}\rho_0^2 x_{\!\scri}^2(x_{\!\scri} D_{x_{\!\scri}}-2\rho_0 D_{\rho_0})$ by~\eqref{EqGAGeb0VF} and compare with~\eqref{EqGAWEdgeN}.) The term $r^{-1}\tilde p_1\pa_{t_*}$ can be thought of as a weak damping term (when the eigenvalues of $\tilde p_1$ are nonnegative), and plays a key role in \cite{HintzVasyMink4,HintzMink4Gauge}. We shall not spell out examples of other first order terms here; these would be arbitrary contributions to the subprincipal part of $\tilde P_2$ (or equivalently to the terms of $\tilde P_1$ which vanish at $\scri^+$ to leading order).
\end{example}

Finally, in the notation of Definition~\ref{DefGAW}\eqref{ItGAWStruct}, we proceed to explain how to measure the size of perturbations $\tilde P$ of $P_0$. First, we make the splitting $\tilde P=\tilde P_1+\tilde P_2$ in~\eqref{EqGAWStruct} unique in the following manner: we fix a collar neighborhood of $\scri^+$ and take $\tilde P_1=\chi_\scri\tilde P_\scri$ where $\chi_\scri\in\CI(M)$ is a fixed cutoff function which is $1$ near $\scri^+$ and has support in the collar neighborhood, and $\tilde P_\scri$ is the (unique) operator, homogeneous of degree $2$ with respect to dilations in the defining function of $\scri^+$ in the collar neighborhood, so that $\tilde P-\chi_\scri\tilde P_\scri\in\rho_0^2 x_{\!\scri}^2\rho_+^2\cA^{(\ell_0,2\ell_{\!\scri},\ell_+,\ell_\cT)}\Diffbetb^2(M;\upbeta^*E)$. (In terms of~\eqref{EqGAWEdgeN}, and further splitting $\tilde P_\scri$ into the sum $\tilde P_\scri=\tilde P_{\scri,0}+\tilde P_{\scri,+}$ of two terms with support disjoint from $\scri^+\cap I^0$, resp.\ $\scri^+\cap\iota^+$ by means of a partition of unity $1=\chi_0+\chi_+$ on $\scri^+$, we can take $2\rho_0^{-2}x_{\!\scri}^{-2}\rho_+^{-2}\tilde P_{\scri,0}=\chi_0(2 i^{-1}p_1(x_{\!\scri} D_{x_{\!\scri}}-2\rho_0 D_{\rho_0})+p_0)$ and $2\rho_0^{-2}x_{\!\scri}^{-2}\rho_+^{-2}\tilde P_{\scri,+}=\chi_+(-2 i^{-1}p_1(x_{\!\scri} D_{x_{\!\scri}}-2\rho_+ D_{\rho_+})+p_0)$.) Next, cover $\{t_*+r\geq-\half r-1\}$ by a finite number of coordinate charts on which $E$ is trivial, and let $\chi_i$ be a partition of unity subordinate to this cover. Fix moreover a finite number of b-vector fields $V_\alpha$ on $M$ which at each point of $M$ span $\Tb M$. Letting $\tilde p_{2,i}\in\rho_0^2 x_{\!\scri}^2\rho_+^2\cA^{(\ell_0,2\ell_{\!\scri},\ell_+,\ell_\cT)}(M)$ denote a coefficient (or matrix element if $E$ is not the trivial bundle) of $\chi_i\tilde P_2$ in such a chart (and trivialization), we let $|\tilde p_{2,i}|_k$ denote the sum of the $L^\infty$-norm of $\rho_0^{-2}x_{\!\scri}^{-2}\rho_+^{-2}\rho_0^{-\ell_0}x_{\!\scri}^{-2\ell_{\!\scri}}\rho_+^{-\ell_+}\rho_\cT^{-\ell_\cT}\tilde p_{2,i}$ and of all its derivatives along $k$-fold compositions of the $V_\alpha$, and $|\chi_i\tilde P_2|_k$ is the sum of all these norms of the coefficients. We similarly define $|\chi_i\tilde P_\scri|_k$ as the $L^\infty$-norm of the restrictions of the coefficients of $\rho_0^{-2}x_{\!\scri}^{-2}\rho_+^{-2}\rho_0^{-\ell_0}\rho_+^{-\ell_+}\tilde P_\scri$ to $\scri^+$, plus analogous norms for up to $k$-fold derivatives of these coefficients along $V_\alpha|_{\scri^+}$. Finally, we let
\begin{equation}
\label{EqGAWNorm}
  |\tilde P|_k,\qquad k\in\N_0,
\end{equation}
denote the (finite) sum of these norms over all $i$.

\section{Analysis of stationary wave type operators}
\label{SSt}

We use the notation of~\S\ref{SsGS} and fix a stationary wave type operator $P$ (Definition~\ref{DefGSO}) with respect to a stationary and asymptotically flat metric $g$ (Definition~\ref{DefGSG}) acting on sections of a stationary bundle $E=\pi_0^*(E_X)\to M_0=\ol{\R_{t_*}}\times X$; we assume that $P$ is spectrally admissible with indicial gap $(\beta^-,\beta^+)$ (Definition~\ref{DefGSOSpec}). We fix a positive definite fiber inner product on $E$, and work with the volume densities $|\dd g|$ and $|\dd g_X|$ on $M_0^\circ$ and $X^\circ$, respectively (Lemma~\ref{LemmaGSGVol}). We use polar coordinates $x=r\omega$ on $X^\circ=\R^n$, and inverse polar coordinates $\rho=r^{-1}$ near $\rho^{-1}(0)\subset X=\ol{\R^n}$.

In~\S\ref{SsStEst}, we upgrade the assumptions made in Definition~\ref{DefGSOSpec} to quantitative mapping and invertibility properties of the spectral family near low or bounded energies, and also show how the nontrapping assumption in Definition~\ref{DefGSG}\eqref{ItGSGNontrap} gives high energy estimates. In~\S\ref{SsStCo}, we use these estimates to prove basic decay results for forward solutions of $P$; see Theorem~\ref{ThmStCo} for the main result in this regard. In~\S\ref{SAS}, we extract the leading order asymptotic profiles for forward solutions of $P$ at $\iota^+$ and $\cT^+$ under an appropriate genericity assumption on $P$, leading to the second main result (Theorem~\ref{ThmAS}) of this paper for stationary wave type operators.

\subsection{Estimates for the spectral family}
\label{SsStEst}

We shall need two types of estimates for the spectral family:
\begin{enumerate}
\item\label{ItStEstb} estimates on weighted b-Sobolev spaces, which are used in~\S\ref{SsStCo} to prove the conormality of the output of the resolvent;
\item\label{ItStEstsc} estimates on scattering (or scattering-b-transition) Sobolev spaces with variable decay order, which are used as normal operator estimates at the translation face in the 3b-analysis in~\S\ref{SsWM}.
\end{enumerate}
These are closely related: the b-estimates are consequences of second microlocal scattering-b-estimates, first obtained by Vasy \cite{VasyLowEnergyLag}, and these can be thought of as sharp variable order scattering estimates where the decay order jumps right at the outgoing radial set (a first instance of which is described in~\S\ref{SssStEstNonzero}). The variable order estimates appeared in the special case of the scalar wave operator on Kerr spacetimes already in \cite[\S3.5]{HintzKdSMS}, but the same arguments go through in the present general setting with only minor modifications. The proofs in this section are relatively small modifications of arguments already appearing in the literature, in particular \cite{MelroseEuclideanSpectralTheory,VasyZworskiScl,VasyMicroKerrdS,VasyLAPLag,VasyLowEnergyLag,HaefnerHintzVasyKerr,HintzPrice}; nonetheless, we give a fair amount of details in the present setting for completeness.

We begin with the analysis of the zero energy operator in~\S\ref{SssStEst0}, followed by bounded nonzero energies in~\S\ref{SssStEstNonzero} and high energies in~\S\ref{SssStEstHi}. Uniform low energy resolvent estimates are proved in~\S\ref{SssStEstLo}.

\subsubsection{Zero frequency}
\label{SssStEst0}

The principal symbol of $\hat P(\sigma)$ as a differential operator on $X^\circ$ is independent of $\sigma\in\C$, and indeed equal to the restriction of the dual metric function $G(\zeta)=g^{-1}(\zeta,\zeta)$ to $\ker\pa_{t_*}$. Since $\pa_{t_*}$ is timelike, this implies that $\hat P(\sigma)$ is elliptic on $X^\circ$.

\begin{lemma}[Zero energy operator]
\label{LemmaStEst0}
  Let $\sfs\in\CI(\Sb^*X)$.\footnote{Here and below, we allow for variable orders as they arise from the conversion of the (necessarily) variable order 3b-orders on the spacetime; see~\S\ref{SsWM}, specifically the proof Proposition~\ref{PropWMFred}. However, unless threshold conditions are required for variable orders (the first instance being Proposition~\ref{PropStEstNz}\eqref{ItStEstNzsc}), the reader may safely assume at first reading that the orders are constant.} Then the operator
  \begin{equation}
  \label{EqStEst0}
    \hat P(0)\colon\Hb^{\sfs,-\frac{n}{2}+\beta}(X;E_X)\to\Hb^{\sfs-2,-\frac{n}{2}+\beta+2}(X;E_X)
  \end{equation}
  is invertible for all $\beta\in(\beta^-,\beta^+)$. Moreover, $\beta^-,\beta^+\in\Re\specb(\rho^{-2}\hat P(0))$, and~\eqref{EqStEst0} is not invertible for $\beta\notin(\beta^-,\beta^+)$.
\end{lemma}

\begin{rmk}[Zero energy weights: II]
\label{RmkStEst0Weights}
  This is the first indication that $(\beta^-,\beta^+)$ is the largest interval of weights which is disjoint from the boundary spectrum and contains one value of $\beta$ for which~\eqref{EqStEst0} is invertible. This is then confirmed in Lemma~\ref{LemmaStEsttf} below, which shows that the restrictions on the values of $\beta$ in condition~\eqref{ItGSOSpectf} of Definition~\ref{DefGSOSpec} arise from the boundary spectrum of $\hat P(0)$ in the same manner.
\end{rmk}

\begin{proof}[Proof of Lemma~\usref{LemmaStEst0}]
  Let $0<\nu_\bop\in\CI(X;\Omegab X)$; thus, near $\rho=0$ we have $\nu_\bop=a\rho^n|\dd g_X|$ where $0<a\in\CI(X)+\cA^\delta(X)$. Due to~\eqref{EqGSGMetricCoeff}, the b-principal symbol of $\rho^{-2}\hat P(0)$ is given, in terms of the coordinates $\xi\frac{\dd\rho}{\rho}+\eta$, $\eta\in T^*\Sph^{n-1}$, by $\xi^2+|\eta|_{\slg^{-1}}^2\bmod\cA^\delta P^2(\Tb^*X)$ (the space $\cA^\delta P^2(\Tb^*X)$ consisting of smooth functions which in local coordinates are homogeneous quadratic polynomials in the fibers of $\Tb^*X$ with coefficients in $\cA^\delta(X)$). Therefore, $\rho^{-2}\hat P(0)$ is an elliptic b-differential operator. Thus, when $\beta\in\R$ is such that $\beta\notin\cS:=\Re\specb(\rho^{-2}\hat P(0))$, the operator
  \begin{equation}
  \label{EqStEst0Op}
    \hat P(0)\colon\Hb^{\sfs,\beta}(X,\nu_\bop;E_X)\to\Hb^{\sfs-2,\beta+2}(X,\nu_\bop;E_X)
  \end{equation}
  is Fredholm; its index is thus constant when $\beta$ varies in a connected component of $\R\setminus\cS$, whereas the index jumps by a nonzero amount when $\beta$ crosses a value in $\cS$ by the relative index theorem \cite[\S6.1]{MelroseAPS}. Moreover, by elliptic regularity, elements of its kernel are automatically conormal, i.e.\ they are elements of $\Hb^{\infty,\beta}(X,\nu_\bop;E_X)\subset\cA^\beta(X;E_X)$ (using Sobolev embedding for this inclusion).

  Now, since $\Hb^{\sfs,\beta}(X,\nu_\bop;E_X)=\Hb^{\sfs,-\frac{n}{2}+\beta}(X,|\dd g_X|;E_X)$, condition~\eqref{ItGSOSpec0} in Definition~\ref{DefGSO} implies that for $\beta\in(\beta^-,\beta^+)$, the operator~\eqref{EqStEst0Op} has trivial kernel. The dual space of its codomain $\Hb^{\sfs-2,-\frac{n}{2}+\beta+2}(X,|\dd g_X|;E_X)$ with respect to the $L^2(X,|\dd g_X|;E_X)$-inner product is
  \[
    \Hb^{-\sfs+2,\frac{n}{2}-\beta-2}(X,|\dd g_X|;E_X)=\Hb^{-\sfs+2,n-2-\beta}(X,\nu_\bop;E_X).
  \]
  Condition~\eqref{ItGSOSpec0} thus also implies that~\eqref{EqStEst0Op} has trivial cokernel for $\beta\in(\beta^-,\beta^+)$. Therefore, we necessarily have $(\beta^-,\beta^+)\subset\R\setminus\cS$. Therefore,~\eqref{EqStEst0Op} is Fredholm for $\beta\in(\beta^-,\beta^+)$ (in particular has closed range), and hence is invertible.
\end{proof}

\subsubsection{Bounded nonzero frequencies}
\label{SssStEstNonzero}

In order to analyze the spectral family
\begin{equation}
\label{EqStEstFamily}
  \hat P(\sigma)=2 i\sigma\rho\Bigl(\rho\pa_\rho-\frac{n-1}{2}-S\Bigr)+\hat P(0)-i\sigma Q+\sigma^2 g^{0 0}
\end{equation}
(see~\eqref{EqGSOSpecFam}) for nonzero $\sigma$, we note first that in the coordinates
\[
  \xi\frac{\dd\rho}{\rho^2} + \frac{\eta}{\rho},\qquad \eta\in T^*\Sph^{n-1},
\]
on $\Tsc^*X$, we have
\begin{equation}
\label{EqStEstNonzeroSymb}
  p := \sigmasc^{2,0}(\hat P(\sigma)) \equiv -2\sigma\xi + \xi^2 + |\eta|_{\slg^{-1}}^2 \bmod \cA^\delta S^2(\Tsc^*X)
\end{equation}
in view of Definition~\ref{DefGSO}\eqref{ItGSOSymb} and the memberships~\eqref{EqGSGMetricCoeff}. Here and below, we write $\cA^\delta$ to indicate the space of symbols (or below: vector fields) which have coefficients of class $\cA^\delta(X)$. The scattering characteristic set of $\hat P(\sigma)$, a subset of $\Tsc^*_{\pa X}X$ in view of the ellipticity in the interior noted previously, is thus
\[
  \Char_\sigma =
  \begin{cases}
    \{ (\omega;\xi,\eta) \in \Tsc^*_{\pa X}X \colon (\sigma-\xi)^2 + |\eta|_{\slg^{-1}}^2 = \sigma^2 \},& 0\neq\sigma\in\R, \\
    \{ (\omega;\xi,\eta) \in \Tsc^*_{\pa X}X \colon \xi=\eta=0 \}, & \sigma\notin\R.
  \end{cases}
\]
We call the zero section $\cR_{\rm out}\subset\Char_\sigma$ of $\Tsc^*_{\pa X}X$ the \emph{outgoing radial set}. For $0\neq\sigma\in\R$, we call $\cR_{\sigma,\rm in}=\{(\omega;\xi,\eta)\colon \xi=2\sigma,\ \eta=0\}\subset\Char_\sigma$ (i.e.\ the graph of $2\sigma\frac{\dd\rho}{\rho^2}=-2\sigma\,\dd r$) the \emph{incoming radial set}. The Hamiltonian vector field $H_p$ of $p$ is (in terms of local coordinates $\omega\in\R^{n-1}$ on $\Sph^{n-1}$, and writing covectors on $\Sph^{n-1}$ as $\eta\,\dd\omega$)
\begin{equation}
\label{EqStEstHam}
\begin{split}
  \sfH_p := \rho^{-1}H_p &= (\pa_\xi p)(\rho\pa_\rho+\eta\pa_\eta) + (\pa_\eta p)\pa_\omega - \bigl((\rho\pa_\rho+\eta\pa_\eta)p\bigr)\pa_\xi - (\pa_\omega p)\pa_\eta \\
    &\equiv 2(\xi-\sigma)(\rho\pa_\rho+\eta\pa_\eta) - 2|\eta|_{\slg^{-1}}^2\pa_\xi + H_{|\eta|_{\slg^{-1}}^2} \bmod \cA^\delta\Vb(\Tsc^*X).
\end{split}
\end{equation}
(The expression on the first line is obtained as follows: if $\tilde\xi\in\R$ and $\tilde\eta\in T^*\Sph^{n-1}$ are the standard momentum coordinates, i.e.\ covectors are $\tilde\xi\,\dd\rho+\tilde\eta$, then $H_p=(\pa_{\tilde\xi}p)\pa_\rho+(\pa_{\tilde\eta}p)\pa_\omega-(\pa_\rho p)\pa_{\tilde\xi}-(\pa_\omega p)\pa_{\tilde\eta}$. Changing variables to $\xi=\rho^2\tilde\xi$ and $\eta=\rho\tilde\eta$ gives the stated expression.) Thus, for $\pm\sigma>0$, the flow of $\pm\sfH_p$ goes from the source $\cR_{\sigma,\rm in}$ to the sink $\cR_{\rm out}$. (This is the conjugated version \cite{VasyLAPLag} of the computations in \cite[\S8]{MelroseEuclideanSpectralTheory}.)

We furthermore introduce the following quantities:
\begin{definition}[Thresholds]
\label{DefStEstThr}
  In terms of the expression~\eqref{EqGSOStruct} for $P$, let
  \[
    \ubar S:=\inf_{p\in\pa X}(\min\Re\spec S(p)),
  \]
  where $\spec S(p)\subset\C$ is the spectrum (in the linear algebra sense) of the linear map $S(p)\colon E_p\to E_p$. Furthermore, we let
  \[
    \ubar S_{\rm in} = \sup_h \inf_{p\in\pa X} \Biggl( \min\spec \frac12\biggl(\bigl(S(p)+S(p)^{*h}\bigr)-\sigmasc^{1,0}\Bigl(\frac{\hat P(0)-\hat P(0)^{*h}}{2 i\rho}\Bigr)\Big|_{2\frac{\dd\rho}{\rho^2}}\biggr) \Biggr),
  \]
  where the supremum is taken over all smooth positive definite fiber inner products on $E_X$ defined near $\pa X$.
\end{definition}

\begin{prop}[Spectral family at nonzero energies]
\label{PropStEstNz}
  Let $\ubar S$, $\ubar S_{\rm in}$ be as in Definition~\usref{DefStEstThr}.
  \begin{enumerate}
  \item\label{ItStEstNzsc}{\rm (Variable order estimates.)} Let $\sfs\in\CI(\Ssc^*X)$, and $\sfr\in\CI(\ol{\Tsc_{\pa X}^*}X)$. Let $\pm\sigma>0$, and suppose that $\pm\sfH_p\sfr\leq 0$ on $\Char_\sigma$; suppose moreover that $\sfr<-\half+\ubar S$ at $\cR_{\rm out}$ and $\sfr>-\half-\ubar S_{\rm in}$ at $\cR_{\sigma,\rm in}$, with $\sfr$ constant near $\cR_{\rm out}$ and $\cR_{\sigma,\rm in}$. Then the operator
    \begin{equation}
    \label{EqStEstNzsc}
      \hat P(\sigma) \colon \bigl\{ u\in\Hsc^{\sfs,\sfr}(X;E_X) \colon \hat P(\sigma)u\in\Hsc^{\sfs-2,\sfr+1}(X;E_X) \bigr\} \to \Hsc^{\sfs-2,\sfr+1}(X;E_X)
    \end{equation}
    is Fredholm and has trivial kernel. If $\sfr$ is a variable order function satisfying these assumptions for all $\sigma$ in a compact subset $I\subset\R\setminus\{0\}$ simultaneously, then we have a uniform estimate
      \begin{equation}
      \label{EqStEstNzscEst}
        \|u\|_{\Hsc^{\sfs,\sfr}(X;E_X)} \leq C \|\hat P(\sigma)u\|_{\Hsc^{\sfs-2,\sfr+1}(X;E_X)},\qquad \sigma\in I.
      \end{equation}
  \item\label{ItStEstNzb}{\rm (Estimates on b-spaces.)} Let $s,\ell\in\R$, and suppose that $\ell<-\half+\ubar S$ and $s+\ell>-\half-\ubar S_{\rm in}$. Then for $\sigma\neq 0$, $\Im\sigma\geq 0$, the operator
    \begin{equation}
    \label{EqStEstNzb}
      \hat P(\sigma) \colon \bigl\{ u\in\Hb^{s,\ell}(X;E_X) \colon \hat P(\sigma)u\in\Hb^{s,\ell+1}(X;E_X) \bigr\} \to \Hb^{s,\ell+1}(X;E_X)
    \end{equation}
    is Fredholm and has trivial kernel. If $\Omega\subset\{\sigma\in\C\colon\sigma\neq 0,\ \Im\sigma\geq 0\}$ is compact, then we have a uniform estimate
    \begin{equation}
    \label{EqStEstNzbEst}
      \|u\|_{\Hb^{s,\ell}(X;E_X)} \leq C\|\hat P(\sigma)u\|_{\Hb^{s,\ell+1}(X;E_X)},\qquad \sigma\in\Omega.
    \end{equation}
  \end{enumerate}
\end{prop}

We shall later show (see Corollary~\ref{CorStEstInv}) that $\hat P(\sigma)$ is in fact invertible in both cases.

\begin{proof}[Proof of Proposition~\usref{PropStEstNz}]
  \pfstep{Part~\eqref{ItStEstNzsc}.} This is essentially standard, see e.g.\ \cite[Proposition~4.13]{VasyMinicourse}: it follows from elliptic estimates, radial point estimates at $\cR_{\rm out}$ and $\cR_{\sigma,\rm in}$, and real principal type propagation in between (which requires the monotonicity of $\sfr$). We shall merely briefly discuss the threshold conditions on $\sfr$; we only consider the case $\sigma>0$.

  The radial point estimate at $\cR_{\rm out}$ is proved using a positive commutator argument which uses a commutant $A=A^*\in\Psisc^{-\infty,2\sfr+1}(X;E_X)$ with principal symbol
  \begin{equation}
  \label{EqStEstscCommA}
    a = \rho^{-2\sfr-1}\chi(\xi^2+|\eta|_{\slg^{-1}}^2)^2\chi(\rho)^2
  \end{equation}
  where $\chi\in\CIc([0,\sigma^2))$ is identically $1$ near $0$ and satisfies $\chi'\leq 0$. One then evaluates the $L^2(X,|\dd g_X|;E_X)$ inner product
  \begin{equation}
  \label{EqStEstscComm}
    2\Im\la\hat P(\sigma)u,A u\ra = \la \sC u,u\ra,\qquad \sC := i\bigl(\hat P(\sigma)^*A-A\hat P(\sigma)\bigr).
  \end{equation}
  Writing\footnote{We use here that $\hat P(\sigma)$ is the sum of a second order scattering differential operator with real scalar principal symbol and scattering differential operators of order $\leq 1$ with coefficients in $\rho\CI+\cA^{1+\delta}$, as follows from~\eqref{EqGSOSpecFam}, Definition~\ref{DefGSO}\eqref{ItGSOStruct}, and the fact that $\rho^j\Diffb^j\subset\Diffsc^j$, $j\in\N_0$.} $\Im\hat P(\sigma)=(2 i)^{-1}(\hat P(\sigma)-\hat P(\sigma)^*)\in\Diffsc^{1,-1}(X;E_X)$ and
  \[
    q_\sigma=\sigmasc^{1,-1}\bigl(\Im\hat P(\sigma)\bigr),
  \]
  the principal symbol of $\sC=i[\hat P(\sigma),A]+2(\Im\hat P(\sigma))A\in\Psisc^{-\infty,2\sfr}(X;E_X)$ is $H_p a+2 q_\sigma a$. At $\cR_{\rm out}$, we have $\rho^{2\sfr}H_p a=2\sigma(2\sfr+1)$ by~\eqref{EqStEstHam}. Denoting by $\pi\colon\cR_{\rm out}\to\pa X$ the projection, we claim that for any $\eps>0$ we can choose a positive definite fiber inner product on $E_X$ so that
  \begin{equation}
  \label{EqStEstscLower}
    2\rho^{2\sfr}q_\sigma a|_{\cR_{\rm out}} = (2\rho^{-1}q_\sigma)|_{\cR_{\rm out}}\in\CI(\cR_{\rm out};\pi^*\End(E_X))\leq -4\sigma(\ubar S-\eps)
  \end{equation}
  in the sense of self-adjoint endomorphisms; this implies that $\rho^{2\sfr}\,\sigmasc^{-\infty,2\sfr}(\sC)\leq 2\sigma(2\sfr+1-2(\ubar S-\eps))$ at $\cR_{\rm out}$, which is negative under the stated threshold condition on $\sfr$ provided we take $\eps>0$ to be small enough. The derivative of $\chi(\xi^2+|\eta|_{\slg^{-1}}^2)$ along $H_p$ has the opposite sign if the support of $\chi$ is sufficiently localized near $0$, and therefore we can propagate scattering decay from a punctured neighborhood of $\cR_{\rm out}$ into $\cR_{\rm out}$.

  To prove~\eqref{EqStEstscLower}, note first that the terms involving $Q\in\cA^{1+\delta}\Diffsc^1(X;E_X)$, $g^{0 0}\in\cA^{1+\delta}(X)$ in~\eqref{EqStEstFamily} do not contribute to $q_\sigma$. The scalar operator $2 i\sigma\rho(\rho\pa_\rho-\frac{n-1}{2})=-2\sigma(i\pa_r+\frac{i(n-1)}{2 r})$ (where $r=\rho^{-1}$) is symmetric with respect to the volume density $r^{n-1}|\dd r\,\dd\slg|$. Hence, by Lemma~\ref{LemmaGSGVol}, its imaginary part with respect to $|\dd g_X|$ and any fiber inner product on $E_X$ lies in $\cA^{1+\delta}(X;\End(E_X))$ and therefore does not contribute to $q_\sigma$ either. Note next that the scattering principal symbol of
  \begin{equation}
  \label{EqStEstImP0}
    \rho^{-1}\Im\hat P(0)\in\rho^{-1}(\rho^2\CI+\cA^{2+\delta})\Diffb^1(X;E_X)=(\rho\CI+\cA^{1+\delta})\Diffb^1(X;E_X)
  \end{equation}
  at $\cR_{\rm out}$ vanishes; indeed, in local coordinates, this operator can be written as a linear combination of $\rho^2 D_\rho$ and $\rho D_\omega$ (both of which have vanishing scattering principal symbol at the zero section) with coefficients in $\End(E_X)$, plus a subprincipal term in $(\rho\CI+\cA^{1+\delta})(X;\End(E_X))$. Finally, choose on $E_X$ a fiber inner product so that $\Re S(p)=\half(S(p)+S(p)^*)\geq\ubar S-\eps$ for all $p\in\pa X$ (see Proposition~\ref{PropLA}); then the only contribution to $q_\sigma$ comes from $\Im(-2 i\sigma\rho S)=-2\sigma\rho\Re S$, which gives~\eqref{EqStEstscLower}.

  Turning to the incoming radial set, one again uses a commutant with main term $\rho^{-2\sfr-1}$; one then has $\rho^{2\sfr}H_p(\rho^{-2\sfr-1})=-2\sigma(2\sfr+1)$ at $\cR_{\sigma,\rm in}$ by~\eqref{EqStEstHam}. For the computation of $(\rho^{-1}q_\sigma)_{\cR_{\sigma,\rm in}}$, the only terms of $\hat P(\sigma)$ that contribute are $-2 i\sigma\rho S$ and $\hat P(0)$. With respect to a fiber inner product $h$ on $E_X$, we thus have
  \begin{equation}
  \label{EqStEstscqSigma}
  \begin{split}
    \bigl(\rho^{-1}q_\sigma\bigr)\big|_{\cR_{\sigma,\rm in}} &= -2\sigma\Re S + \sigmasc^{1,0}\bigl(\rho^{-1}\Im\hat P(0)\bigr)\big|_{\cR_{\sigma,\rm in}} \\
      &= -2\sigma \Bigl( \Re S - \half\,\sigmasc^{1,0}\bigl(\rho^{-1}\Im \hat P(0)\bigr)\big|_{2\frac{\dd\rho}{\rho^2}} \Bigr).
  \end{split}
  \end{equation}
  (We use here that $\sigmasc^{1,0}(\rho^{-1}\Im\hat P(0))$ is linear in the fibers of $\Tsc^*_{\pa X}X$, as follows from~\eqref{EqStEstImP0} and the subsequent discussion.) By definition of $\ubar S_{\rm in}$, we can choose $h$ so that, as a self-adjoint bundle endomorphism on the pullback of $E|_{\pa X}$ to $\cR_{\sigma,\rm in}$, this is $\leq -2\sigma(\ubar S_{\rm in}-\eps)$. Altogether then, we obtain a negative commutator at $\cR_{\sigma,\rm in}$ if $-2\sigma(2\sfr+1+2(\ubar S_{\rm in}-\eps))<0$, which is true by assumption on $\sfr$ when $\eps>0$ is sufficiently small. The derivative of the cutoff term $\chi((\xi-2\sigma)^2+|\eta|_{\slg^{-1}}^2)^2$ along $\sfH_p$ has the same sign, and hence we obtain an estimate for the propagation of scattering decay out of $\cR_{\rm in,\sigma}$.

  Altogether, one obtains an estimate
  \begin{equation}
  \label{EqStEstNzscEstPf}
    \|u\|_{\Hsc^{\sfs,\sfr}} \leq C\bigl( \|\hat P(\sigma)u\|_{\Hsc^{\sfs-2,\sfr+1}} + \|u\|_{\Hsc^{-N,\sfr_0}}\bigr)
  \end{equation}
  for any fixed $N$, where $\sfr_0<\sfr$ still satisfies $\sfr_0>-\half-\ubar S_{\rm in}$ at $\cR_{\sigma,\rm in}$. (This estimate holds in the strong form that $u\in\Hsc^{\sfs,\sfr}(X;E_X)$ if the right hand side is finite.) Analogous arguments for $\hat P(\sigma)^*$ (defined with respect to any fiber inner product on $E_X$---note that adjoints with respect to two different choices differ merely via conjugation by a bundle isomorphism of $E_X$) give the dual estimate
  \begin{equation}
  \label{EqStEstNzscEstAdj}
    \|\tilde u\|_{\Hsc^{-\sfs+2,-\sfr-1}} \leq C\bigl( \|\hat P(\sigma)^*\tilde u\|_{\Hsc^{-\sfs,-\sfr}} + \|\tilde u\|_{\Hsc^{-N,\sfr_1}}\bigr),
  \end{equation}
  where $\sfr_1<-\sfr-1$ satisfies $\sfr_1>-\half-\ubar S$ at $\cR_{\rm out}$. Together, these two estimates imply that~\eqref{EqStEstNzsc} is Fredholm.

  Any element $u$ in the nullspace of~\eqref{EqStEstNzsc} automatically satisfies $u\in\Hsc^{\infty,\sfr}(X;E_X)$, and by the incoming radial point estimate we moreover have $u\in\Hsc^{\infty,\sfr'}(X;E_X)$ for all variable order functions $\sfr'$ satisfying $\sfr'<-\half+\ubar S$ at $\cR_{\rm out}$. Moreover, one can prove iterated regularity of $u$ under application of any number of ps.d.o.\ $\rho^{-1}A$ where the principal symbol of $A\in\Psisc^{1,0}(X;E_X)$ vanishes at $\cR_{\rm out}$; see e.g.\ \cite[\S2]{GellRedmanHassellShapiroZhangHelmholtz} for a detailed discussion of such module regularity (albeit in the non-conjugated perspective) which originated in \cite{HassellMelroseVasySymbolicOrderZero}. Since the set of such $A$ includes all scattering vector fields (acting on sections of $E_X$ using an arbitrary connection), this module regularity is the same as b-regularity. We thus have $u\in\Hb^{\infty,\ell}(X;E_X)\subset\cA^{\frac{n}{2}+\ell}(X;E_X)$ for any $\ell<-\half+\ubar S$. By assumption~\eqref{ItGSOSpec} in Definition~\ref{DefGSOSpec}, this implies that $u=0$, as desired.

  By a standard functional analytic argument, the error term $\|u\|_{\Hsc^{-N,\sfr_0}}$ in~\eqref{EqStEstNzscEstPf} can then be dropped upon increasing the constant $C$. The claim regarding~\eqref{EqStEstNzscEst} follows similarly from the fact that under the stated assumptions, the estimate~\eqref{EqStEstNzscEstPf} holds uniformly for $\sigma\in I$, and the error term can be dropped by essentially the same functional analytic argument (using now also that $I$ is compact).

  \pfstep{Intermezzo: variable order estimates in $\Im\sigma\geq 0$.} Before discussing the b-setting, we prove that~\eqref{EqStEstNzsc} is Fredholm (with trivial kernel) also for $\Im\sigma>0$, with $\sfr$ now arbitrary except for the requirement that $\sfr<-\half+\ubar S$ at (and constant near) $\cR_{\rm out}$. For $\Im\sigma>0$, the scattering principal symbol of $\hat P(\sigma)$ is elliptic except at $\cR_{\rm out}$. To get an estimate at $\cR_{\rm out}$, we again use a positive commutator argument involving the commutant~\eqref{EqStEstscCommA}, which we write as $a=\check a^2$ with $\check a=\rho^{-\sfr-\frac12}\chi(\xi^2+|\eta|_{\slg^{-1}}^2)\chi(\rho)$, and the operator $A=\check A^*\check A$ where $\check A\in\Psisc^{-\infty,\sfr+\frac12}(X;E_X)$ has principal symbol $\check a$. Since now $\sigmasc^{2,0}(\hat P(\sigma))$ is not real, the imaginary part $\Im\hat P(\sigma)$ gives the leading contribution (with one decay order less than the commutator term) to $\sC$ at $\pa X$. Since this leads to somewhat cumbersome expressions, we instead consider the operator
  \begin{equation}
  \label{EqStEstNzscOp}
    \sigma^{-1}\hat P(\sigma) = 2 i\rho\Bigl(\rho\pa_\rho-\frac{n-1}{2}-S\Bigr) + \sigma^{-1}\hat P(0) - i Q + \sigma g^{0 0}
  \end{equation}
  and evaluate the commutator
  \begin{equation}
  \label{EqStEstNzscComm2}
    2 \Im\la\sigma^{-1}\hat P(\sigma)u,A u\ra = \la\sC u,u\ra,\qquad \sC:=\Re\Bigl(\sigma^{-1}i[\hat P(\sigma),A] + 2\bigl(\Im(\sigma^{-1}\hat P(\sigma))\bigr)A\Bigr).
  \end{equation}
  The first summand of $\sC$ lies in $\Psisc^{-\infty,2\sfr}$, with $\rho^{2\sfr}$ times its principal symbol at $\cR_{\rm out}$ equal to $2(2\sfr+1)$. In the second summand, the terms in~\eqref{EqStEstNzscOp} involving $Q\in\cA^{1+\delta}\Diffsc^1$ and $g^{0 0}\in\cA^{1+\delta}$ contribute subprincipal terms ($\Psisc^{-\infty,2\sfr-\delta}$) to $\sC$ as before, and so does the term $2 i\rho(\rho\pa_\rho-\frac{n-1}{2})$, while $-2 i\rho S$ contributes a term with principal symbol $\leq -4(\ubar S_{\rm in}-\eps)$ at $\cR_{\rm out}$ (omitting the factor of $\rho^{-2\sfr}$) for any fixed $\eps>0$ upon choosing an appropriate fiber inner product. (The terms discussed thus far sum up to be a negative multiple of $\rho^{-2\sfr}$ at $\cR_{\rm out}$ in light of the assumed upper bound on $\sfr$ there.) Lastly, expanding
  \[
    \Im(\sigma^{-1}\hat P(0))=|\sigma|^{-2}\bigl(-\Im\sigma \Re\hat P(0)+\Re\sigma \Im\hat P(0)\bigr),
  \]
  the term involving $\Im\hat P(0)$ contributes a term in $\Psisc^{-\infty,2\sfr}$ to $\sC$ but with vanishing principal symbol at $\cR_{\rm out}$, as before; it remains to write the contribution of the other term to $\sC$ as
  \begin{equation}
  \label{EqStEstNzscP0}
    -(\Im\sigma)|\sigma|^{-2}\Re\Bigl(\bigl(\hat P(0)+\hat P(0)^*\bigr)\check A^2\Bigr).
  \end{equation}
  But since by~\eqref{EqGSGMetricCoeff} the b-principal symbol of $\rho^{-2}\hat P(0)$ is a positive definite quadratic form near $\pa X$, we can write (cf.\ \cite[Lemma~3.3]{VasyLAPLag})\footnote{The coefficients of $\nabla$ and $R$ lie in $\CI(X)+\cA^\delta(X)$. We shall not explicitly write this anymore, unless the presence of conormal coefficients requires additional care beyond the smooth coefficient case.}
  \[
    \hat P(0) = (\rho\nabla)^*(\rho\nabla) + \rho^2 R,\qquad R\in \Diffb^1(X;E_X),
  \]
  where $\nabla\in\Diffb^1(X;E_X,\Tb^*X\otimes E_X)$, and $\nabla^*$ is its adjoint with respect to the (already chosen) fiber inner product on $E_X$ and the dual b-metric $\rho^{-2}g^{X X}$ on $X$. Thus $\rho^2 R\in\rho\Diffsc^1$ contributes to~\eqref{EqStEstNzscP0} a term in $\Psisc^{-\infty,2\sfr}$ whose scattering principal symbol vanishes at $\cR_{\rm out}$ (since that of $\rho R\in\rho\Diffb^1$ does). Finally, setting $P_0=(\rho\nabla)^*(\rho\nabla)=P_0^*\in\rho^2\Diffb^2(X;E_X)$, we rewrite the contribution $-(P_0+P_0^*)\check A^2$ (which lies in $\Psisc^{-\infty,2\sfr+1}$, i.e.\ gives the leading order contribution to $\sC$) as $(\Im\sigma)|\sigma|^{-2}$ times
  \begin{align*}
    -\Re\bigl(P_0\check A^2\bigr) &= -\check A P_0\check A - \Re\bigl([P_0,\check A]\check A\bigr) \\
      &= -(\rho\nabla\circ\check A)^*(\rho\nabla\circ\check A) - \half\bigl([P_0,\check A]\check A+\check A[\check A,P_0]\bigr) \\
      &= -(\rho\nabla\circ\check A)^*(\rho\nabla\circ\check A) - \half[[P_0,\check A],\check A].
  \end{align*}
  (This is where the explicit insertion of the real part in~\eqref{EqStEstNzscComm2} is useful.) The first term is $\leq 0$ \emph{as an operator}, i.e.\ when acting on $u$ followed by pairing with $u$ (as in~\eqref{EqStEstNzscComm2}), and hence has the same sign as the main term in the above symbolic computations; it can thus be dropped. The second, double commutator, term on the other hand is an element of $\Psisc^{-\infty,2\sfr-1}$, and therefore subprincipal in the symbolic commutator calculation. See \cite[\S2.5]{VasyMicroKerrdS} for similar arguments in the closed manifold setting, and \cite{NonnenmacherZworskiQuantumDecay,WunschZworskiNormHypResolvent} for further background on \emph{complex absorbing potentials} which $i \Im(\sigma^{-1}\hat P(\sigma))$ can be thought of.

  Altogether, we again obtain the estimate~\eqref{EqStEstNzscEstPf} with arbitrary $\sfr_0$. Together with a matching dual estimate~\eqref{EqStEstNzscEstAdj}, this gives the Fredholm property of~\eqref{EqStEstNzsc}, and by the same arguments as in the real $\sigma$ setting the conormality of elements of its kernel---which thus is trivial due to the spectral admissibility of $\hat P(\sigma)$.

  We also note that a slight extension of the above arguments gives a uniform estimate~\eqref{EqStEstNzscEst} when $\sigma\in\Omega\subset\C$, with $\Omega$ a compact subset of the punctured upper half plane as in part~\eqref{ItStEstNzb}, if in the case $\Omega\cap\R\neq\emptyset$ the same assumptions on $\sfr$ as in part~\eqref{ItStEstNzsc} of the Lemma are satisfied for all $\sigma\in\Omega\cap\R$. Only the case $\Omega\cap\R\neq\emptyset$ requires an elaboration. The above arguments at the outgoing radial set apply uniformly down to $\Omega\cap\R$: apart from the symbolic positive commutator computation, this uses that the term~\eqref{EqStEstNzscP0} is nonpositive as an operator, modulo terms in $\Psisc^{-\infty,2\sfr}$ whose principal symbols vanish at $\cR_{\rm out}$. Moreover, since $\hat P(\sigma)$ is not uniformly elliptic on $\Char_{\Re\sigma}$, one now needs to keep track also of the terms in the symbolic commutator computation which involve derivatives of the cutoffs which localize to $\cR_{\rm out}$: rather than being controlled by elliptic regularity, they give rise to a priori control terms in a punctured neighborhood of $\cR_{\rm out}$ as in the case of real $\sigma$. Uniform (down to $\Omega\cap\R$) versions of the incoming radial point estimate and real principal type propagation estimates, discussed previously for real $\sigma$, can be proved in an analogous manner: the imaginary part of $\sigma^{-1}\hat P(\sigma)$ acts as complex absorption (with the correct sign for propagation from the incoming to the outgoing radial set) just as in the outgoing radial point estimate. Together with elliptic regularity outside of $\Char_{\Re\sigma}$, this completes the proof of~\eqref{EqStEstNzscEst} uniformly for $\sigma\in\Omega$.

  \pfstep{Part~\eqref{ItStEstNzb}.} Second microlocal refinements required for the proof of~\eqref{EqStEstNzb} are proved in \cite{VasyLAPLag} using the second microlocal spaces $H_{\scop,\bop}^{s,r,\ell}$ introduced in~\cite{VasyLowEnergyLag,VasyLAPLag}; these spaces make precise the notion of a scattering Sobolev space whose scattering decay order is equal to $r$ (constant), except right at the zero section where it jumps to $\ell$ (which is then regarded as the b-decay order). One can then prove the Fredholm property of
  \[
    \hat P(\sigma) \colon \bigl\{ u\in H_{\scop,\bop}^{s,r,\ell}(X;E_X) \colon \hat P(\sigma)u\in H_{\scop,\bop}^{s-2,r+1,\ell+1}(X;E_X) \bigr\} \to H_{\scop,\bop}^{s-2,r+1,\ell+1}(X;E_X)
  \]
  for $\sigma\neq 0$, $\Im\sigma\geq 0$ using radial point and propagation estimates and complex absorption type arguments as before; this requires $r>-\half-\ubar S_{\rm in}$, resp.\ $\ell<-\half+\ubar S$ in order for the incoming, resp.\ outgoing radial point estimate to work. The only additional ingredient is now an estimate for the b-normal operator of $\hat P(\sigma)$, given by the ordinary differential operator $2 i\sigma\rho(\rho\pa_\rho-\frac{n-1}{2}-S|_{\pa X})$, between b-Sobolev spaces on $X$ with weights $\ell$ and $\ell+1$ (see \cite[Lemma~4.13]{VasyLAPLag}). The injectivity of $\hat P(\sigma)$ in this second microlocal setting is particularly easy to show, since for $u\in H_{\scop,\bop}^{s,r,\ell}$ with $\hat P(\sigma)u=0$, one gets $u\in H_{\scop,\bop}^{\infty,\infty,\ell}=\Hb^{\infty,\ell}\subset\cA^{\frac{n}{2}+\ell}$ directly from the incoming radial point estimate and propagation all the way down to the resolved zero section.

  The passage to b-Sobolev spaces is accomplished by noting that for $r=s+\ell$ we have $\Hb^{s,\ell}=H_{\bop,\scop}^{s,r,\ell}$ and $H_{\bop,\scop}^{s-2,r+1,\ell+1}\subset H_{\bop,\scop}^{s,r+1,\ell+1}=\Hb^{s,\ell+1}$. This proves~\eqref{EqStEstNzb}--\eqref{EqStEstNzbEst} and completes the proof of the Proposition.
\end{proof}

\subsubsection{High frequencies}
\label{SssStEstHi}

In the high energy regime, where $\pm\Re\sigma\to\infty$ while $\Im\sigma\in[0,\infty)$ remains bounded, we will exploit the nontrapping assumption of Definition~\ref{DefGSG}\eqref{ItGSGNontrap}. More generally, we need to study the regime $|\sigma|\to\infty$ in $\Im\sigma\geq 0$. We thus pass to the semiclassical rescaling
\begin{align*}
  P_{h,z} &:= h^2 \hat P(h^{-1}z) = 2 i z\rho\Bigl(h\rho\pa_\rho-h\frac{n-1}{2}-h S\Bigr) + h^2\hat P(0) - i z h Q + z^2 g^{0 0}, \\
  &\qquad h:=|\sigma|^{-1},\quad z:=\frac{\sigma}{|\sigma|}.
\end{align*}
This is a semiclassical scattering operator of class $\Diffsch^{2,0,0}(X;E_X)$, with principal symbol
\begin{equation}
\label{EqStEstSymbh}
  p_z(x,\mu_\semi) := \sigmasch^{2,0,0}(P_{h,z})(x,\mu_\semi) = h^2 G(-h^{-1}z\,\dd t_*+\mu_\semi)=G(-z\,\dd t_*+h\mu_\semi)
\end{equation}
where $G(\zeta)=g^{-1}(\zeta,\zeta)$ is the dual metric function and $\mu_\semi\in\Tsch^*_x X$. If we write semiclassical scattering covectors near $\pa X$ as
\begin{equation}
\label{EqStEstTschCoord}
  \mu_\semi = \xi_\semi\frac{\dd\rho}{h\rho^2} + \frac{\eta_\semi}{h\rho},\qquad \eta_\semi\in T^*\Sph^{n-1},
\end{equation}
then Definition~\ref{DefGSG} gives
\begin{equation}
\label{EqStEstpz}
  p_z \equiv -2 z\xi_\semi+\xi_\semi^2+|\eta_\semi|_{\slg^{-1}}^2 \bmod \cA^\delta S^2(\Tsch^*X).
\end{equation}
The characteristic set over $h=0$ is
\[
  \Char_{\semi,z} = p_z^{-1}(0) \subset \ol{\Tsch^*_{h^{-1}(0)}}X,
\]
and the semiclassical versions of the radial sets over $\pa X$ are the sets
\begin{align*}
  \cR_{\semi,\pm 1,\rm in} &= \{ (\omega;\xi_\semi,\eta_\semi) \in \Tsch^*_{\pa X}X \colon \xi_\semi=\pm 2,\ \eta_\semi=0 \}, \\
  \cR_{\semi,\rm out} &= \{ (\omega;0,0) \},
\end{align*}
with the incoming radial set defined only for $z=1+\cO(h)$. Finally, we note that the Hamiltonian vector field at a point~\eqref{EqStEstTschCoord} is
\begin{align}
\label{EqStEstHamh}
  \sfH_{p,z} := h^{-1}\rho^{-1}H_{p_z} &= \rho^{-1}\pi_*\Bigl(H_G|_{-z\,\dd t_*+\xi_\semi\frac{\dd\rho}{\rho^2}+\frac{\eta_\semi}{\rho}}\Bigr) \\
\label{EqStEstHamhRad}
    &\equiv 2(\xi_\semi-z)(\rho\pa_\rho+\eta_\semi\pa_{\eta_\semi}) - 2|\eta_\semi|_{\slg^{-1}}^2\pa_{\xi_\semi} + H_{|\eta_\semi|^2_{\slg^{-1}}} \bmod \cA^\delta\Vb(\Tsch^*X);
\end{align}
in~\eqref{EqStEstHamh}, the map $\pi$ is a rescaling of the map $T^*M_0\to T^*X$ induced by the inclusion $X\ni x\mapsto(0,x)\in M_0$, in that it maps $\xi\frac{\dd\rho}{\rho^2}+\frac{\eta}{\rho}\mapsto\xi\frac{\dd\rho}{h\rho^2}+\frac{\eta}{h\rho}$ in the fibers.

Solely using the dynamical assumptions on the metric $g$ and the structure of the operator $P$, but not requiring mode stability (Definition~\ref{DefGSOSpec}\eqref{ItGSOSpec}), we now show:

\begin{prop}[High energy estimates]
\label{PropStEstHi}
  Let $\ubar S$, $\ubar S_{\rm in}$ be as in Definition~\usref{DefStEstThr}.
  \begin{enumerate}
  \item\label{ItStEstHisc}{\rm (Variable order estimates.)} Let $\sfs\in\CI({}^\schop S^*X)$, $\sfr\in\CI(\ol{\Tsch^*_{\pa X}}X)$, and $\sfb\in\CI(\ol{\Tsch^*_{h^{-1}(0)}}X)$. Suppose that $\sfs$, $\sfr$, and $\sfb$ are nonincreasing along the flow of $\pm\sfH_{p_{\pm 1}}$ inside $\Char_{\semi,\pm 1}$, and that moreover $\sfr>-\half-\ubar S_{\rm in}$ at $\cR_{\semi,\pm 1,\rm in}$ and $\sfr<-\half+\ubar S$ at $\cR_{\semi,\rm out}$, with $\sfr$ constant near $\cR_{\semi,\pm 1,\rm in}$ and $\cR_{\semi,\rm out}$. Then there exists $C'>0$ so that for $\sigma\in\R$ with $\pm\sigma\geq C'$, the operator $\hat P(\sigma)$ in~\eqref{EqStEstNzsc} is invertible. Moreover, there exists a constant $C''$ so that
    \begin{equation}
    \label{EqStEstHisc}
      \|u\|_{H_{\scop,|\sigma|^{-1}}^{\sfs,\sfr,\sfb}(X;E_X)} \leq C''\| \hat P(\sigma)u \|_{H_{\scop,|\sigma|^{-1}}^{\sfs-2,\sfr+1,\sfb-1}},\qquad \sigma\in\R,\ \pm\sigma\geq C'.
    \end{equation}
  \item\label{ItStEstHib}{\rm (Estimates on b-spaces.)} Let $s,\ell\in\R$ be as in Proposition~\usref{PropStEstNz}\eqref{ItStEstNzb}. There exists $C'>0$ so that for $\sigma\in\C$ with $\Im\sigma\geq 0$ and $|\sigma|\geq C$, the operator $\hat P(\sigma)$ in~\eqref{EqStEstNzb} is invertible. Moreover, there exists a constant $C$ so that
    \begin{equation}
    \label{EqStEstHib}
      \|u\|_{H_{\bop,|\sigma|^{-1}}^{s,\ell}(X;E_X)} \leq C|\sigma|^{-1} \| \hat P(\sigma)u \|_{H_{\bop,|\sigma|^{-1}}^{s,\ell+1}(X;E_X)},\qquad \Im\sigma\geq 0,\ |\sigma|\geq C'.
    \end{equation}
  \end{enumerate}
\end{prop}
\begin{proof}
  \pfstep{Part~\eqref{ItStEstHisc}.} We only consider the `$+$' sign, so $z=1$. Since $\pa_{t_*}$ is timelike over $X^\circ$, the principal symbol $p_z$ is elliptic for large frequencies, i.e.\ near fiber infinity of $\ol{\Tsch^*}X$. By the nontrapping assumption on $g$, the expression~\eqref{EqStEstHamh} for $\sfH_{p,z}$ (and noting that $-\dd t_*$ is future timelike) implies that every maximally extended integral curve $\gamma\colon I\to\Tsch^*_{h^{-1}(0)}X^\circ$ of $\sfH_{p,z}$ in $\Char_z$ remains in any fixed neighborhood $\rho<\eps$ for arguments sufficiently close to $\inf I$ and $\sup I$. By the source, resp.\ sink nature of $\cR_{\semi,1,\rm in}$, resp.\ $\cR_{\semi,\rm out}$, this implies that $\gamma(s)\to\cR_{\semi,1,\rm in}$, resp.\ $\gamma(s)\to\cR_{\semi,\rm out}$ as $s\searrow\inf I$, resp.\ $s\nearrow\sup I$. (We use here that integral curves of $\sfH_{p,z}$ in the characteristic set over $\pa X$ have the same property by direct computation, unless they are contained in one of the radial sets and thus constant.) Along $\gamma$, we thus have real principal type estimates, possibly with monotonically decreasing semiclassical order $\sfb$. (See \cite[Appendix~A]{HintzVasyCauchyHorizon} for the case of semiclassical orders which only depend on the base point---though the propagation result, \cite[Proposition~A.5]{HintzVasyCauchyHorizon}, does not require this restriction; see also \cite[Appendix~A]{BaskinVasyWunschRadMink} for the homogeneous setting. See also~\cite[\S2.3]{GalkowskiThinBarriers}.)

  Near the radial sets, the positive commutator arguments used in the proof of Proposition~\ref{PropStEstNz} apply here as well upon switching to semiclassical fiber-linear coordinates and quantizations, except now also the derivative falling on the cutoff $\chi(\rho)$ localizing to an $\eps$-neighborhood of $\rho=0$ contributes to the principal symbol of the commutator (called $\sC$ there); for propagation near $\cR_{\semi,1,\rm in}$, resp.\ $\cR_{\semi,\rm out}$, this contribution has the same sign as, resp.\ the opposite sign of the main term (which is negative) at the respective radial set, as follows from the positivity, resp.\ negativity (up to $\cO(\rho^\delta)$ errors) of the $\rho\pa_\rho$ coefficient in~\eqref{EqStEstHamhRad}. Altogether, one thus obtains estimates for $P_{h,z}$ and its adjoint,
  \begin{align*}
    \|u\|_{H_{\scop,h}^{\sfs,\sfr,\sfb}} \leq C h^{-1}\Bigl( \|P_{h,z}u\|_{H_{\scop,h}^{\sfs-2,\sfr+1,\sfb}} + \|u\|_{H_{\scop,h}^{-N,-N,-N}} \Bigr), \\
    \|\tilde u\|_{H_{\scop,h}^{-\sfs+2,-\sfr-1,-\sfb}} \leq C h^{-1}\Bigl( \|P_{h,z}^*\tilde u\|_{H_{\scop,h}^{-\sfs,-\sfr,-\sfb}} + \|\tilde u\|_{H_{\scop,h}^{-N,-N,-N}} \Bigr),
  \end{align*}
  for any fixed $N$, which hold for all $u$ and $\tilde u$ for which all norms are finite. For sufficiently small $h>0$ (i.e.\ for sufficiently large $|\Re\sigma|$), the second terms on the right are less than $\half$ times the left hand side and can thus be absorbed. Since $h^{-1}\|P_{h,z}u\|_{H_{\scop,h}^{\sfs-2,\sfr+1,\sfb}}=\|h^2\hat P(h^{-1}z)u\|_{H_{\scop,h}^{\sfs-2,\sfr+1,\sfb+1}}=\|\hat P(h^{-1}z)u\|_{H_{\scop,h}^{\sfs-2,\sfr+1,\sfb-1}}$, this gives~\eqref{EqStEstHisc}.

  \pfstep{Intermezzo: estimates in the closed upper half plane.} The estimate~\eqref{EqStEstHisc} holds more generally for $\pm\Re\sigma\geq C'$ when $\sigma\in\C$, $\Im\sigma\in[0,C]$ for any fixed $C$, with $C',C''$ depending on $C$; this is most easily proved by considering the operator $z^{-1}P_{h,z}$, analogously to the arguments starting with~\eqref{EqStEstNzscOp}. If $\sigma$ is nonreal, the imaginary part of $P_{h,z}$ contributes (via $z^{-1}h^2\hat P(0)$, with $\Im(z^{-1})=-\Im z$) terms which are one order stronger (in the scattering decay order sense) at the radial sets than the main terms in the real $\sigma$ case; but by a simple modification of the arguments after~\eqref{EqStEstNzscComm2}, these terms have the correct sign as operators modulo terms which can be controlled by the symbolic commutator calculation.

  More generally, for $\sigma\in\C$, $|\sigma|>1$, with $\Im\sigma\geq 0$ now unbounded, proofs of semiclassical estimates for $z^{-1}P_{h,z}$ can be carried out along similar lines, with some modifications which we proceed to explain. First, regarding outgoing radial point estimates, we note the following: while $z^{-1}\sfH_{p,z}$ is no longer a real vector field, its main term $-2\rho\pa_\rho$ yields the same contribution to the principal symbol of the commutator
  \[
    i\bigl( (z^{-1}P_{h,z})^* A - A z^{-1}P_{h,z} \bigr)
  \]
  (called $\sC$ above) upon differentiating the weight $\rho^{-2\sfr-1}$ of the commutant as before, whereas all other contributions can be made arbitrarily small at $\cR_{\semi,\rm out}$ upon localizing to a sufficiently small neighborhood thereof. Moreover, the skew-adjoint part of $z^{-1}\hat P(0)$ contributes terms which are of higher order than the main term of the symbolic computation not only in the scattering decay, but also in the semiclassical order sense when $\Im z$ is not of size $\cO(h)$; these terms are however still nonpositive as operators, modulo terms that can be absorbed in the symbolic commutator calculation near $\cR_{\semi,\rm out}$.

  Next, the incoming radial set only plays a role when $0\leq\Im z\leq\eta$ where $\eta>0$ is arbitrary. The deviation of $z^{-1}\sfH_{p,z}$ from $\pm\sfH_{p,\pm 1}$ is of size $\cO(\eta)$ near $\cR_{\semi,\pm 1,\rm in}$; note moreover that the skew-adjoint part of $z^{-1}\hat P(0)$ is the sum of a term which is nonpositive as an operator and the term $(\pm 1+\cO(\eta))\Im\hat P(0)$ which contributes to the main symbolic term as in~\eqref{EqStEstscqSigma}.

  Regarding the replacement for real principal type propagation estimates, we work with $P_{h,z}$ and consider for suitable commutants $A=\check A^*\check A$ with $\check A=\check A^*$ (which are semiclassical scattering operators) the commutator
  \begin{align*}
    i\bigl( P_{h,z}^* A - A P_{h,z} \bigr) &= i[\Re P_{h,z},A] + (\Im P_{h,z})A + A(\Im P_{h,z}) \\
      &= i[\Re P_{h,z},A] + 2\check A(\Im P_{h,z})\check A + [\check A,[\check A,\Im P_{h,z}]]
  \end{align*}
  Following \cite[\S\S3.2 and 7.2]{VasyMicroKerrdS}, we claim that for $0\leq\Im z\leq\eta\ll 1$ and $z=1+\cO(\eta)$ (the case $z=-1+\cO(\eta)$ being analogous), and on the characteristic set of $\Re P_{h,z}$, the principal symbol of $\Im P_{h,z}$ is nonpositive, which ensures that the second term on the right here has a nonpositive principal symbol as required for forward propagation along the Hamiltonian vector field of $\Re P_{h,z}$. To verify the claim, recall~\eqref{EqStEstSymbh} and note that (switching to standard scattering covectors for notational simplicity) for $\mu\in\Tsc^*_x X$ we have
  \begin{equation}
  \label{EqStEstHiscGImz}
  \begin{split}
    G(-z\,\dd t_*+\mu) &= \Bigl(-(\Im z)^2 G(\dd t_*) + G(-(\Re z)\dd t_*+\mu)\Bigr) \\
      &\quad\qquad + 2 i(\Im z)g^{-1}(-(\Re z)\dd t_*+\mu,-\dd t_*).
  \end{split}
  \end{equation}
  If $\Re G(-z\,\dd t_*+\mu)=0$, we therefore have $G(-(\Re z)\dd t_*+\mu)=(\Im z)^2 G(\dd t_*)\leq 0$ since $\dd t_*$ is timelike. Therefore, the covector $-(\Re z)\dd t_*+\mu$ is causal, and hence $g^{-1}(-(\Re z)\dd t_*+\mu,-\dd t_*)\neq 0$. But by assumption~\eqref{ItGSGTime2} in Definition~\ref{DefGSG}, $\Tsc^*_x X$ is spacelike, and therefore the set of $\mu\in\Tsc^*_x X$ for which $-(\Re z)\dd t_*+\mu$ is causal is convex, and it contains $\mu=0$; therefore $g^{-1}(-(\Re z)\dd t_*+\mu,-\dd t_*)$ has the same sign as $g^{-1}(-(\Re z)\dd t_*,-\dd t_*)<0$, as desired. For later use, we note that these arguments show that $-(\Re z)\dd t_*+\mu$ is future causal when $\Re G(-z\,\dd t_*+\mu)=0$, $\Im z\geq 0$, and $\Re z>0$, and indeed future timelike when $\Im z>0$.

  Lastly, consider the characteristic set $\Char_{\semi,z}$ of $P_{h,z}$ when $\Im z\geq\eta>0$. Over $\pa X$, this is equal to $\cR_{\semi,\rm out}$ by~\eqref{EqStEstpz}. Over $X^\circ$ on the other hand, $\Re p_z=\Re G(-z\,\dd t_*+\mu)=0$ implies that $-(\Re z)\dd t_*+\mu$ is timelike, and thus $\Im G(-z\,\dd t_*+\mu)\neq 0$. This proves the absence of characteristic set over $X^\circ$ at finite semiclassical frequencies; at fiber infinity, there is no characteristic set since $\Tsc^*X$ is spacelike. Thus, $\Char_{\semi,z}=\cR_{\semi,\rm out}$ for $\Im z\geq\eta>0$.

  \pfstep{Part~\eqref{ItStEstHib}.} For bounded $\Im\sigma$, this is proved in \cite[Theorem~1.5]{VasyLAPLag}. In short, semiclassical versions of the radial point and propagation estimates now give estimates on semiclassical second microlocal b/scattering Sobolev spaces,
  \[
    \|u\|_{H_{\scop,\bop,h}^{s,r,\ell}} \leq C h^{-1}\|P_{h,z}u\|_{H_{\scop,\bop,h}^{s-2,r+1,\ell+1}}
  \]
  (where we already absorbed the $\cO(h^N)$ error term into the left hand side, and with the semiclassical order taken to be $0$ in both norms), and a corresponding estimate on dual spaces. As in the reference, this implies the estimate~\eqref{EqStEstHib} on semiclassical b-Sobolev spaces when $\Im\sigma\in[0,C'']$ (for any fixed $C''$) and $|\Re\sigma|$ is sufficiently large. The estimates in $\Im\sigma\geq 0$, $|\sigma|\gg 1$, follow similarly to the variable order case.
\end{proof}

\begin{cor}[Invertibility for nonzero frequencies]
\label{CorStEstInv}
  For $0\neq\sigma\in\C$ with $\Im\sigma\geq 0$, the operator $\hat P(\sigma)$ is invertible as a map~\eqref{EqStEstNzsc} on variable order scattering Sobolev spaces, and also as a map~\eqref{EqStEstNzb} on b-Sobolev spaces.
\end{cor}
\begin{proof}
  Consider the case of variable scattering orders. Proposition~\ref{PropStEstHi} implies that $\hat P(\sigma)$, as a map~\eqref{EqStEstNzsc}, is invertible for sufficiently large $|\Re\sigma|$. Let us write $\cX_\sigma^{\sfs,\sfr}$ for the domain of the map~\eqref{EqStEstNzsc}. Since $\cX_\sigma^{\sfs,\sfr}$ depends on $\sigma$, proving the independence of the Fredholm index of $\hat P(\sigma)\colon\cX_\sigma^{\sfs,\sfr}\to\Hsc^{\sfs-2,\sfr+1}(X;E_X)$ requires an argument; in the case at hand, we shall prove that the set
  \[
    \cA:=\bigl\{\sigma\in\R\setminus\{0\}\colon \ker_{\Hsc^{-\sfs+2,-\sfr-1}(X;E_X)}\hat P(\sigma)^*=\{0\} \bigr\}
  \]
  (where $\sfr$, depending on $\sigma$, satisfies the conditions of Proposition~\ref{PropStEstNz}\eqref{ItStEstNzsc}, with the kernel however being independent of the particular choice of $\sfr$ subject to this condition) is open and closed in $\R\setminus\{0\}$. To prove the openness, we consider $\sigma_0\in\cA$ and fix a variable order function $\sfr$ satisfying the hypotheses of Proposition~\ref{PropStEstNz}\eqref{ItStEstNzsc} for all $\sigma\in[\sigma_0-\eps,\sigma_0+\eps]$ where we fix $\eps<|\sigma_0|$. We then exploit that we have a uniform estimate~\eqref{EqStEstNzscEstAdj} for $\hat P(\sigma)^*$: if there existed a sequence $\sigma_j\in[\sigma_0-\eps,\sigma_0+\eps]$ with $\sigma_j\to\sigma_0$, and $v_j\in\Hsc^{-\sfs+2,-\sfr-1}(X;E_X)$ with norm $1$ and so that $\hat P(\sigma_j)v_j=0$, then~\eqref{EqStEstNzscEstAdj} would give a lower bound $\|v_j\|_{\Hsc^{-N,\sfr_1}(X;E_X)}\geq C^{-1}$; passing to a subsequence of the $v_j$ which converges weakly in $\Hsc^{-\sfs+2,-\sfr-1}$ and thus strongly in $\Hsc^{-N,\sfr_1}$ to a necessarily nonzero limit $v_0\in\Hsc^{-\sfs+2,-\sfr-1}(X;E_X)$, the limit $v_0$ (by virtue merely of being the distributional limit of the $v_j$) satisfies $\hat P(\sigma_0)^*v_0=0$. This contradicts $\sigma_0\in\cA$.

  To prove that $\cA$ is closed, consider a point $\sigma_0\in\R\setminus\{0\}$ in the boundary of $\cA$. Assuming that $\sigma_0\notin\cA$, then setting $d=\dim\ker\hat P(\sigma_0)^*\geq 1$, we could define a map
  \[
    \tilde P(\sigma) \colon \cX_\sigma^{\sfs,\sfr} \oplus \C^d \ni (u,a) \mapsto \hat P(\sigma)u + \sum_{j=1}^d a_j f_j \in \Hsc^{\sfs-2,\sfr+1}(X;E_X),
  \]
  where $f_1,\ldots,f_d\in\Hsc^{\sfs-2,\sfr+1}(X;E_X)$ spans a complementary subspace to $\hat P(\sigma_0)(\cX_{\sigma_0}^{\sfs,\sfr})$. By definition, $\tilde P(\sigma_0)$ is invertible. Since also $\tilde P(\sigma)$ satisfies uniform (for $\sigma$ near $\sigma_0$) Fredholm estimates similar to~\eqref{EqStEstNzscEstPf}--\eqref{EqStEstNzscEstAdj}, repeating the above arguments shows that $\tilde P(\sigma)$ is invertible also for $\sigma$ sufficiently close to $\sigma_0$. Taking such a $\sigma$ which moreover lies in $\cA$, we have $f_1\in\hat P(\sigma)(\cX_\sigma^{\sfs,\sfr})$; but this contradicts the injectivity of $\tilde P(\sigma)$.

  In summary, we conclude that $\cA=\R\setminus\{0\}$ since $\cA$ contains both positive and negative numbers, as noted at the beginning of the proof. The proof in the case of b-Sobolev spaces is completely analogous.
\end{proof}

\begin{cor}[Continuity down to the real line]
\label{CorStEstC0}
  For $s,\ell\in\R$ with $\ell<-\half+\ubar S$ and $s+\ell>-\half-\ubar S_{\rm in}$, the map $\hat P(\sigma)^{-1}\colon\Hb^{s,\ell+1}(X;E_X)\to\Hb^{s,\ell}(X;E_X)$ is continuous in $\sigma\in\C$, $\Im\sigma\geq 0$, $\sigma\neq 0$, in the weak operator topology. Moreover, for $\eps>0$. The map $\hat P(\sigma)^{-1}\colon\Hb^{s,\ell+1}(X;E_X)\to\Hb^{s-\eps,\ell-\eps}(X;E_X)$ is continuous in the norm topology.
\end{cor}
\begin{proof}
  Given the uniform estimate~\eqref{EqStEstNzb}, this is a standard functional analytic argument, see e.g.\ \cite[\S2.7]{VasyMicroKerrdS}: consider a sequence $\sigma_j$ in the punctured upper half plane which converges to $\sigma_0\neq 0$. Given $f\in\Hb^{s,\ell+1}(X;E_X)$, the sequence $u_j=\hat P(\sigma_j)^{-1}f\in\Hb^{s,\ell}(X;E_X)$ is uniformly bounded, and hence upon passing to a subsequence converges weakly, $u_j\weakto u_0\in\Hb^{s,\ell}$. Therefore, $\hat P(\sigma_0)u_0=f$, and thus $u_0=\hat P(\sigma_0)^{-1}f$ due to the invertibility of $\hat P(\sigma_0)$. This proves that $\hat P(\sigma_0)^{-1}f$ is the weak limit of $\hat P(\sigma_j)^{-1}f$, as desired.

  For the final claim, assume the contrary. Then there exist $\eps,\delta>0$, a sequence $\sigma_j$ in the punctured upper half plane with $\sigma_0=\lim\sigma_j$, $\sigma_0\neq 0$, and $f_j\in\Hb^{s,\ell+1}$ with norm $1$ so that for the uniformly bounded sequence $\hat P(\sigma_j)^{-1}f_j\in\Hb^{s,\ell}$ we have
  \begin{equation}
  \label{EqStEstC0Pf}
    \|\hat P(\sigma_j)^{-1}f_j - \hat P(\sigma_0)^{-1}f_j\|_{\Hb^{s-\eps,\ell-\eps}} \geq \delta.
  \end{equation}
  But upon passing to a subsequence, we may assume that $u_j:=\hat P(\sigma_j)^{-1}f_j$ converges weakly to some $u_0\in\Hb^{s,\ell}$, and thus $u_j\to u_0$ in $\Hb^{s-\eps,\ell-\eps}$. We may also assume that $f_j\weakto f_0\in\Hb^{s,\ell+1}$. But then $f_j=\hat P(\sigma_j)u_j\to\hat P(\sigma_0)u_0$ strongly in $\Hb^{s-2-\eps,\ell-\eps}$, and therefore $f_0=\hat P(\sigma_0)u_0$ and thus $u_0=\hat P(\sigma_0)^{-1}f_0$. Therefore, writing
  \[
    \|\hat P(\sigma_j)^{-1}f_j-\hat P(\sigma_0)^{-1}f_j\|_{\Hb^{s-\eps,\ell-\eps}} \leq \|u_j-u_0\|_{\Hb^{s-\eps,\ell-\eps}} + \|\hat P(\sigma_0)^{-1}(f_0-f_j)\|_{\Hb^{s-\eps,\ell-\eps}},
  \]
  the first summand on the right is $<\delta/2$ for sufficiently large $j$, and likewise for the second summand (when $\eps$ is chosen so that $s+\ell-2\eps>-\half-\ubar S_{\rm in}$ still) since $f_j\to f_0$ in $\Hb^{s-\eps,\ell+1-\eps}$. This contradicts~\eqref{EqStEstC0Pf} and thus finishes the proof.
\end{proof}

\subsubsection{Uniform low frequency estimates} 
\label{SssStEstLo}

While we now know that $\hat P(\sigma)$ is invertible for all $\sigma$ with $\Im\sigma\geq 0$, it remains to prove \emph{uniform} estimates at low frequencies. We follow \cite{VasyLowEnergyLag} for the b-estimates, and \cite[\S3.5]{HintzKdSMS} (where the scalar wave operator on subextremal Kerr spacetimes is considered) for the variable order scattering estimates. Recall from~\eqref{EqGSONtf} that if we write $\hat P(0)=P_{(0)}(\rho,\omega,\rho D_\rho,D_\omega)$, then
\[
  N_\tface^\pm(P) = \pm 2 i\hat\rho\Bigl(\hat\rho\pa_{\hat\rho}-\frac{n-1}{2}-S|_{\pa X}\Bigr) + \hat\rho^2 P_{(0)}(0,\omega,\hat\rho D_{\hat\rho},D_\omega).
\]
As a scattering-b-operator on $\tface=[0,\infty]_{\hat\rho}\times\Sph^{n-1}$ (with weight $2$ at the b-end $\zface=\hat\rho^{-1}(\infty)$), this is elliptic in the differential order sense, but has a non-trivial scattering characteristic set over the scattering end $\sctface=\hat\rho^{-1}(0)$,
\[
  \Char_{\tface,\pm} = \{ (\hat\rho,\omega;\hat\xi,\hat\eta) \colon \hat\rho=0,\ (\hat\xi\mp 1)^2+|\hat\eta|_{\slg^{-1}}^2 = 1 \} \subset \Tsc^*_\sctface\tface,
\]
where we write scattering covectors over $\hat\rho\in[0,\infty)$ as $\hat\xi\frac{\dd\hat\rho}{\hat\rho^2}+\frac{\hat\eta}{\hat\rho}$, $\hat\eta\in T^*\Sph^{n-1}$. The rescaled Hamiltonian vector field $\hat\sfH_p=\hat\rho^{-1}H_p$ takes the form~\eqref{EqStEstHam} with $\sigma=\pm 1$ and $\hat\xi,\hat\eta,\hat\rho$ in place of $\xi,\eta,\rho$, and in particular it has the same sink, resp.\ source structure at the outgoing and incoming radial sets
\begin{equation}
\label{EqSttfRad}
  \cR_{\tface,\rm out}=\{(0,\omega;0,0)\},\qquad \cR_{\tface,\pm,\rm in}=\{(0,\omega;\pm 2,0)\} \subset \Char_{\tface,\pm}.
\end{equation}

\begin{lemma}[Transition face normal operator]
\label{LemmaStEsttf}
  Recall the weights $\beta^-<\beta^+$ from Definition~\usref{DefGSOSpec}, and let $\beta\in(\beta^-,\beta^+)$. Fix on $\tface$ the density $|\hat r^{n-1}\dd\hat r\,\dd\slg|$ where $\hat r=\hat\rho^{-1}$.
  \begin{enumerate}
  \item\label{ItStEsttfsc}{\rm (Variable order estimates.)} Let $\sfs\in\CI({}^{\scop,\bop}S^*\tface)$ and $\sfr\in\CI(\ol{\Tsc^*_\sctface}\tface)$, and suppose that $\pm\hat\sfH_p\sfr\leq 0$ on $\Char_{\tface,\pm}$; suppose moreover that $\sfr<-\half+\ubar S$ at $\cR_{\tface,\rm out}$ and $\sfr>-\half-\ubar S_{\rm in}$ at $\cR_{\tface,\pm,\rm in}$, with $\sfr$ constant near $\cR_{\tface,\rm out}$ and $\cR_{\tface,\pm,\rm in}$. Then the operator
    \begin{equation}
    \label{EqStEsttfsc}
    \begin{split}
      &N_\tface^\pm(P) \colon \bigl\{ u\in H_{\scop,\bop}^{\sfs,\sfr,\frac{n}{2}-\beta}(\tface;\upbeta_\tface^*E|_{\pa X}) \colon N_\tface^\pm(P)u \in H_{\scop,\bop}^{\sfs-2,\sfr+1,\frac{n}{2}-\beta-2}(\tface;\upbeta_\tface^*E|_{\pa X}) \bigr\} \\
      &\hspace{21em} \to H_{\scop,\bop}^{\sfs-2,\sfr+1,-\frac{n}{2}-\beta-2}(\tface;\upbeta_\tface^*E|_{\pa X})
    \end{split}
    \end{equation}
    is invertible. (Here, $\upbeta_\tface\colon\tface\to\pa X$ is the blow-down map as in Definition~\usref{DefGSOtf}.)
  \item\label{ItStEsttfb}{\rm (Estimates on b-spaces.)} Let $s,\ell\in\R$, and suppose that $\ell<-\half+\ubar S$ and $s+\ell>-\half-\ubar S_{\rm in}$. Then the operator
    \begin{equation}
    \label{EqStEsttfb}
    \begin{split}
      &N_\tface^\theta(P) \colon \bigl\{ u\in \Hb^{s,\ell,\frac{n}{2}-\beta}(\tface;\upbeta_\tface^*E|_{\pa X}) \colon N_\tface^\theta(P)u \in \Hb^{s,\ell+1,\frac{n}{2}-\beta-2}(\tface;\upbeta_\tface^*E|_{\pa X}) \bigr\} \\
      &\hspace{21em} \to \Hb^{s,\ell+1,\frac{n}{2}-\beta-2}(\tface;\upbeta_\tface^*E|_{\pa X})
    \end{split}
    \end{equation}
    is invertible for all $\theta\in[0,\pi]$ (with uniformly bounded inverse).
  \end{enumerate}
\end{lemma}
\begin{proof}
  \pfstep{Part~\eqref{ItStEsttfsc}.} Note that $N_\tface^\pm(P)$ is elliptic for $\hat\rho\in(0,\infty]$ and also near fiber infinity of $\ol{{}^{\scop,\bop}T^*}\tface$; moreover, radial point (and real principal type propagation) estimates in $\Char_{\tface,\pm}$ can be proved by repeating the arguments of the proof of Proposition~\ref{PropStEstNz}\eqref{ItStEstNzsc}. At $\ztface$, we pass to $\hat r=\hat\rho^{-1}$ and note that
  \[
    N_\tface^\pm(P) = \hat r^{-2}P_{(0)}(0,\omega,-\hat r D_{\hat r},D_\omega) \bmod \hat r^{-1}\Diffb^1([0,1)_{\hat r}\times\Sph^{n-1};\upbeta_\tface^*E|_{\pa X}).
  \]
  Moreover, if $0<\nu_\bop\in\CI(\tface;\Omegab\tface)$, then in $\hat r<1$ we have $|\hat r^{n-1}\dd\hat r\,\dd\slg|=a\hat r^n\nu_\bop$ for some smooth $a>0$, and therefore $\Hb^{s,\frac{n}{2}-\beta}([0,1)_{\hat r}\times\Sph^{n-1},|\hat r^{n-1}\dd\hat r\,\dd\slg|)=\Hb^{s,-\beta}([0,1)_{\hat r},\nu_\bop)$. By Lemma~\ref{LemmaStEst0}, we have $-\beta\notin\Re\specb(\hat r^{-2}P_{(0)}(0,\omega,-\hat r D_{\hat r},D_\omega))$. (The sign switch arises from $\hat\rho=\hat r^{-1}$.) Combining elliptic b-estimates near $\ztface$ with the aforementioned symbolic estimates, and arguing similarly for the adjoint, gives
  \begin{align}
  \label{EqStEstfscPf}
    \|u\|_{H_{\scop,\bop}^{\sfs,\sfr,\frac{n}{2}-\beta}} &\leq C\Bigl( \|N_\tface^\pm(P)u\|_{H_{\scop,\bop}^{\sfs-2,\sfr+1,\frac{n}{2}-\beta-2}} + \|u\|_{H_{\scop,\bop}^{-N,-N,-N}}\Bigr), \\
    \|\tilde u\|_{H_{\scop,\bop}^{-\sfs+2,-\sfr-1,-\frac{n}{2}+\beta+2}} &\leq C\Bigl( \|N_\tface^\pm(P)^*\tilde u\|_{H_{\scop,\bop}^{-\sfs,-\sfr,-\frac{n}{2}+\beta}} + \|\tilde u\|_{H_{\scop,\bop}^{-N,-N,-N}}\Bigr) \nonumber
  \end{align}
  for any fixed $N$. For later use, we note that these estimates also hold uniformly for $N_\tface^\theta(P)$ when $\theta$ is close to $0$, resp.\ $\pi$, under the assumptions for the sign `$-$', resp.\ `$+$', and for all $\theta\in(0,\pi)$ under only the assumption that $\sfr$ is a constant $<-\half+\ubar S$ near $\cR_{\tface,\rm out}$.

  As a consequence, the operator~\eqref{EqStEsttfsc} is Fredholm. As in the proof of Proposition~\ref{PropStEstNz}, one can show that every element $u$ in the nullspace of $N_\tface^\pm(P)$ is conormal at $\sctface$ (with weight $-\half+\ubar S-\eps$ for all $\eps>0$); symbolic ellipticity (including in the b-sense near $\hat r=0$) implies the conormality of $u$ also at $\hat r=0$, so overall $u\in\cA^{(\alpha,-\beta)}(\tface;\upbeta_\tface^*E|_{\pa X})$. Due to the assumption of Definition~\ref{DefGSOSpec}\eqref{ItGSOSpectf}, we conclude that $u=0$, proving the injectivity of $N_\tface^\pm(P)$.

  It remains to prove the injectivity of $N_\tface^\pm(P)^*$. Any
  \[
    v\in H_{\scop,\bop}^{-\sfs+2,-\sfr-1,-\frac{n}{2}+\beta+2}(\tface;\upbeta_\tface^*E|_{\pa X})\cap\ker N_\tface^\pm(P)^*
  \]
  has infinite scattering/b-regularity (by ellipticity at fiber infinity), and moreover infinite scattering decay outside of $\cR_{\tface,\pm,\rm in}$ (by radial point estimates propagating \emph{out} of $\cR_{\tface,\rm out}$, followed by real principal type propagation). Thus, near $\sctface$ we have $\WFsc(v)\subset\cR_{\tface,\pm,\rm in}$, and therefore
  \[
    \tilde N_\tface^\pm(P)w=0,\qquad
      \tilde N_\tface^\pm(P):=e^{\pm 2 i/\hat\rho}N_\tface^\pm(P)^*e^{\mp 2 i/\hat\rho},\quad
      w := e^{\pm 2 i/\hat\rho}v,
  \]
  with $\WFsc(w)\subset\cR_{\tface,\pm,\rm in}\pm 2\dd(\hat\rho^{-1})=o_\scop$ (the zero section of $\Tsc^*\tface$ over $\sctface$), and with $w$ indeed having a scattering decay order at $o_\scop$ of at least $-\half+\ubar S_{\rm in}-\eps$ for any $\eps>0$. Since the scattering principal symbol of $\tilde N_\tface^\pm(P)$ is the pullback of that of $N_\tface^\pm(P)$ along the translation by $\pm 2\,\dd(\hat\rho^{-1})$, one can propagate module regularity (cf.\ the proof of Proposition~\ref{PropStEstNz}) into the radial set $o_\scop$ and conclude that $w$ is conormal at $\hat\rho=0$. Therefore, $v=e^{\mp 2 i/\hat\rho}w\in\exp(\mp 2 i/\hat\rho)\cA^\alpha(\tface;\upbeta_\tface^*E|_{\pa X})$ near $\hat\rho=0$, and thus
  \[
    v\in\exp\Bigl(\mp\frac{2 i}{\hat\rho/(1+\hat\rho)}\Bigr)\cA^{(\alpha,-n+2+\beta)}(\tface;\upbeta_\tface^*E|_{\pa X})
  \]
  globally. (Note that the exponential prefactor is smooth and nonvanishing down to $\hat\rho^{-1}=0$.) Since $P$ is spectrally admissible, $v=0$.

  \pfstep{Part~\eqref{ItStEsttfb}.} The arguments are completely analogous to those in the first part. The injectivity of $N_\tface^\pm(P)^*$ is most cleanly proved by passing to second microlocal scattering-b-spaces at $\sctface$, as a bridge from the b-spaces of current interest to the scattering spaces already discussed above; a radial point estimate at the lift of $\cR_{\tface,\rm out}$ now implies infinite decay on the lift of $\Tsc^*_\sctface\tface$ except at $\cR_{\tface,\pm,\rm in}$, whereas the infinite b-decay (i.e.\ decay at the scattering zero section) follows from a normal operator argument involving the inversion of the b-normal operator $\pm 2 i\hat\rho(\hat\rho\pa_{\hat\rho}-\frac{n-1}{2}-S|_{\pa X})$. (See the proof of \cite[Theorem~6.1]{HaefnerHintzVasyKerr} for a similar argument.) At this point, one can pass fully to scattering Sobolev spaces and thus to a functional setting in which the injectivity of $N_\tface^\pm(P)^*$ was proved above. Since the operator~\eqref{EqStEsttfb} is Fredholm for all $\theta\in[0,\pi]$ and invertible for $\theta=0,\pi$, its invertibility for all $\theta\in[0,\pi]$ follows from its injectivity (Definition~\ref{DefGSOSpec}\eqref{ItGSOSpectf}) via a Fredholm index argument exactly as in the proof of Corollary~\ref{CorStEstInv}.
\end{proof}

The radial sets $\cR_{\sigma,\rm in}$ and $\cR_{\rm out}$ associated with the spectral family $\hat P(\sigma)$ for fixed nonzero real $\sigma$ with $\pm\sigma>0$ can be assembled into
\[
  \cR_{\pm,\rm in},\qquad \cR_{\pm,\rm out} \subset \Tscbt^*_\scface X,
\]
defined as the closures of $\bigcup_{\sigma\in\pm(0,1)}\{\sigma\}\times\cR_{\sigma,\rm in}$ and $\bigcup_{\sigma\in\pm(0,1)}\{\sigma\}\times\cR_{\rm out}$ (the latter simply being the zero section of the scattering-b-transition cotangent bundle over the closure of $\rho=0$, $\pm\sigma>0$). In the region $\hat\rho=\frac{\rho}{|\sigma|}\lesssim 1$, we write sc-b-transition covectors as $\xi_\scbtop\frac{\dd\hat\rho}{\hat\rho^2} + \frac{\eta_\scbtop}{\hat\rho}$, and then the sc-b-transition principal symbol of the family $\pm[0,1)\ni\sigma\mapsto\sigma^{-2}\hat P(\sigma)$ is $-2\xi_\scbtop+\xi_\scbtop^2+|\eta_\scbtop|_{\slg^{-1}}^2$ modulo elements of $S^2(\Tscbt^*X)$ with coefficients which are conormal of order $\delta$ at $\sctface$, cf.\ Lemma~\ref{LemmaGSOLo} and equation~\eqref{EqStEstNonzeroSymb}.

\begin{prop}[Uniform estimates at low energy]
\label{PropStEstLo}
  Use the notation of Definition~\usref{DefGSOSpec}.
  \begin{enumerate}
  \item\label{ItStEstLosc}{\rm (Variable order estimates.)} Fix $\beta\in(\beta^-,\beta^+)$. Let $\sfs\in\CI(\Sscbt^*X)$, and $\sfr\in\CI(\ol{\Tscbt^*_\scface}X)$. Suppose that $\sfr<-\half+\ubar S$ at $\cR_{\pm,\rm out}$ and $\sfr>-\half-\ubar S_{\rm in}$ at $\cR_{\pm,\rm in}$, with $\pm\sfH_p\sfr\leq 0$ and with $\sfr$ constant near $\cR_{\pm,\rm out}$ and $\cR_{\pm,\rm in}$. Then there exists a constant $C$ so that
    \begin{equation}
    \label{EqStEstLosc}
      \|u\|_{H_{\scbtop,\sigma}^{\sfs,\sfr,-\frac{n}{2}+\beta,0}(X;E_X)} \leq C\|\hat P(\sigma)u\|_{H_{\scbtop,\sigma}^{\sfs-2,\sfr+1,-\frac{n}{2}+\beta+2,0}(X;E_X)},\qquad \sigma\in\pm[0,1).
    \end{equation}
  \item\label{ItStEstLob}{\rm (Estimates on b-spaces.)} Let $s,\ell,\nu\in\R$, and suppose that $\ell<-\half+\ubar S$, $s+\ell>-\half-\ubar S_{\rm in}$, and $\ell-\nu\in(-\frac{n}{2}+\beta^-,-\frac{n}{2}+\beta^+)$. Write $\rho\in\CI(X)$ for a boundary defining function.\footnote{Thus, we may take $\rho=r^{-1}$ near $\pa X$ indeed, for consistency of notation, and merely need to smooth this out near $r=0$ to obtain a valid choice of $\rho$.} Then there exists a constant $C$ so that
    \begin{equation}
    \label{EqStEstLob}
      \| (\rho+|\sigma|)^\nu u \|_{\Hb^{s,\ell}(X;E_X)} \leq C \| (\rho+|\sigma|)^{\nu-1}\hat P(\sigma)u \|_{\Hb^{s,\ell+1}(X;E_X)},\qquad \sigma\in[-1,1]+i[0,1].
    \end{equation}
  \end{enumerate}
\end{prop}
\begin{proof}
  \pfstep{Part~\eqref{ItStEstLosc}.} We follow the proof of~\cite[Proposition~3.21]{HintzKdSMS}, and hence shall be brief; the proof proceeds via upgrading symbolic estimates using the invertibility of the transition face and zero energy operators. To wit, symbolic estimates (i.e.\ only using the sc-b-transition principal symbol, and the subprincipal symbol at radial sets) for the weighted scattering-b-transition operator $\hat P(\sigma)$ give
  \begin{equation}
  \label{EqStEstLoscPf}
    \|u\|_{H_{\scbtop,\sigma}^{\sfs,\sfr,-\frac{n}{2}+\beta,0}} \leq C\Bigl( \|\hat P(\sigma)u\|_{H_{\scbtop,\sigma}^{\sfs-2,\sfr+1,-\frac{n}{2}+\beta+2,0}} + \|u\|_{H_{\scbtop,\sigma}^{-N,\sfr_0,-\frac{n}{2}+\beta,0}}\Bigr),
  \end{equation}
  for any $N$ and $\sfr_0\in\CI(\ol{\Tscbt^*_\scface}X)$, chosen so that $\sfr_0<\sfr$ and so that it satisfies the same conditions as $\sfr$ still. Identify $E|_X$ in a collar neighborhood $[0,\rho_0)_\rho\times\pa X$ of $\pa X$ with the pullback of $E|_{\pa X}$; fix a cutoff $\chi\in\CI(X)$, $\chi\equiv 1$ near $\pa X$, with support in this collar neighborhood. Write further $\hat\rho=\frac{\rho}{|\sigma|}$ and denote by $\phi_\sigma(\hat\rho,\omega)=(|\sigma|\hat\rho,\omega)$ the change of coordinates map. We then have
  \[
    \|\chi u\|_{H_{\scbtop,\sigma}^{-N,\sfr_0,-\frac{n}{2}+\beta,0}\bigish(X;E_X;\rho^{-n}|\frac{\dd\rho}{\rho}\dd\slg|\bigish)} \sim |\sigma|^{-\beta}\| \phi_\sigma^*(\chi u) \|_{H_{\scop,\bop}^{-N,\sfr_0,\frac{n}{2}-\beta}\bigish(\tface;\upbeta_\tface^*E|_{\pa X},\hat\rho^{-n}|\frac{\dd\hat\rho}{\hat\rho}\dd\slg|\bigish)}\,,
  \]
  see \cite[Equation~(A.6b)]{HintzKdSMS} or \cite[Proposition~2.21]{Hintz3b}. We can estimate the right hand side in terms of $N_\tface^\pm(P)(\phi_\sigma^*(\chi u))$ using Lemma~\ref{LemmaStEsttf}; using that in the collar neighborhood, the difference $\hat P(\sigma)-\sigma^2 N_\tface^\pm(P)\in\Diffscbt^{2,-1,-3,0}$ has one more order of decay at $\tface$ than $\hat P(\sigma)$ itself (and vanishes at $\scface$, which makes up for the loss of one scattering decay order), we can then pass back to $\hat P(\sigma)$ and obtain
  \[
    \|\chi u\|_{H_{\scbtop,\sigma}^{-N,\sfr_0,-\frac{n}{2}+\beta,0}} \leq C\Bigl( \|\hat P(\sigma)u\|_{H_{\scbtop,\sigma}^{-N-2,\sfr_0+1,-\frac{n}{2}+\beta+2,0}} + \|u\|_{H_{\scbtop,\sigma}^{-N,\sfr_0,-\frac{n}{2}+\beta+\delta,0}}\Bigr)
  \]
  for $\delta=1$ and, a fortiori, any smaller value of $\delta$. We choose $\delta>0$ small enough so that $\beta+\delta\in(\beta^-,\beta^+)$ still. Noting that $\supp(1-\chi)u\cap\tface=\emptyset$, writing the error term in~\eqref{EqStEstLoscPf} as $\chi u+(1-\chi)u)$ and applying the triangle inequality, we have now obtained the improved estimate
  \[
    \|u\|_{H_{\scbtop,\sigma}^{\sfs,\sfr,-\frac{n}{2}+\beta,0}} \leq C\Bigl( \|\hat P(\sigma)u\|_{H_{\scbtop,\sigma}^{\sfs-2,\sfr+1,-\frac{n}{2}+\beta+2,0}} + \|u\|_{H_{\scbtop,\sigma}^{-N,\sfr_0,-\frac{n}{2}+\beta+\delta,0}}\Bigr).
  \]

  Applying a similar argument to the new error term here, now localizing to a neighborhood of $\zface$ using a cutoff $\psi(\frac{|\sigma|}{\rho})$ with $\psi\in\CIc([0,1))$ identically $1$ near $0$, and exploiting Lemma~\ref{LemmaStEst0}, allows one to relax the error term further to $\|u\|_{H_{\scbtop,\sigma}^{-N,\sfr_0,-\frac{n}{2}+\beta+\delta,-1}}\leq C|\sigma|^\delta\|u\|_{H_{\scbtop,\sigma}^{-N,\sfr_0,-\frac{n}{2}+\beta,0}}$, which for small $|\sigma|$ is bounded by $\half\|u\|_{H_{\scbtop,\sigma}^{\sfs,\sfr,-\frac{n}{2}+\beta,0}}$. This proves the uniform estimate~\eqref{EqStEstLosc} for small enough $|\sigma|$. For $\sigma$ away from $0$, the estimate~\eqref{EqStEstLosc} is equivalent to~\eqref{EqStEstNzscEst}.

  \pfstep{Part~\eqref{ItStEstLob}.} This is the content of \cite[Theorem~1.1]{VasyLowEnergyLag} in the form given in \cite[Theorem~2.11]{HintzPrice}, albeit in the more general setting of \cite[Theorem~2.5]{VasyLowEnergyLag}. If one second microlocalizes at the zero section of $\Tscbt^*X$ over $\scface$, thereby introducing the b-decay order $\ell$ (i.e.\ decay at zero scattering frequency), one has the uniform estimate
  \[
    \|u\|_{H_{\scbtop,2}^{s,r,\ell,-\frac{n}{2}+\beta,0}} \leq C\|\hat P(\sigma)u\|_{H_{\scbtop,2}^{s-2,r+1,\ell+1,-\frac{n}{2}+\beta+2,0}},\qquad \sigma\in[-1,1]+i[0,1],
  \]
  where (for simpler comparison with the first part of the proof) we use an ad hoc notation for the corresponding scale of Sobolev spaces in which the scattering decay order is split into two orders: the constant scattering decay order $r$ at nonzero frequencies and the b-decay order $\ell$. (In the notation of \cite{VasyLowEnergyLag}, one has $\|u\|_{H_{\scbtop,2}^{s,r,\ell,\ell-\nu,0}}=\|(\rho+|\sigma|)^\nu u\|_{H_{\scop,\bop,\res}^{s,r,\ell}}$.) This estimate uses the uniform bound $\|u\|_{\Hb^{s,\ell,\frac{n}{2}-\beta}}\lesssim\|N_\tface^\theta(P)u\|_{\Hb^{s,\ell+1,\frac{n}{2}-\beta-2}}$ of Lemma~\ref{LemmaStEsttf}\eqref{ItStEsttfb}. To deduce~\eqref{EqStEstLob} from this, one notes that
  \[
    \|(\rho+|\sigma|)^\nu u\|_{\Hb^{s,\ell}} \sim \|u\|_{H_{\scbtop,2}^{s,r,\ell,-\frac{n}{2}+\beta,0}},\qquad r=s+\ell,\quad \beta=\ell-\nu+\tfrac{n}{2},
  \]
  and that for these $r,\beta$, one also has
  \[
    \|\hat P(\sigma)u\|_{H_{\scbtop,2}^{s-2,r+1,\ell+1,-\frac{n}{2}+\beta+2,0}} \lesssim \|\hat P(\sigma)u\|_{H_{\scbtop,2}^{s,s+\ell+1,\ell+1,\ell-\nu+2,0}} \sim \|(\rho+|\sigma|)^{\nu-1}\hat P(\sigma)u\|_{\Hb^{s,\ell+1}}.\qedhere
  \]
\end{proof}

\subsection{Conormality of the resolvent; decay and regularity of forward solutions}
\label{SsStCo}

In this section, we only use the b-estimates from~\S\ref{SsStEst}. The first result is closely related to \cite[\S\S2.2.1--2.2.2]{HintzPrice} and \cite[Propositions~12.4 and 12.12]{HaefnerHintzVasyKerr}:

\begin{lemma}[Conormality of the resolvent]
\label{LemmaStCo}
  Let $s,\ell\in\R$ and $k\in\N_0$, and suppose that $\ell<-\half+\ubar S$ and $s+\ell>-\half-\ubar S_{\rm in}$ in the notation of Definition~\usref{DefStEstThr}.
  \begin{enumerate}
  \item\label{ItStCoLo}{\rm (Low frequencies.)} Let $\nu\in\R$ with $\ell-\nu\in(-\frac{n}{2}+\beta^-,-\frac{n}{2}+\beta^+)$. Then the operator
    \begin{equation}
    \label{EqStCoLo}
      (\sigma\pa_\sigma)^k \hat P(\sigma)^{-1} \colon (\rho+|\sigma|)^{-\nu+1}\Hb^{s+k,\ell+1}(X;E_X) \to (\rho+|\sigma|)^{-\nu}\Hb^{s,\ell}(X;E_X)
    \end{equation}
    is uniformly bounded for $\sigma\in[-1,1]+i[0,1]$.
  \item\label{ItStCoMed}{\rm (Bounded frequencies.)} Let $0<C_0<C_1$ and $C_2>0$. Then the operator
    \begin{equation}
    \label{EqStCoMed}
      \pa_\sigma^k \hat P(\sigma)^{-1} \colon \Hb^{s+k,\ell+1}(X;E_X) \to \Hb^{s,\ell}(X;E_X)
    \end{equation}
    is uniformly bounded for $\sigma\in\pm[C_0,C_1]+i[0,C_2]$.
  \item\label{ItStCoHi}{\rm (High frequencies.)} Let $C_0>0$. Then the operator
    \begin{equation}
    \label{EqStCoHi}
      (\sigma\pa_\sigma)^k \hat P(\sigma)^{-1} \colon H_{\bop,|\sigma|^{-1}}^{s+k,\ell+1}(X;E_X) \to |\sigma|^{-1+k}H_{\bop,|\sigma|^{-1}}^{s,\ell}(X;E_X)
    \end{equation}
    is uniformly bounded for $\sigma\in\C$, $\Im\sigma\geq 0$, $|\sigma|\geq C_0$.
  \end{enumerate}
\end{lemma}
\begin{proof}
  We begin with part~\eqref{ItStCoMed}. For $k=0$, and using Corollary~\ref{CorStEstInv}, the estimate~\eqref{EqStCoMed} is the same as~\eqref{EqStEstNzbEst}. Next,
  \[
    \pa_\sigma \hat P(\sigma)^{-1} = -\hat P(\sigma)^{-1}\circ\pa_\sigma\hat P(\sigma)\circ\hat P(\sigma)^{-1}
  \]
  is uniformly bounded as a map
  \[
    \Hb^{s+1,\ell+1} \xra{\hat P(\sigma)^{-1}} \Hb^{s+1,\ell} \xra{\pa_\sigma\hat P(\sigma)} \Hb^{s,\ell+1} \xra{\hat P(\sigma)^{-1}} \Hb^{s,\ell},
  \]
  where we use $\pa_\sigma\hat P(\sigma)\in\rho\Diffb^1(X;E_X)$. The estimate~\eqref{EqStCoMed} for $k\geq 2$ follows inductively in an analogous manner, using also that $\pa_\sigma^2\hat P(\sigma)\in\cA^{1+\delta}(X;\End(E_X))$.

  The proof of~\eqref{EqStCoHi} is completely analogous, with the estimate~\eqref{EqStEstHib} giving the case $k=0$; note now that with $h=|\sigma|^{-1}$ denoting the semiclassical parameter, we have $\sigma\pa_\sigma\hat P(\sigma)\in\sigma\rho\Diffb^1(X;E_X)\subset h^{-2}\rho\Diff_{\bop,h}^1(X;E_X)$. Similarly, part~\eqref{ItStCoLo} follows by repeated application of the estimate~\eqref{EqStEstLob} and the observation that for $\sigma\in[-1,1]+i[0,1])$, the map
  \[
    \sigma\pa_\sigma\hat P(\sigma) = 2 i\sigma\rho\Bigl(\rho\pa_\rho-\frac{n-1}{2}-S\Bigr) - i\sigma Q + 2\sigma^2 g^{0 0},\qquad Q\in\cA^{2+\delta}\Diffb^1,\ \ g^{0 0}\in\cA^{1+\delta}\CI,
  \]
  satisfies uniform bounds
  \[
    \|(\rho+|\sigma|)^{\nu-1}\sigma\pa_\sigma\hat P(\sigma)u\|_{\Hb^{s,\ell+1}} \leq \| (\rho+|\sigma|)^\nu u\|_{\Hb^{s+1,\ell}}
  \]
  for any $\nu\in\R$.
\end{proof}

The uniform low energy estimates of Lemma~\ref{LemmaStCo}\eqref{ItStCoLo} can be improved to a certain amount of smoothness at the lift
\[
  \zface\subset X_\res:=\bigl[[-1,1]\times X;\{0\}\times\pa X\bigr]
\]
of $\sigma^{-1}(0)\subset[-1,1]\times X$ to $X_\res$. For present purposes, we only need a result for a simple class of conormal inputs (in fact, the case $\ell=\infty$ in Proposition~\ref{PropStCoReg} below is sufficient for later). More precise results can be obtained with more careful bookkeeping, see \cite{MorganWunschPrice} for a concrete example. As in \cite[\S3.2.2]{HintzPrice}, we write the sc-b-transition single spaces of $X$ for positive and negative frequencies as
\[
  X_\res^\pm = [ \pm[0,1]\times X; \{0\}\times\pa X].
\]
We write $\cA^{(\alpha,\beta,\gamma)}(X_\res^\pm)$ for the space of conormal functions on $X_\res^\pm$ with weights $\alpha$, $\beta$, $\gamma$ at $\scface^\pm$ (the lift of $\pm[0,1]\times\pa X$), $\tface^\pm$ (the front face), $\zface$ (the lift of $\{0\}\times X$), respectively. We let $\cA^{(\alpha-,\beta,\gamma)}(X_\res^\pm):=\bigcup_{\eps>0}\cA^{(\alpha-\eps,\beta,\gamma)}(X_\res^\pm)$ etc. We furthermore write $\cA^{(\alpha,\beta,(j,0))}(X_\res^\pm)$ for functions $u$ which are polyhomogeneous with index set $(j,0)$ down to $\zface$; this means that in a collar neighborhood $[0,1)_{\hat r}\times[0,\rho_0)_\rho\times\Sph^{n-1}$, $\hat r:=\frac{|\sigma|}{\rho}$, of $\zface\cap\tface^\pm\subset X_\res^\pm$, one has $u(\hat r,\rho,\omega)\in\hat r^j\CI([0,1)_{\hat r};\cA^\beta([0,1)_\rho\times\Sph^{n-1}_\omega))$. We furthermore set $\cA^{(\alpha,\beta,((j,0),\gamma))}(X_\res^\pm)=\cA^{(\alpha,\beta,(j,0))}(X_\res^\pm)+\cA^{(\alpha,\beta,\gamma)}(X_\res^\pm)$. Spaces encoding partial expansions at several boundary hypersurfaces are denoted $\cA^{(\alpha,((k,0),\beta),((j,0),\gamma))}(X_\res^\pm)$; see also \cite[Definition~2.13]{HintzPrice}.

\begin{prop}[Improved regularity of the low energy resolvent]
\label{PropStCoReg}
  Define $k:=\lceil\beta^+-\beta^-\rceil\in\N$. Fix a boundary defining function $\rho\in\CI(X)$, and let $\psi\in\CIc(\R)$ be identically $1$ near $0$. Let $\hat f\in\CI([-1,1]_\sigma;\cA^{\ell+2}(X;E_X))$, where $\ell\in\R$ satisfies $\ell+1\geq\max(\ubar S+\frac{n-1}{2},\beta^++1)$, Then there exist
  \begin{align*}
    u_j&\in\cA^{\beta^+-}(X;E_X),\qquad j=0,\ldots,k-1, \\
    \tilde u^\pm&\in\cA^{(\frac{n-1}{2}+\ubar S-,\;\beta^+-,\;\beta^+-\beta^--)}(X_\res^\pm;E_X)
  \end{align*}
  so that
  \[
    \hat P(\sigma)^{-1}\hat f(\sigma) = \sum_{j=0}^{k-1} \psi\Bigl(\frac{\sigma}{\rho}\Bigr)\Bigl(\frac{\sigma}{\rho}\Bigr)^j u_j + \tilde u^+ + \tilde u^-,
  \]
  where we write $\tilde u^\pm$ also for the extension by $0$ to $\pm\sigma<0$.
\end{prop}
\begin{proof}
  By Lemma~\ref{LemmaStEst0} and Sobolev embedding, we have
  \[
    u_0:=\hat P(0)^{-1}\hat f(0)\in\cA^{\beta^+-}(X;E_X).
  \]
  Note now that $\hat f_0:=\hat f\in\cA^{(\ell+2,\ell+2,(0,0))}(X_\res^\pm;E_X)$. Thus, writing $\psi=\psi(\frac{\sigma}{\rho})$, the function $\hat f_1$ defined by
  \begin{equation}
  \label{EqStCoRegfu1}
  \begin{split}
    \hat f^1(\sigma) &:= \hat f(\sigma) - \hat P(\sigma) ( \psi u_0 ) \\
      &= (1-\psi)\hat f(\sigma) + \psi\bigl(\hat f(\sigma)-\hat f(0)\bigr) + \psi(\hat f(0)-\hat P(0)u_0) \\
      &\qquad - \psi(\hat P(\sigma)-\hat P(0))u_0 - [\hat P(\sigma),\psi]u_0
  \end{split}
  \end{equation}
  satisfies $\hat f^1\in\cA^{(\ell+2,\beta^++2-,(1,0))}(X_\res^\pm;E_X)$ since each one of the five terms on the right in~\eqref{EqStCoRegfu1} lies in this space; this uses in particular that the difference $\hat P(\sigma)-\hat P(0)$ maps $\cA^{(\alpha,\beta,((k,0),\gamma))}(X_\res^\pm;E_X)\to\cA^{(\alpha+1,\beta+2,((k+1,0),\gamma+1))}(X_\res^\pm;E_X)$. Therefore,
  \begin{equation}
  \label{EqStCoRegfl1}
    \hat f_1(\sigma) := \sigma^{-1}\hat f^1(\sigma) \in \cA^{(\ell+2,\beta^++1-,(0,0))}(X_\res^\pm;E_X).
  \end{equation}
  We stress that the boundary value $\hat f_1(0)$ at $\zface$ is the same on $X_\res^+$ and on $X_\res^-$.

  If $\beta^+-1>\beta^-$, we can put $u_1:=\hat P(0)^{-1}\hat f_1(0)\in\cA^{\beta^+-1-}(X;E_X)$, and then
  \[
    \hat f_2(\sigma) := \sigma^{-1}\bigl(\hat f_1(\sigma)-\hat P(\sigma)(\psi u_1)\bigr) \in \cA^{(\ell+2,\beta^+-,(0,0))}(X_\res^\pm;E_X).
  \]
  Continuing in this manner, we obtain elements $u_j\in\cA^{\beta^+-j-}(X;E_X)$, $0\leq j<k$, with the property that
  \begin{equation}
  \label{EqStCoRegfj}
    \hat f_j = \sigma^{-j}\biggl( \hat f(\sigma) - \hat P(\sigma)\Bigl(\psi\sum_{m=0}^{j-1} \sigma^m u_m\Bigr) \biggr) \in \cA^{(\ell+2,\beta^++2-j-,(0,0))}(X_\res^\pm;E_X).
  \end{equation}
  Indeed, when $j\leq k-1$, then the decay rate of $\hat f_j(0)\in\cA^{\beta^++2-j-}(X;E_X)$ lies in $(\beta^-+2,\beta^++2)$, and hence $\hat P(0)^{-1}$ from Lemma~\ref{LemmaStEst0} can be applied to it.

  Finally, for the $k$-th error term we record merely the conormal membership
  \[
    \hat f_k \in \cA^{(\ell+2,\beta^++2-k-\eps,0)}(X_\res^\pm;E_X) \subset |\sigma|^{\beta^+-\beta^--k-2\eps}\cA^{(\ell+2,\beta^-+2+\eps,0)}(X_\res^\pm;E_X)
  \]
  for any $\eps>0$. A fortiori, this implies that for $\sigma\in[-1,1]$, the rescaling
  \[
    \tilde f_k(\sigma) := \sigma^k |\sigma|^{-(\beta^+-\beta^-)+2\eps}\hat f_k(\sigma)
  \]
  is uniformly bounded as an element
  \[
    \tilde f_k(\sigma) \in (\rho+|\sigma|)^{-\nu+1}\Hb^{s,\ell+1}(X;E_X),\qquad \nu:=-\beta^--\eps+\ell+\tfrac{n}{2},
  \]
  for all $s\in\R$ and $\eps>0$. (Note that $\frac{n}{2}+(\ell-\nu)=\beta^-+\eps$, with the summand $\frac{n}{2}$ arising from passing between $|\dd g_X|$ and positive b-densities.) Thus, Lemma~\ref{LemmaStCo}\eqref{ItStCoLo} (where we can take $k$ arbitrarily large) implies that
  \[
    \tilde u_k(\sigma) := \hat P(\sigma)^{-1}\tilde f_k(\sigma) \in \cA^{(\frac{n-1}{2}+\ubar S-,\beta^-+\eps-,0)}(X_\res^\pm;E_X)
  \]
  Altogether, we have
  \[
    \hat P(\sigma)^{-1}\hat f = \sum_{m=0}^{k-1} \psi\sigma^m u_m + |\sigma|^{(\beta^+-\beta^-)-2\eps}\tilde u_k(\sigma).
  \]
  Since $\eps>0$ is arbitrary, this completes the proof upon writing $\sigma^m u_m=(\frac{\sigma}{\rho})^m\rho^m u_m$ and renaming $\rho^m u_m\in\cA^{\beta^+-}(X;E_X)$ into $u_m$.
\end{proof}

We can now state and prove our first main result for forward solutions of the stationary model problem:

\begin{thm}[Forward solutions for simple inputs]
\label{ThmStCo}
  Let $P$ be a stationary wave type operator (Definition~\usref{DefGSO}) with respect to a stationary and asymptotically flat Lorentzian metric (Definition~\usref{DefGSG}) on $M_0=\ol\R\times X$, $X=\ol{\R^n}$. Suppose $P$ is spectrally admissible with zero energy weights in $(\beta^-,\beta^+)$ (Definition~\usref{DefGSOSpec}). Define $\ubar S\in\R$ as in Definition~\usref{DefStEstThr}. Let $f\in\CIdot(M_0;E)$, with support in $t_*\geq 0$. (More generally, one can allow $f\in\sS(\R_{t_*};\cA^{\ell+2}(X;E_X))$ where $\ell+1\geq\max(\ubar S+\frac{n-1}{2},\beta^++1)$, still assuming $t_*\geq 0$ on $\supp f$.) Then the unique solution of
  \[
    P u = f,\qquad u|_{t_*<0}=0,
  \]
  is conormal on $M_0$, and on the resolution $M_1=[M_0;\pa\ol\R\times\pa X]$ satisfies
  \begin{equation}
  \label{EqStCoSol}
    u \in \cA^{(\frac{n-1}{2}+\ubar S-,\beta^++1-,\beta^+-\beta^-+1-)}(M_1;E),
  \end{equation}
  where the weights refer to the boundary hypersurfaces $\scri^+$ (lift of $\ol\R\times\pa X$), $\iota^+$ (lift of $\{\infty\}\times\pa X$), $\cT^+\subset M_1$ (lift of $\{\infty\}\times X$), respectively. (See Definition~\usref{DefGAMfd} and Figure~\usref{FigGAMfd}.) That is, fixing a positive definite fiber inner product on $E\to X$, we have the pointwise estimate
  \[
    |Z^J u| \lesssim \rho_{\!\scri}^{\frac{n-1}{2}+\ubar S-\eps}\rho_+^{\beta^++1-\eps}\rho_\cT^{\beta^+-\beta^-+1-\eps}
  \]
  for all $\eps>0$ and all multiindices $J$, where $Z$ is a spanning set of the set of b-vector fields on $M_1$. (One may take $Z=\{\la t_*\ra\pa_{t_*},\,\la x\ra\pa_{x^1},\ldots,\la x\ra\pa_{x^n}\}$.)
\end{thm}
\begin{proof}
  Due to the support assumption on $f$, we have $\hat f(\cdot+i\tau)\in\sS(\R_\sigma;\cA^{\ell+2}(X;E_X))$ for $\tau\geq 0$, with all seminorms bounded uniformly in $\tau$; and $\hat f(\sigma)$ is holomorphic in $\sigma$. Likewise then,
  \[
    \hat u(\sigma) := \hat P(\sigma)^{-1}\hat f(\sigma),\qquad \Im\sigma>0,
  \]
  is holomorphic in $\sigma$, and Schwartz in $\Re\sigma$ (uniformly for $\Im\sigma\geq\eps>0$ for any $\eps>0$) with values in $\cA^{\frac{n-1}{2}+\ubar S-}(X;E_X)$ by Proposition~\ref{PropStEstHi}\eqref{ItStEstHib}. By the Paley--Wiener theorem, we conclude that
  \begin{equation}
  \label{EqStCoContour}
    u(t_*,x) := \cF^{-1}\bigl(\hat P(\sigma)^{-1}\hat f(\sigma,x)\bigr) = (2\pi)^{-1}\int_{\Im\sigma=C} e^{-i\sigma t_*}\hat P(\sigma)^{-1}\hat f(\sigma,x)\,\dd\sigma
  \end{equation}
  is supported in $t_*\geq 0$, where $C>0$ is arbitrary, as follows via contour shifting (justified using Proposition~\ref{PropStEstHi}\eqref{ItStEstHib}). The uniform estimates on $\hat P(\sigma)$ for $0\leq\Im\sigma\leq C$ following from Proposition~\ref{PropStEstNz}\eqref{ItStEstNzb} (bounded frequencies), Proposition~\ref{PropStEstHi}\eqref{ItStEstHib} (high frequencies), and Proposition~\ref{PropStEstLo}\eqref{ItStEstLob} (low frequencies) together with Corollary~\ref{CorStEstC0} imply that we can take $C=0$ in~\eqref{EqStCoContour}.

  Let $\chi\in\CIc((-1,1))$ be identically $1$ near $0$. Setting
  \begin{align*}
    u_{\rm lo}(t_*,x) &= \cF^{-1}\bigl(\chi(\sigma)\hat P(\sigma)^{-1}\hat f(\sigma,x)\bigr), \\
    u_{\rm hi}(t_*,x) &= \cF^{-1}\bigl((1-\chi(\sigma))\hat P(\sigma)^{-1}\hat f(\sigma,x)\bigr),
  \end{align*}
  we can then control $u_{\rm hi}$ using Lemma~\ref{LemmaStCo}\eqref{ItStCoHi} and conclude that
  \begin{equation}
  \label{EqStCouhi}
    u_{\rm hi} \in \sS\bigl(\R_{t_*};\cA^{\frac{n-1}{2}+\ubar S-}(X;E_X)\bigr) \subset \cA^{(\frac{n-1}{2}+\ubar S-,\infty,\infty)}(M_1;E).
  \end{equation}
  Here, $u_{\rm hi}$ typically has full support in $t_*$, and the weights refer to $\scri^+$, the lift of $\pa\ol\R\times\pa X$, and the lift of $\pa\ol\R\times X$ (in this order).

  On the other hand, the low energy contribution $u_{\rm lo}$ is described by Proposition~\ref{PropStCoReg}. In the notation of that Proposition, the inverse Fourier transform of $\psi(\frac{\sigma}{\rho})(\frac{\sigma}{\rho})^j u_j(x)$ for $j=0,\ldots,k-1$, is given by
  \begin{equation}
  \label{EqStCoCOV}
  \begin{split}
    u_j(x) \rho (2\pi)^{-1} \int_\R e^{-i\hat r(\rho t_*)} \psi(\hat r)\hat r^j\,\dd\hat r &\in \rho\cA^{\beta^+-}(X;E_X)\cdot \sS(\R_{\rho t_*}) \\
      &\subset \cA^{(\beta^++1-,\beta^++1-,\infty)}(M_1;E);
  \end{split}
  \end{equation}
  note here that in $|\rho t_*|>1$, the reciprocal $(\rho t_*)^{-1}=\frac{r}{t_*}$ is a defining function of the lift of $\pa\ol\R\times X$. The inverse Fourier transform of $(1-\psi(\hat r))\tilde u^\pm\in\cA^{(\frac{n-1}{2}+\ubar S-,\beta^+-,\infty)}(X_\res^\pm;E_X)$ can be controlled using \cite[Proposition~2.29(1)]{Hintz3b}, with the result
  \begin{equation}
  \label{EqStCoTildeu1}
    \cF^{-1}\bigl((1-\psi)\tilde u^\pm\bigr) \in \cA^{(\min(\frac{n-1}{2}+\ubar S,\beta^++1)-,\beta^++1-,\infty)}(M_1;E).
  \end{equation}
  For the computation of $\cF^{-1}(\psi\tilde u^\pm)$ on the other hand we change coordinates as in~\eqref{EqStCoCOV} and use \cite[Lemma~2.25(2)]{Hintz3b} to obtain
  \begin{equation}
  \label{EqStCoTildeu2}
    \cF^{-1}\bigl(\psi\tilde u^\pm\bigr) \in \cA^{\beta^+-\beta^-+1-}(\ol{\R_{\rho t_*}})\cdot\rho\cA^{\beta^+-}(X;E_X) \subset \cA^{(\beta^++1-,\beta^++1-,\beta^+-\beta^-+1-)}(M_1;E).
  \end{equation}

  With $u_{\rm lo}$ being the sum of~\eqref{EqStCoCOV}--\eqref{EqStCoTildeu2}, and with $u_{\rm hi}$ given by~\eqref{EqStCouhi}, we conclude that the forward solution $u$ satisfies
  \[
    u = u_{\rm lo}+u_{\rm hi} \in \cA^{(\beta_{\!\scri}-,\beta^++1-,\beta^+-\beta^-+1-)}(M_1;E),\qquad
    \beta_{\!\scri} := \min\Bigl(\frac{n-1}{2}+\ubar S,\beta^++1\Bigr).
  \]
  If $\beta^++1\geq\frac{n-1}{2}+\ubar S$, then $\beta_{\!\scri}=\frac{n-1}{2}+\ubar S$, and thus we are done.

  If, on the other hand, $\beta^++1<\frac{n-1}{2}+\ubar S$ and thus $\beta_{\!\scri}=\beta^++1$, then the decay rate at $\scri^+$ can be sharpened as follows. For brevity, we only record weights at $\scri^+$ and $\iota^+$; suppose we have already established $u\in\cA^{(\beta'_{\!\scri}-,\beta^++1-)}$ where $\beta'_{\!\scri}\geq\beta^++1$. Using the leading order behavior~\eqref{EqGAWEdgeN} of $P$ (with $p_1=S|_{\pa X}$ and $p_0=0$, cf.\ Lemma~\ref{LemmaGAWStat}), and using $\rho_{\!\scri}=x_{\!\scri}^2$, we then have
  \[
    \Bigl(\rho_{\!\scri}\pa_{\rho_{\!\scri}}-\Bigl(\frac{n-1}{2}+S|_{\pa X}\Bigr)\Bigr)(\rho_{\!\scri}\pa_{\rho_{\!\scri}}-\rho_+\pa_{\rho_+})u \in \cA^{(\beta'_{\!\scri}+\frac12-,\beta^++1-)}.
  \]
  Since $S|_{\pa X}$ commutes with the vector field $\rho_{\!\scri}\pa_{\rho_{\!\scri}}-\rho_+\pa_{\rho_+}$, we can first integrate this vector field starting from $t_*=0$ where $u$ vanishes; since $\beta^++1\leq\beta_{\!\scri}'<\beta_{\!\scri}'+\half$, this gives
  \[
    \Bigl(\rho_{\!\scri}\pa_{\rho_{\!\scri}}-\frac{n-1}{2}-S|_{\pa X}\Bigr)u \in \cA^{(\beta'_{\!\scri}+\frac12-,\beta^++1-)}.
  \]
  (See also \cite[Lemma~7.7(1)]{HintzVasyMink4}.) Upon integrating this towards $\rho_{\!\scri}=0$, we obtain
  \[
    u \in \cA^{(\beta''_{\!\scri}-,\beta^++1-)},\qquad \beta''_{\!\scri}=\min\Bigl(\beta'_{\!\scri}+\frac12,\ \frac{n-1}{2}+\ubar S\Bigr).
  \]
  Thus, we can improve the order of $u$ at $\scri^+$ by (at most) half a power of $\rho_{\!\scri}$ until we reach the desired order $\frac{n-1}{2}+\ubar S-$.
\end{proof}

\subsection{Sharpness of decay rates}
\label{SAS}

In some cases, the pointwise decay of forward solutions with Schwartz forcing obtained in Theorem~\ref{ThmStCo} (namely, $t_*^{-\beta^++\beta^--1}\la r\ra^{-\beta^-}$ in $t_*\geq 1$, $r<q t$, $q\in(0,1)$, up to a $t_*^\eps$ loss) can be improved. A dramatic example is given by Minkowski spacetimes with even spacetime dimensions $n+1\geq 4$: the sharp Huygens principle gives Schwartz decay of $u$ at $\cT^+\cup\iota^+$, even though $(\beta^-,\beta^+)=(0,n-2)$ (which we verify in~\S\ref{SE}). In $3+1$ dimensions, the wave operator of a stationary and asymptotically flat (with mass $\bhm\neq 0$) metric which is spectral admissible in the sense of \cite[Definitions~2.3 and 2.9]{HintzPrice} has $(\beta^-,\beta^+)=(0,1)$ (cf.\ \cite[Equation~(2.9)]{HintzPrice}), but the pointwise decay rate at $\cT^+\cup\iota^+$ is $t_*^{-3}$ by \cite[Theorem~3.9]{HintzPrice}. We proceed to the analyze the way in which these examples are exceptional in the class of stationary operators considered here, and how Theorem~\ref{ThmStCo} is (conjecturally) generically sharp up to the arbitrarily small $t_*^\eps$ loss.\footnote{This section can be skipped at first reading, as the material developed here is not used elsewhere in the paper.}

Thus, let $P$ be a stationary wave type operator which is spectrally admissible with indicial gap $(\beta^-,\beta^+)$, and define $S,\hat P(0)$ via~\eqref{EqGSOStruct}. Put
\begin{equation}
\label{EqASSpace}
  \cF :=\{f\in\CIdot(M_0;E)\colon t_*\geq 0\ \text{on}\ \supp f\},
\end{equation}
and let $f\in\cF$. Our aim is to find conditions on $P$ and $f$ so that the forward solution $u$ of $P u=f$ obeys a lower bound matching the upper bound (without the arbitrarily small loss) in Theorem~\ref{ThmStCo}. In fact, we shall find conditions so that $u$ has a leading order term at $\cT^+$ and $\iota^+$ (with decay at $\scri^+$ matching Theorem~\ref{ThmStCo}). See Theorem~\ref{ThmAS} below for the final result. To this end, we shall use the algorithmic procedure introduced in \cite{HintzPrice} to produce an expansion of the low energy resolvent $\hat P(\sigma)^{-1}\hat f(\sigma)$ with explicit leading singular terms.

\textbf{Taylor expansion at $\zface$.} To study when inverse polynomial lower bounds hold as $t_*\to\infty$, it suffices by~\eqref{EqStCouhi} to analyze the low energy resolvent. Let
\[
  u_0 = \hat P(0)^{-1}\hat f(0) \in \cA^{\beta^+-}(X;E_X).
\]
We first claim that there is a proper subspace
\begin{equation}
\label{EqASSubspace}
  \cF_0\subsetneq\cF
\end{equation}
so that for $f\in\cF\setminus\cF_0$, the solution $u_0$ does not lie in $\cA^{\beta^++\eta}(X;E_X)$ for any $\eta>0$. Indeed, for all $\beta>\beta^+$, the index of the Fredholm operator $\hat P(0)\colon\Hb^{s,-\frac{n}{2}+\beta}(X;E_X)\to\Hb^{s-2,-\frac{n}{2}+\beta+2}(X;E_X)$ is nonzero by \cite[Theorem~6.5]{MelroseAPS} (see the proof of Lemma~\ref{LemmaStEst0}); but this operator is injective, and therefore it cannot be surjective. In particular, if $\eps>0$ is so small that no element of $\specb(\rho^{-2}\hat P(0))$ has real part in $(\beta^+,\beta^++2\eps)$, then for $\beta=\beta^++\eps$, the range of $\hat P(0)$ (as an operator between the above spaces)---which is a closed subspace of $\Hb^{s-2,-\frac{n}{2}+\beta+2}(X;E_X)$---cannot contain the dense subspace $\CIdot(X;E_X)\subset\Hb^{s-2,-\frac{n}{2}+\beta+2}(X;E_X)$. A fortiori, $\CIdot(X;E_X)\not\subset\bigcup_{\eta>0}\hat P(0)(\cA^{\beta^++\eta}(X;E_X))$.

Since $\hat f(0)\in\CIdot(X;E_X)$, the solution $u_0$ has a polyhomogeneous expansion at $\pa X$. As regards its leading order term, we make the following simplifying assumption:
\begin{equation}
\label{EqASAssmSimple}
  \parbox{0.8\textwidth}{There is only one pole $\lambda=\lambda^\pm$ of $\wh N(\rho^{-2}\hat P(0),-i\lambda)^{-1}$ with $\Re\lambda^\pm=\beta^\pm$, and it is simple. Moreover, $\ker\wh N(\rho^{-2}\hat P(0),-i\lambda^\pm)$ is 1-dimensional and spanned by $v^\pm\in\CI(\pa X;E|_{\pa X})$.}
\end{equation}
For now, we only need this for the `$+$' sign. For $f\in\cF\setminus\cF_0$, we then have
\begin{equation}
\label{EqASu0}
  u_0 = u_{0 0}\chi \rho^{\lambda^+}v^+ + \tilde u_0,\qquad 0\neq u_{0 0}\in\C,\quad\tilde u_0\in\cA^{\beta^++\eps}(X;E_X),
\end{equation}
where $\chi\in\CIc([0,1)_\rho)$ is a cutoff which is identically $1$ near $0$, and $\eps\in(0,1)$ is small. (Assumption~\eqref{EqASAssmSimple} can easily be relaxed: if there are several poles with real part $\beta^+$, then there are several leading order terms of $u_0$ here; and if the poles have higher multiplicity, one needs to allow for additional factors of $|\log\rho|^j$, $j\in\N_0$.)

We next modify~\eqref{EqStCoRegfu1}--\eqref{EqStCoRegfl1} slightly and set
\begin{equation}
\label{EqASf1}
\begin{split}
  \hat f_1(\sigma) := \sigma^{-1}\bigl(\hat f(\sigma)-\hat P(\sigma)u_0\bigr) &= -\pa_\sigma\hat P(0)u_0 + \hat f^{(1)}(\sigma) - \half\sigma\pa_\sigma^2\hat P(0)u_0 \\
    & \in \CI\bigl(\pm[0,1)_\sigma;\cA^{(\beta^++1,0),\beta^++1+\eps}(X;E_X)\bigr)
\end{split}
\end{equation}
where $\hat f^{(1)}(\sigma):=\sigma^{-1}(\hat f(\sigma)-\hat f(0))$.\footnote{Even though at present we are constructing $\hat P(\sigma)^{-1}\hat f(\sigma)$ in Taylor series at $\zface$, we omit the cutoff $\psi$ for two reasons: firstly, we can afford arbitrary imprecisions or decay losses at $\scface^\pm\subset X_\res^\pm$; secondly, without $\psi$, a certain condition on a $\tface^\pm$-model problem below will take a simpler form.} We note that $\hat f^{(1)}\in\sS(\R_\sigma;\CIdot(X;E_X))$ vanishes rapidly at $\tface^\pm$; and the final term in~\eqref{EqASf1} has decay order $\beta^++2+\delta$ at $\tface^\pm$. Moreover, we have $\hat f_1(0)=u_{0 0}\chi\rho^{\lambda^++1}f_{1,0}+\tilde f_1$ where $\tilde f_1\in\cA^{\beta^++1+\eps}(X;E_X)$ and
\begin{equation}
\label{EqASf10}
  f_{1,0} = -2 i\Bigl(\lambda^+-\frac{n-1}{2}-S|_{\pa X}\Bigr)v^+ \in \CI(\pa X;E|_{\pa X}).
\end{equation}
In order to proceed, we need to assume that $f_{1,0}\neq 0$. This holds if and only if there exists $p\in\pa X$ so that $\lambda^+\notin\frac{n-1}{2}+\spec S(p)$; indeed, this condition holds for an open set of $p$, and $v^+$ does not vanish on any nonempty open set by unique continuation.

\begin{rmk}[Counterexamples: I]
\label{RmkASCounterex1}
  By contrast, one \emph{always} has $f_{1,0}=0$ for the wave operator on $(3+1)$-dimensional Minkowski spacetimes (or Minkowski type spacetimes as in Example~\ref{ExGSGMink}), and more generally on stationary asymptotically flat (with any mass $\bhm\in\R$): for such spacetimes, we have $\frac{n-1}{2}=1$ and $S|_{\pa X}=0$. The vanishing of $f_{1,0}$ is due to the coincidence of the $\rho^{\beta^+}=\rho^1=r^{-1}$ (with $\beta^+=1$) decay produced by the zero energy operator inverse and the $\rho^{\frac{n-1}{2}+S|_{\pa X}}=\rho^1=r^{-1}$ decay of outgoing spherical waves.
\end{rmk}

Setting $k=\lceil\beta^+-\beta^-\rceil$, we can iteratively construct $u_j\in\cA^{\beta^+-j}(X;E_X)$, $0\leq j<k$, as in the proof of Proposition~\ref{PropStCoReg}; for $j\geq 1$, we have
\begin{equation}
\label{EqASuj}
  u_j = u_{0 0}\chi\rho^{\lambda^+-j}v_j + \tilde u_j,\qquad
  v_j=\wh N(\rho^{-2}\hat P(0),-i(\lambda^+-j))^{-1}f_{j,0},\quad
  \tilde u_j\in\cA^{\beta^+-j+\eps}(X;E_X),
\end{equation}
and then $\hat f_{j+1}=\sigma^{-1}(\hat f_j(\sigma)-\hat P(\sigma)u_j)$ (which matches~\eqref{EqStCoRegfj} with $j$ increased by $1$, and with $\psi$ absent) satisfies
\[
  \hat f_{j+1} \in \CI\bigl(\pm[0,1)_\sigma;\cA^{(\beta^+-(j+1)+2,0),\beta^+-(j+1)+2+\eps}(X;E_X)\bigr)
\]
and
\begin{align*}
  &\hat f_{j+1}(0) = u_{0 0}\chi\rho^{\lambda^+-(j+1)+2}f_{j+1,0} + \tilde f_{j+1}, \\
  &\qquad f_{j+1,0}=-2 i\Bigl(\lambda^+-j-\frac{n-1}{2}-S|_{\pa X}\Bigr)v_j,\quad \tilde f_{j+1}\in\cA^{\beta^+-(j+1)+2+\eps}(X;E_X).
\end{align*}
We require that the leading order term $f_{j+1,0}$ be nonzero for $j=k-1$ (and thus necessarily for all lower $j$ as well); that is,
\begin{equation}
\label{EqASAssmfLead}
  \parbox{0.8\textwidth}{$f_{k,0}\neq 0\in\CI(\pa X;E|_{\pa X})$, where $f_{j,0}$ is defined for $j=1$ by~\eqref{EqASf10}, and for $j=1,\ldots,k-1$ inductively by\[f_{j+1,0}=-2 i\Bigl(\lambda^+-j-\frac{n-1}{2}-S|_{\pa X}\Bigr)\wh N(\rho^{-2}\hat P(0),-i(\lambda^+-j))^{-1}f_{j,0}.\]}
\end{equation}
We note that this assumption holds if there exists $p\in\pa X$ so that $\lambda^+-j\notin\frac{n-1}{2}+\spec S(p)$ for $j=0,\ldots,k-1$: the nonvanishing of $f_{1,0}$ under this assumption (for $j=0$) was already discussed above, and the nonvanishing of $f_{j+1,0}$ follows inductively.

\begin{rmk}[Counterexamples: II]
\label{RmkASCounterex2}
  On $(n+1)$-dimensional Minkowski (type) spacetimes with odd $n\geq 3$ (in which case $(\beta^-,\beta^+)=(0,n-2)$ and $\lambda^+=n-2$, so $k=\lceil\beta^+-\beta^-\rceil=n-2$ by Proposition~\ref{PropES} below), with the zero energy inverse producing $\rho^{n-2}$ decay and outgoing spherical waves having $\rho^{\frac{n-1}{2}}$ decay, one always has $f_{j,0}=0$ for $j=\frac{n-1}{2}$ (which is the value of $j$ for which $(n-2)-(j-1)-\frac{n-1}{2}=0$). In general, whenever $f_{j,0}=0$ for some $j\in\{0,\ldots,k\}$, one can construct the Taylor series of $\hat P(\sigma)^{-1}\hat f(\sigma)$ at $\zface$ to at least one higher order than done here, which results in higher regularity at $\zface$ and thus in faster pointwise decay of $u$. In such cases, the extraction of a leading asymptotic term of $u$ requires keeping track also of subleading terms of some of the $u_j$, as done in \cite{TataruDecayAsympFlat,HintzPrice,MorganDecay,MorganWunschPrice}.
\end{rmk}

\textbf{Transition face model problem.} Under the assumption~\eqref{EqASAssmfLead}, the error term at zero energy $\hat f_k(0)\equiv u_{0 0}\chi\rho^{\lambda^+-k+2}f_{k,0}\bmod\cA^{\lambda^+-k+2+\eps}(X;E_X)$ does not lie in the codomain of the zero energy operator in Lemma~\ref{LemmaStEst0} for any $\beta\in(\beta^-,\beta^+)$; the next step in the resolvent construction thus requires the inversion of the $\tface^\pm$-normal operators. To this end, we note that in $\pm\sigma\geq 0$, and with $\hat r=\hat\rho^{-1}=\frac{|\sigma|}{\rho}$, we have
\begin{equation}
\label{EqASfppm}
  f'_\pm := u_{0 0}^{-1}\bigl(|\sigma|^{-(\lambda^+-k+2)}\hat f_k(\sigma)\bigr)\big|_{\tface^\pm} = \hat r^{-(\lambda^+-k)-2}f_{k,0} \in \cA^{((\beta^+-k)+2,-(\beta^+-k)+2)}(\tface^\pm;E|_{\pa X}),
\end{equation}
where the orders refer to decay at $\sctface^\pm=\tface^\pm\cap\scface^\pm=\hat\rho^{-1}(0)$ and $\ztface^\pm=\tface^\pm\cap\zface=\hat r^{-1}(0)$, respectively. Using Lemma~\ref{LemmaStEsttf}\eqref{ItStEsttfb} and Sobolev embedding, and noting that $\beta^+-k\leq\beta^-$, there exists a unique $u'_\pm\in\cA^{(\ell,-\beta^--)}(\tface^\pm;E|_{\pa X})$ (with $\ell\in\R$ sufficiently negative) with
\begin{equation}
\label{EqASNtf}
  N_\tface^\pm(P) u'_\pm=f'_\pm.
\end{equation}
In order to describe the asymptotic behavior of $u'_\pm$ at $\hat r=0$, we now take advantage of assumption~\eqref{EqASAssmSimple} for the `$-$' sign. In the case of strict inequality $-(\beta^+-k)>-\beta^-$ (which happens if and only if $\beta^+-\beta^-\notin\N$), we have
\begin{subequations}
\begin{equation}
\label{EqASuprime}
  u'_\pm = u'_{\pm,0}\psi(\hat r)\hat r^{-\lambda^-} v_- + \tilde u_\pm',\qquad \tilde u_\pm'\in\cA^{(\ell,-\beta^-+\eps)}(\tface^\pm;E|_{\pa X}),
\end{equation}
for some $u'_{\pm,0}\in\C$ and with $\eps>0$ small enough so that no element of $\specb(\rho^{-2}\hat P(0))$ has real part in $(\beta^--2\eps,\beta^-)$. Here, $\psi\in\CIc(\R)$ is identically $1$ near $0$. If on the other hand $-(\beta^+-k)=-\beta^-$, then there are two possibilities. Either we have $-(\lambda^+-k)\neq-\lambda^-$, in which case $u'_\pm$ has such an expansion with an additional term $\psi(\hat r)\hat r^{-(\lambda^+-k)}u'_{\pm,1}$ for some $u'_{\pm,1}\in\CI(\pa X;E|_{\pa X})$; we shall not permit this situation here for simplicity of presentation. Or we have $-(\lambda^+-k)=-\lambda^-$, in which case
\begin{equation}
\label{EqASuprimeLog}
  u'_\pm = u'_{\pm,0}\psi(\hat r)\hat r^{-\lambda^-}(\log\hat r)v_- + \psi(\hat r)\hat r^{-\lambda^-}u'_{\pm,1} + \tilde u'_\pm,
\end{equation}
\end{subequations}
again with $u'_{\pm,0}\in\C$, $u'_{\pm,1}\in\CI(\pa X;E|_{\pa X})$, and with $\tilde u'_\pm$ as in~\eqref{EqASuprime}. The  final assumption (which is stronger than~\eqref{EqASAssmfLead}) then is:
\begin{equation}
\label{EqASAssmtfLead}
  \parbox{0.8\textwidth}{Assume that either $\beta^+-\beta^-\notin\N$ or $\lambda^+-k=\lambda^-$. For $f'_\pm=f'_\pm(\hat r,\omega)=\hat r^{-(\lambda^+-k)-2}f_{k,0}$, the conormal solutions $u'_\pm$ of $N_\tface^\pm(P)u'_\pm=f'_\pm$ on $\tface^\pm$, which are of the form~\eqref{EqASuprime} or \eqref{EqASuprimeLog}, have leading coefficient $u'_{\pm,0}\neq 0$.}
\end{equation}
(In the case $\lambda^+-k=\lambda^-$, the coefficients $u'_{\pm,0}$ can be computed easily; they are the unique complex numbers so that $f_{k,0}-u'_{\pm,0}\hat r^{\lambda^-+2}[N_\tface^\pm(P),\log\hat r](\hat r^{-\lambda^-}v_-)\in\CI(\pa X;E|_{\pa X})$ lies in the range of $\wh N(\rho^{-2}\hat P(0),-i\lambda^-)$. See also \cite[Lemma~2.23]{HintzPrice}.) If~\eqref{EqASAssmtfLead} is violated, then the low energy resolvent is more regular at $\zface$ than what we shall prove below under the hypothesis~\eqref{EqASAssmtfLead}, and this leads to stronger decay of waves at $\cT^+$.

Note now that there exists a unique $u_-\in\cA^{\beta^-}(X;E_X)$ so that
\begin{equation}
\label{EqASLargeZero}
  \hat P(0)u_-=0,\qquad \exists\,\eps>0\ \ \text{so that}\ \ u_-=\chi\rho^{\lambda^-}v_-+\tilde u_-,\quad \tilde u_-\in\cA^{\beta^-+\eps}(X;E_X).
\end{equation}
Indeed, one can set $\tilde u_-=-\hat P(0)^{-1}f_-$ where $f_-=\hat P(0)(\chi\rho^{\lambda^-}v_-)\in\cA^{\lambda^-+2+\delta}(X;E_X)$; note that the normal operator of $\hat P(0)$ annihilates $v_-$, which gives the extra order of decay of $f_-$ at $\pa X$. (One might reasonably call $u_-$ a \emph{large zero energy state}.)

\textbf{Final error.} Let us write
\[
  N(\hat P(0))=\rho^2 N(\rho^{-2}\hat P(0)),\qquad
  N(\pa_\sigma\hat P(0))=\rho N(\rho^{-1}\pa_\sigma\hat P(0));
\]
then we have $N(\hat P(0))\equiv\hat P(0)\bmod\cA^{2+\delta}\Diffb^2(X;E_X)$ and $N(\pa_\sigma\hat P(0))\equiv\pa_\sigma\hat P(0)\bmod\cA^{1+\delta}\Diffb^1(X;E_X)$. Moreover, we have $\sigma^2 N_\tface^\pm(P)=N(\hat P(0))+\sigma N(\pa_\sigma\hat P(0))$ as differential operators on $(-1,1)_\sigma\times[0,1)_\rho\times\pa X$.

Consider first the case that $\beta^+-\beta^-\notin\N$, so $u'_\pm$ is given by~\eqref{EqASuprime}. Regarding the error term $\tilde u'_\pm$ as a function of $(\sigma,\rho,\omega)$ for $\rho<1$, we now define an extension of $u'_\pm$ off $\tface^\pm$ by replacing $v_-$ (up to weight factors) by the large zero energy state $u_-$ from~\eqref{EqASLargeZero}, to wit,
\[
  u''_\pm := \chi u'_{\pm,0}\psi\hat r^{-\lambda^-}v_-+\chi\tilde u'_\pm + |\sigma|^{-\lambda^-}\psi u'_{\pm,0}\tilde u_- = |\sigma|^{-\lambda^-}\bigl(\psi u'_{\pm,0}u_- + |\sigma|^{\lambda^-}\chi\tilde u'_\pm \bigr).
\]
Thus, $u''_\pm\in\cA^{(\ell,\ ((0,0),\eps),\ ((-\lambda^-,0),-\lambda^-+\eps))}(X_\res^\pm)$ for some $\ell\in\R$ and $\eps>0$. We then set
\begin{align}
\label{EqASfkp1}
  \hat f_{k+1}(\sigma) &:= \hat f_k(\sigma) - u_{0 0}|\sigma|^{\lambda^+-\lambda^--k} \hat P(\sigma)\bigl( \psi u'_{\pm,0} u_- + |\sigma|^{\lambda^-}\chi\tilde u'_\pm\bigr) \\
    & = u_{0 0}|\sigma|^{\lambda^+-k} \bigl( u_{0 0}^{-1}\sigma^2|\sigma|^{-\lambda^++k-2}\hat f_k(\sigma) - \hat P(\sigma)u''_\pm\bigr) \nonumber\\
    & = u_{0 0}|\sigma|^{\lambda^+-k} \Bigl[ \sigma^2\Bigl( u_{0 0}^{-1}|\sigma|^{-\lambda^++k-2}\hat f_k(\sigma)-\chi N_\tface^\pm(P)u''_\pm\Bigr) \nonumber\\
    &\quad\hspace{5.5em} + \bigl(\chi N(\hat P(0))-\hat P(0)\bigr)u''_\pm + \bigl(\chi \sigma N(\pa_\sigma\hat P(0))-\sigma\pa_\sigma\hat P(0)\bigr)u''_\pm \nonumber\\
    &\quad\hspace{5.5em} - \half\sigma^2\pa_\sigma^2\hat P(0)u''_\pm \Bigr]. \nonumber
\end{align}
The first summand in square parentheses vanishes to leading order at $\tface^\pm$ by construction, and thus lies in $\cA^{(\ell,2+\eps,((-\lambda^-,0),-\lambda^-+\eps))}(X_\res^\pm)$; from now on $\ell$ may vary from line to line. The second and third summands gain $\delta$ orders of decay at $\tface^\pm$ (and $\scface^\pm$) and thus lie in $\cA^{(\ell,2+\delta,((-\lambda^-,0),-\lambda^-+\eps))}(X_\res^\pm)$. The final term lies in $\cA^{(\ell,3+\delta,-\lambda^-+2)}(X_\res^\pm)$. Thus far, we have shown that $\hat f_{k+1}\in|\sigma|^{\lambda^+-k}\cA^{(\ell,\ 2+\eps,\ ((-\lambda^-,0),-\lambda^-+\eps))}(X_\res^\pm)$ (if we choose $\eps<\delta$, as we may). However, the restriction of $|\sigma|^{-\lambda^++\lambda^-+k}\hat f_{k+1}(\sigma)$ to $\zface$ is given by $-u_{0 0}u'_{\pm,0}\hat P(0)u_-=0$ (using that $-\lambda^++\lambda^-+k$ has positive real part). Therefore,
\[
  \hat f_{k+1}(\sigma) \in |\sigma|^{\lambda^+-\lambda^--k+\frac{\eps}{2}}\cA^{(\ell,\lambda^-+2+\frac{\eps}{2},0)}(X_\res^\pm) = \cA^{(\ell,\ \lambda^+-k+2+\eps,\ \lambda^+-\lambda^--k+\frac{\eps}{2})}(X_\res^\pm).
\]
By Lemma~\ref{LemmaStCo}\eqref{ItStCoLo}, and noting that $\Re(\lambda^-+2+\frac{\eps}{2})\in(\beta^-+2,\beta^++2)$, we have
\begin{equation}
\label{EqASukp1}
  u_{k+1}(\sigma) := \hat P(\sigma)^{-1}\hat f_{k+1}(\sigma) \in \cA^{(\ell,\ \lambda^+-k+\eps,\ \lambda^+-\lambda^--k+\frac{\eps}{2})}(X_\res^\pm).
\end{equation}

Putting the pieces~\eqref{EqASu0}, \eqref{EqASuj}, \eqref{EqASfkp1}, and \eqref{EqASukp1} together, we have
\begin{align}
  \hat u_{\rm lo}(\sigma) &:= \chi(\sigma)\hat P(\sigma)^{-1}\hat f(\sigma) \nonumber\\
\label{EqAShatu}
\begin{split}
    &= \chi(\sigma)\Biggl(\;\sum_{j=0}^{k-1} \Bigl(\frac{\sigma}{\rho}\Bigr)^j \rho^j u_j \\
    &\hspace{4em} + u_{0 0}\sum_\pm H(\pm\sigma)(\pm 1)^k|\sigma|^{\lambda^+-\lambda^-}\bigl(\psi u'_{\pm,0}u_- + |\sigma|^{\lambda^-}\chi\tilde u'_\pm\bigr) + \sigma^k u_{k+1}\Biggr),
\end{split}
\end{align}
where $H$ is the Heaviside function. (Here $(\pm 1)^k=|\sigma|^{-k}\sigma^k$.) In particular,
\[
  \hat u_{\rm lo}|_{\pm\sigma\geq 0} \in \cA^{(\ell,\ ((\lambda^+,0),\beta^++\eps),\ ((0,0),\beta^+-\beta^-))}(X_\res^\pm).
\]
The strongest singularity at $\zface$ arises from the term $\sum_\pm H(\pm\sigma)(\pm 1)^k|\sigma|^{\lambda^+-\lambda^-}u'_{\pm ,0}u_-$ (since $|\sigma|^{\lambda^+-\lambda^-}|\sigma|^{\lambda^-}\chi\tilde u'_\pm\in\cA^{(\ell,\beta^+,\beta^+-\beta^-+\eps)}(X_\res^\pm)$ is more regular at $\zface$); thus, this term produces the leading order asymptotics of $u_{\rm lo}=\cF^{-1}\hat u_{\rm lo}$ and thus of $u$ at $\cT^+$. Writing $|\sigma|^{\lambda^+-\lambda^-}$ as a linear combination of $(\sigma\pm i 0)^{\lambda^+-\lambda^-}$, and noting that the $|t_*|^{-(\lambda^+-\lambda^-+1)}$ lower bound of $\cF^{-1}(\sigma-i 0)^{\lambda^+-\lambda^-}$ in $t_*\leq -1$ cannot be canceled by the inverse Fourier transform of any other term in~\eqref{EqAShatu}, we conclude that the strongest singular term is in fact $(\sigma+i 0)^{\lambda^+-\lambda^-}u'_{+,0}u_-$. (This argument is taken from \cite[Remark~3.5]{HintzPrice}.) Its inverse Fourier transform is a nonzero multiple of $t_*^{-(\lambda^+-\lambda^-+1)}u_-$ for $t_*>0$, which is thus the leading order asymptotic profile of $u$ at $\cT^+$.

In order to control $u$ also at $\iota^+$, we argue in a manner similar to \cite[\S3.2.2]{HintzPrice} and note that $\hat u_{\rm lo}$ is partially polyhomogeneous on $X_\res^\pm$. We just discussed its leading order term at $\zface$; for control of $u$ at $\iota^+$, we now turn to its leading order term at $\tface$, which in terms of the function $\hat r_{\rm tot}:=\frac{\sigma}{\rho}=\pm\hat r$ on $\tface=\bigcup_\pm\tface^\pm$ is $u_{0 0}$ times
\[
  \hat u_\tface=\sum_{j=0}^{k-1}\Big(\frac{\sigma}{\rho}\Bigr)^j \rho^{\lambda^+}v_j + |\sigma|^{\lambda^+}\sum_\pm H(\pm\hat r_{\rm tot})(\pm 1)^k u'_\pm=:\rho^{\lambda^+}\hat u_{\tface,0};
\]
here, we set $v_0=v^+$ and recall $v_j\in\CI(\pa X;E|_{\pa X})$ for $j\geq 1$ from~\eqref{EqASuj}, and we recall $u'_\pm$ from~\eqref{EqASfppm}--\eqref{EqASNtf}. Thus, the distribution $\hat u_{\tface,0}=\hat u_{\tface,0}(\hat r_{\rm tot},\omega)$ on $\tface$ is given by
\begin{equation}
\label{EqASutf0}
\begin{split}
  \hat u_{\tface,0}(\hat r_{\rm tot},\omega) &= \sum_{j=0}^{k-1} \hat r_{\rm tot}^j v_j + \sum_\pm H(\pm\hat r_{\rm tot})(\pm 1)^k |\hat r_{\rm tot}|^{\lambda^+}u'_\pm \\
    &= \sum_{j=0}^{k-1} \hat r_{\rm tot}^j v_j + \sum_\pm \bigl(\hat r_{\rm tot}^{\lambda^++2}N_\tface^\pm(P)\hat r_{\rm tot}^{-\lambda^+}\bigr)^{-1}(\hat r_{\rm tot}^k f_{k,0}),
\end{split}
\end{equation}
where we use $\hat r=\pm\hat r_{\rm tot}$ in the expression~\eqref{EqGSONtf} of $N_\tface^\pm(P)$. Thus, $\hat u_{\tface,0}$ is the sum of a polynomial in $\hat r_{\rm tot}$ and a term which is conormal at $\hat r_{\rm tot}^{-1}(0)\subset\tface\subset X_\res$ and has a nonvanishing leading order term at $\hat r_{\rm tot}=0$ with decay rate $\lambda^+-\lambda^-$ in view of~\eqref{EqASuprime} and assumption~\eqref{EqASAssmtfLead}. The inverse Fourier transform of $\hat u_\tface$ in $\sigma=\rho\hat r_{\rm tot}$ is
\[
  \rho^{\lambda^++1}(\cF^{-1}_{v\to\hat r_{\rm tot}}\hat u_{\tface,0})(v,\omega)|_{v=\rho t_*}=t_*^{-(\lambda^+-\lambda^-+1)}r^{-\lambda^-}\cdot v^{\lambda^+-\lambda^-+1}(\cF^{-1}_{v\to\hat r_{\rm tot}}\hat u_{\tface,0})(v,\omega)|_{v=\rho t_*}.
\]
Since $\iota^+=[0,\infty]_v\times\pa X\subset M_1$ for $v=\rho t_*$, this (and recalling~\eqref{EqASfppm}--\eqref{EqASNtf}) identifies the asymptotic profile $a_+:=(t_*^{\lambda^+-\lambda^-+1}r^{\lambda^-}u)|_{\iota^+}$ of $u$ at $\iota^+$, as a distribution on $(0,\infty)_v\times\pa X$, as $u_{0 0}$ times
\begin{equation}
\label{EqASiplus}
\begin{split}
  &a_+ = v^{\lambda^+-\lambda^-+1}a_{+,0}, \\
  &\qquad a_{+,0}=\bigl(\cF^{-1}_{v\to\hat r_{\rm tot}}\hat u_{\tface,0}\bigr)\big|_{v>0} = \cF_{v\to\hat r_{\rm tot}}^{-1}\biggl( \sum_\pm \bigl(\hat r_{\rm tot}^{\lambda^++2}N_\tface^\pm(P)\hat r_{\rm tot}^{-\lambda^+}\bigr)^{-1}(\hat r_{\rm tot}^k f_{k,0})\biggr).
\end{split}
\end{equation}
Note here that the inverse Fourier transform of the first, polynomial, term in~\eqref{EqASutf0} is a sum of differentiated $\delta$-distributions supported at $v=0$ (which is the boundary of $\iota^+$ at null infinity), and thus its restriction to $v>0$ vanishes. The function $a_+(v,\omega)$ is a nonzero function of $v=\rho t_*\in(0,\infty)$ since $\hat u_{\tface,0}$ is not a polynomial in $\hat r_{\rm tot}$ (cf.\ its behavior near $\hat r_{\rm tot}=0$).

We leave a detailed discussion of the case $k=\lambda^+-\lambda^-\in\N$ (and thus $u'_\pm$ is of the form~\eqref{EqASuprimeLog}) to the reader; in a special case, this is discussed in \cite[\S3.1]{HintzPrice}. The main difference is that now the strongest singularity of $\hat u_{\rm lo}(\sigma)$ at $\zface$ is logarithmic, and indeed equal to $(\sigma+i 0)^k\log(\sigma+i 0)u_{+,0}u_-$ with $u_{+,0}\in\C$, and $u_-$ is the large zero energy state $u_-$ from~\eqref{EqASLargeZero}.

\begin{rmk}[Support of the profile at $\iota^+$]
\label{RmkASipSupp}
  The distribution $\cF_{v\to\hat r_{\rm tot}}^{-1}\hat u_{\tface,0}$ appearing in~\eqref{EqASiplus} is defined on $\R_v\times\pa X$; but it necessarily vanishes in $v<0$ since $u$, being a forward solution with forcing supported in $t_*\geq 0$, must vanish for $t_*<0$, and in particular near $I^0\supset\{\rho=0,v<0\}$. (Note also that the inverse Fourier transform of $\hat u_{\rm lo}$ after subtraction of its leading order part at $\tface$ and upon localizing near $\scface^\pm$ is the inverse Fourier transform, from $\sigma$ to $t_*$, of an element of $\cA^{(\ell,\beta^++\eps,\infty)}(X_\res^\pm)$, which has $\rho^{\beta^++1+\eps}=o(\rho^{\beta^++1})$ decay at $\rho=0$ for $v\neq 0$ by \cite[Proposition~2.29(1)]{Hintz3b}, and thus does not contribute to the profile $a_+$.)
\end{rmk}

\begin{rmk}[Alternative characterization of the profile at $\iota^+$]
\label{RmkASipProfile}
  Using~\eqref{EqGSONtf}, the operator on $\tface$ given by $\hat r^2N_\tface^\pm(P)$ on $\tface^\pm$ takes the form
  \[
    \hat r_{\rm tot}^2 N_\tface(P) = -2 i\hat r_{\rm tot}\Bigl( \hat r_{\rm tot}\pa_{\hat r_{\rm tot}}+\frac{n-1}{2}+S|_{\pa X}\Bigr) + P_{(0)}(0,\omega,-\hat r_{\rm tot} D_{\hat r_{\rm tot}},D_\omega).
  \]
  Conjugation by $\hat r^{\lambda^+}$ transforms $\hat r_{\rm tot}\pa_{\hat r_{\rm tot}}$ to $\hat r_{\rm tot}\pa_{\hat r_{\rm tot}}-\lambda^+$. Thus, by~\eqref{EqASiplus}, and noting that $\cF_{v\to\hat r_{\rm tot}}^{-1}$ intertwines $\hat r_{\rm tot}$ and $i\pa_v$, and $\pa_{\hat r_{\rm tot}}$ and $i v$ (see~\eqref{EqWipmFTvr} below for the normalizations used here), we conclude that $a_+=v^{\lambda^+-\lambda^-+1}a_{+,0}$ where $a_{+,0}$ is a distribution on $\R_v\times\pa X$ which for $k=\lceil\Re(\lambda^+-\lambda^-)\rceil$ solves the equation
  \begin{subequations}
  \begin{equation}
  \label{EqASipProfile}
  \begin{split}
    &\biggl(-2\pa_v \Bigl( v\pa_v + \lambda^+ + 1 - \frac{n-1}{2} - S|_{\pa X}\Bigr) \\
    &\hspace{5em} + P_{(0)}(0,\omega,v D_v-i(\lambda^++1),D_\omega)\biggr)a_{+,0}(v,\omega) = f_{k,0} i^k\delta^{(k)}(v).
  \end{split}
  \end{equation}
  Moreover, $a_{+,0}$ is the \emph{unique} such solution which has the additional property that
  \begin{equation}
  \label{EqASipProfile2}
    \parbox{0.8\textwidth}{$a_{+,0}$ is conormal at $v=0,\infty$, has an upper bound $C|v^{-\lambda^++\lambda^--1}|$ as $v\to\infty$ for some $C>0$, is supported in $v\geq 0$, and satisfies $\la a_{+,0},1\ra=0$,}
  \end{equation}
  \end{subequations}
  where $\la-,-\ra$ is the distributional pairing. Indeed, the upper bound and support property of $a_{+,0}$ ensure that the Fourier transform $\wh{a_{+,0}}(\hat r_{\rm tot},\omega)$ is conormal (in fact, polyhomogeneous by virtue of solving~\eqref{EqASipProfile}) at $\hat r_{\rm tot}=0$ and everywhere continuous, and the final condition then enforces the vanishing of $\wh{a_{+,0}}(\hat r_{\rm tot},\omega)$ at $\hat r_{\rm tot}=0$. Thus, on $\tface^\pm$, the rescaling $\hat r^{-\lambda^+}\wh{a_{+,0}}(\hat r,\omega)$ is $\cO(\hat r^{-\lambda^++\eps})$ for some $\eps>0$ as $\hat r\searrow 0$ and therefore lies in a space on which $N_\tface^\pm(P)$ is injective.\footnote{We caution the reader not to conflate the distribution $a_{+,0}$ on $\R_v\times\pa X$ with $\cF_{v\to\hat r_{\rm tot}}^{-1}\hat u_{\tface,0}$: these two distributions differ by the sum of differentiated $\delta$-distributions at $v=0$ caused by the first summand in~\eqref{EqASutf0}.} --- Lastly, assumption~\eqref{EqASAssmtfLead} is equivalent to $a_+=v^{\lambda^+-\lambda^-+1}a_{+,0}$ having a non-trivial restriction to $v=\infty$, i.e.\ there exists $c\neq 0$ so that $a_+(v,\omega)-c v_-(\omega)$ has $v^{-\eps}$ decay as $v\to\infty$ (or equivalently $a_{+,0}(v,\omega)-c v^{-\lambda^++\lambda^--1}v_-(\omega)=\cO(v^{-\lambda^++\lambda^--1-\eps})$).
\end{rmk}

We summarize our analysis as follows:

\begin{thm}[Asymptotic profiles]
\label{ThmAS}
  Let $P$ be a stationary wave type operator which is spectrally admissible with indicial gap $(\beta^-,\beta^+)$, and define $S\in\CI(X;\End(E_X))$ via~\eqref{EqGSOStruct}. Define $\ubar S\in\R$ as in Definition~\usref{DefStEstThr}. Suppose that $P$ satisfies the assumption~\eqref{EqASAssmSimple}. Let $f\in\cF\setminus\cF_0$ (see~\eqref{EqASSpace}--\eqref{EqASSubspace}), and suppose that assumption~\eqref{EqASAssmtfLead} (which uses the notation of~\eqref{EqASAssmfLead}) is satisfied. Then there exist $\eps>0$ and
  \[
    a\in\cA^{(\frac{n-1}{2}+\ubar S-\lambda^--,\ ((0,0),\eps),\ ((0,0),\eps))}(M_1;E)
  \]
  so that neither $a_\cT=a|_{\cT^+}\in\CI(X;E_X)$ nor $a_+=a|_{\iota^+}\in\CI(\iota_+;E)$ are identically $0$, and so that for all $\eta>0$ there exists $C>0$ so that (in the notation of Definition~\usref{DefGAMfd})
  \begin{equation}
  \label{EqAS}
    \Bigl|u(t_*,x) - a(t_*,x) t_*^{-(\lambda^+-\lambda^-+1)}\la x\ra^{-\lambda^-}\Bigr|\leq C\rho_{\!\scri}^{\frac{n-1}{2}+\ubar S-\eta}\rho_+^{\beta^++1+\eps}\rho_\cT^{\beta^+-\beta^-+1+\eps}
  \end{equation}
  for $t_*\geq 1$; such estimates also hold for derivatives of the term in absolute values on the left along any b-differential operator on $M_1$ (i.e.\ along any number of $t_*\pa_{t_*}$, $\la x\ra\pa_x$). Here, $a_\cT$ is a constant multiple of the rescaled large zero energy state $\la x\ra^{\lambda^-}u_-$, with $u_-$ as in~\eqref{EqASLargeZero}, and $a_+$ is given by~\eqref{EqASiplus}. If assumption~\eqref{EqASAssmtfLead} is violated, then~\eqref{EqAS} holds with $a\equiv 0$.
\end{thm}

Thus, Theorem~\ref{ThmAS} gives conditions under which the upper bounds in Theorem~\ref{ThmStCo} are sharp up to $\eps$-losses.

\begin{proof}[Proof of Theorem~\usref{ThmAS}]
  Only the decay at $\scri^+$ was not discussed above. This can be proved by similar arguments as at the end of the proof of Theorem~\ref{ThmStCo}, using \cite[Lemmas~7.6 and 7.7]{HintzVasyMink4} to handle the leading order term of $u$ at $\iota^+$.
\end{proof}

The conditions in Theorem~\ref{ThmAS} can easily be checked in a number of settings; see~\S\ref{SE} for examples.

Note that for any spectrally admissible stationary wave type operator $P$ there exists an open neighborhood $\cU$ of $P$ in the space of stationary wave type operators so that all $P'\in\cU$ are also spectrally admissible (with indicial gaps depending on $P'$). We then conjecture that inside this open set of spectrally admissible stationary wave type operators, there is an open dense subset of operators for which the conditions of Theorem~\ref{ThmAS} are satisfied for all $f\in\cF$ outside a subspace of positive codimension; in this sense, the upper bounds given in Theorem~\ref{ThmStCo} are, conjecturally, generically sharp up to the arbitrarily small losses.

\section{Analysis of admissible wave type operators}
\label{SW}

We now turn to the analysis of non-stationary wave type operators. We use the notation of~\S\ref{SsGA}, and fix an $(\ell_0,2\ell_{\!\scri},\ell_+,\ell_\cT)$-admissible asymptotically flat metric $g$ (relative to a stationary and asymptotically flat Lorentzian metric $g_0$ as in Definition~\ref{DefGSG}); see Definition~\ref{DefGAG}. We moreover fix an admissible wave type operator $P$ (with respect to $g$ and a stationary wave type operator $P_0$ relative to $g_0$---see Definition~\ref{DefGSO}) acting on sections of the pullback $\upbeta^*E\to M$ of a stationary vector bundle $E\to M_0$ (see~\S\ref{SsGS}); see Definition~\ref{DefGAW}.

We first establish some basic properties of the characteristic set and the null-bi\-char\-ac\-ter\-is\-tic flow of $P$. In~\S\ref{SsWT}, we translate the saddle point structure of two radial sets (incoming and outgoing) over $\pa\cT^+$ into radial point estimates (which are microlocal estimates in weighted 3b-Sobolev spaces near $\cT^+\cap\iota^+$). In~\S\ref{SsWR}, we combine these estimates with the microlocal propagation estimates near $\scri^+$ from~\cite{HintzVasyScrieb} to prove a global regularity estimate on edge-3b-spaces (edge-b near $\scri^+$, 3b near $\cT^+$); see Proposition~\ref{PropWR}. Following an analysis of the (Mellin-transformed) normal operator at $\iota^+$ in~\S\ref{SsWip}, we complement the regularity estimate from~\S\ref{SsWR} with the invertibility of the various normal operators (at $\scri^+$ from \cite{HintzVasyScrieb}, $\iota^+$ from~\S\ref{SsWip}, and $\cT^+$ from~\S\ref{SsStEst}) to obtain Fredholm and invertibility properties of admissible wave type operators (see Theorem~\ref{ThmWM}). In~\S\ref{SsWb} finally, we prove higher (b-)regularity for forward solutions of admissible wave type equations.

\begin{notation}[Phase space]
  As already done in parts of~\S\ref{SsGA}, we shall, when working in (local coordinates in) a neighborhood of $\cT^+$, denote the $\betbop$-phase space by $\Ttb^*M$ simply. Similarly, over a neighborhood of $\scri^+$, we shall write $\Teb^*M$; and on a set which contains $\iota^+$ but is disjoint from $I^0$, we shall write $\Tetb^*M$. This is a suggestive (and somewhat imprecise) notation which is easier to parse than $\Tbetb^*_\cU M$ with specifications of the set $\cU$, and it makes explicit the Lie algebra structure on $M$ of main interest at any given step of the argument.
\end{notation}

By definition, the principal symbol $p\colon\zeta\mapsto g^{-1}(\zeta,\zeta)$ of $P$ satisfies
\[
  p \in \rho_0^2 x_{\!\scri}^2\rho_+^2(\CI+\cA^{(\ell_0,2\ell_{\!\scri},\ell_+,\ell_\cT)})P^2(\Tbetb^*M)
\]
on $\{t_*+r>-\half r-1\}$, where we recall that $P^2(\Tbetb^*M)\subset S^2(\Tbetb^*M)$ consists of fiber-wise homogeneous quadratic polynomials. We denote by $\Sigma\subset\Tbetb^*M\setminus o$ the characteristic set of $P$, which is the closure of $p^{-1}(0)\cap T^*M^\circ\setminus o$. Note that $\Sigma$ has two components over $\ft_*\geq -1$,
\[
  \Sigma\cap\Tetb^*_{\{\ft_*\geq -1\}} M = \Sigma^+ \cup \Sigma^-,
\]
where $\Sigma^\pm$ is the closure of the set of $\zeta\in T^*M^\circ$ for which $\pm g^{-1}(\zeta,-\dd\ft_*)<0$. Note that Definition~\ref{DefGAG} and Lemma~\ref{LemmaGAGbetb} imply that in the interior of $\cT^+$, this definition is equivalent to the requirement that $\zeta\in\Sigma\cap\Ttb^*_{(\cT^+)^\circ}M$ satisfy $\pm g^{-1}(\zeta,-\dd t_*)<0$, while near $\iota^+$ and $\scri^+$ the requirement reads $\pm g^{-1}(\zeta,-r\,\dd t)<0$ where $t=t_*+r$.

Near null infinity, $P$ is an admissible operator in the sense of \cite[Definition~3.5]{HintzVasyScrieb}; see Remark~\ref{RmkGAWComments}\eqref{ItGAWCommentsScriEB}. Therefore, the results on the structure of the null-bicharacteristic flow of $P$ and microlocal regularity results near $\scri^+$ from \cite[\S4]{HintzVasyScrieb} apply here. Concretely, let us use the coordinates $x_{\!\scri},\rho_+$ from~\eqref{EqGAGebp} (with $T=-2$) in a neighborhood $\cU_\scri$ of $\scri^+\cap\{\ft_*\geq -1\}$, and write $\omega\in\R^{n-1}$ for local coordinates on $\Sph^{n-1}$. Writing covectors as
\begin{equation}
\label{EqWHebCovec}
  \xi_\ebop\frac{\dd x_{\!\scri}}{x_{\!\scri}} + \eta_\ebop\frac{\dd\omega}{x_{\!\scri}} + \zeta_\ebop\frac{\dd\rho_+}{\rho_+},
\end{equation}
we then have (only recording weights at $\scri^+$ and $\iota^+$)
\begin{equation}
\label{EqWHebG}
  G_\ebop^+ := 2 x_{\!\scri}^{-2}\rho_+^{-2}p \equiv \xi_\ebop(\xi_\ebop-2\zeta_\ebop)+\half x_{\!\scri}^2\xi_\ebop^2+2|\eta_\ebop|_{\slg^{-1}}^2 \bmod \cA^{(2\ell_{\!\scri},\ell_+)}S^2(\Teb^*_{\cU_\scri}M)
\end{equation}
by Lemma~\ref{LemmaGAGeb}\eqref{ItGAGebp} and the membership~\eqref{EqGAGMetric}. The Hamiltonian vector field
\begin{equation}
\label{EqWHameb}
\begin{split}
  H_{G_\ebop^+} &= 2\bigl((1+\half x_{\!\scri}^2)\xi_\ebop-\zeta_\ebop\bigr)(x_{\!\scri}\pa_{x_{\!\scri}}+\eta_\ebop\pa_{\eta_\ebop}) - 2\xi_\ebop\rho_+\pa_{\rho_+} \\
  &\qquad - (x_{\!\scri}^2\xi_\ebop^2+4|\eta_\ebop|_{\slg^{-1}}^2)\pa_{\xi_\ebop} + 4 x_{\!\scri}\slg^{i j}\eta_{\ebop,i}\pa_{\omega^j} \bmod \cA^{(2\ell_{\!\scri},\ell_+)}\Veb(\Teb^*_{\cU_\scri}M),
\end{split}
\end{equation}
computed using \cite[Lemma~2.2]{HintzVasyScrieb}, is fiber-radial only at the sets
\begin{equation}
\label{EqWRadScri}
\begin{split}
  \cR_{\scri,\rm out}^\pm &:= \{ (x_{\!\scri},\rho_+,\omega;\xi_\ebop,\zeta_\ebop,\eta_\ebop) \colon x_{\!\scri}=0,\ \xi_\ebop=\eta_\ebop=0,\ \pm\zeta_\ebop>0 \}, \\
  \cR_{\scri,\rm in,+}^\pm &:= \{ (x_{\!\scri},\rho_+,\omega;\xi_\ebop,\zeta_\ebop,\eta_\ebop) \colon x_{\!\scri}=\rho_+=0,\ \eta_\ebop=0,\ \xi_\ebop=2\zeta_\ebop,\ \pm\xi_\ebop>0 \},
\end{split}
\end{equation}
with the superscript `$\pm$' indicating containment in $\Sigma^\pm$. (In a sufficiently small neighborhood of $\scri^+$, i.e.\ for small $x_{\!\scri}$, this is shown in \cite[\S4.1]{HintzVasyScrieb}. For general $x_{\!\scri}$, this remains true in view of the fact that $g$ is the Minkowski metric to leading order along all of $\iota^+$, and can be checked using~\eqref{EqWHameb}; it is also a consequence of Lemma~\ref{LemmaWFlowT} below.) We recall furthermore from \cite[\S4.1]{HintzVasyScrieb} that $\cR_{\scri,\rm out}^\pm$ is a (local) sink for the $\pm H_{G_\ebop^+}$-flow, whereas $\cR_{\scri,\rm in,+}^\pm$ is a saddle point, with stable manifold given by $\{(x_{\!\scri},\rho_+,\omega;\xi_\ebop,\zeta_\ebop,\eta_\ebop)\colon x_{\!\scri}=0,\ \eta_\ebop=0,\ \xi_\ebop=2\zeta_\ebop,\ \pm\xi_\ebop>0\}$, and unstable manifold equal to $\Sigma^\pm\cap\Teb^*_{\iota^+}M$ in a neighborhood of $\cR_{\scri,\rm in,+}^\pm$.

Next, in an open neighborhood $\cU_\iota$ of $\iota^+\setminus\scri^+$, consider the coordinates $\rho_+=r^{-1}$, $\rho_\cT=\frac{r}{t_*}$, $\omega\in\R^{n-1}$, and write 3b-covectors as
\begin{equation}
\label{EqWtbCoords}
  \sigma_\tbop\frac{\dd t_*}{r} + \xi_\tbop\frac{\dd r}{r} + \eta_\tbop\,\dd\omega.
\end{equation}
Then Lemma~\ref{LemmaGAGbetb} implies that
\begin{equation}
\label{EqWGtb}
  G_\tbop := \rho_+^{-2}p \equiv \xi_\tbop(\xi_\tbop-2\sigma_\tbop) + |\eta_\tbop|_{\slg^{-1}}^2 \bmod \cA^{(\ell_+,\ ((0,0),\ell_\cT))}P^2(\Ttb_{\cU_\iota}^*M),
\end{equation}
where we record the weights at $\iota^+$ and $\cT^+$. Since, in general,
\begin{align*}
  H_\sfp &= (\pa_{\sigma_\tbop}\sfp) r\pa_{t_*} + (\pa_{\xi_\tbop}\sfp) (r\pa_r+\sigma_\tbop\pa_{\sigma_\tbop}) + (\pa_{\eta_\tbop}\sfp)\pa_\omega \\
    &\qquad - (r\pa_{t_*}\sfp)\pa_{\sigma_\tbop} - \bigl((r\pa_r+\sigma_\tbop\pa_{\sigma_\tbop})\sfp\bigr) \pa_{\xi_\tbop} - (\pa_\omega\sfp)\pa_{\eta_\tbop} \\
      &= -(\pa_{\sigma_\tbop}\sfp) \rho_\cT^2\pa_{\rho_\cT} + (\pa_{\xi_\tbop}\sfp) (\rho_\cT\pa_{\rho_\cT}-\rho_+\pa_{\rho_+}+\sigma_\tbop\pa_{\sigma_\tbop}) + (\pa_{\eta_\tbop}\sfp)\pa_\omega \\
      &\qquad + (\rho_\cT^2\pa_{\rho_\cT}\sfp)\pa_{\sigma_\tbop} - \bigl((\rho_\cT\pa_{\rho_\cT}-\rho_+\pa_{\rho_+}+\sigma_\tbop\pa_{\sigma_\tbop})\sfp\bigr)\pa_{\xi_\tbop} - (\pa_\omega\sfp)\pa_{\eta_\tbop},
\end{align*}
we therefore have
\begin{equation}
\label{EqWHam3b}
\begin{split}
  H_{G_\tbop} &\equiv 2\xi_\tbop\rho_\cT^2\pa_{\rho_\cT} + 2(\xi_\tbop-\sigma_\tbop)(\rho_\cT\pa_{\rho_\cT}-\rho_+\pa_{\rho_+}+\sigma_\tbop\pa_{\sigma_\tbop}) + 2\sigma_\tbop\xi_\tbop\pa_{\xi_\tbop} \\
    &\qquad + 2\slg^{i j}\eta_{\tbop,i}\pa_{\omega^j} - \pa_\omega(|\eta_\tbop|^2_{\slg_\omega^{-1}})\pa_{\eta_\tbop} \bmod \cA^{(\ell_+,\ ((0,0),\ell_\cT))}\Vtb(\Ttb^*_{\cU_\iota}M).
\end{split}
\end{equation}
This is fiber-radial at the sets $\cR_{\cT,\rm in}=\bigcup_\pm\cR_{\cT,\rm in}^\pm$ and $\cR_{\cT,\rm out}=\bigcup_\pm\cR_{\cT,\rm out}^\pm$, where
\begin{equation}
\label{EqWRadT}
\begin{split}
  \cR_{\cT,\rm in}^\pm &:= \{ (\rho_\cT,\rho_+,\omega;\sigma_\tbop,\xi_\tbop,\eta_\tbop) \colon \rho_\cT=\rho_+=0,\ \eta_\tbop=0,\ \xi_\tbop=2\sigma_\tbop,\ \pm\sigma_\tbop<0 \}, \\
  \cR_{\cT,\rm out}^\pm &:= \{ (\rho_\cT,\rho_+,\omega;\sigma_\tbop,\xi_\tbop,\eta_\tbop) \colon \rho_\cT=\rho_+=0,\ \xi_\tbop=\eta_\tbop=0,\ \pm\sigma_\tbop<0 \}.
\end{split}
\end{equation}

We shall write $\pa\cR_{\cT,\rm out}^\pm$ for the intersection of the closure of $\cR_{\cT,\rm out}^\pm$ in $\ol{\Ttb^*_{\cU_\iota}}M$ with fiber infinity $\Stb^*_{\cU_\iota}M$; we define $\pa\cR_{\cT,\rm in}^\pm$, $\pa\Sigma^\pm$, etc.\ similarly. In the local coordinates
\begin{equation}
\label{EqWRadTCoord}
  \rho_\infty := \frac{1}{|\sigma_\tbop|}=\mp\frac{1}{\sigma_\tbop},\qquad
  \hat\xi_\tbop := \frac{\xi_\tbop}{\sigma_\tbop},\qquad
  \hat\eta_\tbop := \frac{\eta_\tbop}{\sigma_\tbop},
\end{equation}
we then compute the linearization of $\rho_\infty H_{G_\tbop}$ at fiber infinity of the radial sets~\eqref{EqWRadT}. For a submanifold $S\subset\Stb^*M$, we write $\cI_S\subset\cA^{(((0,0),\ell_+),((0,0),\ell_\cT))}(\ol{\Ttb^*}M)$ for the ideal of functions vanishing at $S$; changing variables in~\eqref{EqWHam3b} then gives
\begin{align}
\label{EqWRadTHamLinIn}
\begin{split}
  \pm\rho_\infty H_{G_\tbop} &\equiv 2(-\rho_\cT\pa_{\rho_\cT}+\rho_+\pa_{\rho_+}+\rho_\infty\pa_{\rho_\infty}+\hat\eta_\tbop\pa_{\hat\eta_\tbop} - (2-\hat\xi_\tbop)\pa_{\hat\xi_\tbop}) \\
    &\quad\hspace{16em} \bmod \cI_{\pa\cR_{\cT,\rm in}^\pm}\Vb(\ol{\Ttb^*}M),
\end{split} \\
\label{EqWRadTHamLinOut}
\begin{split}
  \pm\rho_\infty H_{G_\tbop} &\equiv 2(\rho_\cT\pa_{\rho_\cT}-\rho_+\pa_{\rho_+}-\rho_\infty\pa_{\rho_\infty}-\hat\eta_\tbop\pa_{\hat\eta_\tbop}-2\hat\xi_\tbop\pa_{\hat\xi_\tbop}) \\
    &\quad\hspace{16em} \bmod \cI_{\pa\cR_{\cT,\rm out}^\pm}\Vb(\ol{\Ttb^*}M).
\end{split}
\end{align}
Thus, $\pa\cR_{\cT,\rm in}^\pm$ and $\pa\cR_{\cT,\rm out}^\pm$ are saddle points.

For later use, we record the relationship between edge-b- and 3b-momenta near $(\iota^+)^\circ$:
\begin{equation}
\label{EqWeb3bRel}
\begin{alignedat}{3}
  \xi_\tbop &= -\half\xi_\ebop,&\qquad
  \eta_\tbop&=x_{\!\scri}^{-1}\eta_\ebop,&\qquad
  \sigma_\tbop&=x_{\!\scri}^{-2}\bigl(\half\xi_\ebop-\zeta_\ebop\bigr), \\
  \xi_\ebop &= -2\xi_\tbop, &\qquad
  \eta_\ebop &= x_{\!\scri}\eta_\tbop, &\qquad
  \zeta_\ebop &= -x_{\!\scri}^2\sigma_\tbop - \xi_\tbop.
\end{alignedat}
\end{equation}
This follows from $\xi_\tbop=r\pa_r(\cdot)$, $\eta_\tbop=\pa_\omega(\cdot)$, $\sigma_\tbop=r\pa_{t_*}(\cdot)$ upon changing variables via~\eqref{EqGAGebpVF}.

We now work in a full neighborhood of $\iota^+\cup\cT^+$. Denoting by $\rho_\infty$ a homogeneous degree $-1$ function on $\Tetb^*M\setminus o$, we set
\begin{equation}
\label{EqWRescHamVF}
  \sfH:=\rho_\infty H_{x_{\!\scri}^{-2}\rho_+^{-2}p}.
\end{equation}
This is a b-vector field on the radial compactification of $\Tetb^*M$ over $\ft_*^{-1}([-1,\infty])$; we denote its restriction to $\Setb^*M$ by the same symbol (but caution the reader that this restriction loses the information about the fiber-radial behavior of $\sfH$ at the radial sets).

\begin{lemma}[Null-bicharacteristic flow over $\cT^+$]
\label{LemmaWFlowT}
  Let $\gamma\colon I\to\pa\Sigma^\pm\cap\Stb^*_{\cT^+}M$ denote a maximally extended integral curve of $\pm\sfH$ over $\cT^+$. Then $I=\R$, and either $\gamma\subset\pa\cR_{\cT,\rm out}^\pm\cup\pa\cR_{\cT,\rm in}^\pm$ (in which case $\gamma$ is constant), or $\gamma(s)$ converges to a point in $\pa\cR_{\cT,\rm in}^\pm$, resp.\ $\pa\cR_{\cT,\rm out}^\pm$, as $s\to-\infty$, resp.\ $s\to+\infty$.
\end{lemma}
\begin{proof}
  Let $\zeta\in\Sigma\cap\Ttb^*_{\cT^+}M\setminus o$, and suppose $H_{G_\tbop}$ is radial at $\zeta$. Then $\zeta$ necessarily lies over $\pa\cT^+$ in view of the nontrapping condition in Definition~\ref{DefGSG}\eqref{ItGSGNontrap}. Over $\pa\cT^+$, we then deduce from~\eqref{EqWHam3b} and the vanishing of the $\pa_{\omega^j}$-components that $\eta_\tbop=0$; since then $0=G_\tbop=\xi_\tbop(\xi_\tbop-2\sigma_\tbop)$ by~\eqref{EqWGtb}, we conclude that $\zeta$ lies in one of the sets~\eqref{EqWRadT}.

  Suppose now $\gamma$ is not contained in, and thus disjoint from, $\pa\cR_{\cT,\rm out}^\pm\cup\pa\cR_{\cT,\rm in}^\pm$. On $\Sigma^\pm$ and over $\pa\cT^+$, we compute $\pm H_{G_\tbop}(\xi_\tbop/\sigma_\tbop)=\pm 2\sigma_\tbop^{-1}|\eta_\tbop|^2<0$; therefore, if $\gamma$ lies over $\pa\cT^+$, it tends to $\pa\cR_{\cT,\rm in}^\pm$, resp.\ $\pa\cR_{\cT,\rm out}^\pm$ as $s\to-\infty$, resp.\ $s\to+\infty$.

  Finally, if $\gamma$ does not lie over $\pa\cT^+$, then the nontrapping assumption in Definition~\ref{DefGSG}\eqref{ItGSGNontrap} implies that $\gamma$ enters any neighborhood of $\pa\cT^+$ both as $s\to-\infty$ and as $s\to+\infty$; thus, in both directions it has accumulation points in $\pa\Sigma^\pm\cap\Stb^*_{\pa\cT^+}M$. Since the sets of accumulations points are $\pm\sfH$-invariant, and since $\pa\cR_{\cT,\rm in}^\pm$, resp.\ $\pa\cR_{\cT,\rm out}^\pm$ is a source, resp.\ sink for the $\pm\sfH$-flow over $\cT^+$, we conclude that $\gamma(s)$ must tend to a point in $\pa\cR_{\cT,\rm out}^\pm$ as $s\to\infty$, and to a point in $\pa\cR_{\cT,\rm in}^\pm$ as $s\to-\infty$. The proof is complete.
\end{proof}

\begin{lemma}[Null-bicharacteristic flow over $\iota^+$]
\label{LemmaWFlowI}
  Let $\gamma\colon I\to\pa\Sigma^\pm\cap\Setb^*_{\iota^+}M$ denote a maximally extended integral curve of $\pm\sfH$ over $\iota^+$. Then $I=\R$ and the limits $p_\pm=\lim_{s\to\pm\infty}\gamma(s)$ exist. There are four possibilities:
  \begin{enumerate}
  \item\label{ItWFlowI1} $\gamma\subset\pa\cR_{\scri,\rm out}^\pm\cup\pa\cR_{\scri,\rm in,+}^\pm\cup\pa\cR_{\cT,\rm out}^\pm\cup\pa\cR_{\cT,\rm in}^\pm$; or
  \item\label{ItWFlowI2} $p_-\in\pa\cR_{\scri,\rm in,+}^\pm\cup\pa\cR_{\cT,\rm out}^\pm$ and $p_+\in\pa\cR_{\scri,\rm out}^\pm$; or
  \item\label{ItWFlowI3} $p_-\in\pa\cR_{\scri,\rm in,+}^\pm$ and $p_+\in\pa\cR_{\cT,\rm in}^\pm$; or
  \item\label{ItWFlowI4} $p_-\in\pa\cR_{\cT,\rm in}^\pm$ and $p_+\in\pa\cR_{\cT,\rm out}^\pm$.
  \end{enumerate}
  Moreover, the level sets of $\frac{r}{t_*}=\rho_\cT\in(0,\infty)$ inside $\iota^+$ are strictly convex for the $\pm\sfH$-flow. In the coordinates used in~\eqref{EqWRadT} and~\eqref{EqWRadScri}, the unstable manifold of $\cR_{\cT,\rm out}^\pm$ is given by
  \begin{equation}
  \label{EqWFlowIUnstable}
  \begin{split}
    \cW_{\rm out}^\pm &:= \{ (\rho_\cT,\rho_+,\omega;\sigma_\tbop,\xi_\tbop,\eta_\tbop) \colon \rho_+=0,\ \xi_\tbop=\eta_\tbop=0,\ \pm\sigma_\tbop<0 \} \\
      &= \{ (x_{\!\scri},\rho_+,\omega;\xi_\ebop,\zeta_\ebop,\eta_\ebop) \colon \rho_+=0,\ \xi_\ebop=\eta_\ebop=0,\ \pm\zeta_\ebop>0 \}.
  \end{split}
  \end{equation}
\end{lemma}
\begin{proof}
  We first prove the convexity statement using the expression~\eqref{EqWHam3b} for $H_{G_\tbop}$. If $H_{G_\tbop}\rho_\cT=2\rho_\cT(\rho_\cT\xi_\tbop+(\xi_\tbop-\sigma_\tbop))=0$ at a point in $\Sigma^\pm$ with $\rho_+=0$ implies $H_{G_\tbop}^2\rho_\cT=2\rho_\cT(2\rho_\cT\sigma_\tbop\xi_\tbop+2\sigma_\tbop^2)>0$ since $2\sigma_\tbop\xi_\tbop=\xi_\tbop^2+|\eta_\tbop|_{\slg^{-1}}^2$ by~\eqref{EqWGtb}.

  The rest of the Lemma can be proved by explicitly computing the flow of $\mp\sigma_\tbop^{-1}H_{G_\tbop}$; we leave this to the reader, and instead give here a more qualitative argument. Note first that the flow over $\iota^+\cap\cT^+$ and $\iota^+\cap\scri^+$ was already discussed in Lemma~\ref{LemmaWFlowT} and \cite[\S4.1]{HintzVasyScrieb}. Possibilities~\eqref{ItWFlowI1} and \eqref{ItWFlowI4} can only occur (and they do occur) when $\gamma$ lies over $\iota^+\cap\cT^+$ or $\iota^+\cap\scri^+$. Also possibility~\eqref{ItWFlowI2} may occur for such $\gamma$, with $p_-\in\pa\cR_{\scri,\rm in,+}^\pm$.

  It thus suffices to consider $\gamma\subset\Stb^*_{(\iota^+)^\circ}M$. Suppose $\rho_{\cT,+}:=\lim_{s\to\infty}\rho_\cT(\gamma(s))=0$ (the limit existing due to the convexity of $\rho_\cT\circ\gamma$, which implies that $\rho_\cT(\gamma(s))$ is monotone for large $|s|$), then $\gamma(s)$ has an accumulation point over $\iota^+\cap\cT^+$ as $s\to\infty$ and therefore must tend to a point in $\pa\cR_{\cT,\rm in}^\pm$, since by~\eqref{EqWRadTHamLinOut} $\rho_\cT$ is increasing over $(\iota^+)^\circ$ near the other radial set $\pa\cR_{\cT,\rm out}^\pm$. Note moreover that in this case, $\rho_{\cT,-}:=\lim_{s\to-\infty}\rho_\cT(\gamma(s))$ must be equal to $+\infty$ by convexity, and an analogous argument then proves that necessarily $p_-\in\pa\cR_{\scri,\rm in,+}^\pm$. This is possibility~\eqref{ItWFlowI3}.

  We cannot have $\rho_{\cT,+}\in(0,\infty)$ by convexity. If $\rho_{\cT,+}=+\infty$, then the source-to-sink structure of the flow over $\iota^+\cap\scri^+$ implies that $p_+\in\pa\cR_{\scri,\rm out}^\pm$; if $\rho_{\cT,-}=0$, we must have $p_-\in\pa\cR_{\cT,\rm out}^\pm$, whereas in the remaining possibility $\rho_{\cT,-}=\infty$ we conclude that $p_-\in\pa\cR_{\scri,\rm in,+}^\pm$. This is the remaining part of possibility~\eqref{ItWFlowI2}.

  Regarding~\eqref{EqWFlowIUnstable}, it suffices to note that the linearization of $\sfH$ has, at each point of $\pa\cR_{\cT,\rm out}^\pm$, a one-dimensional unstable eigenspace by~\eqref{EqWRadTHamLinOut}; therefore, the $\sfH$-invariant submanifold $\cW_{\rm out}^\pm\subset\Sigma^\pm$, which has one-dimensional fibers over $\cR_{\cT,\rm out}^\pm$, must be the unstable manifold indeed.
\end{proof}

\begin{prop}[Dynamics of the null-bicharacteristic flow]
\label{PropWFlow}
  Let $\gamma\colon I\to\pa\Sigma^\pm$ denote a maximally extended integral curve of $\pm\sfH$ over the domain $\ft_*^{-1}([-1,\infty])$; write $s_-=\inf I$ and $s_+=\sup I$, and $p_\pm=\lim_{s\to s_\pm}\gamma(s)$ if these limits exist. Then $s_+=\infty$, and moreover exactly one of the following possibilities must occur:
  \begin{enumerate}
  \item\label{ItWFlowContained} $\gamma\subset\pa\cR_{\scri,\rm out}^\pm\cup\pa\cR_{\scri,\rm in,+}^\pm\cup\pa\cR_{\cT,\rm out}^\pm\cup\pa\cR_{\cT,\rm in}^\pm$; or
  \item\label{ItWFlow2} $s_-=-\infty$ and $p_-\in\pa\cR_{\scri,\rm in,+}^\pm\cup\pa\cR_{\cT,\rm out}^\pm$, $p_+\in\pa\cR_{\scri,\rm out}^\pm$; or
  \item\label{ItWFlow3} $s_-=-\infty$ and $p_-\in\pa\cR_{\scri,\rm in,+}^\pm$, $p_+\in\pa\cR_{\cT,\rm in}^\pm$; or
  \item\label{ItWFlow4} $s_-=-\infty$ and $p_-\in\pa\cR_{\cT,\rm in}^\pm$, $p_+\in\pa\cR_{\cT,\rm out}^\pm$; or
  \item\label{ItWFlow5} $s_->-\infty$ and $\ft_*(\gamma(s_-))=-1$ as well as $p_+\in\pa\cR_{\scri,\rm out}^\pm\cup\pa\cR_{\scri,\rm in,+}^\pm$.
  \end{enumerate}
\end{prop}

\begin{figure}[!ht]
\centering
\includegraphics{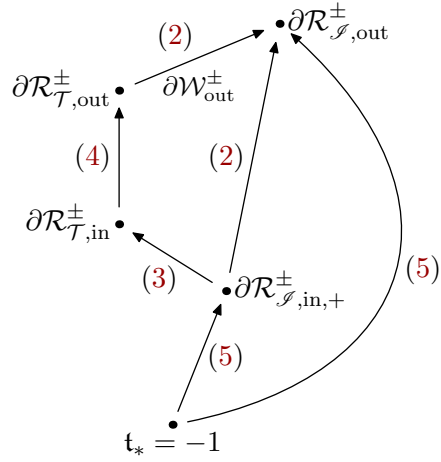}
\caption{Illustration of the dynamics of the (rescaled) null-bicharacteristic flow of $P$ inside $\pa\Sigma^\pm$. The numbers refer to the enumeration in Proposition~\usref{PropWFlow}. Possibility~\eqref{ItWFlowContained} is not drawn here. The boundary at fiber infinity of the unstable manifold $\cW_{\rm out}^\pm$ of $\pa\cR_{\cT,\rm out}^\pm$ (see~\eqref{EqWFlowIUnstable}) is also indicated.}
\label{FigWFlow}
\end{figure}

See Figure~\ref{FigWFlow}. Thus, starting from a point in $\pa\Sigma^+\cap\ft_*^{-1}(-1)$, every point in $\pa\Sigma^\pm$ can be reached by following the null-bicharacteristic flow of $P$ in the future causal direction and passing through the sets of saddle points $\pa\cR_{\scri,\rm in,+}^\pm$, $\pa\cR_{\cT,\rm in}^\pm$, $\pa\cR_{\cT,\rm out}^\pm$ until one reaches the global sink $\pa\cR_{\scri,\rm out}^\pm$.

\begin{proof}[Proof of Proposition~\usref{PropWFlow}]
  A curve $\gamma$ as in part~\eqref{ItWFlowContained} is constant, and therefore $s_\pm=\pm\infty$ and $p_+=p_-$. Thus, from now on, we only consider $\gamma$ which are not of the type~\eqref{ItWFlowContained}. Since either $\gamma$ lies over $\ft_*^{-1}(\infty)=\iota^+\cup\cT^+$ or $\ft_*\circ\gamma$ is monotonically increasing, the fact that $\Sigma^\pm\cap\ft_*^{-1}([-1,\infty])$ is compact implies that $s_+=\infty$. If $\gamma$ lies over $\iota^+\cup\cT^+$, then Lemmas~\ref{LemmaWFlowT} and \ref{LemmaWFlowI} give the desired result; if $\gamma$ lies over $(\scri^+)^\circ$, then \cite[\S4.1]{HintzVasyScrieb} implies that $s_->-\infty$ and $\ft_*(\gamma(s_-))=-1$ and $p_+\in\pa\cR_{\scri,\rm out}^\pm\cup\pa\cR_{\scri,\rm in,+}^\pm$.

  It remains to consider the case that $\gamma\subset S^*M^\circ$ and $\ft_*(\gamma(0))<\infty$ where we assume $0\in I$. Since $\ft_*\circ\gamma$ is monotonically increasing, it has a past limit $\ft_{*,-}:=\lim_{s\to s_-}\ft_*(\gamma(s))\geq -1$. Let $\tilde p_-$ be an accumulation point of $\gamma(s)$ as $s\to s_-$. If $\tilde p_-\in S^*M^\circ$, then $\tilde p_-$ is not a critical point of $\sfH$; thus if $\ft_{*,-}>-1$, then $\gamma$ could be continued past $s=s_-$, contradicting the definition of $\ft_{*,-}$, and we conclude that $\ft_{*,-}=-1$ and $\tilde p_-=p_-$. The possibility that $\tilde p_-$ lies over $\scri^+$ cannot occur, since then $\gamma$ would lie entirely over $\scri^+$.

  Let $\tilde p_+$ be an accumulation point of $\gamma(s)$ as $s\to\infty$. If $\ft_*(\tilde p_+)<\infty$, then using the monotonicity of $\ft_*\circ\gamma$ similarly to before, one concludes that necessarily $\tilde p_+\in\pa\cR_{\scri,\rm out}^\pm$, and the sink nature of this radial set implies that in fact $\gamma$ converges to $p_+$, i.e.\ $\tilde p_+=p_+$. If $\ft_*(\tilde p_+)=\infty$, we claim that necessarily $\tilde p_+=p_+\in\pa\cR_{\scri,\rm out}^\pm$; this follows once we show that the set $\Omega_+$ of accumulation points of $\gamma(s)$ as $s\to\infty$ contains a point in $\pa\cR_{\scri,\rm out}^\pm$. By Lemmas~\ref{LemmaWFlowT} and \ref{LemmaWFlowI}, the set $\Omega_+\cap(\pa\cR_{\scri,\rm in,+}^\pm\cup\pa\cR_{\cT,\rm in}^\pm\cup\pa\cR_{\cT,\rm out}^\pm\cup\pa\cR_{\scri,\rm out}^\pm)$ is nonempty. But since the stable and unstable manifolds of $\pa\cR_{\scri,\rm in,+}^\pm$ lie over $\pa M$ (and are thus disjoint from $M^\circ$), the existence of a point in $\Omega_+\cap\pa\cR_{\scri,\rm in,+}^\pm$ implies that an integral curve inside the unstable manifold of the saddle point $\pa\cR_{\scri,\rm in,+}^\pm$ is contained in $\Omega^+$, and therefore so is its endpoint which either lies in $\pa\cR_{\scri,\rm out}^\pm$ or in $\pa\cR_{\cT,\rm in}^\pm$. In the latter case, similar reasoning implies that $\Omega^+$ contains a curve from $\pa\cR_{\cT,\rm in}^\pm$ to $\pa\cR_{\cT,\rm out}^\pm$ and thus a point in this latter radial set; and by Lemma~\ref{LemmaWFlowI}\eqref{ItWFlowI2} then finally a point in $\pa\cR_{\scri,\rm out}^\pm$. Thus, $\Omega^+$ contains a point in $\pa\cR_{\scri,\rm out}^\pm$ (and by the source nature of $\pa\cR_{\scri,\rm out}^\pm$ we must then in fact have $\Omega^+=\{p_+\}$). This concludes the proof.
\end{proof}

\subsection{Microlocal analysis near the asymptotically translation-invariant regime}
\label{SsWT}

In a collar neighborhood of $\cT^+\subset M$, the microlocal analysis of $P$ utilizes the 3b-pseu\-do\-dif\-fer\-en\-tial algebra. Microlocal elliptic and real principal type propagation estimates on 3b-Sobolev spaces hold by the usual (symbolic) arguments; we thus restrict our attention to neighborhoods of the two radial sets $\cR_{\cT,\rm in}^\pm$ and $\cR_{\cT,\rm out}^\pm$. There, we use the coordinates
\[
  \rho_+=r^{-1}\geq 0,\quad\rho_\cT=\frac{r}{t_*}\geq 0,\quad\omega\in\R^{n-1},\quad
  \sigma_\tbop\in\R,\quad\xi_\tbop\in\R,\quad\eta_\tbop\in\R^{n-1},
\]
as in~\eqref{EqWtbCoords}, with the radial sets being defined in~\eqref{EqWRadT}. We shall also use the coordinates $\rho_\infty,\hat\xi_\tbop,\hat\eta_\tbop$ near fiber infinity defined in~\eqref{EqWRadTCoord}. Lastly, we define the rescaled Hamiltonian vector field
\[
  \sfH=\rho_\infty H_{\rho_+^{-2}p}\in(\CI+\cA^{(\ell_+,\ell_\cT)})\Vb(\ol{\Ttb^*}M).
\]
(This definition is consistent with~\eqref{EqWRescHamVF} in a neighborhood of $\cT^+$.)

\begin{prop}[Propagation through $\cR_{\cT,\rm in}^\pm$]
\label{PropWTIn}
  Let $\alpha_+,\alpha_\cT\in\R$. Let
  \begin{equation}
  \label{EqWTInufuf}
  \begin{alignedat}{2}
    u&\in\Htb^{-\infty,(\alpha_+,\alpha_\cT)}(M;\upbeta^*E),&\qquad f&:=P u\in\Htb^{-\infty,(\alpha_++2,\alpha_\cT)}(M;\upbeta^*E), \\
    \tilde u&\in\Htb^{-\infty,(-\alpha_+-2,-\alpha_\cT)}(M;\upbeta^*E),&\qquad \tilde f&:=P^*\tilde u\in\Htb^{-\infty,(-\alpha_+,-\alpha_\cT)}(M;\upbeta^*E).
  \end{alignedat}
  \end{equation}
  Recall the quantity $\ubar S_{\rm in}$ from Definition~\usref{DefStEstThr}.
  \begin{enumerate}
  \item\label{ItWTInDir}{\rm (Direct problem.)} Let $s>s_0$ where
    \begin{equation}
    \label{EqWTInDirThr}
      s_0+\alpha_+-\alpha_\cT>-\half-\ubar S_{\rm in}.
    \end{equation}
    Suppose that $\WFtb^{s_0,(\alpha_+,\alpha_\cT)}(u)$ and $\WFtb^{s-1,(\alpha_++2,\alpha_\cT)}(f)$ are disjoint from $\pa\cR_{\cT,\rm in}^\pm$, and $\WFtb^{s,(\alpha_+,\alpha_\cT)}(u)\cap(\cU\cap\Stb^*_{(\iota^+)^\circ}M)=\emptyset$ for some neighborhood $\cU\subset\Stb^*_{\iota^+}M$ of $\pa\cR_{\cT,\rm in}^\pm$ over $\iota^+$. Then
    \[
      \WFtb^{s,(\alpha_+,\alpha_\cT)}(u)\cap\pa\cR_{\cT,\rm in}^\pm=\emptyset.
    \]
  \item\label{ItWTInAdj}{\rm (Adjoint problem.)} Let $s\in\R$ with $s+\alpha_+-\alpha_\cT>-\half-\ubar S_{\rm in}$. Suppose that $\WFtb^{-s,(-\alpha_+,-\alpha_\cT)}(\tilde f)\cap\pa\cR_{\cT,\rm in}^\pm=\emptyset$, $\WFtb^{-s+1,(-\alpha_+-2,-\alpha_\cT)}(\tilde u)\cap(\cV\setminus\pa\cR_{\cT,\rm in}^\pm)=\emptyset$ for some neighborhood $\cV\subset\Stb^*_{\cT^+}M$ of $\pa\cR_{\cT,\rm in}^\pm$ over $\cT^+$. Then
  \[
    \WFtb^{-s+1,(-\alpha_+-2,-\alpha_\cT)}(\tilde u)\cap\pa\cR_{\cT,\rm in}^\pm=\emptyset.
  \]
  \end{enumerate}
\end{prop}

In the threshold condition~\eqref{EqWTInDirThr}, one should think of the relative order $\ell:=\alpha_+-\alpha_\cT$ as the effective spatial weight; condition~\eqref{EqWTInDirThr} then reads $s_0+\ell>-\half-\ubar S_{\rm in}$ and in this form (with $s$ in place of $s_0$) already appeared in Proposition~\ref{PropStEstNz}\eqref{ItStEstNzb}. See Remark~\ref{RmkWTOutThr} for an analogous discussion of the second threshold assumption in Proposition~\ref{PropStEstNz}\eqref{ItStEstNzb}.

Our proof of Proposition~\ref{PropWTIn} will proceed via a standard (regularized) positive commutator argument, as explained in detail e.g.\ in \cite[\S4.4]{VasyMinicourse}. Such an argument gives a quantitative estimate which implies the stated qualitative result on the wave front set of $u$. For part~\eqref{ItWTInDir}, the quantitative estimate takes the form
\begin{equation}
\label{EqWTDirEst}
\begin{split}
  \|B u\|_{\Htb^{s,(\alpha_+,\alpha_\cT)}(M;\upbeta^*E)} &\leq C\Bigl( \| G P u \|_{\Htb^{s-1,(\alpha_++2,\alpha_\cT)}(M;\upbeta^*E)} + \| E u \|_{\Htb^{s,(\alpha_+,\alpha_\cT)}(M;\upbeta^*E)} \\
    &\quad\hspace{14em} + \|\chi u\|_{\Htb^{s_0,(\alpha_+,\alpha_\cT)}(M;\upbeta^*E)}\Bigr),
\end{split}
\end{equation}
where $B,E,G\in\Psitb^0(M)$ and $\chi\in\CI(M)$ are as follows:
\begin{itemize}
\item $\WFtb'(B)\subset\Elltb(G)$ and $\pa\cR_{\cT,\rm in}^\pm\subset\Elltb(G)$;
\item all backwards integral curves of $\pm\sfH$ from $\WFtb'(B)\cap\Sigma^\pm$ either tend to $\pa\cR_{\cT,\rm in}^\pm$ or enter $\Elltb(E)$ in finite time, all while remaining in $\Elltb(G)$;
\item the Schwartz kernels of $B,E,G$ are supported in the interior of $\supp\chi\times\supp\chi$.
\end{itemize}
Moreover, the constant $C$ depends only on $s,s_0,\alpha_+,\alpha_\cT$ as well as on a suitable seminorm (which depends on $s$ as well) of $P$; thus, the estimate~\eqref{EqWTDirEst} is locally uniform for continuous families of admissible wave type operators $P$. The qualitative propagation result amounts to the statement that~\eqref{EqWTDirEst} holds in the strong sense that the left hand side is finite provided the right hand side is (and then the estimate holds). We leave statements of quantitative versions of the other wave front set results in this section to the reader.

\begin{proof}[Proof of Proposition~\usref{PropWTIn}]
  For notational simplicity, we only prove estimates near $\cR_{\cT,\rm in}^+$, where $\sigma_\tbop<0$ and $\rho_\infty=-\frac{1}{\sigma_\tbop}$, $\hat\xi_\tbop=\frac{\xi_\tbop}{\sigma_\tbop}$, $\hat\eta_\tbop=\frac{\eta_\tbop}{\sigma_\tbop}$; the arguments near $\cR_{\cT,\rm in}^-$ are completely analogous. Postponing issues of regularization for now, we consider a commutant
  \begin{equation}
  \label{EqWTInCommutant}
  \begin{split}
    &A = \check A^*\check A,\qquad \check A:=\Op_\tbop(\check a), \\
    &\qquad \check a=\rho_\infty^{-s+\frac12}\rho_+^{-\alpha_+-1}\rho_\cT^{-\alpha_\cT}\chi_+(\rho_+)\chi_\cT(\rho_\cT)\chi_1((\hat\xi_\tbop-2)^2)\chi_2(|\hat\eta_\tbop|^2),
  \end{split}
  \end{equation}
  where the cutoffs $\chi_+,\chi_\cT,\chi_1,\chi_2\in\CIc([0,\delta))$ (with $\delta>0$ chosen later) are identically $1$ near $0$; we furthermore arrange $(-\chi_+\chi_+')^{1/2}\in\CIc([0,\delta))$, likewise for $\chi_\cT$, $\chi_1$, and $\chi_2$. Consider then the $L^2$-inner product
  \begin{equation}
  \label{EqWTInComm}
    \Im\la P u,A u\ra = \la\sC u,u\ra,\qquad
    \sC = \frac{i}{2}(P^*A-A P)=\frac{i}{2}[P,A] + \frac{P-P^*}{2 i}A.
  \end{equation}
  Thus, $\sC\in\Psitb^{2 s,2\alpha_+,2\alpha_\cT}(M;\upbeta^*E)$. The only contribution to the principal symbol of the first summand which is elliptic at $\cR_{\cT,\rm in}^+$ arises from differentiation of the weights of $\check a^2$ along $\half$ times the Hamiltonian vector field of $p$. In view of~\eqref{EqWRadTHamLinIn} (and noting that $G_\tbop=\rho_+^{-2}p$ gives $\rho_+^2 H_{G_\tbop}=H_p+p\rho_+^2 H_{\rho_+^{-2}}$, with the second term vanishing on the characteristic set) this contribution is
  \begin{equation}
  \label{EqWTInCommMain}
    2(-s-\alpha_++\alpha_\cT-\half)\rho_\infty^{-2 s}\rho_+^{-2\alpha_+}\rho_\cT^{-2\alpha_\cT}\chi_+^2\chi_\cT^2\chi_1^2\chi_2^2.
  \end{equation}
  Differentiation of $\chi_\cT$ on the other hand gives, for sufficiently small $\delta>0$, a nonnegative contribution (and positive where $\chi'_\cT\neq 0$); differentiation of the other cutoffs gives nonpositive contributions.

  We proceed to compute the principal symbol of the second summand in~\eqref{EqWTInComm} at $\cR_{\cT,\rm in}^+$ and relate this to the analogous computations in the proof of Proposition~\ref{PropStEstNz}\eqref{ItStEstNzsc}. Firstly, we can replace $P$ by $P_0$ (up to an error term in $\Psitb^{2 s,2\alpha_++\ell_+,2\alpha_\cT+\ell_\cT}$), which is thus of the form~\eqref{EqGSOStruct}. As in the proof of Proposition~\ref{PropStEstNz}, we can furthermore work with the volume density $r^2|\dd t_*\,\dd r\,\dd\slg|$ modulo error terms which do not contribute to the principal symbol of $\sC$ at $\iota^+$. Note moreover that the coordinates $(t_*,\rho,\omega)$ used in~\eqref{EqGSOStruct} are related to $(\rho_\cT,\rho_+,\omega)$ via $\rho_\cT=\rho^{-1}t_*^{-1}$, $\rho_+=\rho$, and thus
  \[
    \pa_{t_*}=-\rho_+\rho_\cT^2\pa_{\rho_\cT},\qquad
    \rho\pa_\rho=\rho_+\pa_{\rho_+}-\rho_\cT\pa_{\rho_\cT}.
  \]
  Since $r\pa_{t_*}$, $r\pa_r$, $\pa_\omega$ are smooth 3b-vector fields near $\iota^+\cap\cT^+$, we have $Q\pa_{t_*},g^{0 0}\pa_{t_*}^2\in\cA^{(3+\delta,0)}\Difftb^2(M;\upbeta^*E)$. We thus obtain
  \begin{equation}
  \label{EqWTInImPart}
  \begin{split}
    \frac{P_0-P_0^*}{2 i} &\equiv 2\rho\,\Re S\ D_{t_*} + \Im\hat P(0) \bmod \cA^{(2+\delta,0)}\Difftb^1(M;\upbeta^*E) \\
      &\equiv -2 \rho_+^2\Bigl( \Re S\ \rho_\cT^2 D_{\rho_\cT} - \half\rho^{-2}\Im\hat P(0)\Bigr) \bmod \cA^{(2+\delta,0)}\Difftb^1(M;\upbeta^*E);
  \end{split}
  \end{equation}
  see also~\eqref{EqStEstscqSigma}. Let $p_1=\sigmatb^1(\Im P_0)\in(\rho_+^2\CI+\cA^{(2+\delta,0)})(M;\upbeta^*\End(E))$. In view of $\sigmatb^1(\rho_\cT^2 D_{\rho_\cT})=\sigmatb^1(-r D_{t_*})=-\sigma_\tbop>0$ and
  \[
    \half\,\sigmatb^1\bigl(\rho^{-2}\Im\hat P(0)\bigr)\bigr|_{\sigma_\tbop\frac{\dd t_*}{r}+2\sigma_\tbop\frac{\dd r}{r}} = -\sigma_\tbop\cdot\half\,\sigmasc^{1,0}\Bigl(\frac{\hat P(0)-\hat P(0)^*}{2 i\rho}\Bigr)\Big|_{2\frac{\dd\rho}{\rho^2}},
  \]
  we conclude from Definition~\ref{DefStEstThr} (and recalling $\rho_\infty=(-\sigma_\tbop)^{-1}$) that one can choose a fiber inner product on $E$ near $\iota^+\cap\cT^+$ so that, for any fixed $\eps>0$, we have
  \begin{equation}
  \label{EqWTInImPart2}
    \rho_\infty\rho^{-2}p_1 \leq -2(\ubar S_{\rm in} - \eps)
  \end{equation}
  in a $\delta$-neighborhood of $\pa\cR_{\rm in}^+$ for sufficiently small $\delta>0$.

  We conclude that we can write
  \begin{equation}
  \label{EqWTInCommSq}
    \sigmatb^{2 s,2\alpha_+,2\alpha_\cT}(\sC) = -\eta\rho_\infty^{-1}\rho_+^2\check a^2 - b_0^2 - b_+^2 - b_1^2 - b_2^2 + b_\cT^2
  \end{equation}
  for sufficiently small $\eta>0$, where
  \[
    b_0 = \rho_\infty^{-s}\rho_+^{-\alpha_+}\rho_\cT^{-\alpha_\cT}\chi_+\chi_\cT\chi_1\chi_2 \Bigl( -\half\rho_\infty^{2 s}\rho_+^{2\alpha_+}\rho_\cT^{2\alpha_\cT}H_p\bigl(\rho_\infty^{-2 s+1}\rho_+^{-2\alpha_+-2}\rho_\cT^{-2\alpha_\cT}\bigr) - p_1 - \eta \Bigr)^{1/2}
  \]
  is a weighted symbol which is elliptic at $\cR_{\cT,\rm in}^+$; note that at $\cR_{\cT,\rm in}^+$, the term in parentheses is is bounded from below by $2(s+\alpha_+-\alpha_\cT+\half)+2(\ubar S_{\rm in}-\eps)-\eta$ by~\eqref{EqWTInCommMain} and~\eqref{EqWTInImPart2}, which upon choosing $\eps,\eta>0$ small is indeed positive at $\cR_{\cT,\rm in}^+$ and hence near $\supp\check a$ when $\delta>0$ is sufficiently small. The terms $b_\bullet$ in~\eqref{EqWTInCommSq} for $\bullet=+,1,2,\cT$ arise from differentiation of the cutoff $\chi_\bullet$ in the definition of $\check a$; the signs were already discussed above.

  Assuming that $u\in\Htb^{s,(\alpha_+,\alpha_\cT)}(M;\upbeta^*E)$, we can now prove a quantitative estimate on $u$ in the usual fashion by quantizing the relationship~\eqref{EqWTInCommSq} to an equality of ps.d.o.s,
  \[
    \sC = -\eta(\Lambda\check A)^*(\Lambda\check A) - B_0^*B_0 - B_+^*B_+ - B_1^*B_1 - B_2^*B_2 + B_\cT^*B_\cT + R,
  \]
  where $\Lambda=\Op_\tbop(\rho_\infty^{-1/2}\rho_+)$ and $R\in\Psitb^{2 s-1,2\alpha_+,2\alpha_\cT}$. Plugging this into~\eqref{EqWTInComm} and using a parametrix $\Lambda^-$ of $\Lambda$ with $\Lambda^-\Lambda=I-R'$ where $R'\in\Psitb^{-\infty}$, one obtains the $L^2$-estimate
  \[
    -\eta\|\Lambda\check A u\|^2 - C_\eta\|(\Lambda^-)^*\check A P u\|^2 \geq -\eta\|\Lambda\check A u\|^2 - \|B_0 u\|^2 + \|B_\cT u\|^2 - \| R'' u \|_{\Htb^{s-\frac12,\alpha_+,\alpha_\cT}},
  \]
  where $R''\in\Psitb^0$ is used to bound the terms $|\la R u,u\ra|$ and $|\la \check A P u,R'\check A u\ra|$, and where we dropped the contributions from $B_+$, $B_1$, and $B_2$ on the right hand side. (The integration by parts here is justified using a regularization argument as in \cite[Proof of Proposition~4.11]{VasyMinicourse}.) Upon canceling the first terms on both sides and rearranging, we thus obtain an estimate for $B_0 u$ in terms of $P u$ and an a priori control term $B_\cT u$; the remainder term can be improved to $\|R'''u\|_{\Htb^{s_0,\alpha_+,\alpha_\cT}}$ by an iterative application of this estimate (improving control by half a 3b-differential order at a time).

  In order to prove the propagation of microlocal 3b-regularity, one uses a regularized version of the above argument. Thus, one quantizes $\check a_\gamma=\phi_\gamma(\rho_\infty)\check a$ where $\phi_\gamma(\rho_\infty)=(1+\gamma\rho_\infty^{-1})^{-(s-s_0)}$; for $\gamma>0$, the symbol $\check a_\gamma$ thus has differential order $-s_0+\half$, and therefore the above argument goes through (with $B_0$ replaced by correspondingly regularized versions $B_{0,\gamma}$ etc.) due to the a priori regularity assumption on $u$; indeed one obtains a uniform $L^2$-estimate for $B_{0,\gamma}u$. As $\gamma\searrow 0$, a standard functional analytic argument (using the weak compactness of the unit ball in $L^2$) then implies $B_0 u\in L^2$, and finishes the proof of part~\eqref{ItWTInDir}. (See \cite[\S4.4]{VasyMinicourse} for details.)

  The proof of part~\eqref{ItWTInAdj} proceeds along the same lines, with $s$, $\alpha_+$, and $\alpha_\cT$ in~\eqref{EqWTInCommutant} replaced by $-s+1$, $-\alpha_+-2$, and $-\alpha_\cT$, respectively. Since passing from the operator $P$ to its adjoint $P^*$ also switches the sign of the imaginary part, the operator $\sC$ now has a \emph{positive} principal symbol at $\cR_{\cT,\rm in}^+$ under the stated condition on $s$. Therefore, the propagation estimate has a priori control terms arising from differentiation of the cutoffs $\chi_+$, $\chi_1$, $\chi_2$, which are controlled by the a priori regularity assumption on $u$ when the support of $\check a$ is chosen small enough. We also note that the amount of regularization is now arbitrary (since upon decreasing the regularity $-s+1$ of $\tilde u$, i.e.\ increasing $s$, the threshold condition on $s$ remains valid), which explains why there is no a priori regularity requirement on $\tilde u$ \emph{at} the radial set.
\end{proof}

\begin{prop}[Propagation through $\cR_{\cT,\rm out}^\pm$]
\label{PropWTOut}
  Let $\alpha_+,\alpha_\cT\in\R$. Let $u,f,\tilde u,\tilde f$ be as in~\eqref{EqWTInufuf}. Recall the quantity $\ubar S$ from Definition~\usref{DefStEstThr}.
  \begin{enumerate}
  \item\label{ItWTOutDir}{\rm (Direct problem.)} Let $s\in\R$ with
    \begin{equation}
    \label{EqWTOutDirThr}
      s + \alpha_+ - \alpha_\cT < -\half+\ubar S.
    \end{equation}
    Suppose that $\WFtb^{s-1,(\alpha_++2,\alpha_\cT)}(f)\cap\pa\cR_{\cT,\rm out}^\pm=\emptyset$, $\WFtb^{s,(\alpha_+,\alpha_\cT)}(u)\cap(\cU\setminus\pa\cR_{\cT,\rm out}^\pm)=\emptyset$ for some neighborhood $\cU\subset\Stb^*_{\cT^+}M$ of $\pa\cR_{\cT,\rm out}^\pm$ over $\cT^+$. Then
    \[
      \WFtb^{s,(\alpha_+,\alpha_\cT)}(u)\cap\pa\cR_{\cT,\rm out}^\pm=\emptyset.
    \]
  \item\label{ItWTOutAdj}{\rm (Adjoint problem.)} Let $s<s_0$ where $s_0+\alpha_+-\alpha_\cT<-\half+\ubar S$. Suppose that $\WFtb^{-s_0+1,(-\alpha_+-2,-\alpha_\cT)}(\tilde u)$ and $\WFtb^{-s,(-\alpha_+,-\alpha_\cT)}(\tilde f)$ are disjoint from $\pa\cR_{\cT,\rm out}^\pm$, and $\WFtb^{-s+1,(-\alpha_+-2,-\alpha_\cT)}(\tilde u)\cap(\cV\setminus\pa\cR_{\cT,\rm out}^\pm)=\emptyset$ for some neighborhood $\cV\subset\Stb^*_{\iota^+}M$ of $\pa\cR_{\cT,\rm out}^\pm$ over $\iota^+$. Then
  \[
    \WFtb^{-s+1,(-\alpha_+-2,-\alpha_\cT)}(\tilde u)\cap\pa\cR_{\cT,\rm out}^\pm=\emptyset.
  \]
  \end{enumerate}
\end{prop}

\begin{rmk}[Threshold condition]
\label{RmkWTOutThr}
  Upon introducing the effective weight $\ell=\alpha_+-\alpha_\cT$, condition~\eqref{EqWTOutDirThr} reads $s+\ell<-\half+\ubar S$. Note then that upon embedding $\Tsc^*\cT^+$ as the level set $\sigma_\tbop=\sigma r$ for fixed frequency $\sigma\neq 0$, the quantity $s+\ell$ (with $s$ the 3b-differential order) is the scattering decay order at $\Tsc^*_{\pa\cT^+}\cT^+$, and $\cR_{\cT,\rm out}^\pm$ corresponds to the scattering zero section. In this fashion, condition~\eqref{EqWTOutDirThr} corresponds to the condition in Proposition~\ref{PropStEstNz}\eqref{ItStEstNzb} regarding the upper bound on the scattering decay order at zero scattering frequency by $-\half+\ubar S$.
\end{rmk}

\begin{proof}[Proof of Proposition~\usref{PropWTOut}]
  The proof is a minor modification of that of Proposition~\ref{PropWTIn}. We only comment on part~\eqref{ItWTOutDir}. We use a commutant as in~\eqref{EqWTInCommutant} but with the term $\chi_1$ replaced by $\chi_1(\hat\xi_\tbop^2)$ so as to localize near $\cR_{\cT,\rm out}^+$. Using~\eqref{EqWRadTHamLinOut}, the contribution to the principal symbol of $\sC$ arising from differentiating the weights is now
  \begin{equation}
  \label{EqWTOutCalc1}
    2(s+\alpha_+-\alpha_\cT+\half)\rho_\infty^{-2 s}\rho_+^{-2\alpha_+}\rho_\cT^{-2\alpha_\cT}\chi_+^2\chi_\cT^2\chi_1^2\chi_2^2.
  \end{equation}
  this replaces the expression~\eqref{EqWTInCommMain}. The contribution at $\cR_{\cT,\rm out}^+$ from the imaginary part of $P$ can be estimated using~\eqref{EqWTInImPart} and in the notation of~\eqref{EqWTInImPart2} (and in view of Definition~\ref{DefStEstThr}) satisfies
  \begin{equation}
  \label{EqWTOutCalc2}
    \rho_\infty\rho^{-2}p_1 \leq -2(\ubar S-\eps)
  \end{equation}
  for any fixed $\eps>0$ upon choosing an appropriate fiber inner product on $E$. Differentiation of $\chi_+$, $\chi_1$, $\chi_2$ gives nonnegative terms (which can be written as squares of smooth symbols), which are the origin of the a priori control assumption on $u$ in a punctured neighborhood of $\pa\cR_{\cT,\rm out}^+$ over $\cT^+$; differentiation of $\chi_\cT$ on the other hand gives a nonpositive term. Thus, we can propagate 3b-regularity estimates through the radial set $\pa\cR_{\cT,\rm out}^+$ provided the sum of $2(s+\alpha_+-\alpha_\cT+\half)$ from~\eqref{EqWTOutCalc1} and $-2\ubar S$ from~\eqref{EqWTOutCalc2} is negative at $\pa\cR_{\cT,\rm out}^+$ for some sufficiently small $\eps>0$; this is precisely condition~\eqref{EqWTOutDirThr}. To prove the propagation of regularity, one needs to regularize the computation as before; there is no limit to the amount of regularization since~\eqref{EqWTOutDirThr} remains valid if one decreases $s$.
\end{proof}

\subsection{Global 3b-regularity estimate}
\label{SsWR}

Unless $-\ubar S_{\rm in}>\ubar S$, there is no constant 3b-dif\-fer\-en\-tial order (i.e.\ real number) $s$ for which the threshold conditions~\eqref{EqWTInDirThr} (with $s$ in place of $s_0$) and~\eqref{EqWTOutDirThr} hold simultaneously.\footnote{This is a common situation; for example, $\ubar S=\ubar S_{\rm in}=0$ for scalar wave operators associated with stationary and asymptotically metrics (see the comment after Proposition~\ref{PropES}) and also for admissible asymptotically flat metrics (see the comment after Proposition~\ref{PropES2}). The discussion of the threshold conditions shows that the need for variable orders in such settings has exactly the same origin as the analogous need in stationary scattering theory as described in \cite[Proposition~4.13]{VasyMinicourse}.} Both threshold conditions can however always be satisfied if we use a variable 3b-differential order instead. For real principal type propagation, it is in addition necessary that this order be monotonically decreasing along the null-bicharacteristic flow. We shall work in the region $\ft_*\geq -1$ (using the notation of Definition~\ref{DefGAG}), and thus in the edge-3b-cotangent bundle; by $\sfH$ we denote the rescaled Hamiltonian vector field~\eqref{EqWRescHamVF}.

\begin{lemma}[Existence of variable 3b-differential order functions]
\label{LemmaWROrder}
  Let $\alpha_+,\alpha_\cT\in\R$, and recall the quantities $\ubar S,\ubar S_{\rm in}$ associated with the wave type operator $P$ from Definition~\usref{DefStEstThr}. Fix $s_{\rm in}\geq s_{\rm out}$ with $s_{\rm in}+\alpha_+-\alpha_\cT>-\half-\ubar S_{\rm in}$ and $s_{\rm out}+\alpha_+-\alpha_\cT<-\half+\ubar S$. Then there exists a function $\sfs\in\CI(\Setb^*M\cap\{\ft_*\geq-1\})$, with $s_{\rm out}\leq\sfs\leq s_{\rm in}$, with the following properties.
  \begin{enumerate}
  \item\label{ItWRAbove}{\rm (Above the incoming threshold.)} Near $\pa\cR_{\cT,\rm in}^\pm$, we have $\sfs=s_{\rm in}$.
  \item\label{ItWRBelow}{\rm (Below the outgoing threshold.)} Near $\pa\cR_{\cT,\rm out}^\pm$, we have $\sfs=s_{\rm out}$.
  \item\label{ItWRConst}{\rm (Constancy near radial sets.)} The function $\sfs$ is constant near $\pa\cR_{\cT,\rm in}^\pm$, $\pa\cR_{\cT,\rm out}^\pm$, $\pa\cR_{\scri,\rm in,+}^\pm$, and $\pa\cR_{\scri,\rm out}^\pm$.
  \item\label{ItWRMono}{\rm (Monotonicity.)} We have $\pm\sfH\sfs\leq 0$ on $\pa\Sigma^\pm$.
  \item\label{ItWRStat}{\rm (Properties of the stationary model.)} Properties~\eqref{ItWRAbove}--\eqref{ItWRMono} hold at the same time also for the stationary model $P_0$ (see Definition~\usref{DefGAW}) of $P$, with the monotonicity holding on the characteristic set of $P_0$.
  \end{enumerate}
  Moreover, we can arrange for $\sfs$ to be equal to $s_{\rm in}$ in the complement of any fixed neighborhood of $\bigcup_\pm(\pa\cR_{\scri,\rm out}^\pm\cup\pa\cW_{\rm out}^\pm)$ where $\cW^\pm_{\rm out}\supset\cR_{\cT,\rm out}^\pm$ is the unstable manifold of $\cR_{\cT,\rm out}^\pm$ defined in~\eqref{EqWFlowIUnstable}.
\end{lemma}

\begin{rmk}[Modifications of a given variable order function]
\label{RmkWRMod}
  For later use, we note that given $\sfs$ and $k\geq 0$, the function $\sfs+k$ satisfies all requirements (with $s_{\rm in}+k$ in place of $s_{\rm in}$) except~\eqref{ItWRBelow} (which of course does hold for $s_{\rm out}+k$ in place of $s_{\rm out}$, but $s_{\rm out}+k$ may exceed the required upper bound by $-\half+\ubar S-\alpha_++\alpha_\cT$). Moreover, given any $\lambda>0$ and $q\in\R$ so that $s'_{\rm in}=\lambda s_{\rm in}+q$ and $s'_{\rm out}=\lambda s_{\rm out}+q$ satisfy the required threshold conditions, the function $\sfs'=\lambda\sfs+q$ satisfies the same properties; in particular, if $s'_{\rm in}\geq s_{\rm in}+\delta$ and $s'_{\rm out}\geq s_{\rm out}+\delta$ for some $\delta\in\R$, then also $\sfs'\geq\sfs+\delta$.
\end{rmk}

\begin{rmk}[Same order function for other operators]
\label{RmkWROrderOther}
  The proof gives a stronger result: if $\sP_0$ is a bounded set of stationary wave type operators (with respect to the smooth, or even just the $\cC^1$ topology on the coefficients of $P$ in Definition~\ref{DefGSO}), and $\tilde\sP$ is a set of operators which are bounded in the $|\cdot|_k$ seminorm for $k=1$ defined in~\eqref{EqGAWNorm}, and so that $P'=P_0'+\tilde P'$ is an admissible wave type operator for $P_0'\in\sP_0$ and $\tilde P'\in\tilde\sP$, then one can choose $\sfs$ so that all conditions in Lemma~\ref{LemmaWROrder} hold for all such operators $P'$. The significance of requiring control of up to one derivative on the coefficients of the operators is that this gives pointwise bounds on the Hamiltonian vector field.
\end{rmk}

The idea of the proof is to take $\sfs$ to be constant and above the incoming threshold throughout most of the phase space, except near $\pa\cW_{\rm out}^\pm\cup\pa\cR_{\scri,\rm out}^\pm$ where it decreases to be below the outgoing threshold. A rigorous construction (which the reader may safely skip at first reading) requires some care; it is given in Appendix~\ref{SWROrder}.

\begin{rmk}[Localizing the set where $\sfs$ is not constant]
  For the subsequent analysis, it would be sufficient to use a less regular variable order function $\sfs$ which only satisfies $\sfs\in(\CI+\cA^{(2\ell_{\!\scri},\ell_+,\infty)})(\Setb^*M)$. The gain would be that one could make $\sfs$ constant in $\ft_*\leq C$ for any desired $C$; indeed, one could use the time function $\ft_*$ in place of the function $f_\ft$ in~\eqref{EqWRft}, and the rest of the proof would go through with minor modifications. Smooth $\sfs$ which are constant in $\ft_*\leq C$ typically fail the monotonicity requirement~\eqref{ItWRMono} (the issue being that outgoing null geodesics near $\scri^+$ may lie inside of a hypersurface of the form $\rho_+=C_1-C_2 x_{\!\scri}^{2\ell_{\!\scri}}$, which is not $\cC^{0,1}$ down to $x_{\!\scri}=0$).
\end{rmk}

\begin{prop}[Global control of 3b-regularity]
\label{PropWR}
  In the notation of Definition~\usref{DefGAW}\eqref{ItGAWEdgeN}, and in analogy to Definition~\usref{DefStEstThr}, define
  \begin{equation}
  \label{EqWRubarp1}
    \ubar p_1 := \inf_{p\in\scri^+} (\min\Re\spec p_1(p)).
  \end{equation}
  Let $\alpha_{\!\scri},\alpha_+,\alpha_\cT\in\R$ be such that
  \begin{equation}
  \label{EqWRThr}
    \alpha_+<\alpha_{\!\scri}-\frac12,\qquad
    \alpha_{\!\scri}<-\frac12+\ubar p_1,\qquad
  \end{equation}
  Working in $\ft_*\geq -1$, let $\sfs,\sfs_0,\sfs_1\in\CI(\Setb^*M\cap\{\ft_*\geq-1\})$ be three variable order functions as in Lemma~\usref{LemmaWROrder}, with $\sfs_1>\sfs>\sfs_0$; in particular, $\sfs_0+\alpha_+-\alpha_\cT>-\half-\ubar S_{\rm in}$ at $\pa\cR_{\cT,\rm in}^\pm$ and $\sfs_1+\alpha_+-\alpha_\cT<-\half+\ubar S$ at $\pa\cR_{\cT,\rm out}^\pm$. Let $-1<T_\flat<T_{-1}<T_0<T_1<T_2$, and let $\chi_0,\chi_1\in\CI(\R_{\ft_*})$ be cutoffs with support in $\ft_*\geq T_{-1}$, with $\chi_0=1$ on $[T_0,\infty)$, and $\chi_1=1$ on $[T_1,\infty)$, and let $\chi_\flat\in\CIc((-\infty,T_2))$ be a cutoff with $\chi_\flat=1$ on $[T_\flat,T_1]$. Then for any $N$, there exists a constant $C>0$ so that the estimate
  \begin{equation}
  \label{EqWRDir}
  \begin{split}
    &\|\chi_1 u\|_{H_\etbop^{\sfs,(2\alpha_{\!\scri},\alpha_+,\alpha_\cT)}(M;\upbeta^*E)} \\
    &\qquad \leq C\Bigl( \|\chi_0 P u\|_{H_\etbop^{\sfs-1,(2\alpha_{\!\scri}+2,\alpha_++2,\alpha_\cT)}(M;\upbeta^*E)} + \|\chi_0 u\|_{H_\etbop^{\sfs_0,(2\alpha_{\!\scri},\alpha_+,\alpha_\cT)}(M;\upbeta^*E)} \\
    &\qquad \hspace{18em} + \|\chi_\flat u\|_{H_\etbop^{\sfs,(2\alpha_{\!\scri},-N,-N)}(M;\upbeta^*E)}\Bigr).
  \end{split}
  \end{equation}
  holds in the strong sense that if the right hand side is finite, then so is the left hand side and the estimate holds. We moreover have the following adjoint estimate, which holds in the same strong sense:
  \begin{equation}
  \label{EqWRAdj}
  \begin{split}
    &\|\chi_1\tilde u\|_{H_\etbop^{-\sfs+1,(-2\alpha_{\!\scri}-2,-\alpha_+-2,-\alpha_\cT)}(M;\upbeta^*E)} \\
    &\quad \leq C\Bigl( \|\chi_0 P^*\tilde u\|_{H_\etbop^{-\sfs,(-2\alpha_{\!\scri},-\alpha_+,-\alpha_\cT)}(M;\upbeta^*E)} + \|\chi_0\tilde u\|_{H_\etbop^{-\sfs_1+1,(-2\alpha_{\!\scri}-2,-\alpha_+-2,-\alpha_\cT)}(M;\upbeta^*E)}\Bigr).
  \end{split}
  \end{equation}
\end{prop}

The existence of order functions $\sfs_0,\sfs_1,\sfs$ with the required properties follows from the construction in Remark~\ref{RmkWRMod}. Moreover, the orders on $\chi_\flat u$ in \eqref{EqWRDir} at $\iota^+$ and $\cT^+$ are arbitrary since $\supp\chi_\flat$ is disjoint from these boundary hypersurfaces.

\begin{rmk}[Discussion]
\fakephantomsection
\label{RmkWRDiscuss}
  \begin{enumerate}
  \item In the estimate~\eqref{EqWRDir}, the term involving $\chi_\flat u$ is supported in a bounded `initial' time interval. (This term can easily be controlled from Cauchy data or using a support condition; see~\S\ref{SsWM}.) Thus, the estimate~\eqref{EqWRDir} allows us to propagate edge-3b-regularity of degree $\sfs$ from this initial time interval towards the entire causal future, assuming that $u$ is known to have a minimal amount of regularity as required by the threshold condition at $\pa\cR_{\cT,\rm in}^\pm$.
  \item In the adjoint estimate~\eqref{EqWRAdj} on the other hand, where the propagation direction is reversed, the analogue of the initial control term $\chi_\flat u$ in~\eqref{EqWRDir} are the a priori regularity and boundedness assumptions encoded in the finiteness of the norm of $\chi_0\tilde u$ on the right hand side.
  \item Ignoring the mismatch of cutoff functions on both sides of the estimates~\eqref{EqWRDir}--\eqref{EqWRAdj}, these estimate provide global control of edge-3b-regularity. However, the weights in the norms on $\chi_j u$ or $\chi_j\tilde u$, $j=0,1$, on both sides are still the same. Thus, Proposition~\ref{PropWR} is far from providing a Fredholm setting for $P$ or $P^*$ (see also the comments after~\eqref{EqMFHbRellich}); obtaining this in addition requires control of $u$ to leading order also at the boundary hypersurfaces of $M$, which will be accomplished via the inversion of normal operators in~\S\S\ref{SsWip}--\ref{SsWM}.
  \end{enumerate}
\end{rmk}

\begin{proof}[Proof of Proposition~\usref{PropWR}]
  If we replace $\chi_1 u$ on the left hand side of~\eqref{EqWRDir} by $\chi_1 A u$ where $A\in\Psietb^0(M;\upbeta^*E)$ has operator wave front set disjoint from the characteristic set $\Sigma$ of $P$, the resulting estimate holds true (in the strong sense) by elliptic regularity in the edge-3b-calculus. It thus suffices to consider $\chi_1 u$ microlocally on $\Sigma$.

  By Proposition~\ref{PropWFlow} (and using the timelike nature of $\ft_*$), we can propagate $H_\etbop^{\sfs,(2\alpha_{\!\scri},\alpha_+,\alpha_\cT)}$ regularity of $u$ from $\ft_*^{-1}([T_\flat,T_1])$ (cf.\ the last term in~\eqref{EqWRDir}) in the future causal direction and thus obtain $H_\etbop^{\sfs,(2\alpha_{\!\scri},\alpha_+,\alpha_\cT)}$ control on $u$ in $\Sigma\subset\Setb^*M$ except at $\pa\cR_{\scri,\rm out}$ and over $\iota^+\cup\cT^+$. Using \cite[Lemma~4.10(1)]{HintzVasyScrieb}, which uses the condition $\alpha_+<-\half+\alpha_{\!\scri}$, we can propagate this regularity into $\pa\cR_{\scri,\rm in,+}^\pm$; following this with real principal type propagation over $\iota^+$, we then get control on $u$ in a punctured neighborhood of $\pa\cR_{\cT,\rm in}^\pm$ over $\iota^+$. Proposition~\ref{PropWTIn}\eqref{ItWTInDir} now applies to give control on $u$ at $\pa\cR_{\cT,\rm in}^\pm$; this step uses the incoming threshold condition on $\sfs$, and it is also the reason for the presence of the second term on the right in~\eqref{EqWRDir} whose finiteness guarantees the required a priori regularity requirement~\eqref{EqWTInDirThr} of Proposition~\ref{PropWTIn}\eqref{ItWTInDir}. Real principal type propagation over $\cT^+$ gives control in a punctured neighborhood of $\pa\cR_{\cT,\rm out}^\pm$ over $\cT^+$; this uses the nontrapping assumption (cf.\ Proposition~\ref{PropWFlow}\eqref{ItWFlow4}). Proposition~\ref{PropWTOut}\eqref{ItWTOutDir} then gives control also at $\pa\cR_{\cT,\rm out}^\pm$. From there, following propagation along the flow-out $\pa\cW_{\rm out}^\pm$ of $\pa\cR_{\cT,\rm out}^\pm$ over $\iota^+$, we have control of $\chi_1 u$ in a punctured neighborhood of $\pa\cR_{\scri,\rm out}^\pm$. An application of \cite[Lemma~4.11(1)]{HintzVasyScrieb} to $\chi_1 u$ (which satisfies $P(\chi_1 u)=\chi_1 P u+[P,\chi_1]u=\chi_1(\chi_0 P u)+[P,\chi_1](\chi_\flat u)$, both terms of which are controlled by the right hand side of~\eqref{EqWRDir}), which uses $\alpha_{\!\scri}<-\half+\ubar p_1$, gives global control on $\chi_1 u$ and finishes the proof of~\eqref{EqWRDir}.

  The proof of the adjoint estimate~\eqref{EqWRAdj} proceeds in the reverse order and uses the adjoint versions of the propagation results just used. Thus, we first use \cite[Lemma~4.11(4)]{HintzVasyScrieb} to get $H_\etbop^{-\sfs+1,(-2\alpha_{\!\scri}-2,-\alpha_+-2,-\alpha_\cT)}$ control on $\chi_1\tilde u$ at $\pa\cR_{\scri,\rm out}^\pm$, which we propagate backwards using real principal type propagation. Over $\iota^+$, we thereby get control in a punctured neighborhood of $\pa\cR_{\cT,\rm out}^\pm$ over $\iota^+$, which Proposition~\ref{PropWTOut}\eqref{ItWTOutAdj} lets propagate into $\pa\cR_{\cT,\rm out}^\pm$. Upon propagating further backwards through $\cT^+$, we get control in a punctured neighborhood of $\pa\cR_{\cT,\rm in}^\pm$ over $\cT^+$, which we propagate into $\pa\cR_{\cT,\rm in}^\pm$ using Proposition~\ref{PropWTIn}\eqref{ItWTInAdj}. Propagation from there over $\iota^+$ (and from $\pa\cR_{\scri,\rm out}^\pm\cap\Setb^*_{\scri^+\cap\iota^+}M$ over $\scri^+\cap\iota^+$) then gives control in a punctured neighborhood of $\pa\cR_{\scri,\rm in}^\pm$. An application of \cite[Lemma~4.10(2)]{HintzVasyScrieb} gives control at $\pa\cR_{\scri,\rm in}^\pm$. Together with elliptic regularity, we have thus obtained quantitative control on $\tilde u$ near $\iota^+\cap\cT^+$ and near $\pa\cR_{\scri,\rm out}^\pm$. It remains to propagate this over $M\setminus(\iota^+\cap\cT^+)$ in the backwards direction for finite time until all of $\{\ft_*\geq T_\flat\}$ is covered.
\end{proof}

\subsection{Analysis of the \texorpdfstring{$\iota^+$-normal operator}{normal operator at punctured future timelike infinity}}
\label{SsWip}

We continue to work near $\iota^+\cup\cT^+$. Correspondingly, we work only with edge-3b-structures and drop $I^0$ from the notation.

Since in the notation of Definition~\ref{DefGAW}, the difference $P-P_0=\tilde P$ has decaying coefficients relative to $P_0\in(x_{\!\scri}^2\rho_+^2\CI+\cA^{(2+2\delta,2+\delta,(0,0))})\Diffetb^2(M;\upbeta^*E)$ (see Lemma~\ref{LemmaGAWStatbetb}) as an edge-3b-differential operator near $\iota^+\cup\cT^+$, the $\iota^+$-normal operator of $P$ (i.e.\ the b-normal operator at $\iota^+$ with respect to a total defining function of $\iota^+\cup\cT^+$) is equal to that of $P_0$. Furthermore, in the notation of~\eqref{EqGSOStruct}, and recalling~\eqref{EqGAWStatbetbPf}, the operators $Q\pa_{t_*}\in x_{\!\scri}^{-1}\rho_+\cA^{(4+2\delta,2+\delta,(0,0))}\Diffetb^2(M;\upbeta^*E)=\cA^{(3+2\delta,3+\delta,(0,0))}\Diffetb^2$  and $g^{0 0}\pa_{t_*}^2\in\rho_+^2\cA^{(2+2\delta,1+\delta,(0,0))}(M)=\cA^{(2+2\delta,3+\delta,(0,0))}(M)$ are of lower order compared to
\[
  N_\iota(P) := -2 \pa_{t_*}\rho\Bigl(\rho\pa_\rho-\frac{n-1}{2}-S|_{\pa X}\Bigr) + \rho^2 P_{(0)}(0,\omega,\rho D_\rho,D_\omega),
\]
where the normal operator $P_{(0)}$ of $\hat P(0)$ was introduced in~\eqref{EqGSOP0}. The operator $N_\iota(P)$ is an operator on $\R_{t_*}\times[0,\infty)_\rho\times\pa X$; we have shown:

\begin{lemma}[$N_\iota(P)$ as a normal operator of $P$]
\label{LemmaWip}
  Let $\chi\in\CI(M)$ be $1$ near $\iota^+$ and have support in a collar neighborhood of $\iota^+$. Then
  \[
    P - \chi N_\iota(P)\chi \in x_{\!\scri}^2\rho_+^2\bigl(\rho_+\CI+\cA^{(((0,0),2\ell_{\!\scri}),\ \ell_+,\ ((0,0),\ell_\cT))}\bigr)\Diffetb^2(M;\upbeta^*E).
  \]
\end{lemma}

In order to explicitly reveal the degree $-2$ homogeneity of $N_\iota(P)$ with respect to dilations (i.e.\ with respect to the b-normal vector field at $\iota^+$), we shall use two different sets of coordinates: first, we use
\[
  T := t_*^{-1},\qquad R := \frac{r}{t_*}=\frac{1}{\rho t_*}.
\]
Thus, $T$ a total defining function of $\iota^+\cup\cT^+$, and $R$ is a projective radial coordinate along $\iota^+$. We then put
\begin{equation}
\label{EqWipNorm0}
\begin{split}
  N^0_\iota(P) &:= T^{-2}N_\iota(P) \\
    &= -2 R^{-1}(R\pa_R+T\pa_T)\Bigl(R\pa_R+\frac{n-1}{2}+S|_{\pa X}\Bigr) + R^{-2}P_{(0)}(0,\omega,-R D_R,D_\omega).
\end{split}
\end{equation}
This is dilation-invariant in $T$, and its Mellin-transformed normal operator family (with respect to the total defining function $T$, as appropriate in 3b-analysis, cf.\ \eqref{EqM3bOp}--\eqref{EqM3bNDMT}) is
\begin{equation}
\label{EqWipNormMT}
  \wh{N_\iota^0}(P,\lambda) = -2 R^{-1}(R\pa_R + i\lambda)\Bigl(R\pa_R+\frac{n-1}{2}+S|_{\pa X}\Bigr) + R^{-2}P_{(0)}(0,\omega,-R D_R,D_\omega).
\end{equation}
Another useful set of coordinates is
\[
  \rho,\qquad v := R^{-1} = \rho t_* = \frac{\rho}{T},
\]
with $\rho$ a defining function of $\iota^+\setminus\scri^+$, and we then put
\[
  N^1_\iota(P) := \rho^{-2}N_\iota(P) = -2\pa_v\Bigl( \rho\pa_\rho+v\pa_v - \frac{n-1}{2} - S|_{\pa X}\Bigr) + P_{(0)}(0,\omega,\rho D_\rho+v D_v,D_\omega).
\]
For a function $u$ which depends only on $(v,\omega)$ (or, equivalently, only on $(R,\omega)$), we record the relationship
\begin{equation}
\label{EqWipNormRel}
  \wh{N_\iota^0}(P,\lambda)u = T^{-i\lambda}N_\iota^0(P)(T^{i\lambda}u) = T^{-i\lambda-2}\rho^2 N_\iota^1(P)\bigl(T^{i\lambda}u\bigr) = v^{i\lambda+2}\wh{N_\iota^1}(P,\lambda)\bigl(v^{-i\lambda}u\bigr),
\end{equation}
where the Mellin-transformed normal operator family of $N_\iota^1$ (now with respect to $\rho=T v$) is
\begin{equation}
\label{EqWiprhovNopMT}
  \wh{N_\iota^1}(P,\lambda) = -2\pa_v\Bigl(i(v D_v+\lambda)-\frac{n-1}{2}-S|_{\pa X}\Bigr) + P_{(0)}(0,\omega,v D_v+\lambda,D_\omega)
\end{equation}

Returning to $\wh{N_\iota^0}(P,\lambda)$, note that since $T^{-2}P_0\in x_{\!\scri}^2\rho_\cT^{-2}\Diffetb^2(M;\upbeta^*E)$ (modulo terms with decaying coefficients at $\iota^+$) is a (weighted at $\scri^+$ and $\cT^+$) edge-3b-operator, the structural results from \S\ref{SsMeb} (near $\scri^+$) and \S\ref{SsM3b} (near $\cT^+$) combine to give
\begin{equation}
\label{EqWipNormMem}
  \wh{N^0_\iota}(P,\lambda) \in x_{\!\scri}^2\rho_\cT^{-2}\Diff_{0,\bop}^2(\iota^+;\upbeta^*E|_{\pa X}) = \Diff_{0,\bop}^{2,(-2,2)}(\iota^+;\upbeta^*E|_{\pa X}),
\end{equation}
where we abuse notation and write $\upbeta$ also for the restriction $\iota^+\to\pa X$ of the blow-down map $M\to M_0$.

\begin{rmk}[Structure of $\iota^+$ and $\wh{N_\iota^0}(P,\lambda)$]
\label{RmkWipNormStruct}
We stress that we have $\iota^+=[0,\infty]_{R,1/2}\times\pa X$, where we write $[0,\infty]_{R,1/2}$ is equal to $[0,\infty]$ as a set, but the smooth structure is such that $R^{-1/2}$ is a smooth local boundary defining function of $R^{-1}(0)=\iota^+\cap\scri^+\subset\iota^+$, cf.\ the square root blow-up in Definition~\ref{DefGAMfd}, while a local defining function of $R^{-1}(0)=\iota^+\cap\cT^+\subset\iota^+$ is $R$. The membership~\eqref{EqWipNormMem} then follows also directly from~\eqref{EqWipNormMT} (see also~\eqref{EqWipBd0} below).
\end{rmk}

Furthermore, for any $\gamma_+\in\R$, the families
\begin{equation}
\label{EqWipNormHi}
  (0,1)\ni h \mapsto \wh{N_\iota^0}(P,-i\gamma_+\pm h^{-1})
\end{equation}
define elements of $\Diff_{0,\cop,\semi}^{2,(-2,2,2,2)}(\iota^+;\upbeta^*E|_{\pa X})$; here, we combine the notation from the semiclassical 0- and cone-calculi and write
\begin{equation}
\label{EqWipDiff0ch}
  \Diff_{0,\cop,\semi}^{m,(\gamma,l,\alpha,b)}(\iota^+) = x_{\!\scri}^{-\gamma}\Bigl(\frac{\rho_\cT}{\rho_\cT+h}\Bigr)^{-l}(\rho_\cT+h)^{-\alpha}\Bigl(\frac{h}{\rho_\cT+h}\Bigr)^{-b}\Diff_{0,\cop,\semi}^m(\iota^+),
\end{equation}
with semiclassical cone behavior near $\iota^+\cap\cT^+$ and semiclassical 0-behavior near $\iota^+\cap\scri^+$.

\begin{prop}[Estimates for $\wh{N_\iota^0}(P,\lambda)$: bounded frequencies]
\label{PropWipBd}
  Fix a positive b-density on $\iota^+$ to define $L^2$- and Sobolev spaces. Let $s\in\R$ or $s\in\CI({}^{0,\bop}S^*\iota^+)$. Let $\gamma_{\!\scri},\gamma_\cT\in\R$, $\lambda\in\C$, and suppose that
  \begin{equation}
  \label{EqWipBdWeights}
    -\Im\lambda<1+\beta^+,\qquad
    -\Im\lambda < \gamma_{\!\scri} < \frac{n-1}{2}+\ubar S,\qquad
    \gamma_\cT\in(-\beta^+,-\beta^-),
  \end{equation}
  where $\beta^\pm$ are as in Definition~\usref{DefGSOSpec}, and $\ubar S$ is as in Definition~\usref{DefStEstThr}. Then the operator
  \begin{equation}
  \label{EqWipBdOp}
    \wh{N_\iota^0}(P,\lambda) \colon H_{0,\bop}^{s,(2\gamma_{\!\scri},\gamma_\cT)}(\iota^+;\upbeta^*E|_{\pa X}) \to H_{0,\bop}^{s-2,(2\gamma_{\!\scri}+2,\gamma_\cT-2)}(\iota^+;\upbeta^*E|_{\pa X})
  \end{equation}
  is invertible, with inverse uniformly bounded when $\lambda$ varies in a compact set of $\lambda\in\C$ subject to~\eqref{EqWipBdWeights}. Here, the space $H_{0,\bop}^{s,(2\gamma_{\!\scri},\gamma_\cT)}(\iota^+)$ consists of all distributions on $(\iota^+)^\circ$ which lie in $H_{0,\loc}^{s,2\gamma_{\!\scri}}(\iota^+\setminus\cT^+)$, resp.\ $H_{\bop,\loc}^{s,\gamma_\cT}(\iota^+\setminus\scri^+)$, upon multiplication by a smooth function on $\iota^+$ which vanishes near $\cT^+$, resp.\ $\scri^+$.
\end{prop}

The first upper bound on $-\Im\lambda$ in~\eqref{EqWipBdWeights} (which corresponds to the $L^\infty$-decay rate at $(\iota^+)^\circ$) matches the second ($\iota^+$-decay) order in~\eqref{EqStCoSol}. In a similar manner, the second upper bound in~\eqref{EqWipBdWeights} is directly related to analogous bounds on the $\scri^+$- and $\iota^+$-weights under which one can propagate edge-b-regularity through the corner $\scri^+\cap\iota^+$; see \cite[Lemmas~4.10, 4.11 and Remark~8.3]{HintzVasyScrieb}. In this context, note also that if $\mu_\bop$ is a positive b-density on $M$, then
\begin{equation}
\label{EqWipBdbDensity}
  H_\etbop^{s,(2\alpha_{\!\scri},\alpha_+,\alpha_\cT)}(M) = H_\etbop^{s,\bigish(2(\alpha_{\!\scri}+\frac{n}{2}),\;\alpha_++\frac{n+1}{2},\;\alpha_\cT+\frac12\bigish)}(M,\mu_\bop).
\end{equation}
Correspondingly, we will use Proposition~\ref{PropWipBd} with
\begin{equation}
\label{EqWipBdRel}
  -\Im\lambda=\alpha_++\frac{n+1}{2},\qquad
  \gamma_{\!\scri}=\alpha_{\!\scri}+\frac{n}{2},\qquad
  \gamma_\cT=\alpha_\cT-\alpha_+-\frac{n}{2},
\end{equation}
for various choices of $\alpha_+,\alpha_{\!\scri},\alpha_\cT$, with the bounds~\eqref{EqWipBdWeights} translating into
\begin{equation}
\label{EqWipBdWeights2}
\begin{alignedat}{2}
  \alpha_+&<-\frac{n-1}{2}+\beta^+,&\qquad
  \alpha_+&<-\frac12+\alpha_{\!\scri}, \\
  \alpha_{\!\scri}&<-\frac12+\ubar S,&\qquad
  \alpha_+-\alpha_\cT&\in\Bigl(-\frac{n}{2}+\beta^-,-\frac{n}{2}+\beta^+\Bigr).
\end{alignedat}
\end{equation}
In particular, $\alpha_\cT+\half<\alpha_+-(-\frac{n}{2}+\beta^-)+\half<1+\beta^+-\beta^-$, which matches the $L^\infty$-decay rate at $(\cT^+)^\circ$ in~\eqref{EqStCoSol}.

\begin{proof}[Proof of Proposition~\usref{PropWipBd}]
  \pfstep{Fredholm property of $\wh{N_\iota^0}(P,\lambda)$.} In $R=\frac{r}{t_*}>1$, we pass to the local defining function $x_{\!\scri}:=R^{-1/2}$ of $\scri^+\cap\iota^+$ in the expression~\eqref{EqWipNormMT} and compute
  \begin{equation}
  \label{EqWipBd0}
    x_{\!\scri}^{-2}\wh{N_\iota^0}(P,\lambda) = -\half(x_{\!\scri}\pa_{x_{\!\scri}}-2 i\lambda)\Bigl[x_{\!\scri}\pa_{x_{\!\scri}}-2\Bigl(\frac{n-1}{2}+S|_{\pa X}\Bigr)\Bigr] + x_{\!\scri}^2 P_{(0)}(0,\omega,\half x_{\!\scri} D_{x_{\!\scri}},D_\omega).
  \end{equation}
  Since the principal part of $P_{(0)}$ at $x_{\!\scri}=0$ is $(\half x_{\!\scri} D_{x_{\!\scri}})^2+\slDelta$ (with $\slDelta=\Delta_\slg$), the operator $x_{\!\scri}^{-2}\wh{N_\iota^0}(P,\lambda)$ is elliptic as a 0-operator. (This can also be deduced from the fact that the principal symbol of $\wh{N_\iota^0}(P,\lambda)$ is the pullback of that of $T^{-2}P$ along the inclusion map ${}^{0,\bop}T^*\iota^+\hra\Tetb^*_{\iota^+}M$.) A local parametrix can be constructed in the 0-calculus \cite{MazzeoMelroseHyp,Hintz0Px}; this was done in (the proof of) \cite[Theorem~8.2]{HintzVasyScrieb}. For any two cutoff functions $\chi_1,\chi_2\in\CIc(\iota^+\setminus\cT^+)$ which are identically $1$ near $\iota^+\cap\scri^+$ and satisfy $\supp\chi_1\subset\chi_2^{-1}(1)$, \cite[Theorem~8.2]{HintzVasyScrieb} then gives the following a priori estimate (recording only the weight at $\scri^+$): for $\gamma_{\!\scri}$ as in~\eqref{EqWipBdWeights}, and for any fixed $N$, we have, dropping the bundle $\upbeta^*E|_{\pa X}$ from the notation,\footnote{The operator $x_{\!\scri}^{-2}\wh{N_\iota^0}(P,\lambda)$ is denoted $\wh{P_+}(\lambda)$ in the reference, and $\ubar p_{1,+}$ in the reference is equal to $\ubar S$ in present notation (cf.\ \cite[Definition~4.3]{HintzVasyScrieb} and Definition~\ref{DefGAW}\eqref{ItGAWEdgeN}).}
  \begin{equation}
  \label{EqWipBd0Est}
    \|\chi_1 u\|_{H_0^{s,2\gamma_{\!\scri}}(\iota^+)} \leq C\Bigl( \|\chi_2 x_{\!\scri}^{-2}\wh{N_\iota^0}(P,\lambda)u\|_{H_0^{s-2,2\gamma_{\!\scri}}(\iota^+)} + \|\chi_2 u\|_{H_0^{-N,-N}(\iota^+)}\Bigr).
  \end{equation}
  Moreover, any distribution $u$ with $\chi_1 u\in H_0^{s,2\gamma_{\!\scri}}(\iota^+)$ and $\chi_2\wh{N_\iota^0}(P,\lambda)u=0$ is automatically conormal at $\scri^+$ and satisfies
  \begin{equation}
  \label{EqWipBdCon}
    \chi_1 u\in\cA^{n-1+2\ubar S-\eps}(\iota^+\setminus\cT^+) = (x_{\!\scri}^2)^{\frac{n-1}{2}+\ubar S-\eps}\cA^0(\iota^+\setminus\cT^+)
  \end{equation}
  for all $\eps>0$.

  Turning to the conic point $R=0$, inspection of~\eqref{EqWipNormMT} shows that $R^2\wh{N_\iota^0}(P,\lambda)$ is an elliptic b-operator for $R\in[0,\infty)$, and its b-normal operator at $R=0$ is $P_{(0)}(0,\omega,-R D_R,D_\omega)$. By Lemma~\ref{LemmaStEst0}, the interval $(-\beta^+,-\beta^-)$ is disjoint from its boundary spectrum. Elliptic b-theory thus implies that for cutoffs $\chi_3,\chi_4\in\CIc(\iota^+\setminus\scri^+)$ which are $1$ near $\iota^+\cap\cT^+$ and satisfy $\supp\chi_3\subset\chi_4^{-1}(1)$, we have the a priori estimate
  \[
    \|\chi_3 u\|_{\Hb^{s,\gamma_\cT}(\iota^+)} \leq C\Bigl( \|\chi_4 R^2\wh{N_\iota^0}(P,\lambda)u\|_{\Hb^{s-2,\gamma_\cT}(\iota^+)} + \|\chi_4 u\|_{\Hb^{-N,-N}(\iota^+)}\Bigr).
  \]
  We may take $\chi_3=1-\chi_1$. If $\chi_3 u\in\Hb^{s,\gamma_\cT}(\iota^+)$ and $\chi_4\wh{N_\iota^0}(P,\lambda)u=0$, then
  \begin{equation}
  \label{EqWipBdCon2}
    \chi_3 u \in \cA^{-\beta^--\eps}(\iota^+\setminus\scri^+)
  \end{equation}
  for all $\eps>0$ by b-elliptic regularity and using the above information about the boundary spectrum of $R^2\wh{N_\iota^0}(P,\lambda)$. Putting this estimate together with~\eqref{EqWipBd0Est} then implies the semi-Fredholm estimate
  \[
    \|u\|_{H_{0,\bop}^{s,(2\gamma_{\!\scri},\gamma_\cT)}(\iota^+)} \leq C\Bigl( \|\wh{N_\iota^0}(P,\lambda)u\|_{H_{0,\bop}^{s-2,(2\gamma_{\!\scri}+2,\gamma_\cT-2)}(\iota^+)} + \|u\|_{H_{0,\bop}^{-N,(-N,-N)}(\iota^+)} \Bigr).
  \]
  Similarly, one can prove an adjoint estimate
  \[
    \|\tilde u\|_{H_{0,\bop}^{-s+2,(-2\gamma_{\!\scri}-2,-\gamma_\cT+2)}(\iota^+)} \leq C\Bigl( \|\wh{N_\iota^0}(P,\lambda)^*\tilde u\|_{H_{0,\bop}^{-s,(-2\gamma_{\!\scri},-\gamma_\cT)}(\iota^+)} + \|\tilde u\|_{H_{0,\bop}^{-N,(-N,-N)}(\iota^+)} \Bigr),
  \]
  where the adjoint is defined with respect to a positive b-density and any smooth fiber inner product on $E|_{\pa X}$. Together, these estimate prove that the operator~\eqref{EqWipBdOp} is Fredholm.

  \pfstep{Injectivity.} We next prove the injectivity of~\eqref{EqWipBdOp} by reducing it to the injectivity of the transition face normal operator, which is assumption~\eqref{ItGSOSpectf} in Definition~\ref{DefGSOSpec}. This reduction is possible since the $\iota^+$-normal operator of $P$, beyond being homogeneous with respect to scaling, is also time-translation invariant; see also the algebraic considerations in \cite[Remark~3.27]{Hintz3b}, which we supplement here with corresponding considerations on the level of function spaces. Thus, given $u$ in the kernel of~\eqref{EqWipBdOp}, we have $u\in\cA^{(n-1+2\ubar S-\eps,-\beta^--\eps)}(\iota^+)$ for all $\eps>0$; near $\scri^+$, this is~\eqref{EqWipBdCon}, and near $\cT^+$, this follows from the above discussion of the boundary spectrum. Since $\iota^+$ is the square root blow-up of $[0,\infty]_v\times\pa X$ at $v=0$, we have
  \[
    u\in\cA^{(\frac{n-1}{2}+\ubar S-\eps,-\beta^--\eps)}([0,\infty]_v\times\pa X)
  \]
  by~\eqref{EqWipBdCon} and \eqref{EqWipBdCon2}. By~\eqref{EqWipNormRel}, this implies
  \begin{equation}
  \label{EqWipuflatMem}
    \wh{N_\iota^1}(P,\lambda)u^\flat = 0,\qquad u^\flat:=v^{-i\lambda}u \in \cA^{(\frac{n-1}{2}+\ubar S+\Im\lambda-\eps,-\beta^--\Im\lambda-\eps)}([0,\infty]_v\times\pa X)\quad\forall\,\eps>0.
  \end{equation}
  In view of~\eqref{EqWipBdWeights}, the weight of $u^\flat$ at $v=0$ is positive for sufficiently small $\eps>0$. The extension of $u^\flat$ by $0$ to $v<0$,
  \begin{equation}
  \label{EqWipsharp}
    u^\sharp(v,\omega) := \begin{cases} u^\flat(v,\omega), & v>0, \\ 0, & v\leq 0, \end{cases}
  \end{equation}
  therefore satisfies $u^\sharp,v D_v u^\sharp\in H_\loc^{\frac12+\delta}((-1,1)_v\times\pa X)$ for all $\delta\in(0,\frac{n-1}{2}+\ubar S+\Im\lambda)$. Regarding $\wh{N_\iota^1}(P,\lambda)$ as a smooth coefficient differential operator on $\R_v\times\pa X$, given by~\eqref{EqWiprhovNopMT} also for $v\leq 0$, we infer that
  \begin{equation}
  \label{EqWipInjEqn}
    \wh{N_\iota^1}(P,\lambda)u^\sharp=0\quad\text{on}\quad \R_v\times\pa X
  \end{equation}
  in the sense of distributions; indeed, setting $f^\sharp:=\wh{N_\iota^1}(P,\lambda)u^\sharp$, we have $\supp f^\sharp\subset\{0\}$ by construction, so $f^\sharp$ is a sum of differentiated $\delta$-distributions at $v=0$, and in view of $f^\sharp\in H_\loc^{-\frac12+\delta}((-1,1)_v\times\pa X)$ must therefore vanish.

  The advantage of working with the extended equation~\eqref{EqWipInjEqn} is that the Fourier transform in $v$ becomes available; that this is a useful tool can be seen from the fact that the $v$-Fourier transform is a rescaling of the $t_*$-Fourier transform. Denoting by $\hat r$ the Fourier-dual variable to $v$, we then note that~\eqref{EqWipuflatMem} implies\footnote{This is a standard result on Fourier transforms of symbols or distributions conormal at $0$, see e.g.\ \cite[Lemma~2.25]{Hintz3b}.}
  \begin{equation}
  \label{EqWiwpm}
    w_\pm(\hat r,\omega) := \wh{u^\sharp}(\pm\hat r,\omega)|_{(0,\infty)\times\pa X} \in \cA^{\bigish(((0,0),-\beta^--\Im\lambda-1-\eps),\frac{n-1}{2}+\ubar S+\Im\lambda+1-\eps\bigish)}([0,\infty]\times\pa X).
  \end{equation}
  Moreover, with the convention $\cF_{v\to\hat r}u^\sharp=\int e^{i\hat r v}u^\sharp(v,\omega)\,\dd v$ of the Fourier transform (consistently with our convention~\eqref{EqIFT} for the $t_*$-Fourier transform), we have
  \begin{equation}
  \label{EqWipmFTvr}
    \cF_{v\to\hat r}\circ\pa_v=-i\hat r\circ\cF_{v\to\hat r},\qquad
    \cF_{v\to\hat r}\circ v=-i\pa_{\hat r}\circ\cF_{v\to\hat r};
  \end{equation}
  transforming~\eqref{EqWiprhovNopMT} accordingly, we find that $w_\pm$ satisfies the equation
  \[
    \Bigl[\pm 2 i\hat r\Bigl( i(-D_{\hat r}\hat r+\lambda)-\frac{n-1}{2}-S|_{\pa X}\Bigr) + P_{(0)}(0,\omega,-D_{\hat r}\hat r+\lambda,D_\omega)\Bigr] w_\pm = 0.
  \]
  Finally introducing
  \begin{equation}
  \label{EqWipmflat}
    w_\pm^\flat := \hat r^{-i\lambda+1}w_\pm \in \cA^{(((1-i\lambda,0),-\beta^--\eps),\alpha)}([0,\infty]\times\pa X),\qquad \alpha=\frac{n-1}{2}+\ubar S-\eps,
  \end{equation}
  this implies, upon division by $\hat r^2$,
  \begin{equation}
  \label{EqWiwpmEq}
    \Bigl[\pm 2 i\hat r^{-1}\Bigl( -\hat r\pa_{\hat r}-\frac{n-1}{2}-S|_{\pa X}\Bigr) + \hat r^{-2}P_{(0)}(0,\omega,-\hat r D_{\hat r},D_\omega)\Bigr] w_\pm^\flat = 0.
  \end{equation}
  Upon passing to $\hat\rho=\hat r^{-1}$, this is precisely the statement that $N_\tface^\pm(P)w_\pm^\flat=0$ (cf.\ \eqref{EqGSONtf}). By assumption, we have $1+\Im\lambda>-\beta^+$, and therefore $w_\pm^\flat$, regarded now as a function of $(\hat\rho,\omega)$, lies $\cA^{(\alpha,-\beta)}([0,\infty]_{\hat\rho}\times\pa X)$ for some $\beta\in(\beta^-,\beta^+)$; by Definition~\ref{DefGSO}\eqref{ItGSOSpectf}, this implies $w_\pm^\flat\equiv 0$. We now work our way back: we infer $w_\pm\equiv 0$, which in view of the definition~\eqref{EqWiwpm} implies $\supp\wh{u^\sharp}\subset\{0\}$. Therefore, $u^\sharp(v,\omega)$ is a polynomial in $v$; since it is at the same time supported in $v\geq 0$ by~\eqref{EqWipsharp}, this forces $u^\sharp\equiv 0$. Therefore, $u^\flat\equiv 0$, and by~\eqref{EqWipuflatMem}, we finally deduce $u\equiv 0$, finishing the proof of the injectivity of~\eqref{EqWipBdOp}.

  \pfstep{Surjectivity.} There are two ways to proceed: one way is to show the triviality of $\ker\wh{N_\iota^0}(P,\lambda)^*$ on $H_{0,\bop}^{-s+2,(-2\gamma_{\!\scri}-2,-\gamma_\cT+2)}(\iota^+)$, which can be reduced to the triviality of the kernel of the adjoint of the transition face normal operator in Definition~\ref{DefGSO}\eqref{ItGSOSpectf} by arguments similar to the above; we leave the details to the reader. Another way, which we follow, is to show that $\wh{N_\iota^0}(P,\lambda)$ is invertible for some $\lambda$ satisfying~\eqref{EqWipBdWeights}, which implies that~\eqref{EqWipBdOp} has index $0$ and implies its surjectivity in view of its injectivity. We prove the desired invertibility by semiclassical means in Proposition~\ref{PropWipHi} below.
\end{proof}

We now turn to the invertibility and uniform estimates for the semiclassical family~\eqref{EqWipNormHi} (which also completes the proof of Proposition~\ref{PropWipBd}). We use spaces of distributions on $\iota^+$ which lie in semiclassical 0-Sobolev spaces near $\scri^+$ and in semiclassical cone Sobolev spaces near $\cT^+$; matching~\eqref{EqWipDiff0ch}, we write
\[
  H_{0,\cop,h}^{s,(\gamma,l,\alpha,b)}(\iota^+)
\]
for the function space with differential order $s$, decay order $\gamma$ at $\scri^+$, decay order $l$ at the cone face over $\cT^+$, decay order $\alpha$ at the transition face over $\cT^+$, and semiclassical order $b$. (This amalgamates the spaces~\eqref{EqMFch} and \eqref{EqMF0}.)  That is, localization, via multiplication with a cutoff function, to a neighborhood of $\scri^+$, resp.\ $\cT^+$ gives an element of $H_{0,h}^{s,\gamma,b}$, resp.\ $H_{\cop,h}^{s,l,\alpha,b}$.

\begin{prop}[Estimates for $\wh{N_\iota^0}(P,\lambda)$: high frequencies]
\label{PropWipHi}
  Let $\gamma_{\!\scri},\gamma_+,\gamma_\cT\in\R$, and suppose that\footnote{There is no condition $\gamma_+<1+\beta^+$ here; this was used in the proof of Proposition~\ref{PropWipBd} only in the final step of the proof of injectivity of $\wh{N_\iota^0}(P,\lambda)$. In all concrete settings considered in~\S\ref{SE}, we happen to have $\frac{n-1}{2}+\ubar S\leq 1+\beta^+$ in any case, in which case the condition $\gamma^+<1+\beta^+$ follows from the assumptions of Proposition~\ref{PropWipHi}.}
  \[
    \gamma_+<\gamma_{\!\scri}<\frac{n-1}{2}+\ubar S,\qquad
    \gamma_\cT\in(-\beta^+,-\beta^-).
  \]
  Let $\sfs,\sfs_0$ be two variable order functions as in Lemma~\usref{LemmaWROrder} where $\alpha_+=\gamma_+-\frac{n+1}{2}$ and $\alpha_\cT=\gamma_++\gamma_\cT-\half$. Denote by $\sfs_\infty\in\CI({}^{0\chop}S^*\iota^+)$ and $\sfs_\semi\in\CI(\ol{{}^{0\chop}T^*_\sface}\iota^+)$ the orders at fiber infinity and the semiclassical face induced by $\sfs|_{\Setb^*_{\iota^+}M}$ via pullback along~\eqref{EqM3bSclMaps} (near $\iota^+\cap\cT^+$) and~\eqref{EqMebSclMaps} (near $\iota^+\cap\scri^+$). Then there exists $h_0>0$ so that for $0<h<h_0$, the operator $\wh{N_\iota^0}(P,-i\gamma_+ \pm h^{-1})$ is invertible. (The same $h_0$ works also for all nearby values of $\gamma_+$.) Its inverse satisfies the uniform bound
  \begin{equation}
  \label{EqWipHi}
  \begin{split}
    &\|u\|_{H_{0,\cop,h}^{\sfs_\infty,(2\gamma_{\!\scri},\gamma_\cT,\gamma_\cT,\sfs_\semi)}(\iota^+;\upbeta^*E|_{\pa X})} \\
      &\qquad \leq C\|\wh{N_\iota^0}(P,-i\gamma_+ \pm h^{-1})u\|_{H_{0,\cop,h}^{\sfs_\infty-2,(2\gamma_{\!\scri}+2,\gamma_\cT-2,\gamma_\cT-2,\sfs_\semi-1)}(\iota^+;\upbeta^*E|_{\pa X})},\qquad 0<h<h_0.
  \end{split}
  \end{equation}
\end{prop}
\begin{proof}
  \pfstep{Symbolic estimate.} In view of the relationship between edge-b-analysis on $M$ near $\iota^+\cap\scri^+$ and semiclassical 0-analysis on $\iota^+$ described on the level of phase spaces and thus of principal symbols and Hamiltonian vector fields in and after~\eqref{EqMebSclMap} (and also \eqref{EqMebSclMaps}), and on the level of subprincipal symbols in~\eqref{EqMebSclSubpr}, the proofs (via positive commutator arguments) of symbolic edge-b-propagation estimates near $\scri^+$ in \cite[\S4.2]{HintzVasyScrieb} give, with minor modifications, symbolic semiclassical 0-propagation estimates near $\scri^+\cap\iota^+$; this is discussed in detail in \cite[\S8.2]{HintzVasyScrieb}.

  Similarly, the positive commutator proofs of symbolic 3b-estimates near $\cT^+\cap\iota^+$, presented in~\S\ref{SsWT}, give symbolic semiclassical cone propagation estimates near $\cT^+\cap\iota^+$; the relationships between the orders $\gamma_{\!\scri},\gamma_+,\gamma_\cT$ here and the orders $\alpha_{\!\scri},\alpha_+,\alpha_\cT$ in these propagation estimates are given by~\eqref{EqWipBdRel} with $\gamma_+$ in place of $-\Im\lambda$. We illustrate this for the outgoing radial point estimate over $\pa\cT^+$. Recall that if we write 3b-covectors as $\sigma_\tbop\frac{\dd t_*}{r}+\xi_\tbop\frac{\dd r}{r}+\eta_\tbop\,\dd\omega$ as in~\eqref{EqWtbCoords}, then $\cR_{\cT,\rm out}^\pm=\{\rho_\cT=\rho_+=\xi_\tbop=\eta_\tbop=0,\pm\sigma_\tbop<0\}$ from~\eqref{EqWRadT}. Write $\chop$-covectors near the semiclassical face $\sface\subset\iota^+_\chop$ of the semiclassical cone single space of $\iota^+$, and near $\cT^+\cap\iota^+$, as
  \[
    \xi_\chop\frac{R}{h}\frac{\dd R}{R} + \eta_\chop\frac{R}{h}\dd\omega.
  \]
  Recall $\rho_\cT=R=\frac{r}{t_*}$, $\rho_+=\frac{1}{r}$, $\omega\in\pa X$, and the total defining function $\rho_0=T=t_*^{-1}$ of $\iota^+\cup\cT^+$. The map~\eqref{EqM3bSclMap} then has the local coordinate expression
  \begin{align}
  \label{EqWipHiIotach}
    {}^\chop\iota \colon (h,R,\omega,\xi_\chop,\eta_\chop) &\mapsto \mp h^{-1}\frac{\dd t_*}{t_*} + \xi_\chop\frac{R}{h}\frac{\dd R}{R} + \eta_\chop\frac{R}{h}\dd\omega \\
      & = \Bigl(\frac{h}{R}\Bigr)^{-1} \Bigl( (\mp 1-R\xi_\chop)\frac{\dd t_*}{r} + \xi_\chop\frac{\dd r}{r}+\eta_\chop\,\dd\omega\Bigr) \in \Ttb^*_{(R,0,\omega)}M, \nonumber
  \end{align}
  where we use the local coordinates $R,\omega$ (with $(\rho_\cT,\rho_+,\omega)=(R,0,\omega)$) on $\iota^+\subset M$. Thus, ${}^\chop\iota^*(\sigma_\tbop,\xi_\tbop,\eta_\tbop)=\tilde h^{-1}(\mp 1-R\xi_\chop,\xi_\chop,\eta_\chop)$, where $\tilde h:=\frac{h}{R}$ is a local defining function (in $h\lesssim R$) of $\sface\subset(\iota^+\setminus\scri^+)_\chop$. Therefore, the outgoing radial set $\cR_{\cT,\rm out}^\pm$ lies in the range of ${}^\chop\iota$ for the matching choice of the sign of $h^{-1}$ in~\eqref{EqM3bSclMap} (and here in~\eqref{EqWipHiIotach}); the preimage of $\pa\cR_{\cT,\rm out}^\pm$ under the induced map $\ol{\Tch^*_\sface(\iota^+\setminus\scri^+)}\to\Stb^*_{\iota^+\setminus\scri^+}M$ (see~\eqref{EqM3bSclMaps}) is then
  \[
    {}^\chop\cR_{\cT,\rm out} = \{ (\tilde h,R,\omega;\xi_\chop,\eta_\chop) \colon \tilde h=R=0,\ \xi_\chop=\eta_\chop=0 \}.
  \]

  For the sake of definiteness, let us consider the top sign. Since in the 3b-characteristic set over $\iota^+$, the set $\pa\cR_{\cT,\rm out}^+$ is a source for a sign-preserving rescaling of $H_{T^{-2}p}$, the same is true for ${}^\chop\cR_{\cT,\rm out}$ with respect to the rescaled Hamiltonian flow of the $\chop$-principal symbol of $\wh{N_\iota^0}(P,-i\gamma_+\pm h^{-1})$ (which is the pullback of $T^{-2}p$ along ${}^\chop\iota$). Now, the proof of Proposition~\ref{PropWTOut}\eqref{ItWTOutDir} with weights $\alpha_+,\alpha_\cT$ with respect to the metric density $|\dd g|$, i.e.\ with weights $\gamma_+=\alpha_++\frac{n+1}{2}$, $\tilde\gamma_\cT=\alpha_\cT+\frac{1}{2}$ with respect to a positive b-density on $M$, is based on a positive commutator estimate; for the operator $T^{-2}P=\rho_+^{-2}\rho_\cT^{-2}P\in\Difftb^{2,(0,2)}(M;\upbeta^*E)$ and working with b-densities on $M$, this estimate involves a commutant $\check a^2$ where $\check a$ is the product of the weight
  \[
    \rho_\infty^{-\sfs+\frac12}\rho_+^{-\gamma_+}\rho_\cT^{-\tilde\gamma_\cT+1} = (-\sigma_\tbop)^{s-\frac12} \rho_0^{-\gamma_+} \rho_\cT^{-\gamma_\cT+1},\qquad \gamma_\cT:=\tilde\gamma_\cT-\gamma_+,
  \]
  (cf.\ \eqref{EqWTInCommutant}, with the shifted orders reflecting the re-weighting of the operator and of the volume density used) with appropriate cutoff functions; the definition of $\gamma_\cT$ is consistent with~\eqref{EqWipBdRel}.

  In order to translate this to the $\chop$-setting, we can use one of two equivalent methods. In the first method, one reduces to the case $\gamma_+=0$ by replacing $T^{-2}P$ by $\rho_0^{-\gamma_+}T^{-2}P\rho_0^{\gamma_+}=T^{-2}P + \rho_0^{-\gamma_+}[T^{-2}P,\rho_0^{\gamma_+}]$; this leaves the principal symbol unchanged, but the second summand (which is a first order 3b-operator) has principal symbol $-i\gamma_+ \rho_0^{-1}H_{T^{-2}p}\rho_0$. After this reduction, we can use as a commutant for the $\chop$-estimate at ${}^\chop\cR_{\cT,\rm out}$ the same expression as in the 3b-setting, except without the presence of a localization in $\rho_+$, pulled back along ${}^\chop\iota$. The upshot is that one can propagate semiclassical cone regularity out of ${}^\chop\cR_{\cT,\rm out}$ provided $2(-\gamma_\cT+1+s-\frac12)+(-2\ubar S-n)<0$ there, where $s=\sfs|_{\pa\cR_{\cT,\rm out}^+}$; note here that $\rho_\cT^2\rho_0^{-1}H_{T^{-2}p}\rho_0=0$ at $\cR_{\cT,\rm out}^\pm$ by~\eqref{EqWRadTHamLinOut} (since $\rho_0=\rho_+\rho_\cT$), and the additional shift by $n$ arises from computing the imaginary part of $T^{-2}P$ (or $\wh{N_\iota^0}(P,\pm h^{-1})$) with respect to the b-density $\rho_+^{n+1}\rho_\cT|\dd g|=\rho_0\rho_+^n|\dd g|$ on $M$ instead of that of $P$ with respect to $|\dd g|$ (the constant $-n$ then arising from $\half$ times the derivative of the weight $\rho_0\rho_+^n$ along~\eqref{EqWRadTHamLinOut}). This condition is equivalent to $\sfs-\gamma_\cT<\frac{n-1}{2}+\ubar S$ at $\pa\cR_{\cT,\rm out}^+$, or equivalently to
\begin{equation}
\label{EqWipHiSymbThr}
  \sfs_\semi-\gamma_\cT<\frac{n-1}{2}+\ubar S\quad\text{at}\quad {}^\chop\cR_{\cT,\rm out}.
\end{equation}
But since $\sfs+\alpha_+-\alpha_\cT=(\sfs+\gamma_+-\frac{n+1}{2})-(\gamma_++\gamma_\cT-\half)=\sfs-\gamma_\cT-\frac{n}{2}$, this condition is guaranteed by the choice of $\sfs$ in Lemma~\ref{LemmaWROrder}.

  The other, equivalent, method is to leave $\gamma_+$ unchanged, and still use as the commutant for the semiclassical cone estimate the same commutant as in the 3b-setting but without the weight $\rho_0^{-\gamma_+}$ and without localization in $\rho_+$; in this case, the imaginary part of $\wh{N_\iota^0}(P,-i\gamma_++h^{-1})$ gains a contribution via the second term in~\eqref{EqM3bSubpr}, which, as already observed above, vanishes at the outgoing radial set.

  Analogous arguments apply at the radial point estimate corresponding to the incoming radial set $\cR_{\cT,\rm in}^\pm$.

  Of course, and most simply, one can also prove these radial point estimates directly using positive commutator arguments relying on the explicit form~\eqref{EqWipNormMT} (with $\lambda=-i\gamma_+\pm h^{-1}$) of the operator; this would lead to the same conclusions and threshold conditions, but it would not explain, on a conceptual level, why the threshold conditions here are the same as for the 3b-propagation results in Propositions~\ref{PropWTIn} and \ref{PropWTOut}. We refer the reader to \cite[\S4]{HintzConicProp} for detailed proofs of radial point estimates in the $\chop$-setting.

  Combining these radial point estimates with real principal type propagation and elliptic estimates (and using the global structure of the null-bicharacteristic flow over $\iota^+$ recorded in Lemma~\ref{LemmaWFlowI}, see also the top left quadrangle in Figure~\ref{FigWFlow}), we obtain the estimate (dropping the bundle $\upbeta^*E|_{\pa X}$ from the notation)
  \begin{equation}
  \label{EqWipHiSymb}
  \begin{split}
    \|u\|_{H_{0,\cop,h}^{\sfs_\infty,(2\gamma_{\!\scri},\gamma_\cT,\gamma_\cT,\sfs_\semi)}(\iota^+)} &\leq C\Bigl(\|\wh{N_\iota^0}(P,-i\gamma_+ \pm h^{-1})u\|_{H_{0,\cop,h}^{\sfs_\infty-2,(2\gamma_{\!\scri}+2,\gamma_\cT-2,\gamma_\cT-2,\sfs_\semi-1)}(\iota^+)} \\
      &\hspace{13em} + \|u\|_{H_{0,\cop,h}^{-N,(-N,\gamma_\cT,\gamma_\cT,\sfs_\semi^\flat)}(\iota^+)}\Bigr).
  \end{split}
  \end{equation}
  Here, $\sfs_\semi^\flat\in\CI(\ol{{}^{0\chop}T^*_\sface}\iota^+)$ is arbitrary but fixed. We make the following choice: first, fix the variable 3b-order function $\sfs^\flat$ according to Lemma~\ref{LemmaWROrder} but with $\sfs^\flat\leq\sfs-\eps$ for some small $\eps>0$; then, denote by $\sfs_\semi^\flat$ the pullback of $\sfs^\flat$ along the second map in~\eqref{EqM3bSclMaps}.

  \pfstep{Transition face normal operator estimate.} The error term in the estimate~\eqref{EqWipHiSymb} is not yet small when $h$ is small since in the norms on $u$ on both sides, the orders at the semiclassical cone transition face (i.e.\ the third order in parentheses) are the same. We improve (i.e.\ weaken) it using an estimate for the $\tface$-normal operator
  \begin{equation}
  \label{EqWipHiNDtf}
  \begin{split}
    N_{\cD,\tface}^\pm(P) &= \mp 2 i\tilde R^{-1}\Bigl(\tilde R\pa_{\tilde R}+\frac{n-1}{2}+S|_{\pa X}\Bigr) + \tilde R^{-2}P_{(0)}(0,\omega,-\tilde R D_{\tilde R},D_\omega) \\
      &\in\Diff_{\bop,\scop}^{2,(2,0)}(\tface;E|_{\pa X}),\qquad \tface=[0,\infty]_{\tilde R}\times\pa X,
  \end{split}
  \end{equation}
  of the $\chop$-operator $h^2\wh{N_\iota^0}(P,-i\gamma_+\pm h^{-1})$ (using the notation of~\eqref{EqWipNormMT} and setting $\tilde R=h R$) by following the arguments in \cite[\S4.4]{HintzConicProp}. Note that the change of coordinates $\tilde R=\hat\rho^{-1}$ identifies $N_{\cD,\tface}^\pm(P)$ and $N_\tface^\pm(P)$ in~\eqref{EqGSONtf}; this is an instance of a general phenomenon recalled after~\eqref{EqMtbTtf} (where we write $N_{\cT,\tface}^\pm(P)$ for $N_\tface^\pm(P)$). Thus, from Lemma~\ref{LemmaStEsttf} and working with a positive b-density $|\frac{\dd\tilde R}{\tilde R}\dd\omega|$ on $\tface$, we have the estimate\footnote{The function spaces are b-Sobolev spaces near $\tilde R=0$ and scattering Sobolev spaces near $\tilde R=\infty$. The aforementioned relationship $\hat\rho=\tilde R^{-1}$ explains the switch in position of `b' and `sc' in the notation for the function spaces here.}
  \begin{equation}
  \label{EqWipHiNtf}
    \| u' \|_{H_{\bop,\scop}^{-N,\gamma_\cT,\sfr}(\tface)} \leq C\|N_{\cD,\tface}^\pm(P)u'\|_{H_{\bop,\scop}^{-N-2,\gamma_\cT-2,\sfr+1}(\tface)}
  \end{equation}
  for the choice $\sfr:=\sfs_\semi^\flat|_{\ol{\Tch^*_{\sface\cap\tface}}\iota^+}-\gamma_\cT$ of scattering decay order; note that in view of~\eqref{EqWipHiSymbThr}, we have $\sfr<\frac{n-1}{2}+\ubar S$ at the image of the outgoing radial set $\cR_{\tface,\rm out}$ in~\eqref{EqSttfRad} under the coordinate change $\tilde R=\hat\rho^{-1}$, or equivalently at the pullback of ${}^\chop\cR_{\cT,\rm out}\subset\Tch^*_\tface\iota^+$ under the bundle identification ${}^{\bop,\scop}T^*\tface\cong\Tch^*_\tface\iota^+$; and similarly we have $\sfr>\frac{n-1}{2}-\ubar S_{\rm in}$ at the incoming radial set. (The scattering decay order is $\sfr-\frac{n}{2}$ if we use the density used in Lemma~\ref{LemmaStEsttf}, and thus the threshold conditions in that Lemma are indeed satisfied.)

  We now rewrite the second term on the right in~\eqref{EqWipHiSymb} as follows. Let $\chi\in\CIc([0,\delta)_h\times[0,\delta)_R)$ be a cutoff function, identically $1$ near $h=R=0$, with $\delta>0$ chosen momentarily. When applying the triangle inequality to the norm of $u=\chi u+(1-\chi)u$, the $\tface$-order on the term $(1-\chi)u$ can be taken to be arbitrary. In the term $\chi u$ on the other hand, given $\eps>0$, we can shrink $\delta>0$ so that (omitting the weight at $\iota^+\cap\scri^+$ from the notation)
  \[
    \|\chi u\|_{H_{\cop,h}^{-N,\gamma_\cT,\gamma_\cT,\sfs_\semi^\flat}(\iota^+\setminus\scri^+)} \leq C h^{-\gamma_\cT} \| \chi u \|_{H_{\bop,\scop}^{-N,\gamma_\cT,\sfr+\eps}(\tface)}
  \]
  for an $h$-independent constant $C$ by \cite[Corollary~3.7]{HintzConicProp}; here, we implicitly pass between the coordinates $R,\omega$ and $\tilde R,\omega$. Using~\eqref{EqWipHiNtf}, and for sufficiently small $\eps>0$ (so that $\sfr+\eps$ satisfies the same threshold conditions as $\sfr$ still), we can estimate this by
  \begin{align*}
     &C h^{-\gamma_\cT} \| N_{\cD,\tface}^\pm(P)(\chi u) \|_{H_{\bop,\scop}^{-N-2,\gamma_\cT-2,\sfr+1+\eps}} \\
     &\qquad \leq C\| h^{-2}N_{\cD,\tface}^\pm(P)(\chi u)\|_{H_{\cop,h}^{-N-2,\gamma_\cT-2,\gamma_\cT-2,\sfs_\semi^\flat+1+2\eps}(\iota^+\setminus\scri^+)},
  \end{align*}
  where on the right we write $N_{\cD,\tface}^\pm(P)$ for the semiclassical cone operator (defined near $\supp\chi\subset(\iota^+\setminus\scri^+)_\chop$) given by the expression~\eqref{EqWipHiNDtf} upon passing to $\tilde R=h R$. We then commute $h^{-2}N_{\cD,\tface}^\pm(P)$ through $\chi$ (the commutator having support disjoint from $\tface$) and then replace $h^{-2}N_{\cD,\tface}^\pm(P)$ by $\wh{N_\iota^0}(P,-i\gamma_+\pm h^{-1})\chi\equiv h^{-2}N_{\cD,\tface}^\pm(P)\chi\bmod\Diffch^{2,2,-1,0}(\iota^+\setminus\scri^+)$. Altogether, we thus conclude that~\eqref{EqWipHiSymb} can be strengthened by replacing the second term on the right by
  \begin{equation}
  \label{EqWipHiAlmost}
    C\|u\|_{H_{0,\cop,h}^{-N,(-N,\gamma_\cT,\gamma_\cT-1,\sfs_\semi^\flat+1+2\eps)}(\iota^+)}
  \end{equation}
  for any fixed $\eps>0$ (with the constant $C$ depending on $\eps$).

  If $\sfs_\semi^\flat<\sfs-1$ (which given the requirements on $\sfs$ and $\sfs^\flat$ is possible if and only if $\sfs$ exceeds the required threshold value at $\cR_{\cT,\rm in}$ by more than $1$), the term~\eqref{EqWipHiAlmost} is small compared to the left hand side of~\eqref{EqWipHiSymb} when $h>0$ is small enough, and we obtain~\eqref{EqWipHi} in this case. Even without this strong requirement on $\sfs$, one can prove the estimate~\eqref{EqWipHi} by using the idea at the end of the proof of \cite[Theorem~4.10]{HintzConicProp}: one takes $\chi$ above to be a cutoff, conormal on $(\iota^+\setminus\scri^+)_\chop$, which near $\sface\cap\tface$ is $1$ for $R\leq h$ and $0$ for $R\geq 2 h$, and one can then show that the error term in~\eqref{EqWipHiSymb} can be weakened to $\|u\|_{H_{0,\cop,h}^{-N,(-N,\gamma_\cT,\gamma_\cT-1+\eta,\sfs_\semi^\flat+1+2\eps-\eta)}(\iota^+)}$ for $\eta\in[0,1]$. For $\eta$ close to $1$, the $\tface$- and $\sface$-orders of this error term are still less than those of the left hand side of~\eqref{EqWipHiSymb}, allowing one to conclude as before.
\end{proof}

\subsection{Sharp mapping properties of admissible wave operators}
\label{SsWM}

For $\ft_0>-1$, we let
\begin{equation}
\label{EqWMOmega}
  \Omega_{\ft_0} := \ft_*^{-1}((\ft_0,\infty]),\qquad \ol{\Omega_{\ft_0}} = \ft_*^{-1}([\ft_0,\infty]) \subset M,\qquad
  \Omega := \Omega_0,\qquad
  \ol\Omega := \ol{\Omega_0}.
\end{equation}
Since $\Omega_{\ft_0}$ has nonempty intersection only with the boundary hypersurfaces $\scri^+$, $\iota^+$, and $\cT^+$ of $M$, we shall work with the edge-3b-phase space $\Tetb^*M$ over $\Omega_{\ft_0}$. We moreover work with the space
\[
  \Hsupp_\etbop^{s,(2\alpha_{\!\scri},\alpha_+,\alpha_\cT)}(\ol\Omega)=x_{\!\scri}^{2\alpha_{\!\scri}}\rho_+^{\alpha_+}\rho_\cT^{\alpha_\cT}\Hsupp_\etbop^s(\ol\Omega)
\]
of distributions lying in a weighted edge-3b-Sobolev space on $\Omega_{-1/2}$ with support in $\ft_*\geq 0$, and with the space $\Hext_\etbop^{s,(2\alpha_{\!\scri},\alpha_+,\alpha_\cT)}(\Omega)$ of restrictions to $\Omega$ of such distributions on $\Omega_{-1/2}$. (Cf.\ \eqref{EqMFHbsupp}--\eqref{EqMFHbext} for the notation used here.)

We recall that $P$ is a $(\ell_0,2\ell_{\!\scri},\ell_+,\ell_\cT)$-admissible wave type operator (with respect to $g$ and a stationary wave type operator $P_0$ relative to $g_0$).

\begin{definition}[Admissible orders]
\label{DefWMOrders}
  We call the orders
  \[
    \sfs\in\CI(\Setb^*_{\ol\Omega}M),\quad
    \alpha_{\!\scri},\quad\alpha_+,\quad\alpha_\cT\in\R,
  \]
  \emph{$P$-admissible} if $\sfs$ satisfies the conditions in Lemma~\usref{LemmaWROrder}---that is, $\sfs$ is constant near the radial sets, satisfies $\sfs+\alpha_+-\alpha_\cT>-\half-\ubar S_{\rm in}$ at $\pa\cR_{\cT,\rm in}^\pm$ and $\sfs+\alpha_+-\alpha_\cT<-\half+\ubar S$ at $\pa\cR_{\cT,\rm out}^\pm$, is monotonically decreasing along the null-bicharacteristic flow, and satisfies the same properties also for the stationary model $P_0$ of $P$---, and if moreover
  \begin{align}
  \label{EqWMOrdersAlphaI}
    \alpha_{\!\scri}&<-\frac12+\ubar p_1, \\
    \alpha_+&<-\frac12+\alpha_{\!\scri},\qquad
    \alpha_+<-\frac{n-1}{2}+\beta^+, \nonumber\\
    \alpha_\cT&\in\Bigl(\alpha_++\frac{n}{2}-\beta^+,\alpha_++\frac{n}{2}-\beta^-\Bigr). \nonumber
  \end{align}
  Here, $\ubar p_1$ is given by~\eqref{EqWRubarp1} in the notation of Definition~\usref{DefGAW}\eqref{ItGAWEdgeN}, and $(\beta^-,\beta^+)$ is the indicial gap from Definition~\usref{DefGSOSpec}.
\end{definition}

When studying the stationary operator $P_0$ instead of $P$, the quantity $\ubar p_1$ in~\eqref{EqWMOrdersAlphaI} gets replaced by $\ubar S$ from Definition~\ref{DefStEstThr}. (In general, we have $\ubar p_1\leq\ubar S$.) The following result on a priori estimates for $P$ is the culmination of the global quantitative microlocal regularity theory developed in~\S\ref{SsWR} and the normal operator estimates of \cite[\S7]{HintzVasyScrieb} and \S\S\ref{SsStCo}, \ref{SsWip}:

\begin{prop}[Fredholm estimates]
\label{PropWMFred}
  Let $\sfs,\sfs_0,\sfs_1\in\CI(\Setb^*_{\ol\Omega}M)$, $\alpha_{\!\scri}$, $\alpha_+$, $\alpha_\cT\in\R$, with $\sfs_0<\sfs-1$ and $\sfs_1>\sfs+1$. Suppose that $\sfs_0$, $\sfs_1$, and $\sfs$ are $P$-admissible orders for the weights $\alpha_{\!\scri}$, $\alpha_+$, $\alpha_\cT$. Fix on $M$ the volume density $|\dd g|$ for the definition of Sobolev spaces and $L^2$-adjoints. Then there exist $\delta>0$ and $C>0$ so that the adjoint estimate
  \begin{equation}
  \label{EqWMFredDir}
  \begin{split}
    &\|u\|_{\dot H_{\eop,\tbop}^{\sfs,(2\alpha_{\!\scri},\alpha_+,\alpha_\cT)}(\ol\Omega;\upbeta^*E)} \\
    &\qquad \leq C\Bigl( \| P u \|_{\dot H_{\eop,\tbop}^{\sfs-1,(2\alpha_{\!\scri}+2,\alpha_++2,\alpha_\cT)}(\ol\Omega;\upbeta^*E)} + \|u\|_{\dot H_{\eop,\tbop}^{\sfs_0,(2(\alpha_{\!\scri}-\delta),\alpha_+-\delta,\alpha_\cT-\delta)}(\ol\Omega;\upbeta^*E)} \Bigr)
  \end{split}
  \end{equation}
  holds for all $u$ for which all stated norms are finite. Similarly, we have the estimate
  \begin{equation}
  \label{EqWMFredAdj}
  \begin{split}
    &\|\tilde u\|_{\bar H_{\eop,\tbop}^{-\sfs+1,(-2\alpha_{\!\scri}-2,-\alpha_+-2,-\alpha_\cT)}(\Omega;\upbeta^*E)} \\
    &\qquad \leq C\Bigl( \|P^*\tilde u\|_{\bar H_{\eop,\tbop}^{-\sfs,(-2\alpha_{\!\scri},-\alpha_+,-\alpha_\cT)}(\Omega;\upbeta^*E)} + \|\tilde u\|_{\bar H_{\eop,\tbop}^{-\sfs_1+1,(-2\alpha_{\!\scri}-2-2\delta,-\alpha_+-2-\delta,-\alpha_\cT-\delta)}(\Omega;\upbeta^*E)}\Bigr)
  \end{split}
  \end{equation}
  for all $\tilde u$ for which all stated norms are finite.
\end{prop}

The point in~\eqref{EqWMFredDir} is that all orders in the weak norm on the second, error, term on the right are smaller than the corresponding orders of the strong norm the left, and thus the inclusion of the strong into the weak Sobolev space is compact; likewise for~\eqref{EqWMFredAdj}. See the discussion after~\eqref{EqMFHbRellich}. We refer the reader to Figure~\ref{FigWMFred} for an illustration of the various microlocal estimates and normal operators involved in the Fredholm control of $P$.

\begin{figure}[!ht]
\centering
\includegraphics{FigWMFred}
\caption{Illustration of the key structures of the null-bicharacteristic flow of $P$ and of the normal operators of $P$ (and their interrelationships), as used in the proof of the Fredholm estimates~\eqref{EqWMFredDir}--\eqref{EqWMFredAdj} for $P$. The radial sets (with a sketch of the flow) are defined in~\eqref{EqWRadScri} and \eqref{EqWRadT}; the normal operator ${}^\eop N_{\scri^+}(P)$ at $\scri^+$ is defined in \cite[\S7.2]{HintzVasyScrieb}; the normal operator $N_\iota^0(P)$ and its (conjugated) Mellin-transforms $\wh{N_\iota^0}(P,\lambda)$, $\wh{N_\iota^1}(P,\lambda)$ are defined in~\eqref{EqWipNorm0} and \eqref{EqWipNormMT}--\eqref{EqWiprhovNopMT}. The transition face models $N_\tface^\pm(P_0)$, resp.\ $N_{\cD,\tface}^\pm(P_0)$ of the spectral family $\wh{P_0}(\sigma)$ of the stationary model $N_{\cT^+}(P)=P_0$, resp.\ of the Mellin-transformed normal operator family at $\iota^+$, are defined in~\eqref{EqGSONtf}, resp.\ \eqref{EqWipHiNDtf}. Recall $v=\frac{t_*}{r}$ and $\hat r=\pm\frac{\sigma}{\rho}=\pm\sigma r$, and see the proofs of Propositions~\ref{PropWipBd} and \ref{PropWipHi} for the relationships between these transition face models. Not explicitly noted here are the relationships between the radial sets in the spacetime edge-3b-cotangent bundle (illustrated here) on the one hand and the radial sets in the (semiclassical) scattering or scattering-b-transition cotangent bundles in the analysis of the spectral family in~\S\ref{SsStEst} and the radial sets in the semiclassical cone and 0-cotangent bundles in the analysis of the Mellin-transformed normal operator family in~\S\ref{SsWip} on the other hand.}
\label{FigWMFred}
\end{figure}

\begin{proof}[Proof of Proposition~\usref{PropWMFred}]
  \pfstep{Direct estimate.} We drop the bundle $\upbeta^*E$ from the notation. Proposition~\ref{PropWR} (with $T_\flat<T_{-1}<T_0<T_1<T_2<0$) then implies the estimate
  \begin{equation}
  \label{EqWMFredDir0}
    \|u\|_{\dot H_{\eop,\tbop}^{\sfs,(2\alpha_{\!\scri},\alpha_+,\alpha_\cT)}(\ol\Omega)} \leq C\Bigl( \| P u \|_{\dot H_{\eop,\tbop}^{\sfs-1,(2\alpha_{\!\scri}+2,\alpha_++2,\alpha_\cT)}(\ol\Omega)} + \|u\|_{\dot H_{\eop,\tbop}^{\sfs_0,(2\alpha_{\!\scri},\alpha_+,\alpha_\cT)}(\ol\Omega)}\Bigr).
  \end{equation}
  (For later use, we point out that this estimate, which only uses microlocal regularity estimates, holds in the strong sense that if the right hand side is finite, then so is the left, and the estimate holds.) We further improve (i.e.\ weaken the norm on) the second norm on the right successively at the boundary hypersurfaces of $\ol\Omega$ at infinity.

  \pfsubstep{(i)}{Improvement at $\scri^+$.} We use \cite[Theorem~7.3]{HintzVasyScrieb} for domains $\Omega_{\delta,\rho}^+$ (in the notation of the reference) which contain $\scri^+\cap\{\ft_*\geq 0\}$ but are contained in a neighborhood of $\scri^+\cap\{\ft_*>-1\}$. With $\chi_\scri\in\CI(M)$ supported in a small neighborhood of $\scri^+$, we then estimate the $\dot H_{\eop,\tbop}^{\sfs_0,(2\alpha_{\!\scri},\alpha_+,\alpha_\cT)}(\ol\Omega)$ norm of $\chi_\scri u$ using~\cite[Theorem~7.3(1)]{HintzVasyScrieb}; for any fixed $\eps>0$, we obtain
  \begin{align}
    \|u\|_{\dot H_{\eop,\tbop}^{\sfs_0,(2\alpha_{\!\scri},\alpha_+,\alpha_\cT)}(\ol\Omega)} &\leq C'\|u\|_{H_{\eop,\bop}^{\sfs_0,(2\alpha_{\!\scri},\alpha_+)}(\Omega_{\delta,\rho}^+)^{\bullet,-}} + \|(1-\chi_\scri)u\|_{\dot H_{\eop,\tbop}^{\sfs_0,(2\alpha_{\!\scri},\alpha_+,\alpha_\cT)}(\ol\Omega)} \nonumber\\
      &\leq C'\Bigl(\|P u\|_{H_{\eop,\bop}^{\sfs_0-1,(2\alpha_{\!\scri}+2,\alpha_++2)}(\Omega_{\delta,\rho}^+)^{\bullet,-}} + C_\eps\|u\|_{H_{\eop,\bop}^{\sfs_0+1,(2\alpha_{\!\scri}-2\ell_{\!\scri},\alpha_+)}(\Omega_{\delta,\rho}^+)^{\bullet,-}} \nonumber\\
      &\quad\qquad + \eps\|u\|_{H_{\eop,\bop}^{\sfs_0+1,(2\alpha_{\!\scri},\alpha_+)}(\Omega_{\delta,\rho}^+)^{\bullet,-}}\Bigr) + \|(1-\chi_\scri)u\|_{\dot H_{\eop,\tbop}^{\sfs_0,(2\alpha_{\!\scri},\alpha_+,\alpha_\cT)}(\ol\Omega)} \nonumber\\
  \label{EqWMFredDir1}
  \begin{split}
      &\leq C'\|P u\|_{\dot H_{\eop,\tbop}^{\sfs_0-1,(2\alpha_{\!\scri}+2,\alpha_++2,\alpha_\cT)}(\ol\Omega)} + C'_\eps\|u\|_{\dot H_{\eop,\tbop}^{\sfs_0+1,(2\alpha_{\!\scri}-2\ell_{\!\scri},\alpha_+,\alpha_\cT)}(\ol\Omega)} \\
      &\qquad + C'\eps\|u\|_{\dot H_{\eop,\tbop}^{\sfs_0+1,(2\alpha_{\!\scri},\alpha_+,\alpha_\cT)}(\ol\Omega)},
  \end{split}
  \end{align}
  where we use that $\supp(1-\chi_\scri)\cap\scri^+=\emptyset$ to change the $\scri^+$-decay order in the norm of $(1-\chi_\scri)u$ to $2\alpha_{\!\scri}-2\ell_{\!\scri}$. We fix $\eps>0$ sufficiently small so that
  \[
    C C'\eps\|u\|_{\dot H_{\eop,\tbop}^{\sfs_0+1,(2\alpha_{\!\scri},\alpha_+,\alpha_\cT)}(\ol\Omega)} \leq \half\|u\|_{\dot H_{\eop,\tbop}^{\sfs,(2\alpha_{\!\scri},\alpha_+,\alpha_\cT)}(\ol\Omega)};
  \]
  this uses that $\sfs_0+1\leq\sfs$. We can then absorb this term into the left hand side of~\eqref{EqWMFredDir0}.

  Fix now $\delta\in(0,\ell_{\!\scri}]$ so that\footnote{The assumptions on $\sfs_0$ and $\sfs$ imply that $\sfs_0+1$, $\alpha_{\!\scri}$, $\alpha_+$, $\alpha_\cT$ are $P$-admissible, which implies the existence of $\delta$.} $\sfs_0+1$, $\alpha_{\!\scri}-\delta$, $\alpha_+$, $\alpha_\cT$ are $P$-admissible orders. From~\eqref{EqWMFredDir1} and \eqref{EqWMFredDir0}, we then obtain
  \begin{equation}
  \label{EqWMFredDir2}
    \|u\|_{\dot H_{\eop,\tbop}^{\sfs,(2\alpha_{\!\scri},\alpha_+,\alpha_\cT)}(\ol\Omega)} \leq C\Bigl( \| P u \|_{\dot H_{\eop,\tbop}^{\sfs-1,(2\alpha_{\!\scri}+2,\alpha_++2,\alpha_\cT)}(\ol\Omega)} + \|u\|_{\dot H_{\eop,\tbop}^{\sfs_0,(2(\alpha_{\!\scri}-\delta),\alpha_+,\alpha_\cT)}(\ol\Omega)}\Bigr)
  \end{equation}
  except the differential order in the error term is $\sfs_0+1$ at first; it can however be reduced to $\sfs_0$ via microlocal edge-3b-propagation estimates, i.e.\ concretely using~\eqref{EqWMFredDir0} with $\sfs_0+1<\sfs$ and $\alpha_{\!\scri}-\delta$ in place of $\sfs$ and $\alpha_{\!\scri}$.

  \pfsubstep{(ii)}{Improvement at $\iota^+$.} Let now $\chi_\iota\in\CIc(\Omega)$ be identically $1$ in a collar neighborhood of $\iota^+$. Then $(1-\chi_\iota)u$ has support disjoint from $\iota^+$, and hence the $\iota^+$-decay order in its edge-3b-Sobolev norms can be taken to be arbitrary. On the other hand, we can estimate $\chi_\iota u$ in terms of the $\iota^+$-normal operator of $P$ via the Mellin transform in the total boundary defining function $T=t_*^{-1}$ of $\iota^+\cup\cT^+$. To facilitate this, we first pass to a positive b-density $\mu_\bop$ on $M$ and recall from~\eqref{EqWipBdbDensity}--\eqref{EqWipBdRel} that
  \begin{equation}
  \label{EqWMFredDensitySwitch}
    \|\chi_\iota u\|_{H_{\eop,\tbop}^{\sfs_0,(2(\alpha_{\!\scri}-\delta),\alpha_+,\alpha_\cT)}(M)} \sim \|\chi_\iota u\|_{H_{\eop,\tbop}^{\sfs_0,(2(\gamma_{\!\scri}-\delta),\gamma_+,\gamma_\cT+\gamma_+)}(M,\mu_\bop)}
  \end{equation}
  where $\gamma_{\!\scri}:=\alpha_{\!\scri}+\frac{n}{2}$, $\gamma_+:=\alpha_++\frac{n+1}{2}$, and $\gamma_\cT:=\alpha_\cT+\frac12-\gamma_+=\alpha_\cT-\alpha_+-\frac{n}{2}$. Then Proposition~\ref{PropMFebEq0} (near $\scri^+$) and \cite[Proposition~4.29(2)]{Hintz3b} (near $\cT^+$, with $\mu_\cD=0$, $\mu_\cT=1$, $\hat\mu=-1$, and with $\gamma_+,\gamma_\cT+\gamma_+$ in place of $\alpha_\cD,\alpha_\cT$---corresponding to working with positive b-densities on $M$ and $\iota^+$, denoted $\cD$ in \cite{Hintz3b}) imply for any fixed $\eps>0$, and provided $\chi_\iota$ is supported in a sufficiently small neighborhood (depending on $\eps$) of $\iota^+$, the estimate
  \begin{align*}
    &\|\chi_\iota u\|_{H_{\eop,\tbop}^{\sfs_0,(2(\gamma_{\!\scri}-\delta),\gamma_+,\gamma_\cT+\gamma_+)}(M,\mu_\bop)}^2 \\
    &\qquad \leq C\Biggl( \int_{[-1,1]} \| \wh{\chi_\iota u}(-i\gamma_++\lambda) \|_{H_{0,\bop}^{(\sfs_0)_\infty+\eps,(2(\gamma_{\!\scri}-\delta),\gamma_\cT)}(\iota^+)}^2\,\dd\lambda \\
    &\qquad\hspace{4em} + \sum_\pm \int_{[1,\infty)} \|\wh{\chi_\iota u}(-i\gamma_+\pm\lambda)\|_{H_{0,\cop,|\lambda|^{-1}}^{(\sfs_0)_\infty+\eps,(2(\gamma_{\!\scri}-\delta),\gamma_\cT,\gamma_\cT,(\sfs_0)_\semi+\eps)}(\iota^+)}^2\,\dd\lambda\Biggr).
  \end{align*}
  Here, $\wh{\chi_\iota u}(\lambda)$ (a function on $\iota^+$ depending on $\lambda\in\C$) is the Mellin transform of $\chi_\iota u$ in $T$; and $(\sfs_0)_\infty$ and $(\sfs_0)_\semi$ are the differential and semiclassical orders induced by $\sfs_0$, defined near $\scri^+$ as in Proposition~\ref{PropMFebEq0} and near $\cT^+$ as in \cite[Proposition~4.29(2)]{Hintz3b}. The norms on $\iota^+$ are defined with respect to a positive b-density on $\iota^+$. In view of the $P$-admissibility of $\sfs_0+\eps,\alpha_{\!\scri}-\delta,\alpha_+,\alpha_\cT$ for small $\delta$ and $\eps$, the orders on the function spaces on the right are such that for sufficiently small $\eps>0$, the estimates for the $\iota^+$-normal operator $\wh{N_\iota^0}(P)$ of $T^{-2}P$ in Proposition~\ref{PropWipBd} (for the first integral) and Proposition~\ref{PropWipHi} (for the second integral) apply. Using Proposition~\ref{PropMFebEq0} and \cite[Proposition~4.29(2)]{Hintz3b} again to pass back to spacetime Sobolev spaces, we obtain
  \begin{equation}
  \label{EqWMFredDir3}
    \|\chi_\iota u\|_{H_{\eop,\tbop}^{\sfs_0,(2(\gamma_{\!\scri}-\delta),\gamma_+,\gamma_\cT+\gamma_+)}(M,\mu_\bop)} \leq C\|N_\iota^0(P)(\chi_\iota u)\|_{H_{\eop,\tbop}^{\sfs_0-1+2\eps,(2(\gamma_{\!\scri}-\delta)+2,\gamma_+,\gamma_\cT+\gamma_+-2)}(M,\mu_\bop)},
  \end{equation}
  where in addition to the relaxation of the orders by $\eps$ (to accommodate for the possible failure of $\sfs_0$ to be dilation-invariant near $\iota^+$) we also relax the differential order in~\eqref{EqWipHi} (and likewise for bounded frequencies) from $(\sfs_0)_\infty-2+2\eps$ to $(\sfs_0)_\infty-1+2\eps$, so that the differential order $(\sfs_0)_\infty-1+2\eps$ and the semiclassical order $(\sfs_0)_\semi-1+2\eps$ are the orders induced by the single spacetime edge-3b-order $\sfs_0-1+2\eps$. Moreover, we abuse notation and write $N_\iota^0(P)$ for any element of $\Diff_{\eop,\tbop}^{2,(-2,0,2)}(M;\upbeta^*E)$ with $\iota^+$-normal operator given by $N_\iota^0(P)$ (see~\eqref{EqWipNorm0}).

  Lemma~\ref{LemmaWip} then enables us to pass from $N_\iota^0(P)=T^{-2}N_\iota(P)$ to $P$; moreover, the commutator of $P$ with $\chi_\iota$ is an operator with coefficients having support disjoint from $\iota^+$. Thus, from~\eqref{EqWMFredDir3} we obtain
  \begin{align*}
    &\|\chi_\iota u\|_{H_{\eop,\tbop}^{\sfs_0,(2(\gamma_{\!\scri}-\delta),\gamma_+,\gamma_\cT+\gamma_+)}(M,\mu_\bop)} \\
    &\qquad \leq C\Bigl( \|\chi_\iota P u\|_{H_{\eop,\tbop}^{\sfs_0-1+2\eps,(2(\gamma_{\!\scri}-\delta)+2,\gamma_++2,\gamma_\cT+\gamma_+)}(M,\mu_\bop)} + \|u\|_{\dot H_{\eop,\tbop}^{\sfs_0+1+2\eps,(2(\gamma_{\!\scri}-\delta),\gamma_+-\ell_+,\gamma_\cT+\gamma_+)}(\ol\Omega)}\Bigr).
  \end{align*}
  Passing back to the metric density as in~\eqref{EqWMFredDensitySwitch}, we use this estimate on the error term~\eqref{EqWMFredDir2} (split via $u=\chi_\iota u+(1-\chi_\iota)u$ and the triangle inequality) and obtain the improvement
  \begin{equation}
  \label{EqWMFredDir4}
    \|u\|_{\dot H_{\eop,\tbop}^{\sfs,(2\alpha_{\!\scri},\alpha_+,\alpha_\cT)}(\ol\Omega)} \leq C\Bigl( \| P u \|_{\dot H_{\eop,\tbop}^{\sfs-1,(2\alpha_{\!\scri}+2,\alpha_++2,\alpha_\cT)}(\ol\Omega)} + \|u\|_{\dot H_{\eop,\tbop}^{\sfs_0,(2(\alpha_{\!\scri}-\delta),\alpha_+-\delta,\alpha_\cT)}(\ol\Omega)}\Bigr),
  \end{equation}
  where we shrink $\delta>0$ further if necessary so as to satisfy $\delta\leq\ell_+$, and so that $\sfs_0+1+2\eps$ and $\sfs_0$ are $P$-admissible with weights $\alpha_{\!\scri}-\delta$, $\alpha_+-\delta$, $\alpha_\cT$. As above, the differential order of the last term on the right is initially the stronger $\sfs_0+1+2\eps$, which however can be reduced back to $\sfs_0$ as before using~\eqref{EqWMFredDir0} for $\sfs_0+1+2\eps\leq\sfs$ (where we fix $\eps>0$ sufficiently small), $\alpha_{\!\scri}-\delta$, $\alpha_+-\delta$ in place of $\sfs$, $\alpha_{\!\scri}$, $\alpha_+$.

  \pfsubstep{(iii)}{Improvement at $\cT^+$.} In the final step, we weaken the $\cT^+$-decay order of the second, error, term in~\eqref{EqWMFredDir4}. Thus, let $\chi_\cT\in\CIc(\Omega)$ be identically $1$ in a collar neighborhood of $\cT^+$. Write moreover
  \[
    u' := T^{-\alpha_\cT}u.
  \]
  We only record orders at $\iota^+$ and $\cT^+$ for now. Using \cite[Proposition~4.29(1)]{Hintz3b}, we obtain, for any $\eps>0$ and provided the support of $\chi_\cT$ is contained in a sufficiently small (depending on $\eps>0$) neighborhood of $\cT^+$,
  \begin{equation}
  \label{EqWMFredDir5}
  \begin{split}
    \|\chi_\cT u\|_{\Htb^{\sfs_0,(\alpha_+-\delta,\alpha_\cT)}}^2 &= \|\chi_\cT u'\|_{\Htb^{\sfs_0,(\alpha_+-\alpha_\cT-\delta,0)}}^2 \\
      &\leq C \sum_\pm  \int_{\pm[0,1]} \bigl\|\wh{\chi_\cT u'}(\sigma)\bigr\|_{H_{\scbtop,\sigma}^{(\sfs_0)_\infty+\eps,(\sfs_0)_\scop+\eps+\alpha_+-\alpha_\cT-\delta,\alpha_+-\alpha_\cT-\delta,0}(X)}^2\,\dd\sigma \\
      &\qquad + C\sum_\pm \int_{\pm[1,\infty)} \bigl\|\wh{\chi_\cT u'}(\sigma)\bigr\|_{H_{\scop,|\sigma|^{-1}}^{(\sfs_0)_\infty+\eps,(\sfs_0)_\scop+\eps+\alpha_+-\alpha_\cT-\delta,(\sfs_0)_\semi+\eps}(X)}^2\,\dd\sigma.
  \end{split}
  \end{equation}
  Here, we use the metric volume density $|\dd g|$ on $M$, and the density $|\dd g_X|$ on $X$ (see Lemma~\ref{LemmaGSGVol}). Moreover, $(\sfs_0)_\infty$ and $(\sfs_0)_\scop$ in the second line are the pullbacks of $\sfs_0$ under the maps~\eqref{Eq3bLoMaps}, and $(\sfs_0)_\infty$, $(\sfs_0)_\scop$, and $(\sfs_0)_\semi$ in the third line are the pullbacks of $\sfs_0$ under the maps~\eqref{Eq3bHiMaps}.

  Note then that for sufficiently small $\delta,\eps>0$, the scattering decay order of the spaces in which $\wh{\chi_\cT u'}(\sigma)$ is estimated, i.e.\ $(\sfs_0)_\scop+\alpha_+-\alpha_\cT+\eps-\delta$, is $>-\half-\ubar S_{\rm in}$ at the incoming, and $<-\half+\ubar S$ at the outgoing radial sets in view of the admissibility of $\sfs_0,\alpha_{\!\scri},\alpha_+,\alpha_\cT$ and the relationships~\eqref{Eq3bLoMaps}--\eqref{Eq3bHiMaps} between the 3b-phase space over $\cT^+\cap\iota^+$ on the one hand and the $\scbtop$- and semiclassical scattering phase spaces over $\pa\cT^+\subset\cT^+$ on the other hand. Thus, the low energy estimates from Proposition~\ref{PropStEstLo}\eqref{ItStEstLosc} are applicable, as are the high energy estimates from Proposition~\ref{PropStEstHi}\eqref{ItStEstHisc}; thus, with $\wh{P_0}(\sigma)$ denoting the spectral family of the time-translation invariant model $P_0$ of $P$, we can further estimate~\eqref{EqWMFredDir5} by
  \begin{align*}
    &\|\chi_\cT u\|_{\Htb^{\sfs_0,(\alpha_+-\delta,\alpha_\cT)}}^2 \\
    &\qquad \leq C\sum_\pm\int_{\pm[0,1]} \| \wh{P_0}(\sigma)\wh{\chi_\cT u'}(\sigma) \|_{H_{\scbtop,\sigma}^{(\sfs_0)_\infty-2+\eps,(\sfs_0)_\scop+1+\alpha_+-\alpha_\cT-\delta+\eps,\alpha_+-\alpha_\cT+2-\delta,0}(X)}^2\,\dd\sigma \\
    &\qquad\qquad + C\sum_\pm\int_{\pm[1,\infty)} \| \wh{P_0}(\sigma)\wh{\chi_\cT u'}(\sigma) \|_{H_{\scop,|\sigma|^{-1}}^{(\sfs_0)_\infty-2+\eps,(\sfs_0)_\scop+1+\alpha_+-\alpha_\cT-\delta+\eps,(\sfs_0)_\semi-1+\eps}(X)}^2\,\dd\sigma \\
    &\qquad \leq C \| P_0(\chi_\cT u') \|_{\Htb^{\sfs_0-1+2\eps,(\alpha_+-\alpha_\cT+2-\delta,0)}}^2 \\
    &\qquad = C \| T^{\alpha_\cT}P_0(\chi_\cT T^{-\alpha_\cT}u) \|_{\Htb^{\sfs_0-1+2\eps,(\alpha_++2-\delta,\alpha_\cT)}}^2;
  \end{align*}
  here, in the passage to the penultimate line we applied \cite[Proposition~4.29(1)]{Hintz3b} again (and increased the differential order). Now,
  \begin{align*}
    T^{\alpha_\cT}P_0 \chi_\cT T^{-\alpha_\cT} &= \chi_\cT P + [P,\chi_\cT] - (P-P_0)\chi_\cT + T^{\alpha_\cT}[P_0,T^{-\alpha_\cT}]\chi_\cT \\
      &\equiv \chi_\cT P \bmod \cA^{(2,\min(\ell_\cT,1))}\Difftb^2(M;\upbeta^*E);
  \end{align*}
  indeed, $[P,\chi_\cT]$ vanishes near $\cT^+$, while $P-P_0$ is controlled by Definition~\ref{DefGAW}\eqref{ItGAWStruct}, and lastly $T^{\alpha_\cT}[P_0,T^{-\alpha_\cT}]\chi_\cT\in\Difftb^{2,(-2,-1)}(M;\upbeta^*E)$ by~\eqref{Eq3bConjT}. We therefore obtain (for $\delta<\min(\ell_\cT,1)$)
  \begin{equation}
  \label{EqWMFredDir6}
    \|\chi_\cT u\|_{\Htb^{\sfs_0,(\alpha_+-\delta,\alpha_\cT)}} \leq C\Bigl(\|\chi_\cT P u\|_{\Htb^{\sfs_0-1+2\eps,(\alpha_++2-\delta,\alpha_\cT)}} + \|u\|_{\bar H_{\eop,\tbop}^{\sfs_0+1+2\eps,(-N,\alpha_+-\delta,\alpha_\cT-\delta)}(\ol\Omega)}\Bigr).
  \end{equation}

  We now plug the estimate~\eqref{EqWMFredDir6} into the error term in~\eqref{EqWMFredDir4} which we split using $u=\chi_\cT u+(1-\chi_\cT)u$ and the triangle inequality. Since $\supp(1-\chi_\cT)$ is disjoint from $\cT^+$, we finally obtain (upon fixing $\eps>0$ to be small, and relaxing the differential order of the error term using propagation estimates as before)
  \[
    \|u\|_{\dot H_{\eop,\tbop}^{\sfs,(2\alpha_{\!\scri},\alpha_+,\alpha_\cT)}(\ol\Omega)} \leq C\Bigl( \| P u \|_{\dot H_{\eop,\tbop}^{\sfs-1,(2\alpha_{\!\scri}+2,\alpha_++2,\alpha_\cT)}(\ol\Omega)} + \|u\|_{\dot H_{\eop,\tbop}^{\sfs_0,(2(\alpha_{\!\scri}-\delta),\alpha_+-\delta,\alpha_\cT-\delta)}(\ol\Omega)}\Bigr).
  \]
  This is the desired semi-Fredholm estimate~\eqref{EqWMFredDir}.

  \pfstep{Adjoint estimate.} The starting point is the edge-3b-regularity estimate for $P^*$ given in Proposition~\ref{PropWR}. Concretely, taking $T_{-1}=0$, $T_0=1$, $T_1=2$ in Proposition~\ref{PropWR}, we apply the estimate~\eqref{EqWRAdj}, which, roughly speaking, gives control on a strong norm (with differential order $-\sfs+1$) on $\tilde u$ in $\ft_*\geq 2$ in terms of a strong norm on $P^*\tilde u$ and a weak norm on $\tilde u$ in $\ft_*\geq 1$. Write $\tilde f:=P^*\tilde u\in\bar H_{\eop,\tbop}^{-\sfs,(-2\alpha_{\!\scri},-\alpha_+,-\alpha_\cT)}(M;\upbeta^*E)$, and let $\chi_2\in\CIc([0,4))$ be identically $1$ on $[0,3]$; then $\tilde u_2:=\chi_2\tilde u$ is supported in $\ft_*<4$ and satisfies
  \begin{equation}
  \label{EqWMFredAdjWave}
    P^*\tilde u_2 = \tilde f_2 := \chi_2\tilde f+[P^*,\chi_2]\tilde u.
  \end{equation}
  We can then use \cite[Theorem~6.4(2)]{HintzVasyScrieb} (keeping in mind \cite[Remark~4.7]{HintzVasyScrieb} regarding the numerology for the orders) to estimate $\tilde u_2$, as the unique backwards solution of the equation~\eqref{EqWMFredAdjWave} with weight $-2\alpha_{\!\scri}-2$ at $\scri^+$, in terms of $\tilde f_2$; concretely, keeping only the decay order at $\scri^+$ in the notation, we have
  \[
    \|\tilde u_2\|_{H_\eop^{-\sfs+1,-2\alpha_{\!\scri}-2}(\Omega\setminus\Omega_4)^{-,\bullet}} \leq C\|\tilde f_2\|_{H_\eop^{-\sfs,-2\alpha_{\!\scri}}(\Omega\setminus\Omega_4)^{-,\bullet}},
  \]
  where `$-$' indicates the extendible character at $\ft_*=0$, and `$\bullet$' indicates the supported character at $\ft_*=4$. But since $[P^*,\chi_2]$ is supported in $\ft_*^{-1}([3,4])$ and thus in particular in $\chi_1^{-1}(1)$ (in the notation of the estimate~\eqref{EqWRAdj}), we can further estimate
  \[
    \|\tilde f_2\|_{H_\eop^{-\sfs,-2\alpha_{\!\scri}}(\Omega\setminus\Omega_4)^{-,\bullet}} \leq C\Bigl(\|\tilde f\|_{\Hext_{\eop,\tbop}^{-\sfs,(-2\alpha_{\!\scri},-\alpha_+,-\alpha_\cT)}} + \|\chi_1\tilde u\|_{\Hext_{\eop,\tbop}^{-\sfs+1,(-2\alpha_{\!\scri}-2,-\alpha_+-2,-\alpha_\cT)}}\Bigr).
  \]
  Altogether, we have thus established the estimate
  \begin{equation}
  \label{EqWMFredAdj0}
  \begin{split}
    &\|u\|_{\Hext_{\eop,\tbop}^{-\sfs+1,(-2\alpha_{\!\scri}-2,-\alpha_+-2,-\alpha_\cT)}} \\
    &\qquad \leq C\Bigl( \|P^*\tilde u\|_{\Hext_{\eop,\tbop}^{-\sfs,(-2\alpha_{\!\scri},-\alpha_+,-\alpha_\cT)}} + \|\chi_0\tilde u\|_{H_{\eop,\tbop}^{-\sfs_1+1,(-2\alpha_{\!\scri}-2,-\alpha_+-2,-\alpha_\cT)}}\Bigr),
  \end{split}
  \end{equation}
  where $\chi_0\in\CI(\R_{\ft_*})$ can be taken to be supported in $[1,\infty)$ and identically $1$ on $[2,\infty)$.

  The proof then proceeds similarly to before: we estimate $\tilde u$ near $\scri^+\cap\{\ft_*\geq 1\}$ by means of the normal operator estimate given in \cite[Theorem~7.3(2)]{HintzVasyScrieb}, which allows us to weaken the $\scri^+$-decay order on the error term in~\eqref{EqWMFredAdj0} to $-2\alpha_{\!\scri}-2-2\delta$ for small $\delta>0$. One then estimates $\chi_0\tilde u$ near $\iota^+$ in terms of $P^*\tilde u$ using the (inverse) Mellin transform and the adjoint versions of the estimates of Propositions~\ref{PropWipBd} and \ref{PropWipHi}; this allows one to weaken the $\iota^+$-decay order on the error term to $-\alpha_+-2-\delta$. Finally, one uses the adjoint versions of the estimates for the spectral family given in Proposition~\ref{PropStEstHi} and Proposition~\ref{PropStEstLo} to estimate $\chi_0\tilde u$ near $\cT^+$ in terms of $P^*\tilde u$. This gives~\eqref{EqWMFredAdj} and finishes the proof.
\end{proof}

\begin{cor}[Injectivity estimate]
\label{CorWMInj}
  Under the assumptions of Proposition~\usref{PropWMFred} (with the requirement on the order $\sfs$ being that $\sfs\pm 1$ be $P$-admissible), there exists a constant $C$ so that
  \begin{equation}
  \label{EqWMInj}
    \|u\|_{\dot H_{\eop,\tbop}^{\sfs,(2\alpha_{\!\scri},\alpha_+,\alpha_\cT)}(\ol\Omega;\upbeta^*E)} \leq C \| P u \|_{\dot H_{\eop,\tbop}^{\sfs-1,(2\alpha_{\!\scri}+2,\alpha_++2,\alpha_\cT)}(\ol\Omega;\upbeta^*E)}
  \end{equation}
  for all $u$ for which both norms are finite.
\end{cor}
\begin{proof}
  Since $P$ is a wave operator, it is automatically injective on any space of distributions which are supported in $\ft_*\geq 0$. By a standard functional analytic argument, we can thus drop the error term in~\eqref{EqWMFredDir} upon increasing the constant $C$. Indeed, if this were not possible, there would exist a sequence $u_j\in\dot H_{\eop,\tbop}^{\sfs,(2\alpha_{\!\scri},\alpha_+,\alpha_\cT)}(\ol\Omega;\upbeta^*E)$ with norm $1$ so that $P u_j\to 0$ in $\dot H_{\eop,\tbop}^{\sfs-1,(2\alpha_{\!\scri}+2,\alpha_++2,\alpha_\cT)}$; upon passing to a subsequence, we may assume that $u_j\weakto u_0$. The estimate~\eqref{EqWMFredDir} gives a uniform positive lower bound on the norm of $u_j$ in $\dot H_{\eop,\tbop}^{\sfs-\delta,(2(\alpha_{\!\scri}-\delta),\alpha_+-\delta,\alpha_\cT-\delta)}$, and thus on the norm of $u_0$ (which is the strong limit of $u_j$ in this weaker space); thus $u_0\neq 0$, but $P u_0=0$, which gives the desired contradiction.
\end{proof}

Note that the estimates~\eqref{EqWMFredDir} and~\eqref{EqWMInj} apply also to $P_0$, and the estimate~\eqref{EqWMFredAdj} applies to $P_0^*$.

\begin{cor}[Invertibility of $P_0$]
\label{CorWMP0Inv}
  Under the assumptions of Proposition~\usref{PropWMFred} for $P_0$ (with the requirement on the order $\sfs$ being that $\sfs\pm 1$ be $P_0$-admissible), and defining the domain $\Omega\subset M$ in terms of $t_*$, i.e.\ $\Omega=t_*^{-1}((0,\infty])$, the operators
  \begin{align}
  \label{EqWMP0InvDir}
  \begin{split}
    P_0 &\colon \bigl\{ u\in\dot H_{\eop,\tbop}^{\sfs,(2\alpha_{\!\scri},\alpha_+,\alpha_\cT)}(\ol\Omega;\upbeta^*E) \colon P_0 u\in\dot H_{\eop,\tbop}^{\sfs-1,(2\alpha_{\!\scri}+2,\alpha_++2,\alpha_\cT)}(\ol\Omega;\upbeta^*E) \bigr\} \\
      &\hspace{18em} \to \dot H_{\eop,\tbop}^{\sfs-1,(2\alpha_{\!\scri}+2,\alpha_++2,\alpha_\cT)}(\ol\Omega;\upbeta^*E),
  \end{split} \\
  \label{EqWMP0InvAdj}
  \begin{split}
    P_0^* &\colon \bigl\{ \tilde u\in\bar H_{\eop,\tbop}^{-\sfs+1,(-2\alpha_{\!\scri}-2,-\alpha_+-2,-\alpha_\cT)}(\Omega;\upbeta^*E) \colon P_0^*\tilde u\in\bar H_{\eop,\tbop}^{-\sfs,(-2\alpha_{\!\scri},-\alpha_+,-\alpha_\cT)}(\Omega;\upbeta^*E) \bigr\} \\
      &\hspace{18em} \to \bar H_{\eop,\tbop}^{-\sfs,(-2\alpha_{\!\scri},-\alpha_+,-\alpha_\cT)}(\Omega;\upbeta^*E)
  \end{split}
  \end{align}
  are invertible.
\end{cor}
\begin{proof}
  Corollary~\ref{CorWMInj}, applied to $P_0$, implies that~\eqref{EqWMP0InvDir} is injective and has closed range. We claim that the range contains $\CIdot(\ol\Omega;\upbeta^*E)$; since this is a dense subspace, this implies the surjectivity of~\eqref{EqWMP0InvDir}. But by Theorem~\ref{ThmStCo}, the forward solution $u$ of $P_0 u=f\in\CIdot(\ol\Omega;\upbeta^*E)$ satisfies
  \[
    u \in \dot\cA^{(\frac{n-1}{2}+\ubar S-,\beta^++1-,\beta^+-\beta^-+1-)}(\ol\Omega;\upbeta^*E),
  \]
  where the dot indicates the supported character of $u$ at $t_*=0$. If $\mu_\bop$ is a positive b-density on $M$, this implies
  \begin{align*}
    u &\in \dot H_\bop^{\infty,(\frac{n-1}{2}+\ubar S-,\beta^++1-,\beta^+-\beta^-+1-)}(\ol\Omega,\mu_\bop;\upbeta^*E) \\
      &\subset \dot H_\bop^{\infty,(-\frac12+\ubar S-,-\frac{n-1}{2}+\beta^+-,\beta^+-\beta^-+\frac12-)}(\ol\Omega,|\dd g_0|;\upbeta^*E) \\
      &\subset \dot H_\bop^{\infty,(\alpha_{\!\scri},\alpha_+,\alpha_\cT)}(\ol\Omega,|\dd g_0|;\upbeta^*E),
  \end{align*}
  where in the last inclusion we use~\eqref{EqWMOrdersAlphaI}, which in particular gives $\alpha_\cT<\alpha_++\frac{n}{2}-\beta^-<(-\frac{n-1}{2}+\beta^+)+\frac{n}{2}-\beta^-<\half+\beta^+-\beta^-$. The last space is contained in the domain of $P_0$ in~\eqref{EqWMP0InvDir}, as required.

  The surjectivity of $P_0$ implies the injectivity of $P_0^*$ in~\eqref{EqWMP0InvAdj}. The injectivity of $P_0$ implies that~\eqref{EqWMP0InvAdj} has dense range; but the semi-Fredholm estimate~\eqref{EqWMFredAdj} implies that the range of $P_0^*$ is also closed, and the proof is complete.
\end{proof}

\begin{rmk}[Norm of the inverse]
\label{RmkWMStatConst}
  Since the invertibility of $P_0$ in Corollary~\ref{CorWMP0Inv} was proved using a Fredholm/compactness argument (in the proofs of Corollaries~\ref{CorWMInj} and \ref{CorWMP0Inv}, and already earlier e.g.\ in the estimates for the spectral family at $\cT^+$ and for the Mellin-transformed normal operator family at $\iota^+$), we do not obtain quantitative control on the operator norms of $P_0^{-1}$ and $(P_0^*)^{-1}$ here. However, if $P_0$ varies continuously inside a compact family of spectrally admissible stationary wave type operators, then these operator norms (between the stated spaces, with differential orders and weights which themselves vary over compact sets of admissible orders) are uniformly bounded; this follows from a functional analytic argument much as in the proof of Corollary~\ref{CorWMInj}.
\end{rmk}

We revert to the notation~\eqref{EqWMOmega}. We are now in a position to prove the main result of this paper for non-stationary operators:

\begin{thm}[Forward solutions for $P$]
\label{ThmWM}
  Let $\sfs\in\CI(\Setb^*_{\ol\Omega}M)$, $\alpha_{\!\scri}$, $\alpha_+$, $\alpha_\cT\in\R$. Suppose that both $\sfs-1$ and $\sfs+2$ are $P$-admissible orders for the weights $\alpha_{\!\scri}$, $\alpha_+$, $\alpha_\cT$ (see Definition~\usref{DefWMOrders}). Fix on $M$ the volume density $|\dd g|$ for the definition of Sobolev spaces and $L^2$-adjoints. Then there exists $C>0$ so that the following holds: for any $f\in\dot H_{\eop,\tbop}^{\sfs-1,(2\alpha_{\!\scri}+2,\alpha_++2,\alpha_\cT)}(\ol\Omega;\upbeta^*E)$, the unique forward solution $u$ of $P u=f$ (i.e.\ $u=0$ for $\ft_*<0$) satisfies $u\in\dot H_{\eop,\tbop}^{\sfs,(2\alpha_{\!\scri},\alpha_+,\alpha_\cT)}(\ol\Omega;\upbeta^*E)$ and the estimate
  \begin{equation}
  \label{EqWM}
    \|u\|_{\dot H_{\eop,\tbop}^{\sfs,(2\alpha_{\!\scri},\alpha_+,\alpha_\cT)}(\ol\Omega;\upbeta^*E)} \leq C\|f\|_{\dot H_{\eop,\tbop}^{\sfs-1,(2\alpha_{\!\scri}+2,\alpha_++2,\alpha_\cT)}(\ol\Omega;\upbeta^*E)}.
  \end{equation}
\end{thm}

\begin{rmk}[Mapping properties of $P$ and $P^*$]
\label{RmkWMMap}
  An equivalent formulation (using the uniqueness of forward solutions of $P$) is that
  \begin{equation}
  \label{EqWMMapDir}
  \begin{split}
    P &\colon \bigl\{ u\in\dot H_{\eop,\tbop}^{\sfs,(2\alpha_{\!\scri},\alpha_+,\alpha_\cT)}(\ol\Omega;\upbeta^*E) \colon P u\in\dot H_{\eop,\tbop}^{\sfs-1,(2\alpha_{\!\scri}+2,\alpha_++2,\alpha_\cT)}(\ol\Omega;\upbeta^*E) \bigr\} \\
      &\hspace{18em} \to \dot H_{\eop,\tbop}^{\sfs-1,(2\alpha_{\!\scri}+2,\alpha_++2,\alpha_\cT)}(\ol\Omega;\upbeta^*E).
  \end{split}
  \end{equation}
  is invertible. This also implies that
  \begin{align*}
    P^* &\colon \bigl\{ \tilde u\in\bar H_{\eop,\tbop}^{-\sfs+1,(-2\alpha_{\!\scri}-2,-\alpha_+-2,-\alpha_\cT)}(\Omega;\upbeta^*E) \colon P^*\tilde u\in\bar H_{\eop,\tbop}^{-\sfs,(-2\alpha_{\!\scri},-\alpha_+,-\alpha_\cT)}(\Omega;\upbeta^*E) \bigr\} \\
      &\hspace{18em} \to \bar H_{\eop,\tbop}^{-\sfs,(-2\alpha_{\!\scri},-\alpha_+,-\alpha_\cT)}(\Omega;\upbeta^*E)
  \end{align*}
  is invertible, since $P^*$ is Fredholm by Proposition~\ref{PropWMFred} and has trivial kernel and cokernel in view of the invertibility of~\eqref{EqWMMapDir}.
\end{rmk}

\begin{rmk}[The constant $C$]
\label{RmkWMConst}
  The proof shows also that $C$ in~\eqref{EqWM} can be taken to be a function of $|\tilde P|_k$ for some large but finite $k$ which depends only on $\|\sfs\|_{L^\infty}$. That is, the estimate~\eqref{EqWM} is uniform for bounded families of admissible wave type operators $P$. This continues to hold if we allow $P$ to be admissible with respect to a stationary wave type operator $P'_0$ where $P'_0$ varies over a compact family; see Remark~\ref{RmkWMStatConst}.
\end{rmk}

\begin{proof}[Proof of Theorem~\usref{ThmWM}]
  We claim that there exist $\digamma>1$ and an admissible wave type operator $P_\digamma$ (with respect to an admissible metric $g_\digamma$) so that $P_\digamma=P$ on $\Omega_\digamma=\ft_*^{-1}((\digamma,\infty])$, and so that the estimate~\eqref{EqWMFredAdj} holds for $P_\digamma$ \emph{without} an error term, i.e.
  \begin{equation}
  \label{EqWMdigamma}
    \|\tilde u\|_{\bar H_{\eop,\tbop}^{-\sfs+1,(-2\alpha_{\!\scri}-2,-\alpha_+-2,-\alpha_\cT)}(\Omega;\upbeta^*E)} \leq C \|P_\digamma^*\tilde u\|_{\bar H_{\eop,\tbop}^{-\sfs,(-2\alpha_{\!\scri},-\alpha_+,-\alpha_\cT)}(\Omega;\upbeta^*E)}.
  \end{equation}
  This implies the Theorem for $P_\digamma$ in place of $P$. To prove the Theorem for the original operator $P$, one notes that one can solve the forward problem $P u=f$ on $\Omega_{\digamma+1}\setminus\Omega$ (which is disjoint from $\iota^+\cup\cT^+$), using \cite[Theorem~6.4(1)]{HintzVasyScrieb}; denoting by $u_{\rm in}$ this local in time solution (which is quantitatively controlled by $f$) and by $\chi\in\CI(\R)$ a cutoff with $\chi=0$ on $(-\infty,\digamma]$ and $\chi=1$ on $[\digamma+1,\infty)$, the global solution $u$ on $\Omega$ is then $u=\chi u_{\rm in}+u'$ where $u'=(1-\chi)u$ solves the forward problem
  \[
    P_\digamma u'=P u'=(1-\chi)f-[P,\chi]u=(1-\chi)f-[P,\chi]u_{\rm in}=:f'.
  \]
  Since we have the estimate~\eqref{EqWM} for $u'$ and $f'$, we then obtain the estimate~\eqref{EqWM} also for $u$ and $f$.

  The guiding idea for the construction of $P_\digamma$ and the proof of~\eqref{EqWMdigamma} is that $P$ is very close to $P_0$ for late times $\ft_*\geq\digamma\gg 1$; we shall thus define $P_\digamma$ to be equal to $P$ for late times, and to be globally close to $P_0$. We shall then be able to estimate the error term in~\eqref{EqWMFredAdj} in terms of $P_0^*\tilde u=P_\digamma^*\tilde u+(P_0^*-P_\digamma^*)\tilde u$, where the second summand is small when $\digamma$ is large. Concretely, recalling $T=t_*^{-1}$, let $\phi\in\CIc([0,1))$ be identically $1$ near $0$, and set
  \[
    P_\digamma = \phi(\digamma T)P + \bigl(1-\phi(\digamma T)\bigr)P_0.
  \]
  This is an admissible wave type operator with respect to $g_\digamma$ and $P_0$, where
  \[
    g_\digamma^{-1}=\phi(\digamma T)g^{-1}+(1-\phi(\digamma T))g_0^{-1}.
  \]
  Since $\phi(\digamma T)$, as a function on $M$, is uniformly bounded in $\cA^0(M)$, the family $\{P_\digamma\colon\digamma>1\}$ is bounded in the sense of Remark~\ref{RmkWROrderOther}, i.e.\ any finite number of b-derivatives of the coefficients of
  \begin{equation}
  \label{EqWMPDiff}
    P_\digamma-P_0 = \phi(\digamma T)(P-P_0)
  \end{equation}
  is uniformly bounded. Note now that the microlocal estimates entering in the proof of Proposition~\ref{PropWMFred}, for bounded orders $\sfs_0,\sfs_1,\sfs$, only require some finite number of b-regularity of the coefficients of the wave type operator. Therefore, there exists a constant $C_0>0$ so that we have
  \begin{equation}
  \label{EqWMPdigammaStart}
  \begin{split}
    &\|\tilde u\|_{\bar H_{\eop,\tbop}^{-\sfs+1,(-2\alpha_{\!\scri}-2,-\alpha_+-2,-\alpha_\cT)}(\Omega)} \\
    &\qquad \leq C_0\Bigl( \|P_\digamma^*\tilde u\|_{\bar H_{\eop,\tbop}^{-\sfs,(-2\alpha_{\!\scri},-\alpha_+,-\alpha_\cT)}(\Omega)} + \|\tilde u\|_{\bar H_{\eop,\tbop}^{-\sfs_1+1,(-2\alpha_{\!\scri}-2-2\delta,-\alpha_+-2-\delta,-\alpha_\cT-\delta)}(\Omega)}\Bigr)
  \end{split}
  \end{equation}
  for all $\digamma>1$. (We drop the bundle $\upbeta^*E$ from the notation.) Concretely, in the notation of~\eqref{EqGAWNorm}, the constant $C_0$ can be taken to be a function of $|\tilde P|_k$ for some large but finite $k$.

  We then estimate the error term here by means of Corollary~\ref{CorWMP0Inv} (for $\sfs_1>\sfs+1$ in place of $\sfs$, and choosing $\sfs_1$ so that moreover $\sfs_1\pm 1$, $\alpha_{\!\scri}$, $\alpha_+$, $\alpha_\cT$ are admissible orders), which gives
  \begin{align*}
    &\|\tilde u\|_{\bar H_{\eop,\tbop}^{-\sfs_1+1,(-2\alpha_{\!\scri}-2,-\alpha_+-2,-\alpha_\cT)}(\Omega)} \\
    &\qquad \leq C_0'\|P_0^*\tilde u\|_{\bar H_{\eop,\tbop}^{-\sfs_1,(-2\alpha_{\!\scri},-\alpha_+,-\alpha_\cT)}(\Omega)} \\
    &\qquad \leq C_0'\|P_\digamma^*\tilde u\|_{\bar H_{\eop,\tbop}^{-\sfs_1,(-2\alpha_{\!\scri},-\alpha_+,-\alpha_\cT)}(\Omega)} + C_0'\|(P_\digamma-P_0)^*\tilde u\|_{\bar H_{\eop,\tbop}^{-\sfs_1,(-2\alpha_{\!\scri},-\alpha_+,-\alpha_\cT)}(\Omega)}.
  \end{align*}
  But in view of~\eqref{EqWMPDiff} and using the notation of Definition~\ref{DefGAW}\eqref{ItGAWStruct}, we have $P_\digamma-P_0=\phi(\digamma T)(\tilde P_1+\tilde P_2)$, and therefore
  \[
    P_\digamma-P_0 \to 0\quad\text{in}\quad x_{\!\scri}^2\rho_+^2 \cA^{(0,\ell_+-\eps,\ell_\cT-\eps)}\Diffetb^2(M),\qquad \digamma\to\infty,
  \]
  on $\ft_*\geq -1$ for any fixed $\eps>0$. Therefore, $P_\digamma-P_0\to 0$ in $\cA^{(2,2,0)}\Diffetb^2(M)$. Given any $\delta>0$, we can then choose $\digamma>1$ so large that
  \[
    C_0 C_0'\|(P_\digamma-P_0)^*\tilde u\|_{\bar H_{\eop,\tbop}^{-\sfs_1,(-2\alpha_{\!\scri},-\alpha_+,-\alpha_\cT)}(\Omega)} \leq \delta\|\tilde u\|_{\bar H_{\eop,\tbop}^{-\sfs_1+2,(-2\alpha_{\!\scri}-2,-\alpha_+-2,-\alpha_\cT)}(\Omega)}.
  \]
  Fixing $\delta>0$ small enough so that
  \[
    \delta\|\tilde u\|_{\bar H_{\eop,\tbop}^{-\sfs_1+2,(-2\alpha_{\!\scri}-2,-\alpha_+-2,-\alpha_\cT)}(\Omega)} \leq \frac12\|\tilde u\|_{\bar H_{\eop,\tbop}^{-\sfs+1,(-2\alpha_{\!\scri}-2,-\alpha_+-2,-\alpha_\cT)}(\Omega)},
  \]
  we can then absorb this term into the left hand side of~\eqref{EqWMPdigammaStart}. This establishes~\eqref{EqWMdigamma} and finishes the proof.
\end{proof}

\begin{rmk}[Small non-decaying perturbations]
\label{RmkWNonDecay}
  Theorem~\ref{ThmWM} implies the a priori estimate
  \begin{equation}
  \label{EqWNonDecay}
    \|\tilde u\|_{\Hext_\etbop^{-\sfs+1,(-2\alpha_{\!\scri}-2,-\alpha_+-2,-\alpha_\cT)}(\Omega;\upbeta^*E)} \leq C\|P^*\tilde u\|_{\Hext_\etbop^{-\sfs,(-2\alpha_{\!\scri},-\alpha_+,-\alpha_\cT)}(\Omega;\upbeta^*E)}.
  \end{equation}
  (Indeed, this follows by duality from the solvability for $P$ and the estimate~\eqref{EqWM}.) Consider an operator
  \begin{equation}
  \label{EqWNonDecayR}
    R\in\cA^{(2,2,0)}\Diffetb^1(M;\upbeta^*E);
  \end{equation}
  we stress that we do \emph{not} require $R$ to have decaying coefficients as $t_*\to\infty$ in spatially compact sets. Then there exists $\eps_0>0$ so that the estimate~\eqref{EqWNonDecay} holds for $P+\eps R$ in place of $P$ for all $\eps\in\C$ with $|\eps|<\eps_0$; indeed, it suffices to choose $\eps_0$ so that the operator norm of $R^*\colon\Hext_\etbop^{-\sfs+1,(-2\alpha_{\!\scri}-2,-\alpha_+-2,-\alpha_\cT)}(\Omega;\upbeta^*E)\to\Hext_\etbop^{-\sfs,(-2\alpha_{\!\scri},-\alpha_+,-\alpha_\cT)}(\Omega;\upbeta^*E)$ is less than $\frac{1}{2 C\eps_0}$. By duality, this implies the (unique) solvability of the forward problem
  \[
    (P+\eps R)u=f \in \Hsupp_\etbop^{\sfs-1,(2\alpha_{\!\scri}+2,\alpha_++2,\alpha_\cT)}(\ol\Omega;\upbeta^*E)
  \]
  for all $|\eps|<\eps_0$, with the solution $u\in\Hsupp_\etbop^{\sfs,(2\alpha_{\!\scri},\alpha_+,\alpha_\cT)}(\ol\Omega;\upbeta^*E)$ obeying the estimate~\eqref{EqWM}. One can in fact relax~\eqref{EqWNonDecayR} further to
  \begin{equation}
  \label{EqWNonDecayR2}
    R \in \cA_\etbop^{(2,2,0)}\Diffetb^1(M;\upbeta^*E),
  \end{equation}
  where we denote by $\cA_\etbop^{(2,2,0)}(M)$ the space of all $a\in x_{\!\scri}^2\rho_+^2 L^\infty(M)$ so that $B a\in x_{\!\scri}^2\rho_+^2 L^\infty(M)$ for all $B\in\Diffetb(M)$; the reason is that multiplication by elements of $\cA_\etbop^{(0,0,0)}(M)$ defines bounded linear maps on every weighted edge-3b-Sobolev space. One class of operators $R$ satisfying~\eqref{EqWNonDecayR2} includes\footnote{This entails b-regularity at $\scri^+$, whereas the even weaker edge-regularity would still be sufficient.}
  \[
    R=V(t_*,x),\qquad \Bigl|\rho_+^{-j}\pa_{t_*}^j\la x\ra^{|\alpha|}\pa_x^\alpha V\Bigr| \leq C_{j\alpha}\la x\ra^{-2},
  \]
  where we take $\rho_+=\frac{t}{\la r\ra\la t-r\ra}$ as in~\eqref{EqIRhos} with $t=t_*+r$. (In spatially compact regions, this amounts to uniform bounds on $V$ and all its $(t,x)$-derivatives.) This is more permissive than the assumptions in \cite{MetcalfeSterbenzTataruLED} as far as spatial decay is concerned, though we require more regularity here.
\end{rmk}

Theorem~\ref{ThmWM} result can be combined with the solvability theory near $I^0\cap\scri^+$ from \cite[\S9]{HintzVasyScrieb}; we present a concrete example:

\begin{cor}[Forward solutions for $P$ on a larger spacetime domain]
\label{CorWt}
  Let $t\in\CI(M^\circ)$ be a function with $t=t_*+r$ for large $r$, and suppose $\dd t$ is past timelike on $M^\circ$. Moreover:
  \begin{enumerate}
  \item suppose $\tau\in\CI(M^\circ)$ has the following properties: $\tau=t/r$ near $I^0$; we have $\{\tau\geq 0\}\subset\{t_*+r\geq-\half r-1\}$; the differential $\dd\tau$ is past timelike in $\{\tau>0\}$; and $\ft_*^{-1}(0)\subset\{\tau>0\}$. Denote by $\ol\cD\subset M$ the closure of $\tau^{-1}([0,\infty])$;
  \item let $\sfs\in\CI({}^{\bop,\eop,\tbop}S^*M\cap\ol\cD)$, $\alpha_0$, $\alpha_{\!\scri}$, $\alpha_+$, $\alpha_\cT\in\R$ be $P$-admissible in the following sense: the conditions of Definition~\usref{DefWMOrders} are verified; $\sfs$ is constant near $\pa\cR_{\rm c}^\pm$, $\pa\cR_{\rm in,-}^\pm$, $\pa\cR_{\rm out}^\pm$, which are the boundaries at fiber infinity of the radial sets $\cR_{\rm c}^\pm$, $\cR_{\rm in,-}^\pm\subset\Seb^*_{I^0\cap\scri^+}M$ and $\cR_{\rm out}^\pm\subset\Seb^*_{\scri^+}M$ defined in \cite[Lemma~4.1]{HintzVasyScrieb}; furthermore $\pm\sfH\sfs\leq 0$ where $\sfH=\rho_\infty H_{\rho_0^{-2}x_{\!\scri}^{-2}\rho_+^{-2}p}$ is the rescaled Hamiltonian vector field where $\rho_\infty$ is a classical elliptic symbol on $\Tbetb^*M$ of order $-1$; and finally
  \begin{equation}
  \label{EqWtOrders}
    \alpha_{\!\scri}<\alpha_0+\frac12,\qquad
    \sfs|_{\pa\cR_{\rm c}}>-(\alpha_0-\alpha_{\!\scri})-\Bigl(-\frac12+\ubar p_1-\alpha_{\!\scri}\Bigr).
  \end{equation}
  Suppose, in fact, that $\sfs-1$ and $\sfs+2$ are $P$-admissible orders for the weights $\alpha_0,\alpha_{\!\scri},\alpha_+,\alpha_\cT$.
  \end{enumerate}
  Then there exists $C>0$ so that for any $f\in\dot H_{\bop,\eop,\tbop}^{\sfs-1,(\alpha_0+2,2\alpha_{\!\scri}+2,\alpha_++2,\alpha_\cT)}(\ol\cD;\upbeta^*E)$, the unique forward solution $u$ of $P u=f$ satisfies $u\in\dot H_{\bop,\eop,\tbop}^{\sfs,(\alpha_0,2\alpha_{\!\scri},\alpha_+,\alpha_\cT)}(\ol\cD;\upbeta^*E)$ and
  \[
    \|u\|_{\dot H_{\bop,\eop,\tbop}^{\sfs,(\alpha_0,2\alpha_{\!\scri},\alpha_+,\alpha_\cT)}(\ol\cD;\upbeta^*E)} \leq C\|f\|_{\dot H_{\bop,\eop,\tbop}^{\sfs-1,(\alpha_0+2,2\alpha_{\!\scri}+2,\alpha_++2,\alpha_\cT)}(\ol\cD;\upbeta^*E)}.
  \]
\end{cor}

The existence of an order function $\sfs$ satisfying all these requirements follows from a simple extension of the proof of Lemma~\ref{LemmaWROrder} using the expressions for the edge-b-Hamiltonian vector field near $I^0\cap\scri^+$ from \cite[\S4.1]{HintzVasyScrieb}; we omit the details here.

\begin{proof}[Proof of Corollary~\usref{CorWt}]
  The assumptions on $\tau$ and $t$ merely serve to guarantee the solvability of the forward problem in the appropriate weighted edge-b-Sobolev space up to $\ft_*^{-1}(0)$, from where Theorem~\ref{ThmWM} can be used to continue the solution. Concretely, under the stated assumptions on $\tau$, one can solve $P u=f$ in $\{\tau\leq\tau_0\}\subset M\setminus(\scri^+\cup\iota^+\cup\cT^+)$ for any fixed $\tau_0\in[0,1)$; see also the first step of the proof of \cite[Theorem~9.2]{HintzVasyScrieb}. For $\tau_0$ sufficiently close to $1$, one can then continue the solution (with control in the edge-b-Sobolev space with orders $\sfs$, $\alpha_0$, $2\alpha_{\!\scri}$) by means of \cite[Theorem~6.4(1)]{HintzVasyScrieb} to a neighborhood of $I^0\cap\scri^+$, and indeed by exploiting the timelike nature of $\tau$ and $t$ to a neighborhood of $\ft_*^{-1}(0)\subset M$. (The threshold conditions~\eqref{EqWtOrders} are used in this step.) From there, $u$ can be extended (with quantitative control) by means of Theorem~\ref{ThmWM}.
\end{proof}

\subsection{Higher b-regularity; pointwise decay}
\label{SsWb}

Even in the case that the forcing $f$ in Theorem~\ref{ThmWM} is smooth and compactly supported in $M^\circ\cap\Omega$, the control on the solution $u$ provided by~\eqref{EqWM} is rather weak in the sense that the degree of edge-3b-regularity of $u$ we prove near the flow-out $\cW_{\rm out}=\bigcup_\pm\cW_{\rm out}^\pm$ (see~\eqref{EqWFlowIUnstable}) is less than $-\half+\ubar S-\alpha_++\alpha_\cT$ (cf.\ Lemma~\ref{LemmaWROrder}). As a corollary of the results in this section, one can improve this control to infinite order b-regularity of $u$.

We work in the domain $\Omega=\ft_*^{-1}((0,\infty])$, see~\eqref{EqWMOmega}.

\begin{definition}[Mixed function spaces]
\label{DefWbMix}
  For $\sfs\in\CI(\Setb^*_{\ol\Omega}M)$, $\alpha_{\!\scri},\alpha_+,\alpha_\cT\in\R$ and $k\in\N_0$, we denote by
  \[
    \dot H_{\eop,\tbop;\bop}^{(\sfs;k),(2\alpha_{\!\scri},\alpha_+,\alpha_\cT)}(\ol\Omega)
  \]
  the space of all $u$ so that $A u\in\dot H_{\eop,\tbop}^{\sfs,(2\alpha_{\!\scri},\alpha_+,\alpha_\cT)}(\ol\Omega)$ for all $A\in\Diffb^k(M)$. Spaces of sections of vector bundles over $M$ are defined analogously.
\end{definition}

We can give these spaces the structure of a Hilbert space by defining a squared norm as the sum of squared norms of $A_i u$ where $\{A_i\}\subset\Diffb^k(M)$ is a finite set of operators generating $\Diffb^k(M)$ over $\CI(M)$. On such mixed function spaces, we have $\ebop$-microlocal elliptic regularity and real principal type propagation away from $\cT^+$, and also radial point estimates near $\scri^+$, by the results of \cite[\S\S5.2--5.3]{HintzVasyScrieb}. We need to complement these results with $\tbop$-microlocal estimates near $\cT^+$.

In the following, we work near $\cT^+$, only record weights at $\iota^+$ and $\cT^+$, and only keep the 3b-structure in the notation; the function spaces we use are thus the $(\tbop;\bop)$-spaces introduced in~\eqref{EqMFtbb}. Let $\tilde\chi_\cT,\chi_\cT\in\CI(M)$ be identically $1$ in a collar neighborhood of $\cT^+$, supported in $\{t_*+r>-\half r-1\}$, and so that $\tilde\chi_\cT=1$ near $\supp\chi_\cT$. The Schwartz kernels of all ps.d.o.s below are assumed to be supported in the interior of $\tilde\chi_\cT^{-1}(1)\times\tilde\chi_\cT^{-1}(1)$.

\begin{lemma}[$\tbop;\bop$-elliptic estimates]
\label{LemmaWbEll}
  Let $\sfs\in\CI(\Stb^*M)$, $m\in\R$, $k\in\N_0$, $\alpha=(\alpha_+,\alpha_\cT)\in\R^2$. Let $A\in\Psitb^m(M)$ and $B,G\in\Psitb^0(M)$, and suppose that $\WFtb'(B)\subset\Elltb(A)\cap\Elltb(G)\cap\Ttb^*_{\supp\chi_\cT}M$. Let $N\in\R$. Then there exists $C>0$ so that the estimate
  \begin{equation}
  \label{EqWbEll}
    \|B u\|_{H_{\tbop;\bop}^{(\sfs;k),\alpha}(M)} \leq C\Bigl( \|G A u\|_{H_{\tbop;\bop}^{(\sfs-m;k),\alpha}(M)} + \|\tilde\chi_\cT u\|_{H_{\tbop;\bop}^{(-N;k),\alpha}(M)}\Bigr)
  \end{equation}
  holds in the strong sense that the left hand side is finite if all terms on the right hand side are, and the estimate holds.
\end{lemma}

More generally, one can allow $A$ to have conormal coefficients, i.e.\ we can allow $A\in\cA^{(((0,0),\eps),((0,0),\eps))}\Psitb^m(M)$, $\eps>0$; moreover, one can consider operators $A$ which act on sections of a vector bundle $E\to M$. We leave it to the reader to make the required notational changes here and in Lemma~\ref{LemmaWbProp} below.

\begin{proof}[Proof of Lemma~\usref{LemmaWbEll}]
  Let $Q\in\Psitb^{-m}(M)$ be a microlocal elliptic parametrix of $G A$ near $\WFtb'(B)$, so $I=Q G A+R$ where $R\in\Psitb^0(M)$ with $\WFtb'(R)\cap\WFtb'(B)=\emptyset$. Then $B u=B Q(G A u)+B R u$. Since $B Q\in\Psitb^{-m}(M)$ defines a bounded map $H_{\tbop;\bop}^{(\sfs-m;k),\alpha}(M)\to H_{\tbop;\bop}^{(\sfs;k),\alpha}(M)$, and since $B R=B R\tilde\chi_\cT$ is bounded as a map $H_{\tbop;\bop}^{(-N;k),\alpha}(M)\to H_{\tbop;\bop}^{(\sfs;k),\alpha}(M)$ for all $N$, the estimate~\eqref{EqWbEll} follows.
\end{proof}

\begin{lemma}[Real principal type $\tbop;\bop$-propagation estimate]
\label{LemmaWbProp}
  Let $\sfs\in\CI(\Stb^*M)$, $m\in\R$, $k\in\N_0$, $\alpha=(\alpha_+,\alpha_\cT)\in\R^2$. Suppose $A\in\Psitb^m(M)$ has a real homogeneous principal symbol, and $\sfH\sfs\leq 0$ where $\sfH\in\Vb(\Stb^*M)$ is the restriction of $\rho_\infty^{m-1}H_{\sigmatb^m(A)}$ to fiber infinity $\Stb^*M\subset\ol{\Ttb^*}M$, with $\rho_\infty$ denoting a defining function of fiber infinity. Suppose $B,E,G\in\Psitb^0(M)$ are such that $\WFtb'(B)\subset\Elltb(G)$, and so that all backward null-bicharacteristics from $\WFtb'(B)\cap\Char(A)$ reach $\Elltb(E)$ in finite time while remaining in $\Elltb(G)$. Then for any fixed $N\in\R$, we have the estimate
  \begin{equation}
  \label{EqWbProp}
    \|B u\|_{H_{\tbop;\bop}^{(\sfs;k),\alpha}(M)} \leq C\Bigl( \|G A u\|_{H_{\tbop;\bop}^{(\sfs-m;k),\alpha}(M)} + \|E u\|_{H_{\tbop;\bop}^{(\sfs;k),\alpha}(M)} + \|\tilde\chi_\cT u\|_{H_{\tbop;\bop}^{(-N;k),\alpha}(M)}\Bigr).
  \end{equation}
  This holds in the strong sense as in Lemma~\usref{LemmaWbEll}.
\end{lemma}
\begin{proof}
  We argue by induction on $k$, the case $k=0$ being a standard propagation estimate which is proved using a positive commutator argument (and which was tacitly used already in the proof of Proposition~\ref{PropWR}). The inductive step (i.e.\ proving~\eqref{EqWbProp} for $k+1$ in place of $k$), in which for fixed $B$ we shall need to enlarge the elliptic sets of $G,E$ slightly, can be proved in a manner similar to \cite[Proposition~5.15]{HintzVasyScrieb}. First of all, since $H_{\tbop;\bop}^{(\sfs;k+1),\alpha}(M)\subset H_{\tbop;\bop}^{(\sfs+1;k),\alpha}(M)$, we may apply the estimate~\eqref{EqWbProp} with $\sfs+1$ in place of $\sfs$ in order to obtain an estimate on $\|B u\|_{H_{\tbop;\bop}^{(\sfs+1;k),\alpha}(M)}$. Let $V\in\Vb(M)$ be a time-dilation operator; since over $\supp\tilde\chi_\cT$ any b-vector field is the sum of a 3b-vector field and a smooth multiple of $V$, we have
  \[
    \|B u\|_{H_{\tbop;\bop}^{(\sfs;k+1),\alpha}(M)} \leq C\Bigl( \|B u\|_{H_{\tbop;\bop}^{(\sfs+1;k),\alpha}(M)} + \|V B u\|_{H_{\tbop;\bop}^{(\sfs;k),\alpha}(M)} \Bigr),
  \]
  with the first term on the right hand side already controlled. Moreover, since $V B u=B V u+[V,B]u$ where $[V,B]\in\Psitb^0(M)$ by Lemma~\ref{LemmaMOMComm}, this estimate remains valid (with a different $C$) if we replace $B u$, resp.\ $V B u$ in the first, resp.\ second term by $\tilde B u$, resp.\ $B V u$, where $\tilde B\in\Psitb^0(M)$ is elliptic on the operator wave front set of $B$. But we can estimate $\|B V u\|_{H_{\tbop;\bop}^{(\sfs;k),\alpha}(M)}$ by plugging $V u$ into~\eqref{EqWbProp} instead of $u$. Indeed, writing
  \[
    G A V u = V G A u + [G A,V]u,
  \]
  one applies the inductive hypothesis (with $B$ replaced by an operator $\tilde B$ whose elliptic set contains $\WFtb'(G A)$) to estimate the term $\|[G A,V]u\|_{H_{\tbop;\bop}^{(\sfs-m;k),\alpha}(M)}$ by $\|\tilde B u\|_{H_{\tbop;\bop}^{(\sfs;k),\alpha}(M)}$; we note here that $[G A,V]\in\Psitb^m(M)$ by Lemma~\ref{LemmaMOMComm}. Similarly, we have $E V u=V E u+[E,V]u$ with $[E,V]\in\Psitb^0$, and therefore $\|[E,V]u\|_{H_{\tbop;\bop}^{(\sfs;k),\alpha}(M)}\lesssim\|\tilde E u\|_{H_{\tbop;\bop}^{(\sfs;k),\alpha}(M)}$ if $\tilde E\in\Psitb^0(M)$ is chosen to satisfy $\Elltb(\tilde E)\supset\WFtb(E)$. Finally, due to the support assumptions on the Schwartz kernels, we can replace $\tilde\chi_\cT u$ in the final term in~\eqref{EqWbProp} by $\tilde\chi_\cT' u$ for a suitable cutoff $\tilde\chi_\cT'$ with support in $\tilde\chi_\cT^{-1}(1)$. But then $\tilde\chi_\cT' V u=V\tilde\chi_\cT' u+[V,\tilde\chi_\cT']u$, and the $H_{\tbop;\bop}^{(-N;k),\alpha}(M)$-norm of the second summand is bounded by the norm of $\tilde\chi_\cT u$. This completes the inductive step.
\end{proof}

Lemmas~\ref{LemmaWbEll} is a quantitative version of the qualitative wave front set statement
\[
  \WF_{\tbop;\bop}^{(\sfs;k),\alpha}(u)\subset\WF_{\tbop;\bop}^{(\sfs-m;k),\alpha}(A u)\cup\Char(A),
\]
where we use the notation~\eqref{EqMFtbbWF}. Similarly, Lemma~\ref{LemmaWbProp} is a quantitative version of the statement that if $\gamma\colon[0,1]\to\Stb^*M$ is an integral curve of $\sfH$ inside the set $\Char(A)\setminus\WF_{\tbop;\bop}^{(\sfs-m+1;k),\alpha}(A u)$ then $\gamma(0)\notin\WF_{\tbop;\bop}^{(\sfs;k),\alpha}(u)$ implies $\gamma(t)\notin\WF_{\tbop;\bop}^{(\sfs;k),\alpha}(u)$ for $t\in[0,1]$.

The radial point estimates for admissible wave type operators $P$ over $\pa\cT^+$ are more delicate. Since we are interested in additional b-regularity only for solutions of the forward problem (as opposed to the adjoint problem), we shall only prove estimates for $P$ here. Recall from Proposition~\ref{PropWTOut} that propagation through $\cR_{\cT,\rm out}^\pm$ imposes an \emph{upper} bound on the 3b-regularity order $\sfs$ of the wave $u$. Thus, if one wishes to add $k$ degrees of b-regularity on top of 3b-regularity, one might expect that the sum $\sfs+k$ needs to satisfy the same upper bound; this would impose a fixed upper bound on the total b-regularity of $u$ near $\cW_{\rm out}\cap\Tb^*_{(\iota^+)^\circ}M$. We will demonstrate, to the contrary, that a threshold condition is only required for $\sfs$, while $k$ can be arbitrarily large. This is related to the notion of \emph{module regularity}---here at the flowout $\cW_{\rm out}$---as introduced in \cite{HassellMelroseVasySymbolicOrderZero} and used e.g.\ in \cite{HaberVasyPropagation,BaskinVasyWunschRadMink,HintzVasySemilinear,GellRedmanHassellShapiroZhangHelmholtz,HintzConicProp}. There are two key differences, however. First, the extra b-regularity is of a strictly stronger character than the background 3b-regularity. Second, we shall need to \emph{assume} a priori that $u$ has the amount $k$ of b-regularity; under this assumption, we show that 3b-regularity (relative to $k$ b-derivatives) propagates microlocally. We explain in the proof of Theorem~\ref{ThmWb} below how to obtain this a priori assumption. The main input in the proofs of the microlocal 3b-propagation results relative to a fixed amount of b-regularity is the fact that time dilation vector fields can be commuted through 3b-equations with well-controlled errors by Lemmas~\ref{LemmaM3bbComm2} and \ref{LemmaMOMComm}.

Recall $\rho_+=\rho=r^{-1}$, which is a defining function of $\iota^+$ on $\supp\chi_\cT$. Let $\chi\in\CIc([0,1)_\rho)$ be identically $1$ near $\rho=0$. In order to test for b-regularity, we shall use the vector fields
\begin{equation}
\label{EqWbVF}
\begin{alignedat}{2}
  V_0 &= -t_*\pa_{t_*},&\qquad
  V_1 &= \chi\rho\pa_\rho, \\
  V_j &= \chi\Omega_j\quad (j=2,\ldots,N'), &\qquad
  V_k &= (1-\chi)W_k\quad (k=N'+1,\ldots,N),
\end{alignedat}
\end{equation}
on $\R_{t_*}\times\R^n$, where the $\Omega_j$ span $\cV(\Sph^{n-1})$ over $\CI(\Sph^{n-1})$, and the $W_k$ span $\cV(X^\circ)$ over $\CI(X^\circ)$; thus, the $V_i$ span $\Vb(M_0)$ over $\CI(M_0)$, and their lifts to $M$, which we denote by $V_i$ still, span $\Vb(M)$ over $\CI(M)$. Note that $V_1,\ldots,V_N\in\Vtb(M)$, and only $V_0$ (which is a time dilation vector field, see Definition~\ref{DefM3bb}) does not lie in $\Vtb(M)$. If the vector bundle $E$ is trivial, we set $X_i=V_i$ for all $i$; otherwise, we take $X_0\in\Diffb^1(M;\upbeta^*E)$ to be a time dilation operator, and $X_i\in\Difftb^1(M;\upbeta^*E)$ for $i=1,\ldots,N$ to be operators with scalar principal symbol equal to that of $V_i$.

\begin{lemma}[Commutators of test operators]
\label{LemmaWbComm}
  Recall the definitions~\eqref{EqWRadT}, and write $\cR_{\cT,\rm out}=\bigcup_\pm\cR_{\cT,\rm out}^\pm$ and $\cR_{\cT,\rm in}=\bigcup_\pm\cR_{\cT,\rm in}^\pm$. Then for $j=0,\ldots,N$, we can write
  \begin{equation}
  \label{EqWbComm}
    \chi_\cT[P,X_j] = y_j P + \rho_+^2\biggl(R_j + \sum_{\ell=0}^N Y_{j,\ell}X_\ell\biggr),
  \end{equation}
  where, recording weights only at $\iota^+$ and $\cT^+$, we have
  \[
    y_j\in\CI(M)+\cA^{(\ell_+,\ell_\cT)}(M), \qquad R_j,\ Y_{j,\ell}\in(\CI+\cA^{(\ell_+,\ell_\cT)})\Difftb^1(M;\upbeta^*E),
  \]
  with coefficients supported near $\supp\chi_\cT$, and with $\sigmatb^1(Y_{j,\ell})=0$ at $\cR_{\cT,\rm out}\cup\cR_{\cT,\rm in}$ for $\ell=0,\ldots,N$.
\end{lemma}
\begin{proof}
  We omit the bundle $\upbeta^*E$ from the notation. Recall from Definition~\ref{DefGAW} that $\tilde\chi_\cT(P-P_0)\in\cA^{(2+\ell_+,\ell_\cT)}\Difftb^2(M)$. We will repeatedly use the observation that for any $Q\in\Difftb^k(M)$ we can write
  \[
    \chi_\cT Q = \sum_{\ell=0}^N Q_\ell X_\ell + R,\qquad Q_i,R\in\Difftb^{k-1}(M),
  \]
  and in fact one can take $Q_0=\rho_\cT Q'_0$ with $Q_0\in\Difftb^{k-1}(M)$ where $\rho_\cT=\frac{r}{t_*}$. This follows by considering the 3b-principal symbol of $\chi_\cT Q$ and noting that any fiber-linear function on $\Ttb^*M$ with support on $\supp\tilde\chi_\cT$ can be written as a linear combination, with $\CI(M)$ coefficients, of $\sigmatb^1(\rho_\cT X_0)$, $\sigmatb^1(X_\ell)$ ($\ell=1,\ldots,N$). By multiplying $Q_i,R$ on the left by $\tilde\chi_\cT$, one can moreover arrange for the coefficients of $Q_i$ for $i=0,\ldots,N$ to be supported in $\supp\tilde\chi_\cT$.

  \pfstep{The case $j=0$.} Lemma~\ref{LemmaM3bbComm2} and the fact that $\rho_+=r^{-1}$ commutes with $X_0$ give a decomposition
  \begin{equation}
  \label{EqWbCommX0}
    \chi_\cT[P,X_0] = \rho_\cT P^\flat X_0 + \rho_\cT P_0^\sharp+\tilde P^\sharp,
  \end{equation}
  where $P^\flat\in\rho_+^2\Difftb^1(M)$, $P_0^\sharp\in\rho_+^2\Difftb^2(M)$, and $\tilde P^\sharp\in\cA^{(2+\ell_+,\ell_\cT)}\Difftb^2(M)$. We may multiply $P^\flat$, $P_0^\sharp$, and $\tilde P^\sharp$ with $\tilde\chi_\cT$ on the left; this preserves the validity of $\chi_\cT[P,X_0]=\rho_\cT P^\flat X_0+\rho_\cT P_0^\sharp+\tilde P^\sharp$ but ensures that the coefficients of $P^\flat$, $P_0^\sharp$, and $\tilde P^\sharp$ are supported near $\supp\chi_\cT$. Using the above observation, we can write $\tilde P^\sharp\in\cA^{(2+\ell_+,\ell_\cT)}\Difftb^2(M)$ as
  \[
    \tilde P^\sharp = \tilde C_0^\sharp \rho_\cT X_0 + \sum_{\ell=1}^N \tilde C_\ell^\sharp X_\ell + \tilde R^\sharp,
  \]
  where $\tilde C_i^\sharp$, $\tilde R^\sharp\in\cA^{(2+\ell_+,\ell_\cT)}\Difftb^1(M)$. We can write $\rho_\cT P_0^\sharp$ in a similar manner but with $C_i^\sharp,R^\sharp\in\rho_+^2\rho_\cT\Difftb^1(M)$. Combining this with~\eqref{EqWbCommX0} gives~\eqref{EqWbComm} for $j=0$, with $y_0=0$, $R_0=\rho_+^{-2}(\tilde R^\sharp+R^\sharp)$, $Y_{0,0}=\rho_\cT(P^\flat+\rho_+^{-2}\tilde C_0^\sharp+\rho_+^{-2}C_0^\sharp)$, and $Y_{0,\ell}=\rho_+^{-2}(\tilde C_\ell^\sharp+C_\ell^\sharp)$ for $\ell=1,\ldots,N$.

  \pfstep{The case $j=1$.} Consider next the commutator of $P$ with $X_1$. We use the notation of~\eqref{EqWtbCoords}--\eqref{EqWHam3b}. Note that $p:=\sigmatb^2(P)=\rho_+^2 G_\tbop$, and therefore
  \[
    H_p=\rho_+^2 H_{G_\tbop}+p \rho_+^{-2}H_{\rho_+^2}.
  \]
  Since $\sigmatb^1(i^{-1}X_1)=-\chi(\rho_+)\xi_\tbop$, the 3b-principal symbol of $-i[P,i^{-1}X_1]\in\rho_+^2(\CI+\cA^{(\ell_+,\ell_\cT)})\Difftb^2$ is
  \[
    H_p(\chi\xi_\tbop) = (H_p\chi)\xi_\tbop + \chi H_p\xi_\tbop,
  \]
  where $H_p\chi\in(\CI+\cA^{(0,\ell_\cT)})P^1(\Ttb^*M)$ vanishes near $\iota^+$ and thus is a linear combination of the symbols of $\rho_\cT V_0$, $V_1$, $\ldots$, $V_N$ with coefficients vanishing near $\iota^+$. Write $\cE=\rho_+^2(\rho_+\CI+\cA^{(\ell_+,\ell_\cT)})P^2(\Ttb^*M)$. We moreover have
  \begin{equation}
  \label{EqWbCommX1}
    \rho_+^2 H_{G_\tbop}\xi_\tbop \equiv 2\rho_+^2\sigma_\tbop\xi_\tbop \equiv -2 i^{-1}\rho_+^2\rho_\cT\xi_\tbop\cdot\sigmatb^1(X_0) \bmod \cE.
  \end{equation}
  Since $\rho_\cT=0$ at $\cR_{\cT,\rm out}\cup\cR_{\cT,\rm in}$, this implies that we can write
  \[
    \chi_\cT[P,X_1] \equiv y_1 P + \rho_+^2\biggl( Y_{1,0} X_0 + \sum_{\ell=1}^N Y_{1,\ell}X_\ell\biggr) \bmod \rho_+^2(\CI+\cA^{(\ell_+,\ell_\cT)})\Difftb^1
  \]
  where the principal symbol of $Y_{1,0}$ vanishes over $\cT^+$, and the principal symbol of $Y_{1,\ell}$ for $\ell=1,\ldots,N$ vanishes at $\iota^+$. Taking $R_1$ to be $\rho_+^{-2}$ times the error term here gives~\eqref{EqWbComm}.

  \pfstep{The case $j=2,\ldots,N'$.} For spherical derivatives, we compute, with $\cE$ as above and for smooth functions $f^i$ on $\Sph^{n-1}$,
  \[
    \rho_+^2 H_{G_\tbop}\bigl(f^i(\omega)\eta_{\tbop,i}\bigr) \equiv 2\slg^{j k}(\omega)(\pa_{\omega^k}f^i)\eta_{\tbop,j}\eta_{\tbop,i} -f^i\pa_{\omega^i}\slg^{j k}(\omega)\eta_{\tbop,j}\eta_{\tbop,k} \bmod \cE.
  \]
  Since $\eta_\tbop=0$ at $\cR_{\cT,\rm out}\cup\cR_{\cT,\rm in}$, this implies that we can arrange~\eqref{EqWbComm}.

  \pfstep{The case $j=N'+1,\ldots,N$.} Since $V_j=0$ near $\iota^+$, we can multiply any decomposition~\eqref{EqWbComm} with a cutoff $\chi_0\in\CI(M)$ which is supported in $\tilde\chi_\cT^{-1}(1)$ and vanishes outside a neighborhood of $\supp(1-\chi)$; this ensures that the coefficients of $y_j,R_j,Y_{j,\ell}$ vanish near $\iota^+$. The proof is complete.
\end{proof}

We can now generalize Proposition~\ref{PropWTOut}\eqref{ItWTOutDir} to the case of relative b-regularity:

\begin{prop}[Propagation through $\cR_{\cT,\rm out}^\pm$ with relative b-regularity]
\label{PropWbOut}
  Let $s\in\R$, $k\in\N_0$, $\alpha_+,\alpha_\cT\in\R$. Let
  \[
    u \in H_{\tbop;\bop}^{(-\infty;k),(\alpha_+,\alpha_\cT)}(M;\upbeta^*E),\qquad
    f := P u\in H_{\tbop;\bop}^{(-\infty;k),(\alpha_++2,\alpha_\cT)}(M;\upbeta^*E).
  \]
  Assume that the threshold condition~\eqref{EqWTOutDirThr} holds. Suppose that $\WF_{\tbop;\bop}^{(s-1;k),(\alpha_++2,\alpha_\cT)}(f)\cap\pa\cR_{\cT,\rm out}^\pm=\emptyset$, and $\WF_{\tbop;\bop}^{(s;k),(\alpha_+,\alpha_\cT)}(u)\cap(\cU\setminus\pa\cR_{\cT,\rm out}^\pm)=\emptyset$ for some neighborhood $\cU\subset\Stb^*_{\cT^+}M$ of $\pa\cR_{\cT,\rm out}^\pm$ over $\cT^+$. Then $\WF_{\tbop;\bop}^{(s;k),(\alpha_+,\alpha_\cT)}(u)\cap\pa\cR_{\cT,\rm out}^\pm=\emptyset$.
\end{prop}

Similarly, Proposition~\ref{PropWTIn}\eqref{ItWTInDir} continues to hold on spaces with $k$ degrees of b-regularity, and again we stress that the proof gives quantitative estimates as in~\eqref{EqWTDirEst}. We leave the purely notational adaptation to the reader.

\begin{proof}[Proof of Proposition~\usref{PropWbOut}]
  We argue by induction on $k$. The case $k=0$ is the content of Proposition~\ref{PropWTOut}\eqref{ItWTOutDir}. For the inductive step, we assume that the Proposition has been proved for the value $k\in\N_0$, and for all principally scalar operators with the same principal symbol as $P$ for which the threshold condition is~\eqref{EqWTOutDirThr} still; we shall proceed to prove the Proposition for the value $k+1$. Let
  \[
    u'=(X_0 u,\ldots,X_N u)\in H_{\tbop;\bop}^{(-\infty;k),(\alpha_+,\alpha_\cT)}\bigl(M;(\upbeta^*E)^{N+1}\bigr).
  \]
  By Lemma~\ref{LemmaWbComm}, we can write
  \begin{equation}
  \label{EqWbOutPf}
    \chi_\cT(P' - Y')u' = R' u + f',
  \end{equation}
  where $P'=P\otimes\Id_{N+1}=\diag(P,P,\ldots,P)$ is a wave type operator on the bundle $(\upbeta^*E)^{N+1}$, and
  \begin{align*}
    Y'=(\rho_+^2 Y_{j,\ell})_{j,\ell=0,\ldots,N}, \quad
    R'=\diag\bigl(\rho_+^2 R_j\colon j=0,\ldots,N \bigr), \quad
    f'=(y_j f)_{j=0,\ldots,N}.
  \end{align*}
  The inductive hypothesis implies $\WF_{\tbop;\bop}^{(s-1;k),(\alpha_++2,\alpha_\cT)}(R' u)\cap\pa\cR_{\cT,\rm out}^\pm=\emptyset$, and the a priori regularity assumption on $u$ implies $\WF_{\tbop;\bop}^{(s;k),(\alpha_+,\alpha_\cT)}(u')\cap(\cU\setminus\pa\cR_{\cT,\rm out}^\pm)=\emptyset$. Since $P'-Y'$ is principally scalar, with the same principal symbol as $P$, we can apply the inductive hypothesis to equation~\eqref{EqWbOutPf}; since $P'$ acts component-wise as $P$, and since the principal symbol of the subprincipal term $Y'$ vanishes at $\cR_{\cT,\rm out}$, the threshold condition on the 3b-differential order for propagation through $\cR_{\cT,\rm out}$ for the operator $P'-Y'$ is the same as the threshold condition for $P$ itself. We conclude that $\WF_{\tbop;\bop}^{(s;k),(\alpha_+,\alpha_\cT)}(u')\cap\pa\cR_{\cT,\rm out}^\pm=\emptyset$. Since the operators $\chi_\cT X_0,\ldots,\chi_\cT X_N$ span $\chi_\cT(\Difftb^1/\Difftb^0)(M;\upbeta^*E)$, this gives $\WF_{\tbop;\bop}^{(s;k+1),(\alpha_+,\alpha_\cT)}(u)\cap\pa\cR_{\cT,\rm out}^\pm=\emptyset$. This completes the inductive step and thus finishes the proof.
\end{proof}

We can now strengthen Theorem~\ref{ThmWM}:

\begin{thm}[Forward solutions of $P$ with additional b-regularity]
\label{ThmWb}
  Let $\sfs$, $\alpha_{\!\scri}$, $\alpha_+$, $\alpha_\cT\in\R$, let $k\in\N_0$. Suppose that $\sfs-2$ and $\sfs+2$ are $P$-admissible orders for the weights $\alpha_{\!\scri},\alpha_+,\alpha_\cT$ (see Definition~\usref{DefWMOrders}). Define function spaces on $M$ using the metric density $|\dd g|$. Then there exists $C>0$ so that the following holds: for $f\in\dot H_{\eop,\tbop;\bop}^{(\sfs-1;k),(2\alpha_{\!\scri}+2,\alpha_++2,\alpha_\cT)}(\ol\Omega;\upbeta^*E)$, the unique forward solution $u$ of $P u=f$ satisfies $\dot H_{\eop,\tbop;\bop}^{(\sfs;k),(2\alpha_{\!\scri},\alpha_+,\alpha_\cT)}(\ol\Omega;\upbeta^*E)$ and the estimate
  \[
    \|u\|_{\dot H_{\eop,\tbop;\bop}^{(\sfs;k),(2\alpha_{\!\scri},\alpha_+,\alpha_\cT)}(\ol\Omega;\upbeta^*E)} \leq C\|f\|_{\dot H_{\eop,\tbop;\bop}^{(\sfs-1;k),(2\alpha_{\!\scri}+2,\alpha_++2,\alpha_\cT)}(\ol\Omega;\upbeta^*E)}.
  \]
\end{thm}
\begin{proof}
  Let $X\in\Diffb^1(M;\upbeta^*E)$ be an operator with scalar principal symbol equal to that of one of the following vector fields: an edge-3b-vector field on $M$, the vector field $t_*\pa_{t_*}$, a vector field $\chi_\scri\Omega$ where $\chi_\scri\in\CI(M)$ is equal to $1$ near $\scri^+$ and supported in a collar neighborhood of $\scri^+$ and $\Omega\in\cV(\Sph^{n-1})$. Note that any element of $\Diffb^1(M;\upbeta^*E)$ with support in $\ft_*\geq -1$ can be written as a linear combination of such operators $X$ with coefficients in $\CI(M)$. Moreover, by Lemma~\ref{LemmaMOMComm} and \cite[Lemma~5.6]{HintzVasyScrieb}, we have $[X,L]\in\Psi_{\eop,\tbop}^s(M)$ for any $L\in\Psi_{\eop,\tbop}^s(M)$ with Schwartz kernel supported in both factors in $\ft_*\geq -1$, similarly for $L\in\cA^{(-2\beta_{\!\scri},-\beta_+,-\beta_\cT)}\Psi_{\eop,\tbop}^s(M)$ and variable orders $s$.

  Theorem~\ref{ThmWM} is the case $k=0$. To prove the result for $k+1$ b-derivatives, so for $f\in\dot H_{\eop,\tbop;\bop}^{(\sfs-1;k+1),(2\alpha_{\!\scri}+2,\alpha_++2,\alpha_\cT)}\subset\dot H_{\eop,\tbop;\bop}^{(\sfs-1,k),(2\alpha_{\!\scri}+2,\alpha_++2,\alpha_\cT)}$ and with the inductive hypothesis giving $u\in\dot H_{\eop,\tbop;\bop}^{(\sfs;k),(2\alpha_{\!\scri}+2,\alpha_++2,\alpha_\cT)}$, note that with $X$ as above, we have
  \begin{align*}
    P(X u) = X f + [P,X]u &\in \dot H_{\eop,\tbop;\bop}^{(\sfs-1;k),(2\alpha_{\!\scri}+2,\alpha_++2,\alpha_\cT)} + \dot H_{\eop,\tbop;\bop}^{(\sfs-2;k),(2\alpha_{\!\scri}+2,\alpha_++2,\alpha_\cT)}, \\
      &\subset \dot H_{\eop,\tbop;\bop}^{(\sfs-2;k),(2\alpha_{\!\scri}+2,\alpha_++2,\alpha_\cT)}.
  \end{align*}
  For $k=0$, we apply Theorem~\ref{ThmWM} and use that $\sfs-1$ satisfies its assumptions; for $k\geq 1$, we use the inductive hypothesis. In both cases, we get $X u\in\dot H_{\eop,\tbop;\bop}^{(\sfs-1;k),(2\alpha_{\!\scri},\alpha_+,\alpha_\cT)}$. Since $X$ was arbitrary, this implies $u\in\dot H_{\eop,\tbop;\bop}^{(\sfs-1;k+1),(2\alpha_{\!\scri},\alpha_+,\alpha_\cT)}$. But since $u$ satisfies $P u\in\dot H_{\eop,\tbop;\bop}^{(\sfs-1;k+1),(2\alpha_{\!\scri},\alpha_+,\alpha_\cT)}$, we can apply $(\eop,\tbop)$-microlocal elliptic, real principal type, and radial point estimates: near $\cT^+$, these are Lemmas~\ref{LemmaWbEll}, \ref{LemmaWbProp}, and Proposition~\ref{PropWbOut} (and its analogue at $\cR_{\cT,\rm in}$), and near $\scri^+$, these are \cite[Propositions~5.14--5.16]{HintzVasyScrieb}, while near $(\iota^+)^\circ$ these are standard b-estimates. Therefore, we can recover the edge-3b-regularity order $\sfs$ for $u$, i.e.\ $u\in\dot H_{\eop,\tbop;\bop}^{(\sfs;k+1),(2\alpha_{\!\scri},\alpha_+,\alpha_\cT)}$. This completes the proof.
\end{proof}

We leave the extension to the setting considered in Corollary~\ref{CorWt} to the interested reader; the only additional ingredient is \cite[Corollary~6.6]{HintzVasyScrieb}, which provides the additional b-regularity near $I^0\cup\scri^+\setminus\iota^+$.

\begin{rmk}[Small non-decaying perturbations: b-regularity]
\label{RmkWbNonDecay}
  Continuing Remark~\ref{RmkWNonDecay}, one can extend Theorem~\ref{ThmWb} to operators $P+\eps R$ where $R$ is of class~\eqref{EqWNonDecayR} and $|\eps|<\eps_0$, with $\eps_0>0$ depending on the orders $\sfs,\alpha_{\!\scri},\alpha_+,\alpha_\cT,k$. Indeed, when $\eps$ is sufficiently small, the operator $\eps R$ produces small, and thus absorbable, error terms in the microlocal estimates used in the proof of Theorem~\ref{ThmWb}. Note carefully that b-regularity (more precisely, $k$ degrees of b-regularity, in addition to a large degree of edge-3b-regularity) of the coefficients of $R$ is now crucial to ensure that commutators of $R$ with operators $X$ as in the proof of Theorem~\ref{ThmWb} are bounded maps on weighted edge-3b;b-Sobolev spaces.
\end{rmk}

\begin{cor}[Forward solutions of $P$ with conormal forcing]
\label{CorWbCon}
  Define function spaces on $M$ using the volume density $|\dd g|$. Let $\alpha_{\!\scri},\alpha_+,\alpha_\cT\in\R$, and suppose that in the notation of~\eqref{EqWRubarp1}, Definition~\usref{DefGAW}\eqref{ItGAWEdgeN}, and Definition~\usref{DefGSOSpec}, we have
  \[
    \alpha_{\!\scri} < -\frac12 + \ubar p_1,\ \ 
    \alpha_+ < \min\Bigl(-\frac12+\alpha_{\!\scri},\ -\frac{n-1}{2}+\beta^+),\ \ 
    \alpha_\cT \in \Bigl(\alpha_++\frac{n}{2}-\beta^+,\alpha_++\frac{n}{2}-\beta^-\Bigr).
  \]
  Let $f\in\dot H_\bop^{\infty,(2\alpha_{\!\scri}+2,\alpha_++2,\alpha_\cT)}(\ol\Omega;\upbeta^*E)$.\footnote{We recall that the space $\Hbsupp^{\infty,(2\alpha_{\!\scri},\alpha_+,\alpha_\cT)}(\ol\Omega;\upbeta^*E)$ consists of all smooth sections $v$ of $E$ over $M^\circ=\R_{t_*}\times\R_x^n$ with support in $\ft_*\geq 0$ so that every derivative of $v$ along any finite number of vector fields in the set $Z$ in~\eqref{EqIZ} lies in $\rho_{\!\scri}^{\alpha_{\!\scri}}\rho_+^{\alpha_+}\rho_\cT^{\alpha_\cT}L^2(M,|\dd t_*\,\dd x|)$, where $\rho_{\!\scri},\rho_+$, and $\rho_\cT$ are as in~\eqref{EqIRhos}.} Then the forward solution $u$ of $P u=f$ satisfies $u\in\dot H_\bop^{\infty,(2\alpha_{\!\scri},\alpha_+,\alpha_\cT)}(\ol\Omega;\upbeta^*E)$.
\end{cor}
\begin{proof}
  This follows from Theorem~\ref{ThmWb} upon fixing any admissible variable order function $\sfs$ and taking $k\in\N_0$ to be arbitrary.
\end{proof}

\begin{cor}[Pointwise decay]
\label{CorWbPointwise}
  Let $\gamma_{\!\scri},\gamma_+,\gamma_\cT\in\R$, and suppose that in the notation of~\eqref{EqWRubarp1}, Definition~\usref{DefGAW}\eqref{ItGAWEdgeN}, and Definition~\usref{DefGSOSpec}, we have
  \[
    \gamma_{\!\scri} < \frac{n-1}{2}+\ubar p_1,\qquad
    \gamma_+ < \min(\gamma_{\!\scri},\ 1+\beta^+),\qquad
    \gamma_\cT \in (\gamma_+-\beta^+,\gamma_+-\beta^-).
  \]
  Let $f\in\dot\cA^{(2\gamma_{\!\scri}+2,\gamma_++2,\gamma_\cT)}(\ol\Omega;\upbeta^*E)$.\footnote{Here, $\dot\cA^{(2\gamma_{\!\scri},\gamma_+,\gamma_\cT)}(\ol\Omega;\upbeta^*E)$ is the space of all smooth sections $v$ of $E$ over $M^\circ$ with support in $\ft_*\geq 0$ so that any finite number of derivatives of $v$ along the vector fields in the set $Z$ in~\eqref{EqIZ} lies in $\rho_{\!\scri}^{\gamma_{\!\scri}}\rho_+^{\gamma_+}\rho_\cT^{\gamma_\cT}L^\infty(M;E)$, where $\rho_{\!\scri},\rho_+,\rho_\cT$ are boundary defining functions of $\scri^+,\iota^+,\cT^+\subset M_1$ as in~\eqref{EqIRhos}.} Then the forward solution $u$ of $P u=f$ satisfies $u\in\dot\cA^{(2(\gamma_{\!\scri}-\eps),\gamma_+-\eps,\gamma_\cT-\eps)}(\ol\Omega;\upbeta^*E)$ for all $\eps>0$.
\end{cor}
\begin{proof}
  Let $\alpha_{\!\scri}=\gamma_{\!\scri}-\frac{n}{2}$, $\alpha_+=\gamma_+-\frac{n+1}{2}$, and $\alpha_\cT=\gamma_\cT-\frac12$. If $\mu_\bop$ is a positive b-density on $M$ (such as $\mu_\bop=\la t_*\ra^{-1}\la x\ra^{-n}|\dd g|$), then
  \[
    \Hbsupp^{\infty,(2\alpha_{\!\scri},\alpha_+,\alpha_\cT)}(\ol\Omega,|\dd g|) = \Hbsupp^{\infty,(2\gamma_{\!\scri},\gamma_+,\gamma_\cT)}(\ol\Omega,\mu_\bop).
  \]
  Sobolev embedding moreover gives the inclusions
  \[
    \Hbsupp^{\infty,(2\gamma_{\!\scri},\gamma_+,\gamma_\cT)}(\ol\Omega,\mu_\bop) \subset \dot\cA^{(2\gamma_{\!\scri},\gamma_+,\gamma_\cT)}(\ol\Omega) \subset \Hbsupp^{\infty,(2(\gamma_{\!\scri}-\eps),\gamma_+-\eps,\gamma_\cT-\eps)}(\ol\Omega,\mu_\bop)
  \]
  for all $\eps>0$. An application of Corollary~\ref{CorWbCon} then concludes the proof.
\end{proof}

If the forcing term $f$ in Corollaries~\ref{CorWbCon} and \ref{CorWbPointwise} has more than $2(\frac{n-1}{2}+\ubar p_1)+2$ orders of decay (pointwise, or in $L^2$ relative to a b-density) at $\scri^+$, the solution $u$ typically does not have more than $2(\frac{n-1}{2}+\ubar p_1)$ orders of decay. For example, when the spectrum of $p_1$ in Definition~\ref{DefGAW} consists of a single real eigenvalue (which is thus equal to $\ubar p_1$), then $u$ has a $x_{\!\scri}^{2(\frac{n-1}{2}+\ubar p_1)}$ leading order term (i.e.\ \emph{radiation field}) at $\scri^+$, as follows by direct integration using~\eqref{EqGAWEdgeN}; see also \cite[\S\S1.1.1, 5.1]{HintzVasyMink4} and \cite[Proof of Theorem~3.9]{HintzPrice}.

\begin{rmk}[Comparison with the stationary case]
\label{RmkWbCompStat}
  In the case that $P$ is a \emph{stationary} wave type operator (thus $\ubar p_1=\ubar S$), we distinguish two cases, depending on whether the inequality $1+\beta^+\leq\frac{n-1}{2}+\ubar S$ holds or not.
  \begin{enumerate}
  \item If this inequality does hold, we may take $\gamma_{\!\scri},\gamma_+$, and $\gamma_\cT$ arbitrarily close to the upper bounds $\frac{n-1}{2}+\ubar S$, $1+\beta^+$, and $\beta^+-\beta^-+1$, respectively. Thus, Corollary~\ref{CorWbCon} is significantly more precise than Theorem~\ref{ThmStCo}, since the forcing may have non-trivial decay orders at $\scri^+$, $\iota^+$, and $\cT^+$, and the solution has precisely matching orders. (We recall, however, that Theorem~\ref{ThmStCo} was used in the course of the proof of Corollary~\ref{CorWbCon}, namely in the proof of Corollary~\ref{CorWMP0Inv}.) In light of Theorem~\ref{ThmAS}, we expect the decay rates at $\iota^+$ and $\cT^+$ to be sharp (up to the arbitrarily small loss) for generic $P_0$ and all perturbations $P-P_0$, even for Schwartz forcing. This can likely be proved by writing $P u=f$ as $P_0 u=f-(P-P_0)u$ and using a version of Theorem~\ref{ThmStCo} that applies to conormal inputs with definite decay rates; we do not address this problem in the present paper.
  \item If this inequality does not hold, then the strongest choices of weights allowed in Corollary~\ref{CorWbCon} are arbitrarily close to $\frac{n-1}{2}+\ubar S$ for $\gamma_{\!\scri},\gamma_+$, and $\frac{n-1}{2}+\ubar S-\beta^-$ for $\gamma_\cT$. Thus, the decay at $\iota^+$ and $\cT^+$ falls short of that in Theorem~\ref{ThmStCo} by the amount $(1+\beta^+)-(\frac{n-1}{2}+\ubar S)>0$. For now, we leave it as an open problem to determine, in the case that $f$ has strong decay (e.g.\ Schwartz), what the sharp decay of $u$ is. We mention here only work in progress by Luk--Oh \cite{LukOhTwoTails} which shows that Price's law for angular momenta $\geq 1$ of a massless scalar field on exact Schwarzschild spacetimes \cite{HintzPrice,AngelopoulosAretakisGajicRNPrice} needs to be weakened by a full order on a large class of non-stationary but asymptotically (as $t_*\to\infty$) Schwarzschild spacetimes.
  \end{enumerate}
\end{rmk}

\section{Examples}
\label{SE}

We conclude this paper by demonstrating how our main results---Theorems~\ref{ThmStCo} and \ref{ThmAS} for stationary wave type operators, and~\ref{ThmWM} and \ref{ThmWb} and their Corollaries~\ref{CorWbCon} and \ref{CorWbPointwise} for non-stationary wave type operators---apply in simple concrete settings. In~\S\ref{SsEW} we consider scalar wave operators (and mild perturbations thereof) associated with stationary and non-stationary asymptotically flat metrics. In~\S\ref{SsEV} we couple such operators to asymptotically inverse square potentials. A proof of the nonlinear stability of Minkowski space is expected to be possible by applying the results of the present paper to the linearization of a suitable gauge-fixed Einstein equation around asymptotically Minkowskian metrics; this will be pursued elsewhere.

\subsection{Wave operators and short range potentials}
\label{SsEW}

A natural class of examples of stationary wave type operators is the following:

\begin{prop}[Scalar wave operators of stationary and asymptotically flat metrics]
\label{PropES}
  Let $n\geq 3$, and let $g_0$ be a stationary and asymptotically flat metric (see Definition~\usref{DefGSG}) on $M_0=\ol{\R_{t_*}}\times X$, where $X=\ol{\R^n}$. Denote by $\Box_{g_0}$ the scalar wave operator.
  \begin{enumerate}
  \item\label{ItESMassless} Assume that mode stability at frequencies $\sigma\in\R\setminus\{0\}$ holds for $\Box_{g_0}$. Then $\Box_{g_0}$ is spectrally admissible with indicial gap $(\beta^-,\beta^+)=(0,n-2)$ in the sense of Definition~\usref{DefGSOSpec}.
  \item\label{ItESPotential} Let $\delta>0$ and $V_0\in\cA^{2+\delta}(X)$. If $P_0:=\Box_{g_0}+V_0$ satisfies mode stability (for all $0\neq\sigma\in\C$, $\Im\sigma\geq 0$) and has no zero energy resonances, then $P_0$ is spectrally admissible.
  \end{enumerate}
  In both cases, we thus conclude that the forward solution of $P_0 u=f$, where $f\in\sS(\R_{t_*,x}^{1+n})$ is supported in $t_*\geq 0$, satisfies
  \[
    |Z^J u|\lesssim\rho_{\!\scri}^{\frac{n-1}{2}-\eps}(\rho_+\rho_\cT)^{n-1-\eps}
  \]
  for all $\eps>0$ and multi-indices $J$. Here, $Z$ is the set of vector fields from~\eqref{EqIZ}, and $\rho_{\!\scri}=\frac{\la t_*\ra}{t_*+r}$ and $\rho_+\rho_\cT=\la t_*\ra^{-1}$ as in~\eqref{EqIRhos}.
\end{prop}

\begin{rmk}[Mode stability]
\label{RmkESModeStab}
  The mode stability on the real axis assumed in part~\eqref{ItESMassless} is standard when $g_0$ is ultra-static, i.e.\ $g_0=-\dd t^2+h$ where $h=h(x,\dd x)$ is equal to the Euclidean metric $\dd x^2$ modulo errors in $\cA^{1+\delta}(\ol{\R^n})$ (which implies that $g_0$ is stationary and asymptotically flat in the sense of Definition~\ref{DefGSG}), or more generally a warped product metric of the form $g_0=-a^2\,\dd t^2+h$ where $a\equiv 1\bmod\cA^{1+\delta}(\ol{\R^n})$. (Indeed, one can prove using a boundary pairing argument \cite[\S2.3]{MelroseGeometricScattering} that every conormal element in $\ker\wh{P_0}(\sigma)$, $\sigma\in\R\setminus\{0\}$, must be Schwartz; and then a unique continuation result at infinity such as \cite[Theorem~17.2.8]{HormanderAnalysisPDE3} finishes the proof. Since a more general unique continuation result which would cover general metrics $g_0$ does not appear to exist in the literature, we simply \emph{assume} mode stability in the present paper.) In this setting, the spectral admissibility hypothesis in part~\eqref{ItESPotential} is easy to check when $V_0$ is nonnegative and real-valued.
\end{rmk}

The basic example to which Proposition~\ref{PropES}\eqref{ItESMassless} applies is the wave operator on $(n+1)$-dimensional Minkowski space. We stress that the metric $\slg$ in Definition~\ref{DefGSG} does not need to be the standard metric $g_{\Sph^{n-1}}$ on $\Sph^{n-1}$; for example, we can allow for general cone angles at infinity, so $\slg=a^2 g_{\Sph^{n-1}}$, $a>0$.

\begin{rmk}[Tensor wave operator]
\label{RmkESTensor}
  More generally, for any $p,q\in\N_0$, the tensor wave operator
  \[
    \Box_{g_0} = -\tr_{g_0}\nabla^{g_0}\nabla^{g_0} \in \Diff^2(M_0^\circ; (T M_0^\circ)^{\otimes p} \otimes (T^* M_0^\circ)^{\otimes q})
  \]
  on $(M_0^\circ,g_0)=(\R_{t_*}\times\R^n_x,g_0)$ is a stationary wave type operator for the bundle $E=(\wt T^*M_0)^{\otimes p}\otimes(\wt T^*M_0)^{\otimes q}$. Given Proposition~\ref{PropES}, this follows from the fact that the principal symbol of $\Box_{g_0}$ is scalar and equal to the dual metric function of $g_0$, and hence one has~\eqref{EqGSOStruct} by principal symbol considerations; see also Remark~\ref{RmkGSOStruct}. Condition~\eqref{ItGSOSpectf} with $(\beta^-,\beta^+)=(0,n-2)$ follows by using the trivialization $\wt\Tsc{}^*_{\pa X}X\cong\pa X\times\R^n$ induced by the differentials of the standard coordinates on $\R^n$, in which the transition face normal operator of $\Box_{g_0}$ consists of $(p+q)n$ copies of the scalar one. Conditions~\eqref{ItGSOSpec}--\eqref{ItGSOSpec0} in Definition~\ref{DefGSOSpec} can easily be verified in the ultra-static case mentioned in Remark~\ref{RmkESModeStab}.
\end{rmk}

\begin{proof}[Proof of Proposition~\usref{PropES}]
  \pfstep{Structure of the wave operator.} We first show that $\Box_{g_0}$ is a stationary wave type operator in the sense of Definition~\ref{DefGSO}. We only need to verify the structure~\eqref{EqGSOStruct} of $\Box_{g_0}$ in a collar neighborhood $[0,1)_\rho\times\pa X$ of $\pa X$. To this end, we recall from Lemma~\ref{LemmaGSGVol} that
  \begin{equation}
  \label{EqESVolume}
    |\dd g_0|\in (1+\cA^\delta(X))r^{n-1}|\dd t_*\,\dd r\,\dd\slg|,
  \end{equation}
  and note that $\pa_r=-\rho^2\pa_\rho$. In view of the structure~\eqref{EqGSGMetric}--\eqref{EqGSGMetricCoeff} of $g_0^{-1}$, the terms of $\Box_{g_0}$ involving $g_0^{t_* r}=g_0^{-1}(\dd t_*,\dd r)$ are
  \begin{align*}
    &r^{-(n-1)}(1+\cA^\delta)\pa_{t_*} r^{n-1}(1+\cA^\delta)\pa_r + r^{-(n-1)}(1+\cA^\delta)\pa_r r^{n-1}(1+\cA^\delta)\pa_{t_*} \\
    &\qquad \in \rho\pa_{t_*} \bigl( -2(1+\cA^\delta)\rho\pa_\rho + (n-1) + \cA^\delta \bigr),
  \end{align*}
  where we use that $\pa_r\colon\cA^\delta\to\cA^{1+\delta}$. Moreover, the coefficient of the principal term $\pa_{t_*}\pa_r$ is equal to $2+\cA^{1+\delta}$ by~\eqref{EqGSGMetricCoeff}, which implies that the coefficient of $\rho\pa_\rho$ in parentheses here in fact lies in $-2(1+\cA^{1+\delta})$. We can thus absorb the term $\rho\cdot\cA^{1+\delta}(X)\rho\pa_\rho\in\cA^{2+\delta}\Diffb^1$ into $Q$ in~\eqref{EqGSOStruct}, and the final $\cA^\delta$ term in parentheses above into $S$. Mixed derivatives in $t_*$ and the spherical variables contribute to $Q$ as well, noting that the simple order of vanishing of elements of $\Vsc(X)=\rho\Vb(X)$ times the $\cA^{1+\delta}$ remainder of $g_0^{0 X}$ in~\eqref{EqGSGMetricCoeff} produces an element of $\cA^{2+\delta}\Diffb^1$ indeed. The $\pa_{t_*}^2$ coefficient is necessarily $-g_0^{0 0}\in\cA^{1+\delta}$.

  It remains to compute $\wh{\Box_{g_0}}(0)$; we claim that
  \begin{equation}
  \label{EqESP0}
    \wh{\Box_{g_0}}(0) = \rho^2 P_{(0)} + \tilde P,\qquad P_{(0)}=-(\rho\pa_\rho)^2+(n-2)\rho\pa_\rho + \slDelta,\quad \tilde P\in\cA^{2+\delta}\Diffb^2(X),
  \end{equation}
  where $\slDelta=\Delta_\slg$ is the nonnegative Laplacian on $(\Sph^{n-1},\slg)$. This follows from the expression for $g^{X X}$ in~\eqref{EqGSGMetricCoeff}; the leading order term gives rise to the Laplacian for the metric $g_X:=\dd r^2+r^2\slg$ (which in the case $\slg=g_{\Sph^{n-1}}$ is the Euclidean metric) on $\R^n$, and additional $\cA^{2+\delta}\Diffb^2$ error terms arise from the remainder term in~\eqref{EqGSGMetric} and the $\cA^\delta$-error term in~\eqref{EqESVolume}. This completes the verification of~\eqref{EqGSOStruct} with $S\in\cA^\delta(X)$ and $\wh{\Box_{g_0}}(0)$ given by~\eqref{EqESP0}.

  \pfstep{Spectral admissibility.} We now verify the spectral admissibility of $\Box_{g_0}$ (Definition~\ref{DefGSOSpec}). Consider first the case that $\slg$ is the standard metric, so $g_X$ is Euclidean.

  \pfsubstep{(i)}{Zero energy operator.} The b-analysis of $\Delta_{\R^n}$ is then standard, see \cite[Theorem~4.5]{GellRedmanHaberVasyFeynman} and \cite[Lemma~3.2]{CarronCoulhonHassellRiesz}. Concretely, restricted to degree $\ell\in\N_0$ spherical harmonics on $\Sph^{n-1}$ (with eigenvalue $\ell(\ell+n-2)$), the b-normal operator $\wh N(P_{(0)},\lambda)$ is multiplication by $-(i\lambda)^2+(n-2)i\lambda+\ell(\ell+n-2)$; thus,
  \begin{equation}
  \label{EqESSpecb}
    \specb(\rho^{-2}\wh{\Box_{g_0}}(0))= -\N_0 \cup (n-2+\N_0).
  \end{equation}
  By Lemma~\ref{LemmaStEst0}, the operator $\wh{\Box_{g_0}}(0)\colon\Hb^{s,-\frac{n}{2}+\beta}(X,|\dd g_0|)\to\Hb^{s-2,-\frac{n}{2}+\beta+2}(X,|\dd g_0|)$ is Fredholm for $\beta\in(0,n-2)$. Since any element of its nullspace necessarily lies in $\cA^{n-2}(X)$ and thus decays at infinity, the maximum principle implies that $\wh{\Box_{g_0}}(0)$ is injective; arguing similarly for $\wh{\Box_{g_0}}(0)^*$ implies the invertibility of $\wh{\Box_{g_0}}(0)$. Thus, the indicial gap of $\wh{\Box_{g_0}}(0)$ is $(\beta^-,\beta^+)=(0,n-2)$, as claimed.

  \pfsubstep{(ii)}{Transition face operators.} Condition~\eqref{ItGSOSpectf} in Definition~\ref{DefGSOSpec} is verified for $\theta=0,\pi$ (i.e.\ for $N_\tface^\pm(P)$) in \cite[Lemma~5.10]{HintzConicProp} (with $\sfZ=0$) for the operator obtained from $N_\tface^\pm(\Box_{g_0})$ by switching from $\hat\rho$ to $\hat r=\hat\rho^{-1}$ and conjugating by $e^{-i\hat r}$; see also Remark~\ref{RmkGSOSpecOut}. We shall similarly prove the injectivity of
  \begin{equation}
  \label{EqESNtf}
    N_\tface^\theta(\Box_{g_0}) = -2 i e^{i\theta}\hat r^{-1}\Bigl(\hat r\pa_{\hat r}+\frac{n-1}{2}\Bigr) + \Bigl( -\pa_{\hat r}^2 - \frac{n-1}{\hat r}\pa_{\hat r} + \hat r^{-2}\slDelta \Bigr)
  \end{equation}
  on $\cA^{(\alpha,-\beta)}(\tface)$, $\alpha\in\R$, $\beta\in(\beta^-,\beta^+)$, for $\theta\in(0,\pi)$: if $u$ solves $N_\tface^\theta(\Box_{g_0})u=0$, then the projection $v_\ell=v_\ell(\hat r)$ of $\exp(i e^{i\theta}\hat r)u$ in the angular variables onto the space of spherical harmonics of degree $\ell$ solves the Bessel ODE
  \begin{equation}
  \label{EqESBesselODE}
    -v_\ell'' - \frac{n-1}{\hat r}v_\ell' + \frac{\lambda_\ell}{\hat r^2}v_\ell - e^{2 i\theta} v_\ell = 0,
  \end{equation}
  where $\lambda_\ell=\ell(\ell+n-2)$. The general solution is a linear combination of
  \begin{equation}
  \label{EqESBesselODESol}
    \hat r^{-\frac{n-2}{2}}H^{(1)}_{\nu_\ell}(e^{i\theta}\hat r),\qquad
    \hat r^{-\frac{n-2}{2}}H^{(2)}_{\nu_\ell}(e^{i\theta}\hat r),
  \end{equation}
  where $\nu_\ell=\frac{n-2}{2}+\ell$. Now if $u\in\cA^{(\alpha,-\beta)}(\tface)$, then $v_\ell$ is exponentially decaying as $\hat r\to\infty$; since $H_{\nu_\ell}^{(j)}(e^{i\theta}\hat r)$ is exponentially decaying for $j=1$ and exponentially growing for $j=2$, this means that $v_\ell=c\hat r^{-\frac{n-2}{2}}H_{\nu_\ell}^{(1)}(e^{i\theta}\hat r)$ for some $c\in\C$. But since $|\hat r^{-\frac{n-2}{2}}H_{\nu_\ell}^{(1)}(e^{i\theta}\hat r)|\gtrsim\hat r^{-\frac{n-2}{2}-\nu_\ell}\gtrsim\hat r^{-\frac{n-2}{2}-\nu_0}=\hat r^{-(n-2)}$ does not lie in $\cA^{-\beta}([0,1)_{\hat r})$ for $\beta\in(0,n-2)$, we must have $c=0$ and thus $v_\ell=0$.

  In the case that $\slg$ is an arbitrary Riemannian metric, one merely needs to replace $\ell(\ell+n-2)$ above by the eigenvalues of $\Delta_\slg$. The boundary spectrum of $\rho^{-2}\wh{\Box_{g_0}}(0)$ still has $(0,n-2)$ as the indicial gap (corresponding to the eigenvalue $0$ of $\Delta_\slg$, with eigenspace spanned by constants). The proof of the invertibility, resp.\ injectivity of $N_\tface^\theta(P_0)$ for $\theta=0,\pi$, resp.\ $\theta\in[0,\pi]$, goes through, with minor modifications, in this generality as well.

  \pfsubstep{(iii)}{Mode stability for $\Im\sigma>0$.} Consider a conormal element $u\in\ker\wh{P_0}(\sigma)$. Let $\ft$ denote a smooth function on $M^\circ$ so that $\ft-t_*=F(x)$ is smooth, stationary, and equal to $|x|=r$ for large $r$, and so that $\dd\ft$ is spacelike. (Using Lemma~\ref{LemmaGSGTime}, this can be accomplished by gluing $t_*+r$ for large $r$ to $t_*+C$ in a compact spatial region for an appropriate constant $C$.) Then $P_0(e^{-i\sigma\ft}\tilde u)=0$ where $\tilde u(x)=e^{i\sigma(\ft-t_*)}u(x)=e^{i\sigma F(x)}u(x)$. Thus, $\tilde u$ decays exponentially (together with all derivatives) as $|x|\to\infty$. In particular, the energy of $U(\ft,x):=e^{-i\sigma\ft}\tilde u(x)$,
  \[
    E(s) := \int_{\ft^{-1}(s)} T[U](\pa_\ft,(\dd\ft)^\flat) \,\dd\mu,
  \]
  is finite; here the volume density $\dd\mu$ on $\R^3_x$ is defined by $|\dd g_0|=|\dd\ft\,\dd\mu|$, and $T[U](X,Y)=X(U)Y(U)-\half g_0(X,Y)|\nabla^{g_0}U|_{g_0}^2$ is the stress-energy-momentum tensor. Since $\pa_\ft=\pa_{t_*}$ is a timelike Killing vector field, $E(s)$ is coercive and independent of $s$; but since $e^{-i\sigma\ft}|_{\ft=s}$ grows exponentially as $s\to\infty$, we must have $\nabla^{g_0}\tilde u=0$, thus $\tilde u=0$ (since $\tilde u$ tends to $0$ at infinity), and finally $u\equiv 0$.

  \pfstep{Coupling with stationary potentials.} Lastly, the addition of a potential $V$ with better than inverse quadratic decay at $\pa X$ only requires notational modifications; the absence of zero energy resonances is now an \emph{assumption}, whereas the transition face normal operators are independent of $V$.

  \pfstep{Decay.} The final claim follows immediately from Theorem~\ref{ThmStCo}.
\end{proof}

The proof also shows that $S|_{\pa X}=0$ in the notation of Definition~\ref{DefGSO}, which together with~\eqref{EqESP0} implies that the threshold quantities from Definition~\ref{DefStEstThr} are $\ubar S=\ubar S_{\rm in}=0$.

The following result gives sufficient conditions for the applicability of our second main result (Theorem~\ref{ThmAS}) for stationary wave type operators.

\begin{prop}[Sharp decay in odd spacetime dimensions]
\label{PropESSharp}
  Let $g_0$ be as in Proposition~\usref{PropES} on an $(n+1)$-dimensional spacetime, where $n\geq 4$ is \emph{even}. Let $f\in\sS(\R_{t_*,x}^{1+n})$ have support in $t_*\geq 0$.
  \begin{enumerate}
  \item\label{EqESSharpWave} The forward solution $u$ of the scalar wave equation $\Box_{g_0}u=f$ satisfies
    \begin{equation}
    \label{EqESSharp1}
      \Bigl| u(t_*,x) - c\bigl(t_*(t_*+2\la x\ra)\bigr)^{-\frac{n-1}{2}} \Bigr| \leq C t_*^{-(n-1)-\eps}\left(\frac{t_*}{t_*+\la x\ra}\right)^{\frac{n-1}{2}-\eta}
    \end{equation}
    for some\footnote{A careful inspection of the proof shows that one take $\eps$ arbitrarily close to $1$ (due to the fact that the next indicial root of the zero energy operator after $n-2$ is $n-1=(n-2)+1$).} $\eps>0$ and all $\eta>0$ (with $C$ depending on $\eta$); here, the constant $c=c(f)\in\C$ is nonzero for generic\footnote{more precisely: $f$ does not lie in a positive codimension subspace of the space of Schwartz functions with support in $t_*\geq 0$} $f$.
  \item\label{ItESSharpPot} For $V_0\in\cA^{2+\delta}(X)$ for which $\Box_{g_0}+V_0$ is spectrally admissible, the forward solution $u$ of $(\Box_{g_0}+V_0)u=f$ satisfies
    \begin{equation}
    \label{EqESSharpPot}
      \Bigl| u(t_*,x) - c a_\cT(x)\bigl(t_*(t_*+2\la x\ra)\bigr)^{-\frac{n-1}{2}} \Bigr| \leq C t_*^{-(n-1)-\eps}\left(\frac{t_*}{t_*+\la x\ra}\right)^{\frac{n-1}{2}-\eta}
    \end{equation}
    for some $\eps>0$ and all $\eta>0$, where $a_\cT=a_\cT(x)$ is the unique stationary solution of $(\Box_{g_0}+V_0)a_\cT=0$ with $a_\cT(x)\to 1$ as $|x|\to\infty$; here, $c$ is nonzero for generic $f$.
  \end{enumerate}
\end{prop}

Thus, the $t_*^{-(n-1)+\eps}\la x\ra^\eps$ bounds (in $r<q t$, $q\in(0,1)$) on forward solutions $u$ of $P u=f$ with Schwartz forcing $f$ provided by Proposition~\ref{PropES} are sharp up to the arbitrarily small loss $\eps>0$. We remark that~\eqref{EqESSharp1} matches the well-known $t^{-(n-1)}$ decay rate of linear waves on Minkowski space with odd spacetime dimension $n+1$, see e.g.\ \cite[\S6.2]{HormanderNonlinearLectures}. For any \emph{odd} $n\geq 3$ on the other hand, the decay rates can be improved by at least $\min(\delta,1)$ powers of $t_*^{-1}$ by following the arguments in~\S\ref{SAS}.

\begin{proof}[Proof of Proposition~\usref{PropESSharp}]
  We only need to verify the assumptions of Theorem~\ref{ThmAS}. We begin with part~\eqref{EqESSharpWave}. The discussion leading to~\eqref{EqESSpecb} implies that the only elements of $\specb(\rho^{-2}\wh{\Box_{g_0}}(0))$ with real part $0$ or $n-2$ are $\lambda^-=0$ and $\lambda^+=n-2$, and they arise from the restriction of $N(\rho^{-2}\wh{\Box_{g_0}}(0))$ to spherically symmetric functions, i.e.\ to constants. On constants, $\wh N(\rho^{-2}\wh{\Box_{g_0}}(0),\lambda)$ is multiplication by
  \[
    p(\lambda)=-(i\lambda)^2+(n-2)i\lambda;
  \]
  the inverse of $p(\lambda)$ has a simple pole at $i\lambda=i\lambda^\pm$. In the notation of~\eqref{EqASAssmSimple}, we can take $v^\pm=1\in\CI(\pa X)$. According to~\eqref{EqASf10} (with $S|_{\pa X}=0$), we have $f_{1,0}=-i(n-3)$, and then~\eqref{EqASAssmfLead} gives
  \[
    f_{j+1,0} = -i\frac{n-3-2 j}{p(n-2-j)}f_{j,0},\qquad j=1,\ldots,n-3.
  \]
  Since $n-3$ is odd, all $f_{j+1,0}$ are nonzero, and in particular $f_{n-2,0}\neq 0$. Finally, since $\beta^+-\beta^-=n-2$ is an integer, the solutions of the $\tface^\pm$-model problems~\eqref{EqASNtf} take the form~\eqref{EqASuprimeLog}. Concretely, acting on spherically symmetric functions on $[0,\infty]_{\hat r}\times\Sph^{n-1}$, we have
  \begin{equation}
  \label{EqESSharpNtf}
    N_\tface^\pm(\Box_{g_0}) = \mp 2 i\hat r^{-1}\Bigl(\hat r\pa_{\hat r}+\frac{n-1}{2}\Bigr) + \hat r^{-2}\bigl( -(\hat r\pa_{\hat r})^2 - (n-2)\hat r\pa_{\hat r} \bigr)
  \end{equation}
  by~\eqref{EqESNtf} for $\theta=0,\pi$. Upon multiplying this by $\hat r^2$, this is a regular singular ordinary differential operator at $\hat r=0$, and a solution of $N_\tface^\pm(\Box_{g_0})u'_\pm=\hat r^{-2}f_{n-2,0}$ with $u'_\pm\in\cA^{0-}([0,1)_{\hat r})$ necessarily satisfies $u'_\pm=-\frac{f_{n-2,0}}{n-2}\log\hat r+\cA^{((0,0),1-)}([0,1)_{\hat r})$; therefore,~\eqref{EqASuprimeLog} holds with $u'_{\pm,0}=-\frac{f_{n-2,0}}{n-2}\neq 0$. This verifies condition~\eqref{EqASAssmtfLead}.

  The asymptotic profile at $\cT^+$ is constant since constants are large zero energy states with decay rate $\rho^{\lambda^-}=\rho^0$; that is, $\Box_{g_0}1=0$. The asymptotic profile at $\iota^+$ is the same as that of the wave equation on Minkowski space, and can thus easily be determined by considering the fundamental solution of the wave equation \cite[\S6.2]{HormanderNonlinearLectures}. The latter is a constant multiple of $(t^2-r^2)_+^{-\frac{n-1}{2}}=(t_*(t_*+2 r))_+^{-\frac{n-1}{2}}$, multiplication of which by $t_*^{n-1}$ gives $a_+=(\frac{v}{2+v})_+^{\frac{n-1}{2}}$ up to a multiplicative constant.

  Alternatively, one can find $a_+$ using the description given in~\eqref{EqASipProfile}--\eqref{EqASipProfile2}. Explicitly, denoting by $L=-2\pa_v(v\pa_v+\frac{n-1}{2})-(v\pa_v+n-1)^2+(n-2)(v\pa_v+n-1)$ the operator on the left in~\eqref{EqASipProfile}, one first finds $L u_\delta=f_{n-2,0}i^{n-2}\delta^{(n-2)}(v)$ for $u_\delta=\sum_{j=0}^{n-3} c_j\delta^{(j)}(v)$, where $c_{n-3}=\frac{f_{n-2,0}i^{n-2}}{2(n-2)}$ and $c_j=c_{j+1}\frac{(j+3-n)(j+1)}{2(j+1-\frac{n-1}{2})}$ for $j=n-4,\ldots,0$; in particular, $\la u_\delta,1\ra=c_0\neq 0$. On the other hand, letting $u_1=(v(2+v))_+^{-\frac{n-1}{2}}=0$, we have $L u_1=0$,\footnote{This follows from the fact that this holds in $v>0$ (and trivially in $v<0$), and thus $L u_1$ must equal a sum of differentiated $\delta$-distributions at $v=0$; however, $u_1$ is an asymptotic sum of terms with homogeneities $-\frac{n-1}{2}+j$, $j\in\N_0$, none of which are negative integers. Therefore, $L u_1$ must in fact vanish.} and one finds $\la u_1,1\ra=(-1)^{\frac{n-2}{2}}\frac{2^{n-4}(\frac{n-4}{2})!^2}{(n-3)!}\neq 0$.\footnote{More generally, for $p,q\in\N_0$ with $p+q\geq 1$ (the case $p=q=\frac{n-2}{2}$ being of interest here), one has
  \begin{align*}
    \la v_+^{-\frac12-p},(2+v)^{-\frac12-q}\ra&=\frac{2^{2 p}p!}{(2 p)!}\la v_+^{-\frac12},\pa_v^p(2+v)^{-\frac12-q}\ra \\
      &= \frac{(-1)^p(2(p+q))!p!q!}{(2 p)!(2 q)!(p+q)!}\la v_+^{-\frac12},(2+v)^{-\frac12-p-q}\ra = \frac{(-1)^p 2^{p+q}p!q!(p+q-1)!}{(2 p)!(2 q)!}.
  \end{align*}} Therefore, $a_{+,0}=u_\delta-c_0\la u_1,1\ra^{-1} u_1$ up to a multiplicative constant; and $a_+=v^{n-1}a_{+,0}$.

  For part~\eqref{ItESSharpPot}, the only modification is that the large zero energy state $a_\cT$ is no longer the constant $1$ (unless $V_0=0$). Regarding the relationship of the description of the asymptotic profile in~\eqref{EqESSharpPot} with the expression in~\eqref{EqAS}, define the function
  \[
    a(t_*,x) := t_*^{n-1} \cdot a_\cT(x)\bigl(t_*(t_*+2\la x\ra)\bigr)^{-\frac{n-1}{2}} = a_\cT(x) \Bigl(1+\frac{2\la x\ra}{t_*}\Bigr)^{-\frac{n-1}{2}}.
  \]
  Since $\cT^+$ is (the closure of) the set of limits of $(t_*,x)$, $x\in\R^n$, as $t_*\to\infty$, and noting that $\frac{\la x\ra}{t_*}=0$ at $\cT^+$, we have $a|_{\cT^+}=a_\cT$. Moreover, taking the limit $t_*\to\infty$ while keeping $X=\frac{x}{t_*}\neq 0$ fixed gives $a|_{\iota^+}(X)=(1+2|X|)^{-\frac{n-1}{2}}$, which upon writing $v=|X|^{-1}$ gives the $\iota^+$-profile $a_+=(\frac{v}{2+v})_+^{\frac{n-1}{2}}$ from before.
\end{proof}

In the non-stationary setting, we have:

\begin{prop}[Scalar wave operators of admissible asymptotically flat metrics]
\label{PropES2}
  Let $n\geq 3$, let $g_0$ be a stationary and asymptotically flat metric on $M_0$ (with $\dim M_0=n+1$), and let $g$ be an $(\ell_0,2\ell_{\!\scri},\ell_+,\ell_\cT)$-admissible asymptotically flat metric on $M$ relative to $g_0$ (see Definition~\usref{DefGAG}).
  \begin{enumerate}
  \item  If $\Box_{g_0}$ is spectrally admissible,\footnote{See Proposition~\ref{PropES}\eqref{ItESMassless}.} then the scalar wave operator $\Box_g$ is an admissible wave type operator in the sense of Definition~\usref{DefGAW}.
  \item Let $\delta>0$. Let $V_0\in\cA^{2+\delta}(X)$, and suppose $\Box_{g_0}+V_0$ is spectrally admissible. Let $\tilde V\in\cA^{(2+2\delta,2+\delta,\delta)}(M)$, and let $V=V_0+\tilde V$. Then $\Box_g+V$ is an admissible wave type operator.
  \end{enumerate}
  The main results of this paper for non-stationary wave type operators (Theorems~\usref{ThmWM} and \usref{ThmWb}, Corollaries~\usref{CorWbCon} and~\usref{CorWbPointwise}) thus apply to $P=\Box_g$, resp.\ $\Box_g+V$. Concretely, let $f$ be supported in $\ft_*\geq 0$, and let $u$ be the unique forward solution of $P u=f$. Let $\alpha_{\!\scri},\alpha_+,\alpha_\cT\in\R$, and suppose that
  \[
    \alpha_{\!\scri} < -\frac12,\qquad
    \alpha_+ < -\frac12 + \alpha_{\!\scri},\qquad
    \alpha_\cT \in \Bigl(\alpha_+-\frac{n}{2}+2,\alpha_++\frac{n}{2}\Bigr).
  \]
  Let $\ol\Omega=\ft_*^{-1}([0,\infty])\subset M$. Then, for differential orders $\sfs\in\CI(\Setb^*_{\ol\Omega}M)$ so that $\sfs-2$ and $\sfs+2$ are $P$-admissible for the weights $\alpha_{\!\scri},\alpha_+,\alpha_\cT$ (see Definition~\usref{DefWMOrders}), the following results hold.
  \begin{enumerate}
  \myitem{ItES2Var}{\rm(a)}{a}{\rm (Variable order edge-3b-regularity, and additional b-regularity.)}
  \[
    f\in\Hsupp_\etbop^{\sfs-1,(2\alpha_{\!\scri}+2,\alpha_++2,\alpha_\cT)}(\ol\Omega)=\rho_{\!\scri}^{\alpha_{\!\scri}+1}\rho_+^{\alpha_++2}\rho_\cT^{\alpha_\cT}\Hsupp_\etbop^{\sfs-1}(\ol\Omega) \implies u\in\Hsupp_\etbop^{\sfs,(2\alpha_{\!\scri},\alpha_+,\alpha_\cT)}(\ol\Omega);
  \]
  more generally, $f\in\Hsupp_{\etbop;\bop}^{(\sfs-1;k),(2\alpha_{\!\scri}+2,\alpha_++2,\alpha_\cT)}(\ol\Omega)$ implies $u\in\Hsupp_{\etbop;\bop}^{(\sfs;k),(2\alpha_{\!\scri},\alpha_+,\alpha_\cT)}(\ol\Omega)$ for all $k\in\N_0\cup\{\infty\}$.
  \myitem{ItES2Hb}{\rm(b)}{b}{\rm ($L^2$-control with infinite order b-regularity.)}
  \[
    f\in\Hbsupp^{\infty,(2\alpha_{\!\scri}+2,\alpha_++2,\alpha_\cT)}(\ol\Omega) \implies u\in\Hbsupp^{\infty,(2\alpha_{\!\scri},\alpha_+,\alpha_\cT)}(\ol\Omega).
  \]
  \myitem{ItES2Pw}{\rm(c)}{c}{\rm (Pointwise control with infinite order b-regularity.)} If $\gamma_+-(n-2)<\gamma_\cT<\gamma_+<\gamma_{\!\scri}<\frac{n-1}{2}$, then given $f\in\dot\cA^{(2\gamma_{\!\scri}+2,\gamma_++2,\gamma_\cT)}(\ol\Omega)=\rho_{\!\scri}^{\gamma_{\!\scri}+1}\rho_+^{\gamma_++2}\rho_\cT^{\gamma_\cT}\dot\cA^{(0,0,0)}(\ol\Omega)$, we have $u\in\dot\cA^{(2(\gamma_{\!\scri}-\eps),\gamma_+-\eps,\gamma_\cT-\eps)}(\ol\Omega)$ for all $\eps>0$.
  \end{enumerate}
\end{prop}

We refer the reader to Lemmas~\ref{LemmaGAGPert} and \ref{LemmaGAGPertMet} for a description of the class of metrics $g$ which Proposition~\ref{PropES2} can handle.

\begin{proof}[Proof of Proposition~\usref{PropES2}]
  By \cite[Example~3.13]{HintzVasyScrieb}, $P=\Box_g$ has the required form near $\scri^+$. The reference also gives the vanishing of $p_0$ and $p_1$; and $S|_{\pa X}=0$, so $\ubar S=\ubar S_{\rm in}=0$. The only remaining thing to check is that we have $\chi_\cT\Box_g\equiv\chi_\cT\Box_{g_0}\bmod\cA^{(2+\ell_+,\ell_\cT)}\Difftb^2(M)$ where $\chi_\cT\in\CI(M)$ is $1$ near $\cT^+$ and supported in a collar neighborhood of $\cT^+$; the orders $2+\ell_+$ and $\ell_\cT$ refer to decay orders at $\iota^+$ and $\cT^+$, respectively. This follows, similarly to \cite[Example~3.8]{HintzVasyScrieb}, from the following observation about the Levi-Civita connection of $g$: keeping only the weights at $\iota^+$ and $\cT^+$, the Koszul formula implies
  \[
    \nabla^g \in (\CI+\cA^{(1+\ell_+,\ell_\cT)})\Difftb^1(M;\Ttb M,\Ttb^*M\otimes\Ttb M),\qquad
    \nabla^g-\nabla^{g_0}\in\cA^{(1+\ell_+,\ell_\cT)}\Difftb^1,
  \]
  similarly for the connection acting on sections of tensor products of $\Ttb M$ and $\Ttb^*M$.

  The regularity and decay estimates for forward solutions now follow from Theorems~\ref{ThmWM} and \ref{ThmWb} and Corollaries~\ref{CorWbCon} and~\ref{CorWbPointwise}
\end{proof}

\begin{rmk}[Introductory theorem]
\label{RmkEWIntro}
  Parts~\eqref{ItIVSchwartz}--\eqref{ItIVLinfty} of Theorem~\ref{ThmIV} follow Proposition~\ref{PropES2}, and part~\eqref{ItIVStat} follows from Proposition~\ref{PropESSharp}.
\end{rmk}

\subsection{Waves coupled to inverse square potentials}
\label{SsEV}

Consider again a stationary and asymptotically flat metric $g_0$ on $M_0=\ol{\R_{t_*}}\times X$ with $X=\ol{\R^n}$; we now only assume that $n\geq 1$. Let $g$ be an admissible asymptotically flat metric on $M$ relative to $g_0$. We couple $\Box_{g_0}$ with a stationary potential $V_0$ on $X=\ol{\R^n}$ which near infinity is approximately of inverse square type. Concretely, consider
\[
  V_0 \in \rho^2\CI(X) + \cA^{2+\delta}(X),\qquad \delta>0.
\]
We shall assume that
\[
  \alpha:=\rho^{-2}V_0|_{\pa X}
\]
is a constant; we allow for $\alpha$ to be complex. We state the analogues of Propositions~\ref{PropES}, \ref{PropES2}, and \ref{PropESSharp} all at once:

\begin{prop}[Approximate inverse square potentials]
\label{PropEV}
  Suppose mode stability in the closed upper half plane holds for $P_0=\Box_{g_0}+V_0$ in the sense of Definition~\usref{DefGSOSpec}\eqref{ItGSOSpec}--\eqref{ItGSOSpec0}. Suppose that $\alpha\in\C\setminus(-\infty,-(\frac{n-2}{2})^2]$. Then $P_0$ is spectrally admissible with indicial gap
  \[
    (\beta^-,\beta^+)=(\Re\lambda^-,\Re\lambda^+),\qquad
    \lambda^\pm:=\frac{n-2}{2}\pm\nu_0,\quad \nu_0:=\sqrt{\Bigl(\frac{n-2}{2}\Bigr)^2+\alpha}.
  \]
  The square root here is the unique one with positive real part. Therefore, the solution of $P_0 u_0=f$, with $f$ Schwartz and supported in $t_*\geq 0$, satisfies the pointwise estimate $|Z^J u_0|\lesssim\rho_{\!\scri}^{\frac{n-1}{2}-\eps}\rho_+^{\frac{n}{2}+\Re\nu_0-\eps}\rho_\cT^{2\Re\nu_0+1-\eps}$ for all $\eps>0$ and multi-indices $J$. Moreover:
  \begin{enumerate}
  \item\label{ItEVSharp} Suppose that either $\nu_0\in\R\setminus(\half+\N_0)$, or $\nu_0\in\C\setminus\R$ and $\Re\nu_0\notin\half\N$. Then for generic $f$, we have pointwise lower bounds $|u_0|\gtrsim\rho_{\!\scri}^{\frac{n-1}{2}-\eps}\rho_+^{\frac{n}{2}+\Re\nu_0}\rho_\cT^{2\Re\nu_0+1}$. More precisely, $P_0$ satisfies the assumptions of Theorem~\usref{ThmAS}, and the asymptotic profiles $a_\cT$ and $a_+$ (see~\eqref{EqAS}) of the forward solution $u_0$ at $\cT^+$ and $\iota^+$ are as follows: $a_\cT=a_\cT(x)$ solves $P_0 a_\cT=0$ and asymptotes to a constant as $|x|\to\infty$; and $a_+=a_+(v,\omega)$ (with $v=\rho t_*=\frac{t_*}{r}>0$) is a constant multiple of $(\frac{v}{2+v})^{\frac12+\nu_0}$.
  \item\label{ItEVAdm} Let $\tilde V\in\cA^{(2+2\delta,2+\delta,\delta)}(M)$ and $V=V_0+\tilde V$.\footnote{See also Example~\ref{ExGAWPotential}.} Then the operator $P=\Box_g+V$ is an admissible wave type operator (Definition~\usref{DefGAW}). The conclusions~\eqref{ItES2Var}--\eqref{ItES2Pw} in Proposition~\usref{PropES2} hold with $\alpha_{\!\scri}<-\frac12$, $\alpha_+<-\frac12+\alpha_{\!\scri}$, and $\alpha_\cT\in(\alpha_++1-\Re\nu_0,\alpha_++1+\Re\nu_0)$ (for the $L^2$-bounds) and $\gamma_+<\gamma_{\!\scri}<\frac{n-1}{2}$, and $\gamma_\cT\in(\gamma_+-\frac{n-2}{2}-\Re\nu_0,\gamma_+-\frac{n-2}{2}+\Re\nu_0)$ (for the $L^\infty$-bounds).
  \end{enumerate}
\end{prop}

Allowing for $\alpha$ (thus $\nu_0$, and thus the potential $V_0$) to be complex-valued is only of limited interest in applications. We work in this generality here, however, to underline the fact that our analysis does not require any reality or symmetry conditions on the wave type operators. Checking the mode stability assumption in Proposition~\ref{PropEV} on the other hand is typically a delicate task, unless $g_0$ is ultra-static or a warped product metric as in Remark~\ref{RmkESModeStab} and the potential $V_0$ is real and nonnegative in which case it is easy to show, or if $V_0$ is a perturbation of such a potential. See also \cite[\S3.1]{GajicInverseSquare}.

For $\alpha\in\R$ (and thus $\nu_0>0$), the condition in part~\eqref{ItEVSharp} is equivalent to
\begin{equation}
\label{EqEValpha}
  \alpha\neq
  \begin{cases}
     \ell(\ell+n-2) \ \text{for some $\ell\in\N_0$}, & n\ \text{odd}, \\
    (\ell+\half)(\ell+n-\frac32),\ \text{for some $\ell\in\N_0$}, & n\ \text{even}.
  \end{cases}
\end{equation}
In the special case that $n=3$ and $\alpha\in\R$, and restricting to spherically symmetric $g$, this reproves \cite[Theorem~1.1]{GajicInverseSquare}, with the caveat that unlike \cite{GajicInverseSquare} (where $g$ is the Schwarzschild metric), we do not allow for $g$ to have trapping in the present paper. (Trapping has no bearing on the low energy resolvent behavior, however, and indeed the nontrapping assumption can be removed easily. See Remark~\ref{RmkITrapping} and also \cite[\S4]{HintzPrice}.)

We note that the decay estimates we obtain here for $P$ (Corollaries~\ref{CorWbCon} and~\ref{CorWbPointwise}) are strictly weaker than those for $P_0$, cf.\ Remark~\ref{RmkWbCompStat} (with $\ubar S=0$).

\begin{proof}[Proof of Proposition~\usref{PropEV}]
  The transition face normal operator of $P_0$ is
  \[
    N_\tface^\theta(P_0) = -2 i e^{i\theta}\hat r^{-1}\Bigl(\hat r\pa_{\hat r}+\frac{n-1}{2}\Bigr) + \Bigl( -\pa_{\hat r}^2 - \frac{n-1}{\hat r}\pa_{\hat r} + \hat r^{-2}\slDelta + \hat r^{-2}\alpha \Bigr).
  \]
  Its invertibility for $\theta=0,\pi$ is checked in \cite[Lemma~5.10]{HintzConicProp} (with $\sfZ=\alpha$, and noting that the proof goes through also for $n=1$), with the boundary spectrum of $\rho^{-2}\hat P_0(0)$ restricted to degree $\ell\in\N_0$ spherical harmonics given by\footnote{The definition of the boundary spectrum in the reference does not include the factor $i$ in~\eqref{EqMspecb}, and has an additional overall minus sign since it is defined with respect to $\hat r=\hat\rho^{-1}$ there, instead of $\hat\rho$ here.}
  \begin{equation}
  \label{EqEVnuell}
    \frac{n-2}{2} \pm \nu_\ell,\qquad \nu_\ell:=\sqrt{\Bigl(\frac{n-2}{2}+\ell\Bigr)^2 + \alpha}.
  \end{equation}
  Since for $\alpha\in\C\setminus(-\infty,-(\frac{n-2}{2})^2]$, the real part of $\nu_\ell$ is a strictly increasing function of $\ell\geq 0$, this completes the determination of the indicial gap. The injectivity of $N_\tface^\theta(P_0)$ for $\theta\in(0,\pi)$ follows by noting that if $u\in\cA^{(\alpha,-\beta)}(\tface)$, with $\beta$ in the indicial gap, lies in its kernel, then the projection $v_\ell(\hat r)$ of $\exp(i e^{i\theta}\hat r)u$ to degree $\ell\in\N_0$ spherical harmonics satisfies the Bessel ODE~\eqref{EqESBesselODE} with $\lambda_\ell=\ell(\ell+n-2)+\alpha$, and is thus a linear combination of the weighted Hankel functions~\eqref{EqESBesselODESol} with $\nu_\ell$ as in~\eqref{EqEVnuell}. As in the arguments after~\eqref{EqESBesselODESol}, the exponential decay of $v_\ell$ as $\hat r\to\infty$ selects the Hankel function of the first kind, and the $\hat r^{-\beta}$ upper bound on $v_\ell$ near $\hat r=0$ then forces $v_\ell=0$. Thus, $P_0$ is spectrally admissible in the sense of Definition~\ref{DefGSOSpec} (with $S|_{\pa X}=0$), and therefore $P$ is admissible in the sense of Definition~\ref{DefGAW} (with $p_0,p_1=0$, and with the threshold quantities $\ubar S$ and $\ubar S_{\rm in}$ in Definition~\ref{DefStEstThr} being $0$). This finishes the proof of part~\eqref{ItEVAdm}.

  Turning to part~\eqref{ItEVSharp}, note that the strict monotonicity of $\Re\nu_\ell$ in $\ell$ implies that assumption~\eqref{EqASAssmSimple} is satisfied with $v^\pm=1\in\CI(\pa X)$. According to~\eqref{EqASf10} (with $S|_{\pa X}=0$), we have $f_{1,0}=-i(-1+2\nu_0)$, and then~\eqref{EqASAssmfLead} gives
  \[
    f_{j+1,0} = -i\frac{-1+2\nu_0-2 j}{p\left(\frac{n-2}{2}+\nu_0-j\right)}f_{j,0},\qquad p(\lambda)=-\lambda^2+(n-2)\lambda+\alpha,
  \]
  for $j=1,\ldots,k-1$ where $k=\lceil\beta^+-\beta^-\rceil=\lceil 2\nu_0\rceil$. Thus, $f_{k,0}\neq 0$ provided $\nu_0\notin\half+\N_0$.

  We next need to study the solutions of the transition face model problems~\eqref{EqASNtf}. We directly consider~\eqref{EqASiplus} and note that the conjugation of $N_\tface^+(P_0)$ entering in~\eqref{EqASiplus}, restricted to spherically symmetric inputs, takes the form
  \[
    \hat r^{\lambda^+}N_\tface^+(P_0)\hat r^{-\lambda^+} = -2 i\hat r^{-1}\Bigl(\hat r\pa_{\hat r}-\lambda^++\frac{n-1}{2}\Bigr) + \hat r^{-2}\bigl( -(\hat r\pa_{\hat r}-\lambda^+)^2 - (n-2)(\hat r\pa_{\hat r}-\lambda^+) + \alpha \bigr)
  \]
  and $f'_\pm=\hat r^{-\lambda^++k-2}f_{k,0}$. Suppose first that $2\nu_0=\lambda^+-\lambda^-\in\N$ (which in particular forces $\alpha$ to be real). Then $u'_\pm$ necessarily has a nonvanishing logarithmic term $\hat r^{-\lambda^-}\log\hat r$ at $\hat r=0$ by a slight modification of the argument following~\eqref{EqESSharpNtf}. If $\nu_0$ is complex, we only consider the case that $\Re(2\nu_0)\notin\N$; and also for real $\nu_0$ we can assume that this is the case, since the case $2\nu_0\in\N$ was just discussed. Now for $\Re(2\nu_0)\notin\N$, we have $\lambda^+-\lambda^-\notin\N$. One then finds that the function $\tilde a_{+,0}(v,\omega):=(v(2+v))_+^{-\frac12-\nu_0}$ solves equation~\eqref{EqASipProfile} with right hand side $0$. Moreover, using $\Re\nu_0>0$, the pairing $\la\tilde a_{+,0},1\ra$ is well-defined and nonzero (as is easily seen upon writing $v_+^{-\frac12-\nu_0}$ as a nonzero constant times the $p$-th order derivative of $v_+^{-\frac12-\nu_0+p}$ where $p=\lfloor\Re\nu_0+\frac12\rfloor$). Then $a_{+,0}$ is the sum of $\tilde a_{+,0}$ and appropriate differentiated $\delta$-distributions at $0$ by arguments similar to those in the proof of Proposition~\ref{PropESSharp}. Since $a_{+,0}(v,\omega)$ therefore has a non-trivial $v^{-1-2\nu_0}$ leading order term as $v\to\infty$, this completes the verification of assumption~\eqref{EqASAssmtfLead}.
\end{proof}

\appendix

\section{Analysis of non-stationary ODEs}
\label{SODE}

We consider an ordinary differential operator
\[
  P(t,D_t) := \sum_{j=0}^m p_j(t)D_t^j,\qquad p_j\in\CI([0,\infty)),
\]
of order $m\geq 1$ on $[0,\infty)$ which is \emph{asymptotically stationary} in the following sense: the limits $p_{0,j}=\lim_{t\to\infty}p_j(t)\in\C$ exist, and there exists $\delta>0$ so that $p_j(t)=p_{0,j}+\cO(t^{-\delta})$ as $t\to\infty$ in the precise sense that
\[
  \tilde p_j := p_j - p_{0,j} \in \cA^\delta,
\]
by which we mean $|\la t\ra^k\pa_t^k\tilde p_j(t)|\leq C_{j k}\la t\ra^{-\delta}$ for all $k\in\N_0$ and $t\geq 0$. The \emph{stationary model} of $P$ is thus the constant coefficient operator
\[
  P_0(D_t) := \sum_{j=0}^m p_{0,j}D_t^j.
\]
We write
\begin{equation}
\label{EqODEDiff}
  \tilde P:=P-P_0\in\cA^\delta\Diff^m
\end{equation}
to mean that the coefficients of $\tilde P$ (the non-stationary perturbation) are of class $\cA^\delta$. We make the following assumptions.
\begin{enumerate}
\item{\rm (Principal symbol.)} $P$ is uniformly elliptic, i.e.\ $p_{0,m}\neq 0$ and $p_m(t)\neq 0$ for all $t\in[0,\infty)$. (Equivalently, $\inf_{t\in[0,\infty)}|p_m(t)|>0$.)
\item{\rm (Mode stability.)} Define the spectral family of $P_0$ by $\wh{P_0}(\sigma):=\sum_{j=0}^m p_{0,j}(-\sigma)^j$, $\sigma\in\C$. Then $\wh{P_0}(\sigma)$ is nonzero for $\Im\sigma\geq 0$.
\end{enumerate}

We write $H^{s,\alpha}(\R)=\la t\ra^{-\alpha}H^s(\R)$ for standard weighted Sobolev spaces, and
\begin{align*}
  \dot H^{s,\alpha}([0,\infty))&:=\{u\in H^{s,\alpha}(\R)\colon \supp u\subset[0,\infty) \}, \\
  \bar H^{s,\alpha}((0,\infty))&:=\{u|_{(0,\infty)} \colon u\in H^{s,\alpha}(\R) \}
\end{align*}
for the corresponding spaces of supported, resp.\ extendible distributions on the positive half line (cf.\ \cite[Appendix~B]{HormanderAnalysisPDE3}). The norms on these spaces are the restriction of, resp.\ the quotient norm induced by, the norm on $H^{s,\alpha}(\R)$. The following (elementary) result is the ODE analogue of Theorem~\ref{ThmWM}.

\begin{thm}[Forward solutions for $P$]
\label{ThmODE}
  Let $s,\alpha\in\R$. Then there exists $C>0$ so that the following holds: given $f\in\Hsupp^{s-m,\alpha}([0,\infty))$, the unique forward solution $u$ of $P u=f$ (i.e.\ $u=0$ for $t<0$) satisfies $u\in\Hsupp^{s,\alpha}([0,\infty))$ and the estimate
  \[
    \|u\|_{\Hsupp^{s,\alpha}([0,\infty))} \leq C\|f\|_{\Hsupp^{s-m,\alpha}([0,\infty))}.
  \]
\end{thm}
\begin{proof}
  We give a proof mirroring our approach to Theorem~\ref{ThmWM} as closely as possible (which is thus certainly more complicated than needed if one only wishes to solve the ODE). Unlike in the wave equation setting, there are no propagation phenomena here (since $P$ is elliptic); there do not exist analogues of null infinity or punctured future timelike infinity (since compactifying $\R$ at $+\infty$ merely adds a single point); and the low energy resolvent analysis is trivial (and in particular the resolvent is holomorphic in $\Im\sigma\geq 0$, with zero energy playing no special role---ultimately because $\wh{P_0}(\sigma)$ acts on a compact space, viz.\ a point) and thus there are no threshold conditions on weights or differentiability orders.

  \pfstep{Step 1. Decay of forward solutions for $P_0$.} (Cf.\ Theorem~\ref{ThmStCo}.) We first consider $f\in\dot\sS([0,\infty))$, i.e.\ Schwartz functions on $\R$ which vanish on $(-\infty,0)$. By the Paley--Wiener theorem, this is equivalent to $\hat f(\sigma)=\int_\R e^{i\sigma t}f(t)\,\dd t$ being holomorphic in $\Im\sigma>0$ and satisfying
  \begin{equation}
  \label{EqODEPW}
    |\hat f(\sigma)| \leq C_N\la\Re\sigma\ra^{-N},\qquad \sigma\in\C,\ \Im\sigma\geq 0,
  \end{equation}
  for all $N$. By ellipticity of $P_0$, we have $|\wh{P_0}(\sigma)|\geq c\la\sigma\ra^m$ for some $c>0$ when $|\sigma|$ is sufficiently large (this is thus a high energy estimate); moreover, $\wh{P_0}(\sigma)$ is nonzero for $\Im\sigma\geq 0$. Therefore, $\wh{u_0}(\sigma):=\wh{P_0}(\sigma)^{-1}\hat f(\sigma)$ satisfies the bounds~\eqref{EqODEPW} as well, and hence $u_0=\cF^{-1}\wh{u_0}\in\dot\sS([0,\infty))$ is the desired forward solution of $P_0 u=f$.\footnote{The Schwartz decay of $u_0$ should be contrasted with the definite decay rates of Theorems~\ref{ThmStCo} and \ref{ThmAS}, which are entirely due to the limited regularity of the low energy resolvent.}
  
  \pfstep{Step 2. Global regularity estimate for $P$ and $P^*$.} (Cf.\ Proposition~\ref{PropWR}.) Let $\alpha$, $s$, $N\in\R$. Then we have
  \begin{equation}
  \label{EqODEEll}
    \|u\|_{\dot H^{s,\alpha}([0,\infty))} \leq C\Bigl(\|P u\|_{\dot H^{s-m,\alpha}([0,\infty))} + \|u\|_{\dot H^{-N,\alpha}([0,\infty))}\Bigr)
  \end{equation}
  since $P$ is uniformly elliptic.\footnote{We caution the reader that the error term $\|u\|_{\dot H^{-N,\alpha}([0,\infty))}$ cannot be dropped without using spectral information on $P$ (as we are about to show), even though $P$ is evidently injective on spaces of distributions supported in $t\geq 0$. For example, the range of $P=D_t\colon\dot H^{s,\alpha}([0,\infty))\to H^{s-1,\alpha}([0,\infty))$ (which violates the mode stability hypothesis for $\sigma=0$) is not closed for any $s,\alpha\in\R$.}

 Let $\chi,\tilde\chi\in\CI([0,\infty))$ be identically $0$ on $[0,1]$ and $1$ on $[2,\infty)$, with $\tilde\chi=1$ on $\supp\chi$. For the adjoint, we then have\footnote{Elliptic estimates for partial differential operators require a priori control on a slightly larger region, hence the larger cutoff $\tilde\chi$ on the right hand side. In the present ODE setting, a sharper estimate is of course possible, but we proceed with the stated estimate for better comparison with the wave equation setting in the main part of the paper.}
  \begin{equation}
  \label{EqODEAdj}
    \|\chi\tilde u\|_{H^{-s+m,-\alpha}(\R)} \leq C\Bigl(\|\tilde\chi P^*\tilde u\|_{H^{-s,-\alpha}(\R)} + \|\tilde\chi\tilde u\|_{H^{-N,-\alpha}(\R)}\Bigr).
  \end{equation}

  \pfstep{Step 3. Fredholm estimate.} (Cf.\ Proposition~\ref{PropWMFred}.)

  \pfsubstep{(3.i)}{Direct estimate.} We have
  \[
    \|u\|_{\dot H^{-N,\alpha}([0,\infty))} \leq C\|(1-\chi)u\|_{\dot H^{-N,-N}([0,\infty))} + \|\chi u\|_{H^{-N,\alpha}(\R)}.
  \]
  Let $T=t^{-1}$ and $u':=T^{-\alpha}u$. Then $\|\chi u\|_{H^{-N,\alpha}(\R)}\sim\|\chi u'\|_{H^{-N}(\R)}$, and we estimate
  \[
    \|\chi u'\|_{H^{-N}(\R)} = \bigl\|\wh{\chi u'}(\sigma)\bigr\|_{\la\sigma\ra^N L^2(\R)} \leq C\bigl\|\wh{P_0}(\sigma)\wh{\chi u'}(\sigma)\bigr\|_{\la\sigma\ra^{N+m}L^2(\R)} = C\|P_0(\chi u')\|_{H^{-N-m}(\R)},
  \]
  where we used the mode stability and high energy estimates in the first inequality. This is then bounded by a constant times
  \begin{align*}
    &\|\chi P u'\|_{H^{-N-m}(\R)} + \|(P-P_0)\chi u'\|_{H^{-N-m}(\R)} + \|[P_0,\chi]u'\|_{H^{-N-m}(\R)} \\
    &\qquad \leq C\Bigl(\|\chi P u'\|_{H^{-N-m}(\R)} + \|\tilde\chi u'\|_{H^{-N,-\delta}(\R)} + \|\tilde\chi u'\|_{H^{-N-1,-N}(\R)}\Bigr),
  \end{align*}
  where we used~\eqref{EqODEDiff} to estimate the second term, and the compact support in $t$ of the coefficients of $[P_0,\chi]$ to estimate the third term. We now pass back to $u$ and use that $P u'=T^{-\alpha}P u+[P,T^{-\alpha}]u$; but $T^\alpha[P,T^{-\alpha}]\in\cA^1\Diff^{m-1}$ since $t^{-\alpha}[\pa_t,t^\alpha]=\alpha t^{-1}$ (which is a special case of the commutation properties of time dilation operators $t\pa_t$, cf.\ Lemma~\ref{LemmaM3bbComm}), and therefore
  \[
    \|\chi P u'\|_{H^{-N-m}(\R)} \leq C\Bigl(\|\chi P u\|_{H^{-N-m,\alpha}(\R)} + \|\tilde\chi u\|_{H^{-N-1,\alpha-1}(\R)}\Bigr).
  \]

  Altogether, we can now improve~\eqref{EqODEEll} to the estimate
  \begin{equation}
  \label{EqODESemiFred}
    \|u\|_{\dot H^{s,\alpha}([0,\infty))} \leq C\Bigl( \|P u\|_{\dot H^{s-m,\alpha}([0,\infty))} + \|u\|_{\dot H^{-N,\alpha-\eta}([0,\infty))}\Bigr)
  \end{equation}
  for $\eta=\min(\delta,1)>0$. Taking $-N<s$ here, this is a semi-Fredholm estimate (so $P$ has finite-dimensional kernel and closed range) since the inclusion $\dot H^{s,\alpha}([0,\infty))\hra\dot H^{-N,-\alpha-\eta}([0,\infty))$ is compact by the Rellich--Kondrakhov theorem. Since $P$ is injective on spaces of distributions vanishing for $t<0$, a standard functional analytic argument now implies that the second, error, term on the right in~\eqref{EqODESemiFred} can be dropped (cf.\ Corollary~\ref{CorWMInj}).

  \pfsubstep{(3.ii)}{Adjoint estimate.} Using similar arguments involving the stationary model $P_0$, we can improve~\eqref{EqODEAdj} to
  \begin{equation}
  \label{EqODESemiFredAdj0}
    \|\chi\tilde u\|_{H^{-s+m,-\alpha}(\R)} \leq C\Bigl(\|\tilde\chi P^*\tilde u\|_{H^{-s,-\alpha}(\R)} + \|\tilde\chi\tilde u\|_{H^{-N,-\alpha-\eta}(\R)}\Bigr)
  \end{equation}
  for the adjoint of $P$. Combining this with local-in-time theory near $t=0$, this implies the global estimate
  \begin{equation}
  \label{EqODESemiFredAdj}
    \|\tilde u\|_{\bar H^{-s+m,-\alpha}((0,\infty))} \leq C\Bigl(\|P^*\tilde u\|_{\bar H^{-s,-\alpha}((0,\infty))} + \|\tilde u\|_{\bar H^{-N,-\alpha-\eta}((0,\infty))}\Bigr).
  \end{equation}
  In more detail, let $\psi\in\CIc([0,3))$ be identically $1$ on $[0,2]$; then
  \[
    \|\psi\tilde u\|_{H^{-s+m}((0,3])^{-,\bullet}} \leq \|P^*(\psi\tilde u)\|_{H^{-s}((0,3])^{-,\bullet}},
  \]
  where $H^r((0,3])^{-,\bullet}=\{u|_{(0,\infty)}\colon u\in H^r(\R),\ \supp u\subset(-\infty,3]\}$. This follows from local-in-time solvability (in the direction of decreasing $t$) for $P^*$. But then $P^*(\psi\tilde u)=\psi P^*\tilde u+[P^*,\psi]\tilde u$, and $\|[P^*,\psi]\tilde u\|_{H^{-s}}\leq C\|\chi\tilde u\|_{H^{-s+m-1,\alpha}(\R)}$ since the coefficients of $[P^*,\psi]$ are supported in the compact set $\supp\dd\psi\subset\chi^{-1}(1)$. In combination with~\eqref{EqODESemiFredAdj0}, this gives~\eqref{EqODESemiFredAdj}.

  The estimates~\eqref{EqODESemiFred} and \eqref{EqODESemiFredAdj} imply that $P\colon\dot H^{s,\alpha}([0,\infty))\to\dot H^{s-m,\alpha}([0,\infty))$ is Fredholm (and in fact injective). This of course also applies to $P_0$.

  \pfstep{Step 4. Invertibility of the stationary model.} (Cf.\ Corollary~\ref{CorWMP0Inv}.) The range of
  \[
    P_0\colon\dot H^{s,\alpha}([0,\infty))\to\dot H^{s-m,\alpha}([0,\infty))
  \]
  is closed. But it contains the dense subspace $\dot\sS([0,\infty))$ by Step 1; therefore, $P_0$ is surjective, and in view of its injectivity bijective. As a consequence, we have a quantitative injectivity estimate for its adjoint,
  \begin{equation}
  \label{EqODEAdjInj}
    \|\tilde u\|_{\bar H^{-s+m,-\alpha}((0,\infty))} \leq C'\|\tilde P_0^*\tilde u\|_{\bar H^{-s,-\alpha}((0,\infty))}.
  \end{equation}

  \pfstep{Step 5. Proof of the Theorem.} (Cf.\ the proof of Theorem~\ref{ThmWM}.) For $\digamma>1$, and fixing $\phi\in\CIc([0,2))$ to be identically $1$ on $[0,1]$, we define
  \[
    P_\digamma := \phi(\digamma T)P + (1-\phi(\digamma T))P_0
  \]
  where $T=t^{-1}$. Thus $P_\digamma$ is uniformly bounded in $\cA^0\Diff^m$, and therefore~\eqref{EqODESemiFredAdj} applies to $P_\digamma^*$, with the constant being uniform for all $\digamma>1$. Using~\eqref{EqODEAdjInj}, we can estimate the error term of~\eqref{EqODESemiFredAdj} by
  \begin{align*}
    C\|\tilde u\|_{\bar H^{-N,-\alpha-\eta}((0,\infty))} &\leq C C'\|\tilde P_0^*\tilde u\|_{\bar H^{-N-m,-\alpha-\eta}((0,\infty))} \\
      &\leq C C'\Bigl(\|\tilde P_\digamma^*\tilde u\|_{\bar H^{-N-m,-\alpha-\eta}((0,\infty))} + \|(\tilde P_\digamma-\tilde P_0)^*\tilde u\|_{\bar H^{-N-m,-\alpha-\eta}((0,\infty))}\Bigr).
  \end{align*}
  But $\tilde P_\digamma-\tilde P_0\to 0$ in $\cA^{-\eps}\Diff^m$ as $\digamma\to\infty$ for all $\eps>0$; taking $\eps=\eta$, this implies
  \[
    C C'\|(\tilde P_\digamma-\tilde P_0)^*\tilde u\|_{\bar H^{-N-m,-\alpha-\eta}((0,\infty))} \leq \frac12\|\tilde u\|_{\bar H^{-N,-\alpha}((0,\infty))}
  \]
  for all sufficiently large $\digamma>1$. Altogether, we obtain
  \[
    \|\tilde u\|_{\bar H^{-s+m,-\alpha}((0,\infty))} \leq C\|P_\digamma^*\tilde u\|_{\bar H^{-s,-\alpha}((0,\infty))}.
  \]
  This implies the surjectivity, and thus invertibility, of $P_\digamma\colon\dot H^{s,\alpha}([0,\infty))\to H^{s-m,\alpha}([0,\infty))$.

  In order to deduce the invertibility of $P$, note that $P_\digamma=P$ for $t\geq\digamma$. We can uniquely solve $P u_{\rm in}=f\in\dot H^{s-m,\alpha}([0,\infty))$ on $[0,2\digamma)$ using local-in-time theory with $u_{\rm in}\in H^s([0,2\digamma))^{\bullet,-}$. We claim that we can solve $P u=f$ globally by writing
  \[
    u=\chi_{\rm in}u_{\rm in}+u'
  \]
  where $\chi_{\rm in}\in\CIc([0,2\digamma))$ equals $1$ on $[0,\digamma]$, and $u'\in\dot H^{s,\alpha}([0,\infty))$ is the unique solution of
  \[
    P_\digamma u'=f':=(1-\chi_{\rm in})f-[P,\chi_{\rm in}]u_{\rm in} \in \dot H^{s-m,\alpha}([0,\infty)).
  \]
  Indeed, $u'$ is supported in $[\digamma,\infty)$ (since this is true for $f'$), and therefore $P u'=P_\digamma u'=f'$. This implies $P(\chi_{\rm in}u_{\rm in}+u')=\chi_{\rm in}f+[P,\chi_{\rm in}]u_{\rm in}+P u'=f$, as desired.
\end{proof}

Similarly to Theorem~\ref{ThmWb}, we can extend Theorem~\ref{ThmODE} to spaces encoding higher degrees of b-regularity:

\begin{thm}[Forward solutions of $P$ with additional b-regularity]
\label{ThmODEb}
  Write $\dot H_{;\bop}^{(s;k),\alpha}([0,\infty))$ for the space of all $u\in\dot H^{s,\alpha}([0,\infty))$ so that $(\la t\ra\pa_t)^j u\in\dot H^{s,\alpha}([0,\infty))$ for all $j=0,\ldots,k$. Given $f\in\dot H_{;\bop}^{(s-m;k),\alpha}([0,\infty))$, the unique forward solution $u$ of $P u=f$ satisfies $u\in\dot H_{;\bop}^{(s;k),\alpha}([0,\infty))$, with norm bounded by a constant times the norm of $f$.
\end{thm}

Pointwise bounds analogous to Corollary~\ref{CorWbPointwise} follow from this via Sobolev embedding.

\begin{proof}[Proof of Theorem~\usref{ThmODEb}]
  Letting $X=\la t\ra\pa_t$, we note that
  \[
    P(X u)=X f+[P,X]u,
  \]
  where $[P,X]\in\cA^0\Diff^m$. (This uses the conormality, i.e.\ iterated regularity under application of $X$, of the coefficients of $P$, and moreover that $[\pa_t,\la t\ra\pa_t]=(\pa_t\la t\ra)\pa_t$ has bounded coefficients, despite the commutant $\la t\ra\pa_t$ having a strong weight at infinity; cf.\ Lemma~\ref{LemmaM3bbComm}.) If we have already proved the Theorem for $k-1$ degrees of b-regularity, then we already know $u\in\dot H_{;\bop}^{(s;k-1),\alpha}([0,\infty))$ when $f\in\dot H_{;\bop}^{(s-m;k),\alpha}([0,\infty))$. Thus, $X f+[P,X]u\in\dot H_{;\bop}^{(s-m;k-1),\alpha}([0,\infty))$. Applying the inductive hypothesis again $X u\in\dot H_{;\bop}^{(s;k-1),\alpha}([0,\infty))$. This completes the inductive step and finishes the proof of the Theorem.
\end{proof}

\section{Near-optimal inner products}
\label{SLA}

The goal is to prove the following statement:

\begin{prop}[Inner products]
\label{PropLA}
  Let $X$ be a smooth manifold, and let $E\to X$ be a complex vector bundle of finite rank. Let $S\in\CI(X;\End(E))$. For each $x\in X$, let
  \[
    \bar S(x):=\max\{\Re\mu\colon \mu\in\spec S(x)\},\quad
    \ubar S(x):=\min\{\Re\mu\colon \mu\in\spec S(x)\}.
  \]
  Then for all $\eps>0$, there exists a smooth positive definite Hermitian fiber inner product on $E$ with respect to which
  \begin{equation}
  \label{EqLA}
    \ubar S(x) - \eps < \half(S(x) + S(x)^*) < \bar S(x) + \eps\quad\forall\,x\in X.
  \end{equation}
  That is, for $e\in E_x$ with $\|e\|=1$, we have $\Re\la S(x)e,e\ra\in(\ubar S(x)-\eps,\bar S(x)+\eps)$.
\end{prop}

Note that
\begin{equation}
\label{EqLATriv}
  \min\spec\Bigl(\frac12(S(x)+S(x)^*)\Bigr)\leq \min\Re\spec S(x):=\min_{\mu\in\spec S(x)}\Re\mu,
\end{equation}
with strict inequality in general; this follows by applying the operator on the left to an eigenvector of $S(x)$ corresponding to the eigenvalue with smallest real part and taking the inner product with that eigenvector.

We first prove the following result:

\begin{lemma}[Inner products on $\C^N$]
\label{LemmaLAPw}
  Let $S\in\C^{N\times N}$, and let $\ubar S$, resp.\ $\bar S$ denote the minimum, resp.\ maximum of the real parts of the eigenvalues of $S$. Let $\eps>0$, and let $\cB(\eps)$ denote the set of Hermitian inner products on $\C^N$ with respect to which the estimate $\ubar S-\eps<\half(S+S^*)<\bar S+\eps$ holds. Then $\cB(\eps)$ is a nonempty open convex cone.
\end{lemma}
\begin{proof}
  Denote by $\la-,-\ra$ the standard Hermitian inner product on $\C^N$. Any positive definite inner product on $\C^N$ can be written in the form $(u,v)\mapsto\la u,B v\ra$ for some positive definite linear map $B\colon\C^N\to\C^N$. Such an inner product is an element of $\cB(\eps)$ if and only if
  \begin{equation}
  \label{EqLAPwIneq}
    2(\ubar S-\eps)\la u,B u\ra < \la S u,B u\ra + \la u,B S u\ra < 2(\bar S+\eps)\la u,B u\ra,\qquad u\in\C^N,\ u\neq 0.
  \end{equation}
  It is then clear that $\cB(\eps)$ is a convex cone. Moreover, since this chain of inequalities is homogeneous with respect to dilations of $u$, it suffices to check it for $u$ lying in the compact set $\{u\in\C^N\colon\la u,u\ra=1\}$; thus $\cB(\eps)$ is indeed open.

  It remains to prove that $\cB(\eps)$ is nonempty. This is discussed in \cite[\S3.4]{HintzPsdoInner}, but we give a detailed argument here for completeness. Pass to the Jordan normal form of $S$. It suffices to consider the case that $S$ only has a single Jordan block, as in the case of several Jordan blocks one may declare the basis vectors corresponding to different blocks to be orthogonal. Upon subtracting $\ubar S=\bar S$ times the identity from $S$, we thus have
  \[
    S = \begin{pmatrix} 0 & 1 & \dots & 0 \\ \vdots & \vdots & \ddots & \vdots \\ 0 & 0 & \dots & 1 \\ 0 & 0 & \dots & 0 \end{pmatrix}
  \]
  and $\bar S=\ubar S=0$. We then take $B=\diag(\delta^{N-1},\delta^{N-2},\ldots,1)$. For $u=(u_1,\ldots,u_N)$, the second inequality in~\eqref{EqLAPwIneq} holds for $\delta<\eps^2/4$ in view of
  \[
    2\sum_{j=1}^{N-1} \delta^{N-j}\Re(\ol{u_j}u_{j+1}) \leq \frac{\eps}{2}\sum_{j=1}^{N-1} \delta^{N-j}|u_j|^2 + \frac{2\delta}{\eps}\sum_{j=2}^N \delta^{N-j}|u_j|^2 \leq \eps\sum_{j=1}^N \delta^{N-j}|u_j|^2.
  \]
  The first inequality is proved similarly.
\end{proof}

\begin{proof}[Proof of Proposition~\usref{PropLA}]
  Fix any smooth positive definite Hermitian fiber inner product $\la-,-\ra$ on $E$. For each $x_0\in X$ then, we can choose a positive definite $B'_{x_0}\in\End(E_{x_0})$ so that~\eqref{EqLA} holds for the inner product $(e,f)\mapsto\la e,B'_{x_0}f\ra$ at $x_0$. Choose any smooth extension $B_{x_0}\in\CI(\cU_{x_0};E)$ of $B'_{x_0}$ to an open neighborhood $\cU_{x_0}$ of $x_0$. Then for a sufficiently small open set $\cV_{x_0}\subset\cU_{x_0}$ containing $x_0$, the estimate~\eqref{EqLA} continues to hold for all $x\in\cV_{x_0}$. The proof is then completed by using a partition of unity on $X$ subordinate to the cover $\{\cV_{x_0}\colon x_0\in X\}$, in view of the convexity statement of Lemma~\ref{LemmaLAPw}.
\end{proof}

\section{Construction of variable 3b-differential order functions}
\label{SWROrder}

\begin{proof}[Proof of Lemma~\usref{LemmaWROrder}]
  We begin by constructing local auxiliary functions with monotonicity properties along the $\pm\sfH$-flow.
  
  \pfstep{Constructions near $\scri^+$.} We first consider a neighborhood
  \[
    \Omega_\scri = \{ \ft_*\geq -1,\ x_{\!\scri}\leq\eps_\scri \}
  \]
  of $\scri^+$, with $\eps_\scri>0$ to be chosen; here we use the coordinates $x_{\!\scri}=\sqrt{(t_*+2)\rho}$, $\rho=r^{-1}$, and $\rho_+=(t_*+2)^{-1}$ from~\eqref{EqGAGebp} (with $T=-2$). Let $w'=x_{\!\scri}+x_{\!\scri}^{2\ell_{\!\scri}}\rho_+^{\ell_+}$. On $\Omega_\scri\cap\Sigma^\pm$, we have
  \begin{equation}
  \label{EqWROrderScriBd}
    \pm\zeta_\ebop\geq\frac13|\xi_\ebop|+|\eta_\ebop|
  \end{equation}
  when $\eps_\scri$ is sufficiently small; indeed, the expression~\eqref{EqWHebG} implies $\zeta_\ebop^2\geq(1-C w')(\xi_\ebop-\zeta_\ebop)^2+(2-C w')|\eta_\ebop|^2-C w'\zeta_\ebop^2$ on $\Sigma^\pm$, so $|\zeta_\ebop|\geq|\eta_\ebop|$ for small $w'$; and if $|\zeta_\ebop|<\half|\xi_\ebop|$, then still $(1-C w')(\xi_\ebop-\zeta_\ebop)^2\geq\frac14(1-C w')\xi_\ebop^2$, so $\zeta_\ebop^2\geq\frac18\xi_\ebop^2$ (and thus $|\zeta_\ebop|\geq\frac{1}{2\sqrt{2}}|\xi_\ebop|$) for small $w'$. This proves~\eqref{EqWROrderScriBd} since $\pm\zeta_\ebop>0$ on the characteristic set over $\scri^+$.
  
  On $\Omega_\scri\cap\Sigma^\pm$ then, we measure proximity to $\cR_{\scri,\rm out}^\pm$ (see~\eqref{EqWRadT}) by means of the (nonnegative) function
  \begin{equation}
  \label{EqWROrderScrif}
    f_{\scri,\rm out} := \frac{\xi_\ebop}{\zeta_\ebop} =: \hat\xi_\ebop.
  \end{equation}
  Using~\eqref{EqWHameb}, we compute, with $\hat\eta_\ebop=\frac{\eta_\ebop}{\zeta_\ebop}$,
  \begin{equation}
  \label{EqWRLocebComp}
    \pm|\zeta_\ebop|^{-1}H_{G_\ebop^+}f_{\scri,\rm out} = \zeta_\ebop^{-1}H_{G_\ebop^+}f_{\scri,\rm out} \leq -4|\hat\eta_\ebop|^2 + C w'.
  \end{equation}
  This is therefore negative on $\Omega_\scri\cap\{|\hat\eta_\ebop|^2\geq\frac{1}{32}\}$ provided $\eps_\scri>0$ is sufficiently small. Upon shrinking $\eps_\scri$ further, we can moreover ensure that points in $\Sigma^\pm\cap\Omega_\scri$ with $|\hat\xi_\ebop|\leq\frac12$ satisfy
  \begin{equation}
  \label{EqWRLoceb}
    \hat\xi_\ebop\geq-\frac18,\qquad
    |\hat\eta_\ebop|^2 \geq \frac12|\hat\xi_\ebop|-\frac{1}{16},\qquad
    |\hat\eta_\ebop|^2\leq 2|\hat\xi_\ebop|+\frac{1}{2} \leq \frac32;
  \end{equation}
  this is possible since by~\eqref{EqWHebG} we have $|\hat\eta_\ebop|^2=\half\hat\xi_\ebop(2-\hat\xi_\ebop)+\cO(w'(1+|\hat\eta_\ebop|^2))$ on the characteristic set over $\Omega_\scri$. Define then on $\Omega_\scri$ the function
  \begin{equation}
  \label{EqWRchiscriout}
    \chi_{\scri,\rm out} = \phi(f_{\scri,\rm out})\psi(|\hat\eta_\ebop|^2)
  \end{equation}
  where $\phi\in\CIc((-\frac12,\frac12))$ is identically $1$ on $[-\frac14,\frac14]$ and satisfies $\phi'\leq 0$ for positive arguments, and $\psi\in\CIc((-3,3))$ is identically $1$ on $[2,2]$. By~\eqref{EqWRLoceb}, we have $\supp(\phi\,\dd\psi)\cap\Omega_\scri\cap\Sigma^\pm=\emptyset$. Moreover, $\chi_{\scri,\rm out}$ is equal to $0$, resp.\ $1$ near $\pa\cR_{\scri,\rm in,+}^\pm$, resp.\ $\pa\cR_{\scri,\rm out}^\pm$. Furthermore, on $\Omega_\scri\cap\Sigma^\pm$, we have
  \[
    \pm|\zeta_\ebop|^{-1}H_{G_\ebop^+}\chi_{\scri,\rm out} = \phi'(f_{\scri,\rm out})\cdot(\pm|\zeta_\ebop|^{-1}H_{G_\ebop^+}f_{\scri,\rm out}) \geq 0,
  \]
  since on $\supp\dd\phi$ we have $|\hat\xi_\ebop|\in[\frac14,\frac12]$, which on $\Sigma^\pm$ forces $\hat\xi_\ebop\in[\frac14,\frac12]$ and thus $|\hat\eta_\ebop|^2\geq\frac{1}{16}$ by~\eqref{EqWRLoceb}, and hence~\eqref{EqWRLocebComp} and the line following it give the desired conclusion.
  
  In order to effect localization to a neighborhood of $\scri^+\cup\iota^+$, we employ the (nonnegative) monotone inverse time function
  \begin{equation}
  \label{EqWRft}
    f_\ft := \frac{1}{t_*+r} = (\rho_+^{-1}+\rho_+^{-1}x_{\!\scri}^{-2}+T)^{-1}
  \end{equation}
  with $T=-2$. In $\ft_*\geq -1$, we then find from~\eqref{EqWHameb}
  \[
    \pm|\zeta_\ebop|^{-1}H_{G_\ebop^+}f_\ft = -f_\ft^{-2} \zeta_\ebop^{-1}H_{G_\ebop^+}f_\ft^{-1} \leq -2 f_\ft^{-2} (2-\hat\xi_\ebop-C w')\rho_+^{-1}x_{\!\scri}^{-2}
  \]
  This is thus negative for $|\hat\xi_\ebop|\leq\frac12$ in $\Omega_\scri$ upon shrinking $\eps_\scri>0$ further if necessary. Furthermore, for small $x_{\!\scri}>0$ and for $\ft_*\geq -1$ we have $t_*\geq-\frac32$ and therefore $\rho_+=(t_*+2)^{-1}\leq 2$. Given $\rho_{+,\rm lo}>0$, we then have $x_{\!\scri}^2 f_\ft^{-1}=\rho_+^{-1}(1+x_{\!\scri}^2)+T x_{\!\scri}^2\leq 2\rho_{+,\rm lo}^{-1}$ for $x_{\!\scri}\leq 1$ and $\rho_+\in[\rho_{+,\rm lo},2]$, and therefore $x_{\!\scri}\leq\sqrt{2\rho_{+,\rm lo}^{-1}f_\ft}$. With the choice of $\rho_{+,\rm lo}$ deferred to later, we fix $\phi_0\in\CIc([0,1))$ to be identically $1$ on $[0,\half]$ and with $\phi_0'\leq 0$ on $[0,1)$, and consider the function
  \begin{equation}
  \label{EqWRchit}
    \chi_\ft := \phi_0\Bigl( 2\eps_\scri^{-1}\sqrt{2\rho_{+,\rm lo}^{-1}f_\ft}\,\Bigr),
  \end{equation}
  which for $\rho_+\in[\rho_{+,\rm lo},2]$ is thus supported in $x_{\!\scri}\leq\half\eps_\scri$, and which is equal to $1$ in a neighborhood of $\scri^+$. Moreover, the function $\chi_\ft\chi_{\scri,\rm out}$ is monotonically increasing along the $\pm\sfH$-flow over $\Omega_\scri$ since $|\hat\xi_\ebop|=|f_{\scri,\rm out}|\leq\half$ on $\supp\chi_{\scri,\rm out}$.
  
  \pfstep{Constructions near $\iota^+$.} With $C_\iota>0$ chosen momentarily, and for appropriate $\eps_\iota>0$ (depending on $C_\iota$), we now work with the coordinates $\rho_\cT=\frac{r}{t_*}$ and $\rho_+=\frac{1}{r}$ near $\iota^+\setminus\scri^+$ and in the domain
  \[
    \Omega_\iota := \Sigma^\pm\cap\{\rho_\cT\leq C_\iota,\ \rho_+\leq\eps_\iota\}.
  \]
  Note that for $\rho_\cT=C_\iota$ and $\rho_+\in[0,\eps_\iota]$, we have $x_{\!\scri}=\sqrt{(t_*+2)/r}=\rho_\cT^{-1}\sqrt{1+2\rho_+\rho_\cT^2}\leq\rho_\cT^{-1}=C_\iota^{-1}$. We shall thus fix $C_\iota\geq(\half\eps_\scri)^{-1}$. In the 3b-momentum coordinates~\eqref{EqWtbCoords}, there then exists a number $\eps_\iota>0$ so that $|\sigma_\tbop|\geq\frac13|\xi_\tbop|+\half|\eta_\tbop|$ on $\Omega_\iota\cap\Sigma^\pm$; this follows from~\eqref{EqWGtb} by arguments similar to the earlier analogous edge-b-computation. Furthermore, $\mp\sigma_\tbop>0$ on $\Omega_\iota\cap\Sigma^\pm$. We introduce $\hat\xi_\tbop=\frac{\xi_\tbop}{\sigma_\tbop}$ and $\hat\eta_\tbop=\frac{\eta_\tbop}{\sigma_\tbop}$ as in~\eqref{EqWRadTCoord}.
  
  For the null-bicharacteristic flow on Minkowski space with metric $-\dd t^2+\dd r^2+r^2\slg$, the momentum $\sigma=\zeta(\pa_t)$ (where $\zeta=\sigma\,\dd t+\xi\,\dd r+\eta$ is a point in the characteristic set) dual to $t$ is conserved and negative in the future characteristic set; moreover, $r\xi=r\zeta(\pa_r)$ is monotonically increasing in the future direction. Therefore, $-r\xi/\sigma=-r\zeta(\pa_r)/\zeta(\pa_t)$ is monotonically increasing and in addition homogeneous of degree $0$ in the fibers of the cotangent bundle. Translated into the 3b-fiber coordinates (and recalling that when passing from the coordinates $t$ and $r$ to $t_*=t-r$ and $r$, the vector fields $\pa_t$, $\pa_r$ become $\pa_{t_*}$, $\pa_r-\pa_{t_*}$), we thus consider in $\Omega_\iota\cap\{\rho_+>0\}$ the function
  \[
    f_\cT := -\rho_+^{-1}\frac{\xi_\tbop-\sigma_\tbop}{\sigma_\tbop} = \rho_+^{-1}(1-\hat\xi_\tbop).
  \]
  Let $w=\rho_+^{\ell_+}$. Then a computation using~\eqref{EqWHam3b} gives
  \begin{equation}
  \label{EqWRfTMono}
    \pm|\sigma_\tbop|^{-1}H_{G_\tbop}f_\cT  = -\sigma_\tbop^{-1}H_{G_\tbop}f_\cT \geq 2\rho_+^{-1}(1-C w) \geq \rho_+^{-1}\quad\text{on}\quad\Omega_\iota,
  \end{equation}
  with the final inequality holding for sufficiently small $\eps_\iota>0$. Note that $f_\cT\to+\infty$ as one approaches $\cR_{\cT,\rm out}^\pm$ (where $\rho_+=0$ and $\hat\xi_\tbop=0$) from $(\cT^+)^\circ$; note also that when $|\hat\xi_\tbop|<\frac12$ and $f_\cT$ is large, then $\rho_+^{-1}$ is large, i.e.\ $\rho_+$ is small. Thus, we shall use $f_\cT$ to measure proximity to $\pa\cR_{\cT,\rm out}^\pm$. We also note that the restriction of the (bounded) function $\frac{f_\cT}{\la f_\cT\ra}$ to the region $\hat\xi_\tbop\leq\half$ is smooth down to $\rho_+=0$.
  
  Working in $\Omega_\iota$ still, we next consider the function
  \[
    f_{\cW,\rm out} := -\frac{\xi_\tbop}{\sigma_\tbop} = -\hat\xi_\tbop.
  \]
  By~\eqref{EqWGtb} and~\eqref{EqWHam3b},
  \[
    \pm|\sigma_\tbop|^{-1}H_{G_\tbop}f_{\cW,\rm out} \geq 2\bigl(|\hat\eta_\tbop|^2 - C w\bigr)\quad\text{on}\quad\Sigma^\pm.
  \]
  For any fixed $\delta\in(0,1]$, this quantity is thus positive on $\Omega_\iota\cap\{|\hat\eta_\tbop|^2\geq\frac{\delta}{16}\}$ for sufficiently small $\eps_\iota$ (depending on $C_\iota$ and $\delta$). We shall use the function $f_{\cW,\rm out}$ to measure proximity to the flow-out $\pa\cW_{\rm out}^\pm$ of $\pa\cR_{\cT,\rm out}^\pm$, see~\eqref{EqWFlowIUnstable}, where $f_{\cW,\rm out}$ vanishes (whereas on the flow-in of $\pa\cR_{\cT,\rm in}^\pm$, one has $f_{\cW,\rm out}=-2$).
  
  Upon shrinking $\eps_\iota$ further, we can moreover ensure that the momentum coordinates of points in $\Sigma^\pm\cap\Omega_\iota$ with $|\hat\xi_\tbop|\leq\frac{\delta}{2}$ satisfy
  \begin{equation}
  \label{EqWRtbBds}
    \hat\xi_\tbop\geq-\frac{\delta}{8},\qquad
    |\hat\eta_\tbop|^2\geq|\hat\xi_\tbop|-\frac{\delta}{8},\qquad
    |\hat\eta_\tbop|^2\leq 4|\hat\xi_\tbop|+\delta\leq 3\delta;
  \end{equation}
  indeed, by~\eqref{EqWGtb} we have $|\hat\eta_\tbop|^2=\hat\xi_\tbop(2-\hat\xi_\tbop)+\cO(w(1+|\hat\eta_\tbop|^2))$ on the characteristic set over $\Omega_\iota$. Define then on $\Omega_\scri$ the function
  \[
    \chi_{\iota,\delta} = \phi\Bigl(\frac{f_{\cW,\rm out}}{\delta}\Bigr)\tilde\psi(\eps_\iota\delta f_\cT)\psi\Bigl(\frac{|\hat\eta_\tbop|^2}{2\delta}\Bigr),
  \]
  where as before $\phi\in\CIc((-\half,\half))$, $\phi|_{[-\frac14,\frac14]}=1$, $\phi'\leq 0$ on $[0,\half]$, and $\psi\in\CIc((-3,3))$, $\psi|_{[2,2]}=1$; furthermore, $\tilde\psi\in\CI(\R)$ vanishes on $(-\infty,4]$ and is equal to $1$ on $[5,\infty)$, with $\tilde\psi'\geq 0$ on $[4,5]$. (Here $\eps_\iota$ is sufficiently small, depending on the choices of $C_\iota$ and $\delta\in(0,1)$.) Furthermore, on $\supp\chi_{\iota,\delta}$, we have $f_\cT\geq 4(\eps_\iota\delta)^{-1}$ and $|\hat\xi_\tbop|<\frac{\delta}{2}<\half$ and therefore $\rho_+=(1-\hat\xi_\tbop)f_\cT^{-1}<\half\eps_\iota\delta$. In addition, on $\Omega_\iota\cap\Sigma^\pm$, and noting that $\psi=1$ on $\supp(\phi\tilde\psi)\cap\Omega_\iota\cap\Sigma^\pm=\emptyset$ by~\eqref{EqWRtbBds}, we have
  \[
    \pm\delta|\sigma_\tbop|^{-1}H_{G_\tbop}\chi_{\iota,\delta} = \phi'\tilde\psi\cdot(\pm|\sigma_\tbop|^{-1}H_{G_\tbop}f_{\cW,\rm out}) + \eps_\iota\phi\tilde\psi'\cdot(\pm|\sigma_\tbop|^{-1}H_{G_\tbop}f_\cT).
  \]
  For the first term, note that on $\supp\phi'$, we have $|\hat\xi_\tbop|\in[\frac{\delta}{4},\frac{\delta}{2}]$, which on $\supp\Sigma^\pm$ forces $\hat\xi_\tbop\in[\frac{\delta}{4},\frac{\delta}{2}]$ and thus $|\hat\eta_\tbop|^2\geq\frac{\delta}{8}$ by~\eqref{EqWRtbBds}, and therefore the Hamiltonian derivative is positive. The second term is nonnegative by~\eqref{EqWRfTMono}.
  
  \pfstep{Gluing.} Let $\cU=x_{\!\scri}^{-1}([\half\eps_\scri,\eps_\scri])$. We claim that for sufficiently small $\delta>0$, the set $\cU\cap\supp\chi_{\iota,\delta}\cap\Sigma^\pm$ is contained in $\cU\cap\{\chi_\ft\chi_{\scri,\rm out}=1\}\cap\Sigma^\pm$. To show this, we now pass back to the edge-b-coordinates $\rho_+=(t_*+2)^{-1}$ and $x_{\!\scri}=\sqrt{(t_*+2)/r}$ (so in particular $r^{-1}=x_{\!\scri}^2\rho_+$). Note then that on $\cU\cap\supp\chi_{\iota,\delta}$ we have
  \[
    r^{-1} < \frac12\eps_\iota\delta,\qquad
    |\hat\xi_\tbop|\leq\frac{\delta}{2},\qquad
    |\hat\eta_\tbop|^2\leq 6\delta.
  \]
  (We recall that $\eps_\iota$ is small, depending on $\delta$, and can thus in particular be chosen to decrease with $\delta$.) In terms of edge-b-momentum coordinates, which are related to 3b-momentum coordinates via~\eqref{EqWeb3bRel}, this implies
  \begin{equation}
  \label{EqWRGlueeb}
  \begin{alignedat}{3}
    \rho_+ &< \frac12\eps_\iota\delta x_{\!\scri}^{-2} &\;\leq\;& 2\eps_\scri^{-2}\eps_\iota\delta, \\
    \Bigl|\frac{\xi_\ebop}{\xi_\ebop-2\zeta_\ebop}\Bigr| &\leq \half\delta x_{\!\scri}^{-2} &\;\leq\;& 2\eps_\scri^{-2}\delta, \\
    \Bigl|\frac{\eta_\ebop}{\xi_\ebop-2\zeta_\ebop}\Bigr| &\leq 6\delta x_{\!\scri}^{-1} &\;\leq\;& 12\eps_\scri^{-1}\delta.
  \end{alignedat}
  \end{equation}
  Note now that on $\cU\cap\Sigma^\pm$, we have $\hat\xi_\ebop(\hat\xi_\ebop-2)+\half x_{\!\scri}^2\hat\xi_\ebop^2+2|\hat\eta_\ebop|^2\leq C\rho_+^{\ell_+}$ (where we recall $\hat\xi_\ebop=\frac{\xi_\ebop}{\zeta_\ebop}$ and $\hat\eta_\ebop=\frac{\eta_\ebop}{\zeta_\ebop}$) by~\eqref{EqWHebG}; therefore, $|\hat\xi_\ebop|$ and thus $|\hat\xi_\ebop-2|$ is bounded by some constant $\hat C<\infty$ which we can take to be independent of $\eps_\scri,\delta$). Since $|\frac{\xi_\ebop}{\xi_\ebop-2\zeta_\ebop}|=|\frac{\hat\xi_\ebop}{\hat\xi_\ebop-2}|$, we therefore conclude from~\eqref{EqWRGlueeb} that
  \begin{equation}
  \label{EqWRSeteb}
    \rho_+,\ |f_{\scri,\rm out}|=|\hat\xi_\ebop|,\ |\hat\eta_\ebop| \leq C_0\delta\quad\text{on}\quad \cU\cap\supp\chi_{\iota,\delta}\cap\Sigma^\pm
  \end{equation}
  with $C_0=\max(2\eps_\scri^{-2},2\hat C\eps_\scri^{-2},12\hat C\eps_\scri^{-1})$, provided $\delta>0$ is small enough. In particular, $\chi_{\scri,\rm out}=1$ on this set; and recalling the definitions of $f_\ft$ and $\chi_\ft$ from~\eqref{EqWRft} and~\eqref{EqWRchit}, we note that $2\eps_\scri^{-1}\sqrt{2\rho_{+,\rm lo}^{-1}f_\ft}\leq 4\eps_\scri^{-1}C_1 x_{\!\scri}\rho_+^{1/2}\rho_{+,\rm lo}^{-1/2}\leq C_2\rho_{+,\rm lo}^{-1/2}\delta^{1/2}<\half$ (where $C_1,C_2$ are independent of $\delta$) for small enough $\delta$---depending on the value of $\rho_{+,\rm lo}$---, and hence $\chi_\ft=1$ on the set~\eqref{EqWRSeteb}.
  
  It remains to choose $\rho_{+,\rm lo}$. First, note that for $x_{\!\scri}\in[\half\eps_\scri,\eps_\scri]$, we have $f_\ft\geq c\rho_+$ for some $c>0$, and therefore on $\supp\chi_\ft$ we obtain $8\eps_\scri^{-2}f_\ft\rho_{+,\rm lo}^{-1}\leq 1$, which implies $\rho_+\leq C_3\rho_{+,\rm lo}$ for some $C_3$. Now on $\cU\cap\supp\chi_{\scri,\rm out}$ (see~\eqref{EqWROrderScrif} and \eqref{EqWRchiscriout}), so $|\hat\xi_\ebop|\leq\half$ and $|\hat\eta_\ebop|^2\leq 3$, we have
  \[
    \pm|\zeta_\ebop|^{-1}H_{G_\ebop^+}x_{\!\scri} \leq -2\bigl(1-(1+\half\eps_\scri^2)\half\bigr)x_{\!\scri} + C\rho_+^{\ell_+}x_{\!\scri} < 0
  \]
  by~\eqref{EqWHameb} provided $\rho_+$ is sufficiently small; we thus fix $\rho_{+,\rm lo}>0$ so that
  \begin{equation}
  \label{EqWRxIMono}
    \pm|\zeta_\ebop|^{-1}H_{G_\ebop^+}x_{\!\scri}<0\quad\text{on}\quad\cU\cap\supp\chi_\ft\chi_{\scri,\rm out}.
  \end{equation}
  
  In order to glue together $\chi_\ft\chi_{\scri,\rm out}$ and $\chi_{\iota,\delta}$, we will use a function
  \[
    \psi_0 \in \CI([0,\infty)),\qquad \psi_0=1\ \text{on}\ [0,\half\eps_\scri],\qquad\psi_0=0\ \text{on}\ [\eps_\scri,\infty),\qquad \psi_0'\leq 0.
  \]
  Choosing $\rho_{+,\rm lo}>0$ small, and then $\delta>0$ small as above, we now set
  \[
    \chi := \psi_0(x_{\!\scri}) \chi_\ft\chi_{\scri,\rm out} + (1-\psi_0(x_{\!\scri})) \chi_{\iota,\delta} \in \CI(\Setb^*M\cap\{\ft_*\geq -1\}).
  \]
  Away from $\supp\dd\psi_0$, we have $\pm\sfH\chi\geq 0$ on $\pa\Sigma^\pm$. On the other hand, we have arranged $\supp\dd\psi_0\cap\supp\chi\cap\pa\Sigma^\pm\subset\cU\cap\supp(\chi_\ft\chi_{\scri,\rm out})\cap\pa\Sigma^\pm$, and on this set $\pm\sfH x_{\!\scri}<0$ by~\eqref{EqWRxIMono} and therefore $\pm\sfH\psi_0\geq 0$. Therefore, $\pm\sfH\chi\geq 0$ on $\pa\Sigma^\pm$.
  
  \pfstep{Conclusion of the proof.} We now construct the variable order function $\sfs$ by setting
  \[
    \sfs := s_{\rm in} - (s_{\rm out}-s_{\rm in})\chi.
  \]
  Since $\chi$ is identically $0$ near $\pa\cR_{\cT,\rm in}^\pm\cup\pa\cR_{\scri,\rm in,+}^\pm$ and identically $1$ near $\pa\cR_{\cT,\rm out}^\pm\cup\pa\cR_{\scri,\rm out}^\pm$ conditions~\eqref{ItWRAbove}--\eqref{ItWRConst} are satisfied; and since $\pm\sfH\chi\geq 0$ on $\pa\Sigma^\pm$, also condition~\eqref{ItWRMono} holds. This completes the construction of a variable order function satisfying conditions~\eqref{ItWRAbove}--\eqref{ItWRMono}. Condition~\eqref{ItWRStat} holds as well provided in the definition~\eqref{EqWRchiscriout} of $\chi_{\scri,\rm out}$ the supports of the cutoff functions $\phi,\psi$ are chosen sufficiently small, likewise for the support of $\phi_0$ in~\eqref{EqWRchit}. This completes the proof.
\end{proof}

\section{Pictorial summary of function spaces}
\label{SF}

We provide a pictorial quick reference for the most important function spaces used in the paper. In the figures below, we label the boundary hypersurfaces of the manifold (or total space, in the case of function spaces with parameter-dependent norms) with the weight of the function space. Below, $X$ is a manifold with boundary and $\rho\in\CI(X)$ is a boundary defining function (such as $X=\ol{\R^n}$ and $\rho=\la r\ra^{-1}$, or $X=[0,1)_R\times\Sph^{n-1}$ and $\rho=R$), and $M$ is the spacetime manifold of Definition~\ref{DefGAMfd}. For the definitions of the spaces, see~\S\ref{SM}.

\begin{figure}[!ht]
\centering
\includegraphics{FigFb}
\caption{\textit{On the left:} the b-Sobolev space $\Hb^{s,\ell}(X)$. \textit{On the right:} the scattering Sobolev space $\Hsc^{s,r}(X)$.}
\label{FigFb}
\end{figure}

\begin{figure}[!ht]
\centering
\includegraphics{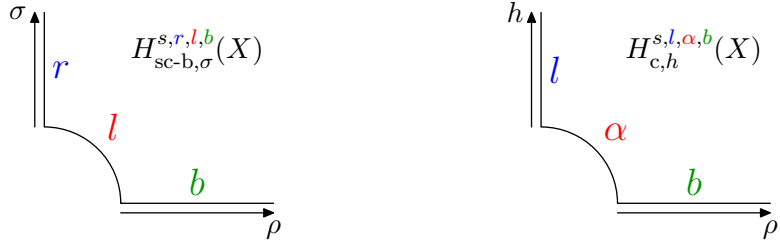}
\caption{\textit{On the left:} the scattering-b-transition Sobolev space $H_{\scbtop,\sigma}^{s,r,l,b}(X)$. \textit{On the right:} the semiclassical cone Sobolev space $H_{\cop,h}^{s,l,\alpha,b}(X)$.}
\label{FigFscbt}
\end{figure}

\begin{figure}[!ht]
\centering
\includegraphics{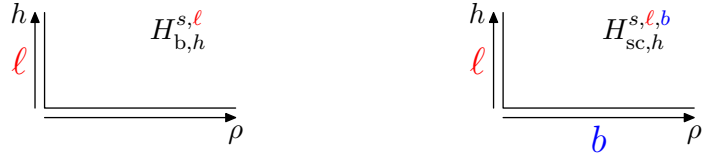}
\caption{\textit{On the left:} the semiclassical b-Sobolev space $H_{\bop,h}^{s,\ell}(X)$ (with the semiclassical order, i.e.\ power of $h$, always $0$ in this paper). \textit{On the right:} the semiclassical scattering Sobolev space $H_{\scop,h}^{s,r,b}(X)$.}
\label{FigFbh}
\end{figure}

\begin{figure}[!ht]
\centering
\includegraphics{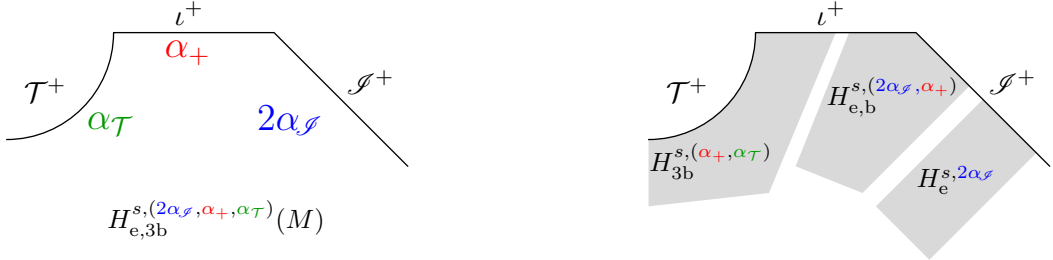}
\caption{\textit{On the left:} the edge-3b-Sobolev space $H_{\eop,\tbop}^{s,(2\alpha_{\!\scri},\alpha_+,\alpha_\cT)}(M)$. \textit{On the right:} Edge, edge-b-, and 3b-Sobolev spaces on $M$. We recall that the factor of $2$ is introduced so that $\alpha_{\!\scri}$ measures decay in terms of powers of $r^{-1}$ (near $(\scri^+)^\circ$), whereas the stated order $2\alpha_{\!\scri}$ refers to powers of $x_{\!\scri}$.}
\label{FigFM}
\end{figure}

\bibliographystyle{alphaurl}
\newcommand{\etalchar}[1]{$^{#1}$}


\end{document}